\documentclass[11pt,a4paper,english]{book}
\usepackage[T1]{fontenc}
\usepackage[latin9]{inputenc}
\usepackage{color}
\usepackage{babel}

\usepackage{varioref}
\usepackage{float}
\usepackage{tipa}
\usepackage{amsmath}
\usepackage{graphicx}
\usepackage{amssymb}
\usepackage{esint}
\usepackage[authoryear]{natbib}
\usepackage[all]{xy}
\usepackage[unicode=true, pdfusetitle,
 bookmarks=false,
 breaklinks=false,pdfborder={0 0 0},backref=false,colorlinks=true]
 {hyperref}

\makeatletter
\newcommand{\lyxaddress}[1]{
\par {\raggedright #1
\vspace{1.4em}
\noindent\par}
}
  \@ifundefined{theoremstyle}{\usepackage{amsthm}}{}
  \theoremstyle{plain}
  \newtheorem{thm}{Theorem}[section]
  \theoremstyle{plain}
  \newtheorem{defn}[thm]{Definition}
  \theoremstyle{definition}
  \newtheorem{rem}[thm]{Remark}
  \theoremstyle{definition}
  \newtheorem*{example*}{Example}
 \theoremstyle{definition}
  \newtheorem{example}[thm]{Example}
  \theoremstyle{plain}
  \newtheorem{cor}[thm]{Corollary}
  \theoremstyle{definition}
  \newtheorem*{rem*}{Remark}
  \theoremstyle{plain}
  \newtheorem{lem}[thm]{Lemma}
  \theoremstyle{plain}
  \newtheorem{assumption}[thm]{Assumption}

\usepackage[all]{xy} 
\usepackage{amscd}
\usepackage{natbib}
\usepackage{babel}
\usepackage{wasysym}
\usepackage{fancyhdr} 
\renewcommand{\chaptermark}[1]{\markboth {\chaptername\ \thechapter. \ #1}{}}

\makeatother 
\makeatletter 
\newcommand{\xyR}[1]{ \makeatletter
\xydef@\xymatrixrowsep@{#1} \makeatother} 
\newcommand{\xyC}[1]{ \makeatletter
\xydef@\xymatrixcolsep@{#1} \makeatother} 
\entrymodifiers={++[ ][F]} 

\numberwithin{equation}{section}
\newcommand{\h}[1]{\hspace{#1 truein}}  
\newcommand{\sss}{\scriptscriptstyle} 

\providecommand{\norm}[1]{\left\lVert#1\right\rVert}

\newcommand{\freccia}{\longrightarrow} 
\newcommand{\xfreccia}[1]{\xrightarrow{\ \ #1\ \ }} 
\newcommand{\xfrecciad}[1]{\xrightarrow{\displaystyle{\ \ \ #1\ \ \ }}}

\newcommand{\eps}{\varepsilon} 
\renewcommand{\phi}{\varphi} 
\newcommand{\lo}{\text{o}} 
\newcommand{\field}[1]{\mathbb{#1}}
\newcommand{\R}{\field{R}}                        
\newcommand{\ER}{{^\bullet\R}}                    
\newcommand{\N}{\field{N}}                        
\newcommand{\hyperR}{{^*\R}}                      
\newcommand{\diff}[1]{\,{\rm d}#1}

\newcommand{\st}[1]{{^\circ #1}} 
\newcommand{\Top}[1]{{\mbox{\Large$\tau$}}_{\sss{#1}}} 
\newcommand{\In}[1]{\in_{_{\sss{#1}}}} 
\newcommand{\InUp}[1]{\in^{\sss{#1}}} 
\newcommand{\Nil}{\mathcal{N}} 
\newcommand{\No}[1]{\mathcal{N}_{#1}} 
\newcommand{\Cc}{\mathcal{C}} 
\newcommand{\D}{\mathcal{D}} 
\newcommand{\F}{\mathcal{F}} 
\newcommand{\Ds}[1]{\mathcal{D}_{\sss #1}} 
\newcommand{\Fs}[1]{\mathcal{F}_{\sss #1}} 
\newcommand{\ext}[1]{{}^\bullet #1} 
\DeclareMathOperator*{\colim}{colim} 
\newcommand{\existsext}[1]{\exists_{\,{}^{^{\bullet}\!}#1}}
\newcommand{\forallext}[1]{\forall_{\,{}^{^{\bullet}\!}#1}}
\newcommand{\corners}[1]{\ulcorner{#1}\urcorner}

\DeclareMathOperator*{\Man}{{\bf Man}\mbox{${}^n$}}
\DeclareMathOperator*{\ManInfty}{{\bf Man}}
\DeclareMathOperator{\Set}{{\bf Set}}
\DeclareMathOperator*{\ORn}{{\bf O}\mbox{$\R^n$}}
\DeclareMathOperator*{\ORInfty}{{\bf O}\mbox{$\R^\infty$}}
\newcommand{\Cn}{\boldsymbol{\mathcal{C}}^n} 
\newcommand{\CInfty}{\boldsymbol{\mathcal{C}}^{\infty}} 
\newcommand{\ECn}{\ext{\Cn}} 
\newcommand{\ECInfty}{\ext{\CInfty}} 
\DeclareMathOperator{\SERn}{{\bf S}\mbox{$\ER^n$}} 
\DeclareMathOperator{\SERInfty}{{\bf S}\mbox{$\ER^{\infty}$}} 

\newcommand{\Lim}[2]{\lim_{t \to 0}\frac{{#1}(t)-{#2}(t)}{t}} 
\newcommand{\Limup}[1]{\lim_{t \to 0}\frac{#1}{t}} 

\newcommand{\qedWithFinalEq}{\phantom{|}\hfill \vrule height 1.5ex width 1.3ex depth-.2ex}
\newcommand{\qedNoNewLine}{\\\phantom{|}\hfill \vrule height 1.5ex width 1.3ex depth-.2ex}
\newcommand{\e}[1]{\text{\h{.1}\  #1 \h{.1}}} 
\newcommand{\ee}[1]{\text{\h{.3}\  #1 \h{.3}}} 
\newcommand{\ptind}{\displaystyle \mathop {\ldots\ldots\,}} 
\newcommand{\pti}{:\;\;\;} 
\newcommand{\then}{\quad \Longrightarrow \quad}
\newcommand{\DIff}{ \quad\;\; :\!\iff \quad } 

\makeatother

\begin{document}
\date{}

\title{\textbf{\Huge Fermat Reals}\textbf{\Large }\\
\textbf{\large Nilpotent Infinitesimals and Infinite Dimensional
Spaces}}

\author{Paolo Giordano}

\maketitle

\lyxaddress{\noindent {\scriptsize Modeling and Applications of Complex Systems
Laboratory (MACS-Lab), Università della Svizzera italiana, via Canavée,
CH-6850, Mendrisio, Switzerland. Email: }\texttt{\scriptsize paolo.giordano@usi.ch}}

\noindent \begin{center}
\textbf{Abstract}
\par\end{center}

\noindent \textbf{\footnotesize F.}{\footnotesize{} Good morning Hermann,
I would like to talk with you about infinitesimals.}{\footnotesize \par}

\noindent \textbf{\footnotesize G.}{\footnotesize{} Tell me Pierre.}{\footnotesize \par}

\noindent \textbf{\footnotesize F.}{\footnotesize{} I'm fed up of all
these slanders about my attitude to be non rigorous, so I've started
to study nonstandard analysis (NSA) and synthetic differential geometry
(SDG).}{\footnotesize \par}

\noindent \textbf{\footnotesize G.}{\footnotesize{} Yes, I've read
something\ldots{}}{\footnotesize \par}

\noindent \textbf{\footnotesize F.}{\footnotesize{} Ok, no problem
about their rigour. But, when I've seen that the sine of an infinite
in NSA is infinitely near to a real number I was astonished: what
is the intuitive meaning of this number, if any? Then, I've seen that
to work in SDG I must learn to work in intuitionistic logic\ldots{}
You know, I love margins of books, and I don't want to loose too much
time, I have many things to do\ldots{}}{\footnotesize \par}

\noindent \textbf{\footnotesize G.}{\footnotesize{} In SDG they also
say that every infinitesimal is at the same time positive and negative,
what is the meaning of all these? And why does the square of a first
order infinitesimal equal zero, whereas the product of two first order
infinitesimals is not necessarily zero? And do you know that from
any single infinitesimal in NSA is possible to construct a non measurable
set? Without using the axiom of choice!}{\footnotesize \par}

\noindent \textbf{\footnotesize F.}{\footnotesize{} Yes, I know, I
know\ldots{} Ok, listen: why cannot we start from standard functions
$x:\R\freccia\R$ and use\ldots{}}{\footnotesize \par}

{\footnotesize \medskip{}
}{\footnotesize \par}

{\footnotesize This work is the ideal continuation of this dialogue:
a theory of actual infinitesimals that do not need a background of
formal logic to be understood, with a clear intuitive meaning and
with non trivial applications to differential geometry of both finite
and infinite dimensional spaces.}{\footnotesize \par}

\begin{center}
{\Large \newpage{}}
\par\end{center}{\Large \par}

\thispagestyle{empty}

\newpage{}

\pagenumbering{roman}

{\Large \tableofcontents{}}{\Large \par}

\newpage{}

\listoffigures

\newpage{}

\thispagestyle{empty}
\begin{verse}
\noindent \emph{If we do not believe in the existence of God, then,
from G$\text{\textipa{ö}}$del's ontological theorem it follows that
an absolute, necessary and not only possible moral system cannot exists.
But I believe that the human kind can achieve this type of moral system,
so I have to believe in God.}

\noindent \emph{If you further assume, suitably formalized with a
rigorous mathematical language, that any good thing has in God its
first cause and that mathematics is a good thing, then you cannot
believe in genius anymore. The simple consequences for the everyday
work of a mathematician and, more in general, for scientific collaboration
are left to the reader.}
\end{verse}
\noindent \begin{flushright}
P. Castelluccia
\par\end{flushright}

\newpage{}

\pagestyle{empty}

\noindent \begin{center}
\textbf{\huge Abstract and structure}
\par\end{center}{\huge \par}

The main aim of the present work is to start a new theory of actual
infinitesimals, called \emph{theory of Fermat reals}. After the work
of A. Robinson on nonstandard analysis (NSA), several theories of
infinitesimals have been developed: synthetic differential geometry,
surreal numbers, Levi-Civita field, Weil functors, to cite only some
of the most studied. We will discuss in details of these theories
and their characteristics, first of all comparing them with our Fermat
reals. One of the most important differences is the philosophical
thread that guided us during all the development of the present work:
we tried to construct a theory with a strong intuitive interpretation
and with non trivial applications to the infinite-dimensional differential
geometry of spaces of mappings. This driving thread tried to develop
a good dialectic between formal properties, proved in the theory,
and their informal interpretations. The dialectic has to be, as far
as possible, in both directions: theorems proved in the theory should
have a clear and useful intuitive interpretation and, on the other
hand, the intuition corresponding to the theory has to be able to
suggest true sentences, i.e. conjectures or sketch of proofs that
can then be converted into rigorous proofs. Almost all the present
theories of actual infinitesimals are either based on formal approaches,
or are not useful in differential geometry. As a meaningful example,
we can say that the Fermat reals can be represented geometrically
(i.e. they can be drawn) respecting the total order relation.

The theory of Fermat reals takes a strong inspiration from synthetic
differential geometry (SDG), a theory of infinitesimals grounded in
Topos theory and incompatible with classical logic. SDG, also called
smooth infinitesimal analysis, originates from the ideas of \citet{La1}
and has been greatly developed by several categorists. The result
is a powerful theory able to develop both finite and infinite dimensional
differential geometry with a formalism that takes great advantage
of the use of infinitesimals. This theory is however incompatible
with classical logic and one is forced to work in intuitionistic logic
and to construct models of SDG using very elaborated topoi. The theory
of Fermat reals is sometimes formally very similar to SDG and indeed,
several proofs are simply a reformulation in our theory of the corresponding
proofs in SDG. However, our theory of Fermat reals is fully compatible
with classical logic. We can thus describe our work as a way to bypass
an impossibility theorem of SDG, i.e. a way considered as impossible
by several researchers. The differences between the two theories are
due to our constraint to have always a good intuitive interpretation,
whereas SDG develops a more formal approach to infinitesimals.

Generally speaking, we have constructed a theory of infinitesimals
which does not need a background of logic to be understood. On the
contrary, nonstandard analysis and SDG need this non trivial background,
and this is a great barrier for potential users like physicists or
engineers or even several mathematicians. This is a goal strongly
searched in NSA, so as to facilitate the diffusion of the theory.

Many parts of our construction are completely constructive and this
result, also considered by several researchers in NSA, opens good
possibilities for a computer implementation of our Fermat reals, with
interesting potential applications in automatic proof theory or in
automatic differentiation theory.

Our infinitesimals $h$, like in SDG, are \emph{nilpotent} so that
we have $h\ne0$, but $h$ is {}``so small'' that for some power
$n\in\N_{>1}$ we have $h^{n}=0$. This permits to obtain an equality
between a function and its tangent straight line in a first order
infinitesimal neighborhood, i.e.\begin{equation}
f(x+h)=f(x)+h\cdot f'(x),\label{eq:exampleTaylorFirst}\end{equation}
where $h^{2}=0$. More generally, we will prove infinitesimal Taylor's
formulas without any rest, so that every smooth functions, in our
framework, is equal to a $k$-th order polynomial in every $k$-th
order infinitesimal neighborhood.

\medskip{}

The second part of the work is devoted to the development of a theory
of smooth infinite dimensional spaces, first of all thinking applications
in differential geometry. Our approach is based on a generalization
of the notion of diffeology (see e.g. \citet{Igl}). This permits
to obtain a cartesian closed complete and cocomplete category in which
the category of smooth manifolds is embedded. Using the above mentioned
generalization we can obtain a category containing the extension of
all smooth manifolds using our new infinitesimal points. We have hence
the category $\CInfty$ of diffeological spaces, which contains all
the smooth manifolds, and a functor $\ext{(-)}:\CInfty\freccia\ECInfty$,
called Fermat functor, which extends every space $X\in\CInfty$ adding
infinitesimal points. E.g. the ring of Fermat reals is $\ER:=\ext{(\R)}$.
The above mentioned categorical properties of these categories permits
to say that we can construct infinite products, spaces of mapping,
infinite sums, quotient spaces and we also have that mappings like
composition, insertion, evaluation, and practically all the interesting
set-theoretical operations are always smooth. We have hence a flexible
framework where infinitesimal methods are also available.

Moreover, the Fermat functor possesses very good properties: it preserves
products of manifolds and intersections, unions, inclusions, counter-images
of open sets, to cite some of them. We will study in general this
preservation properties, discovering some relationships between the
Fermat functor and intuitionistic logic.\medskip{}

In the third part of the work we will present the basis for the whole
development of the differential and integral calculus both for smooth
functions defined on open sets of Fermat reals and on infinitesimal
domains. We also give some first results of differential geometry
using infinitesimal methods, always considering the case of the space
of all the smooth mappings between two manifolds. A very general proof
of the Euler-Lagrange equations is also given, with Lagrangians defined
on spaces of mappings of the general form\[
\ECInfty(\ext{M_{s}},\ECInfty(\ext{M_{s-1}},\cdots,\ECInfty(\ext{M_{2}},\ext{\R})\cdots),\]
where $M_{i}$ are manifolds. The space $\ext{M_{i}}\in\ECInfty$
is the application of the Fermat functor to the manifold $M_{i}$,
so that it can be thought as the manifold with the adding of our new
infinitesimal points. In this section we also give a sketch of some
ideas for a further development of the present work.\medskip{}

The fourth part of the work is composed of an appendix that fixes
common notations for the concepts of category theory that we have
used, and of a detailed study about the comparison of our theory and
other theories of infinitesimals.\bigskip{}

The detailed structure of our work is as follows. After a motivational
Chapter \ref{cha:IntroductionAndGeneralProblem} where we will also
give an explanation for the name \emph{Fermat reals}, in Chapter \ref{cha:Definition-and-algebraic-prop-of-Fermat-reals}
we will define the ring $\ER$ of Fermat reals and the ideals $D_{k}$
of $k$-th order infinitesimals. Having a ring which contains nilpotent
elements, one of the most difficult algebraic problem is the dealing
with products of powers of these nilpotent numbers. In this chapter
we will also prove several effective results that permits to solve
these powers (i.e. to decide whether they are zero or not) in an algorithmic
way.

The derivative $f'(x)$ in a Taylor's formula like \eqref{eq:exampleTaylorFirst}
is determined only up to second order infinitesimals. In Chapter \ref{cha:equalityUpTo_k-thOrder}
we will deeply study these equality up to $k$-th order infinitesimals,
the corresponding cancellation law and its application to Taylor's
formulas.

In Chapter \ref{cha:orderRelation} we will define the total order
relation. We will show that, generally speaking, the order relation
can be total only if the derivative in \eqref{eq:exampleTaylorFirst}
is \emph{not} uniquely determined. In this chapter we will also prove
that the Fermat reals are in bijective correspondence with suitable
curves of the plane $\R^{2}$, i.e. the geometrical representation
of $\ER$.

Chapter \ref{cha:ApproachesToDiffGeomOfInfDimSpaces} starts the second
part of the work, devoted to our approach to infinite dimensional
spaces. In this chapter we review the most studied approaches to infinite
dimensional spaces used in differential geometry: Banach manifolds
and locally convex vector spaces, the convenient vector spaces settings,
diffeological spaces and SDG, presenting some of their positive features
and some possible deficiencies.

In Chapter \ref{cha:theCartesianClosure} we present our generalization
of the notion of diffeological space, which permits to define in the
same framework both the categories $\CInfty$ and $\ECInfty$, respectively
domain and codomain of the Fermat functor. We called this generalization
\emph{the cartesian closure of a given category of figures}.

In Chapter \ref{cha:TheCategoryCn} the cartesian closure is applied
to the category of open sets in spaces of the form $\R^{n}$ and smooth
mappings, obtaining the category $\CInfty$ of diffeological spaces.
We review the embedding of smooth finite dimensional manifolds and
give several examples: infinite dimensional manifolds modeled on convenient
vector spaces (which include manifolds modeled on Banach spaces),
integro-differential operators, set-theoretical operations like compositions
and evaluations, and we prove that the space of all the diffeomorphisms
between two manifolds is a Lie group.

In Chapter \ref{cha:Extending-smooth-spaces} we generalize the construction
of the Fermat ring $\ER$ to any smooth diffeological space $X\in\CInfty$
and we define the category $\ECInfty$ of smooth Fermat spaces, which
includes all the spaces of the form $\ext{X}$.

Chapter \ref{cha:TheFermatFunctor} starts the study of the Fermat
functor $\ext{(-)}:\CInfty\freccia\ECInfty$ that extends every smooth
space $X\in\CInfty$ by adding infinitesimal points. We prove that
this functor preserves products of manifolds and we prove that manifolds
are also embedded in the category $\ECInfty$. In this chapter we
also prove that a standard part functor, right adjoint of the Fermat
functor, does not exists. This correspond to analogous results dealing
with the standard part map in constructive NSA.

We then study, in Chapter \ref{cha:logicalPropertiesOfTheFermatFunctor},
the logical properties of the Fermat functor, i.e. all the logical
operations which are preserved by it. We will see that, even if the
theory of Fermat reals is fully compatible with classical logic, the
best properties of this functor are present in the case of an intuitionistic
interpretation of these logical operations, confirming the good dialectic
between smooth differential geometry and intuitionistic logic.

The third part of the work starts with the study, in Chapter \ref{cha:CalculusOnOpenDomains},
of the development of the basis for the differential and integral
calculus of smooth functions $f:\ext{U}\freccia\ER^{d}$ defined on
an open set $\ext{U}\subseteq\ER^{n}$. These functions generalize
the standard smooth functions and can be expressed, locally, as the
extension of standard smooth functions $\ext{\alpha}(p,-)$ with a
fixed parameter $p\in\ER^{\sf p}$. The differential calculus is based
on the analogous, in SDG, of the Fermat-Reyes property, and formalizes
perfectly the informal methods used originally by P. de Fermat. In
this chapter we also prove the inverse function theorem in $\ER$
and the existence of primitives, which represent a non trivial problem
in non-Archimedean fields.

Due to the connections between total order and nilpotent infinitesimals,
the differential calculus for function defined on infinitesimal sets,
like $D_{k}=\left\{ h\in\ER\,|\, h^{k+1}=0\right\} $, must be developed
using the properties of the equality up to $k$-th order infinitesimals.
This is done in Chapter \ref{cha:CalculusOnInfinitesimalDomains}.

The purpose of Chapter \ref{cha:Infinitesimal-differential-geometry}
is to show the possibilities of the theory of Fermat reals for differential
geometry, in particular for spaces of mappings. We essentially develop
only tangency theory and the existence of integral curves using infinitesimal
methods. We devoted a particular attention to always include in our
results spaces of the form $\ext{M}^{\ext{N}}\in\ECInfty$ for $M$
and $N$ manifolds. In this chapter we also prove the above mentioned
general version of the Euler-Lagrange equations.

In the final Chapter \ref{cha:FurtherDevelopments} we sketch the
ideas of some possible further developments of our work.

\newpage{}

\pagenumbering{arabic}

\pagestyle{fancy}

\part{Algebraic and order properties of Fermat reals}

\chapter{\label{cha:IntroductionAndGeneralProblem}Introduction and general
problem}

Frequently in work by physicists it is possible to find informal calculations
like \begin{equation}
\frac{1}{\sqrt{1-{\displaystyle \frac{v^{2}}{c^{2}}}}}=1+\frac{v^{2}}{2c^{2}}\qquad\qquad\sqrt{1-h_{44}(x)}=1-\frac{1}{2}h_{44}(x)\label{eq:EinsteinInfinitesimal}\end{equation}
 with explicit use of infinitesimals $v/c\ll1$ or $h_{44}(x)\ll1$
such that e.g. $h_{44}(x)^{2}=0$. For example \cite{Ein} (pag. 14)
wrote the formula (using the equality sign and not the approximate
equality sign $\simeq$) \begin{equation}
f(x,t+\tau)=f(x,t)+\tau\cdot\frac{\partial f}{\partial t}(x,t)\label{eq:EinsteinDerivationFormula}\end{equation}
 justifying it with the words {}``\emph{since $\tau$ is very small}'';
the formulas \eqref{eq:EinsteinInfinitesimal} are a particular case
of the general \eqref{eq:EinsteinDerivationFormula}. \cite{Dir}
wrote an analogous equality studying the Newtonian approximation in
general relativity.

Using this type of infinitesimals we can write an \emph{equality},
in some infinitesimal neighborhood, between a smooth function and
its tangent straight line, or, in other words, a Taylor's formula
without remainder. Informal methods based on actual infinitesimals
are sometimes used in differential geometry too. Some classical examples
are the following: a tangent vector is an infinitesimal arc of curve
traced on the manifold and the sum of tangent vectors is made using
infinitesimal parallelograms; tangent vectors to the tangent bundle
are infinitesimal squares on the manifold; a vector field is sometimes
intuitively treated as an {}``infinitesimal transformation'' of
the space into itself and the Lie brackets of two vector fields as
the commutator of the corresponding infinitesimal transformations.

There are obviously many possibilities to formalize this kind of intuitive
reasonings, obtaining a more or less good dialectic between informal
and formal thinking, and indeed there are several theories of actual
infinitesimals (from now on, for simplicity, we will say {}``infinitesimals''
instead of {}``actual infinitesimals'' as opposed to {}``potential
infinitesimals''). Starting from these theories we can see that we
can distinguish between two type of definitions of infinitesimals:
in the first one we have at least a ring $R$ containing the real
field $\R$ and infinitesimals are elements $\eps\in R$ such that
$-r<\eps<r$ for every positive standard real $r\in\R_{>0}$. The
second type of infinitesimal is defined using some algebraic property
of nilpotency, i.e. $\eps^{n}=0$ for some natural number $n\in\N$.
For some ring $R$ these definitions can coincide, but anyway they
lead, of course, only to the trivial infinitesimal $\eps=0$ if $R=\R$.

However these definitions of infinitesimals correspond to theories
which are completely different in nature and underlying ideas. Indeed
these theories can be seen in a more interesting way to belong to
two different classes. In the first one we can put theories that need
a certain amount of non trivial results of mathematical logic, whereas
in the second one we have attempts to define sufficiently strong theories
of infinitesimals without the use of non trivial results of mathematical
logic. In the first class we have Non-Standard Analysis (NSA) and
Synthetic differential geometry (SDG, also called Smooth Infinitesimal
Analysis), in the second one we have, e.g., Weil functors, Levi-Civita
fields, surreal numbers, geometries over rings containing infinitesimals
(see Appendix \ref{app:theoriesOfInfinitesimals} for an introduction
to several approaches to infinitesimals, together with a first comparison
with our approach, and for references). More precisely we can say
that to work in NSA and SDG one needs a formal control deeply stronger
than the one used in {}``standard mathematics''. In NSA one needs
this control to apply the transfer theorem and in SDG one has to be
sufficiently formal to be sure that the proofs can be seen as belonging
to intuitionistic logic. Indeed to use NSA one has to be able to formally
write the sentences one needs to transfer. Whereas SDG does not admit
models in classical logic, but in intuitionistic logic only, and hence
we have to be sure that in our proofs there is no use of the law of
the excluded middle, or e.g. of the classical part of De Morgan's
law or of some form of the axiom of choice or of the implication of
double negation toward affirmation and any other logical principle
which is not valid in intuitionistic logic. Physicists, engineers,
but also the greatest part of mathematicians are not used to have
this strong formal control in their work, and it is for this reason
that there are attempts to present both NSA and SDG reducing as much
as possible the necessary formal control, even if at some level this
is technically impossible (see e.g. \cite{Hen}, and \cite{Be-DN1, Be-DN2}
for NSA; \cite{Bel} and \cite{Lav} for SDG, where using an axiomatic
approach the authors try to postpone the very difficult construction
of an intuitionistic model of a whole set theory using Topos).

On the other hand NSA is essentially the only theory of infinitesimals
with a discrete diffusion and a sufficiently great community of working
mathematicians and published results in several areas of mathematics
and its applications, see e.g. \cite{Al-Fe-Ho-Li}. SDG is the only
theory of infinitesimals with non trivial, new and published results
in differential geometry concerning infinite dimensional spaces like
the space of all the diffeomorphisms of a generic (e.g. non compact)
smooth manifold. In NSA we have only few results concerning differential
geometry (we cite \cite{Sch} and \cite{Ham}, and references therein,
where NSA methods are used in problems of differential geometry).
Other theories of infinitesimals have not, at least up to now, the
same formal strength of NSA or SDG or the same potentiality to be
applied in several different areas of mathematics.

One of the aim of the present work is to find a theory of infinitesimals
within {}``standard mathematics'' (in the precise sense explained
above of a formal control more {}``standard'' and not so strong
as the one needed e.g. in NSA or SDG) with results comparable with
those of SDG, without forcing the reader to learn a strong formal
control of the mathematics he is doing. Because it has to be considered
inside {}``standard mathematics'', our theory of infinitesimals
must be compatible with classical logic. Let us note that this is
not incompatible with the possibility\emph{ }to obtain some results
that need a strong formal control (like, e.g., a transfer theorem),
because they represent a good potential instrument for the reader
that likes such a strong formal control, but they do not force, concretely,
all the readers to have such a formal aptitude. For these reasons,
we think that it is wrong to frame the present work as in opposition
to NSA or SDG.

Concretely, the idea of the present work is to by-pass the impossibility
theorem about the incompatibility of SDG with classical logic that
forces SDG to find models within intuitionistic logic. This by-pass
has to be made, as much as possible, keeping the same properties and
final results. We think that the obtained result is meaningful not
only for differential geometry, but also for other fields, like the
calculus of variations, and we will give a first sketch of results
in this direction.

Another point of view about a powerful theory like NSA is that, in
spite of the fact that frequently it is presented using opposed motivations,
it lacks the intuitive interpretation of what the powerful formalism
permits to do. E.g. what is the intuitive meaning and usefulness of
$\st{\sin(I)}\in\R$, i.e. the standard part of the sine of an infinite
number $I\in\hyperR$? This and the above-mentioned {}``strong formal
control'' needed to work in NSA, together with very strong but scientifically
unjustified cultural reasons, may be some motivations for the not
so high success of the spreading of NSA in mathematics, and consequently
in its didactics.

Analogously in SDG from the intuitive, classical, point of view, it
is a little strange that we cannot exhibit {}``examples'' of infinitesimals
(indeed in SDG it is only possible to prove that $\neg\neg\exists\, d\in D$,
where $D=\{h\in R\,|\, h^{2}=0\}$ is the set of first order infinitesimals).
Because of this, e.g., we cannot construct a physical theory containing
a fixed infinitesimal parameter; another example of a counter intuitive
property is that any $d\in D$ is, at the same time, positive $d>0$
and negative $d<0$. Similar counter intuitive properties can be found
in other theories of infinitesimals that use ideals of rings of polynomials
as a formal scheme to construct particular type of infinitesimals.
Among these theories we can cite {}``Weil functors'' (see \cite{Ko-Mi-Sl}
and \cite{Kr-Mi2} and Appendix \ref{app:theoriesOfInfinitesimals}
of the present work for other references) and {}``differential geometry
over general base fields and rings'' (see \cite{Ber} and Appendix
\ref{app:theoriesOfInfinitesimals}). The final conclusion after the
establishment of this type of counter intuitive examples (even if,
of course, in these theories there are also several intuitively clear
examples and concepts), is that if one wants to work in these types
of frameworks, sometimes one has to follow a completely formal point
of view, loosing the dialectic with the corresponding intuitive meaning.

Another aim of the present work is to construct a new theory of infinitesimals
preserving always a very good dialectic between formal properties
and intuitive interpretation. A first hint to show this positive feature
of our construction is that our is the first theory, as far as we
know, where it is possible to represent geometrically its new type
of numbers%
\footnote{I.e. it is possible to establish a bijective correspondence between
suitable lines of the plane and the numbers belonging to a given infinitesimal
neighborhood.%
}, and it is undeniable that to be able to represent standard real
numbers by a straight line inspired, and it still inspires, several
mathematicians.

More technically we want to show that it is possible to extend the
real field adding nilpotent infinitesimals, arriving at an enlarged
real line $\ER$, by means of a very simple construction completely
inside {}``standard mathematics''. Indeed to define the extension
$\ER\supset\R$ we shall use elementary analysis only. To avoid misunderstandings
is it important to clarify that present work's purpose is not to give
an alternative foundation of differential and integral calculus (like
NSA), but to obtain a theory of nilpotent infinitesimals and to use
it for the foundation of a smooth ($\Cc^{\infty}$) differential geometry,
in particular in the case of infinite dimensional spaces, like the
space of all the smooth functions $\ManInfty(M;N)$ between two generic
manifolds (e.g. without compactness hypothesis on the domain $M$).
This focus on the foundation of differential geometry only, without
including the whole calculus, is more typical of SDG, Weil functors
and geometries over generic rings.

The usefulness of the extension $\ER\supset\R$ can be glimpsed by
saying e.g. that using $\ER$ it is possible to write in a completely
rigorous way that a smooth function is equal to its tangent straight
line in a first order neighborhood; it is possible to use infinitesimal
Taylor's formulas without remainder; to define a tangent vector as
an infinitesimal curve and sum them using infinitesimal parallelograms;
to see a vector field as an infinitesimal transformation, in general
to formalize these and many other non-rigorous methods used in physics
and geometry. This is important both for didactic reasons and because
it was by means of these methods that mathematicians like S. Lie and
E. Cartan were originally led to construct important concepts of differential
geometry.

We can use the infinitesimals of $\ER$ not only as a good language
to reformulate well-known results, but also as a very useful tool
to construct, in a simple and meaningful way, a differential geometry
in classical infinite-dimensional objects like \textbf{Man}$(M,N)$
the space of all the $\Cc^{\infty}$ mapping between two manifolds
$M$, $N$. Here with {}``simple and meaningful'' we mean the idea
to work directly on the geometric object in an intrinsic way without
being forced to use charts, but using infinitesimal points (see \cite{Lav}).
Some important examples of spaces of mappings used in applications
are the space of configurations of a continuum body, groups of diffeomorphisms
used in hydrodynamics as well as in magnetohydrodynamics, electromagnetism,
plasma dynamics, and paths spaces for calculus of variations (see
\cite{Kr-Mi, Ab-Ma-Ra, AJPS, Al-Fe-Ho-Li, Alb1} and references therein).
Interesting applications in classical field theories can also be found
in \cite{Ab-Ma}.

\section{Motivations for the name {}``Fermat reals''}

It is well known that historically two possible reductionist constructions
of the real field starting from the rationals have been made. The
first one is Dedekind's order completion using sections of rationals,
the second one is Cauchy's metric space completion. Of course there
are no historical reason to attribute our extension $\ER\supset\R$
of the real field, to be described below, to Fermat, but there are
strong motivations to say that, probably, he would have liked the
underlying spirit and some properties of our theory. For example: 
\begin{enumerate}
\item we will see that a formalization of Fermat's infinitesimal method
to derive functions is provable in our theory. We recall that Fermat's
idea was, roughly speaking and not on the basis of an accurate historical
analysis which goes beyond the scope of the present work (see e.g.
\cite{Bo-Fr-Ri, Edw, Eve}), to suppose first $h\ne0$, to construct
the incremental ratio\[
\frac{f(x+h)-f(x)}{h}\]
 and, after suitable simplifications (sometimes using infinitesimal
properties), to take in the final result $h=0$. 
\item Fermat's method to find the maximum or minimum of a given function
$f(x)$ at $x=a$ was to take $e$ to be extremely small so that the
value of $f(x+h)$ was approximately equal to that of $f(x)$. In
modern, algebraic language, it can be said that $f(x+h)=f(x)$ only
if $h^{2}=0$, that is if $e$ is a first order infinitesimal. Fermat
was aware that this is not a {}``true'' equality but some kind of
approximation (see e.g. \cite{Bo-Fr-Ri, Edw, Eve}). We will follow
a similar idea to define $\ER$ introducing a suitable equivalence
relation to represent this equality. 
\item Fermat has been described by \cite{EBell} as {}``the king of amateurs''
of mathematics, and hence we can suppose that in its mathematical
work the informal/intuitive part was stronger with respect to the
formal one. For this reason we can think that he would have liked
our idea to obtain a theory of infinitesimals preserving always the
intuitive meaning and without forcing the working mathematician to
be too much formal. 
\end{enumerate}
For these reason we chose the name {}``Fermat reals'' for our ring
$\ER$ (note: without the possessive case, to underline that we are
not attributing our construction of $\ER$ to Fermat).

We already mentioned that the use of nilpotent infinitesimals in the
ring $\ER$ permits to develop many concepts of differential geometry
in an intrinsic way without being forced to use coordinates, as we
shall see in some examples in the course of the present work. In this
way the use of charts becomes specific of stated areas, e.g. where
one strictly needs some solution in a finite neighborhood and not
in an infinitesimal one only (e.g. this is the case for the inverse
function theorem). We can call \emph{infinitesimal differential geometry}
this kind of intrinsic geometry based on the ring $\ER$ (and on extensions
of manifolds $\ext{M}$ and also on more generic object like the exponential
objects $\ext{M}^{\ext{N}}$, see the second part of the present work).

\chapter[Definition of Fermat reals]{\label{cha:Definition-and-algebraic-prop-of-Fermat-reals}Definition
and algebraic properties of Fermat reals}

\section{The basic idea}

We start from the idea that a smooth ($\Cc^{\infty}$) function $f:\ER\freccia\ER$
is actually \emph{equal} to its tangent straight line in the first
order neighborhood e.g. of the point $x=0$, that is

\begin{equation}
\forall h\in D\pti f(h)=f(0)+h\cdot f'(0)\label{eq:FunctionEqualTangent}\end{equation}
 where $D$ is the subset of $\ER$ which defines the above-mentioned
neighborhood of $x=0$. The equality \eqref{eq:FunctionEqualTangent}
can be seen as a first-order Taylor's's formula without remainder
because intuitively we think that $h^{2}=0$ for any $h\in D$ (indeed
the property $h^{2}=0$ defines the first order neighborhood of $x=0$
in $\ER$). These almost trivial considerations lead us to understand
many things: $\ER$ must necessarily be a ring and not a field because
in a field the equation $h^{2}=0$ implies $h=0$; moreover we will
surely have some limitation in the extension of some function from
$\R$ to $\ER$, e.g. the square root, because using this function
with the usual properties, once again the equation $h^{2}=0$ implies
$|h|=0$. On the other hand, we are also led to ask whether \eqref{eq:FunctionEqualTangent}
uniquely determines the derivative $f^{\prime}(0)$: because, even
if it is true that we cannot simplify by $h$, we know that the polynomial
coefficients of a Taylor's's formula are unique in classical analysis.
In fact we will prove that

\begin{equation}
\exists!\, m\in\R\,\,\forall h\in D\pti f(h)=f(0)+h\cdot m\label{eq:IdeaDF}\end{equation}
 that is the slope of the tangent is uniquely determined in case it
is an ordinary real number. We will call formulas like \eqref{eq:IdeaDF}
\emph{derivation formula}s.

If we try to construct a model for \eqref{eq:IdeaDF} a natural idea
is to think our new numbers in $\ER$ as equivalence classes $[h]$
of usual functions $h:\R\freccia\R$. In this way we may hope both
to include the real field using classes generated by constant functions,
and that the class generated by $h(t)=t$ could be a first order infinitesimal
number. To understand how to define this equivalence relation we have
to think at \eqref{eq:FunctionEqualTangent} in the following sense:

\begin{equation}
f(h(t))\sim f(0)+h(t)\cdot f^{\prime}(0),\label{eq:IdeaFunction}\end{equation}
where the idea is that we are going to define $\sim$. If we think
$h(t)$ {}``sufficiently similar to $t$'', we can define $\sim$
so that (\ref{eq:IdeaFunction}) is equivalent to \[
\Limup{f(h(t))-f(0)-h(t)\cdot f^{\prime}(0)}=0,\]
 that is

\begin{equation}
x\sim y\DIff\Lim{x}{y}=0.\label{eq:IdeaRelEq}\end{equation}
 In this way \eqref{eq:IdeaFunction} is very near to the definition
of differentiability for $f$ at 0. 

\noindent It is important to note that, because of de L'Hôpital's
theorem we have the isomorphism \[
\Cc^{1}(\R,\R)/\!\sim\,\,\,\simeq\,\,\R[x]/(x),\]
 the left hand side is (isomorphic to) the usual tangent bundle of
$\R$ and thus we obtain nothing new. It is not easy to understand
what set of functions we have to choose for $x$, $y$ in (\ref{eq:IdeaRelEq})
so as to obtain a non trivial structure. The first idea is to take
continuous functions at $t=0$, instead of more regular ones like
$\mathcal{C}^{1}$-functions, so that e.g. $h_{k}(t)=|t|^{1/k}$ becomes
a $k$-th order nilpotent infinitesimal ($h^{k+1}\sim0$); indeed
for almost all the results presented in this article, continuous functions
at $t=0$ work well. However, only in proving the non-trivial property

\noindent \begin{equation}
\left(\forall x\in\ER\pti x\cdot f(x)=0\right)\then\forall x\in\ER\pti f(x)=0\label{eq:UniquenessIncrementalRatio}\end{equation}
 (here $f:\ER\freccia\ER$ is a smooth function, in a sense we shall
make precise afterwards), we will see that it does not suffice to
take continuous functions at $t=0$. Property (\ref{eq:UniquenessIncrementalRatio})
is useful to prove the uniqueness of smooth incremental ratios, hence
to define the derivative $f^{\prime}:\ER\freccia\ER$ of a smooth
function $f:\ER\freccia\ER$ which, generally speaking, is not the
extension to $\ER$ of an ordinary function defined on $\R$ (like,
e.g., the function $t\mapsto\sin(h\cdot t)$, where $h\in\ER\setminus\R$,
which is used in elementary physics to describe the small oscillations
of the pendulum ). To prove (\ref{eq:UniquenessIncrementalRatio})
the following functions turned out to be very useful: 
\begin{defn}
\noindent \label{def:NilpotentFunctions}If $x:\R_{\ge0}\freccia\R$,
then we say that $x$ \emph{is nilpotent} iff $|x(t)-x(0)|^{k}=o(t)$
as $t\to0^{+}$, for some $k\in\N$. $\Nil$ will denote the set of
all the nilpotent functions.
\end{defn}
\noindent In the previous definition, and we will do it also in the
following, we have used the Landau notation of little-oh functions
(see e.g. \cite{Pro, Sil}). E.g. any H$\text{\textipa{ö}}$lder function
$|x(t)-x(s)|\le c\cdot|t-s|^{\alpha}$ (for some constant $\alpha>0$)
is nilpotent. The choice of nilpotent functions instead of more regular
ones establish a great difference of our approach with respect to
the classical definition of jets (see e.g. \cite{Bro, Go-Gu}), that
\eqref{eq:IdeaRelEq} may recall. Indeed in our approach all the $\mathcal{C}^{1}$-functions
$x$ with the same value and derivative at $t=0$ generate the same
$\sim$-equivalence relation. Only a non differentiable function at
$t=0$ like $x(t)=\sqrt{t}$ generates non trivial nilpotent infinitesimals.

Another problem necessarily connected with the basic idea \eqref{eq:FunctionEqualTangent}
is that the use of nilpotent infinitesimals very frequently leads
to consider terms like $h_{1}^{i_{1}}\cdot\ldots\cdot h_{n}^{i_{n}}$.
For this type of products the first problem is to know whether $h_{1}^{i_{1}}\cdot\ldots\cdot h_{n}^{i_{n}}\ne0$
and what is the order $k$ of this new infinitesimals, that is for
what $k$ we have $(h_{1}^{i_{1}}\cdot\ldots\cdot h_{n}^{i_{n}})^{k}\ne0$
but $(h_{1}^{i_{1}}\cdot\ldots\cdot h_{n}^{i_{n}})^{k+1}=0$. We will
have a good frame if we will be able to solve these problems starting
from the order of each infinitesimal $h_{j}$ and from the value of
the power $i_{j}\in\N$. On the other hand almost all the examples
of nilpotent infinitesimals are of the form $h(t)=t^{\alpha}$, with
$0<\alpha<1$, and their sums; these functions have great properties
both in the treatment of products of powers and, as we will see, in
connection with the order relation. It is for these reasons that we
shall focus our attention on the following family of functions $x:\R_{\ge0}\freccia\R$
in the definition \eqref{eq:IdeaRelEq} of $\sim$:
\begin{defn}
\label{def:LittleOhPolynomials} We say that $x$ \emph{is a little-oh
polynomial}, and we write $x\in\R_{o}[t]$ iff 
\begin{enumerate}
\item $x:\R_{\ge0}\freccia\R$ 
\item We can write\[
x_{t}=r+\sum\limits _{i=1}^{k}\alpha_{i}\cdot t^{a_{i}}+o(t)\quad\text{ as }\quad t\to0^{+}\]
 for suitable\[
k\in\N\]
\[
r,\alpha_{1},\dots,\alpha_{k}\in\R\]
\[
a_{1},\dots,a_{k}\in\R_{\ge0}\]

\end{enumerate}
\end{defn}
\noindent Hence a little-oh polynomial%
\footnote{actually in the following notation the variable $t$ is mute%
} $x\in\R_{o}[t]$ is a polynomial function with real coefficients,
in the real variable $t\ge0$, with generic positive powers of $t$,
and up to a little-oh function as $t\to0^{+}$. 
\begin{rem}
\label{rem:WithLittleOhWeMean}In the following, writing $x_{t}=y_{t}+o(t)$
as $t\to0+$ we will always mean\[
\lim_{t\to0^{+}}\frac{x_{t}-y_{t}}{t}=0\quad\text{and\ensuremath{\quad}}x_{0}=y_{0}.\]
 In other words, every little-oh function we will consider is continuous
as $t\to0^{+}$.\end{rem}
\begin{example*}
Simple examples of little-oh polynomials are the following:\end{example*}
\begin{enumerate}
\item $x_{t}=1+t+t^{1/2}+t^{1/3}+o(t)$ 
\item $x_{t}=r\quad\forall t$. Note that in this example we can take $k=0$,
and hence $\alpha$ and $a$ are the void sequence of reals, that
is the function $\alpha=a:\emptyset\freccia\R$, if we think of an
$n$-tuple $x$ of reals as a function $x:\left\{ 1,\dots,n\right\} \freccia\R$. 
\item $x_{t}=r+o(t)$ 
\end{enumerate}

\section{First properties of little-oh polynomials\label{sec:firstPropertiesOfLittle-ohPolynomials}}

\subsubsection*{Little-oh polynomials are nilpotent:}

First properties of little-oh polynomials are the following: if $x_{t}=r+\sum_{i=1}^{k}\alpha_{i}\cdot t^{a_{i}}+o_{1}(t)$
as $t\to0^{+}$ and $y_{t}=s+\sum_{j=1}^{N}\beta_{j}\cdot t^{b_{j}}+o_{2}(t)$,
then $(x+y)=r+s+\sum_{i=1}^{k}\alpha_{i}\cdot t^{a_{i}}+\sum_{j=1}^{N}\beta_{j}\cdot t^{b_{j}}+o_{3}(t)$
and $(x\cdot y)_{t}=rs+\sum_{i=1}^{k}s\alpha_{i}\cdot t^{a_{i}}+\sum_{j=1}^{N}r\beta_{j}\cdot t^{b_{i}}+\sum_{i=1}^{k}\sum_{j=1}^{N}\alpha_{i}\beta_{j}\cdot t^{a_{i}}t^{b_{j}}+o_{4}(t)$,
hence the set of little-oh polynomials is closed with respect to pointwise
sum and product. Moreover little-oh polynomials are nilpotent (see
Definition \ref{def:NilpotentFunctions}) functions; to prove this
we firstly prove that the set of nilpotent functions $\mathcal{N}$
is a subalgebra of the algebra $\R^{\R}$ of real valued functions.
Indeed, let $x$ and $y$ be two nilpotent functions such that $|x-x(0)|^{k}=o_{1}(t)$
and $|y-y(0)|^{N}=o_{2}(t)$, then we can write $x\cdot y-x(0)\cdot y(0)=x\cdot[y-y(0)]+y(0)\cdot[x-x(0)]$,
so that we can consider $|x\cdot[y-y(0)]|^{k}=|x|^{k}\cdot|y-y(0)|^{k}=|x|^{k}\cdot o_{1}(t)$
and $\frac{|x|^{k}\cdot o_{1}(t)}{t}\to0$ as $t\to0^{+}$ because
$|x|^{k}\to|x(0)|^{k}$, hence $x\cdot[y-y(0)]\in\mathcal{N}$. Analogously
$y(0)\cdot[x-x(0)]\in\mathcal{N}$ and hence the closure of $\mathcal{N}$
with respect to the product follows from the closure with respect
to the sum. The case of the sum follows from the following equalities
(where we use $x_{t}:=x(t)$, $u:=x-x_{0}$, $v:=y-y_{0}$, $|u_{t}|^{k}=o_{1}(t)$
and $|v_{t}|^{N}=o_{2}(t)$ and we have supposed $k\ge N$):\[
u^{k}=o_{1}(t),\,\,\, v^{k}=o_{2}(t)\]
 \[
(u+v)^{k}=\sum_{i=0}^{k}\binom{k}{i}u^{i}\cdot v^{k-i}\]
 \[
\forall i=0,\dots,k\pti\frac{u_{t}^{i}\cdot v_{t}^{k-i}}{t}=\frac{\left(u_{t}^{k}\right)^{\frac{i}{k}}\cdot\left(v_{t}^{k}\right)^{\frac{k-i}{k}}}{t^{\frac{i}{k}}\cdot t^{\frac{k-i}{k}}}=\left(\frac{u_{t}^{k}}{t}\right)^{\frac{i}{k}}\cdot\left(\frac{v_{t}^{k}}{t}\right)^{\frac{k-i}{k}}.\]
 Now we can prove that $\R_{o}[t]$ is a subalgebra of $\mathcal{N}$.
Indeed every constant $r\in\R$ and every power $t^{a_{i}}$ are elements
of $\mathcal{N}$ and hence $r+\sum_{i=1}^{k}\alpha_{i}\cdot t^{a_{i}}\in\mathcal{N}$,
so it remains to prove that if $y\in\mathcal{N}$ and $w=o(t)$, then
$y+w\in\mathcal{N}$, but this is a consequence of the fact that every
little-oh function is trivially nilpotent, and hence it follows from
the closure of $\mathcal{N}$ with respect to the sum.

\subsubsection*{Closure of little-oh polynomials with respect to smooth functions:\label{sub:ClosureOfLittle-ohPolyWRTSmoothFunctions}}

Now we want to prove that little-oh polynomials are preserved by smooth
functions, that is if $x\in\R_{o}[t]$ and $f:\R\freccia\R$ is smooth,
then $f\circ x\in\R_{o}[t]$. Let us fix some notations:\[
x_{t}=r+\sum_{i=1}^{k}\alpha_{i}\cdot t^{a_{i}}+w(t)\quad\text{with}\quad w(t)=o(t)\]
 \[
h(t):=x(t)-x(0)\quad\forall t\in\R_{\ge0}\]
 hence $x_{t}=x(0)+h_{t}=r+h_{t}$. The function $t\mapsto h(t)=\sum_{i=1}^{k}\alpha_{i}\cdot t^{a_{i}}+w(t)$
belongs to $\R_{o}[t]\subseteq\mathcal{N}$ so we can write $|h|^{N}=o(t)$
for some $N\in\N$ and as $t\to0^{+}$. From Taylor's's formula we
have\begin{equation}
f(x_{t})=f(r+h_{t})=f(r)+\sum_{i=1}^{N}\frac{f^{(i)}(r)}{i!}\cdot h_{t}^{i}+o(h_{t}^{N})\label{eq:TaylorLittlePolyClosedSmooth}\end{equation}
 But\[
\frac{|o(h_{t}^{N})|}{|t|}=\frac{|o(h_{t}^{N})|}{|h_{t}^{N}|}\cdot\frac{|h_{t}^{N}|}{|t|}\to0\]
 hence $o(h_{t}^{N})=o(t)\in\R_{o}[t]$. From this, the formula \eqref{eq:TaylorLittlePolyClosedSmooth},
the fact that $h\in\R_{o}[t]$ and using the closure of little-oh
polynomials with respect to ring operations, the conclusion $f\circ x\in\R_{o}[t]$
follows.

\section{Equality and decomposition of Fermat reals}
\begin{defn}
\label{def:equalityInFermatReals}Let $x$, $y\in\R_{o}[t]$, then
we say that $x\sim y$ or that $x=y$ in $\ER$ iff $x(t)=y(t)+o(t)$
as $t\to0^{+}$. Because it is easy to prove that $\sim$ is an equivalence
relation, we can define $\ER:=\R_{o}[t]/\sim$, i.e. $\ER$ is the
quotient set of $\R_{o}[t]$ with respect to the equivalence relation
$\sim$. 
\end{defn}
\noindent The equivalence relation $\sim$ is a congruence with respect
to pointwise operations, hence $\ER$ is a commutative ring. Where
it will be useful to simplify notations we will write {}``$x=y$
in $\ER$'' instead of $x\sim y$, and we will talk directly about
the elements of $\R_{o}[t]$ instead of their equivalence classes;
for example we can say that $x=y$ in $\ER$ and $z=w$ in $\ER$
imply $x+z=y+w$ in $\ER$.\\
 The immersion of $\R$ in $\ER$ is $r\longmapsto\hat{r}$ defined
by $\hat{r}(t):=r$, and in the sequel we will always identify $\hat{\R}$
with $\R$, which is hence a subring of $\ER$. Conversely if $x\in\ER$
then the map $\st(-):x\in\ER\mapsto\st x=x(0)\in\R$, which evaluates
each extended real in $0$, is well defined. We shall call $\st(-)$
the \emph{standard part map}%
\footnote{This denomination should obviously not be confused with the one with
the same name in NSA.%
}. Let us also note that, as a vector space over the field $\R$ we
have $\dim_{\R}\ER=\infty$, and this underlines even more the difference
of our approach with respect to the classical definition of jets (see
e.g. \cite{Bro, Go-Gu}). As we will see, more explicitly later on
in the course of the present work, our idea is more near to NSA, where
standard sets can be extended adding new infinitesimal points, and
this is not the point of view of jet theory.

With the following theorem we will introduce the decomposition of
a Fermat real $x\in\ER$, that is a unique notation for its standard
part and all its infinitesimal parts.
\begin{thm}
\label{thm:existenceUniquenessDecomposition}If $x\in\ER$, then there
exist one and only one sequence\[
(k,r,\alpha_{1},\ldots,\alpha_{k},a_{1},\ldots,a_{k})\]
 such that\[
k\in\N\]
\[
r,\alpha_{1},\dots,\alpha_{k},a_{1},\dots,a_{k}\in\R\]
and 
\begin{enumerate}
\item \label{enu:DecompositionSum}\emph{$x=r+\sum\limits _{i=1}^{k}\alpha_{i}\cdot t^{a_{i}}$}
in \emph{$\ER$} 
\item \emph{\label{enu:DecompositionIncreasingOfOrders}$0<a_{1}<a_{2}<\dots<a_{k}\le1$} 
\item \emph{\label{enu:DecompositionInfinitesimalPartsNotZero}$\alpha_{i}\ne0\quad\forall i=1,\dots,k$} 
\end{enumerate}
\end{thm}
\noindent In this statement we have also to include the void case
$k=0$ and $\alpha=a:\emptyset\freccia\R$. Obviously, as usual, we
use the definition $\sum_{i=1}^{0}b_{i}=0$ for the sum of an empty
set of numbers. As we shall see, this is the case where $x$ is a
standard real, i.e. $x\in\R$.\\
 In the following we will use the notations $t^{a}:=\diff{t_{1/a}}:=[t\in\R_{\ge0}\mapsto t^{a}\in\R]_{\sim}\in\ER$
so that e.g. $\diff{t_{2}}=t^{1/2}$ is a second order infinitesimal%
\footnote{Let us point out that we make hereby an innocuous abuse of language
using the same notation both for the value of the function, $t^{a}\in\R$,
and for the equivalence class, $t^{a}\in\ER$.%
}. In general, as we will see from the definition of order of a generic
infinitesimal, $\diff{t_{a}}$ is an infinitesimal of order $a$.
In other words these two notations for the same object permit to emphasize
the difference between an actual infinitesimal $\diff{t_{a}}$ and
a potential infinitesimal $t^{1/a}$: an actual infinitesimal of order
$a\ge1$ corresponds to a potential infinitesimal of order $\frac{1}{a}\le1$
(with respect to the classical notion of order of an infinitesimal
function from calculus, see e.g. \cite{Pro, Sil}).
\begin{rem}
\label{rem:simplePropertiesOf_diff}Let us note that $\diff{t_{a}}\cdot\diff{t_{b}}=\diff{t_{\frac{ab}{a+b}}}$,
moreover $\diff{t}_{a}^{\alpha}:=(\diff{t_{a}})^{\alpha}=\diff{t_{\frac{a}{\alpha}}}$
for every $\alpha\ge1$ and finally $\diff{t_{a}}=0$ for every $a<1$.
E.g. $\diff{t}_{a}^{[a]+1}=0$ for every $a\in\R_{>0}$, where $[a]\in\N$
is the integer part of $a$, i.e. $[a]\le a<[a]+1$. 
\end{rem}
\noindent \textbf{Existence proof:}

Since $x\in\R_{o}[t]$, we can write $x_{t}=r+\sum_{i=1}^{k}\alpha_{i}\cdot t^{a_{i}}+o(t)$
as $t\to0^{+}$, where $r$, $\alpha_{i}\in\R$, $a_{i}\in\R_{\ge0}$
and $k\in\N$. Hence $x=r+\sum_{i=1}^{k}\alpha_{i}\cdot t^{a_{i}}$
in $\ER$ and our purpose is to pass from this representation of $x$
to another one that satisfies conditions \ref{enu:DecompositionSum},
\ref{enu:DecompositionIncreasingOfOrders} and \ref{enu:DecompositionInfinitesimalPartsNotZero}
of the statement. Since if $a_{i}>1$ then $\alpha_{i}\cdot t^{a_{i}}=0$
in $\ER$, we can suppose that $a_{i}\le1$ for every $i=1,\dots,k$.
Moreover we can also suppose $a_{i}>0$ for every $i$, because otherwise,
if $a_{i}=0$, we can replace $r\in\R$ by $r+\sum\{\alpha_{i}\,|\, a_{i}=0,\, i=1,\dots,k\}$.

\noindent Now we sum all the terms $t^{a_{i}}$ having the same $a_{i}$,
that is we can consider\[
\bar{\alpha_{i}}:=\sum\{\alpha_{j}\,|\, a_{j}=a_{i}\,,\, j=1,\dots,k\}\]
 so that in $\ER$ we have\[
x=r+\sum_{i\in I}\bar{\alpha_{i}}\cdot t^{a_{i}}\]
 where $I\subseteq\left\{ 1,\dots,k\right\} $, $\{a_{i}\,|\, i\in I\}=\{a,\dots,a_{k}\}$
and $a_{i}\ne a_{j}$ for any $i$, $j\in I$ with $i\ne j$. Neglecting
$\alpha_{i}$ if $\alpha_{i}=0$ and renaming $a_{i}$, for $i\in I$,
in such a way that $a_{i}<a_{j}$ if $i$, $j\in I$ with $i<j$,
we obtain the existence result. Note that if $x=r\in\R$, in the final
step of this proof we have $I=\emptyset$.\medskip{}

\noindent \textbf{Uniqueness proof:}

Let us suppose that in $\ER$ we have\begin{equation}
x=r+\sum_{i=1}^{k}\alpha_{i}\cdot t^{a_{i}}=s+\sum_{j=1}^{N}\beta_{j}\cdot t^{b_{j}}\label{eq:uniquenessDecompositionHypothesis}\end{equation}
 where $\alpha_{i}$, $\beta_{j}$, $a_{i}$ and $b_{j}$ verify the
conditions of the statement. First of all $\st{x}=x(0)=r=s$ because
$a_{i}$, $b_{j}>0$. Hence $\alpha_{1}t^{a_{1}}-\beta_{1}t^{b_{1}}+\sum_{i}\alpha_{i}\cdot t^{a_{i}}-\sum_{j}\beta_{j}\cdot t^{b_{j}}=o(t)$.
By reduction to the absurd, if we had $a_{1}<b_{1}$, then collecting
the term $t^{a_{1}}$ we would have\begin{equation}
\alpha_{1}-\beta_{1}t^{b_{1}-a_{1}}+\sum_{i}\alpha_{i}\cdot t^{a_{i}-a_{1}}-\sum_{j}\beta_{j}\cdot t^{b_{j}-a_{1}}=\frac{o(t)}{t}\cdot t^{1-a_{1}}.\label{eq:firstTermUniquenessProofDecomposition}\end{equation}
 In \eqref{eq:firstTermUniquenessProofDecomposition} we have that
$\beta_{1}t^{b_{1}-a_{1}}\to0$ for $t\to0^{+}$ because $a_{1}<b_{1}$
by hypothesis; $\sum_{i}\alpha_{i}\cdot t^{a_{i}-a_{1}}\to0$ because
$a_{1}<a_{i}$ for $i=2,\dots,k$; $\sum_{j}\beta_{j}\cdot t^{b_{j}-a_{1}}\to0$
because $a_{1}<b_{1}<b_{j}$ for $j=2,\dots,N$, and finally $t^{1-a_{1}}$
is limited because $a_{1}\le1$. Hence for $t\to0^{+}$ we obtain
$\alpha_{1}=0$, which conflicts with condition \ref{enu:DecompositionInfinitesimalPartsNotZero}
of the statement. We can argue in a corresponding way if we had $b_{1}<a_{1}$.
In this way we see that we must have $a_{1}=b_{1}$. From this and
from equation \eqref{eq:firstTermUniquenessProofDecomposition} we
obtain\begin{equation}
\alpha_{1}-\beta_{1}+\sum_{i}\alpha_{i}\cdot t^{a_{i}-a_{1}}-\sum_{j}\beta_{j}\cdot t^{b_{j}-a_{1}}=\frac{o(t)}{t}\cdot t^{1-a_{1}}\label{eq:endFirstTerm}\end{equation}
 and hence for $t\to0^{+}$ we obtain $\alpha_{1}=\beta_{1}$. We
can now restart from \eqref{eq:endFirstTerm} to prove, in the same
way, that $a_{2}=b_{2}$, $\alpha_{2}=\beta_{2}$, etc. At the end
we must have $k=N$ because, otherwise, if we had e.g. $k<N$, at
the end of the previous recursive process, we would have\[
\sum_{j=k+1}^{N}\beta_{j}\cdot t^{b_{j}}=o(t).\]
 From this, collecting the terms containing $t^{b_{k+1}}$, we obtain\begin{equation}
t^{b_{k+1}-1}\cdot[\beta_{k+1}+\beta_{k+2}\cdot t^{b_{k+2}-b_{k+1}}+\dots+\beta_{N}\cdot t^{\beta_{N}-\beta_{k+1}}]\to0.\label{eq:uniquenessDecompositionN_equal_k}\end{equation}
 In this sum $\beta_{k+j}\cdot t^{b_{k+j}-b_{k+1}}\to0$ as $t\to0^{+}$,
because $b_{k+1}<b_{k+j}$ for $j>1$ and hence $\beta_{k+1}+\beta_{k+2}\cdot t^{b_{k+2}-b_{k+1}}+\dots+\beta_{N}\cdot t^{\beta_{N}-\beta_{k+1}}\to\beta_{k+1}\ne0$,
so from \eqref{eq:uniquenessDecompositionN_equal_k} we get $t^{b_{k+1}-1}\to0$,
that is $b_{k+1}>1$, in contradiction with the uniqueness hypothesis
$b_{k+1}\le1$.

Let us note explicitly that the uniqueness proof permits also to affirm
that the decomposition is well defined in $\ER$, i.e. that if $x=y$
in $\ER$, then the decomposition of $x$ and the decomposition of
$y$ are equal.\qedNoNewLine

On the basis of this theorem we introduce two notations: the first
one emphasizing the potential nature of an infinitesimal $x\in\ER$,
and the second one emphasizing its actual nature.
\begin{defn}
\label{def:potentialDecomposition}If $x\in\ER$, we say that\begin{equation}
x=r+\sum_{i=1}^{k}\alpha_{i}\cdot t^{a_{i}}\ \text{is the potential decomposition (of }x\text{)}\label{eq:potentialDecomposition}\end{equation}
 iff conditions \ref{enu:DecompositionSum}., \ref{enu:DecompositionIncreasingOfOrders}.,
and \ref{enu:DecompositionInfinitesimalPartsNotZero}. of theorem
\ref{thm:existenceUniquenessDecomposition} are verified. Of course
it is implicit that the symbol of equality in \eqref{eq:potentialDecomposition}
has to be understood in $\ER$. 
\end{defn}
\noindent For example $x=1+t^{1/3}+t^{1/2}+t$ is a decomposition
because we have increasing powers of $t$. The only decomposition
of a standard real $r\in\R$ is the void one, i.e. that with $k=0$
and $\alpha=a:\emptyset\freccia\R$; indeed to see that this is the
case, it suffices to go along the existence proof again with this
case $x=r\in\R$ (or to prove it directly, e.g. by contradiction).
\begin{defn}
\label{def:actualDecomposition}Considering that $t^{a_{i}}=\diff{t_{1/a_{i}}}$
we can also use the following notation, emphasizing more the fact
that $x\in\ER$ is an actual infinitesimal:\begin{equation}
x=\st{x}+\sum_{i=1}^{k}\st{x_{i}}\cdot\diff{t_{b_{i}}}\label{eq:firstActualDecomposition}\end{equation}
 where we have used the notation $\st{x_{i}}:=\alpha_{i}$ and $b_{i}:=1/a_{i}$,
so that the condition that uniquely identifies all $b_{i}$ is $b_{1}>b_{2}>\dots>b_{k}\ge1$.
We call \eqref{eq:firstActualDecomposition} the \emph{actual decomposition}
of $x$ or simply the \emph{decomposition} of $x$. We will also use
the notation $\diff^{i}{x}:=\st{x_{i}}\cdot\diff{t_{b_{i}}}$ (and
simply $\diff{x}:=\diff{^{1}x}$) and we will call $\st{x_{i}}$ the
$i$\emph{-th standard part of} $x$ and $\diff{^{i}x}$ the $i$\emph{-th
infinitesimal part of} $x$ or the $i$-th \emph{differential} of
$x$. So let us note that we can also write\[
x=\st{x}+\sum_{i}\diff{^{i}x}\]
 and in this notation all the addenda are uniquely determined (the
number of them too). Finally, if $k\ge1$ that is if $x\in\ER\setminus\R$,
we set $\omega(x):=b_{1}$ and $\omega_{i}(x):=b_{i}$. The real number
$\omega(x)=b_{1}$ is the greatest order in the actual decomposition
\eqref{eq:firstActualDecomposition}, corresponding to the smallest
in the potential decomposition \eqref{eq:potentialDecomposition},
and is called the \emph{order} of the Fermat real $x\in\ER$. The
number $\omega_{i}(x)=b_{i}$ is called the $i$-th order of $x$.
If $x\in\R$ we set $\omega(x):=0$ and $\diff^{i}x:=0$. Observe
that in general $\omega(x)=\omega(\diff{x})$, $\diff{(\diff{x})}=\diff{x}$
and that, using the notations of the potential decomposition \eqref{def:potentialDecomposition},
we have $\omega(x)=1/a_{1}$. \end{defn}
\begin{example*}
If $x=1+t^{1/3}+t^{1/2}+t$, then $\st{x}=1$, $\diff{x}=\diff{t_{3}}$
and hence $x$ is a third order infinitesimal, i.e. $\omega(x)=3$,
$\diff{^{2}x}=\diff{t_{2}}$ and $\diff{^{3}x}=\diff{t}$; finally
all the standard parts are $\st{x_{i}}=1$. \end{example*}
\begin{rem}
\emph{\label{rem:differentMeaningForStandardAndPotentialOrder} }To
avoid misunderstanding, it is important to underline that there is
an opposite meaning of the word {}``order'' in standard analysis
and in the previous definition. Indeed, in standard analysis if we
say that the infinitesimal function (for $t\to0^{+})$ $t\mapsto x(t)$
is of order greater than the function $t\mapsto y(t)$, we mean that\[
\lim_{t\to0^{+}}\frac{x(t)}{y(t)}=0\]
 Intuitively this implies that we have to think $x(t)$ smaller than
$y(t)$, at least for sufficiently small $t\in(0,\delta)$. Because
the connection between the definition of order given in Definition
\ref{def:actualDecomposition} and the standard definition of order
(with respect to the standard infinitesimal $t\mapsto t$) is given
by $\omega(x)=1/a_{1}$, for Fermat reals the meaning will be the
opposite one: if $x$, $y\in D_{\infty}$ are two infinitesimals,
having every standard part positive $\st{x_{i}}$, $\st{y_{j}}>0$,
and with $\omega(x)>\omega(y)$, then we have to think at $x\in\ER$
as a\emph{ bigger number }with respect to $y\in\ER$. More formally,
in the next section we will see that this will correspond to say that
if $a$, $b\in\N$ are such that $x^{a}\ne0$ and $x^{a+1}=0$, $y^{b}\ne0$
and $y^{b+1}=0$, then $a\ge b$. When we will introduce the order
relation in $\ER$ (see \ref{cha:orderRelation}), we will see that
if for two infinitesimals we have $\omega(x)>\omega(y)$, then $x>y$
iff $\st{x_{1}}>0$. Recalling the Remark \ref{rem:simplePropertiesOf_diff}
we can remember this difference between classical and actual order,
recalling that $\diff{t_{a}}>\diff{t_{b}}$ if $a>b$ and that the
smallest non zero infinitesimal is $\diff{t_{1}}=\diff{t}$, because
$\diff{t_{a}}=0$ if $a<1$.\emph{ }
\end{rem}

\section{The ideals $D_{k}$}

In this section we will introduce the sets of nilpotent infinitesimals
corresponding to a $k$-th order neighborhood of 0. Every smooth function
restricted to this neighborhood becomes a polynomial of order $k$,
obviously given by its $k$-th order Taylor's's formula (without remainder).
We start with a theorem characterizing infinitesimals of order less
than $k$.
\begin{thm}
\label{thm:characterizationOf_x_to_k_equal_0}If $x\in\ER$ and $k\in\N_{>1}$,
then $x^{k}=0$ in $\ER$ if and only if $\st{x}=0$ and $\omega(x)<k$. 
\end{thm}
\noindent \textbf{Proof:} If $x^{k}=0$, then taking the standard
part map of both sides, we have $\st{(x^{k})}=(\st{x})^{k}=0$ and
hence $\st{x}=0$. Moreover $x^{k}=0$ means $x_{t}^{k}=o(t)$ and
hence $\left(\frac{x_{t}}{t^{1/k}}\right)^{k}\to0$ and $\frac{x_{t}}{t^{1/k}}\to0$.
We rewrite this condition using the potential decomposition $x=\sum_{i=1}^{k}\alpha_{i}\cdot t^{a_{i}}$
of $x$ (note that in this way we have $\omega(x)=\frac{1}{a_{1}}$)
obtaining \[
\lim_{t\to0^{+}}\sum_{i}\alpha_{i}\cdot t^{a_{i}-\frac{1}{k}}=0=\lim_{t\to0^{+}}t^{a_{1}-\frac{1}{k}}\cdot\left[\alpha_{1}+\alpha_{2}\cdot t^{a_{2}-a_{1}}+\dots+\alpha_{k}\cdot t^{a_{k}-a_{1}}\right]\]
 But $\alpha_{1}+\alpha_{2}\cdot t^{a_{2}-a_{1}}+\dots+\alpha_{k}\cdot t^{a_{k}-a_{1}}\to\alpha_{1}\ne0$,
hence we must have that $t^{a_{1}-\frac{1}{k}}\to0$, and so $a_{1}>\frac{1}{k}$,
that is $\omega(x)<k$.\\
 Vice versa if $\st{x}=0$ and $\omega(x)<k$, then $x=\sum_{i=1}^{k}\alpha_{i}\cdot t^{a_{i}}+o(t)$,
and\[
\lim_{t\to0^{+}}\frac{x_{t}}{t^{1/k}}=\lim_{t\to0^{+}}\sum_{i}\alpha_{i}\cdot t^{a_{i}-\frac{1}{k}}+\lim_{t\to0^{+}}\frac{o(t)}{t}\cdot t^{1-\frac{1}{k}}\]
 But $t^{1-\frac{1}{k}}\to0$ because $k>1$ and $t^{a_{i}-\frac{1}{k}}\to0^{+}$
because $\frac{1}{a_{i}}\le\frac{1}{a_{1}}=\omega(x)<k$ and hence
$x^{k}=0$ in $\ER$.\qedNoNewLine

If we want that in a $k$-th order infinitesimal neighborhood a smooth
function is equal to its $k$-th Taylor's's formula, i.e.\[
\forall h\in D_{k}:f(x+h)=\sum_{i=0}^{k}\frac{h^{i}}{i!}\cdot f^{(i)}(x)\]
 we need to take infinitesimals which are able to delete the remainder,
that is, such that $h^{k+1}=0$. The previous theorem permits to extend
the definition of the ideal $D_{k}$ to real number subscripts instead
of natural numbers $k$ only.
\begin{defn}
\label{def:idealD_a}If $a\in\R_{>0}\cup\left\{ \infty\right\} $,
then\[
D_{a}:=\left\{ x\in\ER\,|\,\st{x}=0,\ \omega(x)<a+1\right\} \]
 Moreover we will simply denote $D_{1}$ by $D$. \end{defn}
\begin{enumerate}
\item If $x=\diff{t_{3}}$, then $\omega(x)=3$ and $x\in D_{3}$. More
in general $\diff{t_{k}}\in D_{a}$ if and only if $\omega(\diff{t_{k}})=k<a+1$.
E.g. $\diff{t_{k}}\in D$ if and only if $1\le k<2$. 
\item $D_{\infty}=\bigcup_{a}D_{a}=\left\{ x\in\ER\,|\,\st{x}=0\right\} $
is the set of all the infinitesimals of $\ER$. 
\item $D_{0}=\left\{ 0\right\} $ because the only infinitesimal having
order strictly less than 1 is, by definition of order, $x=0$ (see
the Definition \ref{def:actualDecomposition}). 
\end{enumerate}
The following theorem gathers several expected properties of the sets
$D_{a}$ and of the order of an infinitesimal $\omega(x)$:
\begin{thm}
\label{thm:propertyOfIdeal_D_a}Let $a$, $b\in\R_{>0}$ and $x$,
$y\in D_{\infty}$, then 
\begin{enumerate}
\item \label{enu:IncreasingIdeals}$a\le b\then D_{a}\subseteq D_{b}$ 
\item \label{enu:xInItsOwnIdeal}$x\in D_{\omega(x)}$ 
\item \label{enu:consistencyWithIntegerDefinition}$a\in\N\then D_{a}=\{x\in\ER\,|\, x^{a+1}=0\}$ 
\item \label{enu:fromRealIdealToIntegerPower}$x\in D_{a}\then x^{\lceil a\rceil+1}=0$ 
\item \label{enu:integerPartOfOrder}$x\in D_{\infty}\setminus\{0\}$ and
$k=[\omega(x)]\then x\in D_{k}\setminus D_{k-1}$ 
\item \label{enu:diffOfProduct}$\diff{(x\cdot y)}=\diff{x}\cdot\diff{y}$ 
\item \label{enu:orderOfProduct}$x\cdot y\ne0\then{\displaystyle \frac{1}{\omega(x\cdot y)}=\frac{1}{\omega(x)}+\frac{1}{\omega(y)}}$ 
\item \label{enu:orderOfSum}$x+y\ne0\then\omega(x+y)=\omega(x)\vee\omega(y)$ 
\item \label{enu:D_aIsAnIdeal}$D_{a}$ is an ideal 
\end{enumerate}
\end{thm}
\noindent In this statement if $r\in\R$, then $\lceil r\rceil$ is
the \emph{ceiling} of the real $r$, i.e. the unique integer $\lceil r\rceil\in\mathbb{Z}$
such that $\lceil r\rceil-1<r\le\lceil r\rceil$. Moreover if $r$,
$s\in\R$, then $r\vee s:=\max(r,s)$.

\medskip{}

\noindent \textbf{Proof:} Property \emph{\ref{enu:IncreasingIdeals}}.
and \emph{\ref{enu:xInItsOwnIdeal}}. follow directly from Definition
\ref{def:idealD_a} of $D_{a}$, whereas property \emph{\ref{enu:consistencyWithIntegerDefinition}}.
follows from Theorem \ref{thm:characterizationOf_x_to_k_equal_0}.
From \emph{\ref{enu:IncreasingIdeals}}. and \emph{\ref{enu:consistencyWithIntegerDefinition}}.
property \emph{\ref{enu:fromRealIdealToIntegerPower}}. follows: in
fact $x\in D_{a}\subseteq D_{\lceil a\rceil}$ because $a\le\lceil a\rceil$,
hence $x^{\lceil a\rceil+1}=0$ from property \emph{\ref{enu:consistencyWithIntegerDefinition}}.
To prove property \emph{\ref{enu:integerPartOfOrder}}., if $k=[\omega(x)]$,
then $k\le\omega(x)<k+1$, hence directly from Definition \ref{def:idealD_a}
the conclusion follows.

\noindent To prove \emph{\ref{enu:diffOfProduct}}. let\begin{equation}
x=\sum_{i=1}^{k}\st{x_{i}}\cdot\diff{t_{a_{i}}}\quad\text{and}\quad y=\sum_{j=1}^{N}\st{y_{j}}\cdot\diff{t_{b_{j}}}\label{eq:decompositionOfx_y}\end{equation}
 be the decompositions of $x$ and $y$ (considering that they are
infinitesimals, so that $\st{x}=\st{y}=0$). Recall that $\diff{x}=\st{x_{1}}\cdot\diff{t_{a_{1}}}$
and $\diff{y}=\st{y_{1}}\cdot\diff{t_{b_{1}}}$. From \eqref{eq:decompositionOfx_y}
we have \begin{equation}
x\cdot y=\sum_{i=1}^{k}\sum_{j=1}^{N}\st{x_{i}}\st{y_{j}}\diff{t_{a_{i}}}\diff{t_{b_{j}}}=\sum_{i=1}^{k}\sum_{j=1}^{N}\st{x_{i}}\st{y_{j}}\diff{t_{\frac{a_{i}b_{j}}{a_{i}+b_{j}}}}\label{eq:decompositionOfProduct}\end{equation}
 where we have used the Remark \emph{\ref{rem:simplePropertiesOf_diff}}.
But $\omega(x)=a_{1}\ge a_{i}$ and $\omega(y)=b_{1}\ge b_{j}$ from
the Definition \ref{def:actualDecomposition} of decomposition. Hence\[
\frac{1}{a_{1}}+\frac{1}{b_{1}}\le\frac{1}{a_{i}}+\frac{1}{b_{j}}\]
 \[
\frac{a_{1}b_{1}}{a_{1}+b_{1}}\ge\frac{a_{i}b_{j}}{a_{i}+b_{j}}\]
 so that the greatest infinitesimal in the product \eqref{eq:decompositionOfProduct}
is\[
\diff{(x\cdot y)}=\st{x_{1}}\st{y_{1}}\diff{t_{a_{1}}}\diff{t_{b_{1}}}=\diff{x}\cdot\diff{y}\]
 From this proof, property \emph{\ref{enu:orderOfProduct}}. follows,
because $x\cdot y\ne0$ by hypothesis, and hence its order is given
by\[
\omega(x\cdot y)=\frac{a_{1}b_{1}}{a_{1}+b_{1}}=\left(\frac{1}{a_{1}}+\frac{1}{b_{1}}\right)^{-1}=\left(\frac{1}{\omega(x)}+\frac{1}{\omega(y)}\right)^{-1}\]
 From the decompositions \eqref{eq:decompositionOfx_y} we also have\[
x+y=\sum_{i=1}^{k}\st{x_{i}}\diff{t_{a_{i}}}+\sum_{j=1}^{N}\st{y_{j}}\diff{t_{b_{j}}}\]
 and therefore, because by hypothesis $x+y\ne0$, its order is given
by the greatest infinitesimal in this sum, that is\[
\omega(x+y)=a_{1}\vee b_{1}=\omega(x)\vee\omega(y)\]
 It remains to prove property \emph{\ref{enu:D_aIsAnIdeal}}. First
of all $\omega(0)=0<a+1$, hence $0\in D_{a}$. If $x$, $y\in D_{a}$,
then $\omega(x)$ and $\omega(y)$ are strictly less than $a+1$ and
hence $x+y\in D_{a}$ follows from property \emph{\ref{enu:orderOfSum}}.
Finally if $x\in D_{a}$ and $y\in\ER$, then $x\cdot y=x\cdot\st{y}+x\cdot(y-\st{y})$,
so $\omega(x\cdot y)=\omega(x\cdot\st{y})\vee\omega(x\cdot(y-\st{y}))=\omega(x)\vee\omega(x\cdot z)$,
where $z:=y-\st{y}\in D_{\infty}$ is an infinitesimal. If $x\cdot z=0$,
we have $\omega(x\cdot y)=\omega(x)<a+1$, otherwise from property
\emph{\ref{enu:orderOfProduct}}.\[
\frac{1}{\omega(x\cdot z)}=\frac{1}{\omega(x)}+\frac{1}{\omega(z)}\ge\frac{1}{\omega(x)}\]
 and hence $\omega(x\cdot y)\le\omega(x)<a+1$; in any case the conclusion
$x\cdot y\in D_{a}$ follows.\qedNoNewLine

Property \emph{\ref{enu:fromRealIdealToIntegerPower}}. of this theorem
cannot be proved substituting the ceiling $\lceil a\rceil$ with the
integer part $[a]$. In fact if $a=1.2$ and $x=\diff{t_{2.1}}$,
then $\omega(x)=2.1$ and $[a]+1=2$ so that $x^{[a]+1}=x^{2}=\diff{t_{\frac{2.1}{2}}}\ne0$
in $\ER$, whereas $\lceil a\rceil+1=3$ and $x^{3}=\diff{t_{\frac{2.1}{3}}}=0$.

Finally let us note the increasing sequence of ideals/neighborhoods
of zero:\begin{equation}
\{0\}=D_{0}\subset D=D_{1}\subset D_{2}\subset\dots\subset D_{k}\subset\dots\subset D_{\infty}\label{eq:sequenceOfIdeals}\end{equation}
 Because of \eqref{eq:sequenceOfIdeals} and of the property $\diff{t_{a}}=0$
if $a<1$, we can say that $\diff{t}$ is the smallest infinitesimals
and $\diff{t_{2}}$, $\diff{t_{3}}$, etc. are greater infinitesimals.
As we mentioned in the Remark \ref{rem:differentMeaningForStandardAndPotentialOrder},
after the introduction of the order relation in $\ER$, we will see
that this {}``algebraic'' idea of order of magnitude will correspond
to a property of this order relation, so that we will also have $\diff{t}<\diff{t_{2}}<\diff{t_{3}}<\dots$.
Moreover, from the properties \emph{\ref{enu:IncreasingIdeals}}.
and \emph{\ref{enu:integerPartOfOrder}}. of the previous theorem
it follows that if $x^{a}\ne0$ and $x^{a+1}=0$, then $a=[\omega(x)]$,
so that if also $y^{b}\ne0$, $y^{b+1}=0$ and $\omega(x)>\omega(y)$,
then $a\ge b$. This proves what has been stated in Remark \ref{rem:differentMeaningForStandardAndPotentialOrder}.

\section{Products of powers of nilpotent infinitesimals}

In this section we will introduce several simple instruments that
will be very useful to decide whether a product of the form $h_{1}^{i_{1}}\cdot\ldots\cdot h_{n}^{i_{n}}$
, with $h_{k}\in D_{\infty}\setminus\{0\}$, is zero or whether it
belongs to some $D_{k}$.
\begin{thm}
\label{thm:productOfPowers}Let $h_{1},\dots,h_{n}\in D_{\infty}\setminus\{0\}$
and $i_{1},\dots,i_{n}\in\N$, then 
\begin{enumerate}
\item \label{enu:iffForZeroProduct}${\displaystyle {\displaystyle h_{1}^{i_{1}}\cdot\ldots\cdot h_{n}^{i_{n}}=0}\quad\iff\quad\sum_{k=1}^{n}\frac{i_{k}}{\omega(h_{k})}>1}$ 
\item \label{enu:orderOfProductOfPowers}$h_{1}^{i_{1}}{\displaystyle \cdot\ldots\cdot h_{n}^{i_{n}}\ne0\then\frac{1}{\omega(h_{1}^{i_{1}}\cdot\ldots\cdot h_{n}^{i_{n}})}=\sum_{k=1}^{n}\frac{i_{k}}{\omega(h_{k})}}$ 
\end{enumerate}
\end{thm}
\noindent \textbf{Proof:} Let \begin{equation}
h_{k}=\sum_{r=1}^{N_{k}}\alpha_{kr}t^{a_{kr}}\label{eq:decompositionOf_h_k}\end{equation}
 be the potential decomposition of $h_{k}$ for $k=1,\dots,n$. Then
by definition \ref{def:potentialDecomposition} of potential decomposition
and the definition\ref{def:actualDecomposition} of order, we have
$0<a_{k1}<a_{k2}<\dots<a_{kN_{k}}\le1$ and $j_{k}:=\omega(h_{k})=\frac{1}{a_{k1}}$,
hence $\frac{1}{j_{k}}\le a_{kr}$ for every $r=1,\dots,N_{k}$. Therefore
from \eqref{eq:decompositionOf_h_k}, collecting the terms containing
$t^{1/j_{k}}$ we have\[
h_{k}=t^{1/j_{k}}\cdot\left(\alpha_{k1}+\alpha_{k2}t^{a_{k2}-1/j_{k}}+\dots+\alpha_{kN_{k}}t^{a_{kN_{k}-1/j_{k}}}\right)\]
 and hence\begin{align}
h_{1}^{i_{1}}\cdot\ldots\cdot h_{n}^{i_{n}} & =t^{\frac{i_{1}}{j_{1}}+\dots+\frac{i_{n}}{j_{n}}}\cdot\left(\alpha_{11}+\alpha_{12}t^{a_{12}-\frac{1}{j_{1}}}+\dots+\alpha_{1N_{1}}t^{a_{1N_{1}}-\frac{1}{j_{1}}}\right)^{i_{1}}\cdot\ldots\nonumber \\
 & \ldots\cdot\left(\alpha_{n1}+\alpha_{n2}t^{a_{n2}-\frac{1}{j_{n}}}+\dots+\alpha_{nN_{n}}t^{a_{nN_{n}}-\frac{1}{j_{n}}}\right)^{i_{n}}\label{eq:productAfterGathering}\end{align}
 Hence if $\sum_{k}\frac{i_{k}}{j_{k}}>1$ we have that $t^{\frac{i_{1}}{j_{1}}+\dots+\frac{i_{n}}{j_{n}}}=0$
in $\ER$, so also $h_{1}^{i_{1}}\cdot\ldots\cdot h_{n}^{i_{n}}=0$.
Vice versa if $h_{1}^{i_{1}}\cdot\ldots\cdot h_{n}^{i_{n}}=0$, then
the right hand side of \eqref{eq:productAfterGathering} is a $o(t)$
as $t\to0^{+}$, that is\begin{align*}
t^{\frac{i_{1}}{j_{1}}+\dots+\frac{i_{n}}{j_{n}}-1}\cdot\left(\alpha_{11}+\alpha_{12}t^{a_{12}-\frac{1}{j_{1}}}+\dots+\alpha_{1N_{1}}t^{a_{1N_{1}}-\frac{1}{j_{1}}}\right)^{i_{1}}\cdot\ldots\\
\ldots\cdot\left(\alpha_{n1}+\alpha_{n2}t^{a_{n2}-\frac{1}{j_{n}}}+\dots+\alpha_{nN_{n}}t^{a_{nN_{n}}-\frac{1}{j_{n}}}\right)^{i_{n}} & \to0\end{align*}
 But each term $\left(\alpha_{k1}+\alpha_{k2}t^{a_{k2}-\frac{1}{j_{k}}}+\dots+\alpha_{kN_{k}}t^{a_{kN_{k}}-\frac{1}{j_{k}}}\right)^{i_{k}}\to\alpha_{k}^{i_{k}}\ne0$
so, necessarily, we must have $\frac{i_{1}}{j_{1}}+\dots+\frac{i_{n}}{j_{n}}-1>0$,
and this concludes the proof of \emph{\ref{enu:iffForZeroProduct}}.

\noindent To prove \emph{\ref{enu:orderOfProductOfPowers}}. it suffices
to apply recursively property \emph{\ref{enu:orderOfProduct}}. of
Theorem \eqref{thm:propertyOfIdeal_D_a}, in fact\begin{align*}
\frac{1}{\omega(h_{1}^{i_{1}}\cdot\ldots\cdot h_{n}^{i_{n}})} & =\frac{1}{\omega(h_{1}^{i_{1}})}+\frac{1}{\omega(h_{2}^{i_{2}}\cdot\ldots\cdot h_{n}^{i_{n}})}=\\
 & \frac{1}{\omega(h_{1}\cdot\ptind^{i_{1}}\cdot h_{1})}+\frac{1}{\omega(h_{2}^{i_{2}}\cdot\ldots\cdot h_{n}^{i_{n}})}=\dots\\
 & =\frac{i_{1}}{\omega(h_{1})}+\frac{1}{\omega(h_{2}^{i_{2}}\cdot\ldots\cdot h_{n}^{i_{n}})}=\frac{i_{1}}{\omega(h_{1})}+\dots+\frac{i_{n}}{\omega(h_{n})}\end{align*}
 and this concludes the proof.\qedWithFinalEq
\begin{example}
\label{exa:productOfPowers}$\omega(\diff{t_{a_{1}}^{i_{1}}}\cdot\ldots\cdot\diff{t_{a_{n}}^{i_{n}}})^{-1}=\sum_{k}\frac{i_{k}}{\omega(\diff{t_{a_{k}}})}=\sum_{k}\frac{i_{k}}{a_{k}}$
and $\diff{t_{a_{1}}^{i_{1}}}\cdot\ldots\cdot\diff{t_{a_{n}}^{i_{n}}}=0$
if and only if $\sum_{k}\frac{i_{k}}{a_{k}}>1$, so e.g. $\diff{t}\cdot h=0$
for every $h\in D_{\infty}$. 
\end{example}
\noindent From this theorem we can derive four simple corollaries
that will be useful in the course of the present work. Some of these
corollaries are useful because they give properties of powers like
$h_{1}^{i_{1}}\cdot\ldots\cdot h_{n}^{i_{n}}$ in cases where exact
values of the orders $\omega(h_{k})$ are unknown. The first corollary
gives a necessary and sufficient condition to have $h_{1}^{i_{1}}\cdot\ldots\cdot h_{n}^{i_{n}}\in D_{p}\setminus\{0\}$.
\begin{cor}
\label{cor:idealForProductOfPowers}In the hypotheses of the previous
Theorem \ref{thm:productOfPowers} let $p\in\R_{>0}$, then we have\[
h_{1}^{i_{1}}\cdot\ldots\cdot h_{n}^{i_{n}}\in D_{p}\setminus\{0\}\quad\iff\quad\frac{1}{p+1}<\sum_{k=1}^{n}\frac{i_{k}}{\omega(h_{k})}\le1\]

\end{cor}
\noindent \textbf{Proof:} This follows almost directly from Theorem
\ref{thm:productOfPowers}. In fact if $h_{1}^{i_{1}}\cdot\ldots\cdot h_{n}^{i_{n}}\in D_{p}\setminus\{0\}$,
then its order is given by $\omega(h_{1}^{i_{1}}\cdot\ldots\cdot h_{n}^{i_{n}})=\left[\sum_{k}\frac{i_{k}}{\omega(h_{k})}\right]^{-1}=:a$
and moreover $a\ge1$ because $h_{1}^{i_{1}}\cdot\ldots\cdot h_{n}^{i_{n}}\ne0$.
Furthermore, $h_{1}^{i_{1}}\cdot\ldots\cdot h_{n}^{i_{n}}$ being
an element of $D_{p}$, we also have $a<p+1$, from which the conclusion
$\frac{1}{p+1}<\frac{1}{a}\le1$ follows.

Vice versa if $\frac{1}{p+1}<\frac{1}{a}:=\sum_{k}\frac{i_{k}}{\omega(h_{k})}\le1$,
then from Theorem \ref{thm:productOfPowers} we have $h_{1}^{i_{1}}\cdot\ldots\cdot h_{n}^{i_{n}}\ne0$
and $\omega(h_{1}^{i_{1}}\cdot\ldots\cdot h_{n}^{i_{n}})=a$; but
$a<p+1$ by hypothesis, hence $h_{1}^{i_{1}}\cdot\ldots\cdot h_{n}^{i_{n}}\in D_{p}$.\qedNoNewLine

Now we will prove a sufficient condition to have $h_{1}^{i_{1}}\cdot\ldots\cdot h_{n}^{i_{n}}=0$,
starting from the hypotheses $h_{k}\in D_{j_{k}}$ only, that is $\omega(h_{k})<j_{k}+1$.
The typical situation where this applies is for $j_{k}=[\omega(h_{k})]\in\N$.
\begin{cor}
\label{cor:suffCondToHaveProductZeroFrom_D_j_k}Let $h_{k}\in D_{j_{k}}$
for $k=1,\dots,n$ and $i_{1},\dots,i_{n}\in\N$, then\[
\sum_{k=1}^{n}\frac{i_{k}}{j_{k}+1}\ge1\then h_{1}^{i_{1}}\cdot\ldots\cdot h_{n}^{i_{n}}=0\]

\end{cor}
\noindent In fact $\sum_{k=1}^{n}\frac{i_{k}}{\omega(h_{k})}>\sum_{k=1}^{n}\frac{i_{k}}{j_{k}+1}\ge1$
because $\omega(h_{k})<j_{k}+1$, hence the conclusion follows from
Theorem \ref{thm:productOfPowers}.

Let $h$, $k\in D$; we want to see if $h\cdot k=0$. Because in this
case $\sum_{k}\frac{i_{k}}{j_{k}+1}=\frac{1}{2}+\frac{1}{2}=1$ we
always have\begin{equation}
h\cdot k=0\label{eq:productOfFirstOrderInfinitesimalIsZero}\end{equation}
 We will see that this is a great conceptual difference between Fermat
reals and the ring of SDG, where, not necessarily, the product of
two first order infinitesimal is zero. The consequences of this property
of Fermat reals arrive very deeply in the development of the theory
of Fermat reals, forcing us, e.g., to develop several new concepts
if we want to generalize the derivation formula \eqref{eq:IdeaDF}
to functions defined on infinitesimal domains, like $f:D\freccia\ER$
(see \ref{cha:equalityUpTo_k-thOrder}). We will return more extensively
to this difference between Fermat reals and SDG in Chapter \ref{cha:orderRelation}
about order relation on $\ER$. We only mention here that looking
at the simple Definition \ref{def:equalityInFermatReals}, the equality
\eqref{eq:productOfFirstOrderInfinitesimalIsZero} has an intuitively
clear meaning, and it is to preserve this intuition that we keep this
equality instead of changing completely the theory toward a less intuitive
one.

The next corollary solves the same problem of the previous one, but
starting from the hypotheses $h_{k}^{j_{k}}=0$:
\begin{cor}
\label{cor:suffCondToHaveProductZeroFromZeroPowers}If $h_{1},\dots,h_{n}\in D_{\infty}$
and $h_{k}^{j_{k}}=0$ for $j_{1},\dots,j_{n}\in\N$, then if $i_{1},\dots,i_{n}\in\N$,
we have\[
\sum_{k=1}^{n}\frac{i_{k}}{j_{k}}\ge1\then h_{1}^{i_{1}}\cdot\ldots\cdot h_{n}^{i_{n}}=0\]

\end{cor}
\noindent In fact if $h_{k}^{j_{k}}=0$, then $j_{k}>0$ and $h_{k}\in D_{j_{k}-1}$
by Theorem \ref{thm:propertyOfIdeal_D_a}, so the conclusion follows
from the previous corollary.

\noindent Finally, the latter corollary permits e.g. to pass from

\[
\forall h\in D_{p}^{n}:\quad f(h)=\sum_{\substack{i\in\N^{n}\\
|i|\le p}
}h^{i}\cdot a_{i}\]
 to\[
\forall h\in D_{q}^{n}:\quad f(h)=\sum_{\substack{i\in\N^{n}\\
|i|\le q}
}h^{i}\cdot a_{i}\]
 if $q<p$. In the previous formulas $D_{a}^{n}=D_{a}\times\ptind^{n}\times D_{a}$
and we have used the classical multi-indexes notations, e.g. $h^{i}=h_{1}^{i_{1}}\cdot\ldots\cdot h_{n}^{i_{n}}$
and $|i|=\sum_{k=1}^{n}i_{k}$.
\begin{cor}
\label{cor:suffCondToHaveProductZeroMulti-indexes}Let $p\in\N_{>0}$
and $h_{k}\in D_{p}$ for each $k=1,\dots,n$; $i\in\N^{n}$ and $h\in D_{\infty}^{n}$,
then\[
|i|>p\then h^{i}=0\]

\end{cor}
\noindent To prove it, we only have to apply Corollary \ref{cor:suffCondToHaveProductZeroFrom_D_j_k}:\[
\sum_{k=1}^{n}\frac{i_{k}}{p+1}=\frac{\sum_{k}i_{k}}{p+1}=\frac{|i|}{p+1}\ge\frac{p+1}{p+1}=1\]
 Let us note explicitly that the possibility to prove all these results
about products of powers of nilpotent infinitesimals is essentially
tied with the choice of little-oh polynomials in the definition of
the equivalence relation $\sim$ in Definition \ref{def:LittleOhPolynomials}.
Equally effective and useful results are not provable for the more
general family of nilpotent functions (see e.g. \cite{Gio3}).

\section{Identity principle for polynomials}

In this section we want to prove that if a polynomial $a_{0}+a_{1}x+a_{2}x^{2}+\dots+a_{n}x^{n}$
of $\ER$ is identically zero, then $a_{k}=0$ for all $k=0,\ldots,n$.
To prove this conclusion, it suffices to mean {}``identically zero''
as {}``equal to zero for every $x$ belonging to the extension of
an open subset of $\R$''. Therefore we firstly define what is this
extension.
\begin{defn}
\label{def:extensionOfSubsetsOfR}If $U$ is an open subset of $\R^{n}$,
then $\ext{U}:=\{x\in\ER^{n}\,|\,\st{x}\in U\}$. Here with the symbol
$\ext{\R^{n}}$ we mean $\ER^{n}:=\ER\times\ptind^{n}\times\ER$.
\end{defn}
\noindent We shall give further the general definition of the extension
functor $\ext(-)$; in these first chapters we only want to examine
some elementary properties of the ring $\ER$ that will be used later.

The identity principle for polynomials can now be stated in the following
way:
\begin{thm}
\label{thm:PIP}Let $a_{0},\dots,a_{n}\in\ER$ and $U$ be an open
neighborhood of $0$ in $\R$ such that\begin{equation}
a_{0}+a_{1}x+a_{2}x^{2}+\dots+a_{n}x^{n}=0\e{in}\ER\quad\forall x\in\ext{U}\label{eq:polyIdentity}\end{equation}

\noindent Then\[
a_{0}=a_{1}=\dots=a_{n}=0\e{in}\ER\]

\end{thm}
\noindent \textbf{\textit{\emph{Proof:}}} \textit{\emph{Because $U$
is an open neighborhood of $0$ in $\R$, we can always find $x_{1},\dots,x_{n+1}\in U$
such that $x_{i}\ne x_{j}$ for $i$, $j=1,\dots,n+1$ with $i\ne j$.
Hence from hypothesis \eqref{eq:polyIdentity} we have\[
a_{n}x_{k}^{n}+\dots+a_{1}x_{k}+a_{0}=0\e{in}\ER\quad\forall k=1,\dots,n+1\]
 That is, in vectorial form\[
(a_{n},\dots,a_{0})\cdot\left[\begin{array}{cccc}
x_{1}^{n} & x_{2}^{n} & \dots & x_{n+1}^{n}\\
x_{1}^{n-1} & x_{2}^{n-1} & \dots & x_{n+1}^{n-1}\\
\vdots\\
x_{1} & x_{2} & \dots & x_{n+1}\\
1 & 1 & \dots & 1\end{array}\right]=0\e{in}\ER\]
 \[
\]
 This matrix $V$ is a Vandermonde matrix, hence it is invertible\[
(a_{n},\dots,a_{0})\cdot V=\underline{0}\e{in}\ER^{n+1}\]
 \[
(a_{n},\dots,a_{0})\cdot V\cdot V^{-1}=\underline{0}\e{in}\ER^{n+1}\]
 hence $a_{k}=0$ in $\ER$ for every $k=0,\dots,n$.\qedNoNewLine}}

\textit{\emph{This theorem can be extended to polynomials with more
than one variable using recursively the previous theorem, one variable
per time:}}
\begin{cor}
\label{cor:PIPMultivariable}Let $a_{i}\in\ER$ for every $i\in\N^{n}$
with $|i|\le d$. Let $U$ be an open neighborhood of $\underline{0}$
in $\R^{n}$ such that\[
\sum_{\substack{i\in\N^{n}\\
|i|\le d}
}a_{i}x^{i}=0\quad\forall x\in\ext{U}\]
 Then\emph{\[
a_{i}=0\quad\forall i\in\N^{n}:\ \ |i|\le d\]
}
\end{cor}

\section{\textit{\emph{Invertible Fermat reals}}}

We can see more formally that to prove \eqref{eq:FunctionEqualTangent}
we cannot embed the reals $\R$ into a field but only into a ring,
necessarily containing nilpotent element. In fact, applying \eqref{eq:FunctionEqualTangent}
to the function $f(h)=h^{2}$ for $h\in D$, where $D\subseteq\ER$
is a given subset of $\ER$, we have\[
f(h)=h^{2}=f(0)+h\cdot f'(0)=0\quad\forall h\in D\]
 Hereby we have supposed the preservation of the equality $f'(0)=0$
from $\R$ to $\ER$. In other words, if $D$ and $f(h)=h^{2}$ verify
\eqref{eq:FunctionEqualTangent}, then necessarily each element $h\in D$
must be a new type of number whose square is zero. Of course in a
field the only subset $D$ verifying this property is $D=\{0\}$.

Because we cannot have property \eqref{eq:FunctionEqualTangent} and
a field at the same time, we need a sufficiently good family of cancellation
laws as substitutes. We will dedicate a full chapter of this work
to this problem, developing the notion of equality up to a $k$-th
order infinitesimal (see Section \ref{cha:equalityUpTo_k-thOrder}).
At present to prove the uniqueness of \eqref{eq:IdeaDF} we need the
following simplest form of these cancellation laws:
\begin{thm}
\label{thm:firstCancellationLaw} If $x\in\ER$ is a Fermat real and
$r$, $s\in\R$ are standard real numbers, then\[
x\cdot r=x\cdot s\text{ in }\ER\e{and}x\ne0\then r=s\]
\end{thm}
\begin{rem*}
\noindent As a consequence of this result, we can always cancel a
non zero Fermat real in an equality of the form $x\cdot r=x\cdot s$
where $r$, $s$ are \emph{standard} reals. This is obviously tied
with the uniqueness part of \eqref{eq:IdeaDF} and implies that formula
\eqref{eq:IdeaDF} uniquely identifies the first derivative in case
it is a standard real number.\smallskip{}

\end{rem*}
\noindent \textbf{Proof:} From the Definition \ref{def:equalityInFermatReals}
of equality in $\ER$ and from $x\cdot r=x\cdot s$ we have\[
\lim_{t\to0^{+}}\frac{x_{t}\cdot(r-s)}{t}=0\]
 But if we had $r\ne s$ this would implies $\lim_{t\to0^{+}}\frac{x_{t}}{t}=0$,
that is $x=0$ in $\ER$ and this contradicts the hypothesis $x\ne0$.\qedNoNewLine

The last result of this section takes its ideas from similar situations
of formal power series and gives also a formula to compute the inverse
of an invertible Fermat real.
\begin{thm}
\label{thm:formulaForInvertibleFermatdReals}Let $x=\st{x}+\sum_{i=1}^{n}\st{x_{i}}\cdot\diff{t_{a_{i}}}\in\ER$
be a Fermat real and its decomposition. Then\[
x\text{ is invertible}\]

\noindent if and only if $\st{x}\ne0$, and in this case\begin{equation}
\frac{1}{x}=\frac{1}{\st{x}}\cdot\sum_{j=0}^{+\infty}(-1)^{j}\cdot\left(\sum_{i=1}^{n}\frac{\st{x_{i}}}{\st{x}}\cdot\diff{t_{a_{i}}}\right)^{j}\label{eq:formulaforInvetibleFermatReal}\end{equation}

\end{thm}
\noindent In the formula \eqref{eq:formulaforInvetibleFermatReal}
we have to note that the series is actually a finite sum because any
$\diff{t_{a_{i}}}$ is nilpotent.
\begin{enumerate}
\item $(1+\diff{t_{2}})^{-1}=1-\diff{t_{2}}+\diff{t_{2}^{2}}-\diff{t_{2}^{3}}+\dots=1-\diff{t_{2}}+\diff{t}$
because $\diff{t_{2}^{3}}=0$ 
\item $(1+\diff{t_{3}})^{-1}=1-\diff{t_{3}}+\diff{t_{3}^{2}}-\diff{t_{3}^{3}}+\diff{t_{3}^{4}}-\dots=1-\diff{t_{3}}+\diff{t_{3}^{2}}-\diff{t}$ 
\end{enumerate}
\textbf{Proof:} If $x\cdot y=1$ for some $y\in\ER$, then, taking
the standard parts of each side we have $\st{x}\cdot\st{y}=1$ and
hence $\st{x}\ne0$. Vice versa the idea is to start from the series\[
\frac{1}{1+r}=\sum_{j=0}^{+\infty}(-1)^{j}\cdot r^{j}\quad\forall r\in\R:\ \ |r|<1\]
 and, intuitively, to define\begin{align*}
\left(\st{x}+\sum_{i}\st{x_{i}}\diff{t_{a_{i}}}\right)^{-1} & =\st{x}^{-1}\cdot\left(1+\sum_{i}\frac{\st{x_{i}}}{\st{x}}\diff{t_{a_{i}}}\right)^{-1}\\
 & =\st{x}^{-1}\cdot\sum_{j=0}^{+\infty}(-1)^{j}\cdot\left(\sum_{i}\frac{\st{x_{i}}}{\st{x}}\diff{t_{a_{i}}}\right)^{j}\end{align*}
 So let $y:=\st{x}^{-1}\cdot\sum_{j=0}^{+\infty}(-1)^{j}\cdot\left(\sum_{i}\frac{\st{x_{i}}}{\st{x}}\diff{t_{a_{i}}}\right)^{j}$
and $h:=x-\st{x}=\sum_{i}\st{x_{i}}\diff{t_{a_{i}}}\in D_{\infty}$
so that we can also write\[
y=\st{x}^{-1}\cdot\sum_{j=0}^{+\infty}(-1)^{j}\cdot\frac{h^{j}}{\st{x}^{j}}\]
 But $h\in\ER$ is a little-oh polynomial with $h(0)=0$, so it is
also continuous, hence for a sufficiently small $\delta>0$ we have\[
\forall t\in(-\delta,\delta):\ \ \left|\frac{h_{t}}{\st{x}}\right|<1\]
 Therefore\[
\forall t\in(-\delta,\delta):\ \ y_{t}=\frac{1}{\st{x}}\cdot\left(1+\frac{h_{t}}{\st{x}}\right)^{-1}=\frac{1}{\st{x}+h_{t}}=\frac{1}{x_{t}}\]
 From this equality it follows $x\cdot y=1$ in $\ER$ from Definition
\ref{def:equalityInFermatReals}.\qedNoNewLine

\section{The derivation formula}

Even if, in the following of this work, we will see several generalizations
of the derivation formula \eqref{eq:IdeaDF}, we want to give here
a proof of \eqref{eq:IdeaDF} because it has been the principal motivation
for the construction of the ring of Fermat reals $\ER$. Anyhow, before
considering the proof of the derivation formula, we have to understand
how to extend a given smooth function $f:\R\freccia\R$ to a certain
function $\ext f:\ER\freccia\ER$.
\begin{defn}
\label{def:extensionOfFunctions}Let $A$ be an open subset of $\R^{n}$,
$f:A\freccia\R$ a smooth function and $x\in\ext A$ then we define
\[
\ext f(x):=f\circ x\]

\end{defn}
\noindent This definition is well defined because we have seen (see
Subsection \ref{sub:ClosureOfLittle-ohPolyWRTSmoothFunctions}) that
little-oh polynomials are preserved by smooth functions, and because
the function $f$ is locally Lipschitz, so\[
\left|\frac{f(x_{t})-f(y_{t})}{t}\right|\le K\cdot\left|\frac{x_{t}-y_{t}}{t}\right|\quad\forall t\in(-\delta,\delta)\]
 for a sufficiently small $\delta$ and some constant $K$, and hence
if $x=y$ in $\ER$, then also $\ext{f}(x)=\ext{f}(y)$ in $\ER$.

The function $\ext f$ is an extension of $f$, that is \[
\ext f(r)=f(r)\quad{\rm {in}\quad\ER\quad{\rm \forall r\in\R,}}\]
 as it follows directly from the definition of equality in $\ER$
(i.e. Definition \ref{def:equalityInFermatReals}), thus we can still
use the symbol $f(x)$ both for $x\in\ER$ and $x\in\R$ without confusion.
After the introduction of the extension of smooth functions, we can
also state the following useful \emph{elementary transfer theorem}
for equalities, whose proof follows directly from the previous definitions:
\begin{thm}
\label{thm:elementaryTransferTheorem}Let $A$ be an open subset of
$\R^{n}$, and $\tau$, $\sigma:A\freccia\R$ be smooth functions.
Then \[
\forall x\in{\ext A}\pti\ext\tau(x)=\ext\sigma(x)\]
 iff \[
\forall r\in A\pti\tau(r)=\sigma(r).\]

\end{thm}
\noindent Now we will prove the derivation formula \eqref{eq:IdeaDF}.
\begin{thm}
\label{thm:DerivationFormula} Let $A$ be an open set in $\R$, $x\in A$
and $f:A\freccia\R$ a smooth function\emph{,} then\begin{equation}
\exists!\, m\in\R\;\forall h\in D\pti f(x+h)=f(x)+h\cdot m\label{eq:DerivationFormula}\end{equation}

\noindent In this case we have $m=f^{\prime}(x)$, where $f^{\prime}(x)$
is the usual derivative of $f$ at $x$. 
\end{thm}
\noindent \textbf{\textit{\emph{Proof:}}} \textit{\emph{Uniqueness
follows from the previous cancellation law Theorem \ref{thm:firstCancellationLaw},
indeed if $m_{1}\in\R$ and $m_{2}\in\R$ both verify \eqref{eq:DerivationFormula},
then $h\cdot m_{1}=h\cdot m_{2}$ for every $h\in D$. But there exists
a non zero first order infinitesimal, e.g. $\diff{t}\in D$, so from
Theorem \eqref{thm:firstCancellationLaw} it follows $m_{1}=m_{2}$.}}

To prove the existence part, take $h\in D$, so that $h^{2}=0$ in
$\ER$, i.e. $h_{t}^{2}=o(t)$ for $t\to0^{+}.$ But $f$ is smooth,
hence from its second order Taylor's formula we have\[
f(x+h_{t})=f(x)+h_{t}\cdot f'(x)+\frac{h_{t}^{2}}{2}\cdot f^{''}(x)+o(h_{t}^{2})\]
 But\[
\frac{o(h_{t}^{2})}{t}=\frac{o(h_{t}^{2})}{h_{t}^{2}}\cdot\frac{h_{t}^{2}}{t}\to0\e{for}t\to0^{+}\]
 so\[
\frac{h_{t}^{2}}{2}\cdot f^{''}(x)+o(h_{t}^{2})=o_{1}(t)\e{for}t\to0^{+}\]
 and we can write\[
f(x+h_{t})=f(x)+h_{t}\cdot f'(x)+o_{1}(t)\e{for}t\to0^{+}\]
 that is \[
f(x+h)=f(x)+h\cdot f'(x)\e{in}\ER\]
 and this proves the existence part because $f'(x)\in\R$.\qedNoNewLine

\noindent For example $e^{h}=1+h$, $\sin(h)=h$ and $\cos(h)=1$
for every $h\in D$.

Analogously we can prove the following infinitesimal Taylor's formula;
in its statement we use the usual multi-indexes notations (see e.g.
\cite{Pro2, Sil2}) and the notation $D_{n}^{d}:=D_{n}\times\ptind^{d}\times D_{n}$.
\begin{lem}
\label{thm:OrdinaryTaylorFor_nVariables}Let $A$ be an open set in
$\R^{d}$, $x\in A$, $n\in\N_{>0}$ and $f:A\freccia\R$ a smooth
function\emph{,} then\[
\forall h\in D_{n}^{d}\pti f(x+h)=\sum_{\substack{j\in\N^{d}\\
|j|\le n}
}\frac{h^{j}}{j!}\cdot\frac{\partial^{|j|}f}{\partial x^{j}}(x)\]

\end{lem}
\noindent For example $\sin(h)=h-\frac{h^{3}}{6}$ if $h\in D_{3}$
so that $h^{4}=0$.

It is possible to generalize several results of the present work to
functions of class $\mathcal{C}^{n}$ only, instead of smooth ones.
However it is an explicit purpose of this work to simplify statements
of results, definitions and notations, even if, as a result of this
searching for simplicity, its applicability will only hold for a more
restricted class of functions. Some more general results, stated for
$\mathcal{C}^{n}$ functions, but less simple can be found in \cite{Gio3}.

Note that $m=f^{\prime}(x)\in\R$, i.e. the slope is a standard real
number, and that we can use the previous formula with standard real
numbers $x$ only, and not with a generic $x\in\ER$, but we shall
remove this limitation in a subsequent chapter.

In other words we can say that the derivation formula \eqref{eq:IdeaDF}
allows us to differentiate the usual differentiable functions using
a language with infinitesimal numbers and to obtain from this an ordinary
function.

If we apply this theorem to the smooth function $p(r):=\int_{x}^{x+r}f(t)\diff{t}$,
for $f$ smooth, then we immediately obtain the following
\begin{cor}
\label{cor:derivationFormulaForIntegrals_1D}Let $A$ be open in $\R$,
$x\in A$ and $f:A\freccia\R$ smooth. Then \[
\forall h\in D\pti\int_{x}^{x+h}f(t)\diff{t}=h\cdot f(x).\]
 Moreover $f(x)\in\R$ is uniquely determined by this equality. 
\end{cor}
\noindent We close this section by introducing a very simple notation
useful to emphasize some equalities: if $h,k\in\ER$ then we say that
$\exists h/k$ iff $\exists!r\in\R\;:\; h=r\cdot k$, and obviously
we denote this $r\in\R$ with $h/k$. Therefore we can say, e.g.,
that \begin{align*}
f^{\prime}(x) & =\frac{f(x+h)-f(x)}{h} & {}\\
{} & {} & \forall h\in D_{\neq0}\\
f(x) & =\frac{1}{h}\cdot\int_{x}^{x+h}f(t)\diff{t}. & {}\end{align*}
 Moreover we can prove some natural properties of this {}``ratio'',
like the following one \[
\exists\frac{u}{v},\frac{x}{y}\e{and}vy\neq0\Longrightarrow\frac{u}{v}+\frac{x}{y}=\frac{uy+vx}{vy}\]

\begin{example}
\label{exa:ithStandardPartUsingRatio}Consider e.g. $x=1+2\diff{t_{3}}+\diff{t_{2}}+5\diff{t_{4/3}}$,
then using the previous ratio we can find a formula to calculate all
the coefficients of this decomposition. Indeed, let us consider first
the term $2\diff{t_{3}}$: if we multiply both sides by $\diff{t_{3/2}}$,
where\[
\frac{3}{2}=\frac{1}{1-\frac{1}{\omega(\diff{t_{3}})}}\]
 we obtain\[
(x-\st{x})\cdot\diff{t_{3/2}}=2\diff{t_{3}}\diff{t_{3/2}}+\diff{t_{2}}\diff{t_{3/2}}+5\diff{t_{4/3}}\diff{t_{3/2}}\]
 but $\diff{t_{3}}\diff{t_{3/2}}=\diff{t}$ whereas $\diff{t_{a}}\diff{t_{3/2}}=0$
if $a<3$, so\[
\frac{(x-\st{x})\diff{t_{3/2}}}{\diff{t}}=2\]
 Analogously we have\[
\frac{(x-\st{x}-2\diff{t_{3}})\diff{t_{2}}}{\diff{t}}=1\ee{and}\frac{(x-\st{x}-2\diff{t_{3}}-\diff{t_{2}})\diff{t_{4}}}{\diff{t}}=5\]
 where\[
2=\frac{1}{1-\frac{1}{\omega(\diff{t_{2}})}}\ee{and}4=\frac{1}{1-\frac{1}{\omega(\diff{t_{4/3}})}}\]
 Using the same idea we can prove the recursive formula\[
\alpha_{i+1}=\frac{1}{1-\frac{1}{\omega_{i+1}(x)}}\then\frac{\left(x-\st{x}-\sum_{k=1}^{i}x_{i}\diff{t_{\omega_{i}(x)}}\right)\cdot\diff{t_{\alpha_{i+1}}}}{\diff{t}}=x_{i+1}\]
 Finally, directly from the definition of decomposition it follows\[
\alpha\ne\frac{1}{1-\frac{1}{\omega_{i+1}(x)}}\,\,\forall i\then\frac{\left(x-\st{x}-\sum_{k=1}^{i}x_{i}\diff{t_{\omega_{i}(x)}}\right)\cdot\diff{t_{\alpha}}}{\diff{t}}=0\]

\noindent \[
\frac{\left(x-\st{x}-\sum_{k=1}^{i}x_{i}\diff{t_{\omega_{i}(x)}}\right)\cdot\diff{t_{\alpha}}}{\diff{t}}\ne0\then\alpha=\frac{1}{1-\frac{1}{\omega_{i+1}(x)}}\]
 so that all the terms of the decomposition of a Fermat real are uniquely
determined by these recursive formulas.\[
\]

\end{example}

\chapter{Equality up to $k$-th order infinitesimals\label{cha:equalityUpTo_k-thOrder}}

\section{Introduction}

As proved in Theorem \ref{thm:DerivationFormula}, the derivation
formula has several limitations that we are forced to avoid if we
want to obtain results like Stokes's theorem in the space $\ManInfty(M;N)$
of smooth functions between two smooth manifolds $M$, $N$. Let us
analyze the hypotheses of Theorem \ref{thm:DerivationFormula} so
as to motivate some generalizations:
\begin{enumerate}
\item {}``\emph{The point $x\in A$ is a standard real}''. This is the
hypothesis that can be more easily generalized. Indeed we can consider
that any general Fermat real $x\in\ER$ can be written as the sum
of its standard part $\st{x}\in\R$, and of its infinitesimal part
$k_{x}:=x-\st{x}\in D_{\infty}$. The infinitesimal part $k_{x}$
is of course nilpotent and hence, for $h\in D$, we can compute $f(x+h)=f(\st{x}+k_{x}+h)$
using the usual infinitesimal Taylor's formula or arbitrary order
(see Theorem \ref{thm:OrdinaryTaylorFor_nVariables}). We will follow
this idea in this chapter, but another solution is included in the
generalization of the following hypothesis. 
\item {}``\emph{The function $f:A\freccia\R$ is a standard smooth function}''.
As we already mentioned, not every function we are interested in is
of type $\ext{f}$, i.e. is the extension of a classical smooth function.
We already mentioned, as a simple example, the function $t\in\R_{\ge0}\mapsto\sin(h\cdot t)\in\ER$,
where $h\in D_{k}$ is an infinitesimal. More generally any function
of type $x\in\ER\mapsto\ext{g}(h,x)\in\ER$, where $g:\R^{2}\freccia\R$
is a given smooth function and $h\in\ER\setminus\R$ is a non standard
Fermat real, is not of type $\ext{f}$ for some $f$, because it can
happen that $ $$\ext{g}(h,r)\in\ER\setminus\R$ for a standard $r\in\R$
(whereas, of course, $\ext{f}(r)=f(r)\in\R$ for every $r\in\R$).
This implies that, on the one hand, we need a more general notion
of smooth function, surely including domains and codomains of type
$\varphi:\ext{U}\freccia\ER$, where $U\subseteq\R$ is open; on the
other hand we have to define a notion of derivative for this new type
of smooth function.\\
 We will solve this problem introducing the \emph{smooth incremental
ratio} (an idea that is mainly due to G. Reyes, see \cite{Koc}) ,
i.e. for every $x\in\ext{U}$, a function $h\in\ext{(-\delta_{x},\delta_{x})}\mapsto\varphi[x,h]\in\ER$
verifying\begin{equation}
\varphi(x+h)=\varphi(x)+\varphi[x,h]\cdot h\quad\forall h\in\ext{(-\delta_{x},\delta_{x})}\label{eq:smoothIncrementalRatio}\end{equation}
 and formalizing Fermat's method: $\varphi'(x):=\varphi[x,0]$. These
results are not usable for functions of the type $\varphi:D_{k}\freccia\ER$
which are not defined on the extension of a standard open set. This
problem is tied with the next hypothesis analyzed in this list. 
\item \label{enu:hypothesesOnDomain}{}``\emph{The domain of the smooth
function $f:A\freccia\R$ is an open set}''. Especially considering
spaces like spaces of functions, more general than locally flat spaces,
sometimes the more general results can be stated only in infinitesimal
domains like the above $D_{k}$. Examples are the existence and uniqueness
of the flux corresponding to a given vector field or the existence
and uniqueness of the exterior derivative of an $n$-form. For this
reason we will have to define some notion of derivative for functions
of type $\psi:D_{k}\freccia\ER$. At first sight, the definition of
derivative for this type of function may seem an easy goal. In fact,
intuitively, a function of this type can be thought as some type of
polynomial of degree $k\in\N$. The problem is due to the fact that,
in our setting, the derivation formula does not determine uniquely
the coefficients of this polynomial. Indeed we know that in $\ER$
we have $h\cdot k=0$ for every first order infinitesimal $h$, $k\in D$,
so both the coefficients $m_{1}=1+k$, for a fixed $k\in D$, and
$m_{2}=1$, verify, for every $h\in D$, the derivation formula $f(h)=f(0)+h\cdot m_{i}$
if $f(h)=h$.\\
 We want to underline, even if it will be formally clear only later
in this work, that here we do not have a problem of existence, but
of uniqueness only. In other words, e.g. for a function of type $\psi:D\freccia\ER$
there always exists an $m\in\ER$ such that $\psi(h)=\psi(0)+h\cdot m$
for every $h\in D$, but this coefficient $m\in\ER$ is not uniquely
determined by this formula. We can tackle this problem in several
ways. For example we can try to find another formula that uniquely
identifies what we intuitively think of as the derivative of $\psi:D_{k}\freccia\ER$
at $0$. The idea of the smooth incremental ratio \eqref{eq:smoothIncrementalRatio}
goes in this direction. Anyhow, in this work we followed another idea:
because we only have a uniqueness and not an existence problem, we
shall try to define precisely {}``\emph{what is the simplest $m\in\ER$
that verifies the derivation formula,} \emph{and we will call {}``derivative''
this simplest coefficient}''. E.g. among $m_{1}=1+k$ and $m_{2}=1$
in the previous example, the simplest one will surely be $m_{2}=1$
and hence we shall have $\psi'(0)=1$. This chapter is devoted to
the development of these ideas. Indeed $m_{2}=1$ is simpler than
$m_{1}=1+k$ in the sense that $m_{2}$ is $m_{1}$ up to second order
infinitesimals. 
\end{enumerate}
Let us start from the hypothesis\[
m\in\ER\e{and}\forall h\in D\pti h\cdot m=0\]
 and try to derive some necessary condition on $m\in\ER$ based on
the idea that {}``because we have a product with $h\in D$, some
infinitesimal in the decomposition of $m$ will give zero if multiplied
by $h$, so not every infinitesimal in the decomposition of $m$ is
really useful to obtain the final value of the product $h\cdot m$''.
In fact let\begin{equation}
m=\st{m}+\sum_{i=1}^{N}\st{m_{i}}\cdot\diff{t_{a_{i}}}\label{eq:decompositionOf_m}\end{equation}
 be the decomposition of $m$, and\begin{equation}
h=\sum_{j=1}^{k}\st{h_{j}}\cdot\diff{t_{b_{j}}}\label{eq:decompositionOf_h}\end{equation}
 be the decomposition of a generic $h\in D$. Then\begin{equation}
h\cdot m=\sum_{j=1}^{k}\st{m}\st{h_{j}}\diff{t_{b_{j}}}+\sum_{j=1}^{k}\sum_{i=1}^{N}\st{h_{j}}\st{m_{i}}\diff{t_{\frac{a_{i}b_{j}}{a_{i}+b_{j}}}}\label{eq:product_h_m}\end{equation}
 But $h\in D$, hence $b_{j}\le\omega(h)=b_{1}<2$ so that if we had
$a_{i}\le2$ we would have\[
\frac{1}{b_{j}}+\frac{1}{a_{i}}>\frac{1}{2}+\frac{1}{2}=1\]
 and hence $\frac{a_{i}b_{j}}{a_{i}+b_{j}}<1$ and $\diff{t_{\frac{a_{i}b_{j}}{a_{i}+b_{j}}}}=0$.
Therefore we can write

\begin{align}
\forall h & \in D\pti h\cdot m=\sum_{j=1}^{k}\st{m}\st{h_{j}}\diff{t_{b_{j}}}+\sum_{j=1}^{k}\sum_{\substack{i=1\\
a_{i}>2}
}^{N}\st{h_{j}}\st{m_{i}}\diff{t_{\frac{a_{i}b_{j}}{a_{i}+b_{j}}}}=\nonumber \\
= & h\cdot\left(\st{m}+\sum_{\substack{i=1\\
a_{i}>2}
}^{N}\st{m_{i}}\cdot\diff{t_{a_{i}}}\right)\label{eq:remainsTheSimplest}\end{align}
 Looking at \eqref{eq:remainsTheSimplest}, we can say that {}``in
a product of type $h\cdot m$, with $h\in D$, only sufficiently big
infinitesimals ($a_{i}>2$) in the decomposition of $m$ will survive''.
In other words, {}``all the infinitesimals of order less or equal
$2$ are useless to define the value of the product $h\cdot m$''.

The quantity\[
\iota_{2}(m):=\st{m}+\sum_{\substack{i=1\\
a_{i}>2}
}^{N}\st{m_{i}}\cdot\diff{t_{a_{i}}}\]
 is exactly the number $m$ up to second order infinitesimals%
\footnote{Remember that for Fermat reals, the greater is the order and the bigger
the infinitesimal has to be thought, see Remark \ref{rem:differentMeaningForStandardAndPotentialOrder}.%
}, in the sense that it is obtained from $m$ neglecting all the {}``small''
infinitesimals $\diff{^{i}m}=\st{m_{i}}\diff{t_{a_{i}}}$ of order
$\omega(\st{m_{i}}\diff{t_{a_{i}}})=a_{i}\le2$. In the example mentioned
at the item \ref{enu:hypothesesOnDomain}., where $k\in D$, we have
$\iota_{2}(m_{1})=\iota_{2}(1+k)=1$ and, indeed, among all the Fermat
reals $m\in\ER$ that verify the derivation formula, $\iota_{2}(m)$
will be our candidate for the definition of {}``the simplest Fermat
real that verifies the derivation formula''. In fact, the formula
\eqref{eq:remainsTheSimplest} can be written as\[
\forall h\in D\pti h\cdot m=h\cdot\iota_{2}(m)\]
 and it can be interpreted intuitively saying {}``among all the numbers
$m$ that gives the same value of the product $h\cdot m$, the number
$\iota_{2}(m)$ is the simplest one because it contains the minimal
information, neglecting all the useless infinitesimals, i.e. not useful
to define the value of the product $h\cdot m$''.

\section{Equality up to $k$-th order infinitesimals}

The considerations of the previous section give us sufficient heuristic
motivations to define:
\begin{defn}
\label{def:upTo_k-thOrderInfinitesimals} Let $m=\st{m}+\sum_{i=1}^{N}\st{m_{i}}\cdot\diff{t}_{\omega_{i}(m)}$
be the decomposition of $m\in\ER$ and $k\in\R_{\ge0}\cup\{\infty\}$,
then 
\begin{enumerate}
\item ${\displaystyle \iota_{k}m:=\iota_{k}(m):=\st{m}+\sum_{\substack{i=1\\
\omega_{i}(m)>k}
}^{N}\st{m_{i}}\cdot\diff{t_{\omega_{i}(m)}}}$
\item $\ER_{k}:=\left\{ \iota_{k}m\,|\, m\in\ER\right\} $.
\end{enumerate}
Finally if $x$, $y\in\ER$, we will say $x=_{k}y$ iff $\iota_{k}x=\iota_{k}y$
in $\ER$, and we will read it as \emph{$x$ is equal to $y$ up to
$k$-th order infinitesimals.} \end{defn}
\begin{rem}
\noindent \label{rem:compositionsOf_iota_k} Firstly note that if
$0\le k<1$ then the condition $\omega_{i}(m)>k$ is trivial because
we always have that $\omega_{i}(m)\ge1$. Hence\[
\iota_{0}m=m\e{and}\ER_{0}:=\ER\]
 Moreover $\iota_{\infty}m=\st{m}$ and $\ER_{\infty}:=\R$.

The first simple property we can note about $\iota_{k}$ is that $\iota_{j}(\iota_{k}x)=\iota_{j\vee k}(x)$
(recalling that $j\vee k:=\max(j,k))$ so that we have, e.g., $\iota_{j}(\iota_{k}x)=\iota_{k}(\iota_{j}x)$
and $\iota_{k}x=x$ for every $x\in\ER_{k}$. Moreover we have the
following chain of inclusions \begin{equation}
\R=\ER_{\infty}\subseteq\ldots\subseteq\ER_{3}\subseteq\ER_{2}\subseteq\ER_{1}\subseteq\ER_{0}=\ER\label{eq:sequenceOf_ER_k}\end{equation}
 In fact if $z\in\ER_{k}$, we can write $z=\iota_{k}m$ for some
$m\in\ER$; but $\iota_{j}(\iota_{k}m)=\iota_{k}(m)=z$ if $j\le k$
and hence $z$ is also of type $\iota_{j}(m')$ for $m'=\iota_{k}m$
and so $z\in\ER_{j}$. The intuitive meaning of \eqref{eq:sequenceOf_ER_k}
is clear: the more infinitesimals we neglect from a Fermat real $m\in\ER$
and the less terms will remain in the decomposition of $m$; continuing
in this way, only the standard part $\st{m}$ remains.

\noindent We start the study of $\iota_{k}$ considering the relationships
between this operation and the algebraic operations on $\ER$.\end{rem}
\begin{thm}
\label{thm:ER_k_andAlgebraicOperations} Let $x$, $y\in\ER$ and
$k\in\R_{\ge1}$,then 
\begin{enumerate}
\item \label{enu:iotaForGenericSumOfDifferentials}If $x=r+\sum_{h=1}^{M}\gamma_{h}\cdot\diff{t_{c_{h}}}$
in $\ER$ (not necessarily the decomposition of $x$), with $r$,
$\gamma_{h}\in\R$ and $c_{h}\in\R_{\ge1}$, then $\iota_{k}x=r+\sum\limits _{h:c_{h}>k}\gamma_{h}\cdot\diff{t_{c_{h}}}$ 
\item \label{enu:iotaPreservesSum}$\iota_{k}(x+y)=\iota_{k}x+\iota_{k}y$ 
\item \label{enu:iotaPreservesZeros}$\iota_{k}0=0$ 
\item \label{enu:iotaOmogeneous}$\iota_{k}(r\cdot x)=r\cdot\iota_{k}x\quad\forall r\in\R$ 
\item \label{enu:iotaPreservesProductUpTo_k-thOrder}$\iota_{k}(x\cdot y)=\iota_{k}(\iota_{k}x\cdot\iota_{k}y)$,
that is $x\cdot y=_{k}\iota_{k}x\cdot\iota_{k}y$ 
\item \label{enu:ERupTo_k-isARing} The relation $=_{k}$ is an equivalence
relation and the quotient set $\ER_{\sss{=_{k}}}:=\ER/=_{k}$ is a
ring with respect to pointwise operations 
\end{enumerate}
\end{thm}
\noindent \textbf{Proof:} To prove \emph{\ref{enu:iotaForGenericSumOfDifferentials}}.
we can consider that if $x=r+\sum_{h=1}^{M}\gamma_{h}\cdot\diff{t_{c_{h}}}$,
then\begin{align*}
x & =r+\sum_{h=1}^{M}\gamma_{h}\cdot\diff{t_{c_{h}}}=\\
 & =\st{x}+\sum_{q:q\in\{c_{j}\,|\, j=1,\ldots,M\}}\diff{t_{q}}\cdot\sum\{\gamma_{h}\,|\, h=1,\ldots,M\,,\, c_{h}=q\}\end{align*}
 where we have summed all the addends $\gamma_{h}\diff{t_{c_{h}}}$
having the same order $c_{h}=q$. Now call $\{q_{1},\ldots,q_{P}\}:=\{c_{j}\,|\, j=1,\ldots,M\}$
the distinct elements of the set of all the $c_{j}$, and $\bar{\gamma}_{a}:=\sum\{\gamma_{h}\,|\, h=1,\ldots,M\,,\, c_{h}=q_{a}\}$.
Hence $x=\st{x}+\sum_{a=1}^{P}\bar{\gamma}_{a}\diff{t_{q_{a}}}$,
and we can suppose that every $\bar{\gamma}_{a}\ne0$. Recalling the
construction of the decomposition of a Fermat real (see the existence
proof of Theorem \ref{thm:existenceUniquenessDecomposition}) we can
state that $P=N$, where $N$ is the number of addends in the decomposition
of $x$, and that permuting the addends of this sum we obtain the
decomposition of $x$, i.e. for a suitable permutation $\sigma$ of
$\{1,\ldots,N\}$ we have\[
q_{\sigma(i)}=\omega_{i}(x)\ee{and}\bar{\gamma}_{\sigma(i)}=\st{x_{i}}\]
 or, in other words, we can say that\[
x=\st{x}+\sum_{i=1}^{N}\bar{\gamma}_{\sigma(i)}\diff{t_{q_{\sigma(i)}}}\quad\text{is the decomposition of }x.\]
 Therefore, by Definition \ref{def:upTo_k-thOrderInfinitesimals}\begin{align*}
\iota_{k}x & =\st{x}+\sum_{\substack{i=1\\
q_{\sigma(i)}>k}
}^{N}\bar{\gamma}_{\sigma(i)}\diff{t_{q_{\sigma(i)}}}=\st{x}+\sum_{\substack{a=1\\
q_{a}>k}
}^{N}\bar{\gamma}_{a}\diff{t_{q_{a}}}=\\
 & =\st{x}+\sum_{\substack{a=1\\
q_{a}>k}
}^{N}\diff{t_{q_{a}}}\sum\{\gamma_{h}\,|\, h=1,\ldots,M\,,\, c_{h}=q_{a}\}=\\
 & =\st{x}+\sum_{\substack{q:q\in\{c_{j}\,|\, j=1,\ldots,M\}\\
q>k}
}\diff{t_{q}}\cdot\sum\{\gamma_{h}\,|\, h=1,\ldots,M\,,\, c_{h}=q\}=\\
 & =r+\sum_{\substack{h=1\\
c_{h}>k}
}^{M}\gamma_{h}\diff{t_{c_{h}}}.\end{align*}
 \emph{\ref{enu:iotaPreservesSum}}.) \ \ \ We consider the decompositions
of $x$ and $y$, so that we have that\begin{equation}
\iota_{k}x+\iota_{k}y=\st{x}+\sum_{i:\omega_{i}(x)>k}\st{x_{i}}\cdot\diff{t_{\omega_{i}(x)}}+\st{y}+\sum_{j:\omega_{j}(y)>k}\st{y_{j}}\cdot\diff{t_{\omega_{j}(y)}}\label{eq:sumOfIota}\end{equation}
 On the other hand we have\[
x+y=\st{x}+\st{y}+\sum_{i}\st{x_{i}}\diff{t_{\omega_{i}(x)}}+\sum_{j}\st{y_{j}}\diff{t_{\omega_{j}(y)}}\]
 From this and from the previous result \emph{\ref{enu:iotaForGenericSumOfDifferentials}}.
we have that\[
\iota_{k}(x+y)=\st{x}+\st{y}+\sum_{i:\omega_{i}(x)>k}\st{x_{i}}\diff{t_{\omega_{i}(x)}}+\sum_{j:\omega_{j}(y)>k}\st{y_{j}}\diff{t_{\omega_{j}(y)}}=\iota_{k}x+\iota_{k}y\]
 Property \emph{\ref{enu:iotaPreservesZeros}}. is a general consequence
of \emph{\ref{enu:iotaPreservesSum}}. for $x=y=0$.

\smallskip{}

\noindent \emph{\ref{enu:iotaOmogeneous}}.) \ \ \ We multiply
$x$ by $r\in\R$ obtaining\[
r\cdot x=r\cdot\st{x}+\sum_{i=1}^{N}r\cdot\st{x_{i}}\diff{t_{\omega_{i}(x)}}\]
 so that, once again from \emph{\ref{enu:iotaForGenericSumOfDifferentials}}.,
we have\[
\iota_{k}(r\cdot x)=r\cdot\st{x}+\sum_{i:\omega_{i}(x)>k}r\cdot\st{x_{i}}\diff{t_{\omega_{i}(x)}}=r\cdot\iota_{k}x\]
 \emph{\ref{enu:iotaPreservesProductUpTo_k-thOrder}}.) \ \ \ Let
us consider the product of the decompositions of $x$ and $y$ and
let $a_{i}:=\omega_{i}(x)$, $b_{j}:=\omega_{j}(y)$ for simplicity,
then we have\[
x\cdot y=\st{x}\,\st{y}+\sum_{j}\st{x}\,\st{y_{j}}\diff{t_{b_{j}}}+\sum_{i}\st{y}\,\st{x_{i}}\diff{t_{a_{i}}}+\sum_{i,j}\st{x_{i}}\,\st{y_{j}}\diff{t_{\frac{a_{i}b_{j}}{a_{i}+b_{j}}}}\]
 Hence from \emph{\ref{enu:iotaForGenericSumOfDifferentials}}. we
have\begin{align*}
\iota_{k}(x\cdot y) & =\st{x}\,\st{y}+\sum_{j:b_{j}>k}\st{x}\,\st{y_{j}}\diff{t_{b_{j}}}+\sum_{i:a_{i}>k}\st{y}\,\st{x_{i}}\diff{t_{a_{i}}}+\\
 & +\sum\left\{ \st{x_{i}}\,\st{y_{j}}\diff{t_{\frac{a_{i}b_{j}}{a_{i}+b_{j}}}}\,|\,\frac{a_{i}b_{j}}{a_{i}+b_{j}}>k\right\} \end{align*}
 On the other hand we have\begin{align}
\iota_{k}x\cdot\iota_{k}y & =\st{x}\,\st{y}+\sum_{j:b_{j}>k}\st{x}\,\st{y_{j}}\diff{t_{b_{j}}}+\sum_{i:a_{i}>k}\st{y}\,\st{x_{i}}\diff{t_{a_{i}}}+\nonumber \\
 & +\sum_{\substack{i:a_{i}>k\\
j:b_{j}>k}
}\st{x_{i}}\,\st{y_{j}}\diff{t_{\frac{a_{i}b_{j}}{a_{i}+b_{j}}}}\label{eq:productOfIota_k}\end{align}
 and applying \emph{\ref{enu:iotaForGenericSumOfDifferentials}}.
to \eqref{eq:productOfIota_k} we get\begin{align}
\iota_{k}(\iota_{k}x\cdot\iota_{k}y) & =\st{x}\,\st{y}+\sum_{j:b_{j}>k}\st{x}\,\st{y_{j}}\diff{t_{b_{j}}}+\sum_{i:a_{i}>k}\st{y}\,\st{x_{i}}\diff{t_{a_{i}}}+\nonumber \\
 & +\sum\left\{ \st{x_{i}}\,\st{y_{j}}\diff{t_{\frac{a_{i}b_{j}}{a_{i}+b_{j}}}}\,|\, a_{i}>k,\, b_{j}>k,\,\frac{a_{i}b_{j}}{a_{i}+b_{j}}>k\right\} \label{eq:iota_kOfProductOfIota_k}\end{align}
 So it suffices to prove that the set of Fermat reals in the third
summation sign both in \eqref{eq:productOfIota_k} and \eqref{eq:iota_kOfProductOfIota_k}
are equal. But immediately we can see that the set in \eqref{eq:iota_kOfProductOfIota_k}
is a subset of the set of numbers in \eqref{eq:productOfIota_k}.
For the opposite inclusion we have\begin{equation}
\frac{a_{i}b_{j}}{a_{i}+b_{j}}>k\then\frac{1}{a_{i}}+\frac{1}{b_{j}}<\frac{1}{k}\label{eq:a_i-b_j}\end{equation}
 but $\frac{1}{a_{i}}<\frac{1}{a_{i}}+\frac{1}{b_{j}}$ because $b_{j}=\omega_{j}(y)>0$
and hence from \eqref{eq:a_i-b_j} we obtain $a_{i}>k$. Analogously
we can prove that $b_{j}>k$ so that the two sets of Fermat reals
are equal.

\noindent \emph{\ref{enu:ERupTo_k-isARing}}.)\ \ \ We have only
to prove that the ring operations on the quotient set $\ER/=_{k}$,
are well defined, i.e. that\begin{align}
x=_{k}x'\e{and}y=_{k}y'\then x+x'=_{k}y+y'\label{eq:equalityUpTo_kAndSum}\\
x=_{k}x'\e{and}y=_{k}y'\then x\cdot x'=_{k}y\cdot y'\label{eq:equalityUpTo_kAndProduct}\end{align}
 Indeed if $x=_{k}x'$ and $y=_{k}y'$, then $\iota_{k}x=\iota_{k}x'$
and $\iota_{k}y=\iota_{k}y'$ (obviously these equalities have to
be understood in $\ER$), so $\iota_{k}(x)\cdot\iota_{k}(y)=\iota_{k}(x')\cdot\iota_{k}(y')$.
Applying $\iota_{k}$ to both sides we obtain $\iota_{k}\left(\iota_{k}(x)\cdot\iota_{k}(y)\right)=\iota_{k}\left(\iota_{k}(x')\cdot\iota_{k}(y')\right)$
so that from \emph{\ref{enu:iotaPreservesProductUpTo_k-thOrder}}.
we have $\iota_{k}(x\cdot y)=\iota_{k}(x'\cdot y')$, that is $x\cdot y=_{k}x'\cdot y'$.
Analogously, using \emph{\ref{enu:iotaPreservesSum}}., we can prove
\eqref{eq:equalityUpTo_kAndSum}.\qedNoNewLine
\begin{rem}
\label{rem:equalityInER_k}If $m\in\ER_{k}$ and $m=_{k}0$ then we
can write $m=\iota_{k}n$ for some $n\in\ER$; but, on the other hand,
$\iota_{k}m=\iota_{k}0=0$ in $\ER$ because $m=_{k}0$. But $\iota_{k}m=\iota_{k}(\iota_{k}n)=\iota_{k}n=m$,
so we can finally deduce that $m$ must be zero in $\ER$, i.e. $m=0$.
Therefore:\[
m\in\ER_{k}\e{and}m=_{k}0\then m=0\e{in}\ER\]
 This can also be restated saying that the notion of equality up to
$k$-th order infinitesimals, i.e. the equivalence relation $=_{k}$,
is trivial in $\ER_{k}$, i.e. if $m$, $n\in\ER_{k}$ and $m=_{k}n$,
then $m=n$ in $\ER$.\\
 Moreover we can also state \eqref{eq:equalityUpTo_kAndSum} and
\eqref{eq:equalityUpTo_kAndProduct} saying that if we work with equality
up to $k$-th order infinitesimals, that is with the equivalence relation
$=_{k}$, we can always use ring operations sum and product of $\ER$
and this will preserve the equality $=_{k}$. \end{rem}
\begin{example*}
\noindent Whereas property \emph{\ref{enu:iotaPreservesSum}}. says
that $\ER_{k}$ is closed with respect to sums, in general it is not
closed with respect to products. Indeed let

\[
x=\diff{t_{3}}=y\]
 then $\iota_{2}x=\iota_{2}y=\diff{t_{3}}$ and $\iota_{2}x\cdot\iota_{2}y=\left(\diff{t_{3}}\right)^{2}=\diff{t_{3/2}}$.
On the other hand $x\cdot y=\diff{t_{3/2}}$ and $\iota_{2}(x\cdot y)=0$,
so $\iota_{2}(x\cdot y)\ne\iota_{2}(x)\cdot\iota_{2}(y)$. This counterexample
exhibits why we stated the relationships between $\iota_{k}$ and
the product as in \emph{\ref{enu:iotaPreservesProductUpTo_k-thOrder}}
of Theorem \ref{thm:ER_k_andAlgebraicOperations}. \bigskip{}

\end{example*}
We close this section with a theorem that states some properties of
the order of $\iota_{k}x-\iota_{j}x$. The starting idea is roughly
the following: with $\iota_{k}x$ we {}``delete'' in the decomposition
of $x$ all the infinitesimals of order less or equal to $k$; we
do the same with $\iota_{j}x$, so if $j>k$ in the difference $\iota_{k}x-\iota_{j}x$
there will remain only infinitesimals of order between $k$ and $j$.
\begin{thm}
\label{thm:orderAndIota_k}Let $x\in\ER$ and $j$, $k\in\R_{\ge1}$,
with $j>k$, then 
\begin{enumerate}
\item \label{enu:orderOfDifferenceOfIota}$k<\omega(\iota_{k}x-\iota_{j}x)\le j$
and hence $\iota_{k}x-\iota_{j}x\in D_{j}$ 
\item \label{enu:orderOfIota_k}$k<\omega(\iota_{k}x)$ 
\item \label{enu:orderOf_x_MinusIota_x}$\omega(x-\iota_{j}x)\le j$ 
\item \label{enu:differenceOfIotaAndProductInD_p}$\forall h\in D_{\frac{1}{j-1}}\pti h\cdot(\iota_{k}x-\iota_{j}x)=0$ 
\end{enumerate}
\end{thm}
\noindent \textbf{Proof:} To prove \ref{enu:orderOfDifferenceOfIota}\emph{.}
let $x=r+\sum_{i=1}^{N}\alpha_{i}\diff{t_{a_{i}}}$ be the decomposition
of $x$, then\begin{align}
\iota_{k}x-\iota_{j}x & =r+\sum_{i:a_{i}>k}\alpha_{i}\diff{t_{a_{i}}}-r-\sum_{i:a_{i}>j}\alpha_{i}\diff{t_{a_{i}}}=\nonumber \\
 & =\sum_{i:k<a_{i}\le j}\alpha_{i}\diff{t_{a_{i}}}+\sum_{i:a_{i}>j}\alpha_{i}\diff{t_{a_{i}}}-\sum_{i:a_{i}>j}\alpha_{i}\diff{t_{a_{i}}}=\nonumber \\
 & =\sum_{i:k<a_{i}\le j}\alpha_{i}\diff{t_{a_{i}}}\label{eq:decompositionOfDifferenceOfIota}\end{align}
 So $\sum_{i:k<a_{i}\le j}\alpha_{i}\diff{t_{a_{i}}}$ is the decomposition
of $\iota_{k}x-\iota_{j}x$ and its order is given by $\omega(\iota_{k}x-\iota_{j}x)=a_{p}$,
where $p$ is the smallest index $i=1,\ldots,N$ in the decomposition
\eqref{eq:decompositionOfDifferenceOfIota}, i.e. $p:=\min\{i=1,\ldots,N\,|\, k<a_{i}\le j\}$.
Therefore \textit{$k<\omega(\iota_{k}x-\iota_{j}x)=a_{p}\le j$.}

\noindent Property \emph{\ref{enu:orderOfIota_k}}. can be proved
exactly as the previous \emph{\ref{enu:orderOfDifferenceOfIota}}.
but with $j=+\infty$.

\noindent Property \emph{\ref{enu:orderOf_x_MinusIota_x}}. is simply
property \emph{\ref{enu:orderOfDifferenceOfIota}}. with $k=0$.

\noindent \emph{\ref{enu:differenceOfIotaAndProductInD_p}}.)\ \ \ If
$h\in D_{\frac{1}{j-1}}$, then $\omega(h)<\frac{1}{j-1}+1=\frac{j}{j-1}$.
Let us analyze the product $h\cdot(\iota_{k}x-\iota_{j}x)$:\[
\frac{1}{\omega(h)}+\frac{1}{\omega(\iota_{k}x-\iota_{j}x)}>\frac{j-1}{j}+\frac{1}{j}=1\]
 Hence from \ref{thm:productOfPowers} the conclusion follows.\qedNoNewLine

As a consequence of the previous property \emph{\ref{enu:orderOfIota_k}}
of Theorem \ref{thm:orderAndIota_k} we have the following simple
cancellation law.
\begin{cor}
Let $m\in\ER_{k}$, with $k\ge1$, and $h\in D_{\infty}$ with $k+\omega(h)\le k\cdot\omega(h)$,
then\[
m\cdot h=0\then m=0\]

\end{cor}
\noindent \textbf{Proof:} First we note that $h\ne0$, because otherwise
we had $\omega(h)=\omega(0)=0$ and $k\le0$ from the hypothesis $k+\omega(h)\le k\cdot\omega(h)$.
Secondly, from the hypothesis $m\cdot h=0$ we immediately have $\st{m}=0$,
so that $m\in D_{\infty}$. Ad absurdum, if we had $m\ne0$, then
we would also have\[
\frac{1}{\omega(m)}+\frac{1}{\omega(h)}<\frac{1}{k}+\frac{1}{\omega(h)}\]
indeed $\omega(m)>k$ because $m\in\ER_{k}$ and the previous Theorem
\ref{thm:orderAndIota_k}. But $\frac{1}{k}+\frac{1}{\omega(h)}=\frac{k+\omega(h)}{k\cdot\omega(h)}\le1$
by hypothesis and therefore $m\cdot h\ne0$ by Theorem \ref{thm:productOfPowers},
in contradiction with the hypothesis.\qedNoNewLine

For example if we take $m$ as above and consider the infinitesimal
$h=\diff{t_{j}}$ with $2\le j\le k$, then\[
m\cdot\diff{t_{j}}=0\then m=0\]
 indeed $\frac{1}{k}+\frac{1}{\omega(\diff{t_{j}})}=\frac{1}{k}+\frac{1}{j}\le\frac{2}{j}\le1$,
i.e. $k+\omega(\diff{t_{j}})\le k\cdot\omega(\diff{t_{j}})$.

\section{Cancellation laws up to $k$-th order infinitesimals}

The goal of this section is to find for what infinitesimals $h\in D_{\infty}$
and for what power $j\in\N$ and order $k\in\R_{\ge1}$ we have $h^{j}\cdot m=h^{j}\cdot\iota_{k}m$.
We recall that we started this chapter motivating the definition of
$\iota_{k}x$ starting from the property\[
\forall h\in D\pti h\cdot m=h\cdot\iota_{2}m\]
 In this section we want to generalize this property. We will see
that, as a consequence of this generalization, we will obtain a cancellation
law up to $k$-th order infinitesimals of the form\begin{align}
\text{If } & \quad\forall h\in D_{\alpha_{1}}\times\cdots\times D_{\alpha_{n}}\pti h^{j}\cdot m=0\label{eq:cancellationLawUpTo_kOrder}\\
\text{then} & \quad m=_{k}0\nonumber \end{align}
 and hence a general Taylor's formula for smooth functions of the
type $f:\R^{n}\freccia\R$ with \emph{independent infinitesimals increments},
that is a formula useful to compute with a polynomial a term like
$f(x_{1}+h_{1},\ldots,x_{n}+h_{n})$, with $(h_{1},\ldots,h_{n})\in D_{\alpha_{1}}\times\cdots\times D_{\alpha_{n}}$,
i.e. with infinitesimal increments in general of different orders.

We shall use the classical multi-indexes notations (see e.g. \cite{Pro2}
) frequently used in the study of several variables functions. E.g.
in \eqref{eq:cancellationLawUpTo_kOrder} we already used $h^{j}:=h_{1}^{j_{1}}\cdot\ldots\cdot h_{n}^{j_{n}}$.

\noindent We start proving one simple lemma that will be useful in
the following.
\begin{lem}
\label{lem:producOfInfinitesimalAndIota}Let $m\in\ER$, $k\in\R_{\ge1}$
and $h\in D_{\infty}$ such that\begin{equation}
\frac{1}{k}+\frac{1}{\omega(h)}>1\label{eq:firstConditionForSolution}\end{equation}

\noindent then\begin{equation}
h\cdot m=h\cdot\iota_{k}m\label{eq:firstConclusion}\end{equation}

\end{lem}
\noindent Condition \eqref{eq:firstConditionForSolution} is not difficult
to foresee if we want to obtain \eqref{eq:firstConclusion}, because
it implies, as we will see in the following proof, that all the infinitesimals,
in the decomposition of $m$, having an order which is less than or
equal to $k$, multiplied by $h$ will give 0 (compare property \eqref{eq:firstConditionForSolution}
with Theorem \ref{thm:productOfPowers}).

\bigskip{}

\noindent \textbf{Proof:} Let $h=\sum_{p=1}^{N}\beta_{p}\diff{t_{b_{p}}}$
resp. $m=r+\sum_{i=1}^{M}\alpha_{i}\diff{t_{a_{i}}}$ be the decompositions
of $h$ and $m$. Then\begin{equation}
h\cdot m=\sum_{p=1}^{n}r\beta_{p}\diff{t_{b_{p}}}+\sum_{i,p}\alpha_{i}\beta_{p}\diff{t_{\frac{a_{i}b_{p}}{a_{i}+b_{p}}}}\label{eq:product_h_mToHaveIota_k}\end{equation}
 But if $a_{i}\le k$, then\[
\frac{1}{a_{i}}+\frac{1}{b_{p}}\ge\frac{1}{k}+\frac{1}{\omega(h)}>1\]
 So if $a_{i}\le k$, then $\frac{a_{i}b_{p}}{a_{i}+b_{p}}<1$ and
$\diff{t_{\frac{a_{i}b_{p}}{a_{i}+b_{p}}}}=0$, hence we can write
\eqref{eq:product_h_mToHaveIota_k} as\begin{align*}
h\cdot m & =\sum_{p=1}^{n}r\beta_{p}\diff{t_{b_{p}}}+\sum_{p}\sum_{i:a_{i}>k}\alpha_{i}\beta_{p}\diff{t_{\frac{a_{i}b_{p}}{a_{i}+b_{p}}}}=\\
 & =h\cdot\left(r+\sum_{i:a_{i}>k}\alpha_{i}\diff{t_{a_{i}}}\right)=h\cdot\iota_{k}m\end{align*}
 \qedNoNewLine

In the proof of \eqref{eq:cancellationLawUpTo_kOrder} the exponents
$j\in\N^{n}$ will be tied with the ideals $D_{\alpha_{i}}$ through
the following term:
\begin{defn}
\label{def:quotientOfn-tuples}If $j\in\N^{n}$, with $n\in\N_{>0}$,
and $\alpha\in\left(\R_{>0}\cup\{\infty\}\right)^{n}$, then we set
by definition\[
\frac{j}{\alpha+1}:=\sum_{i=1}^{n}\frac{j_{i}}{\alpha_{i}+1}\]

\end{defn}
\noindent Let us note that in the notation $\frac{j}{\alpha+1}$,
the variables $j$ and $\alpha$ are $n$-tuples. In the particular
case $n=1$, we have that $j$ and $\alpha$ are real numbers and
the notation $\frac{j}{\alpha+1}$ has the usual meaning of a fraction.
If $\alpha_{i}=\infty$, then we define $\frac{j_{i}}{\infty+1}:=0$.
Now we can state and prove the main theorem of this section
\begin{thm}
\label{thm:generalCancellationLaw}Let $m\in\ER$, $n\in\N_{>0}$,
$j\in\N^{n}\setminus\{\underline{0}\}$ and $\alpha\in\R_{>0}^{n}$.
Moreover let us consider $k\in\R$ defined by \begin{equation}
\frac{1}{k}+\frac{j}{\alpha+1}=1\label{eq:conditionToFind_k}\end{equation}

\noindent then 
\begin{enumerate}
\item \label{enu:generalThm_hTojAndProductWithIota}$\forall h\in D_{\alpha_{1}}\times\cdots\times D_{\alpha_{n}}\pti h^{j}\cdot m=h^{j}\cdot\iota_{k}m$ 
\item \label{enu:existenceOf_hWithSumOfInverseOfOrders1}$\omega(m)>k\then\exists h\in D_{\alpha_{1}}\times\cdots\times D_{\alpha_{n}}\pti\frac{1}{\omega(m)}+\frac{1}{\omega(h^{j})}=1$ 
\item \textup{\emph{\label{enu:generalCancellationLaw}If $h^{j}\cdot m=0$
for every $h\in D_{\alpha_{1}}\times\cdots\times D_{\alpha_{n}}$,
then $m=_{k}0$}} 
\end{enumerate}
\end{thm}
\noindent The idea of the cancellation law \emph{\ref{enu:generalCancellationLaw}.}
is that if we have $h_{1}^{j_{1}}\cdot\ldots\cdot h_{n}^{j_{n}}\cdot m=0$
for every $(h_{1},\ldots,h_{n})\in D_{\alpha_{1}}\times\cdots\times D_{\alpha_{n}}$,
then condition \eqref{eq:conditionToFind_k} permits to find the best
$k\ge1$ such that $m=_{k}0$. Note that there is no limitation neither
on the exponents $j\in\N^{n}\setminus\{\underline{0}\}$ nor on the
ideals $D_{\alpha_{i}}$, so we can call \emph{\ref{enu:generalCancellationLaw}}.
the \emph{general cancellation law}.

\bigskip{}

\noindent \textbf{Proof of Theorem \ref{thm:generalCancellationLaw}:}

\noindent \emph{\ref{enu:generalThm_hTojAndProductWithIota}}.)\ \ \ By
the definition of $k$ it follows\begin{equation}
\frac{1}{k}+\frac{1}{\omega(h^{j})}=1-\frac{j}{\alpha+1}+\frac{1}{\omega(h^{j})}=1-\sum_{i=1}^{n}\frac{j_{i}}{\alpha_{i}+1}+\sum_{i=1}^{n}\frac{j_{i}}{\omega(h_{i})}\label{eq:firstOneOver_kPlusOneOverOrder_hToj}\end{equation}
 where we have supposed $h^{j}\ne0$, otherwise the conclusion is
trivial, and we have applied Theorem \ref{thm:productOfPowers}.

\noindent But $\omega(h_{i})<\alpha_{i}+1$ because $h_{i}\in D_{\alpha_{i}}$,
so\begin{equation}
\sum_{i=1}^{n}\frac{j_{i}}{\omega(h_{i})}-\sum_{i=1}^{n}\frac{j_{i}}{\alpha_{i}+1}>0\label{eq:addendOfOne}\end{equation}
 Hence from \eqref{eq:firstOneOver_kPlusOneOverOrder_hToj} and \eqref{eq:addendOfOne}
we have $\frac{1}{k}+\frac{1}{\omega(h^{j})}>1$ and the conclusion
follows from Lemma \ref{lem:producOfInfinitesimalAndIota}.

\noindent \emph{\ref{enu:existenceOf_hWithSumOfInverseOfOrders1}}.)\ \ \ For
simplicity let $x_{i}:=(\alpha_{i}+1)\cdot\frac{j}{\alpha+1}$, and
$\frac{1}{a}:=\omega(m)$, then $\frac{1}{a}>k$ from the hypothesis
$\omega(m)>k$ and\[
\frac{1-a}{x_{i}}>\frac{1}{x_{i}}\left(1-\frac{1}{k}\right)=\frac{1}{x_{i}}\cdot\frac{j}{\alpha+1}=\frac{1}{\alpha_{i}+1}\]
 so%
\footnote{Here we are using the usual abuse of notation that consists in indicating
the Fermat real (equivalence class modulo $\sim$, see \ref{def:equalityInFermatReals})
$[t\in\R_{\ge0}\mapsto t^{b}]_{\sim}$ simply by $t^{b}$.%
} $\diff{t_{\frac{x_{i}}{1-a}}}=t^{\frac{1-a}{x_{i}}}\in D_{\alpha_{i}}$,
and if we set\[
h:=\left(\diff{t_{\frac{x_{1}}{1-a}}},\ldots,\diff{t_{\frac{x_{n}}{1-a}}}\right)\in D_{\alpha_{1}}\times\cdots\times D_{\alpha_{n}}\]
 we have\[
h^{j}=t^{\frac{1-a}{x_{1}}j_{1}}\cdot\ldots\cdot t^{\frac{1-a}{x_{n}}j_{n}}=t^{(1-a)\cdot\sum_{i}\frac{j_{i}}{x_{i}}}\]

\noindent But\[
\sum_{i=1}^{n}\frac{j_{i}}{x_{i}}=\sum_{i=1}^{n}\frac{j_{i}}{(\alpha_{i}+1)\cdot\sum_{k=1}^{n}\frac{j_{k}}{\alpha_{k}+1}}=1\]
 hence $h^{j}=t^{1-a}$ and\[
\frac{1}{\omega(m)}+\frac{1}{\omega(h^{j})}=a+(1-a)=1\]

\noindent \emph{\ref{enu:generalCancellationLaw}}.)\ \ \ This
part is essentially the contrapositive of \emph{\ref{enu:existenceOf_hWithSumOfInverseOfOrders1}}.
Indeed from \emph{\ref{enu:existenceOf_hWithSumOfInverseOfOrders1}}.
we have\begin{equation}
\left(\forall h\in D_{\alpha_{1}}\times\cdots\times D_{\alpha_{n}}\pti\frac{1}{\omega(m)}+\frac{1}{\omega(h^{j})}\ne1\right)\then\omega(m)\le k\label{eq:contrapositiveOf2}\end{equation}
 so if we assume that $h^{j}\cdot m=0$ for every $h\in D_{\alpha_{1}}\times\cdots\times D_{\alpha_{n}}$,
then it immediately follows $\st{m}=0$ and from Theorem \ref{thm:productOfPowers}
the equality $h^{j}\cdot m=0$ becomes equivalent to\[
\frac{1}{\omega(m)}+\frac{1}{\omega(h^{j})}>1\]

\noindent Therefore \eqref{eq:contrapositiveOf2} is actually stronger
than the hypothesis of \emph{\ref{enu:generalCancellationLaw}}.

\noindent From \eqref{eq:contrapositiveOf2} it follows $\omega(m)\le k$
and hence $m=_{k}0$.\qedNoNewLine

For example suppose we want to obtain $m=_{2}0$ from a product of
the type $h\cdot m=0$ for every $h\in D_{\alpha}$. What kind of
infinitesimals $D_{\alpha}$ do we have to choose? We have $k=2$,
$j=1$ and $n=1$, hence we must have\[
\frac{1}{k}+\frac{j}{\alpha+1}=\frac{1}{2}+\frac{1}{\alpha+1}=1\]
 hence $\alpha=1$ and from the general cancellation law we have\[
\left(\forall h\in D\pti h\cdot m=0\right)\then m=_{2}0\]

\noindent Analogously if we want $n=2$, then we must have\[
\frac{1}{2}=1-\frac{1}{\alpha_{1}+1}-\frac{1}{\alpha_{2}+1}\]
 so that we must choose $\alpha=(\alpha_{1},\alpha_{2})$ so that
$\frac{1}{\alpha_{1}+1}+\frac{1}{\alpha_{2}+1}=\frac{1}{2}$, e.g.
$\alpha=(3,3)$, i.e.\[
\left(\forall h,k\in D_{3}\pti h\cdot k\cdot m=0\right)\then m=_{2}0\]

\noindent Vice versa now suppose to have $D_{\alpha}$, with $\alpha\in\N_{>0}$
and we want to find $k$:\[
\frac{1}{k}=1-\frac{1}{\alpha+1}\]
 hence $k=\frac{\alpha+1}{\alpha}$ and we obtain\[
\begin{array}{ccc}
\left(\forall h\in D_{\alpha}\pti h\cdot m=0\right) & \then & m=_{\frac{\alpha+1}{\alpha}}0\end{array}\]

\noindent and\[
\left(\forall h\in D_{\infty}\pti h\cdot m=0\right)\then m=_{1}0\]
Let us note explicitly that the best we can obtain from the general
cancellation law is that $m$ is equal to zero up to first order infinitesimals.
As an immediate consequence of the definition of equality in $\ER$,
it follows that $h\cdot\diff{t}=0$ for every infinitesimal $h\in D_{\infty}$,
and because $\diff{t}=_{1}0$ but $\diff{t}\ne0$, this exhibits that
a better result cannot be obtained from this type of cancellation
law.

\subsection*{A counterexample}

The idea to have a cancellation law like the general one \emph{\ref{enu:generalCancellationLaw}}.
of Theorem \ref{thm:generalCancellationLaw} comes from SDG. The particularity
of this law is that it is not of the form

\begin{center}
{}``if \emph{a given number} $h\in\ER$ has the property $\mathcal{P}(h)$
(e.g. $h$ is invertible), and $h\cdot m=0$, then $m=0$'', 
\par\end{center}

\noindent as usual, but it is of the form, e.g.

\noindent \begin{center}
{}``if $h\cdot m=0$ \emph{for every} $h\in D$, then $m=_{2}0$''. 
\par\end{center}

\noindent We can foresee that these differences will not cause any
problem each time we will use infinitesimal Taylor's formulae. Indeed,
as we will see concretely later in the present work, typically these
formulae are used for \emph{generic} infinitesimal increments $h\in D_{\alpha}^{n}$,
i.e. usually we will be able to prove our equalities derived from
Taylor's formulae \emph{for every} $h\in D_{\alpha}^{n}$.

Finally, these cancellation laws do not guarantee a strict equality
but an equality up to infinitesimals of a suitable order $k$. As
we will see, this correspond to have Taylor's formulae with uniqueness
up to some order $k$. If we use these formulae to define derivatives,
this implies that we will have derivatives identified up to infinitesimals
of some order $k$. Roughly speaking, even if this is unusual for
derivatives of smooth functions, it is very common in mathematics;
think e.g. to definite integrals or Radon-Nikodym derivatives, where
certain operators are defined up to a suitable notion of equality
(i.e. an equivalence relation) like {}``up to a constant'' or {}``up
to a set of measure zero''. In the same way, e.g., we will define
first derivatives of smooth functions of the type $f:D_{\alpha}^{n}\freccia\ER$
up to second order infinitesimals. Exactly as for the Radon-Nikodym
derivative, the meaningful properties will be only those {}``up to
second order infinitesimals''.

Now we want to see that it is not possible to avoid the quantifier
{}``\emph{for every} $h$'' in the cancellation law. More precisely
let us suppose to have an infinitesimal $h\in D_{\infty}$ with the
property of being deleted from every product, i.e. such that\begin{equation}
\forall m\in\ER\pti h\cdot m=0\then m=0\label{eq:hypothesisForCounterexample}\end{equation}
 Does such an infinitesimal exist?

Of course $h\ne0$ and hence $\omega(h)\ge1$, but it cannot be that
$\omega(h)=1$ because otherwise\[
\frac{1}{\omega(h)}+\frac{1}{\omega(\diff{t})}=2\]
 and hence $h\cdot\diff{t}=0$ even if $\diff{t}\ne0$, in contradiction
with \eqref{eq:hypothesisForCounterexample}. Hence it must be that
$\omega(h)>1$.

Now we want to find $k\ge1$ such that $\frac{1}{k}+\frac{1}{\omega(h)}>1$,
that is\[
\frac{1}{k}>1-\frac{1}{\omega(h)}=\frac{\omega(h)-1}{\omega(h)}>0\]
 the latter inequality being due to $\omega(h)>1$. Therefore the
number $k$ we are searching for must satisfy\begin{equation}
1\le k<\frac{\omega(h)}{\omega(h)-1}\label{eq:conditionFor_k}\end{equation}
 A $k$ in this interval exists always because $\omega(h)>1$, and
if we set $m:=\diff{t_{k}}$, then $\frac{1}{k}+\frac{1}{\omega(h)}>1$,
and from Lemma \ref{lem:producOfInfinitesimalAndIota} we have\[
h\cdot m=h\cdot\iota_{k}m=h\cdot\iota_{k}(\diff{t_{k}})=0\]
 but $m=\diff{t_{k}}\ne0$. For these reasons we can affirm that an
infinitesimal $h\in D_{\infty}$ with the property of being deleted
form every product, i.e. such that \eqref{eq:hypothesisForCounterexample}
holds, does not exist.

\section{Applications to Taylor's formulae}

\subsection*{General forms of uniqueness in Taylor's formulae}
\begin{cor}
\label{cor:uniquenessInGeneralTaylorFormulae}Let $n\in\N_{>0}$,
$\alpha\in\R_{>0}^{n}$, and for every $j\in\N^{n}$ with $0<\frac{j}{\alpha+1}<1$,
let $m_{j}\in\ER$ and set $k_{j}\in\R$ such that\[
\frac{1}{k_{j}}+\frac{j}{\alpha+1}=1\e{,}k_{\underline{0}}:=0\]

\noindent Then there exists one and only one\[
\bar{m}:\left\{ j\in\N^{n}\,|\,\frac{j}{\alpha+1}<1\right\} \freccia\ER\]

\noindent such that 
\begin{enumerate}
\item \label{enu:coefficientsAreInER_k_j}$\bar{m}_{j}\in\ER_{k_{j}}$ for
every $j\in\N^{n}$ such that $\frac{j}{\alpha+1}<1$
\item \label{enu:equalityBetweenTaylorFormulae}$\forall h\in D_{\alpha_{1}}\times\cdots\times D_{\alpha_{n}}\pti{\displaystyle \sum_{\substack{j\in\N^{n}\\
\frac{j}{\alpha+1}<1}
}\frac{h^{j}}{j!}\cdot m_{j}=\sum_{\substack{j\in\N^{n}\\
\frac{j}{\alpha+1}<1}
}\frac{h^{j}}{j!}\cdot\bar{m}_{j}}$ 
\end{enumerate}
Moreover the unique $\bar{m}_{j}$ is given by $\bar{m}_{j}=\iota_{k_{j}}m_{j}$.
\end{cor}
\noindent To motivate the statement let us observe that if $h\in D_{\alpha_{1}}\times\cdots\times D_{\alpha_{n}}$
and $\frac{j}{\alpha+1}\ge1$, then $\omega(h_{i})<\alpha_{i}+1$,
and so\[
\frac{j_{i}}{\omega(h_{i})}>\frac{j_{i}}{\alpha_{i}+1}\]
 and hence also\[
\sum_{i=1}^{n}\frac{j_{i}}{\omega(h_{i})}>\sum_{i=1}^{n}\frac{j_{i}}{\alpha_{i}+1}=\frac{j}{\alpha+1}\ge1\]
thus  $h^{j}=h_{1}^{j_{1}}\cdot\ldots\cdot h_{n}^{j_{n}}=0$. For
this reason the general Taylor's formula is restricted to $j\in\N^{n}$
such that $\frac{j}{\alpha+1}<1$.

The meaning of this corollary is that if we have an infinitesimal
Taylor's formula like\[
f(h)=f(0)+\sum_{j=1}^{n}\frac{h^{j}}{j!}\cdot m_{j}\quad\forall h\in D_{n}\]
 then we can substitute the coefficients $m_{j}\in\ER$ by $\bar{m}_{j}=\iota_{k_{j}}(m_{j})\in\ER_{k_{j}}$,
that is with $m_{j}$ up to infinitesimals of order $k_{j}$, and
the formula remains unchanged\begin{equation}
f(h)=f(0)+\sum_{j=1}^{n}\frac{h^{j}}{j!}\cdot\bar{m}_{j}\quad\forall h\in D_{n}\label{eq:TaylorWithUniqueCoefficientsComments}\end{equation}
 But now the new coefficients $\bar{m}_{j}\in\ER_{k_{j}}$ are uniquely
determined by \eqref{eq:TaylorWithUniqueCoefficientsComments}.

E.g. this will permit to prove that if $f:D\freccia\ER$, then there
exist one and only one pair\begin{align*}
a\in\ER\\
b\in\ER_{2}\end{align*}
 such that $f(h)=a+h\cdot b$ for every $h\in D$.

\bigskip{}

\noindent \textbf{Proof of Corollary \ref{cor:uniquenessInGeneralTaylorFormulae}:}

\noindent \emph{Existence:} Let $\bar{m}_{j}:=\iota_{k_{j}}(m_{j})$
for every $j\in\N^{n}$ such that $\frac{j}{\alpha+1}<1$. Note that
if $j=\underline{0}$, then $\bar{m}_{j}=m_{j}=m_{\underline{0}}$,
because $k_{\underline{0}}:=0$. Moreover if $j\ne\underline{0}$,
then $0<\frac{j}{\alpha+1}<1$ and hence $k_{j}>1$.

\noindent We have $\bar{m}_{j}\in\ER_{k_{j}}$ and, from Theorem \ref{thm:generalCancellationLaw}
for every $j$ we obtain\[
\forall h\in D_{\alpha_{1}}\times\cdots\times D_{\alpha_{n}}\pti h^{j}\cdot m_{j}=h^{j}\cdot\bar{m}_{j}\]
 and hence also the conclusion\[
\forall h\in D_{\alpha_{1}}\times\cdots\times D_{\alpha_{n}}\pti{\displaystyle \sum_{\substack{j\in\N^{n}\\
\frac{j}{\alpha+1}<1}
}\frac{h^{j}}{j!}\cdot m_{j}=\sum_{\substack{j\in\N^{n}\\
\frac{j}{\alpha+1}<1}
}\frac{h^{j}}{j!}\cdot\bar{m}_{j}}\]
 \emph{Uniqueness:} Let us consider $\hat{m}_{j}\in\ER_{k_{j}}$ that
verify the identity \emph{\ref{enu:equalityBetweenTaylorFormulae}}.,
we shall use the identity principle for polynomials (Theorem \ref{thm:PIP}).
Indeed for each fixed $h\in D_{\alpha_{1}}\times\cdots\times D_{\alpha_{n}}$
and every $r\in\ext{(-1,1)}$ we have $r\cdot h\in D_{\alpha_{1}}\times\cdots\times D_{\alpha_{n}}$,
hence\[
\sum_{\substack{j\in\N^{n}\\
\frac{j}{\alpha+1}<1}
}r^{j}\cdot\frac{h^{j}}{j!}\cdot\left(\bar{m}_{j}-\hat{m}_{j}\right)=0\quad\forall r\in\ext{(-1,1)}\]
 From the identity principle of polynomials every coefficient of this
polynomial in $r$ is zero, i.e.\[
\forall j\pti\frac{h^{j}}{j!}\cdot\left(\bar{m}_{j}-\hat{m}_{j}\right)=0\]
 These equalities are also true for every $h\in D_{\alpha_{1}}\times\cdots\times D_{\alpha_{n}}$,
therefore from Theorem \ref{thm:generalCancellationLaw} we obtain
$\bar{m}_{j}=_{k_{j}}\hat{m}_{j}$, that is $\bar{m}_{j}=\hat{m}_{j}$
because $\bar{m}_{j}$, $\hat{m}_{j}\in\ER_{k_{j}}$ (see Remark \ref{rem:equalityInER_k}).\qedNoNewLine

Using the equalities up to $k_{j}$-th order infinitesimals we can
state this uniqueness in another equivalent form:
\begin{cor}
In the hypotheses of the previous Corollary \ref{cor:uniquenessInGeneralTaylorFormulae},
if $p_{j}\in\ER$ for every $j\in\N^{n}$ with $\frac{j}{\alpha+1}<1$
are such that\[
\forall h\in D_{\alpha_{1}}\times\cdots\times D_{\alpha_{n}}\pti{\displaystyle \sum_{\substack{j\in\N^{n}\\
\frac{j}{\alpha+1}<1}
}\frac{h^{j}}{j!}\cdot m_{j}=\sum_{\substack{j\in\N^{n}\\
\frac{j}{\alpha+1}<1}
}\frac{h^{j}}{j!}\cdot p_{j}}\]

\noindent then $m_{j}=_{k_{j}}p_{j}$ for every $j$. 
\end{cor}
\noindent \textbf{Proof:} In fact we can apply the previous Corollary
\ref{cor:uniquenessInGeneralTaylorFormulae} both with $(m_{j})_{j}$
and with $(p_{j})_{j}$ obtaining that the unique $(\bar{m}_{j})_{j}$
is given by $\bar{m}_{j}=\iota_{k_{j}}(m_{j})=\iota_{k_{j}}(p_{j})$,
so $m_{j}=_{k_{j}}p_{j}$ for every $j$.\qedNoNewLine

\subsection*{Existence in Taylor's formulae for ordinary smooth functions}

The following theorem is a \emph{very simple} evidence that a suitable
and meaningful mathematical language can be useful to extend even
well known classical results. Indeed, using the language of actual
nilpotent infinitesimals we shall see that it is possible to extend
the Taylor's formula for $f(x+h)$ to generic infinitesimal increments
$h\in D_{\alpha_{1}}\times\cdots\times D_{\alpha_{d}}$ (the classical
formulation being for $\alpha_{1}=\dots=\alpha_{d}$):
\begin{thm}
\label{thm:GenericTaylorFormulaForOrdinarySmoothFunctions}Let $f:U\freccia\R^{u}$
be a smooth function, with $U$ open in $\R^{d}$. Take a standard
point $x\in U$ and $\alpha_{1},\ldots,\alpha_{d}\in\R_{>0}$, then
there exist one and only one\[
m:\left\{ j\in\N^{d}\,|\,\frac{j}{\alpha+1}<1\right\} \freccia\R^{u}\]

\noindent such that\[
\forall h\in D_{\alpha_{1}}\times\cdots\times D_{\alpha_{d}}\pti{\displaystyle f(x+h)=\sum_{\substack{j\in\N^{d}\\
\frac{j}{\alpha+1}<1}
}\frac{h^{j}}{j!}\cdot m_{j}}\]

\end{thm}
\noindent \textbf{Proof:} For simplicity let \[
I:=\left\{ j\in\N^{d}\,|\,\frac{j}{\alpha+1}<1\right\} \e{,}n:=\max\left\{ |j|\in\N\,|\, j\in I\right\} \]
 where, of course, $|j|:=j_{1}+\ldots+j_{n}$. Let us take the infinitesimal
Taylor's formula of $f$ of order $n$ (see Theorem \ref{thm:OrdinaryTaylor'sFor_nVariables}):\begin{equation}
\forall h\in D_{n}^{d}\pti f(x+h)=\sum_{\substack{j\in\N^{d}\\
|j|\le n}
}\frac{h^{j}}{j!}\cdot m_{j}\label{eq:TaylorOrder_nForGeneralTaylor}\end{equation}
 where $m_{j}:=\frac{\partial^{|j|}f}{\partial x^{j}}(x)\in\R^{u}$.
Now if we take $h\in D_{\alpha_{1}}\times\cdots\times D_{\alpha_{d}}$,
then $h_{i}\in D_{\alpha_{i}}$ and hence $\omega(h_{i})<\alpha_{i}+1$.
We want to apply \eqref{eq:TaylorOrder_nForGeneralTaylor} with this
$h$, so we have to prove that $h_{i}\in D_{n}$, i.e. that $\omega(h_{i})<n+1$.
But if we set $j:=(0,\ptind^{i-1},0,\alpha_{i},0,\ldots,0)$, then\[
\frac{j}{\alpha+1}=\sum_{k=1}^{d}\frac{j_{k}}{\alpha_{k}+1}=\frac{\alpha_{i}}{\alpha_{i}+1}<1\]
 because $\alpha_{i}>0$, so $j\in I$ and hence $n\ge|j|=\alpha_{i}$.
Therefore $\omega(h_{i})<\alpha_{i}+1\le n+1$, that is $h_{i}\in D_{n}$
and we can apply \eqref{eq:TaylorOrder_nForGeneralTaylor} obtaining:\begin{equation}
f(x+h)=\sum_{\substack{j\in\N^{d}\\
|j|\le n}
}\frac{h^{j}}{j!}\cdot m_{j}\label{eq:TaylorWithTooMuchAddends}\end{equation}
 But we know that if $\frac{j}{\alpha+1}\ge1$, then $h^{j}=0$, so
the sum in \eqref{eq:TaylorWithTooMuchAddends} is extended to $j\in I$
only. This proves the existence part. Uniqueness follows from Corollary
\ref{cor:uniquenessInGeneralTaylorFormulae}.\qedNoNewLine

At present the previous version of the Taylor's formula can be applied
to ordinary smooth functions and to standard points $x\in U$ only.
In the following results we will remove the limitation that the base
point $x$ has to be standard.
\begin{lem}
\label{lem:TaylorForGenericBasePoint}Let $A$ be an open set in $\R^{d}$,
$x\in\ext{A}$, $n\in\N_{>0}$ and $f:A\freccia\R$ a smooth function\emph{,}
then\begin{equation}
\forall h\in D_{n}^{d}\pti f(x+h)=\sum_{\substack{j\in\N^{d}\\
|j|\le n}
}\frac{h^{j}}{j!}\cdot\frac{\partial^{|j|}f}{\partial x^{j}}(x)\label{eq:TaylorWithGenericPoint}\end{equation}

\end{lem}
Note that in \eqref{eq:TaylorWithGenericPoint} we do not have the
problem to define the derivatives of the function $f$ at the non
standard point $x\in\ext{A}$, because we have to mean $\frac{\partial^{|j|}f}{\partial x^{j}}(x)$
as\[
\frac{\partial^{|j|}f}{\partial x^{j}}(x)=^{^{^{^{^{{\scriptstyle \,\,\bullet\!\!\!\!\!}}}}}}\left(\frac{\partial^{|j|}f}{\partial x^{j}}\right)(x),\]
that is as the Fermat extension of the smooth function $\frac{\partial^{|j|}f}{\partial x^{j}}(x)$
applied to the non standard point $x$.

\noindent \textbf{Proof:} We prove the result for $d=1$ only; the
proof for the multivariable case is analogous using the suitable multi-indexes
notations. Let $k:=x-\st{x}$ be the nilpotent part of $x\in\ext{A}$,
then $f(x+h)=f(\st{x}+k+h)$, and we can use the infinitesimal Taylor's
formula (Theorem \ref{thm:OrdinaryTaylorFor_nVariables}) for $f$
at the standard point $\st{x}$ and with infinitesimal increment $k+h$.
Let us firstly suppose $k+h\ne0.$ Then the order of this sum is $\omega(k+h)=\omega(k)\vee\omega(h)$
(see Theorem \ref{thm:propertyOfIdeal_D_a}) and we can write \begin{align*}
f(x+h) & =f(\st{x}+k+h)=\sum_{b=0}^{\omega(h)\vee\omega(k)}\frac{(k+h)^{b}}{b!}\cdot f^{(b)}(\st{x})=\\
 & =\sum_{b=0}^{\omega(h)\vee\omega(k)}\sum_{a=0}^{b}\binom{b}{a}\frac{k^{a}h^{b-a}}{b!}\cdot f^{(b)}(\st{x})=\\
 & =\sum_{b=0}^{\omega(h)\vee\omega(k)}\sum_{a=0}^{b}\frac{k^{a}h^{b-a}}{a!(b-a)!}\cdot f^{(b)}(\st{x}).\end{align*}
But $k^{a}=0$ if $a>\omega(k)$, hence\begin{align*}
f(x+h) & =\sum_{b=0}^{\omega(h)\vee\omega(k)}\sum_{a=0}^{b\wedge\omega(k)}\frac{k^{a}h^{b-a}}{a!(b-a)!}\cdot f^{(b)}(\st{x})=\\
 & =\sum_{s\in I}s.\end{align*}
where\begin{equation}
I:=\left\{ \frac{k^{a}h^{b-a}}{a!(b-a)!}\cdot f^{(b)}(\st{x})\,\left|\,\substack{b=0,\ldots,\omega(h)\vee\omega(k)\\
a=0,\ldots,b\wedge\omega(k)}
\right.\right\} .\label{eq:leftSet}\end{equation}
On the other hand we have \begin{align*}
\sum_{j=0}^{\omega(h)}\frac{h^{j}}{j!}\cdot f^{(j)}(x) & =\sum_{j=0}^{\omega(h)}\frac{h^{j}}{j!}\cdot f^{(j)}(\st{x}+k)=\\
 & =\sum_{j=0}^{\omega(h)}\sum_{i=0}^{\omega(k)}\frac{h^{j}}{j!}\frac{k^{i}}{i!}\cdot f^{(j+i)}(\st{x})=\\
 & =\sum_{t\in J}t\end{align*}
where\begin{equation}
J:=\left\{ \frac{h^{j}}{j!}\frac{k^{i}}{i!}\cdot f^{(j+i)}(\st{x})\,\left|\,\substack{j=0,\ldots,\omega(h)\\
i=0,\ldots,\omega(k)\\
\frac{i}{\omega(k)}+\frac{j}{\omega(h)}\le1}
\right.\right\} \label{eq:rightSet}\end{equation}
in fact $h^{j}k^{i}=0$ if $\frac{i}{\omega(k)}+\frac{j}{\omega(h)}>1$
(see Theorem \ref{thm:productOfPowers}). Now we can prove that the
two sets of addends $I$ and $J$ are equal or they differ at most
for zero, i.e. $I\cup\left\{ 0\right\} =J\cup\left\{ 0\right\} $.
Indeed, take an element $t$ from $J$ \[
t=\frac{h^{j}k^{i}}{j!i!}\cdot f^{(j+i)}(\st{x})\]
\[
j=0,\ldots,\omega(h)\]
\[
i=0,\ldots,\omega(k)\]
\begin{equation}
\frac{i}{\omega(k)}+\frac{j}{\omega(h)}\le1,\label{eq:conditionNoZerokihj}\end{equation}
then setting $a:=i$, $b:=a+j=i+j$ we have\[
t=\frac{h^{b-a}k^{a}}{(b-a)!a!}\cdot f^{(b)}(\st{x}).\]
Moreover, we have that $a=i\le\omega(k)$ and $a=b-j\le b$, so $a=0,\ldots,b\wedge\omega(k)$.
From \ref{eq:conditionNoZerokihj} we have $i\le\omega(k)-j\frac{\omega(k)}{\omega(h)}$
and hence $i+j\le\omega(k)+j\cdot\left(1-\frac{\omega(k)}{\omega(h)}\right)=\omega(k)+j\cdot\frac{\omega(h)-\omega(k)}{\omega(h)}$.
If $\omega(k)\ge\omega(h)$, then $i+j\le\omega(k)+j\cdot\frac{\omega(h)-\omega(k)}{\omega(h)}\le\omega(k)=\omega(h)\vee\omega(k)$.
Vice versa if $\omega(k)<\omega(h)$, then since $\frac{j}{\omega(h)}\le1$
we have that $\omega(k)+j\cdot\frac{\omega(h)-\omega(k)}{\omega(h)}\le\omega(k)+\omega(h)-\omega(k)=\omega(h)=\omega(h)\vee\omega(k)$.
In any case we have proved that $b=i+j\le\omega(h)\vee\omega(k)$,
so the addend $t$ is indeed an element of $J$.

Vice versa, let us consider\[
s=\frac{h^{b-a}k^{a}}{(b-a)!a!}\cdot f^{(b)}(\st{x})\in J\]
\[
b=0,\ldots,\omega(h)\vee\omega(k)\]
\[
a=0,\ldots,b\wedge\omega(k),\]
and set $i:=a$, $j=b-a$, then\[
s=\frac{h^{j}k^{i}}{j!i!}\cdot f^{(j+i)}(\st{x}).\]
Moreover, $i\le b\wedge\omega(k)\le\omega(k)$ and $j\le\omega(h)$
or, in the opposite case, we have $h^{b-a}=h^{j}=0=s$. Analogously
we have $\frac{i}{\omega(k)}+\frac{j}{\omega(h)}\le1$ or, in the
opposite case, we have $h^{b-a}k^{a}=h^{j}k^{i}=0=s$. At the end
we have proved that $s\in I\cup\left\{ 0\right\} $.

It remains to prove the case $h+k=0$. But with the previous deduction
we have proved that\begin{equation}
\sum_{b=0}^{\omega(h)\vee\omega(k)}\frac{(k+h)^{b}}{b!}\cdot f^{(b)}(\st{x})=\sum_{j=0}^{\omega(h)}\frac{h^{j}}{j!}\cdot f^{(j)}(\st{x}+k).\label{eq:equalityOfTwoTaylor}\end{equation}
If $h+k=0$ the right hand side of \eqref{eq:equalityOfTwoTaylor}
gives $f(\st{x})=f(\st{x}+k+h)=f(x+h)$.\qedWithFinalEq

Using this lemma, and the general uniqueness of Corollary \ref{cor:uniquenessInGeneralTaylorFormulae},
we can repeat equal the proof of Theorem \ref{thm:GenericTaylorFormulaForOrdinarySmoothFunctions}
obtaining its generalization to a non standard base point $x\in\ext{U}$:
\begin{thm}
\label{thm:TaylorForStndSmoothAtNonStndPoint}Let $f:U\freccia\R^{u}$
be a smooth function, with $U$ open in $\R^{d}$. Take a point $x\in\ext{U}$
and $\alpha_{1},\ldots,\alpha_{d}\in\R_{>0}$, then there exists one
and only one\[
m:\left\{ j\in\N^{d}\,|\,\frac{j}{\alpha+1}<1\right\} \freccia\ER\]
 such that
\begin{enumerate}
\item $\bar{m}_{j}\in\ER_{k_{j}}$ for every $j\in\N^{d}$ such that $\frac{j}{\alpha+1}<1$
\item $\forall h\in D_{\alpha_{1}}\times\cdots\times D_{\alpha_{d}}\pti{\displaystyle f(x+h)=\sum_{\substack{j\in\N^{d}\\
\frac{j}{\alpha+1}<1}
}\frac{h^{j}}{j!}\cdot m_{j}}$ 
\end{enumerate}
\end{thm}
In the following chapters we shall see how to generalize these theorems
to more general functions with respect to $\ext{f}$, i.e. extension
of standard functions. We have an example of a function which is not
of this type, considering e.g. $f=\ext{g}(p,-)$ for $p\in\ER$ and
$g\in\Cc^{\infty}(\R,\R)$. In this case we will see that in general
the coefficients of the corresponding Taylor's formulae will be generic
elements $m_{j}\in\ER_{k_{j}}$ and not standard reals only.
\begin{example*}
We want to find the Taylor's formula of\begin{equation}
f(h,k)=\frac{\sin(h)}{\cos(k)}\label{eq:exampleForGenericTaylorFormula}\end{equation}
 for $h\in D_{3}$ and $k\in D_{4}$. We can note, using the previous
theorem, that the sum in this Taylor's formula is extended to all
pair $(i,j)\in\N^{2}$ such that\[
\frac{i}{3+1}+\frac{j}{4+1}<1\]
 that is such that $5i+4j<20$. But to find this Taylor's formula
it is simpler to substitute in \eqref{eq:exampleForGenericTaylorFormula}
the Taylor's formulae of $\sin(h)=h-\frac{h^{3}}{6}$ for $h\in D_{3}$
and of $\cos(k)=1-\frac{k^{2}}{2}+\frac{k^{4}}{24}$ for $k\in D_{4}$
and to apply the \emph{algebraic} calculus of nilpotent infinitesimals
we have developed until now:\begin{align}
f(h,k) & =\frac{h-\frac{h^{3}}{6}}{1-\frac{k^{2}}{2}+\frac{k^{4}}{24}}=\left(h-\frac{h^{3}}{6}\right)\cdot\left(1+\frac{k^{2}}{2}-\frac{k^{4}}{24}+\left[\frac{k^{4}}{24}-\frac{k^{2}}{2}\right]^{2}\right)=\nonumber \\
 & =h+\frac{1}{2}hk^{2}-\frac{h^{3}}{6}\quad\forall h\in D_{3}\ \forall k\in D_{4}\label{eq:exampleGeneralTaylor}\end{align}
 For example to obtain this result we have used the equalities $hk^{4}=0$
and $h^{3}k^{2}=0$, easily deducible from Corollary \ref{cor:suffCondToHaveProductZeroFrom_D_j_k}.
\end{example*}
\noindent From \eqref{eq:exampleGeneralTaylor} and from the previous
Theorem \ref{thm:GenericTaylorFormulaForOrdinarySmoothFunctions}
we have\begin{align*}
\frac{\partial f}{\partial h}(0,0) & =1\quad\frac{\partial^{3}f}{\partial h\partial k^{2}}(0,0)=1\quad\frac{\partial^{3}f}{\partial h^{3}}(0,0)=-1\end{align*}
 and for all the other indexes $i$, $j\in\N$ such that $5i+4j<20$
we have

\noindent \[
\frac{\partial^{i+j}f}{\partial h^{i}\partial k^{j}}(0,0)=0\]

Of course this is \emph{only an elementary example}, similar to several
exercises one can find in elementary courses of Calculus. The only
meaningful difference is that we have not used directly neither the
concept of limit nor any rest in the form of suitable little-oh functions.
An easy to use algebraic language of nilpotent infinitesimals have
been used instead. It can also be useful to note that, in comparison
with SDG, for Fermat reals it is very easy to decide if products of
type $h^{3}k^{2}$, with $h\in D_{3}$ and $k\in D_{4}$, are zero
or not; the same easiness is not possible in SDG where starting only
from the belonging to some $D_{k}$ it is not possible to decide products
of this type (see e.g. \cite{Koc} for more details).

\section{Extension of some results to $D_{\infty}$}

In this section we want to extend some of the results of the previous
sections to the ideal $D_{\infty}$ of all the infinitesimals (see
Definition \ref{def:idealD_a}).
\begin{cor}
\label{cor:generalCancellationLawForDInfty}Let $m\in\ER$, $n\in\N_{>0}$,
$j\in\N^{n}\setminus\{\underline{0}\}$

\noindent then 
\begin{enumerate}
\item \label{enu:generalThm_hTojAndProductWithIotaForDInfty}$\forall h\in D_{\infty}^{n}\pti h^{j}\cdot m=h^{j}\cdot\iota_{1}m$ 
\item \textup{\emph{\label{enu:generalCancellationLawForDInfty}If $h^{j}\cdot m=0$
for every $h\in D_{\infty}^{n}$, then $m=_{1}0$, that is $m=\alpha\diff{t}$
for some $\alpha\in\R$.}} 
\end{enumerate}
\end{cor}
\noindent \textbf{Proof:} We recall that\[
\iota_{1}m=\st{m}+\sum_{i:\omega_{i}(m)>1}^{N}\st{m_{i}}\diff{t_{\omega_{i}(m)}}\]
 and because $\omega_{i}(m)\ge1$ for every Fermat real $m$ and every
$i=1,\ldots,N$ we can write\[
m=\alpha\diff{t}+\iota_{1}m\]
 where $\alpha:=\st{m_{\bar{\imath}}}$ if $\omega_{\bar{\imath}}(m)=1$
for some $\bar{\imath}=1,\ldots,N$, otherwise $\alpha:=0$. Therefore
if $h\in D_{\infty}^{n}$, we have $h^{j}\cdot m=h^{j}\cdot\alpha\diff{t}+h^{j}\cdot\iota_{1}m=h^{j}\cdot\iota_{1}m$
because $k\diff{t}=0$ for every infinitesimal $k\in D_{\infty}$.
This proves \emph{\ref{enu:generalThm_hTojAndProductWithIotaForDInfty}}.

\noindent \emph{\ref{enu:generalCancellationLawForDInfty}}.)\ \ \ From
the hypothesis $h^{j}\cdot m=0$ for every $h\in D_{\infty}^{n}$,
because $D_{a}\subset D_{\infty}$, it follows\begin{equation}
\forall h\in D_{a}^{n}\pti h^{j}\cdot m=0\label{eq:hTojTimes_mIsZeroForEveryD_a}\end{equation}
 where $a\ge1$. So we can apply Theorem \ref{thm:generalCancellationLaw}
for each one of these $a$. We set $k_{a}$ with\[
\frac{1}{k_{a}}+\frac{j}{(a,\ptind^{n},a)+1}=1\]
 that is\[
\frac{1}{k_{a}}=1-\sum_{i=1}^{n}\frac{j_{i}}{a+1}=1-\frac{|j|}{a+1}\]

\[
k_{a}=\frac{a+1}{a+1-|j|}\]
 and hence from \eqref{eq:hTojTimes_mIsZeroForEveryD_a} and Theorem
\ref{thm:generalCancellationLaw} we have that $m=_{k_{a}}0$ for
every $a\ge1$, that is\begin{equation}
\st{m}+\sum_{\omega_{i}(m)>k_{a}}\st{m_{i}}\diff{t_{\omega_{i}(m)}}=0\quad\forall a\ge1\label{eq:iota_k_aIsZero}\end{equation}
 This implies $\st{m}=0$ and for each $a\ge1$ the sum in \eqref{eq:iota_k_aIsZero}
does not have addends, i.e.\begin{equation}
\forall a\ge1\,\nexists i=1,\ldots,N:\,\omega_{i}(m)>k_{a}\label{eq:thereIsNotOmega_iGreaterThan_k_a}\end{equation}
 But $\lim_{a\to+\infty}k_{a}=1^{+}$, so if we had\[
\exists\bar{\imath}=1,\ldots,N\pti\omega_{\bar{\imath}}(m)>1\]
 then we could find a suitable $\bar{a}\ge1$ such that $\omega_{\bar{\imath}}(m)>k_{\bar{a}}\ge1$,
in contradiction with \eqref{eq:thereIsNotOmega_iGreaterThan_k_a}.
Therefore $\omega_{i}(m)\le1$ for each $i=1,\ldots,N$, i.e. $m=_{1}0$
and hence $m=\alpha\diff{t}$ for $\alpha=\st{m_{1}}$ or $\alpha=0$
if $m=0$.\qedNoNewLine

Using exactly the same ideas used in the proof of the previous corollary
we can also generalize Taylor's formulae to the case of $D_{\infty}$.
First the uniqueness:
\begin{cor}
Let $(m_{j})_{j\in\N^{d}\setminus\{\underline{0}\}}$ be a sequence
of $\ER$, then there exists one and only one\[
(m_{j})_{j\in\N^{d}\setminus\{\underline{0}\}}\text{ sequence of }\ER_{1}\]

\noindent such that\[
\forall h\in D_{\infty}^{d}\pti\sum_{\substack{j\in\N^{d}\\
j\ne\underline{0}}
}\frac{h^{j}}{j!}\cdot m_{j}=\sum_{\substack{j\in\N^{d}\\
j\ne\underline{0}}
}\frac{h^{j}}{j!}\cdot\bar{m}_{j}.\]

\end{cor}
\qedWithFinalEq

\noindent Secondly, it is also easy to derive the Taylor's formula
for standard functions:
\begin{cor}
Let $f:U\freccia\R^{u}$ be a smooth function, with $U$ open in $\R^{d}$.
Take $x\in\ext{U}$, then there exist one and only one\[
m:\N^{d}\freccia\R^{u}\]

\noindent such that\begin{equation}
\forall h\in D_{\infty}^{d}\pti f(x+h)=\sum_{j\in\N^{d}}\frac{h^{j}}{j!}\cdot m_{j}\label{eq:TaylorInDIinftyForStandardFunctions}\end{equation}

\end{cor}
\qedWithFinalEq

\noindent Let us point out that in formulae like \eqref{eq:TaylorInDIinftyForStandardFunctions}
we do not have a series but a finite sum because every $h_{i}\in D_{\infty}$
is nilpotent.

\section{\label{sec:someElementaryExample}Some elementary examples}

The elementary examples presented in this section want to show, in
a few rows, the simplicity of the algebraic calculus of nilpotent
infinitesimals. Here {}``simplicity'' means that the dialectic with
the corresponding informal calculations, used e.g. in engineering
or in physics, is really faithful. The importance of this dialectic
can be glimpsed both as a proof of the flexibility of the new language,
but also for researches in artificial intelligence like automatic
differentiation theories (see e.g. \citet{Gri} and references therein).
Last but not least, it may also be important for didactic or historical
researches. 
\begin{enumerate}
\item \noindent \textbf{Commutation of differentiation and integration.}
This example derives from \citet{Koc,Lav}. Suppose we want \emph{to
discover} the derivative of the function \[
g(x):=\int_{\alpha(x)}^{\beta(x)}f(x,t)\diff{t}\qquad\forall x\in\R\]
 where $\alpha$, $\beta$ and $f$ are smooth functions. We can see
$g$ as a composition of smooth functions, hence we can apply the
derivation formula, i.e. Theorem \ref{thm:DerivationFormula}: \begin{align*}
g(x+h)= & \int_{\alpha(x+h)}^{\beta(x+h)}f(x+h,t)\diff{t}=\\
= & \int_{\alpha(x)+h\alpha'(x)}^{\alpha(x)}f(x,t)\diff{t}+h\cdot\int_{\alpha(x)+h\alpha'(x)}^{\alpha(x)}\frac{\partial f}{\partial x}(x,t)\diff{t}+\\
{} & +\int_{\alpha(x)}^{\beta(x)}f(x,t)\diff{t}+h\cdot\int_{\alpha(x)}^{\beta(x)}\frac{\partial f}{\partial x}(x,t)\diff{t}+\\
{} & +\int_{\beta(x)}^{\beta(x)+h\beta'(x)}f(x,t)\diff{t}+h\cdot\int_{\beta(x)}^{\beta(x)+h\beta'(x)}\frac{\partial f}{\partial x}(x,t)\diff{t}.\end{align*}
 Now we use $h^{2}=0$ to obtain e.g. (see Corollary \ref{cor:derivationFormulaForIntegrals_1D}):
\[
h\cdot\int_{\alpha(x)+h\alpha'(x)}^{\alpha(x)}\frac{\partial f}{\partial x}(x,t)\diff{t}=-h^{2}\cdot\alpha'(x)\cdot\frac{\partial f}{\partial x}(\alpha(x),t)=0\]
 and \[
\int_{\alpha(x)+h\alpha'(x)}^{\alpha(x)}f(x,t)\diff{t}=-h\cdot\alpha'(x)\cdot f(\alpha(x),t).\]
Calculating in an analogous way similar terms we finally obtain the
well known conclusion. Note that the final formula comes out by itself
so that we have \emph{{}``discovered''} it and not simply we have
proved it. From the point of view of artificial intelligence or from
the didactic point of view, surely this discovering is not a trivial
result.
\item \textbf{Circle of curvature.} A simple application of the infinitesimal
Taylor's formula is the parametric equation for the circle of curvature,
that is the circle with second order osculation with a curve $\gamma:[0,1]\freccia\R^{3}$.
In fact if $r\in(0,1)$ and $\dot{\gamma}_{r}$ is a unit vector,
from the second order infinitesimal Taylor's formula we have \begin{equation}
\forall h\in D_{2}\pti\gamma(r+h)=\gamma_{r}+h\,\dot{\gamma}_{r}+\frac{h^{2}}{2}\,\ddot{\gamma}_{r}=\gamma_{r}+h\,\vec{t}_{r}+\frac{h^{2}}{2}c_{r}\,\vec{n}_{r}\label{eq:_CircleOfCurvature}\end{equation}
 where $\vec{n}$ is the unit normal vector, $\vec{t}$ is the tangent
one and $c_{r}$ the curvature. But once again from Taylor's formula
we have $\sin(ch)=ch$ and $\cos(ch)=1-\frac{c^{2}h^{2}}{2}.$ Now
it suffices to substitute $h$ and $\frac{h^{2}}{2}$ from these formulas
into (\ref{eq:_CircleOfCurvature}) to obtain the conclusion \[
\forall h\in D_{2}\pti\gamma(r+h)=\left(\gamma_{r}+\frac{\vec{n}_{r}}{c_{r}}\right)+\frac{1}{c_{r}}\cdot\left[\sin(c_{r}h)\vec{t}_{r}-\cos(c_{r}h)\vec{n}_{r}\right].\]
 In a similar way we can prove that any $f\in\Cc^{\infty}(\R,\R)$
can be written $\forall h\in D_{k}$ as \[
f(h)=\sum_{n=0}^{k}a_{n}\cdot\cos(nh)+\sum_{n=0}^{k}b_{n}\cdot\sin(nh),\]
so that now the idea of the Fourier series comes out in a natural
way. 
\item \noindent \textbf{Schwarz's theorem.} Using nilpotent infinitesimals
we can obtain a simple and meaningful proof of Schwarz's theorem.
This simple example aims to show how to manage some differences between
our setting and Synthetic Differential Geometry (see \citet{Koc,Lav,Mo-Re}).
Let $f:V\freccia E$ be a $\Cc^{2}$ function between spaces of type
$V=\R^{m}$, $E=\R^{n}$ (in subsequent chapters we will see that
the same proof is still valid for Banach spaces too) and $a\in V$,
we want to prove that ${\rm d}^{2}{f}(a):V\times V\freccia E$ is
symmetric. Take\begin{align*}
{} & k\in D_{2}\\
{} & h,j\in\D_{\infty}\\
{} & jkh\in D_{\ne0}\end{align*}
(e.g. we can take $k_{t}=\diff{t}_{2},h_{t}=j_{t}=\diff{t}_{4}$ so
that $jkh=\diff{t}$, see also Theorem \ref{thm:productOfPowers}).
Using $k\in D_{2}$ and Lemma \ref{lem:TaylorForGenericBasePoint},
we have \begin{equation}
\begin{split}j\cdot f(x & +hu+kv)=\\
 & =j\cdot\left[f(x+hu)+k\,\partial_{v}f(x+hu)+\frac{k^{2}}{2}\partial_{v}^{2}f(x+hu)\right]\\
 & =j\cdot f(x+hu)+jk\cdot\partial_{v}f(x+hu)\end{split}
\label{eq: x+h+k}\end{equation}
where we used the fact that $k^{2}\in D$ and $j$ infinitesimal imply
$jk^{2}=0$. Now we consider that $jkh\in D$ so that any product
of type $jkhi$ is zero for every $i\in D_{\infty}$, so we obtain
\begin{equation}
jk\cdot\partial_{v}f(x+hu)=jk\cdot\partial_{v}f(x)+jkh\cdot\partial_{u}(\partial_{v}f)(x).\label{eq: jkh}\end{equation}
But $k\in D_{2}$ and $jk^{2}=0$ hence \[
j\cdot f(x+kv)-j\cdot f(x)=jk\cdot\partial_{v}f(x).\]
Substituting this in \eqref{eq: jkh} and hence in \eqref{eq: x+h+k}
we obtain \begin{equation}
\begin{split} & j\cdot\left[f(x+hu+kv)-f(x+hu)-f(x+kv)+f(x)\right]=\\
 & =jkh\cdot\partial_{u}(\partial_{v}f)(x).\end{split}
\label{eq: SecondOrderIncrementalRatio}\end{equation}
The left hand side of this equality is symmetric in $u,v$, hence
changing them we have \[
jkh\cdot\partial_{u}(\partial_{v}f)(x)=jkh\cdot\partial_{v}(\partial_{u}f)(x)\]
and thus we obtain the conclusion because $jkh\ne0$ and $\partial_{u}(\partial_{v}f)(x)$,
$\partial_{v}(\partial_{u}f)(x)\in E$. From (\ref{eq: SecondOrderIncrementalRatio})
it follows directly the classical limit relation \[
\lim_{t\to0^{+}}\frac{f(x+h_{t}u+k_{t}v)-f(x+h_{t}u)-f(x+k_{t}v)+f(x)}{h_{t}k_{t}}=\partial_{u}\partial_{v}f(x)\]

\item \textbf{Electric dipole}. In elementary physics, an electric dipole
is usually defined as {}``\emph{a pair of charges with opposite sign
placed at a distance $d$ very less than the distance $r$ from the
observer''}.

Conditions like $r\gg d$ are frequently used in Physics and very
often we obtain a correct formalization if we ask $d\in\ER$ infinitesimal
but $r\in\R\setminus\{0\}$, i.e. $r$ finite. Thus we can define
an electric dipole as a pair $(p_{1},p_{2})$ of electric particles,
with charges of equal intensity but with opposite sign such that their
mutual distance at every time $t$ is a first order infinitesimal:
\begin{equation}
\forall t\pti\vert p_{1}(t)-p_{2}(t)\vert=:\vert\vec{d}_{t}\vert=:d_{t}\in D.\label{eq:_DefinitionOfDipole}\end{equation}
 In this way we can calculate the potential at the point $x$ using
the properties of $D$ and using the hypothesis that $r$ is finite
and not zero. In fact we have \[
\varphi(x)=\frac{q}{4\pi\epsilon_{0}}\cdot\left(\frac{1}{r_{1}}-\frac{1}{r_{2}}\right)\qquad\qquad\vec{r_{i}}:=x-p_{i}\]
 and if $\vec{r}:=\vec{r}_{2}-\frac{\vec{d}}{2}$ then \[
\frac{1}{r_{2}}=\left(r^{2}+\frac{d^{2}}{4}+\vec{r}\boldsymbol{\cdot}\vec{d}\right)^{-1/2}=r^{-1}\cdot\left(1+\frac{\vec{r}\boldsymbol{\cdot}\vec{d}}{r^{2}}\right)^{-1/2}\]
 because for (\ref{eq:_DefinitionOfDipole}) $d^{2}=0$. For our hypotheses
on $d$ and $r$ we have that ${\displaystyle \frac{\vec{r}\boldsymbol{\cdot}\vec{d}}{r^{2}}\in D}$
hence from the derivation formula \[
\left(1+\frac{\vec{r}\boldsymbol{\cdot}\vec{d}}{r^{2}}\right)^{-1/2}=1-\frac{\vec{r}\boldsymbol{\cdot}\vec{d}}{2r^{2}}\]
 In the same way we can proceed for $1/r_{1}$, hence: \begin{align*}
\varphi(x) & =\frac{q}{4\pi\epsilon_{0}}\cdot\frac{1}{r}\cdot\left(1+\frac{\vec{r}\boldsymbol{\cdot}\vec{d}}{2r^{2}}-1+\frac{\vec{r}\boldsymbol{\cdot}\vec{d}}{2r^{2}}\right)=\\
 & =\frac{q}{4\pi\epsilon_{0}}\cdot\frac{\vec{r}\boldsymbol{\cdot}\vec{d}}{r^{3}}\end{align*}
 The property $d^{2}=0$ is also used in the calculus of the electric
field and for the moment of momentum. 

\item \textbf{Newtonian limit in Relativity.} Another example in which we
can formalize a condition like $r\gg d$ using the previous ideas
is the Newtonian limit in Relativity; in it we can suppose to have 

\begin{itemize}
\item $\,\forall t\pti v_{t}\in D_{2}\e{and}c\in\R$ 
\item $\,\forall x\in M_{4}\pti g_{ij}(x)=\eta_{ij}+h_{ij}(x)\e{with}h_{ij}(x)\in D.$ 
\end{itemize}
where $\left(\eta_{ij}\right)_{ij}$ is the matrix of the Minkowski's
metric. This conditions can be interpreted as $v_{t}\ll c$ and $h_{ij}(x)\ll1$
(low speed with respect to the speed of light and weak gravitational
field). In this way we have, e.g. the equalities: \[
\frac{1}{\sqrt{{\displaystyle 1-\frac{v^{2}}{c^{2}}}}}=1+\frac{v^{2}}{2c^{2}}\ee{and}\sqrt{1-h_{44}(x)}=1-\frac{1}{2}\, h_{44}(x).\]

\item \textbf{Linear differential equations.} Let \begin{gather*}
L(y):=A_{\sss0}\frac{\diff{}^{\sss N}y}{\diff{}t^{\sss N}}+\ldots+A_{\sss N-1}\frac{\diff{}y}{\diff{}t}+A_{\sss N}\cdot y=0\end{gather*}
 be a linear differential equation with constant coefficients. Once
again we want \emph{to discover} independent solutions in case the
characteristic polynomial has multiple roots e.g. \[
(r-r_{\sss1})^{2}\cdot(r-r_{\sss3})\cdot\ldots\cdot(r-r_{\sss N})=0.\]
 The idea is that in $\ER$ we have $(r-r_{1})^{2}=0$ also if $r=r_{1}+h$
with $h\in D$. Thus $y(t)={\rm e}^{(r_{1}+h)t}$ is a solution too.
But ${\rm e}^{(r_{1}+h)t}={\rm e}^{r_{1}t}+ht\cdot{\rm e}^{r_{1}t}$,
hence \begin{align*}
L\left[{\rm e}^{(r_{1}+h)t}\right] & =0\\
{} & =L\left[{\rm e}^{r_{1}t}+ht\cdot{\rm e}^{r_{1}t}\right]\\
{} & =L\left[{\rm e}^{r_{1}t}\right]+h\cdot L\left[t\cdot{\rm e}^{r_{1}t}\right]\end{align*}
 We obtain $L\left[t\cdot{\rm e}^{r_{1}t}\right]=0$, that is $y_{1}(t)=t\cdot{\rm e}^{r_{1}t}$
must be a solution. Using $k$-th order infinitesimals we can deal
with other multiple roots in a similar way.
\end{enumerate}
We think that these elementary examples are able to show that some
results that frequently may appear as unnatural in a standard context,
using Fermat reals may be even discovered, even by suitably designed
algorithm.

\chapter{Order relation\label{cha:orderRelation}}

\section{Infinitesimals and order properties\label{sec:Infinitesimals-and-order-properties}}

Like in other disciplines, also in mathematics the layout of a work
reflects the personal philosophical ideas of the authors. In particular
the present work is based on the idea that a good mathematical theory
is able to construct a good dialectic between formal properties, proved
in the theory, and their informal interpretations. The dialectic has
to be, as far as possible, in both directions: theorems proved in
the theory should have a clear and useful intuitive interpretation
and, on the other hand, the intuition corresponding to the theory
has to be able to suggest true sentences, i.e. conjectures or sketch
of proofs that can then be converted into rigorous proofs.

In a theory of new numbers, like the present one about Fermat reals,
the introduction of an order relation can be a hard test of the excellence
of this dialectic between formal properties and their informal interpretations.
Indeed if we introduce a new ring of numbers (like $\ER$) extending
the real field $\R$, we want that the new order relation, defined
on the new ring, will extend the standard one on $\R$. This extension
naturally leads to the wish of findings a geometrical representation
of the new numbers, in accord with the above principle of having a
good formal/informal dialectic.

For example, on the one hand in NSA the order relation on $\hyperR$
has the best formal properties among all the theories of actual infinitesimals.
On the other hand, the dialectic of these properties with the informal
interpretations is not always good, due to the use of, e.g., an ultrafilter
in the construction of $\hyperR$. Indeed, in an ultrafilter on $\N$
we can always find a highly non constructive set $A\subset\N$; any
sequence of reals $x:\N\freccia\R$ which is constant to 1 on $A$
is strictly greater than 0 in $\hyperR$, but it seems not easy to
give neither an intuitive interpretation nor a clear and meaningful
geometric representation of the relation $x>0$ in $\hyperR$. In
fact, it is also for motivations of this type that some approaches
to give a constructive definition of a field similar to $\hyperR$
have been attempted (see e.g. \cite{Pal1, Pal2, Pal3} and references
therein).

\medskip{}

In SDG we have a preorder relation (i.e. a reflexive and transitive
relation, which is not necessarily anti-symmetric) with very poor
properties only. Nevertheless, the works developed in SDG (see e.g.
\cite{Lav}) exhibits that meaningful geometric results can be obtained
in infinite dimensional spaces, even if the order properties of the
ground base ring are not so rich. Once again, the dialectic between
formal properties and their intuitive interpretations represents a
hard test for SDG too. E.g. it seems not so easy to interpret intuitively
that every infinitesimal $h$ in SDG verifies both $h\ge0$ and $h\le0$.
The lack of a total order, i.e. of the trichotomy law\begin{equation}
x<y\e{or}y<x\e{or}x=y\label{eq:trichotomyLaw}\end{equation}
 makes really difficult, or even impossible, to have a geometrical
representation of the infinitesimals of SDG.

\medskip{}

\begin{figure}[h]
\begin{centering}
\includegraphics[scale=0.4]{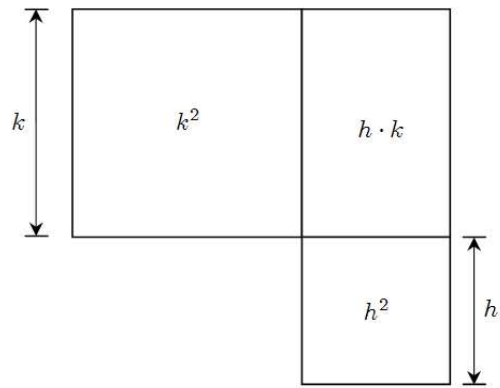} 
\par\end{centering}

\caption{{\small How to guess that $h\cdot k=0$ for two first order infinitesimals
$h$, $k\in D$}}

\label{fig:howToGuessThatTheProductIsZero}
\end{figure}

We want to start this section showing that in our setting there is
a strong connection between some order properties and some algebraic
properties. In particular, we will show that it is not possible to
have good order properties and at the same time a uniqueness without
limitations in the derivation formula (see the discussion starting
Chapter \ref{cha:equalityUpTo_k-thOrder}). We know that in $\ER$
the product of any two first order infinitesimals $h$, $k\in D$
is always zero: $h\cdot k=0$, and a consequence of this property
is that we have some limitations in the uniqueness of the derivation
formula, and for these reasons we introduce the notion of equality
up to $k$-th order infinitesimals (see Chapter \ref{cha:equalityUpTo_k-thOrder}).
In the following theorem we can see that the property $h\cdot k=0$
is a general consequence if we suppose to have a total order on $D$.
The idea of this theorem can be glimpsed at from the Figure \ref{fig:howToGuessThatTheProductIsZero},
where it is represented that if we neglect $h^{2}$ and $k^{2}$ because
we consider them zero, then we have strong reasons to expect that
also $h\cdot k$ will be zero%

\noindent From this picture comes the idea to find a formal demonstration
based on the implication\[
h,k\ge0\e{,}h\le k\then0\le hk\le k^{2}=0\]

\noindent All these ideas conduct toward the following theorem. 
\begin{thm}
\noindent \label{thm:orderImpliesProductOfFirstOrderInfinitesimalsIsZero}Let
$(R,\le)$ be a generic ordered ring and $D\subseteq R$ a subset
of this ring, such that
\begin{enumerate}
\item $0\in D$ 
\item $\forall h\in D\pti h^{2}=0$ and $-h\in D$ 
\item $(D,\le)$ is a total order 
\end{enumerate}
\noindent then\begin{equation}
\forall h,k\in D\pti h\cdot k=0\label{eq:productOfFirstOrderInfinitesimalsIsZero}\end{equation}

\end{thm}
\noindent This theorem implies that if we want a total order in our
theory of infinitesimal numbers, and if in this theory we consider
$D=\{h\,|\, h^{2}=0\}$, then we must accept that the product of any
two elements of $D$ must be zero. For example, if we think that a
geometric representation of infinitesimals cannot be possible if we
do not have, at least, the trichotomy law, then in this theory we
must also have that the product of two first order infinitesimals
is zero. Finally, because in SDG property \eqref{eq:productOfFirstOrderInfinitesimalsIsZero}
is false, this theorem also implies that in SDG it is not possible
to define a total order (and not only a preorder) on the set $D$
of first order infinitesimals compatible with the ring operations.

\noindent \bigskip{}

\noindent \textbf{Proof of Theorem \ref{thm:orderImpliesProductOfFirstOrderInfinitesimalsIsZero}:}

\noindent Let $h$, $k\in D$ be two elements of the subset $D$.
By hypotheses $0,$ $-h$, $-k\in D$, hence all these elements are
comparable with respect to the order relation $\le$, because, by
hypotheses this relation is total (i.e. \eqref{eq:trichotomyLaw}
is true). E.g.\[
h\le k\e{or}k\le h\]
 We will consider only the case $h\le k$, because analogously we
can deal with the case $k\le h$, simply exchanging everywhere $h$
with $k$ and vice versa.

\medskip{}

\noindent \emph{First sub-case}: \emph{$k\ge0$.\ \ \ }By multiplying
both sides of $h\le k$ by $k\ge0$ we obtain \begin{equation}
hk\le k^{2}\label{eq:hk_le_kSquare}\end{equation}
 If $h\ge0$ then, multiplying by $k\ge0$ we have $0\le hk$, so
from \eqref{eq:hk_le_kSquare} we have $0\le hk\le k^{2}=0$, and
hence $hk=0$.

\noindent If $h\le0$ then, multiplying by $k\ge0$ we have\begin{equation}
hk\le0\label{eq:hk_le0}\end{equation}
 If, furthermore, $h\ge-k$, then multiplying by $k\ge0$ we have
$hk\ge-k^{2}$, hence form \eqref{eq:hk_le0} $0\ge hk\ge-k^{2}=0$,
hence $hk=0$.

\noindent If, otherwise, $h\le-k$, then multiplying by $-h\ge0$
we have $-h^{2}=0\le hk\le0$ from \eqref{eq:hk_le0}, hence $hk=0$.
This concludes the discussion of the case $k\ge0$.

\medskip{}

\noindent \emph{Second sub-case}: $k\le0$.\ \ \ In this case we
have $h\le k\le0$. Multiplying both inequalities by $h\le0$ we obtain
$h^{2}=0\ge hk\ge0$ and hence $hk=0$.\qedNoNewLine

Property \eqref{eq:productOfFirstOrderInfinitesimalsIsZero} is incompatible
with the uniqueness in a possible derivation formula like \begin{equation}
\exists!\, m\in R\pti\forall h\in D\pti f(h)=f(0)+h\cdot m\label{eq:derivationFormulaInAGenericRing}\end{equation}
 framed in the ring $R$ of Theorem \ref{thm:orderImpliesProductOfFirstOrderInfinitesimalsIsZero}.
In fact, if $a$, $b\in D$ are two elements of the subset $D\subseteq R$,
then both $a$ and $b$ play the role of $m\in R$ in \eqref{eq:derivationFormulaInAGenericRing}
for the linear function\[
f:h\in D\mapsto h\cdot a=0\in R\]
 So, if the derivation formula \eqref{eq:derivationFormulaInAGenericRing}
applies to linear functions (or less, to constant functions), the
uniqueness part of this formula cannot hold in the ring $R$.

In the next section we will introduce a natural and meaningful total
order relation on $\ER$. Therefore, the previous Theorem \ref{thm:orderImpliesProductOfFirstOrderInfinitesimalsIsZero}
strongly motivate that for the ring of Fermat reals $\ER$ we must
have that the product of two first order infinitesimals must be zero
and hence, that for the derivation formula in $\ER$ the uniqueness
cannot hold in its strongest form. Since we will also see that the
order relation permits to have a geometric representation of Fermat
reals, we can summarize the conclusions of this section saying that
the uniqueness in the derivation formula is incompatible with a natural
geometric interpretation of Fermat reals and hence with a good dialectic
between formal properties and informal interpretations in this theory.

\section{\label{sec:OrderRelation}Order relation}

From the previous sections one can draw the conclusion that the ring
of Fermat reals $\ER$ is essentially {}``the little-oh'' calculus.
But, on the other hand the Fermat reals give us more flexibility than
this calculus: working with $\ER$ we do not have to bother ourselves
with remainders made of {}``little-oh'', but we can neglect them
and use the useful algebraic calculus with nilpotent infinitesimals.
But thinking the elements of $\ER$ as new numbers, and not simply
as {}``little-oh functions'', permits to treat them in a different
and new way, for example to define on them an order relation with
a clear geometrical interpretation%
\footnote{We will see that this order relation is different from the order of
infinite or infinitesimal originally introduced by P. Du Bois-Reymond
(see \cite{Har}).%
}.

First of all, let us introduce the useful notation\[
\forall^{0}t\ge0\pti\mathcal{P}(t)\]

\noindent and we will read the quantifier $\forall^{0}t\ge0$ saying
\emph{{}``for every $t\ge0$ (sufficiently) small}'', to indicate
that the property $\mathcal{P}(t)$ is true for all $t$ in some right%
\footnote{We recall that, by Definition \ref{def:LittleOhPolynomials}, our
little-oh polynomials are always defined on $\R_{\ge0}$%
} neighborhood of $t=0$, i.e.\[
\exists\delta>0\pti\forall t\in[0,\delta)\pti\mathcal{P}(t)\]

\noindent The first heuristic idea to define an order relation is
the following\[
x\le y\iff x-y\le0\iff\exists z\pti z=0\e{in}\ER\e{and}x-y\le z\]
 More precisely, if $x$, $y\in\ER$ are two little-oh polynomials,
we want to ask locally that%
\footnote{We recall that, to simplify the notations, we do not use equivalence
classes as elements of $\ER$ but directly little-oh functions. The
only notion of equality between little-oh functions is, of course,
the equivalence relation defined in Definition \ref{def:equalityInFermatReals}
and, as usual, we must always prove that our relations between little-oh
polynomials are well defined.%
} $x_{t}$ is less than or equal to $y_{t}$, but up to a $o(t)$ for
$t\to0^{+}$, where the little-oh function $o(t)$ depends on $x$
and $y$. Formally: 
\begin{defn}
\noindent \label{def:orderRelation}Let $x$, $y\in\ER$, then we
say\[
x\le y\]

\noindent iff we can find $z\in\ER$ such that $z=0$ in $\ER$ and\[
\forall^{0}t\ge0\pti x_{t}\le y_{t}+z_{t}\]

\end{defn}
\noindent Recall that $z=0$ in $\ER$ is equivalent to $z_{t}=o(t)$
for $t\to0^{+}$. It is immediate to see that we can equivalently
define $x\le y$ if and only if we can find $x'=x$ and $y'=y$ in
$\ER$ such that $x_{t}\le y_{t}$ for every $t$ sufficiently small.
From this it also follows that the relation $\le$ is well defined
on $\ER$, i.e. if $x'=x$ and $y'=y$ in $\ER$ and $x\le y$, then
$x'\le y'$. As usual we will use the notation $x<y$ for $x\le y$
and $x\ne y$.
\begin{thm}
\noindent The relation $\le$ is an order, i.e. is reflexive, transitive
and anti-symmetric; it extends the order relation of $\R$ and with
it $(\ER,\le)$ is an ordered ring. Finally the following sentences
are equivalent:
\begin{enumerate}
\item $h\in D_{\infty}$, i.e. $h$ is an infinitesimal 
\item $\forall r\in\R_{>0}\pti-r<h<r$ 
\end{enumerate}
\end{thm}
\noindent Hence an infinitesimal can be thought of as a number with
standard part zero, or as a number smaller than every standard positive
real number and greater than every standard negative real number (thus
it has in this sense the same property as an infinitesimal both in
NSA and in SDG (in the latter case with real numbers of type $\frac{1}{n}$
($n\in\N_{>0}$) only).

\bigskip{}

\noindent \textbf{Proof:} It is immediate to prove that the relation
is reflexive. To prove transitivity, if $x\le y$ and $y\le w$, then
we have\[
\forall^{0}t\ge0\pti x_{t}\le y_{t}+z_{t}\e{and}\forall^{0}t\ge0\pti y_{t}\le w_{t}+z'_{t}\]
 and these imply\[
\forall^{0}t\ge0\pti x_{t}\le y_{t}+z_{t}\le w_{t}+z_{t}+z'_{t}\]

\noindent showing that $x\le w$. To prove that it is also anti-symmetric,
take $x\le y$ and $y\le x$, then we have\begin{align}
x_{t} & \le y_{t}+z_{t}\quad\forall t\in[0,\delta_{1})\label{eq:antisymmetricDelta1}\\
y_{t} & \le x_{t}+z'_{t}\quad\forall t\in[0,\delta_{2})\label{eq:antisymmetricDelta2}\end{align}
 \[
\lim_{t\to0^{+}}\frac{z_{t}}{t}=0\e{and}\lim_{t\to0^{+}}\frac{z'_{t}}{t}=0\]
 because $z$ and $z'$ are equal to zero in $\ER$, that is are $o(t)$
for $t\to0^{+}$. Hence from \eqref{eq:antisymmetricDelta1} and \eqref{eq:antisymmetricDelta2}
for $\delta:=\min\{\delta_{1},\delta_{2}\}$ we have\[
-\frac{z'_{t}}{t}\le\frac{x_{t}-y_{t}}{t}\le\frac{z_{t}}{t}\quad\forall t\in[0,\delta)\]
 and hence $\lim_{t\to0^{+}}\frac{x_{t}-y_{t}}{t}=0$, that is $x=y$
in $\ER$.

\bigskip{}

If $r$, $s\in\R$ and $r\le s$ as real numbers, then it suffices
to take $z_{t}=0$ for every $t\ge0$ in the Definition \ref{def:orderRelation}
to obtain that $r\le s$ in $\ER$ too. Vice versa if $r\le s$ in
$\ER$, then for some $z=0$ in $\ER$ we have\[
\forall^{0}t\ge0\pti r\le s+z_{t}\]
 and hence for $t=0$ we have $r\le s$ in $\R$ because $z=0$ and
hence $z_{0}=0$. This proves that the order relation $\le$ defined
in $\ER$ extends the order relation on $\R$.

\bigskip{}

The relationships between the ring operations and the order relation
can be stated as\begin{align*}
x & \le y\then x+w\le y+w\\
x & \le y\then-x\ge-y\\
x & \le y\e{and}w\ge0\then x\cdot w\le y\cdot w\end{align*}
 The first two are immediate consequences of the Definition \ref{def:orderRelation}.
To prove the last one, let us suppose that\begin{align}
x_{t} & \le y_{t}+z_{t}\quad\forall^{0}t\ge0\label{eq:xLessThan_yPlus_z}\\
w_{t} & \ge z'_{t}\quad\forall^{0}t\ge0\nonumber \end{align}
 then $w_{t}-z'_{t}\ge0$ for every $t$ small and hence from \eqref{eq:xLessThan_yPlus_z}\[
x_{t}\cdot(w_{t}-z'_{t})\le y_{t}\cdot(w_{t}-z'_{t})+z_{t}\cdot(w_{t}-z'_{t})\quad\forall^{0}t\ge0\]
 from which it follows\[
x_{t}\cdot w_{t}\le y_{t}\cdot w_{t}+(-x_{t}z'_{t}-y_{t}z'_{t}+z_{t}w_{t}-z_{t}z'_{t})\quad\forall^{0}t\ge0\]
 But $-xz'-yz'+zw-zz'=0$ in $\ER$ because $z=0$ and $z'=0$ and
hence the conclusion follows.

\bigskip{}

Finally we know (see Definition \ref{def:idealD_a}) that $h\in D_{\infty}$
if and only if $\st{h}=0$ and this is equivalent to\begin{equation}
\forall r\in\R_{>0}\pti-r<\st{h}<r\label{eq:infinitesimalAndRealsWithStandardPart}\end{equation}
 But if, e.g., $\st{h}<r$, then\[
\forall^{0}t\ge0\pti h_{t}\le r\]
 because the function $t\to h_{t}$ is continuous, and hence we also
have $h\le r$ in $\ER$. Analogously, from \eqref{eq:infinitesimalAndRealsWithStandardPart}we
can prove that $-r\le h$ for all $r\in\R_{>0}$. Of course $r\notin D_{\infty}$
if $r\in\R$, so it cannot be that $h=r$.

Vice versa if \[
\forall r\in\R_{>0}\pti-r<h<r\]
 then, e.g., $h_{t}\le r+z_{t}$ for $t$ small. Hence, for $t=0$
we have $-r\le\st{h}=h_{0}\le r$ for every $r>0$, and so $\st{h}=0$.\qedNoNewLine 
\begin{example*}
We have e.g. $\diff{t}>0$ and $\diff{t_{2}}-3\diff{t}>0$ because
for $t\ge0$ sufficiently small $t^{1/2}>3t$ and hence\[
t^{1/2}-3t>0\quad\forall^{0}t\ge0\]
 From examples like these ones we can guess that our little-oh polynomials
are always locally comparable with respect to pointwise order relation,
and this is the first step to prove that for our order relation the
trichotomy law holds. In the following statement we will use the notation
$\forall^{0}t>0:\mathcal{P}(t)$, that naturally means\[
\forall^{0}t\ge0\pti t\ne0\then\mathcal{P}(t)\]

\noindent where $\mathcal{P}(t)$ is a generic property depending
on $t$.\end{example*}
\begin{lem}
\label{lem:pointwiseComparison}Let $x$, $y\in\ER$, then
\begin{enumerate}
\item \label{enu:standardPartAreDifferent}$\st{x}<\st{y}\then\forall^{0}t\ge0\pti x_{t}<y_{t}$ 
\item \label{enu:standardPartAreEqual}If $\st{x}=\st{y}$, then\[
\left(\forall^{0}t>0\pti x_{t}<y_{t}\right)\ \ \text{or}\ \ \left(\forall^{0}t>0\pti x_{t}>y_{t}\right)\ \ \text{or}\ \ \left(x=y\e{in}\ER\right)\]

\end{enumerate}
\end{lem}
\noindent \textbf{Proof:}

\noindent \emph{\ref{enu:standardPartAreDifferent}}.)\ \ \ Let
us suppose that $\st{x}<\st{y}$, then the continuous function $t\ge0\mapsto y_{t}-x_{t}\in\R$
assumes the value $y_{0}-x_{0}>0$ hence is locally positive, i.e.\[
\forall^{0}t\ge0\pti x_{t}<y_{t}\]
 \emph{\ref{enu:standardPartAreEqual}}.)\ \ \ Now let us suppose
that $\st{x}=\st{y}$, and introduce a notation for the potential
decompositions of $x$ and $y$ (see Definition \ref{def:potentialDecomposition}).
From the definition of equality in $\ER$, we can always write \begin{align*}
x_{t} & =\st{x}+\sum_{i=1}^{N}\alpha_{i}\cdot t^{a_{i}}+z_{t}\quad\forall t\ge0\\
y_{t} & =\st{y}+\sum_{j=1}^{M}\beta_{j}\cdot t^{b_{j}}+w_{t}\quad\forall t\ge0\end{align*}
 where $x=\st{x}+\sum_{i=1}^{N}\alpha_{i}\cdot t^{a_{i}}$ and $y=\st{y}+\sum_{j=1}^{M}\beta_{j}\cdot t^{b_{j}}$
are the potential decompositions of $x$ and $y$ (hence $0<\alpha_{i}<\alpha_{i+1}\le1$
and $0<\beta_{j}<\beta_{j+1}\le1$), whereas $w$ and $z$ are little-oh
polynomials such that $z_{t}=o(t)$ and $w_{t}=o(t)$ for $t\to0^{+}$.

\noindent \smallskip{}

\noindent \emph{Case}: $a_{1}<b_{1}$\ \ \ In this case the least
power in the two decompositions is $\alpha_{1}\cdot t^{a_{1}}$, and
hence we expect that the second alternative of the conclusion is the
true one if $\alpha_{1}>0$, otherwise the first alternative will
be the true one if $\alpha_{1}<0$ (recall that always $\alpha_{i}\ne0$
in a decomposition). Indeed, let us analyze, for $t>0$, the condition
$x_{t}<y_{t}$: the following formulae are all equivalent to it\[
\ \sum_{i=1}^{N}\alpha_{i}\cdot t^{a_{i}}<\sum_{j=1}^{N}\beta_{j}\cdot t^{b_{j}}+w_{t}-z_{t}\]
 \[
t^{a_{1}}\cdot\left[\alpha_{1}+\sum_{i=2}^{N}\alpha_{i}\cdot t^{a_{i}-a_{1}}\right]<\ t^{a_{1}}\cdot\left[\sum_{j=1}^{N}\beta_{j}\cdot t^{b_{j}-a_{1}}+(w_{t}-z_{t})\cdot t^{-a_{1}}\right]\]
 \[
\alpha_{1}+\sum_{i=2}^{N}\alpha_{i}\cdot t^{a_{i}-a_{1}}<\sum_{j=1}^{N}\beta_{j}\cdot t^{b_{j}-a_{1}}+(w_{t}-z_{t})\cdot t^{-a_{1}}\]
 Therefore, let us consider the function\[
f(t):=\sum_{j=1}^{N}\beta_{j}\cdot t^{b_{j}-a_{1}}+(w_{t}-z_{t})\cdot t^{-a_{1}}-\alpha_{1}-\sum_{i=2}^{N}\alpha_{i}\cdot t^{a_{i}-a_{1}}\quad\forall t\ge0\]
 We can write\[
(w_{t}-z_{t})\cdot t^{-a_{1}}=\frac{w_{t}-z_{t}}{t}\cdot t^{1-a_{1}}\]
 and $\frac{w_{t}-z_{t}}{t}\to0$ as $t\to0^{+}$ because $w_{t}=o(t)$
and $z_{t}=o(t)$. Furthermore, $a_{1}\le1$ hence $t^{1-a_{1}}$
is bounded in a right neighborhood of $t=0$. Therefore, $(w_{t}-z_{t})\cdot t^{-a_{1}}\to0$
and the function $f$ is continuous at $t=0$ too, because $a_{i}<a_{i}$
and $a_{1}<b_{1}<b_{j}$. By continuity, the function $f$ is locally
strictly positive if and only if $f(0)=-\alpha_{1}>0$, hence\begin{align*}
\left(\forall^{0}t>0\pti x_{t}<y_{t}\right) & \iff\alpha_{1}<0\\
\left(\forall^{0}t>0\pti x_{t}>y_{t}\right) & \iff\alpha_{1}>0\end{align*}
 \emph{Case}: $a_{1}>b_{1}$\ \ \ We can argue in an analogous
way with $b_{1}$ and $\beta_{1}$ instead of $a_{1}$ and $\alpha_{1}$.

\smallskip{}

\noindent \emph{Case}: $a_{1}=b_{1}$\ \ \ We shall exploit the
same idea used above and analyze the condition $x_{t}<y_{t}$. The
following are equivalent ways to express this condition\[
t^{a_{1}}\cdot\left[\alpha_{1}+\sum_{i=2}^{N}\alpha_{i}\cdot t^{a_{i}-a_{1}}\right]<t^{a_{1}}\cdot\left[\beta_{1}+\sum_{j=2}^{N}\beta_{j}\cdot t^{b_{j}-a_{1}}+(w_{t}-z_{t})\cdot t^{-a_{1}}\right]\]
 \[
\alpha_{1}+\sum_{i=2}^{N}\alpha_{i}\cdot t^{a_{i}-a_{1}}<\beta_{1}+\sum_{j=2}^{N}\beta_{j}\cdot t^{b_{j}-a_{1}}+(w_{t}-z_{t})\cdot t^{-a_{1}}\]
 Hence, exactly as we have demonstrated above, we can state that\begin{align*}
\alpha_{1}<\beta_{1} & \then\forall^{0}t>0\pti x_{t}<y_{t}\\
\alpha_{1}>\beta_{1} & \then\forall^{0}t>0\pti x_{t}>y_{t}\end{align*}

\noindent Otherwise $\alpha_{1}=\beta_{1}$ and we can restart with
the same reasoning using $a_{2}$, $b_{2}$, $\alpha_{2}$, $\beta_{2}$,
etc. If $N=M$, the number of addends in the decompositions, using
this procedure we can prove that\[
\forall t\ge0\pti x_{t}=y_{t}+w_{t}-z_{t},\]
 that is $x=y$ in $\ER$.s

It remains to consider the case, e.g., $N<M$. In this hypotheses,
using the previous procedure we would arrive at the following analysis
of the condition $x_{t}<y_{t}$: \[
0<\sum_{j>N}\beta_{j}\cdot t^{b_{j}}+w_{t}-z_{t}\]
 \[
0<t^{b_{N+1}}\cdot\Bigg[\beta_{N+1}+\sum_{j>N+1}\beta_{j}\cdot t^{b_{j}-b_{N+1}}+(w_{t}-z_{t})\cdot t^{-b_{N+1}}\Bigg]\]
 \[
0<\beta_{N+1}+\sum_{j>N+1}\beta_{j}\cdot t^{b_{j}-}{}^{b_{N+1}}+(w_{t}-z_{t})\cdot t^{-b_{N+1}}\]
 Hence\[
\beta_{N+1}>0\then\forall^{0}t>0\pti x_{t}<y_{t}\]

\[
\beta_{N+1}<0\then\forall^{0}t>0\pti x_{t}>y_{t}\]
 \qedWithFinalEq

\noindent This lemma can be used to find an equivalent formulation
of the order relation. 
\begin{thm}
\noindent \label{thm:equivalentFormulationForOrderRelation}Let $x$,
$y\in\ER$, then
\begin{enumerate}
\item \label{enu:equivalentFormulationFor_x_le_y}$x\le y\iff\left(\forall^{0}t>0\pti x_{t}<y_{t}\right)$
~~or~~ $\ \ (x=y$ in $\ER)$ 
\item \label{enu:equivalentFormulationFor_x_less_y}$x<y\iff\left(\forall^{0}t>0\pti x_{t}<y_{t}\right)$
~~and~~ $(x\ne y$ in $\ER)$ 
\end{enumerate}
\end{thm}
\noindent \textbf{Proof:}

\noindent \emph{\ref{enu:equivalentFormulationFor_x_le_y}}.) $\Rightarrow$\ \ \ If
$\st{x}<\st{y}$ then, from the previous Lemma \ref{lem:pointwiseComparison}
we can derive that the first alternative is true. If $\st{x}=\st{y}$,
then from Lemma \ref{lem:pointwiseComparison} we have\begin{equation}
\left(\forall^{0}t>0\pti x_{t}<y_{t}\right)\e{or}\left(x=y\e{in}\ER\right)\e{or}\left(\forall^{0}t>0\pti x_{t}>y_{t}\right)\label{eq:threeAlternatives}\end{equation}
 In the first two cases we have the conclusion. In the third case,
from $x\le y$ we obtain\begin{equation}
\forall^{0}t\ge0\pti x_{t}\le y_{t}+z_{t}\label{eq:x_le_yWith_z}\end{equation}
 with $z_{t}=o(t)$. Hence from the third alternative of \eqref{eq:threeAlternatives}
we have\[
0<x_{t}-y_{t}\le z_{t}\quad\forall^{0}t>0\]
 and hence $\lim_{t\to0^{+}}\frac{x_{t}-y_{t}}{t}=0$, i.e. $x=y$
in $\ER$.

\medskip{}

\noindent \emph{\ref{enu:equivalentFormulationFor_x_le_y}}.) $\Leftarrow$\ \ \ This
follows immediately from the reflexive property of $\le$ or from
the Definition \ref{def:orderRelation}.

\medskip{}

\noindent \emph{\ref{enu:equivalentFormulationFor_x_less_y}}.) $\Rightarrow$\ \ \ From
$x<y$ we have $x\le y$ and $x\ne y$, so the conclusion follows
from the previous \emph{\ref{enu:equivalentFormulationFor_x_le_y}}.

\medskip{}

\noindent \emph{\ref{enu:equivalentFormulationFor_x_less_y}}.) $\Leftarrow$\ \ \ From
$\forall^{0}t>0:x_{t}<y_{t}$ and from \emph{\ref{enu:equivalentFormulationFor_x_le_y}}.
it follows $x\le y$ and hence $x<y$ from the hypotheses $x\ne y$.\qedNoNewLine

\noindent Now we can prove that our order is total
\begin{cor}
\noindent Let $x$, $y\in\ER$, then in $\ER$ we have
\begin{enumerate}
\item \label{enu:trichotomyWithLessOrEqual}$x\le y\e{or}y\le x\e{or}x=y$ 
\item \label{enu:trichotomyWithLess}$x<y\e{or}y<x\e{or}x=y$ 
\end{enumerate}
\end{cor}
\noindent \textbf{Proof:}

\noindent \emph{\ref{enu:trichotomyWithLessOrEqual}}.)\ \ \ If
$\st{x}<\st{y}$, then from Lemma \ref{lem:pointwiseComparison} we
have $x_{t}<y_{t}$ for $t\ge0$ sufficiently small. Hence from Theorem
\ref{thm:equivalentFormulationForOrderRelation} we have $x\le y$.
We can argue in the same way if $\st{x}>\st{y}$. Also the case $\st{x}=\st{y}$
can be handled in the same way using \emph{\ref{enu:standardPartAreEqual}}.
of Lemma \ref{lem:pointwiseComparison}.

\medskip{}

\noindent \emph{\ref{enu:trichotomyWithLess}}.)\ \ \ This part
is a general consequence of the previous one. Indeed, if we have $x=y$,
then we have the conclusion. Otherwise we have $x\ne y$, and using
the previous \emph{\ref{enu:trichotomyWithLessOrEqual}}. we can deduce
strict inequalities from inequalities because $x\ne y$.\qedNoNewLine

From the proof of Lemma \ref{lem:pointwiseComparison} and from Theorem
\ref{thm:equivalentFormulationForOrderRelation} we can deduce the
following 
\begin{thm}
\label{thm:effectiveCriterionForTheOrder}Let $x$, $y\in\ER$. If
$\st{x}\ne\st{y}$, then\[
x<y\iff\st{x}<\st{y}\]

\noindent Otherwise, if $\st{x}=\st{y}$, then
\begin{enumerate}
\item If $\omega(x)>\omega(y)$, then $x>y$ iff $\st{x_{1}}>0$ 
\item If $\omega(x)=\omega(y)$, then\begin{align*}
\st{x_{1}}>\st{y_{1}} & \then x>y\\
\st{x_{1}}<\st{y_{1}} & \then x<y\end{align*}

\end{enumerate}
\end{thm}
\noindent This Theorem proves also some sentences about the order
relation anticipated in the Remark \ref{rem:differentMeaningForStandardAndPotentialOrder}. 
\begin{example*}
The previous Theorem gives an effective criterion to decide whe\-ther
$x<y$ or not. Indeed, if the two standard parts are different, then
the order relation can be decided on the basis of these standard parts
only. E.g. $2+\diff{t_{2}}>3\diff{t}$ and $1+\diff{t_{2}}<3+\diff{t}$.
\end{example*}
\noindent Otherwise, if the standard parts are equal, we firstly have
to look at the order and at the first standard parts, i.e. $\st{x_{1}}$
and $\st{y_{1}}$, which are the coefficients of the biggest infinitesimals
in the decompositions of $x$ and $y$. E.g. $3\diff{t_{2}}>5\diff{t}$,
and $\diff{t_{2}}>a\diff{t}$ for every $a\in\R$, and $\diff{t}<\diff{t_{2}}<\diff{t_{3}}<\ldots<\diff{t_{k}}$
for every $k>3$, and $\diff{t_{k}}>0$.

\noindent If the orders are equal we have to compare the first standard
parts. E.g. $3\diff{t_{5}}>2\diff{t_{5}}$.

\noindent The other cases fall within the previous ones, because of
the properties of the ordered ring $\ER$. E.g. we have that $\diff{t_{5}}-2\diff{t_{3}}+3\diff{t}<\diff{t_{5}}-2\diff{t_{3}}+\diff{t_{3/2}}$
if and only if $3\diff{t}<\diff{t_{3/2}}$, which is true because
$\omega(\diff{t})=1<\omega(\diff{t_{3/2}})=\frac{3}{2}$. Finally
$\diff{t_{5}}-2\diff{t_{3}}+3\diff{t}>\diff{t_{5}}-2\diff{t_{3}}-\diff{t}$
because $3\diff{t}>-\diff{t}$.

\subsection{Absolute value}

Having a total order we can define the absolute value
\begin{defn}
Let $x\in\ER$, then\[
|x|:=\begin{cases}
x & \text{if }x\ge0\\
-x & \text{if }x<0\end{cases}\]

\end{defn}
\noindent Exactly like for the real field $\R$ we can prove the usual
properties of the absolute value:\begin{align*}
 & |x|\ge0\\
 & |x+y|\le|x|+|y|\\
 & -|x|\le x\le|x|\\
 & ||x|-|y||\le|x-y|\\
 & |x|=0\iff x=0\end{align*}
 Moreover, also the following cancellation law is provable. 
\begin{thm}
\noindent Let $h\in\ER\setminus\{0\}$ and $r$, $s\in\R$, then\[
|h|\cdot r\le|h|\cdot s\then r\le s\]

\end{thm}
\noindent \textbf{Proof:} In fact if $|h|\cdot r\le|h|\cdot s$ then
from Theorem \ref{thm:equivalentFormulationForOrderRelation} we obtain
that either\begin{equation}
\forall^{0}t>0\pti|h_{t}|\cdot r\le|h_{t}|\cdot s\label{eq:absOf_h_tTimes_rLessOrEqualAbs_h_tTimes_s}\end{equation}
 or $|h|\cdot r=|h|\cdot s$. But $h\ne0$ so\[
\left(\forall^{0}t>0\pti h_{t}>0\right)\e{or}\left(\forall^{0}t>0\pti h_{t}<0\right)\]
 hence we can always find a $\bar{t}>0$ such that $|h_{\bar{t}}|\ne0$
and to which \eqref{eq:absOf_h_tTimes_rLessOrEqualAbs_h_tTimes_s}
is applicable. Therefore, in the first case we must have $r\le s$.
In the second one we have\[
|h|\cdot r=|h|\cdot s\]
 but $h\ne0$, hence $|h|\ne0$ and so the conclusion follows from
Theorem \ref{thm:firstCancellationLaw}.\qedNoNewLine

\section{Powers and logarithms}

In this section we will tackle definition and properties of powers
$x^{y}$ and logarithms $\log_{x}y$. Due to the presence of nilpotent
elements in $\ER$, we cannot define these operations without any
limitation. E.g. we cannot define the square root having the usual
properties, like\begin{align}
x\in\ER & \then\sqrt{x}\in\ER\label{eq:SquareRootTakesLittle-ohPolyInLittle-ohPoly}\\
x=y\e{in}\ER & \then\sqrt{x}=\sqrt{y}\e{in}\ER\label{eq:SquareRootIsWellDefined}\\
 & \sqrt{x^{2}}=|x|\nonumber \end{align}
 because they are incompatible with the existence of $h\in D$ such
that $h^{2}=0$, but $h\ne0$. Indeed, the general property stated
in the Subsection \ref{sub:ClosureOfLittle-ohPolyWRTSmoothFunctions}
permits to obtain a property like \eqref{eq:SquareRootTakesLittle-ohPolyInLittle-ohPoly}
(i.e. the closure of $\ER$ with respect to a given operation) only
for smooth functions. Moreover, the Definition \ref{def:extensionOfFunctions}
states that to obtain a well defined operation we need a locally Lipschitz
function. For these reasons, we will limit $x^{y}$ to $x>0$ only,
and $\log_{x}y$ to $x$, $y>0$.
\begin{defn}
Let $x$, $y\in\ER$, with $x>0$, then
\begin{enumerate}
\item $x^{y}:=[t\ge0\mapsto x_{t}^{y_{t}}]_{=\text{ in }\ER}$ 
\item If $y>0$, then $\log_{x}y:=[t\ge0\mapsto log_{x_{t}}y_{t}]_{=\text{ in }\ER}$ 
\end{enumerate}
\end{defn}
\noindent Because of Theorem \ref{thm:equivalentFormulationForOrderRelation}
from $x>0$ we have \[
\forall^{0}t>0\pti x_{t}>0\]
 so that, exactly as we proved in Subsection \ref{sub:ClosureOfLittle-ohPolyWRTSmoothFunctions}
and in Definition \ref{def:extensionOfFunctions}, the previous operations
are well defined in $\ER$.

\noindent From the elementary transfer theorem \ref{thm:elementaryTransferTheorem}
the usual properties follow:\begin{align*}
\left(x^{y}\right)^{z} & =x^{y\cdot z}\\
x^{y}\cdot x^{z} & =x^{y+z}\\
x^{n} & =x\cdot\ptind^{n}\cdot x\e{if}n\in\N\\
\log_{x}\left(x^{y}\right) & =y\\
x^{\log_{x}y} & =y\\
\log(x\cdot y) & =\log x+\log y\\
\log_{x}\left(y^{z}\right) & =z\cdot\log_{x}y\\
x^{\log y} & =y^{\log x}\end{align*}
 About the monotonicity properties, it suffices to use Theorem \ref{thm:equivalentFormulationForOrderRelation}
to prove immediately the usual properties\begin{align*}
z>0\e{and}x\ge y>0 & \then x^{z}\ge y^{z}\\
z<0\e{and}x\ge y>0 & \then x^{z}\le y^{z}\\
z>1\e{and}x\ge y>0 & \then\log_{z}x\ge\log_{z}y\\
0<z<1\e{and}x\ge y>0 & \then\log_{z}x\le\log_{z}y\end{align*}
 Analogous implications, but with strict equalities, are true if we
suppose $x>y$.

Finally, it can be useful to state here the \emph{elementary transfer
theorem for inequalities}, whose proof follows immediately from the
definition of $\le$ and from Theorem \ref{thm:equivalentFormulationForOrderRelation}: 
\begin{thm}
\label{thm:elementaryTransferTheoremForInequalities}Let $A$ be an
open subset of $\R^{n}$, and $\tau$, $\sigma:A\freccia\R$ be smooth
functions. Then \[
\forall x\in{\ext A}\pti\ext\tau(x)\le\ext\sigma(x)\]
 iff \[
\forall r\in A\pti\tau(r)\le\sigma(r).\]

\end{thm}

\section{\label{sec:drawingOfFermatReals}Geometrical representation of Fermat
reals}

At the beginning of this chapter we argued that one of the conducting
idea in the construction of Fermat reals is to maintain always a clear
intuitive meaning. More precisely, we always tried, and we will always
try, to keep a good dialectic between provable formal properties and
their intuitive meaning. In this direction we can see the possibility
to find a geometrical representation of Fermat reals.

\begin{figure}[h]
\begin{centering}
\includegraphics[scale=0.15]{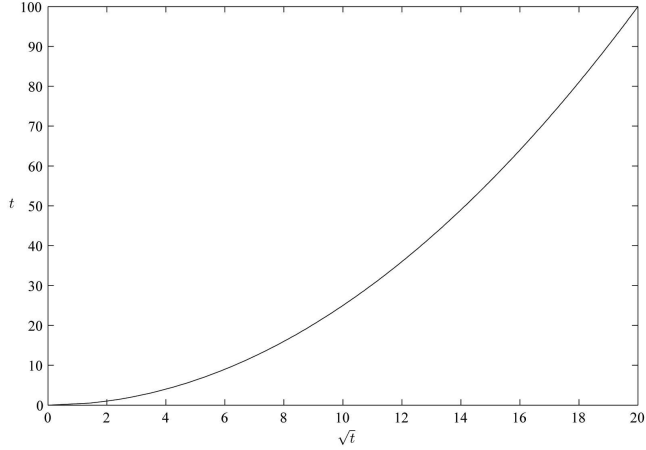} 
\par\end{centering}

\caption{{\small The function representing the Fermat real $\diff{t_{2}}\in D_{3}$}}

\label{fig:representationOf_dt_2} 
\end{figure}

The idea is that to any Fermat real $x\in\ER$ we can associate the
function\begin{equation}
t\in\R_{\ge0}\mapsto\st{x}+\sum_{i=1}^{N}\st{x_{i}}\cdot t^{1/\omega_{i}(x)}\in\R\label{eq:functionsForGeometricalRepresentation}\end{equation}
 where $N$ is, of course, the number of addends in the decomposition
of $x$. Therefore, a geometric representation of this function is
also a geometric representation of the number $x$, because different
Fermat reals have different decompositions, see \ref{thm:existenceUniquenessDecomposition}.
Finally, we can guess that, because the notion of equality in $\ER$
depends only on the germ generated by each little-oh polynomial (see
Definition \ref{def:equalityInFermatReals}), we can represent each
$x\in\ER$ with only the first small part of the function \eqref{eq:functionsForGeometricalRepresentation}.
\begin{defn}
If $x\in\ER$ and $\delta\in\R_{>0}$, then\[
\text{\emph{graph}}_{\delta}(x):=\left\{ (\st{x}+\sum_{i=1}^{N}\st{x_{i}}\cdot t^{1/\omega_{i}(x)},t)\,|\,0\le t<\delta\right\} \]

\noindent where $N$ is the number of addends in the decomposition
of $x$. 
\end{defn}
\noindent Note that the value of the function are placed in the abscissa
position, so that the correct representation of $\text{graph}_{\delta}(x)$
is given by the Figure \ref{fig:representationOf_dt_2}. This inversion
of abscissa and ordinate in the $\text{graph}_{\delta}(x)$ permits
to represent this graph as a line tangent to the classical straight
line $\R$ and hence to have a better graphical picture (see the following
Figures). Finally, note that if $x\in\R$ is a standard real, then
$N=0$ and the $\text{graph}_{\delta}(x)$ is a vertical line passing
through $\st{x}=x$.%

\noindent The following theorem permits to represent geometrically
the Fermat reals
\begin{thm}
\noindent \label{thm:representationTheorem}If $\delta\in\R_{>0}$,
then the function \[
x\in\ER\mapsto\text{\emph{graph}}_{\delta}(x)\subset\R^{2}\]
 is injective. Moreover if $x$, $y\in\ER$, then we can find $\delta\in\R_{>0}$
(depending on $x$ and $y$) such that\[
x<y\]
 if and only if\begin{equation}
\forall p,q,t\pti(p,t)\in\text{\emph{graph}}_{\delta}(x)\e{,}(q,t)\in\text{\emph{graph}}_{\delta}(y)\then p<q\label{eq:orderInTheGeometricalRepresentation}\end{equation}

\end{thm}
\noindent \textbf{Proof:} The application $\rho(x):=\text{graph}_{\delta}(x)$
for $x\in\ER$ is well defined because it depends on the terms $\st{x}$,
$\st{x_{i}}$ and $\omega_{i}(x)$ of the decomposition of $x$ (see
Theorem \ref{thm:existenceUniquenessDecomposition} and Definition
\ref{def:actualDecomposition}). Now, suppose that $\text{graph}_{\delta}(x)=\text{graph}_{\delta}(y)$,
then\begin{equation}
\forall t\in[0,\delta)\pti\st{x}+\sum_{i=1}^{N}\st{x_{i}}\cdot t^{1/\omega_{i}(x)}=\st{y}+\sum_{j=1}^{M}\st{y_{j}}\cdot t^{1/\omega_{j}(y)}\label{eq:representingFunctionsAreLocallyEqual}\end{equation}
 Let us consider the Fermat reals generated by these functions, i.e.\begin{align*}
x': & =\left[t\ge0\mapsto\st{x}+\sum_{i=1}^{N}\st{x_{i}}\cdot t^{1/\omega_{i}(x)}\right]_{=\text{ in }\ER}\\
y': & =\left[t\ge0\mapsto\st{y}+\sum_{j=1}^{M}\st{y_{j}}\cdot t^{1/\omega_{j}(y)}\right]_{=\text{ in }\ER}\end{align*}
 then the decompositions of $x'$ and $y'$ are exactly the decompositions
of $x$ and~$y$\begin{align}
x' & =\st{x}+\sum_{i=1}^{N}\st{x_{i}}\diff{t_{\omega_{i}(x)}}=x\label{eq:xPrimeEqualx}\\
y' & =\st{y}+\sum_{j=1}^{M}\st{y_{j}}\diff{t_{\omega_{j}(y)}}=y\label{eq:yPrimeEqualy}\end{align}
 But from \eqref{eq:representingFunctionsAreLocallyEqual} it follows
$x'=y'$ in $\ER$, and hence also $x=y$ from \eqref{eq:xPrimeEqualx}
and \eqref{eq:yPrimeEqualy}.

\bigskip{}

Now suppose that $x<y$, then, using the same notations of the previous
part of this proof, we have also $x'=x$ and $y'=y$ and hence\[
x'=\st{x}+\sum_{i=1}^{N}\st{x_{i}}\cdot t^{1/\omega_{i}(x)}<\st{y}+\sum_{j=1}^{M}\st{y_{j}}\cdot t^{1/\omega_{j}(y)}=y'\]
 We apply Theorem \ref{thm:equivalentFormulationForOrderRelation}
obtaining that locally $x'_{t}<y'_{t}$, i.e.\[
\exists\delta>0\pti\forall^{0}t\ge0\pti\st{x}+\sum_{i=1}^{N}\st{x_{i}}\cdot t^{1/\omega_{i}(x)}<\st{y}+\sum_{j=1}^{M}\st{y_{j}}\cdot t^{1/\omega_{j}(y)}\]
 This is an equivalent formulation of \eqref{eq:orderInTheGeometricalRepresentation},
and, because of Theorem \ref{thm:equivalentFormulationForOrderRelation}
it is equivalent to $x'=x<y'=y$.\qedNoNewLine 
\begin{example*}
In Figure \ref{fig:firstOrderInfinitesimals} we have the representation
of some first order infinitesimals.%
\begin{figure}[h]
\begin{centering}
\includegraphics{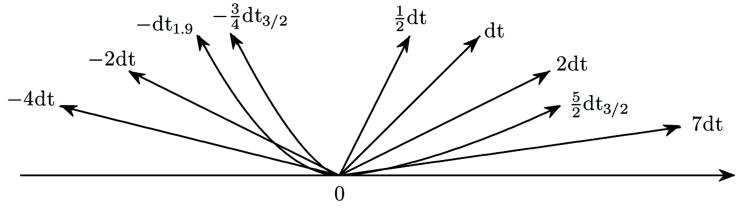} 
\par\end{centering}

\caption{{\small Some first order infinitesimals}}

\label{fig:firstOrderInfinitesimals}
\end{figure}

\end{example*}
\noindent The arrows are justified by the fact that the representing
function \eqref{eq:functionsForGeometricalRepresentation} is defined
on $\R_{\ge0}$ and hence has a clear first point and a direction.
The smaller is $\alpha\in(0,1)$ and the nearer is the representation
of the product $\alpha\diff{t}$, to the vertical line passing through
zero, which is the representation of the standard real $x=0$. Finally,
recall that $\diff{t_{k}}\in D$ if and only if $1\le k<2$.

\begin{figure}[h]
\begin{centering}
\includegraphics[scale=0.7]{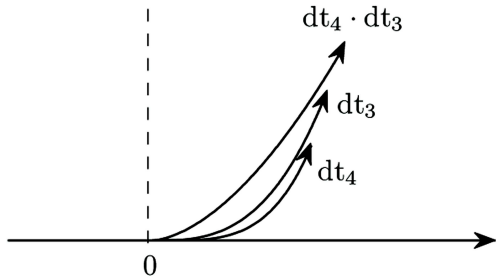} 
\par\end{centering}

\caption{{\small The product of two infinitesimals}}

\label{fig:productOfTwoInfinitesimals}
\end{figure}

\noindent If we multiply two infinitesimals we obtain a smaller number,
hence one whose representation is nearer to the vertical line passing
through zero, as represented in the Figure \ref{fig:productOfTwoInfinitesimals}%

\begin{figure}[h]
\begin{centering}
\includegraphics{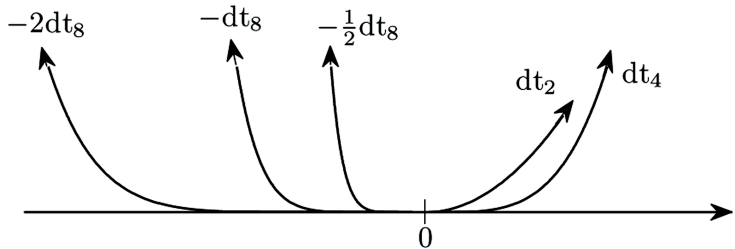} 
\par\end{centering}

\caption{Some higher order infinitesimals}

\label{fig:higherOrderInfinitesimals}
\end{figure}

\noindent In Figure \ref{fig:higherOrderInfinitesimals} we have a
representation of some infinitesimals of order greater than $1$.
We can see that the greater is the infinitesimal $h\in D_{a}$ (with
respect to the order relation $\le$ defined in $\ER$) and the higher
is the order of intersection of the corresponding line $\text{graph}_{\delta}(h)$.%

\begin{figure}[h]
\begin{centering}
\includegraphics[scale=1.2]{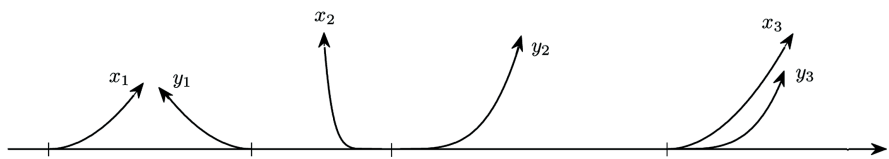} 
\par\end{centering}

\caption{Different cases in which $x_{i}<y_{i}$}

\label{fig:representationOfOrderRelation}
\end{figure}

\noindent Finally, in Figure \ref{fig:representationOfOrderRelation}
we represent the order relation on the basis of Theorem \ref{thm:representationTheorem}.
Intuitively, the method to see if $x<y$ is to look at a suitably
small neighborhood (i.e. at a suitably small $\delta>0$) at $t=0$
of their representing lines $\text{graph}_{\delta}(x)$ and $\text{graph}_{\delta}(y)$:
if, with respect to the horizontal directed straight line, the curve
$\text{graph}_{\delta}(x)$ comes before the curve $\text{graph}_{\delta}(y)$,
then $x$ is less than $y$.%

\part{Infinite dimensional spaces}

\chapter[Approaches to $\infty$-dimensional diff. geom.]{\label{cha:ApproachesToDiffGeomOfInfDimSpaces}Approaches to differential
geometry of infinite dimensional spaces}

\section[Approaches to infinite-dimensional spaces]{Introduction\label{sec:ApproachesToDiffGeomOfInfDimSpaces}}

In this section we want to list some of the most important, i.e. well-established,
approaches that are used to define geometrical structures in infinite
dimensional spaces. One of the most important example we have in mind
is the set $\ManInfty(M,N)$ of all the smooth applications between
two finite dimensional manifolds $M$ and $N$. For the aims of the
present section, we are interested to list some of the most studied
structures on $\ManInfty(M,N)$, and its subspaces, that permit to
develop at least a tangency theory, i.e. the notion of tangent functor
and the notion of differentiability of maps between this type of infinite
dimensional spaces, and have sufficiently good categorical properties.
This is not a trivial goal because, for example, an important example
we can cite is the group $\text{Diff}(M)$ of all the diffeomorphisms
of a manifold $M$. Flows in a compact manifold $M$ can be considered
as $1$-parameter subgroups of $\text{Diff}(M)$, and it would seem
useful to express the smoothness of a flow by means of a suitable
differentiable structure on $\text{Diff}(M)$, which should also behave
like a classical Lie group with respect to this structure.

A typical restriction to distinguish among different approaches to
infinite dimensional spaces is the hypotheses of compactness of the
domain $M$, assumed to obtain some desired property: is this a necessary
hypotheses or are we forced to assume it due to some restrictions
of the chosen approach?

Another interesting property is the possibility to extend the classical
notion of manifold to a more general type of space, so as to get better
categorical properties, like the existence of infinite products or
co-products or a cartesian closed category%
\footnote{For a short introduction, mainly motivated to fix common notations,
of the few notions of category theory used in the present work, see
Appendix \ref{app:someNotionsOfCategoryTheory}.%
}.

Finally, several authors had to tackle the following problem: suppose
we have a new notion of smooth space able to include the space $\ManInfty(M,N)$,
at least for $M$ compact and finite dimensional, and to embed faithfully
(i.e. injectively, see Appendix \ref{app:someNotionsOfCategoryTheory})
the category of smooth finite dimensional manifolds. Even if the extension
of the notion of finite dimensional manifold is faithful, usually
the category $\boldsymbol{\mathcal{C}}$ of these new smooth spaces
includes spaces which are too much general, so that it seems really
hard to generalize for these spaces meaningful results of differential
geometry of finite dimensional manifolds. For this reason, several
authors (see e.g. \citet{Kr-Mi,Fr-Kr,Lav,Mo-Re}) try to select, among
all their new smooth spaces in $\boldsymbol{\mathcal{C}}$, the best
ones having some new more restrictive properties. In this way the
category $\boldsymbol{\mathcal{C}}$ acts as a universe, usually closed
with respect to strong categorical operations (like arbitrary limits,
colimits and cartesian closedness), and the restricted class of smooth
spaces works as a true generalization of the notion of manifold.

For example, in \citet{Kr-Mi} the category of Frölicher spaces acts
as a universe, but indeed the monograph is about manifolds modeled
in convenient vector spaces instead of classical Banach spaces (see
subsection \emph{The convenient vectors spaces settings} \vpageref{sec:TheConvenientVectorSpacesSettings}).
This permits to \citet{Kr-Mi} to generalize as far as possible to
infinite dimensional manifolds the results of finite dimensional spaces,
but as a consequence the class of manifolds modeled in convenient
vector spaces loses some desired categorical properties.

Analogously, in SDG (see e.g. \citet{Lav,Mo-Re,Koc}) the class of
restricted smooth spaces is introduced with the notion of microlinear
space and the universe is a suitable topos, i.e. a whole model for
intuitionistic set theory. In this approach, the infinitesimals are
used to define the properties of this class of restricted, better
behaved, spaces.

Of course, this is not possible in theories that have not an explicit
language of actual infinitesimals, like in the case of diffeological
spaces (see \citet{Igl}). For them we can proceed either as in convenient
vector spaces theory considering the notion of vector space in the
category of smooth diffeological spaces (i.e. smooth diffeological
spaces that are also vector spaces with smooth operations, see \citet{Igl})
and considering manifolds modeled in diffeological vector spaces,
or we can try to develop directly for a generic diffeological space
some notion of differential geometry (see e.g. \citet{Igl,Lau2,Lau1,He-Ma,Hec1,Sou2,Sou1}).
In the following subsections we will return to this problem giving
some more precise definitions.

To understand better some differences between the approaches we are
going to describe shortly in this section, we want to motivate the
notion of cartesian closure, because is one of the basic choice shared
by several authors like \citet{Bas,Bel,Br1,Br2,Br3,Che,Col,Fr-Bu,Fr-Kr,Koc,Kr-Mi,Lav,La1,La2,Mo-Re,Sei,Sou1,Ste,Vog}.
We firstly fix the notations for the notions of \emph{adjoint} of
a map.
\begin{defn}
If $X$, $Y$, $Z$ are sets and $f:X\freccia Z^{Y}$, $g:X\times Y\freccia Z$
are maps, then\begin{align*}
\forall(x,y)\in X\times Y\pti f^{\vee}(x,y) & :=\left[f(x)\right](y)\in Z\\
\forall x\in X\pti g^{\wedge}(x) & :=g(x,-)\in Z^{Y}\end{align*}
 hence\begin{align*}
f^{\vee}: & X\times Y\freccia Z\\
g^{\wedge}: & X\freccia Z^{Y}\end{align*}
 The map $f^{\vee}$ is called the \emph{adjoint} of $f$ and the
map $g^{\wedge}$ is called the \emph{adjoint}%
\footnote{Here we are using the notations of \citet{Ad-He-St}, but some authors,
e.g. \citet{Kr-Mi}, used opposite notations for the adjoint maps.%
} of $g$. 
\end{defn}
\noindent Let us note that $\left(f^{\vee}\right)^{\wedge}=f$ and
$\left(g^{\wedge}\right)^{\vee}=g$, that is the two applications\begin{align*}
(-)^{\vee}: & \left(Z^{Y}\right)^{X}\freccia Z^{X\times Y}\\
(-)^{\wedge}: & Z^{X\times Y}\freccia\left(Z^{Y}\right)^{X}\end{align*}
 are one the inverse of the other and hence represent in explicit
form the bijection of sets $\left(Z^{Y}\right)^{X}\simeq Z^{X\times Y}$
i.e. $\Set(X,\Set(Y,Z))\simeq\Set(X\times Y,Z)$.

One of the main aim of the second part of the present work is to generalize
the notions of smooth manifold and of smooth map between two manifolds
so as to obtain a new category {}``with good properties'' that will
be denoted by $\CInfty$; if we call \emph{smooth maps} the morphisms
of $\CInfty$ and \emph{smooth spaces} its objects, then this category
must be \emph{cartesian closed}, i.e. it has to verify the following
properties for every pair of smooth space $X$, $Y\in\CInfty$: 
\begin{enumerate}
\item $\CInfty(X,Y)$ is a smooth space, i.e. $\CInfty(X,Y)\in\CInfty$\label{enu:expObjects} 
\item \label{enu:defCartesianClosureNaturalIso}The maps $(-)^{\vee}$ and
$(-)^{\wedge}$ are smooth, i.e. they realize in the category $\CInfty$
the bijection $\CInfty(X,\CInfty(Y,Z))\simeq\CInfty(X\times Y,Z)$ 
\end{enumerate}
\noindent Property \ref{enu:expObjects}. is another way to state
that the category we want to construct must contain as objects the
space of all the smooth maps between two generic objects $X$, $Y\in\CInfty$\begin{align*}
\CInfty(X,Y) & =\{f\,|\, X\xfreccia{f}Y\text{ \ is smooth}\}=\\
 & =\{f\,|\, X\xfreccia{f}Y\text{ \ is a morphism of }\CInfty\}.\end{align*}
 Moreover, let us note that as a consequence of \ref{enu:defCartesianClosureNaturalIso}.
we have that\begin{align}
X\xfreccia{f}\CInfty(Y,Z)\ \ \text{is smooth}\ \iff & \ \ X\times Y\xfreccia{f^{\vee}}Z\ \ \text{is smooth}\label{eq:doubleImplicationVee}\\
X\times Y\xfreccia{g}Z\ \ \text{is smooth}\ \iff & \ \ X\xfreccia{g^{\wedge}}\CInfty(Y,Z)\ \ \text{is smooth}.\label{eq:doubleImplicationWedge}\end{align}

\noindent The importance of \eqref{eq:doubleImplicationVee} and \eqref{eq:doubleImplicationWedge}
can be explained saying that if we want to study a smooth map having
values in the space $\CInfty(Y,Z)$, then it suffices to study its
adjoint map $f^{\vee}$. If, e.g., the spaces $X$, $Y$ and $Z$
are finite dimensional manifolds, then $\CInfty(Y,Z)$ is infinite-dimensional,
but $f^{\vee}:X\times Y\freccia Z$ is a standard smooth map between
finite dimensional manifolds, and hence we have a strong simplification.
Conversely, if $g:X\times Y\freccia Z$ is a smooth map, then it generates
a smooth map with values in $\CInfty(Y,Z)$, and all the smooth maps
with values in this type of spaces can be generated in this way. Of
course, this idea is frequently used, even if informally, in the calculus
of variations. Let us note explicitly that the cartesian closure of
the category $\CInfty$, i.e. properties \ref{enu:expObjects}. and
\ref{enu:defCartesianClosureNaturalIso}., does not say anything about
smooth maps with a \emph{domain} of the form $\CInfty(Y,Z)$, but
it reformulates in a convenient way the problem of smoothness of maps
with \emph{codomain} of this type. For a more abstract notion of cartesian
closed category, see e.g. \citet{Mac,Bor,Ar-Ma,Ad-He-St}.

We also want to see a different motivation drawn from \citet{Fr-Kr}.
Let us suppose to have a smooth function $g:\R\times I\freccia\R$,
where $I=[a,b]$, and define the integral function\[
f(t):=\int_{a}^{b}g(t,s)\diff{s}\quad\forall t\in\R.\]
 Then we can look at the function $f$ as the composition of two applications\[
f:t\in\R\mapsto g(t,-)\mapsto\int_{a}^{b}g(t,-)\in\R.\]
 Hence, if we denote\[
i:h\in\mathcal{C}^{\infty}(I,\R)\mapsto\int_{a}^{b}h\in\R,\]
 then\[
f=i\circ g^{\wedge}\ee{i.e.}f(t)=i\left(g^{\wedge}(t,-)\right)\quad\forall t\in\R.\]
 In this way, it is natural to try a proof of the formula for the
derivation under the integral sign in the following way:\begin{multline}
\frac{\diff{f}}{\diff{t}}(t)=\frac{\diff}{\diff{t}}\left(i\circ g^{\wedge}\right)=\diff{i\left(g^{\wedge}(t)\right)}\left[\frac{\diff{g^{\wedge}}}{\diff{t}}(t)\right]=\\
=i\left[\partial_{1}g(t,-)\right]=\int_{a}^{b}\partial_{1}g(t,s)\diff{s}.\label{eq:IdeaForTheDerivationUnderTheIntegralSign}\end{multline}
 Here we have supposed that the following properties hold: 
\begin{itemize}
\item $g^{\wedge}:\R\freccia\mathcal{C}^{\infty}(I,\R)$ is smooth, 
\item $i:\mathcal{C}^{\infty}(I,\R)\freccia\R$ is smooth, 
\item the chain rule for the derivative of the composition of two functions, 
\item the differential of the function $i$ is given by $\diff{i}(h)=i$
for every $h\in\mathcal{C}^{\infty}(I,\R)$, because $i$ is linear, 
\item $\frac{\diff{g^{\wedge}}}{\diff{t}}(t)=\partial_{1}g(t,-)$. 
\end{itemize}
Let us note explicitly that the space $\Cc^{\infty}(I,\R)$ is infinite
dimensional.

Even if in the present work we will be able to prove all these properties,
the aim of \eqref{eq:IdeaForTheDerivationUnderTheIntegralSign} is
not to suggest a new proof, but to hint that a theory where we can
consider the previous properties seems to be very flexible and powerful.

\section{{\large \label{sec:BanachManifoldsAndLocallyConvexVectorSpaces}Banach
manifolds and locally convex vector spaces}}

Banach manifolds is the more natural generalization of finite-dimensional
manifolds if one takes Banach spaces as local model spaces. Even if,
as we will see more precisely in this section, this theory does not
satisfy our condition to present in this chapter only generalized
notions of manifolds able to develop at least a tangency theory and
having sufficiently good categorical properties, Banach manifolds
are the most studied concept in infinite dimensional differential
geometry. Some well known references on Banach manifolds are \citet{Lan,Ab-Ma-Ra}.
Among the most important theorems in this framework we can cite the
implicit and inverse function theorems and the existence and uniqueness
of solutions of Lipschitz ordinary differential equations on such
spaces. The use of charts to prove these fundamental results is indispensable,
so it is not easy to generalize them to more general contexts where
we cannot use the notion of chart having values in some modeling space
with sufficiently good properties.

For the purposes of the present analysis, a typical example of infinite-dimensional
Banach space is the space $\mathcal{C}^{r}(M,E)$ of $\mathcal{C}^{r}$-maps,
where $M$ is a compact manifold and $E$ is a Banach space. The vector
space $\mathcal{C}^{r}(M,E)$ is a Banach space with respect to the
norm\begin{equation}
\left\Vert f\right\Vert _{r}:=\max_{1\le i\le r}\sup_{m\in M}\left\Vert \diff{^{i}f}(m)\right\Vert ,\label{eq:normInBanach}\end{equation}
 but the theory fails for the space $\mathcal{C}^{\infty}(M,E):=\bigcap_{r=1}^{+\infty}\mathcal{C}^{r}(M,E)$
of smooth mappings defined in $M$ and with values in $E$. On the
one hand, even if it is not a formal motivation, but it remains very
important in the real development of mathematics, the hypotheses of
considering $r<+\infty$ and $M$ compact in the previous definition
\ref{eq:normInBanach} are not intrinsic to the problem but are motivated
solely by the limitations of the instrument we are trying to implement,
i.e. a norm in the space $\mathcal{C}^{r}(M,E)$. On the other hand,
more formally, any two different norms $\Vert-\Vert_{r}$ and $\Vert-\Vert_{s}$
are not equivalent, and hence the space $\mathcal{C}^{\infty}(M;,E)$
is not normable with a norm generating the same topology generated
by the family of norms $\left(\Vert-\Vert_{r}\right)_{r=1}^{+\infty}$
(for details, see e.g. \citet{Fri}; in the following, saying that
the space $\mathcal{C}^{\infty}(M,E)$ is not normable, we will always
mean with respect to this topology).

Moreover, $\mathcal{C}^{\infty}(M,E)$ is not a Banach manifold: indeed,
it is separable and metric (see \citet{Fri}), hence if it were a
Banach manifold, then it would be embeddable as an open subset of
an Hilbert space (see \citet{Hend}), and hence it would be normable.

Therefore, the category of Banach manifolds and smooth maps $\mathbf{Ban}$
is not cartesian closed because it is not closed with respect to exponential
objects $\mathbf{Ban}(M,E)=\mathcal{C}^{\infty}(M,E)$, see condition
\ref{enu:expObjects}. in the previous definition of cartesian closed
category, section \ref{sec:ApproachesToDiffGeomOfInfDimSpaces}.

This also proves that the category of Banach manifolds $\mathbf{Ban}$
and smooth maps does not have arbitrary limits: in fact if it had
infinite products (a particular case of limit in a category, see Appendix
\ref{app:someNotionsOfCategoryTheory}), then we would have\[
\prod_{m\in M}E=\mathbf{Ban}(M,E)=\mathcal{C}^{\infty}(M,E),\]
 but we had already seen that this space is not a Banach manifold.

These important counter-examples can conduct us toward the idea of
considering spaces equipped with a family of norms, like $\left(\Vert-\Vert_{r}\right)_{r=1}^{+\infty}$,
or, more generally, of seminorms, i.e. toward the theory of locally
convex vector spaces (see e.g. \citet{Jar}). But any locally convex
topology on the space $\mathcal{C}^{\infty}(M,E)$ is incompatible
with cartesian closure, as stated in the following 
\begin{thm}
\label{thm:cartesianClosednessAndLCVS}Let $F$ be a locally convex
vector space contained in a subcategory $\mathbf{\mathcal{T}}$ of
the category $\mathbf{Top}$ of topological spaces and continuous
functions such that $\mathcal{T}(F,\R)$ always contains all the linear
continuous functionals on the space $F$\[
\text{\emph{Lin}}(F,\R)\subseteq\mathcal{T}(F,\R).\]
 Then we have the following implication\[
\mathbf{\mathcal{T}}\text{ is cartesian closed}\then F\text{ is normable}.\]
 Hence the category $\mathbf{Ban}$ is not cartesian closed because
the space \[
F=\mathcal{C}^{\infty}(\R,\R)\]
is not normable. 
\end{thm}
\noindent \textbf{Proof:} We can argue as in \citet{Kr-Mi}: because
$\mathbf{\mathcal{T}}$ is cartesian closed, every evaluation\[
\text{ev}_{XY}(x,f):=f(x)\quad\forall x\in X\:\forall f\in\mathbf{\mathcal{T}}(X,Y)\]
 is an arrow of $\mathbf{\mathcal{T}}$ (this is a general result
in every cartesian closed category, see e.g. \citet{Mac}) and hence
it is also a continuous function, because $\mathbf{\mathcal{T}}$
is a subcategory of $\mathbf{Top}$ by hypotheses. In this case, we
also have that the restriction of $\text{ev}_{F\R}$ to the subspace
$F^{*}:=\text{Lin}(F,\R)\subseteq\mathbf{\mathcal{T}}(F,\R)$ of linear
continuous functionals on the space $F$ would also be (jointly) continuous:\[
\varepsilon:=\text{ev}_{F\R}|_{F\times F^{*}}:F\times F^{*}\freccia\R.\]
 Then we can find neighborhoods $U\subseteq F$ and $V\subseteq F^{*}$
of zero such that $\varepsilon(U\times V)\subseteq[-1,1]$, that is\[
U\subseteq\left\{ u\in F\,|\,\forall f\in V:\,\,|f(u)|\le1\right\} .\]
 But then, taking a generic functional we can always find $\lambda\in\R_{\ne0}$
such that $\lambda g\in V$, and hence $|g(u)|\le1/\lambda$ for every
$u\in U$. Any continuous functional is thus bounded on $U$, so the
neighborhood $U$ itself is bounded (see e.g. \citet{Jar,Kr-Mi}).
But any locally convex vector space with a bounded neighborhood of
zero is normable (see e.g. \citet{Jar,Do-Sm}). \qedNoNewLine

\noindent This theorem also asserts that notions like Fréchet manifolds
(manifolds modeled in locally convex metrizable and complete vector
spaces) are incompatible with cartesian closedness too.

For a more detailed study about cartesian closedness and Banach manifolds,
see \citet{Br1,Br2,Br3}; for a more detailed study about the relationships
between the topology on spaces of continuous linear functionals $\text{Lin}(F,E)$
and normable spaces, see \citet{Kel,Mai}.

Because one of our aim is to obtain a category $\CInfty$ of {}``smooth''
(and hence topological) spaces embedding the category $\mathbf{Ban}$,
a direct consequence of Theorem \ref{thm:cartesianClosednessAndLCVS}
is that, in general, we will not have a locally convex topology on
spaces of functions like $\CInfty(M,\R)$. Nevertheless, in $\CInfty$
we always have that every arrow (i.e. every smooth function in a generalized
sense) is also continuous and every evaluation is smooth.

Finally, another important problem in the theory of Banach manifolds
is tied with infinite dimensional Lie groups. As it is well known,
they appear in several connections in physics, like in the study of
both compressible and incompressible fluids, in magnetohydrodynamics,
in plasma-dynamics or in electrodynamics (see e.g. \citet{Ab-Ma-Ra}
and references therein). The fundamental results of \citet{Om1} (see
also \citet{Om-dlH,Om2}) show that a Banach Lie group $G$ acting
smoothly, transitively and effectively on a compact manifold $M$
must necessary be finite dimensional. This result strongly underlines
that the space of all the diffeomorphisms $G=\text{Diff}(M)$ of a
compact manifold in itself cannot be a Banach Lie group.

It is important to note that the present work is not in contrast with
the theory of Banach manifolds, but rather it tries to complement
it overpassing some of its defects, like the absence of a calculus
of actual infinitesimals and the lacking of spaces of mappings. On
the one hand, a first aim of the present work is to obtain a category
$\CInfty$ of smooth spaces with better categorical properties (e.g.
we will see that the category $\CInfty$ is cartesian closed and possesses
arbitrary limits and colimits, e.g. infinite products, infinite disjoint
sums or quotient spaces). On the other hand, of course we aim at exploiting
the language of nilpotent infinitesimals. We will see that the category
$\mathbf{Ban}$ of smooth Banach manifolds is faithfully embedded
in our category $\CInfty$ of smooth spaces.

\section{\label{sec:TheConvenientVectorSpacesSettings}The convenient vector
spaces settings}

It is very interesting to note that the original idea to define the
differential of functions $f:\R^{n}\freccia\R^{m}$ reducing it to
the composition $f\circ c$ with differentiable curves $c:\R\freccia\R^{n}$
goes back (for didactic reasons!) to \citet{Had}: in this work a
function $f:\R^{2}\freccia\R$ was called differentiable if all the
compositions $f\circ c$ with differentiable curves $c:\R\freccia\R^{2}$
are again differentiable and satisfy the chain rule. Later (see \citet{Mic})
this notion has been extended to mapping $f:E\freccia F$ between
generic topological vector spaces: $f$ is defined to be differentiable
at $x\in E$ if there exists a continuous linear mapping $l:E\freccia F$
such that $f\circ c:\R\freccia F$ is differentiable at $0$ with
derivative $(l\circ c')(0)$ for each everywhere differentiable curve
$c:\R\freccia E$ with $c(0)=x$. This notion of differentiable function
is really more restrictive that the usual one, but it is equivalent
to the standard notion of smooth function if in it we replace the
word {}``differentiable'' with {}``smooth''. More generally if
we replace {}``differentiable'' with {}``of class $\Cc^{k}$ and
with locally Lipschitz $k$-th derivative'', we obtain an equivalence
with the classical notion. These results have been proved by \citet{Bom}
and all the theory of convenient vector spaces depends strongly on
these non trivial results.

Several theories which detach from the theory of Banach manifolds,
like the convenient vector spaces setting or the following diffeological
spaces, are grounded on generalization of this idea (not necessarily
knowing the cited article \citet{Had}). In particular, the theory
of convenient vector spaces is probably the most developed theory
of infinite dimensional manifolds ables to overpass several problems
of Banach manifolds. Presently, the most complete reference is \citet{Kr-Mi},
even if the theory started with \citet{Fr-Bu} and \citet{Fr-Kr}.

Only to mention few results, in the convenient vector spaces setting
the hard implicit function theorem of Nash and Moser (see \citet{Hami,Kr-Mi})
can be proved, very good results can also be obtained for both holomorphic
and real analytic calculus, the theorem of De Rham can be proved and
the theory of infinite dimensional Lie groups can be well developed.

Because in the present work we will show that any manifold modeled
in convenient vector spaces can be embedded in our category $\CInfty$,
we present very briefly one of the possible equivalent definitions
of this type of spaces and some few notions about smooth manifolds
modeled in convenient vector spaces. 
\begin{defn}
We say that $E$ is a \emph{convenient vector space} iff $E$ is a
locally convex vector space where every smooth curve has a primitive,
i.e.\[
\forall c\in\mathcal{C}^{\infty}(\R,E)\,\exists p\in\mathcal{C}^{\infty}(\R,E):\ p'=c\]

\end{defn}
\noindent Considering the Cauchy-Bochner integral, any Banach space
is hence a convenient vector space, but several non trivial example
directly comes from the cartesian closedness of the category of all
the convenient vector spaces (see \citet{Kr-Mi}).

As mentioned above what type of topology can be considered in a convenient
vector space, due to the cartesian closedness of the related category,
is a non trivial point. The idea to reduce, as far as possible, any
possible notion to the corresponding notion for smooth curves, can
carry us toward the natural idea to consider the final topology for
which any smooth curve is also continuous, i.e. the following
\begin{defn}
Let $E$ be a convenient vector space, then we say that\[
U\text{ is }c^{\infty}\text{-open in }E\]

iff\[
\forall c\in\Cc^{\infty}(\R,E):\ c^{-1}(U)\text{ is open in }\R\]

\end{defn}
\noindent The category of convenient vector spaces is cartesian closed
so that, e.g. $\Cc^{\infty}(\R,\R)$ is again a convenient vector
space. We can now define as usual the notion of chart modeled in a
$c^{\infty}$-open set of a convenient vector space and hence the
corresponding notion of smooth manifold and of smooth map between
two manifolds. So as to avoid confusion with our category $\CInfty$,
in the following we will denote with $\Cc_{\text{cvs}}^{\infty}$
the category of smooth manifolds modeled in convenient vector spaces.
Using suitable generalizations of Boman's theorem (\citet{Bom}),
it is hence possible to prove the following (see \citet{Kr-Mi}) 
\begin{thm}
\noindent Let $M$, $N$ be manifolds modeled on convenient vector
spaces, then we have that $f:M\freccia N$ is smooth iff

\[
\forall c\in\Cc_{\text{\emph{cvs}}}^{\infty}(\R,M)\pti f\circ c\in\Cc_{\text{\emph{cvs}}}^{\infty}(\R,N).\]

\end{thm}
\noindent Using the notion of $c^{\infty}$-open subset of a convenient
vector space and the notion of chart is possible to define a topology
on every manifold considering the final topology in which every chart
is continuous. We have hence the expected result that $W$ is open
in this topology on $M$ if and only if $c^{-1}(W)$ is open in $\R$
for every smooth curve $c\in\Cc_{\text{cvs}}^{\infty}(\R,M)$ (see
\citet{Kr-Mi}).

The notion of Frölicher space provides the possibility to construct
a category with very good properties acting as a universe for the
class of manifolds modeled in convenient vector spaces. We cite here
the definition of Frölicher space only to underline the analogies
with our smooth spaces in $\CInfty$: 
\begin{defn}
A Frölicher space is a triple $(X,\Cc_{X},\mathcal{F}_{X})$ consisting
of a set $X$, a subset $\mathcal{C}_{X}\subseteq X^{\R}$ of curves
on this set, and a subset $\mathcal{F}_{X}\subseteq\R^{X}$ of real
valued functions defined on $X$, with the following properties: 
\begin{enumerate}
\item $\forall f\pti f\in\mathcal{F}_{X}\iff\left[\forall c\in\Cc_{X}:\,\, f\circ c\in\Cc^{\infty}(\R,\R)\right]$ 
\item $\forall c\pti c\in\mathcal{C}_{X}\iff\left[\forall f\in\mathcal{F}_{X}:\,\, f\circ c\in\Cc^{\infty}(\R,\R)\right]$ 
\end{enumerate}
\end{defn}
\noindent The category of Frölicher spaces is cartesian closed and
possesses arbitrary limits and colimits. A locally convex vector space
$E$ is a convenient vector space if and only if it is a Frölicher
space with respect to curves and functions defined as $\Cc_{X}:=\Cc_{\text{cvs}}^{\infty}(\R,E)$
and $\mathcal{F}_{X}:=\Cc_{\text{cvs}}^{\infty}(E,\R)$. Finally,
because of cartesian closedness, it is possible to define a unique
structure of Frölicher space on the set $Y:=\Cc^{\infty}(M,N)$ of
all the smooth maps between two manifolds given by\[
\Cc_{Y}:=\left\{ c:\R\freccia N^{M}\,|\, c^{\vee}:\R\times N\freccia M\text{ is smooth}\right\} \]
 and\[
\mathcal{F}_{Y}:=\left\{ f:N^{M}\freccia\R\,|\,\forall c\in\Cc_{Y}:\,\, f\circ c\in\Cc^{\infty}(\R;,\R)\right\} .\]
 In the following we will use again the symbol $\Cc^{\infty}(M,N)$
to indicate this structure of Frölicher space.

As mentioned at the beginning of this chapter, the notion of manifold
modeled in convenient vector spaces permits to include several infinite
dimensional spaces non ascribable into Banach manifold theory, but,
at the same time, forces us to lose some good categorical property.
In particular the space of all smooth mappings $\Cc^{\infty}(M,N)$
between two manifolds has a manifold structure only for $M$ and $N$
finite dimensional (see \citet{Kr-Mi}, Chapter IX). Moreover, if
$\mathfrak{C}^{\infty}(M,N)$ is this manifold structure%
\footnote{Note that, e.g. if $M=N=\R$, this structure is different from the
structure of convenient vector space (and Frölicher space) $\Cc^{\infty}(\R,\R)$;
for this reason the authors of \citet{Kr-Mi} use a different symbol
$\mathfrak{C}^{\infty}(\R,\R)$.%
} on the set $\Cc^{\infty}(M,N)$, then the exponential law\[
\Cc^{\infty}(M,\mathfrak{C}^{\infty}(N,P))\simeq\Cc^{\infty}(M\times N,P)\]

\noindent holds if and only if $N$ is compact (see \citet{Kr-Mi},
Theorem 42.14).

Using an intuitive interpretation introduced by \citet{La1} we can
say that in the convenient vector spaces settings the fundamental
figure of our spaces is the curve and every notion is reduced to a
corresponding notion about curves. We will use several times later
this intuitive, and fruitfully, interpretations also for other types
of figures. In the notion of Frölicher space there is a particular
stress in the symmetry between curves and functions, but this symmetry
has not been adopted by other authors, like in the following approach
about diffeological spaces.

We will see that both Frölicher spaces and manifolds modeled in convenient
vector spaces are embedded in our category $\CInfty$ of smooth spaces,
so that our approach can supply a language of actual infinitesimals
also to these settings.

\section{\label{sec:DiffeologicalSpaces}Diffeological spaces}

Using the language of the {}``fundamental figures'' given on a general
space $X$ introduced by \citet{La1}, we can describe diffeological
spaces as a natural generalization of the previously seen idea to
take as fundamental figures all the smooth curves $c:\R\freccia X$
on the space $X$. To define the concept of diffeological space, we
first denote with\[
\text{Op}:=\left\{ U\,|\,\exists n\in\N\,:\, U\text{ is open in }\R^{n}\right\} \]

\noindent the set of all the domains of our new figures in the space
$X$. In informal words, the idea of a diffeological space is to say
that the structure on the space $X$ is specified if we give all the
smooth figures $p:U\freccia X$, for $U\in\text{Op}$. More formally,
we have 
\begin{defn}
\noindent We say that $(\mathcal{D},X)$ is a \emph{diffeological
space} iff $X$ is a set and $\mathcal{D}=\left\{ \mathcal{D}_{U}\right\} _{U\in\text{\emph{Op}}}$
is a family of sets of functions\[
\mathcal{D}_{U}\subseteq\Set(U,X)\quad\forall U\in\text{\emph{Op}}\]
 The functions $p\in\mathcal{D}_{U}$ are called \emph{parametrizations}
or \emph{plots} or \emph{figures} \emph{on} $X$ \emph{of type $U$.}
The family $\mathcal{D}$ has to satisfies the following conditions: 
\begin{enumerate}
\item Every point of $X$ is a figure, i.e. for every $U\in\text{\emph{Op}}$
and every constant map $p:U\freccia X$, we must have that $p\in\mathcal{D}_{U}$. 
\item Every set of figures $\mathcal{D}_{U}$ is closed with respect to
re-parametrization, i.e. if $p:U\freccia X$ is a figure in $\mathcal{D}_{U}$,
and $f\in\Cc^{\infty}(V,U)$, where $V\in\text{\emph{Op}}$, then
$p\circ f\in\mathcal{D}_{V}$. 
\item The family $\mathcal{D}=\left\{ \mathcal{D}_{U}\right\} _{U\in\text{\emph{Op}}}$
verifies a sheaf property, i.e. let $V\in\text{\emph{Op}}$, $(U_{i})_{i\in I}$
be an open cover of $V$ and $p:V\freccia X$ a map such that $p|_{U_{i}}\in\mathcal{D}_{U_{i}}$,
then $p\in\mathcal{D}_{V}$. In other words, to be locally a figure
implies to be a figure globally too. 
\end{enumerate}
\noindent Finally a map $f:X\freccia Y$ between two diffeological
spaces $(X,\mathcal{D}^{X})$ and $(Y,\mathcal{D}^{Y})$ is said to
be \emph{smooth} if it takes figures of the domain space in figures
of the codomain space, i.e. if\[
\forall U\in\text{\emph{Op}}\,\,\forall p\in\mathcal{D}_{U}^{X}\pti f\circ p\in\mathcal{D}_{U}^{Y}\]

\end{defn}
If compared with Frölicher spaces, in Diffeology (i.e. the study of
diffeological spaces, see \citet{Igl}) the principal differences
are in the generalization of the types of figures, in the losing of
the symmetry between figures and corresponding functions (i.e. maps
of type $f:X\freccia U$ for $U\in\text{Op}$) and in the fundamental
sheaf property. For example, the generalization to figures of arbitrary
dimension instead of curves only, permits to prove the cartesian closure
of the category of diffeological spaces very easily and without the
use of the non trivial Boman's theorem (see \citet{Fr-Kr,Kr-Mi,Bom}).
The original idea to consider figures of general dimension instead
of curves only, and the fundamental sheaf condition date back to \citet{Che2,Che};
the definition of diffeological space, essentially in the form given
above, is originally of \citet{Sou1,Sou2}.

\noindent The category of diffeological spaces has very good categorical
properties, with arbitrary limits (subspaces, products, pullbacks,
etc.) and colimits (quotient spaces, sums, pushforwards, etc.) and
cartesian closedness (so that set theoretical compositions and evaluations
are always smooth). Classical Fréchet manifolds are fully and faithfully
embedded in this category (see \citet{Los2}).

We can now define a \emph{diffeological vector space} (over $\R$)
any diffeological space $(E,\mathcal{D})$, where $E$ is a vector
space (over $\R$), and such that the addiction and the multiplication
by a scalar\[
(u,v)\in E\times E\mapsto u+v\in E\e{and}(r,u)\in\R\times E\mapsto ru\in E\]

\noindent are smooth (with respect to the suitable product diffeologies
on the domains) and, as usual, the notion of \emph{smooth manifolds
modeled on diffeological vector spaces}.

Anyway, differential geometry on generic diffeological spaces can
be developed surprisingly far as showed e.g. by \citet{Igl}: homotopy
theory, exterior differential calculus, differential forms, Lie derivatives,
integration on chains and Stokes formula, de Rham cohomology, Cartan
formula, generalization of symplectic geometry to diffeological spaces,
etc. As said in \citet{Igl}:
\begin{quotation}
Thanks to the strong stability of diffeology under the most important
categorical operations {[}...{]} every general construction relating
to this theory applies to spaces of functions, differential forms,
fiber bundles, homotopy, etc. without leaving the strict framework
of diffeology. This makes the development of differential geometry
much more easier, much more natural, than usually.
\end{quotation}
\noindent It is also interesting to note that some of these generalizations
(like Stokes formula) are general consequences of this type of extension
of the notion of manifolds, as proved by \citet{Los}, and hence are
not peculiar of Diffeology.

From the point of view of the present work, Diffeology is surely formally
clear, but sometimes lacks from the point of view of the intuitive
geometrical interpretation. To illustrate this assertion, we can consider
the notion of tangent vector as formulated in \citet{Igl}. In the
following we will assume that $(X,\mathcal{D})$ is a diffeological
space and $x\in X$ is a point in the space $X$. The first idea is
that the figures $q:U\freccia X$ of type $U\subseteq\R^{n}$ of the
space $X$ permit to define the notion of smooth $p$-form without
having the notion of tangent vector, but abstracting the properties
of the pullback $q^{*}$ of the figure $q\in\mathcal{D}_{U}$. In
other words, let us suppose that we have already defined what is a
differential $p$-form on $X$, then we would be able to define the
pullback $q^{*}$ of $q$ as a map that associates to each point $u\in U\subseteq\R^{n}$
a $p$-form in $\Lambda^{p}(\R^{n})$. The idea is hence to define
directly a $p$-form as this action on figures through pullback, and
asking the natural condition of composition of pullbacks in case we
take a parametrization $f\in\Cc^{\infty}(V,U)$ of the domain of the
figure $q$: 
\begin{defn}
A \emph{differential $p$-form defined on $X$} is a family of maps
$(\alpha_{U})_{U\in\text{\emph{Op}}}$. Each $\alpha_{U}$, for $U$
open in $\R^{n}$, associates to each figure $q\in\mathcal{D}_{U}$
a smooth $p$-form $\alpha_{U}(q):U\freccia\Lambda^{p}(\R^{n})$ ,
i.e.\[
\alpha_{U}:\mathcal{D}_{U}\freccia\Cc^{\infty}(U,\Lambda^{p}(\R^{n})),\]

\noindent and it has to satisfies the condition\[
\alpha_{V}(q\circ f)=f^{*}(\alpha_{U}(q))\]

\noindent for every plot $q\in\mathcal{D}_{U}$ and for every smooth
parametrization $f\in\Cc^{\infty}(V,U)$ defined on the open set $V\in\text{\emph{Op}}$.
The set of all the differential $p$-forms defined on $X$ will be
denoted by $\Omega^{p}(X)$. 
\end{defn}
\noindent The method used to arrive at this definition is the (frequently
used in mathematics) {}``inversion of the effect with the cause''
in case of bijection between effects and causes. Indeed, if $X=$
is an open set of $\R^{d}$, then it is possible to prove that we
have a natural isomorphism between the new definition and the classical
notion of smooth $p$-form, i.e. $\Omega^{p}(U)\simeq\Cc^{\infty}(U,\Lambda^{p}(U))$,
in other words pullbacks of $p$-forms uniquely determine the $p$-forms
themselves.

\noindent The previous definition satisfy all the properties one needs
from it, like the possibility to define a diffeology on $\Omega^{p}(X)$,
vector space structure, pullbacks, exterior differential, exterior
product, a natural notion of germ generated by a $p$-form so that
two forms are equal if and only if they generate the same germ (that
if they are {}``locally'' equal), etc.

The first intuitive drawback of the definition of $\Omega^{p}(X)$
is that there is no mention to spaces $\Lambda_{x}^{p}(X)$ of $p$-forms
associated to each point $x\in X$ and of the relationships between
these spaces and the whole $\Omega^{p}(X)$. Therefore, to understand
better the following definitions, we introduce the following 
\begin{defn}
We say that two forms $\alpha$, $\beta\in\Omega^{p}(X)$ \emph{have
the same value at $x$}, and we write $\alpha\sim_{x}\beta$, if and
only if for every figure $q\in\mathcal{D}_{U}$ such that\[
0\in U\e{and}q(0)=x\]

\noindent (in this case we will say that $q$ \emph{is centered at
$x$)} we have that\emph{\[
\alpha(q)(0)=\beta(q)(0).\]
} Equivalence classes of $p$-forms by means of the equivalence relation
$\sim_{x}$ are called \emph{values of $\alpha$ at $x$} and we will
denote with \emph{$\Lambda_{x}^{p}(X):=\Omega^{p}(X)/\sim_{x}$ this
quotient set.} 
\end{defn}
\noindent Using these values of 1-forms we can define tangent vectors.
Firstly we introduce the paths on $X$ and the values of a 1-form
on each path with the following 
\begin{defn}
\noindent Let us introduce the space of all the paths on $X$, i.e.\[
\text{\emph{Paths}}(X):=\Cc^{\infty}(\R,X)\]

\noindent and for each path $q\in\text{\emph{Paths}}(X)$, the map
$j(q):\Omega^{1}(X)\freccia\R$ evaluating each 1-form at zero\[
j(q):\alpha\in\Omega^{1}(X)\mapsto\alpha(q)(0)\in\R.\]

\noindent The map $j(q)$ is linear and smooth (because it is an evaluation),
hence\[
j:\text{\emph{Paths}}(X)\freccia L^{\infty}(\Omega^{1}(X),\R),\]

\noindent where $L^{\infty}(\Omega^{1}(X),\R)$ is the space of all
the linear smooth functionals defined on the space of 1-forms of $X$. 
\end{defn}
\noindent Secondly we say that the set of all these values $j(q)$
generates the whole tangent space. The set of these generators is
introduced in the following
\begin{defn}
\noindent The space $C_{x}^{\wedge}(X)$ is the image of all the paths
passing through $x$ under the map $j$: \[
C_{x}^{\wedge}(X):=\left\{ j(q)\,|\, q\in\text{\emph{Paths}}(X)\text{ and }q(0)=x\right\} \subseteq L^{\infty}(\Omega^{1}(X),\R)\]

\end{defn}
\noindent In the space $C_{x}^{\wedge}(X)$ is naturally defined a
multiplication by a scalar $r\in\R$ that formalizes the idea to increase
the speed of going through a given path $q\in\text{Paths}(X)$:\[
r\cdot j(q)=j\left[q(r\cdot(-))\right],\]

\noindent where $q(r\cdot(-))$ is the path $q(r\cdot(-)):s\in\R\freccia q(r\cdot s)\in X$.
But the space $C_{x}^{\wedge}(X)$ is not necessarily a vector space
because is not closed with respect to addiction of these values $j(q)$
of 1-forms on paths $q$ centered at $x$, hence we finally define 
\begin{defn}
\noindent A tangent vector $v\in T_{x}(X)$ is a linear combination
of elements of $C_{x}^{\wedge}(X)$, i.e. \[
v=\sum_{i=1}^{n}s_{i}v_{i}\]
 for some \begin{align*}
n & \in\N\\
(v_{i})_{i=1}^{n} & \text{ sequence of }C_{x}^{\wedge}(X)\\
(s_{i})_{i=1}^{n} & \text{ sequence of }\R.\end{align*}

\end{defn}
\noindent As we said, even if the definitions we have just introduced
are formally correct, their intuitive geometric meaning remains obscure.
In classical manifolds theory, the definition of tangent vector through
1-forms is not geometrically intrinsic unless of Riemannian manifolds,
so it is not clear why passing to a more general space we are able
to obtain this identification in an intrinsic way. Secondly, diffeological
spaces include also spaces with singular points, like $X=\left\{ (x,y)\in\R^{2}\,|\, x\cdot y=0\right\} $.
At the origin $x=(0,0)\in X$ there is no way to define in a geometrically
meaningful way the sum of the two tangent vectors corresponding to
$\boldsymbol{i}=(1,0)$ and $\boldsymbol{j}=(0,1)$ (without using
the superspace $\R^{2}$). This is the principal motivation that conducts
SDG to introduce the notion of \emph{microlinear space} as the spaces
where to each pair of tangent vectors it is possible to associate
an infinitesimal parallelogram, fully contained in the space itself,
whose diagonal represents the sum of these two tangent vectors. The
previous space $X$ is not microlinear exactly at the origin.

As we will see, our category $\CInfty$ is exactly the category of
diffeological spaces and concretely we will only generalize the definition
of diffeological space so as to obtain a more flexible instruments
that will permit us to define the category $\ECInfty$ of spaces extended
with the new infinitesimal points. E.g. we will have that $\R\in\CInfty$
and $\ER\in\ECInfty$. Hence, the theory of Fermat reals naturally
includes diffeological spaces and also provides to them a language
of actual infinitesimals. The use of these infinitesimals opens the
possibility to simplify and clarify some concepts already developed
in the framework of diffeological spaces, e.g. gaining a more clear
geometrical meaning. We will also see that using these infinitesimal
we will also arrive to new results, like the existence of infinitesimals
flows corresponding to a given smooth vector field.

\section{\label{sec:SDG}Synthetic differential geometry}

The fundamental ideas upon which SGD%
\footnote{Frequently SDG is also called smooth infinitesimal analysis.%
} born, originate from the work of \citet{Ehr}, \citet{Wei} and A.
Grothendieck (see \citet{SGA4}). \citet{Ehr} introduced the concept
of $k$-jet at a point $p$ in a manifold $M$ as an important geometric
structure determined by the $k$-th order Taylor's formula of real
valued functions $f$ defined in a neighborhood of $p\in M$. As said
by \citet{Mac2}: 
\begin{quotation}
{[}...{]} the study of jets can be seen as a development of the earlier
idea of studying the {}``infinitely nearby'' points on algebraic
curves on manifolds. Presumably it was Ehresmann's initiative which
stimulated the paper by \citet{Wei}. 
\end{quotation}
\noindent In this work A. Weil introduced the idea to formalize nilpotent
infinitesimals using algebraic methods, more precisely using quotient
rings like $\R[x]/(x^{2})$ or $\R[x,y]/(x^{2},y^{2})$, in general
formal power series in $n$ variables $\R[[x_{1},\dots,x_{n}]]$ modulo
the $(k+1)$-th power of a given ideal $I=(i_{1},\dots,i_{m})$ of
series $i_{1},\dots,i_{m}\in\R[[x_{1},\ldots,x_{n}]]$ with zeros
constant term, i.e. such that $i_{j}(\underbar{0})=0$ for every $j=1,\dots,m$.
These type of objects are now called \emph{Weil algebras}, and C.
Ehresmann's jets are also special cases of Weil algebras. Very roughly,
we can guess the fundamental idea of A. Weil saying that, e.g., an
element $p\in\R[x]/(x^{2})$ can be written as $p=a+x\cdot b$, with
$a$, $b\in\R$, with addiction given in the more obvious way and
multiplication given by $(a+x\cdot b)\cdot(\alpha+x\cdot\beta)=a\alpha+x\cdot(a\beta+b\alpha)$,
that is the same result we would obtain if we multiply the two polynomials
$a+x\cdot b$ and $\alpha+x\cdot\beta$ with the formal rules $x^{2}=0$.
At the end, with a construction as simple as the definition of the
field of complex numbers, we have extended the real field into a ring
with a non-zero element $x$ having zero square, i.e. a first order
infinitesimal (but in this ring there are not infinitesimals of greater
order). Using the same idea, we can see that with the Weil algebra
$\R[x,y]/(x^{2},y^{2})$ we have extended the real field with two
first order infinitesimals $x$, $y$ whose product is not zero%
\footnote{We recall Section \ref{sec:Infinitesimals-and-order-properties} to
underline an important difference with our approach.%
} $x\cdot y\ne0$. Suitably generalized to algebras of germs of smooth
functions defined on manifolds, these two examples, i.e. $\R[x]/(x^{2})$
and $\R[x,y]/(x^{2},y^{2})$, correspond isomorphically to the first
and second tangent bundle respectively (see e.g. \citet{Wei,Kr-Mi,Kr-Mi2,Koc,Lav,Mo-Re,Ber}
for more details). The next fundamental step to obtain a single framework
where all these types of nilpotent infinitesimals are available, has
been performed by A. Grothendieck who tried to use nilpotent infinitesimals
in his theory of schemes to treat infinitesimal structures in algebraic
geometry. The basic idea was to study an algebraic locus like $S^{1}=\left\{ (x,y)\in\R^{2}\,|\, x^{2}+y^{2}=1\right\} $,
not only as a subset of points in the plane, but as the functor $S_{\text{F}}^{1}:\textbf{CRing}\freccia\Set$
from the category $\textbf{CRing}$ of commutative rings with 1 to
the category of sets defined as\begin{align*}
S_{\text{F}}^{1}(A) & :=\left\{ (a,b)\in A^{2}\,|\, a^{2}+b^{2}=0\right\} \\
S_{\text{F}}^{1}(A\xfreccia{f}B) & :=(f\times f)|_{S_{\text{F}(A)}^{1}}:S_{\text{F}}^{1}(A)\freccia S_{\text{F}}^{1}(B)\end{align*}

\noindent (where $f:A\freccia B$ is a ring homomorphism and$ $ $f\times f:(a,b)\in A^{2}\mapsto(f(a),f(b))\in B^{2}$).
Using this approach algebraic geometers started to understand that
the functor corresponding to the trivial locus $\left\{ x\in\R\,|\, x=x\right\} =\R$,
i.e. the functor $R(A):=\left\{ a\in A\,|\, a=a\right\} =A=$ the
underlying set of the ring $A$, behaves like a set of scalars containing
infinitesimals. E.g. $D(A):=\left\{ a\in A\,|\, a^{2}=0\right\} $
is a subfunctor of this functor $R$ and plays the role of the space
of first order infinitesimals. Being a subfunctor, $D$ {}``behaves''
like a subset%
\footnote{In the sense that each Topos is a model of intuitionistic set theory,
so that it is possible to define a formal language for intuitionistic
set theory where sentences like $D\subseteq R$ are rigorous and true
in the model (see \citet{Mo-Re,Koc} for more details).%
} of $R$. These ideas conducted to the notion of Grothendieck topos.
Lawvere found that in the Grothendieck topos, and in other similar
categories that later will originate the general notion of topos (see
\citet{Gra}), an intuitionistic set-theoretic language can be directly
interpreted in any topos. In \citet{La1} he proposes a way to generalize
these construction of algebraic geometry to smooth manifolds theory,
and to use this generalization as a foundation for infinitesimal reasoning
valid both for finite and infinite dimensional manifolds. This proposal
was part of a big project whose objective is to establish an intrinsic
axiomatizaton for continuum mechanics. The inclusion of infinite dimensional
spaces like functions spaces is a natural consequence of the cartesian
closedness of every topos.

The construction of a model for SDG which embeds the category of smooth
finite dimensional manifolds is not a simple task. Classical references
are \citet{Mo-Re,Koc}. Here we only want to sketch some of the fundamental
ideas, first of all to underline the conceptual differences between
SDG and the above mentioned approaches to infinite dimensional differential
geometry.

\noindent The first idea to generalize from the context of algebraic
geometry to manifolds theory is to find a corresponding of the category
of $\textbf{CRing}$ of commutative rings, i.e. to pass from a context
of polynomial operations to more general smooth functions. Indeed,
that category is replaced by that of $\Cc^{\infty}$-rings:
\begin{defn}
\noindent A $\Cc^{\infty}$-ring $(A,+,\cdot,\iota)$ is a ring $(A,+,\cdot)$
together with an interpretation $\iota(f)$ of each possible smooth
map $f:\R^{n}\freccia\R^{m}$, that is a map\[
i(f):A^{n}\freccia A^{m}\]

\noindent such that $\iota$ preserves projections, compositions and
identity maps, i.e.: 
\begin{enumerate}
\item If $p:\R^{m}\freccia\R$ is a projection, then $\iota(p):A^{m}\freccia A$
is a projection. 
\item If $\R^{d}\xfreccia{g}\R^{n}\xfreccia{f}\R^{m}$ are smooth, then
$\iota(f\circ g)=\iota(f)\circ\iota(g)$. 
\item If $1_{\R^{n}}:\R^{n}\freccia\R^{n}$ is the identity map, then $\iota(1_{\R^{n}})=1_{\iota(A^{n})}$. 
\end{enumerate}
A \emph{homomorphism of} $\Cc^{\infty}$\emph{-rings} is a ring homomorphism
which preserves the interpretation of smooth maps, that is such that
\[
\xyR{40pt}\xyC{40pt}\xymatrix{ & A^{n}\ar[d]_{A(f)}\ar[r]^{\varphi^{n}} & B^{n}\ar[d]^{B(f)}\\
 & A^{m}\ar[r]^{\varphi^{m}} & B^{m}}
\]

\end{defn}
\noindent We may define a $\Cc^{\infty}$-ring in an equivalent but
more concise way: let $C^{\infty}$ denote the category whose objects
are the spaces $\R^{d}$, $d\ge0$, and with smooth functions as arrows,
then a $\Cc^{\infty}$-ring is a finite product preserving functor
$A:\Cc^{\infty}\freccia\Set$, and a $\Cc^{\infty}$-homomorphism
is just a natural transformation $\varphi:A\freccia B$. Indeed, given
such a functor, the set $A(\R)$ has the structure of a commutative
ring $(A(\R),+_{A},\cdot_{A})$ given by $+_{A}:=A(\R\times\R\xfreccia{+}\R)$
and $\cdot_{A}:=A(\R\times\R\xfreccia{\cdot}\R)$, where $+:\R\times\R\freccia\R$
and $\cdot:\R\times\R\freccia\R$ are the ring operations on $\R$.

\noindent Here are some examples of $\Cc^{\infty}$-rings 
\begin{example}
\noindent The ring $\Cc^{\infty}(\R^{d},\R)$ of real valued smooth
functions $a:\R^{d}\freccia\R$, with pointwise ring operations, is
a $\Cc^{\infty}$-ring. Usually it is denoted simply with $\Cc^{\infty}(\R^{d})$.
The smooth function $f:\R^{n}\freccia\R^{m}$ is interpreted in the
following way. Let $(h_{1},\dots,h_{n})\in\Cc^{\infty}(\R^{d},\R)^{n}$,
be $n$ elements of the ring $\Cc^{\infty}(\R^{d})$. Their product
\[
(h_{1},\dots,h_{n}):x\in\R^{d}\mapsto(h_{1}(x),\dots,h_{n}(x))\in\R^{n}\]

\noindent can be can be composed with $f:\R^{n}\freccia\R^{m}$ and
projected into its $m$ components obtaining\[
\iota(f):=(p_{1}\circ f\circ(h_{1},\dots,h_{n}),\dots,p_{m}\circ f\circ(h_{1},\dots,h_{n}))\in\Cc^{\infty}(\R^{d},\R),\]

\end{example}
\noindent where $p_{i}:\R^{m}\freccia\R$ are the projections. \medskip{}

\begin{example}
If $M$ is a smooth manifold, the ring of real valued functions defined
on $M,$ i.e. $\Cc^{\infty}(M,\R)$, is a $\Cc^{\infty}$-ring. Here
a smooth function $f:\R^{n}\freccia\R^{m}$ is interpreted using composition,
similarly to the previous example. This ring is also denoted by $\Cc^{\infty}(M)$.
Moreover, it is well known that\[
\Cc^{\infty}(M)=\Cc^{\infty}(N)\then M=N.\]

\noindent If $g:N\freccia M$ is a smooth map between manifolds, then
the $\Cc^{\infty}$-ho\-mo\-mor\-phi\-sm given by\[
\Cc^{\infty}(g):a\in\Cc^{\infty}(M,\R)\mapsto a\circ g\in\Cc^{\infty}(N,\R)\]

\noindent verifies the analogous embedding property:\[
\Cc^{\infty}(g)=\Cc^{\infty}(h)\then g=h.\]

\end{example}
\noindent This means that manifolds can be faithfully considered as
$\Cc^{\infty}$-rings. \medskip{}

\begin{example}
Let $A$ be a $\Cc^{\infty}$-ring and $I$ an ideal of $A$, then
the quotient ring $A/I$ is also a $\Cc^{\infty}$-ring. Indeed, if
$A(f):A^{n}\freccia A^{m}$ is the interpretation of $f:\R^{n}\freccia\R^{m}$,
we can define the interpretation $(A/I)(f):(A/I)^{n}\freccia(A/I)^{m}$
as\begin{multline*}
(A/I)(f)([a_{1}]_{I},\dots,[a_{n}]_{I}):=\\
=([p_{1}(A(f)(a_{1},\dots a_{n}))]_{I},[p_{m}(A(f)(a_{1},\dots a_{n}))]_{I}),\end{multline*}

\end{example}
\noindent where $[a_{i}]_{I}\in A/I$ denotes the equivalent classes
of the quotient ring, and $p_{j}:A^{m}\freccia A$ are the projections
(see e.g. \citet{Mo-Re} for more details). Examples included in this
case are the analogous of the above mentioned $D_{k}:=\Cc^{\infty}(\R)/(x^{k+1})$
and $D(2):=\Cc^{\infty}(\R)/(x^{2},y^{2})$, or the ring $\triangle:=\Cc_{0}^{\infty}(\R^{n})=\Cc^{\infty}(\R^{n})/m_{\left\{ 0\right\} }^{g}$,
where $m_{\left\{ 0\right\} }^{g}$ is the ideal of smooth functions
having zero germs at $0\in\R^{n}$ and finally $\mathbb{I}:=\Cc_{0}^{\infty}(\R^{n}\setminus\left\{ 0\right\} )$.
These $\Cc^{\infty}$-rings will play the role, in the final model,
respectively of \emph{infinitesimals of $k$-th order $D_{k}$,} of
\emph{pairs of infinitesimals of first order whose product is not
necessarily zero} $D(2)$, of the set of all the infinitesimals $\triangle$
and of the set of all the \emph{invertible infinitesimals $\mathbb{I}$.}
\medskip{}

For each subset $X\subseteq\R^{n}$, a function $f:X\freccia\R$ is
said to be \emph{smooth} if there is an open superset $U\supseteq X$
and a smooth function $g:U\freccia\R$ which extends $f$, i.e. $g|_{X}=f$.
We can proceed as in the previous example using composition to define
the $\Cc^{\infty}$-ring $\Cc^{\infty}(X)$ of real valued functions
defined on $X$. An important example that uses this generalization
and the previous example is $\Cc^{\infty}(\N)/K$, where $\Cc^{\infty}(\N)$
is the ring of smooth functions on the natural numbers, and $K$ is
the ideal of eventually vanishing functions. This ring will act, in
the final model, as the set of \emph{infinitely large natural numbers}.
\medskip{}

\begin{example}
A $\Cc^{\infty}$-ring $A$ is called \emph{finitely generated} if
it is isomorphic to one of the form $\Cc^{\infty}(\R^{n})/I$, for
some $n\in\N$ and some finitely generated ideal $I=(i_{1},\dots,i_{m})$.
For example, given an open subset $U\subseteq\R^{n}$ we can find
a smooth function $f:\R^{n}\freccia\R$ such that $f(x)\ne0$ if and
only if $x\ne U$. So $U$ is diffeomorphic to the closed set $\hat{U}=\left\{ (x,y)\,|\, y\cdot f(x)=1\right\} \subseteq\R^{n+1}$.
Hence we have the isomorphism of $\Cc^{\infty}$-rings\[
\Cc^{\infty}(U)\simeq\Cc^{\infty}(\R^{n+1})/(y\cdot f(x)-1).\]

\noindent This proves that $\Cc^{\infty}(U)$ is finitely generated.
Using this result and Whitney's embedding theorem it is possible to
prove that for a manifold $M$, the $\Cc^{\infty}$-ring $\Cc^{\infty}(M)$
is finitely generated too (see \citet{Mo-Re,Koc}). 
\end{example}
\noindent Therefore, the category $\mathbb{L}$ of finitely generated
$\Cc^{\infty}$-rings seems a good step toward the goal to embed finite
dimensional manifolds in a category with infinitesimal objects. However,
function spaces can in general not be constructed in $\mathbb{L}$.
In order to have these function spaces, the first step is to extend
the category $\mathbb{L}$ in the category $\Set^{\mathbb{L}^{\text{op}}}$
of presheaves on $\mathbb{L}$, i.e. of functors $F:\mathbb{L}\freccia\Set$:\[
\textbf{Man}\subseteq\mathbb{L}\subseteq\Set^{\mathbb{L}^{\text{op}}}.\]

\noindent This is a natural step in this context because the embedding
$\mathbb{L}\subseteq\Set^{\mathbb{L}^{\text{op}}}$ is a well know
result in category theory (see Yoneda embedding in Appendix \ref{app:someNotionsOfCategoryTheory}),
and because the category $\Set^{\mathbb{L}^{\text{op}}}$ is a topos.
So we concretely see the possibility to embed the category of smooth
manifolds in a topos containing infinitesimal objects too. Let us
note that manifolds are directly embedded in $\Set^{\mathbb{G^{\text{op}}}}$
without {}``an extension with new infinitesimal points'', so the
approach is very different with respect, e.g., to NSA or to the present
work.

So, what is the ring of scalars representing the geometric line in
the topos $\Set^{\mathbb{L}^{\text{op}}}$? If $A$, $B\in\mathbb{L}$
are finitely generated $\Cc^{\infty}$-rings, and $f:A\freccia B$
is a $\Cc^{\infty}$-homomorphism, this geometric line is represented
by the functor\begin{align}
R(A) & =\mathbb{L}(A,\Cc^{\infty}(\R))\label{eq:lineObjectFunctor}\\
R(A\xfreccia{f}B) & :g\in R(A)\mapsto g\circ f\in R(B)\label{eq:lineObjectFunctor_II}\end{align}

\noindent corresponding, via the Yoneda embedding, to the $\Cc^{\infty}$-ring
$\Cc^{\infty}(\R)$. The set of first order infinitesimal $D$ corresponds
in the topos $\Set^{\mathbb{L}^{\text{op}}}$ to the functor\begin{align}
D(A) & =\mathbb{L}(A,\Cc^{\infty}(\R)/(x^{2}))\label{eq:firstOrderInfinitesimalsFunctor}\\
D(A\xfreccia{f}B) & :g\in D(A)\mapsto g\circ f\in D(B).\label{eq:firstOrderInfinitesimalsFunctor_II}\end{align}

Indeed, the topos $\Set^{\mathbb{L}^{\text{op}}}$ is not the final
model of SDG for several reasons. Among these, we can cite that in
the topos $\Set^{\mathbb{L}^{\text{op}}}$ are not provable properties
like $1\ne0$ or $\forall r\in\R(x\text{ is invertible }\vee\,(1-x)\text{ is invertible})$,
and this is essentially due because the embedding $\textbf{Man}\subseteq\Set^{\mathbb{L}^{\text{op}}}$
does not preserve open covers. A description of the final models is
outside the scopes of the present work. For more details see e.g.
\citet{Mo-Re} and references therein. In the light of the examples
\eqref{eq:lineObjectFunctor}, \eqref{eq:lineObjectFunctor_II} and
\eqref{eq:firstOrderInfinitesimalsFunctor}, \eqref{eq:firstOrderInfinitesimalsFunctor_II}
we can quote \citet{Mo-Re}: 
\begin{quotation}
In recent years, several alternative solutions to the problem of generalizing
manifolds to include function spaces and spaces with singularities
have been proposed in the literature. A particularly appealing one
is the theory of convenient vector spaces {[}...{]}. These structures
are in a way simpler than the sheaves considered in this book, but
one should notice that the theory of convenient vector spaces does
not include an attempt to develop an appropriate framework for infinitesimal
structures, which is one of the main motivations of our approach 
\end{quotation}
\noindent The present work tries to go exactly in the direction to
have a simple generalization of manifolds (indeed, simpler than convenient
vector spaces and as simple as diffeological spaces) and at the same
time infinitesimals structures.

Hence, it is in the opinion of the researchers in SDG that these topos
models are not sufficiently simple, even if, at the same time, they
are very rich and formally powerful. For these reasons smooth infinitesimal
analysis is usually presented in an {}``axiomatic'' way, in the
framework of a naive intuitionistic set theory%
\footnote{Exactly as almost every mathematician works in naive (classical) set
theory. On the other hand to work in SDG, one has to learn to work
in intuitionistic logic, i.e. avoiding the law of the excluded middle,
the proofs by reduction ad absurdum ending with a double negation,
the full De Morgan laws, the equivalence between double negation and
affirmation, the full equivalence between universal and existential
quantifiers through negation, the axiom of choice, etc.%
}, but with explicit introduction of particular axioms useful to deal
with smooth spaces (i.e. objects of $\Set^{\mathbb{L}^{\text{op}}}$
or a better model) and smooth functions (i.e. arrows of $\Set^{\mathbb{L}^{\text{op}}}$
or a better model). This possibility is due to the above mentioned
internal language for a set theory that can be defined in every topos
(that represents its intuitionistic semantics). For example a basic
assumption is the so-called Kock-Lawvere axiom:
\begin{assumption}
$R$ is a ring and we define $D:=\left\{ h\in R\,|\, h^{2}=0\right\} $,
called the set of first order infinitesimal. They satisfy:\begin{equation}
\forall f:D\freccia R\,\,\exists!m\in R:\,\,\forall h\in D\pti f(h)=f(0)+h\cdot m.\label{eq:Kock-Lawvere}\end{equation}

\end{assumption}
\noindent The universal quantifier {}``for every function $f:D\freccia R$''
really means {}``for every set theoretical function from $D$ to
$R$'', but definable using intuitionistic logic. In semantical terms,
this corresponds to {}``for every arrow in the model $\Set^{\mathbb{L}^{\text{op}}}$'',
i.e. for every smooth natural transformation between the functor $D$
(see \eqref{eq:firstOrderInfinitesimalsFunctor} and \eqref{eq:firstOrderInfinitesimalsFunctor_II})
and the functor $R$ (see \eqref{eq:lineObjectFunctor} and \eqref{eq:lineObjectFunctor_II}).
It is not surprising to assert that \eqref{eq:Kock-Lawvere} is incompatible
with classical logic: putting\begin{equation}
f(h)=\begin{cases}
1\text{ if} & h\ne0\\
0\text{ if} & h=0\end{cases}\end{equation}

\noindent then applying the Kock-Lawvere axiom \eqref{eq:Kock-Lawvere}
with this function and considering the hypothesis $\exists h_{0}\in D:\, h_{0}\ne0$,
we obtain\[
1=0+h_{0}\cdot m.\]

\noindent Squaring this equality we obtain $1=0$. Considering this
incompatibility with classical logic a motivation to consider intuitionistic
logic, is a natural passage only in a context of topos theory and
only if one already is thinking to the existence of models like $\Set^{\mathbb{L}^{\text{op}}}$.
But in another context we think that the more natural idea is to criticize
\eqref{eq:Kock-Lawvere} asking some kind of limitation on the functions
to which it can be really applied. Indeed, this was one of the first
motivation to start the present work. Indeed, we will take strong
inspiration from SDG in this work, but we can affirm that these two
theories are very different. Our attention to stress the intuitive
meaning of the new infinitesimals numbers does not find a correspondence
in SDG, where infinitesimal of very different types can be defined,
but sometimes loosing the corresponding intuitive meaning. About this
point of view we can quote \citet{Con}: 
\begin{quotation}
\label{quo:ConwayOnSDG}I think I should point out that {[}SDG{]}
isn't really trying to be a candidate for setting up infinitesimal
analysis. It's just a formal algebraic technique for {}``working
up to any given order in some small variable $ $$s$'' - for instance
if you want to work up to second order in $s$, you just declare that
$s^{3}=0$. 
\end{quotation}
\noindent Even if we do not completely agree with this strong affirmation,
it represents an authoritative opinion that underlines the differences
between SDG and our approach.

Finally we cite that the work of \citet{Wei} has been the base for
several other research tempting to formalize in some way nilpotent
infinitesimal methods (but without getting all the difficulties of
SDG). In this direction we can cite Weil functors (see \citet{Kr-Mi,Ko-Mi-Sl,Kr-Mi2})
and the recent \citet{Ber}.

\chapter{\label{cha:theCartesianClosure}The cartesian closure of a category
of figures}

\section{\label{sec:hypothesisForCartesianClosure}Motivations and basic hypotheses}

In this section we shall define the basic constructions which will
lead us to the category of $\Cn$ \emph{spaces} and $\Cn$ \emph{functions};
we will realize these constructions for a generic $n\in\N_{>0}\cup\left\{ +\infty\right\} $,
even if in the next chapters concerning calculus and differential
geometry we will consider the case $n=+\infty$ only. Any $\Cc^{n}$
manifold is a $\Cn$ space too, and the category $\Cn$ of all $\Cn$
spaces is cartesian closed (see Section \ref{sec:ApproachesToDiffGeomOfInfDimSpaces}),
hence it contains several infinite-dimensional spaces, the first of
which we are interested in is $\Cc^{n}(M,N)$, i.e. the space of all
the usual $\Cc^{n}$ functions between two manifolds $M$ and $N$.
It is important to note that, exactly as in \citet{Kr-Mi} and in
\citet{Mo-Re}, the category $\Cn$ contains many {}``pathological''
spaces; actually $\Cn$ works as a {}``cartesian closed universe''
and we will see that, like in \citet{Koc,Lav,Mo-Re}, the particular
\emph{inf-linear $\Cn$ spaces} have the best properties, and will
work as a good substitute of manifolds (we have already made some
comments about this way of proceeding in Section \ref{sec:ApproachesToDiffGeomOfInfDimSpaces}).

The ideas used in this section arise from analogous ideas about diffeological
spaces and Fr\"olicher spaces (see Section \ref{sec:TheConvenientVectorSpacesSettings}),
in particular our first references are \citet{Che} and \citet{Fr-Kr};
actually $\CInfty$ is the category of diffeological spaces (see Section
\ref{sec:DiffeologicalSpaces}). For these reasons, in this section
we will not present the proofs of the most elementary facts; these
can be indeed easily generalized from analogous proofs of \citet{Che,Fr-Kr,Kr-Mi}
or \citet{Igl}. The results presented in this and the following chapter
have been already published in \citet{Gio3}.

\noindent We present the definition of cartesian closure starting
from a concrete category $\F$ of topological spaces (satisfying few
conditions) and embedding it in a cartesian closed category $\bar{\F}$.
We will call $\bar{\F}$ the \emph{cartesian closure of} $\F$. We
need this generality because we shall use it to define both domain
and codomain of the extension functor $\ext(-)\;:\;\CInfty\freccia\ECInfty$,
that generalizes the construction $\R\mapsto\ER$ associating to each
smooth space $M\in\CInfty$ its extension with our infinitesimal points
$\ext{M}\in\ECInfty$. Indeed, the categories acting as the domain
and the codomain of this functor will be defined starting from two
different categories $\F$ and applying the cartesian closure. 

The problem to generalize the definition of $\ER$ to a functor $\ext(-)$
can also be seen from the following point of view: at this stage of
the present work, it is natural to define a tangent vector to a manifold
$M$ as a map\[
t:D\freccia\ext{M}.\]
 But we have to note that the map $t$ has to be {}``regular'' in
some sense, hence we need some kind of geometric structure both on
the domain of first order infinitesimals $D$ and on the codomain
$\ext{M}$. On the other hand, it is natural to expect that the ideal
$D$ is not of type $\ext{N}$ for some manifold $N$ because the
only standard real number in $D$ is $0$. We shall define suitable
structures on $D$ and $\ext{M}$ so that they will become objects
of the category $\ECInfty$ of extended smooth spaces, i.e. so that
$D$, $\ext{M}\in\ECInfty$. Subsequently we shall define the concept
of tangent vector so that $t\in\ECInfty(D,\ext{M})$, i.e. $t$ will
be an arrow of the category $\ECInfty$ of smooth extended spaces
and smooth extended functions.

In this chapter we will assume the following hypotheses on the category
$\mathcal{F}$:
\begin{assumption}
\noindent \begin{flushleft}
\textbf{\vspace{-8mm}
\label{ass:basicHypothesesOn_F}}
\par\end{flushleft}
\begin{enumerate}
\item \emph{\label{enu:F_SubcategoryOfTop}$\F$ }is a subcategory of the
category of topological spaces\emph{ }\textbf{\emph{$\textbf{Top}$}},
and contains all the constant maps $c:H\freccia X$ and all the open
subspaces $U\subseteq H$ (with the induced topology) of every object
$H\in\F$. The corresponding inclusion $i:U\hookrightarrow H$ is
also an arrow of $\F$, i.e. $i\in\Fs{UH}:=\F(U,H)$.
\end{enumerate}
\noindent \emph{In the following we will denote by $|-|\;:\;\F\freccia\Set$
the forgetful functor which associates to any $H\in F$ its support
set $|H|\in\Set$. Moreover with $\Top{H}$ we will denote the topology
of $H$ and with $(U\prec H)$ the topological subspace of $H$ induced
on the open set $U\in\Top{H}$. The remaining assumptions on $\F$
are the following:}
\begin{enumerate}
\item [2.]\setcounter{enumi}{2}\label{enu:F_closedWRT_Restrictions}The
category $\F$ is closed with respect to restrictions to open sets,
that is if $f\in\Fs{HK}$ and $U$, $V$ are open sets in $H$, $K$
resp. and finally $f(U)\subseteq V$, then $f|_{U}\in\F(U\prec H,V\prec K)$;
\item \label{enu:F_SheafProperty}Every topological space $H\in\F$ has
the following {}``sheaf property'': let $H$, $K\in\F$ be two objects
of $\mathcal{F}$, $(H_{i})_{i\in I}$ an open cover of $H$ and $f:|H|\freccia|K|$
a map such that\[
\forall i\in I:f|_{H_{i}}\in\F(H_{i}\prec H,K),\]
then $f\in\Fs{HK}$.\emph{ }
\end{enumerate}
\end{assumption}
For the construction of the domain of the extension functor $\ext(-)\;:\;\CInfty\freccia\ECInfty$
we want to consider a category $\F$ which permits to embed finite
dimensional manifolds in $\Cn$. For this aim we will set $\F=\ORn$,
the category having as objects open sets $U\subseteq\R^{u}$ (with
the induced topology), for some $u\in\N$ depending on $U$, and with
hom-set the usual space $\Cc^{n}(U,V)$ of $\Cc^{n}$ functions between
the open sets $U\subseteq\R^{u}$ and $V\subseteq\R^{v}$. Thus, $\Cn:=\overline{\ORn}$,
i.e. $\Cn$ is the cartesian closure of the category $\ORn$.

In general, what type of category $\F$ we have to choose depends
on the setting we need: e.g. in case we want to consider manifolds
with boundary we have to take the analogous of the above mentioned
category $\ORn$ but having as objects sets of type $U\subseteq\R_{+}^{u}=\{x\in\R^{u}\,|\, x_{u}\ge0\}$.

\section{\label{sec:DefinitionOfCartesianClosure}The cartesian closure and
its first properties}

The basic idea to define a $\Cn$ space $X$ (which faithfully generalizes
the notion of manifold) is to substitute the notion of chart by a
family of mappings $d:H\freccia X$ with $H\in\F$. Indeed, for $\F=\ORn$
these mappings are of type $d:U\freccia X$ with $U$ open in some
$\R^{u}$, thus they can be thought of as $u$-dimensional figures
on $X$ (see also Sections \ref{sec:DiffeologicalSpaces} and \ref{sec:TheConvenientVectorSpacesSettings}).
Hence, a $\Cn$ space can be thought as a support set together with
the specification of all the finite-dimensional figures on the space
itself. Generally speaking we can think of $\F$ as a category of
\emph{types of figures} (see \cite{La1} for this interpretation).
Always considering the case $\F=\ORn$, we can also think $\F$ as
a category which represents \emph{a well known notion of regular space
and regular function}: with the cartesian closure $\bar{\F}$ we want
to extend this notion to a more general type of spaces (e.g. spaces
of mappings). These are the ideas we have already seen in Section
\ref{sec:DiffeologicalSpaces} in the case of diffeological spaces,
only suitably generalized to a category of topological spaces $\mathcal{F}$
instead of $\F=\ORInfty$, which is the case of diffeology. This generalization
permits to obtain in an easy way the cartesian closedness of $\bar{\F}$,
and thus to have at our disposal a general instrument $\F\mapsto\bar{\F}$
very useful in the construction of the codomain of the extension functor
$\ext(-)$, where we will choose a different category of types of
figures $\F$.
\begin{defn}
In the sequel we will frequently use the notation $f\cdot g:=g\circ f$
for the composition of maps so as to facilitate the lecture of diagrams,
but we will continue to evaluate functions {}``on the right'' hence
$(f\cdot g)(x)=g(f(x))$.
\end{defn}
Objects and arrows of $\bar{\F}$ generalize the same notions of the
diffeological setting (see Section \ref{sec:DiffeologicalSpaces}). 
\begin{defn}
If $X$ is a set, then we say that $(\D,X)$ is an object of $\bar{\F}$
(or simply an $\bar{\mathcal{F}}$-object) if $\D=\{\Ds{H}\}_{\sss H\in\F}$
is a family with \[
\Ds{H}\subseteq\Set(|H|,X)\quad\forall H\in\mathcal{F}.\]
 We indicate by the notation $\Fs{JH}\cdot\Ds{H}$ the set of all
the compositions $f\cdot d$ of functions $f\in\Fs{JH}$ and $d\in\Ds{H}$.
The family $\D$ has finally to satisfy the following conditions:
\begin{enumerate}
\item \label{enu:Def_F_Bar_reparametrization}$\Fs{JH}\cdot\Ds{H}\subseteq\Ds{J}$. 
\item \label{enu:Def_F_Bar_ClosureConstantMaps}$\Ds{H}$ contains all the
constant maps $d:|H|\freccia X$. 
\item \label{enu:Def_F_Bar_SheafProperty}Let $H\in\F$, $(H_{i})_{i\in I}$
an open cover of $H$ and $d:|H|\freccia X$ a map such that $d|_{H_{i}}\in\Ds{(H_{i}\prec H)}$,
then $d\in\Ds{H}$. 
\end{enumerate}
Finally, we set $|(\D,X)|:=X$ to denote the underlying set of the
space $(\D,X)$.
\end{defn}
Because of condition \emph{\ref{enu:Def_F_Bar_reparametrization}.}
we can think of $\Ds{H}$ as the set of \emph{all} the regular functions
defined on the {}``well known'' object $H\in\F$ and with values
in the new space $X$; in fact this condition says that the set of
figures $\Ds{H}$ is closed with respect to re-parametrizations with
a generic $f\in\Fs{JH}$. Condition \ref{enu:Def_F_Bar_SheafProperty}.
is the above mentioned sheaf property and asserts that the property
of being a figure $d\in\mathcal{D}_{H}$ has a local character depending
on $\F$.

We will frequently write $d\In{H}X$ to indicate that $d\in\Ds{H}$
and we can read it%
\footnote{The following are common terminologies used in topos theory, see \cite{La1, Koc, Mo-Re}%
} saying that $d$ \emph{is a figure of }$X$ \emph{of type} $H$ or
$d$ \emph{belong to }$X$ \emph{at the level }$H$ or $d$ \emph{is
a generalized element of }$X$ \emph{of type }$H$.

The definition of arrow $f:X\freccia Y$ (also called \emph{smooth
function in $\bar{\F}$})\emph{ }between two spaces $X$, $Y\in\bar{\F}$
is the usual one for diffeological spaces, that is $f$ takes, through
composition, generalized elements $d\In{H}X$ of type $H$ in the
domain $ $$X$ to generalized elements of the same type in the codomain
$Y$
\begin{defn}
Let $X$, $Y$ be $\bar{\F}$-objects, then we will write\[
f:X\freccia Y\]

\noindent or, more precisely if needed\,%
\footnote{We shall frequently use notations of type $\mathbb{C}\vDash f:A\freccia B$
if we need to specify better the category $\mathbb{C}$ we are considering
(see Appendix \ref{app:someNotionsOfCategoryTheory}).%
}\[
\bar{\F}\vDash f:X\freccia Y\]
iff $f$ maps the support set of $X$ into the support set of $Y$:\[
f:|X|\freccia|Y|\]
and\[
d\cdot f\In{H}Y\]
for every type of figure $H\in\F$ and for every figure $d$ of $X$
of that type, i.e. $d\In{H}X$. In this case, we will also use the
notation $f(d):=d\cdot f$.
\end{defn}
\noindent Note that we have $f:X\freccia Y$ in $\bar{\F}$ iff\[
\forall H\in\F\,\,\forall x\In{H}X\pti f(x)\In{H}Y,\]
moreover $X=Y$ iff\[
\forall H\in\F\,\,\forall d\pti d\In{H}X\iff d\In{H}Y.\]
These and many other properties justify the notation $\In{H}$ and
the name {}``generalized elements''.

With these definitions $\bar{\F}$ becomes a category. Note that it
is, in general, in the second Grothendieck universe (see \cite{SGA4, Ad-He-St})
because $\D$ is a family indexed in the set of objects of $\F$ (this
is not the case for $\F=\ORn$, which is a set and not a class).

\noindent The simplest $\bar{\F}$-object is $\bar{K}:=(\Fs{(-)K},|K|)$
for $K\in\F$, where we recall that $\F_{HK}=\F(H,K)=\left\{ f\,|\, H\xfreccia{f}K\text{ in }\F\right\} $.
For the space $\bar{K}\in\bar{\F}$ we have that\[
\bar{\F}\vDash\bar{K}\xfrecciad{d}X\quad\iff\quad d\In{K}X.\]
Moreover, $\F(H,K)=\bar{\F}(\bar{H},\bar{K})$. Therefore $\F$ is
fully embedded in $\bar{\F}$ if $\bar{H}=\bar{K}$ implies $H=K$;
e.g. this is true if the given category $\F$ verifies the following
hypothesis \[
|H|=|K|=S\text{\ \ and\ \ }H\xfrecciad{1_{S}}K\xfrecciad{1_{S}}H\then H=K.\]
E.g. this is true for $\F=\ORn$.

Moreover, let us note that the composition of two smooth functions
in $\bar{\F}$ of type $d:\bar{H}\freccia X$ and $f:X\freccia\bar{K}$
for $H$, $K\in\F$, gives $d\cdot f\in\bar{\F}(\bar{H},\bar{K})=\F(H,K)$,
which is an arrow in the old category of types of figures $\F$.

Another way to construct an object of $\bar{\F}$ on a given support
set $X$ is to generate it starting from a given family $\mathcal{D}^{0}=(\Ds{H}^{0})_{\sss{H}}$,
with $\Ds{H}^{0}\subseteq\Set(|H|,X)$ for any $H\in\F$, closed with
respect to constant functions, i.e. such that\[
\forall H\in\F\,\,\forall d:|H|\freccia X\text{ is constant}\then d\in\D_{H}^{0}.\]
We will indicate this space by $(\F\cdot\D^{0},X)$. Its figures are,
locally, compositions $f\cdot d$ with $f\in\Fs{HK}$ and $d\in\Ds{K}^{\sss{0}}$.
More precisely $\delta\In{H}(\F\cdot\D^{\sss{0}},X)$ iff $\delta:|H|\freccia X$
and for every $h\in|H|$ there exist an open neighborhood $U$ of
$h$ in $H$, a space $K\in\F$, a figure $d\in\Ds{K}^{0}$ and $f:(U\prec H)\freccia K$
in $\F$ such that $\delta|_{U}=f\cdot d$. Diagrammatically we have:\[
\xyR{45pt}\xyC{55pt}\xymatrix{H\ar[r]^{\delta} & X\\
h\in U\ar[r]_{f}\ar@{_{(}->}[u]\ar[ru]^{\delta|_{U}} & K\ar[u]_{d}}
\]

On each space $X\in\bar{\F}$ we can put the final topology $\Top{X}$
for which any figure $d\In{H}X$ is continuous, that is
\begin{defn}
\label{def:topologyOnTheCartesianClosure}If $X\in\bar{\F}$, then
we say that a subset $U\subseteq|X|$ is \emph{open} in $X$, and
we will write $U\in\Top{X}$ iff $d^{-1}(U)\in\Top{H}$ for any $H\in\F$
and any $d\In{H}X$.
\end{defn}
With respect to this topology any arrow of $\bar{\F}$ is continuous
and we still have the initial $\Top{H}$ in the space $\bar{H}$,
that is $\Top{H}=\Top{\bar{H}}$ (recall that, because of the fundamental
hypotheses \ref{ass:basicHypothesesOn_F}, every type of figure $H\in\F$
is a topological space).

Recalling that in the case $\F=\ORInfty$ we obtain that the cartesian
closure $\bar{\F}$ is the category of diffeological spaces, it can
be useful to cite here \citet{Igl}:
\begin{quotation}
Even if diffeology is a theory which avoids topology on purpose, topology
is not completely absent from its content. But, in contrary to some
approach of standard differential geometry, here the topology is a
byproduct of the main structure, that is diffeology. Locality, through
local smooth maps, or local diffeomorphisms, is introduced without
referring to any topology a priori but will suggest the definition
of a topology a posteriori {[}i.e. $\Top{X}${]}.
\end{quotation}
Ultimately, this choice is due to the necessity to obtain a cartesian
closed category. In fact, if we do not start from a primitive notion
of topology in the definition of $\bar{\F}$-space, we can obtain
cartesian closedness without having the problem to define a topology
in the set of maps $\bar{\F}(X,Y)$. Indeed, this is not an easy problem,
and classical solutions like the compact-open topology (see e.g. \citet{Dug,Kr-Mi}
and references therein) is not applicable to the smooth case. In fact,
the compact-open topology, which essentially coincides with the topology
of uniform convergence, is well suited for continuous maps $f:X\freccia Y$
between locally compact Haussdorff topological spaces $X$ and $Y$
(indeed, the category of these topological spaces is cartesian closed,
see \citet{Mac}). It can be generalized to the case of $\Cc^{k}$-regularity
using $k$-jets ($k\in\N_{>0}$), i.e. using Taylor's formulae up
to $k$-th order (see e.g. \citet{Kr-Mi}), but a generalization including
the smooth case $\Cc^{\infty}$ even for a compact domain $X$ fails.
In fact, for $X$ compact and $Y$ a Banach space, the space $\Cc^{k}(X,Y)$
with the $\Cc^{k}$ compact-open topology is normable, but the space
$\Cc^{\infty}(\R,\R)$ is not normable, so its topology cannot be
the compact-open one (see also Section \ref{sec:BanachManifoldsAndLocallyConvexVectorSpaces}
for more details).

The study of the relationships between different topologies on the
space of maps $\Cc^{\infty}(M,N)$ for $M$, $N$ manifolds, is not
completely solved (see again \citet{Kr-Mi} for some results in this
direction).

\section{\label{sec:CategoricalPropertiesOfTheCartesianClosure}Categorical
properties of the cartesian closure}

We shall now examine subobjects in $\bar{\F}$ and their relationships
with restrictions of functions; after this we will analyze completeness,
co-completeness and cartesian closure of $\bar{\F}$. 
\begin{defn}
\label{def:subspaceIn__F_bar}Let $X\in\bar{\F}$ be a space in the
cartesian closure of $\F$, and $S\subseteq|X|$ a subset, then we
define \[
(S\prec X):=(\D,S)\]
 where, for every type of figure $H\in\F$, we have set \[
d\in\Ds{H}\DIff d\;:\;|H|\freccia S\e{and}d\cdot i\In{H}X.\]
Here $i:S\hookrightarrow|X|$ is the inclusion map. In other words,
we have a figure $d$ of type $H$ in the subspace $S$ iff composing
$d$ with the inclusion map $i$ we obtain a figure of the same type
in the superspace $X$. We will call $(S\prec X)$ \emph{the subspace
induced on }$S$ \emph{by} $X$.
\end{defn}
Using this definition only it is very easy to prove that $(S\prec X)\in\bar{\F}$
and that its topology $\Top{(S\prec X)}$ contains the induced topology
by $\Top{X}$ on the subset $S$. Moreover we have that $\Top{(S\prec X)}\subseteq\Top{X}$
if $S$ is open in $X$, hence in this case we have on $(S\prec X)$
exactly the induced topology.

Finally we can prove that these subspaces have good relationships
with restrictions of maps:
\begin{thm}
\noindent \label{thm:cartesianClosureAndRestrictionOfMaps}Let $f:X\freccia Y$
be an arrow of $\bar{\F}$ and $U$, $V$ be subsets of $|X|$ and
$|Y|$ respectively, such that $f(U)\subseteq V$, then \[
(U\prec X)\xfrecciad{f|_{\sss{U}}}(V\prec Y)\e{in}\bar{\F}.\]
\qedWithFinalEq
\end{thm}
\noindent Using our notation for subobjects we can prove the following
useful and natural properties directly from definition \ref{def:subspaceIn__F_bar}. 
\begin{itemize}
\item $(U\prec\bar{H})=\overline{(U\prec H)}$ for $U$ open in $H\in\F$
(recall the definition of $\bar{H}\in\bar{\F}$, for $H\in\F$, given
in Section \ref{sec:DefinitionOfCartesianClosure} and also recall
that, because of the Hypotheses \ref{ass:basicHypothesesOn_F} the
subspace $(U\prec H)$ is a type of figure, i.e. $(U\prec H)\in\F$,
and we can thus apply the operator $\bar{(-)}:\F\freccia\bar{\F}$
of inclusion of the types of figures $\F$ into the cartesian closure
$\bar{\F}$).
\item $i:(S\prec X)\hookrightarrow X$ is the lifting%
\footnote{For the notion of \emph{lifting} and \emph{co-lifting} see Definition
\ref{def:lifting}%
} of the inclusion $i:S\hookrightarrow|X|$ from $\Set$ to $\bar{\F}$ 
\item $(|X|\prec X)=X$ 
\item $(S\prec(T\prec X))=(S\prec X)$ \ \ if \ \ $S\subseteq T\subseteq|X|$
\item $(S\prec X)\times(T\prec Y)=(S\times T\prec X\times Y)$. 
\end{itemize}
These properties imply that the relation $X\subseteq Y$ iff $|X|\subseteq|Y|$
and $(|X|\prec Y)=X$ is a partial order. Note that this relation
is stronger than saying that the inclusion is an arrow, because it
asserts that $X$ and the inclusion verify the universal property
of $(|X|\prec Y)$, that is $X$ is a subobject of $Y$. A trivial
but useful property of this subobjects notation is the following
\begin{cor}
\label{cor:changeOfSuperSpace}Let $S\subseteq|X'|$ and $X'\subseteq X$
in $\bar{\F}$, then\[
(S\prec X')=(S\prec X),\]
that is in the operator $(S\prec-)$ we can change the superspace
$X$ with any one of its subspaces $X'\subseteq X$ containing $S$.
\end{cor}
\textbf{Proof:} In fact $X'\subseteq X$ means $X'=(|X'|\prec X)$
and hence $(S\prec X')=(S\prec(|X'|\prec X))=(S\prec X)$ because
of the previous properties of the operator $(-\prec-)$.$\qedNoNewLine$

An expected property that transfers from $\F$ to $\bar{\F}$ is the
sheaf property; in other words it states that the property of being
a smooth arrow of the cartesian closure $\bar{\F}$ is a local property.
\begin{thm}
Let $X$, $Y\in\bar{\F}$ be spaces in the cartesian closure, $(U_{i})_{i\in I}$
an open cover of $X$ and $f:|X|\freccia|Y|$ a map from the support
set of $X$ to that of $Y$ such that\[
\bar{\F}\vDash(U_{i}\prec X)\xfrecciad{f|_{U_{i}}}Y\quad\forall i\in I.\]

\noindent Then\[
\bar{\F}\vDash X\xfrecciad{f}Y.\]
\qedWithFinalEq
\end{thm}
Completeness and co-completeness are analyzed in the following theorem.
For its standard proof see e.g. \citet{Fr-Kr} for a similar theorem.
\begin{thm}
\label{thm: limit in Cn} Let $(X_{i})_{i\in I}$ be a family of objects
in $\bar{\F}$ and $p_{i}:|X|\freccia|X_{i}|$ maps for every $i\in I$.
Let us define \[
d\In{H}X\DIff d\;:\;|H|\freccia|X|\e{and}\forall\, i\in I\pti d\cdot p_{i}\In{H}X_{i}\]
 then $\,{\displaystyle (X\xfreccia{p_{i}}X_{i})_{i\in I}\,}$ is
a lifting of $\,{\displaystyle (|X|\xfreccia{p_{i}}|X_{i}|)_{i\in I}\,}$
in $\bar{\F}$.

\noindent Moreover, let $j_{i}:|X_{i}|\freccia|X|$ be maps for every
$i\in I$, and let us suppose that \[
\forall\, x\in|X|\,\,\exists\,\, i\in I\,\,\exists\, x_{i}\in X_{i}\pti x=j_{i}(x_{i}).\]
Let us define $d\In{H}X$ iff $d:|H|\freccia|X|$ and for every $h\in|H|$
there exist an open neighborhood $U$ of $h$ in $H$, an index $i\in I$
and a figure $\delta\In{U}X_{i}$ such that $d|_{U}=\delta\cdot j_{i}$;
then we have that $\,{\displaystyle (X_{i}\xfreccia{j_{i}}X)_{i\in I}\,}$
is a co-lifting of $\,{\displaystyle (|X_{i}|\xfreccia{j_{i}}|X|)_{i\in I}\,}$
in $\bar{\F}$.\qedWithFinalEq
\end{thm}
The category of $ $$\bar{\F}$ spaces is thus complete and co-complete
and we can hence consider spaces like quotient spaces $X/\sim$, disjoint
sums $\sum_{i\in I}X_{i}$, arbitrary products $\prod_{i\in I}X_{i}$,
equalizers, etc. (see Theorem \ref{thm:relationBetweenLimitsAndLifting}
for further details about the connections between limits, co-limits,
lifting and co-lifting).

Directly from the definitions of lifting and co-lifting, it is easy
to prove that on quotient spaces we exactly have the quotient topology
and that on any product we have a topology stronger than the product
topology. We can write this assertion in the following symbolic way:\begin{equation}
\Top{X/\sim}=\Top{X}/\sim\end{equation}

\begin{equation}
\Top{X}\times\Top{Y}\subseteq\Top{X\times Y},\label{eq:productTopologyAndProductSpaces}\end{equation}

\noindent where: $X$ and $Y$ are $\bar{\F}$ spaces, $\sim$ is
an equivalence relation on $|X|$, $(X/\sim)\in\bar{\F}$ is the quotient
space, $\Top{X}/\sim$ is the quotient topology, and $\Top{X}\times\Top{Y}$
is the product topology. Analogously, let $j_{i}:X_{i}\freccia\sum_{i\in I}X_{i}$
be the canonical injections in the disjoint sum of the family of $\bar{\F}$-spaces
$(X_{i})_{i\in I}$, i.e. $j_{i}(x)=(x,i)$. Then we can prove that
$A$ is open in $\sum_{i\in I}X_{i}$ if and only if\begin{equation}
\forall i\in I\pti j_{i}^{-1}(A)\in\Top{X_{i}},\label{eq:openSetsInDisjointSums}\end{equation}

\noindent that is on the disjoint sum we have exactly the colimit
topology. Because any colimit can be obtain as a lifting from $\Set$
of quotient spaces and disjoint sums (see \citet{Mac}), we have the
general result that the topology on the colimit of $\bar{\F}$-spaces
is exactly the colimit topology. In symbolic notations we can write\[
{\mbox{\Large$\tau$}}\left({\displaystyle \colim_{i\in I}}X_{i}\right)=\colim_{i\in I}\Top{X_{i}}.\]

Finally if we define \[
\Ds{H}:=\{d:|H|\freccia\bar{\F}(X,Y)\;|\;\bar{H}\times X\xfrecciad{d^{\vee}}Y\e{in}\bar{\F}\}\quad\forall H\in\F\]
 (we recall that we use the notations $d^{\vee}(h,x):=d(h)(x)$ and
$\mu^{\wedge}(x)(y):=\mu(x,y)$, see Section \ref{sec:ApproachesToDiffGeomOfInfDimSpaces})
then $\langle\D,\bar{\F}(X,Y)\rangle=:Y^{X}$ is an object of $\bar{\F}$.
With this definition, see e.g. \citet{Che} or \citet{Fr-Kr}, it
is easy to prove that $\bar{\F}$ is cartesian closed, i.e. that the
$\bar{\F}$-isomorphism $(-)^{\vee}$ realizes \[
(Y^{X})^{Z}\simeq Y^{Z\times X}.\]

\chapter{\label{cha:TheCategoryCn}The category $\Cn$}

\section{Observables on $\Cn$ spaces and separated spaces}

If our aim is to embed the category of $\Cc^{n}$ manifolds into a
cartesian closed category, the most natural way to apply the results
of the previous Chapter \ref{cha:theCartesianClosure} is to take
as category $\F$ of types of figures $\F=\Man$, that is to consider
directly the cartesian closure of the category of finite dimensional
$\Cc^{n}$ manifolds%
\footnote{We shall not formally assume any hypothesis on the topology of a manifold
because we will never need it in the following; moreover if not differently
specified, with the word {}``manifold'' we will always mean {}``finite
dimensional manifold''.%
}. We shall not follow this idea for several reasons; as we have already
mentioned, we will consider instead $\Cn:=\overline{\ORn}$, that
is the cartesian closure of the category $\ORn$ of open sets and
$\Cc^{n}$ maps. For $n=\infty$ this gives exactly diffeological
spaces. Indeed, as we noted in the previous Chapter \ref{cha:theCartesianClosure},
$\overline{\Man}$ is in the second Grothendieck universe and, essentially
for simplicity, from this point of view the choice $\F=\ORn$ is better.
In spite of this choice, it is natural to expect, and in fact we will
prove it, that the categories of both finite and infinite-dimensional
manifolds are faithfully embedded in the previous $\Cn=\overline{\ORn}$.
Another reason to choose our definition of $\Cn$ is that in this
way the category $\Cn$ is more natural to accept against $\overline{\Man}$;
hence, ones again we are opting for a reason of simplicity. We will
see that manifolds modelled in convenient vector spaces (see Chapter
\ref{cha:ApproachesToDiffGeomOfInfDimSpaces}) are faithfully embedded
in $\Cn$, hence our choice to take finite dimensional objects in
the definition of $\Cn=\overline{\ORn}$ is not restrictive from this
point of view.

Now we pay attention to another type of maps which go {}``in the
opposite direction'' with respect to figures $d:K\freccia X$. They
are important also because we shall use them to introduce new infinitesimal
points for any $X\in\Cn$. We will introduce these notions for a generic
cartesian closure $\bar{\F}$ of a given category if figures $\F$,
because we will use them e.g. also in the category $\ECn$ of extended
spaces. So, in the following $\bar{\F}$ will be a category of figures
(see Hypothesis \ref{ass:basicHypothesesOn_F}).
\begin{defn}
Let $X$ be an $\bar{\F}$ space, then we say that \[
UK\text{\emph{ is a zone (in $X$)}}\]
 iff $U\in\Top{X}$, i.e. $U$ is open in $X$, and $K\in\F$ is a
type of figure. Moreover we say that \[
c\text{\emph{ is an observable on }$UK$}\quad\text{and we will write}\quad c\InUp{UK}X\]
 iff $c\;:\;(U\prec X)\freccia\bar{K}$ is a map of the cartesian
closure $\bar{\F}$.
\end{defn}
\noindent So, observables of a $\Cn$ space $X$ are simply maps of
class $\Cn$ (i.e. are arrows of this category) defined on an open
set of $X$ and with values in an open set $K\subseteq\R^{d}$ for
some $d\in\N$. Recall (see Section \ref{sec:DefinitionOfCartesianClosure})
that for any open set $K\in\ORn$ in the $\Cn$ space $\bar{K}$ we
take as figures of type $H\in\ORn$ all the ordinary $\Cc^{n}$-maps
$\Cc^{n}(H,K)$, i.e. we have \[
\bar{K}=(\Cc^{n}(-,K),K).\]
Therefore, the composition of figures $d\In{H}X$ with observables
$c\InUp{UK}X$ gives ordinary $\Cc^{n}$ maps:\[
d|_{S}\cdot c\in\Cc^{n}(S,K),\e{where}S:=d^{-1}(U),\]
\[
\Cn\vDash\overline{(S\prec H)}\xfrecciad{d|s}(U\prec X)\xfrecciad{c}\bar{K}.\]
From our previous theorems of Chapter \ref{cha:theCartesianClosure},
it follows that $\Cn$ functions $f:X\freccia Y$ take observables
on the codomain to observables on the domain i.e.: \begin{equation}
c\InUp{UK}Y\then f|_{\sss{S}}\cdot c\InUp{SK}X,\label{thm: fPreservesObs}\end{equation}

\noindent where $S:=f^{-1}(U)$:\[
\xymatrix{\xyR{40pt}\xyC{50pt}(S\prec X)\ar[r]^{f|_{S}}\ar[dr]_{f|_{S}\cdot c} & (U\prec Y)\ar[d]^{c}\\
 & \bar{K}}
\]
 Therefore isomorphic $\Cn$ spaces have isomorphic sets of figures
and observables and the isomorphisms are given by suitable simple
compositions.

Generalizing, through observables, the equivalence relation of Definition
\ref{def:equalityInFermatReals} to generic $\Cn$ spaces, we will
have to study the following condition, which is also connected with
the faithfulness of the extension functor.
\begin{defn}
\noindent \label{def: xIdentifiedy} If $X\in\Cn$ is a $\Cn$ space
and $x$, $y\in|X|$ are two points, then we write \[
x\asymp y\]
iff for every zone $UK$ and every observable $c\InUp{UK}X$ we have 
\begin{enumerate}
\item $x\in U\iff y\in U$ 
\item $x\in U\then c(x)=c(y).$
\end{enumerate}
\noindent In this case we will read the relation $x\asymp y$ saying
$x$\emph{ and $y$ are identified in $X$.} Moreover we say that
$X$ \emph{is separated} iff $x\asymp y$ implies $x=y$ for any $x$,
$y\in|X|$.
\end{defn}
We point out that if two points are identified in $X$, then a generic
open set $U\in\Top{X}$ contains the first if and only if it contains
the second too (take a constant observable $c:U\freccia\R$). Furthermore,
from \eqref{thm: fPreservesObs} it follows that $\Cn$ functions
$f:X\freccia Y$ preserve the relation $\asymp$: \[
x\asymp y\text{ in }X\then f(x)\asymp f(y)\text{ in }Y\qquad\forall x,y\in|X|.\]
Trivial examples of separated spaces can be obtained considering the
objects $\bar{U}\in\Cn$ with $U\in\ORn$ (here $\overline{(-)}:\ORn\freccia\Cn$
is the embedding of the types of figures $\ORn$ into $\Cn$, see
\ref{sec:DefinitionOfCartesianClosure}) or taking subobjects of separated
spaces. But the full subcategory of separated $\Cn$ spaces has sufficiently
good properties, as proved in the following
\begin{thm}
\noindent The category of separated $\Cn$ spaces is complete and
admits co-products. Moreover if $X$, $Y$ are separated then $Y^{X}$
is separated too, and hence separated spaces form a cartesian closed
category. 
\end{thm}
\noindent \textbf{Sketch of the proof:} We only do some considerations
about co-products, because from the definition of lifting (see Theorem
\ref{thm: limit in Cn}) it can be directly proved that products and
equalizers of separated spaces are separated too. Let us consider
a family $(\mathcal{X}_{i})_{i\in I}$ of separated spaces with support
sets $X_{i}:=|\mathcal{X}_{i}|$. Constructing their sum in $\Set$
\[
X:=\sum_{i\in I}X_{i}\]
 \[
j_{i}\;:\; x\in X_{i}\longmapsto(x,i)\in X,\]
 from the completeness of $\Cn$ we can lift this co-product of sets
into a co-product ${\displaystyle (\mathcal{X}_{i}\xfreccia{j_{i}}\mathcal{X})_{i\in I}}$
in $\Cn$. To prove that $\mathcal{X}$ is separated we take two points
$x$, $y\in X=|\mathcal{X}|$ identified in $\mathcal{X}$. These
points are of the form $x=(x_{r},r)$ and $y=(y_{s},s)$, with $x_{r}\in X_{r},y_{s}\in X_{s}$
and $r,s\in I$. We want to prove that $r$ and $s$ are necessarily
equal. In fact, from \eqref{eq:openSetsInDisjointSums}, for a generic
$A\subseteq\mathcal{X}$ we have that \[
A\in\Top{\mathcal{X}}\quad\iff\quad\forall\, i\in I\pti j_{i}^{-1}(A)\in\Top{\mathcal{X}_{i}}.\]
and hence $X_{r}\times\{r\}$ is open in $\mathcal{X}$ and $x=(x_{r},r)\asymp y=(y_{s},s)$
implies \[
(x_{r},r)\in X_{r}\times\{r\}\iff(y_{s},s)\in X_{r}\times\{r\}\ee{hence}r=s.\]
 Thus $x=y$ iff $x_{r}$ and $y_{s}=y_{r}$ are identified in $\mathcal{X}_{r}$
and this is a consequence of the following facts: 
\begin{enumerate}
\item if $U$ is open in $\mathcal{X}_{r}$ then $U\times\{r\}$ is open
in $\mathcal{X}$; 
\item if $c\InUp{UK}\mathcal{X}_{r}$, then $\gamma(x,r):=c(x)\,\,\forall\, x\in U$
is an observable of $\mathcal{X}$ defined on $U\times\{r\}$. 
\end{enumerate}
Now let us consider exponential objects. If $f$, $g\in|Y^{X}|$ are
identified, to prove that they are equal is equivalent to prove that
$f(x)$ and $g(x)$ are identified in $Y$ for any $x$. To obtain
this conclusion is sufficient to consider that the evaluation in $x$
i.e. the application $\eps_{x}:\varphi\in|Y^{X}|\longmapsto\varphi(x)\in|Y|$
is a $\Cn$ map and hence from any observable $c\InUp{UK}Y$ we can
always obtain the observable $\eps_{x}|_{\sss U^{\prime}}\cdot c\InUp{U^{\prime}K}Y^{X}$
where $U^{\prime}:=\eps_{x}^{-1}(U)$.\qedNoNewLine 

Finally let us consider two $\Cn$ spaces such that the topology $\Top{X\times Y}$
is equal to the product of the topologies $\Top{X}$ and $\Top{Y}$
(recall that in general we have $\Top{X}\times\Top{Y}\subseteq\Top{X\times Y}$).
Then if $x$, $x^{\prime}\in|X|$ and $y$, $y^{\prime}\in|Y|$ directly
from the definition it is possible to prove that $x\asymp x'$ in
$X$ and $y\asymp y'$ in $Y$ if and only if $(x,y)\asymp(x',y')$
in $X\times Y$.

\section{Manifolds as objects of $\Cn$}

We can associate in a very natural way a $\Cn$ space $\bar{M}$ to
any manifold $M\in\Man$ (the category of $\Cc^{n}$ manifolds and
$\Cc^{n}$ functions) with the following
\begin{defn}
\label{def:embeddingOfManifolds}The underlying set of $\bar{M}$
is the underlying set of the manifold, i.e. $|\bar{M}|:=|M|$, and
for every $H\in\ORn$ the figures $d:H\freccia M$ of type $H$ are
all the ordinary $\Cc^{n}$ maps from $H$ to the manifold $M$, i.e.
\[
d\In{H}\bar{M}\DIff d\in\Man(H,M).\]

\end{defn}
\noindent This definition is only the trivial generalization from
the smooth case to $\Cc^{n}$ of the embedding of manifolds into the
category of diffeological spaces (see e.g. \citet{Igl}).

\noindent With $\bar{M}$ we obtain a $\Cn$ space with the same topology
of the starting manifold. Moreover the observables of $\bar{M}$ are
the most natural ones we could expect. In fact, as a consequence of
the Definition \ref{def:embeddingOfManifolds} it follows that\begin{equation}
c\InUp{UK}\bar{M}\quad\iff\quad c\in\Man(U,K).\label{eq:observablesInManifolds}\end{equation}
Hence it is clear that the space $\bar{M}$ is separated, because
from \eqref{eq:observablesInManifolds} we get that charts are observables
of the space. The following theorem says that the application $M\mapsto\bar{M}$
from $\Man$ to $\Cn$ we are considering is a full embedding, and
therefore it also says that the notion of $\Cn$-space is a non-trivial
generalization of the notion of manifold which includes infinite-dimensional
spaces too.
\begin{thm}
\label{thm: ManifoldsEmbedded} Let $M$ and $N$ be $\Cc^{n}$ manifolds,
then 
\begin{enumerate}
\item \label{enu:embeddingOfManifoldsOnObjects}$\bar{M}=\bar{N}\then M=N$ 
\item \label{enu:embeddingOfManifoldsOnArrows}$\Cn\vDash{\displaystyle \bar{M}\xfrecciad{f}\bar{N}\quad\iff\quad\Man\vDash M\xfrecciad{f}N}$. 
\end{enumerate}
\noindent Hence $\Man$ is fully embedded in $\Cn$.
\end{thm}
\noindent \textbf{Proof:}

\noindent \emph{\ref{enu:embeddingOfManifoldsOnObjects})} If $(U,\varphi)$
is a chart on $M$ and $A:=\varphi(U)$, then $\varphi^{-1}|_{A}:A\freccia M$
is a figure of $\bar{M}$, that is $\varphi^{-1}|_{A}\In{A}\bar{M}=\bar{N}$.
But if $\psi:U\freccia\psi(U)\subseteq\R^{k}$ is a chart of $N$,
then it is also an observable of $\bar{N}$. We have hence obtained
a figure $\varphi^{-1}|_{A}\In{A}\bar{N}$ and an observable $\psi\InUp{U\psi(U)}\bar{N}$
of the space $\bar{N}$. But composition of figures and observables
gives ordinary $\Cc^{n}$ maps, that is the atlases of $M$ and $N$
are compatible.

\medskip{}

\noindent \emph{\ref{enu:embeddingOfManifoldsOnArrows})} For the
implication $\Rightarrow$ we use the same ideas as above and furthermore
that $\varphi^{-1}|_{A}\In{A}\bar{M}$ implies $\varphi^{-1}|_{A}\cdot f\In{A}\bar{N}$.
Finally we can compose this $A$-figure of $\bar{N}$ with a chart
(observable) of $N$ obtaining an ordinary $\Cc^{n}$ map. The implication
$\Leftarrow$ follows directly from the Definition \ref{def:embeddingOfManifolds}.\qedNoNewLine

Directly from these definitions we can prove that for two manifolds
we also have \[
\overline{M\times N}=\bar{M}\times\bar{N}.\]
 This property is useful to prove the properties stated in the following
examples.

\section{\label{sec:ExamplesOfCnSpacesAndFunctions}Examples of $\Cn$ spaces
and functions}
\begin{enumerate}
\item Let $M$ be a $\Cc^{\infty}$ manifold modelled on convenient vector
spaces (see Section \ref{sec:TheConvenientVectorSpacesSettings}).
We can define $\bar{M}$ analogously as above, saying that $d\In{H}\bar{M}$
iff $d:H\freccia M$ is a smooth map from $H$ (open in some $\R^{h}$)
to the manifold $M$. In this way smooth curves on $M$ are exactly
the figures $c\in_{\R}\bar{M}$ of type $\R$ in $\bar{M}$. On $M$
we obviously think of the natural topology, that is the identification
topology with respect to some smooth atlas, which is also the final
topology with respect to all smooth curves and hence is also the final
topology $\Top{\bar{M}}$ with respect to all figures of $\bar{M}$.
More easily with respect to the previous case of finite dimensional
manifolds (due to the results available for manifolds modelled on
convenient vector spaces, see Section \ref{sec:TheConvenientVectorSpacesSettings}),
it is possible to study observables, obtaining that $c\InUp{UK}\bar{M}$
if and only if $c:U\freccia K$ is smooth as a map between manifolds
modelled on convenient vector spaces. Moreover if $(U,\varphi)$ is
a chart of $M$ on the convenient vector space $E$, then $\varphi:(U\prec\bar{M})\freccia(\varphi(U)\prec\bar{E})$
is $\CInfty$. Using these results it is easy to prove the analogous
of Theorem \ref{thm: ManifoldsEmbedded} for the category of manifolds
modelled on convenient vector spaces. Hence also classical smooth
manifolds modelled on Banach spaces are embedded in $\CInfty$. 
\item It is not difficult to prove that the following applications, frequently
used e.g. in calculus of variations, are smooth, that is they are
arrows of $\CInfty$. 

\begin{enumerate}
\item The operator of derivation:\begin{align*}
\partial_{i}:\  & u\in\Cc^{\infty}(\R^{n},\R^{k})\longmapsto\frac{\partial u}{\partial x_{i}}\in\Cc^{\infty}(\R^{n},\R^{k})\end{align*}
To prove that this operator is smooth, i.e. it is an arrow of the
category $\CInfty$, we have to show that it takes figures of type
$H\in\ORInfty$ on its domain to figures of the same type on the codomain.
Figures of type $H$ of the space $\Cc^{\infty}(\R^{n},\R^{k})$ are
maps of type $d:H\freccia\Cc^{\infty}(\R^{n},\R^{k})$, so that we
have to consider the composition $d\cdot\partial_{i}$. Using cartesian
closedness we get that $d^{\vee}:H\times\R^{n}\freccia\R^{k}$ is
an ordinary smooth map. But, always due to cartesian closedeness,
the composition $d\cdot\partial_{i}:H\freccia\Cc^{\infty}(\R^{n},\R^{k})$
is a figure if and only if its adjoint $\left(d\cdot\partial_{i}\right)^{\vee}:H\times\R^{n}\freccia\R^{k}$
is an ordinary smooth map, and by a direct calculation we get that
$\left(d\cdot\partial_{i}\right)^{\vee}=\partial_{u+i}d^{\vee}$,
where $u\in\N$ is the dimension of $H\subseteq\R^{u}$. In fact\begin{align*}
\left(d\cdot\partial_{i}\right)^{\vee}(h,r) & =\partial_{i}(d(h))(r)=\frac{\partial d(h)}{\partial x_{i}}(r)=\\
 & =\lim_{\delta\to0}\frac{d(h)(r+\delta\vec{e_{i}})-d(h)(r)}{\delta}=\\
 & =\lim_{\delta\to0}\frac{d^{\vee}(h,r+\delta\vec{e_{i}})-d^{\vee}(h,r)}{\delta}=\\
 & =\partial_{u+i}d^{\vee}(h,r)\end{align*}
where $\vec{e_{i}}=(0,\ptind^{i-1},0,1,0,\dots,0)\in\R^{n}$. This
equality proves that $d\cdot\partial_{i}$ is a figure and hence that
the operator $\partial_{i}$ is smooth.
\item We can proceed in an analogous way (but here we have to use the derivation
under the integral sign) to prove that the integral operator: \begin{align*}
i:\ \Cc^{\infty}(\R^{2},\R) & \freccia\Cc^{\infty}(\R,\R)\\
u & \longmapsto{\displaystyle \int_{a}^{b}u(-,s)\diff{s}}\end{align*}
 is smooth.
\end{enumerate}
\item Because of cartesian closedness set-theoretical operations like the
following are examples of $\Cn$ arrows (see e.g. \citet{Ad-He-St}): 

\begin{itemize}
\item composition: \[
(f,g)\in B^{A}\times C^{B}\;\;\mapsto\;\; g\circ f\in C^{A}\]
 
\item evaluation: \[
(f,x)\in Y^{X}\times X\;\;\mapsto\;\; f(x)\in Y\]
 
\item insertion: \[
x\in X\;\;\mapsto\;\;(x,-)\in(X\times Y)^{Y}\]

\end{itemize}
\item Using the smoothness of the previous set-theoretical operations and
the smoothness of the derivation and integral operators, we can easily
prove that the classical operator of the calculus of variations is
smooth \[
\mathcal{I}(u)(t):=\int_{a}^{b}F[u(t,s),\partial_{2}u(t,s),s]\diff{s}\]
 \[
\mathcal{I}:\Cc^{\infty}(\R^{2},\R^{k})\freccia\Cc^{\infty}(\R,\R),\]
 where the function $F:\R^{k}\times\R^{k}\times\R\freccia\R$ is smooth.
\item Inversion between smooth manifolds modelled on Banach spaces\[
(-)^{-1}:f\in\text{Diff}(N,M)\;\;\mapsto\;\; f^{-1}\in\text{Diff}(M,N)\]
is a smooth mapping, where $\text{Diff}(M,N)$ is the subspace of
$N^{M}=\CInfty(\bar{M},\bar{N})$ given by the diffeomorphisms between
$M$ and $N$.\\
So $(\text{Diff}(M,M),\circ)$ is a (generalized) Lie group. To
prove that $(-)^{-1}$ is smooth let us consider a figure $d\In{U}\text{Diff}(N,M)$,
then, using cartesian closedness, the map $f:=(d\cdot i)^{\vee}:U\times N\freccia M$,
where $i:\text{Diff}(N,M)\hookrightarrow M^{N}$ is the inclusion,
is an ordinary smooth function between Banach manifolds. We have to
prove that $g:=[d\cdot(-)^{-1}\cdot j]^{\vee}:U\times M\freccia N$
is smooth, where $j:\text{Diff}(M,N)\hookrightarrow N^{M}$ is the
inclusion. But $f[u,g(u,m)]=m$ and $\text{{\bf D}}_{2}f(u,n)=\text{{\bf D}}[d(u)](n)$
hence the conclusion follows from the implicit function theorem because
$d(u)\in\text{Diff}(N,M)$.
\item Since the category $\Cn$ is complete, we can also have $\Cn$ spaces
with singular points like e.g. the equalizer%
\footnote{See the Appendix \ref{app:someNotionsOfCategoryTheory} for the notion
of equalizer.%
} $\{x\in X\,|\, f(x)=g(x)\}$. In this way, any algebraic curve is
a $\CInfty$ separated space too.
\item Another type of space with singular points is the following. Let $\varphi\in\Cc^{n}(\R^{k},\R^{m})$
and consider the subspace $([0,1]^{k}\prec\R^{k})$, then $(\varphi([0,1]^{k})\prec\R^{m})\in\Cn$
is a deformation in $\R^{m}$ of the hypercube $[0,1]^{k}$. 
\item Let $C$ be a continuum body, $I$ the interval for time, and ${\cal E}$
the 3-di\-men\-sional Euclidean space. We can define on $C$ a natural
structure of $\mathbf{\mathcal{C}}^{\infty}$-space. In fact, for
any point $p\in C$ let $p_{r}(t)\in{\cal E}$ be the position of
$p$ at time $t$ in the frame of reference $r$; we define figures
of type $U$ on $C$ ($U\in\ORn$) the functions $d:U\freccia C$
for which the following application \begin{align*}
\tilde{d}:U\times I & \freccia\mathcal{E}\\
(u,t) & \longmapsto d(u)_{r}(t)\end{align*}
is smooth. For example if $U=\R$ then we can think of $d:\R\freccia C$
as a curve traced on the body and parametrized by $u\in\R$. Hence
we are requiring that the position $d(u)_{r}(t)$ of the particle
$d(u)\in C$ in the frame of reference $r$ varies smoothly with the
parameter $u$ and the time $t$. This is a generalization of the
continuity of motion of any point of the body (take $d$ constant).
This smooth (that is diffeological) space will be separated, as an
object of $\mathbf{\mathcal{C}}^{\infty}$, if different points of
the body cannot have the same motion: \[
p_{r}(-)=q_{r}(-)\then p=q\qquad\forall p,q\in C.\]
The configuration space of $C$ can be viewed (see \citet{Wan}) as
a space of type \[
M:=\sum_{t\in I}M_{t}\ee{where}M_{t}\subseteq\mathcal{E}^{C}\]
and so, for the categorical properties of $\mathbf{\mathcal{C}}^{\infty}$
the spaces $\mathcal{E}^{C}$, $M_{t}$ (no matter how we choose these
subspaces $M_{t}$) and $M$ are always objects of $\mathbf{\mathcal{C}}^{\infty}$
as well. With this structure the motion of $C$ in the frame $r$:
\begin{align*}
\mu_{r}:\  & C\times I\freccia\mathcal{E}\\
{} & \h{.055}(p,t)\longmapsto p_{r}(t)\end{align*}
is a smooth map. Note that to obtain these results we need neither
$M_{t}$ nor $C$ to be manifolds, but only the possibility to associate
to any point $p$ of $C$ a motion $p_{r}(-):I\freccia\mathcal{E}$.
If we had the possibility to develop a differential geometry for these
spaces too we would have the possibility to obtain many results of
continuum mechanics for bodies which cannot be naturally represented
using a manifold or having an infinite-dimensional configuration space.
Moreover in the next chapter we will see how to extend any $\mathbf{\mathcal{C}}^{\infty}$
space with infinitesimal points, so that we can also consider infinitesimal
sub-bodies of $C$. 
\end{enumerate}

\chapter{\label{cha:Extending-smooth-spaces}Extending smooth spaces with
infinitesimals}

\section{Introduction}

\noindent The main aim of this chapter is to extend any $\CInfty$
space and any $\CInfty$ function by means of our {}``infinitesimal
points''. First of all, we will extend to a generic space $X\in\CInfty$
the notion of nilpotent path and of little-oh polynomial. The sets
of these paths will be denoted by $\mathcal{N}_{X}$ and $X_{o}[t]$
respectively%
\footnote{\noindent See Definition \ref{def:NilpotentFunctions} and Definition
\ref{def:LittleOhPolynomials} for the case $X=\R$.%
}. Afterward, we shall use the observables $\varphi$ of the space
$X$ to generalize the equivalence relation $\sim$ (i.e. the equality
in $\ER$, see Definition \ref{def:equalityInFermatReals}) using
the following idea \[
\varphi(x_{t})=\varphi(y_{t})+\text{o}(t)\ee{with}\varphi\InUp{UK}X.\]
Using this equivalence relation we will define $\ext{X}:=X_{o}[t]/\sim$,
which will be the generalization of the Definition $\ER:=\R_{o}[t]/\sim$.
Following this idea, the main problem is to understand how to relate
the little-oh polynomials $x$, $y$ with the domain $U$ of $\varphi$.
The second problem is that with this definition, $\ext{X}$ is a set
only, without any kind of structure. Indeed, we will tackle the problem
to define a meaningful category $\ECInfty$ and a suitable structure
on $\ext{X}$ so that $\ext{X}\in\ECInfty$. In the subsequent sections
we will also prove some results that will permit us to prove that
the extension functor $\ext{(-)}:\CInfty\freccia\ECInfty$ preserves
the product of manifolds, i.e. \[
\ext{(M\times N)}\simeq\ext{M}\times\ext{N}\]
for $M,N$ manifolds. The fact that this useful theorem is not proved
for generic $\CInfty$ spaces is due to the fact that the topology
on a product between $\CInfty$ spaces is generally stronger than
the product topology (see \eqref{eq:productTopologyAndProductSpaces},
but recall the final considerations of Section \ref{sec:DefinitionOfCartesianClosure}).

\subsection{Nilpotent paths}

\noindent If $X$ is a $\CInfty$ space, then using the topology $\Top{X}$
we can define the set $\Cc_{0}(X)$ of all the maps $x:\R_{\ge0}\freccia X$
which are continuous at the origin $t=0$. We want to simplify the
notations avoiding the use of germs of continuous functions as equivalent
classes (see \citet{Bou}), but, at the same time, we will keep attention
to consider only local properties $\mathcal{P}(x)$ when we will treat
paths $x\in\Cc_{0}(X)$ continuous at the origin, i.e. we will always
verify that\begin{equation}
\left(x,y\in\Cc_{0}(X)\quad\text{and}\quad x|_{[0,\varepsilon)}=y|_{[0,\varepsilon)}\quad\text{and}\quad\mathcal{P}(x)\right)\then\mathcal{P}(y).\label{eq:localPropertiesForC0Functions}\end{equation}
Following this constraint, it is not important how we extend%
\footnote{To be really rigorous, one has to fix, once and for all, a function
$E=E_{X}^{U,\eps}$ to perform such an extension, but taking into
consideration the fact that the whole construction does not depend
on this extension function. This function is defined on the set $\Cc_{0}([0,\eps),U)$
of function $x:[0,\eps)\freccia U$ continuous at $t=0^{+}$ and with
values in the subset $U\subseteq X$, i.e it is of the type $E:\Cc_{0}([0,\eps),U)\freccia\Cc_{0}(X)$,
and has the property $E(x)|_{[0,\eps)}=x$.%
} to the whole $\R_{\ge0}$ a locally defined function $x:[0,\eps)\freccia X$.

Because any $\CInfty$ function $f:X\freccia Y$ is continuous with
respect to the topologies $\Top{X}$ and $\Top{Y}$, we have that
$f\circ x\in\Cc_{0}(Y)$ if $x\in\Cc_{0}(X)$. More locally, if $U$
is open in $X$ and $x(0)\in U$, then on the subspace $(U\prec X)$
we have the induced topology and from this it follows that $\varphi\circ x\in\Cc_{0}(K)$
if $\varphi\InUp{UK}X$ is an observable of the space $X$. Let us
note explicitly that this is a local property, and hence, on the one
hand, with the notation $\varphi\circ x$ we have to mean a function
$\varphi\circ x:\R_{\ge0}\freccia X$ (because $\varphi\circ x$ is
an element of $\Cc_{0}(X)$). On the other hand, for this function
the only important property is that \[
\exists\,\eps>0\pti[0,\eps)\subseteq\left\{ t\in\R_{\ge0}\,|\, t\in\text{dom}(x)\quad\text{and}\quad x(t)\in\text{dom}(\varphi)\right\} ,\]
i.e. that the set of $t\in\R_{\ge0}$ for which the composition $\varphi(x(t))$
is defined, contains a right neighborhood of the origin.

As many other concepts we will introduce in this chapter, the notion
of nilpotent map is defined by means of the composition with a generic
observable and by a suitable logical implication to relate the starting
value $x(0)$ of a given path $x\in\Cc_{0}(X)$ with the domain of
the observable%
\footnote{Recall that, as usual, we will also use the notation $x_{t}$ for
the evaluation of $x\in Cc_{0}(X)$ at $t\in\text{dom}(x)$ and that
our little-oh functions (always for $t\to0^{+}$) are always continuous
at the origin (see Remark \ref{rem:WithLittleOhWeMean}).%
}.
\begin{defn}
\label{def:nilpotentForCnSpaces}Let $X$ be a $\CInfty$ space and
let $x\in\Cc_{0}(X)$ a path continuous at the origin, then we say
that $x$ \emph{is nilpotent (rel. $X$)} iff for every zone $UK$
of $X$ and every observable $\varphi\InUp{UK}X$ we have that the
following implication is true\[
x(0)\in U\then\exists k\in\N:\norm{\varphi(x_{t})-\varphi(x_{0})}^{k}=\text{{\rm o}}(t).\]
 Moreover we define \[
\No{X}:=\No{}(X):=\{x\in\Cc_{0}(X)\,\,|\,\, x\text{{\rm \ is nilpotent}}\}.\]

\end{defn}
A direct verification proves that the property of a path to be nilpotent
is a local property. Moreover, we will prove later that this definition
generalizes the particular notion expressed in Definition \ref{def:NilpotentFunctions}.

Because every $f\in\CInfty(X,Y)$ preserves the observables (see property
\eqref{thm: fPreservesObs}), if $x\in\No{X}$ then $f\circ x\in\No{Y}$,
that is $\CInfty$ functions preserve nilpotent maps too. In case
of a manifold $M$ (identified with its embedding $M=\bar{M}\in\CInfty$)
we can state the property of being nilpotent with an existential quantifier
instead of an implication
\begin{thm}
\label{thm:NilpotentForManifolds}Let $M$ be a $\Cc^{\infty}$ manifold
and let us consider a map $x:\R_{\ge0}\freccia|M|$, then $x$ is
nilpotent iff we can find a chart $(U,\varphi)$ on $x_{0}$ such
that $\norm{\varphi(x_{t})-\varphi(x_{0})}^{k}=\text{{\rm o}}(t)$
for some $k\in\N$.
\end{thm}
\textbf{Proof:} If we start from the hypothesis $x\in\No{M}$, then
it suffices to take any chart on $x_{0}$ and to use the property
that charts are observables of $M$ to get the conclusion formulated
in the statement.

To prove the opposite implication, let us take an observable $\psi\InUp{VK}\bar{M}$,
where $K$ is open in $\R^{p}$ and with $x_{0}\in V$. Recalling
\eqref{eq:observablesInManifolds} we get that $\psi\in\ManInfty(V,K)$,
i.e. $\psi$ is an ordinary $\Cc^{\infty}$ function. The idea is
to use the equality\[
\forall^{0}t\ge0\pti\psi(x_{t})=\psi\left[\varphi^{-1}(\varphi(x_{t}))\right],\]
which is locally true%
\footnote{Recall the definition of $\forall^{0}t\ge0$ given in Section \ref{sec:OrderRelation}.%
}, and the Lipschitz property of $\psi\circ\varphi^{-1}$. Diagrammatically,
in the category $\ManInfty$ of smooth manifolds, our situation is
the following\[
\xymatrix{\xyR{40pt}\xyC{50pt}U\cap V\ar[r]_{\sim}^{\varphi|_{U\cap V}}\ar[d]^{\ \psi|_{U\cap V}\ } & \varphi(U\cap V)\\
K}
\]
Therefore\[
\gamma:=\left(\varphi|_{U\cap V}\right)^{-1}\cdot\psi|_{U\cap V}\in\ManInfty(\varphi(U\cap V),K),\]
and hence $\gamma$ is locally Lipschitz with respect to some constant
$C>0$. But $x\in\Cc_{0}(\bar{M})$ and $U\cap V\in\Top{M}$, hence\[
\forall^{0}t\ge0\pti x_{t}\in U\cap V,\]
and we can write\begin{align*}
\norm{\psi(x_{t})-\psi(x_{0})}^{k} & =\norm{\gamma\left[\varphi(x_{t})\right]-\gamma\left[\varphi(x_{0})\right]}^{k}\le\\
 & \le C^{k}\cdot\norm{\varphi(x_{t})-\varphi(x_{0})}^{k}=\text{o}(t)\end{align*}
\qedWithFinalEq

This will be a typical idea in several definitions of the present
work: working with generic $\CInfty$ spaces we do not have the possibility
to consider charts on every point, so we require a condition \emph{for
every} \emph{observable} that \emph{potentially} (i.e. by means of
a logical implication) contains the starting point of a given path.
We have already used this idea in the Definition \ref{def: xIdentifiedy}
of {}``$x$ is identified with $y$'', i.e. of the relation $x\asymp y$.
Usually, in case we have charts, like in the previous Theorem \ref{thm:NilpotentForManifolds}
we will be able to transform in an equivalent statement this type
of implications using an existential quantifier. This theorem also
proves that the previous Definition \ref{def:NilpotentFunctions}
is a generalization of the old Definition \ref{def:nilpotentForCnSpaces}.

Finally we consider the relations between the product of two manifolds
$M$, $N$ and nilpotent paths in the following
\begin{thm}
\noindent \label{thm: NilpotentAndProduct} Let $M,N$ be smooth manifolds
and $x:\R_{\ge0}\freccia|M|$, $y:\R_{\ge0}\freccia|N|$ be two maps,
then \[
x\in\No{\bar{M}}\e{and}y\in\No{\bar{N}}\quad\iff\quad(x,y)\in\No{\bar{M}\times\bar{N}},\]
where we set $(x,y)_{t}:=(x_{t},y_{t})$.
\end{thm}
\noindent \textbf{Proof:} 

\noindent $\Leftarrow$ : If $(x,y)\in\No{\bar{M}\times\bar{N}}=\No{\overline{M\times N}}$,
by the previous Theorem \ref{thm:NilpotentForManifolds} we get the
existence of two charts $(U,\varphi)$ of $M$ and $(V,\psi)$ of
$N$ with $x_{0}\in U$ and $y_{0}\in V$ and such that%
\footnote{If it will be clear from the context, we will sometimes omit the parenthesis
in compositions like $fg(x)=f(g(x))$.%
}\[
\exists\, k\in\N\pti\norm{(\varphi x_{t},\psi y_{t})-(\varphi x_{0}t,\psi y_{0})}^{k}=\text{o}(t).\]
Therefore, we also have $\left\{ \norm{\varphi x_{t}-\varphi x_{0}}+\norm{\psi y_{t}-\psi y_{0}}\right\} ^{k}=\lo(t)$.

\noindent But $\norm{\phi x_{t}-\phi x_{0}}\le\norm{\phi x_{t}-\phi x_{0}}+\norm{\psi y_{t}-\psi y_{0}}$
and hence\[
\norm{\phi x_{t}-\phi x_{0}}^{k}\le\left\{ \norm{\phi x_{t}-\phi x_{0}}+\norm{\psi y_{t}-\psi y_{0}}\right\} ^{k}=\lo(t).\]
Therefore also $\norm{\phi x_{t}-\phi x_{0}}^{k}=\lo(t)$, that is
$x\in\No{M}$. Analogously we can proceed for $y$.

\medskip{}

\noindent $\Rightarrow$ : From the hypotheses $x\in\No{\bar{M}}$
and Theorem \ref{thm:NilpotentForManifolds} we get a chart $(U,\phi)$
of $M$ on $x_{0}$ and a $k\in\N$ such that\begin{equation}
\norm{\phi x_{t}-\phi x_{0}}^{k}=\lo(t).\label{eq:chartFrom_x}\end{equation}
Analogously, from $y\in\No{\bar{N}}$ we obtain a chart $(V,\psi)$
on $y_{0}$ and a $k'\in\N$ such that\begin{equation}
\norm{\psi y_{t}-\psi y_{0}}^{k'}=\lo(t).\label{eq:chartFrom_y}\end{equation}
We can suppose $k=k'$. Therefore $(U\times V,\phi\times\psi)$ is
a chart of $M\times N$ on $(x_{0},y_{0})$. Let us try to compute
the term\begin{multline}
\norm{(\phi\times\psi)(x_{t},y_{t})-(\phi\times\psi)(x_{0},y_{0})}^{k}=\norm{(\phi x_{t}-\phi x_{0},\psi y_{t}-\psi y_{0})}^{k}=\\
=\sum_{i=0}^{k}\binom{k}{i}\norm{\phi x_{t}-\phi x_{0}}^{i}\cdot\norm{\psi y_{t}-\psi y_{0}}^{k-i}.\label{eq:normOfProduct}\end{multline}
But\begin{multline*}
\frac{\norm{\phi x_{t}-\phi x_{0}}^{i}\cdot\norm{\psi y_{t}-\psi y_{0}}^{k-i}}{t}=\\
=\frac{\left\{ \norm{\phi x_{t}-\phi x_{0}}^{k}\right\} ^{\frac{i}{k}}}{t^{\frac{i}{k}}}\cdot\frac{\left\{ \norm{\psi y_{t}-\psi y_{0}}^{k}\right\} ^{\frac{k-i}{k}}}{t^{\frac{k-i}{k}}}=\\
=\left\{ \frac{\norm{\phi x_{t}-\phi x_{0}}^{k}}{t}\right\} ^{\frac{i}{k}}\cdot\left\{ \frac{\norm{\psi y_{t}-\psi y_{0}}^{k}}{t}\right\} ^{\frac{k-i}{k}}.\end{multline*}
Each factor of this product goes to zero for $t\to0^{+}$ because
of \eqref{eq:chartFrom_x} and \eqref{eq:chartFrom_y}. Hence also
\eqref{eq:normOfProduct} goes to zero and this, because of Theorem
\ref{thm:NilpotentForManifolds}, proves that\[
(x,y)\in\No{\overline{M\times N}}=\No{\bar{M}\times\bar{N}}.\]
\qedWithFinalEq

\subsection{Little-oh polynomials in $\CInfty$}

We can proceed in a similar way with respect to the generalization
of the notion of little-oh polynomial: at first we will define what
is a little-oh polynomial in $\R^{d}$, and secondly we will generalize
this notion to a generic space $X\in\CInfty$ using observables.
\begin{defn}
\label{def:LittleOhPolynomialsRd} We say that $x$ \emph{is a little-oh
polynomial in $\R^{d}$}, and we write $x\in\R_{o}^{d}[t]$, iff 
\begin{enumerate}
\item $x:\R_{\ge0}\freccia\R^{d}$
\item We can write\[
x_{t}=r+\sum\limits _{i=1}^{k}\alpha_{i}\cdot t^{a_{i}}+o(t)\quad\text{as}\quad t\to0^{+}\]
 for suitable\[
k\in\N\]
\[
r,\alpha_{1},\dots,\alpha_{k}\in\R^{d}\]
\[
a_{1},\dots,a_{k}\in\R_{\ge0}.\]

\end{enumerate}
Now let $X\in\CInfty$ and $x\in\Cc_{0}(X)$, then we say that $x$
\emph{is a little-oh polynomial (of $X$) }iff for every zone $UK$
of $X$, with $K\subseteq\R^{\sf k}$, and every observable $\phi\InUp{UK}X$
we have\[
x_{0}\in U\then\phi\circ x\in\R_{o}^{\sf k}[t].\]
Moreover\[
X_{o}[t]:=X_{o}:=\left\{ x\in\Cc_{0}(X)\,|\, x\text{ is a little-oh polynomial of }X\right\} .\]

\end{defn}
\noindent Let us note that for $d=1$ we have exactly the old Definition
\ref{def:LittleOhPolynomials}. A direct verification proves that
being a little-oh polynomial is a local property. Moreover, we will
prove later that the two parts of this definition (i.e. that of $X_{o}[t]$
and that of $\R_{o}^{d}[t]$ are equivalent if $X=\R^{d}$).

Now we have to prove the analogous for little-oh polynomials of the
previous results stated for nilpotent paths. Once again, because every
$f\in\CInfty(X,Y)$ preserves the observables, we have that $\CInfty$
functions preserve little-oh polynomials too\[
x\in X_{0}[t]\then f\circ x\in Y_{o}[t].\]
The other results we want to prove relate the notion of little-oh
polynomial with that of manifold: at first, as usual, we want to reformulate
the Definition \ref{def:LittleOhPolynomialsRd} for manifolds; secondly
we want to make clear the relationships between little-oh polynomials
and the product of manifolds. For these results we need the following
Lemmas.
\begin{lem}
\label{lem:productOf_Rn_andLittleOhPolynomials}Let $x:\R_{\ge0}\freccia\R^{m}$
and $y:\R_{\ge0}\freccia\R^{n}$ be two maps, then\[
x\in\R_{o}^{m}[t]\e{and}y\in\R_{o}^{n}[t]\quad\iff\quad(x,y)\in\R_{o}^{m+n}[t].\]

\end{lem}
\noindent \textbf{Proof:}

\noindent $\Rightarrow$ : Let us fix the notations for the little-oh
polynomials $x$ and $y$:\[
x_{t}=r+\sum\limits _{i=1}^{K}\alpha_{i}\cdot t^{a_{i}}+o_{1}(t)\]
\[
y_{t}=s+\sum\limits _{j=1}^{N}\beta_{j}\cdot t^{b_{j}}+o_{2}(t),\]
where $r,\alpha_{1},\ldots,\alpha_{K}\in\R^{m}$ and $s,\beta_{1},\ldots,\beta_{N}\in\R^{n}$.
Define $u:=(r,s)\in\R^{m+n}$, $\gamma_{i}:=(\alpha_{i},\underbar{0})\in\R^{m+n}$,
$\gamma_{j+K}:=(\underbar{0},\beta_{j})\in\R^{m+n}$, $c_{i}:=a_{i}$
and $c_{j+K}:=b_{j}$, then\begin{align*}
(x_{t},y_{t}) & =\left(r+\sum\limits _{i=1}^{K}\alpha_{i}\cdot t^{a_{i}}+o_{1}(t),s+\sum\limits _{j=1}^{N}\beta_{j}\cdot t^{b_{j}}+o_{2}(t)\right)=\\
 & =(r,s)+\sum_{i=1}^{K}(\alpha_{i},\underbar{0})\cdot t^{a_{i}}+(o_{1}(t),\underbar{0})+\\
 & \phantom{=}+\sum_{j=1}^{N}(\underbar{0},\beta_{j})\cdot t^{b_{j}}+(\underbar{0},o_{2}(t))=\\
 & =u+\sum_{i=1}^{K}\gamma_{i}\cdot t^{c_{i}}+\sum_{i=K+1}^{K+N}\gamma_{i}\cdot t^{c_{i}}+(o_{1}(t),o_{2}(t)),\end{align*}
and this proves the conclusion because $(o_{1}(t),o_{2}(t))=o(t)$.

\medskip{}

\noindent $\Leftarrow$ : By hypotheses we can write\[
(x_{t},y_{t})=u+\sum_{k=1}^{H}\gamma_{k}\cdot t^{c_{k}}+o(t).\]
We only have to reverse the previous ideas defining:\[
r:=(u_{1},\ldots,u_{m})\qquad s:=(u_{m+1},\ldots,u_{m+n})\]
\[
\alpha_{k}:=(\gamma_{k}^{1},\ldots,\gamma_{k}^{m})\qquad\beta_{k}:=(\gamma_{k}^{m+1},\ldots,\gamma_{k}^{m+n})\]
\[
o_{1}(t):=(o^{1}(t),\ldots,o^{m}(t))\qquad o_{2}(t):=(o^{m+1}(t),\ldots,o^{m+n}(t))\]
\[
a_{i}:=c_{i}\qquad b_{k}:=c_{k}\]
where we have used the notations\[
\gamma_{k}=(\gamma_{k}^{1},\ldots,\gamma_{k}^{m+n})\]
\[
o(t)=(o^{1}(t),\ldots,o^{m+n}(t))\]
for the components. Then\begin{align*}
(x_{t},y_{t}) & =(r,\underbar{0})+\sum_{k=1}^{H}(\alpha_{k},\underbar{0})\cdot t^{c_{k}}+(o_{1}(t),\underbar{0})+\\
 & \phantom{=}(\underbar{0},s)+\sum_{k=1}^{H}(\underbar{0},\beta_{k})\cdot t^{c_{k}}+(\underbar{0},o_{2}(t))=\\
 & =\left(r+\sum_{i=1}^{H}\alpha_{i}\cdot t^{c_{i}}+o_{1}(t),s+\sum_{j=1}^{H}\beta_{j}\cdot t^{c_{j}}+o_{2}(t)\right),\end{align*}
and hence the conclusion follows.\qedNoNewLine

\noindent From this lemma, if $x\in\R_{o}^{d}[t]$, then each component
is a 1-dimensional little-oh polynomial $x_{i}\in\R_{o}[t]$ for $i=1,\ldots,d$.
But we know (see Section \ref{sec:firstPropertiesOfLittle-ohPolynomials})
that each one of these polynomial is nilpotent, i.e. $x_{i}\in\mathcal{N}$.
Therefore, from Theorem \ref{thm: NilpotentAndProduct} it follows
that $x\in\No{\R^{d}}$, i.e.\[
\R_{o}^{d}[t]\subseteq\No{\R^{d}};\]
from this it also follows that $X_{o}[t]\subseteq\No{X}$.
\begin{lem}
\noindent \label{lem:smoothFunctionPreserversLittleOhPolynomials}Let
$x\in\R_{o}^{d}[t]$ and $f\in\Cc^{\infty}(A,\R^{p})$, with $A$
open in $\R^{d}$ and such that, locally, the path $x$ has values
in $A$:\[
\forall^{0}t\ge0\pti x_{t}\in A.\]
Then $f\circ x\in\R_{o}^{p}[t]$.
\end{lem}
\noindent \textbf{Proof:} Let us fix some notations:\[
x_{t}=r+\sum_{i=1}^{k}\alpha_{i}\cdot t^{a_{i}}+w(t)\quad\text{with}\quad w(t)=o(t)\]
 \[
h(t):=x(t)-x(0)\quad\forall t\in\R_{\ge0}\]
 hence $x_{t}=x(0)+h_{t}=r+h_{t}$. The function $t\mapsto h(t)=\sum_{i=1}^{k}\alpha_{i}\cdot t^{a_{i}}+w(t)$
belongs to $\R_{o}^{d}[t]\subseteq\mathcal{N}_{\R^{d}}$, so we can
write $\norm{h_{t}}^{N}=o(t)$ for some $N\in\N$ if we take as observable
of $\R^{d}$ the identity. From Taylor's's formula we have\begin{equation}
f(x_{t})=f(r+h_{t})=f(r)+\sum_{\substack{i\in\N^{d}\\
|i|\le N}
}\frac{\partial^{i}f}{\partial x^{i}}(r)\cdot\frac{h_{t}^{i}}{i!}+o\left(\norm{h_{t}}^{N}\right).\label{eq:TaylorLittleOhPolynomialsClosedSmooth}\end{equation}
 But\[
\frac{o\left(\norm{h_{t}}^{N}\right)}{|t|}=\frac{o\left(\norm{h_{t}}^{N}\right)}{\norm{h_{t}^{N}}}\cdot\frac{\norm{h_{t}^{N}}}{|t|}\to0\]
hence $o\left(\norm{h_{t}}^{N}\right)=o(t)\in\R_{o}^{p}[t]$. Now
we have to note that, for a multi-index $i\in\N^{d}$, it results
that $h_{t}^{i}=h_{1}^{i_{1}}(t)\cdot\ldots\cdot h_{d}^{i_{d}}(t)\in\R_{o}[t]$
because, from the previous Lemma \ref{lem:productOf_Rn_andLittleOhPolynomials},
each function $h_{j}(t)\in\R_{o}[t]$ and because $\R_{o}[t]$ is
an algebra. Moreover, if $\beta\in\R^{p}$ and $h\in\R_{o}[t]$, then
$\beta\cdot h\in\R_{o}^{p}[t]$, so each addend $\frac{\partial^{i}f}{\partial x^{i}}(r)\cdot\frac{h_{t}^{i}}{i!}$
is a little-oh polynomial of $\R_{o}^{p}[t]$. From \eqref{eq:TaylorLittleOhPolynomialsClosedSmooth}
and the closure of little-oh polynomials $\R_{o}^{p}[t]$ with respect
to linear operations, the conclusion $f\circ x\in\R_{o}^{p}[t]$ follows.\qedNoNewLine

Using these lemmas we can prove the above cited results about little-oh
polynomials in manifolds.
\begin{thm}
\label{thm:manifoldsAndLittleOhPolynomials}If $M$ is a $\Cc^{\infty}$
manifold and $x:\R_{\ge0}\freccia|M|$ is a map, then we have that
$x\in\bar{M}_{o}[t]$ if and only if there exists a chart $(U,\phi)$
of $M$ such that:
\begin{enumerate}
\item \label{enu:firstCondManifoldAndNilpotent}$x(0)\in U$
\item \label{enu:secondCondManifoldAndNilpotent}$\phi\circ x\in\R_{o}^{d}[t]$,
where $d:=\dim(M)$.
\end{enumerate}
\end{thm}
\noindent \textbf{Proof:} To prove that the hypotheses $x\in\bar{M}_{o}[t]$
implies conditions \ref{enu:firstCondManifoldAndNilpotent}. and \ref{enu:secondCondManifoldAndNilpotent}.
it suffices to take any chart on $x_{0}$ and to use the property
that charts are observables of $\bar{M}$.

For the opposite implication we start considering that, by the Definition
\ref{def:LittleOhPolynomialsRd} we have that $\phi\circ x\in\R_{o}^{d}[t]$
is continuous at $t=0^{+}$, and hence also $x$ is continuous at
$t=0^{+}$, i.e $x\in\Cc_{0}(M)$. Now, take a generic observable
$\psi\InUp{VK}\bar{M}$, where $K$ is open in $\R^{p}$ and $x_{0}\in V$.
We have\[
\xymatrix{\xyC{65pt}\R^{d}\supseteq\phi(U\cap V)\ar[r]\sp(0.6){\left(\phi|_{U\cap V}\right)^{-1}}\sb(0.6){\sim} & U\cap V\ar[r]^{\psi|_{U\cap V}} & K\subseteq\R^{p}}
,\]
and hence \begin{equation}
\left(\phi|_{U\cap V}\right)^{-1}\cdot\psi|_{U\cap V}\in\Cc^{\infty}(\phi(U\cap V,K).\label{eq:inverseOfChartAndObservableIsCn}\end{equation}
From $x_{0}\in U\cap V$ and from the continuity of the path $x$
at $t=0^{+}$ we get\[
\forall^{0}t\ge0\pti x_{t}\in U\cap V.\]
From this, from \eqref{eq:inverseOfChartAndObservableIsCn}, from
the hypotheses $\phi\circ x\in\R_{o}^{d}[t]$ and from Lemma \ref{lem:smoothFunctionPreserversLittleOhPolynomials}
the conclusion $\left(\phi|_{U\cap V}\right)^{-1}\cdot\psi|_{U\cap V}\circ\phi\circ x=\psi\circ x\in\R_{o}^{p}[t]$
follows.\qedNoNewLine
\begin{thm}
\label{thm:littleOhPolynomialsAndProductOfManifolds}Let $M$, $N$
be $\Cc^{\infty}$ manifolds and $x:\R_{\ge0}\freccia|M|$, $y:\R_{\ge0}\freccia|M|$
two maps. Then\[
x\in\bar{M}_{o}[t]\e{and}y\in\bar{N}_{o}[t]\quad\iff\quad(x,y)\in\left(\bar{M}\times\bar{N}\right)_{o}[t].\]

\end{thm}
\textbf{Proof:} The proof is an almost purely logical consequence
of Theorem \ref{thm:manifoldsAndLittleOhPolynomials} and of Lemma
\ref{lem:productOf_Rn_andLittleOhPolynomials}.

\noindent $\mathbb{\Leftarrow}$ : By hypotheses $(x,y)\in\left(\bar{M}\times\bar{N}\right)_{o}[t]=\left(\overline{M\times N}\right)_{o}[t]$.
Because $M\times N$ is a manifold, from Theorem \ref{thm:manifoldsAndLittleOhPolynomials}
we get the existence of charts $(U,\phi)$ of $M$ and $(V,\psi)$
of $N$ on $x_{0}$ and $y_{0}$ resp. such that\[
\left(\phi\times\psi\right)\circ(x,y)=\left(\phi\circ x,\psi\circ y\right)\in\R_{o}^{m+n}[t].\]
Hence $\phi\circ x\in\R_{o}^{m}[t]$ and $\psi\circ y\in\R_{o}^{n}[t]$
from Lemma \ref{lem:productOf_Rn_andLittleOhPolynomials}.

\medskip{}

\noindent $\Rightarrow$ : Analogously, if $x\in\bar{M}_{o}[t]$ and
$y\in\bar{N}_{o}[t]$, then we can find charts as above, but with
$\phi\circ x\in\R_{o}^{m}[t]$ and $\psi\circ y\in\R_{o}^{n}[t]$.
Once again from Lemma \ref{lem:productOf_Rn_andLittleOhPolynomials}
we obtain\[
\left(\phi\circ x,\psi\circ y\right)=\left(\phi\times\psi\right)\circ(x,y)\in\R_{o}^{m+n}[t],\]
and hence the conclusion follows from Theorem \ref{thm:manifoldsAndLittleOhPolynomials}.\qedNoNewLine

\section{\label{sec: extensions}The Fermat extension of spaces and functions}

Considering the previous definitions of nilpotent and little-oh paths
and the Definition \ref{def: xIdentifiedy} it is now clear how to
generalize the definition of equality in $\ER$ (see Definition \ref{def:equalityInFermatReals})
to a generic $X\in\CInfty$:
\begin{defn}
\label{def:equalityInExt_X}Let $X$ be a $\CInfty$ space and let
$x$, $y\in X_{o}[t]$ be two little-oh polynomials, then we say that
\[
x\sim y\e{in}X\ee{or\ simply}x=y\e{in}\ext{X}\]
 iff for every zone $UK$ of $X$ and every observable $\varphi\InUp{UK}X$
we have
\begin{enumerate}
\item $x_{0}\in U\iff y_{0}\in U$
\item $x_{0}\in U\then\varphi(x_{t})=\varphi(y_{t})+\text{o}(t)$
\end{enumerate}
\end{defn}
Obviously we will write $\ext{X}:=X_{o}[t]/\!\sim$ and $\ext{f}(x):=f\circ x$
if $f\in\CInfty(X,Y)$ and $x\in\ext{X}$ and we will call them the
\emph{Fermat extension of $X$ and of $f$} respectively. As usual,
we will also define the \emph{standard part of }$x\in\ext{X}$ as
$\st{x}:=x(0)\in X$.

We prove the correctness of the definition of $\ext{f}$ in the following:
\begin{thm}
If $f\in\CInfty(X,Y)$ and $x=y$ in $\ext{X}$ then $\ext{f}(x)=\ext{f}(y)$
in $\ext{Y}$.
\end{thm}
\textbf{Proof:} Take a zone $VK$ in $Y$ and an observable $\psi\InUp{VK}Y$,
then from the continuity of $f$, we have $U:=f^{-1}(V)\in\Top{X}$.
We can thus apply hypothesis $x=y$ in $\ext{X}$ with the zone $UK$
and the observable $\varphi:=f|_{U}\cdot\psi\InUp{UK}X$. From this
the conclusion follows considering that $f\circ x,f\circ y\in Y_{o}[t]$
and $x_{0}\in U=f^{-1}(V)$ iff $f(x_{0})\in V$. \qedNoNewLine

Using the continuity of $\phi\circ x$ we can note that $x=y$ in
$\ext{X}$ implies that $x_{0}$ and $y_{0}$ are identified in $X$
(see Definition \ref{def: xIdentifiedy}) and thus using constant
maps $\hat{x}(t):=x$, for $x\in X$, we obtain an injection $\hat{{\scriptstyle (-)}}:|X|\freccia\ext{X}$
if the space $X$ is separated. Therefore, if $Y$ is separated too,
$\ext{f}$ is really an extension of $f$. Finally, note that the
application $\ext{(-)}$ preserves compositions and identities.

Using ideas very similar to the ones used above for similar theorems,
we can prove that if $X=M$ is a $\Cc^{\infty}$ manifold then we
have that $x=y$ in $\ext{M}$ iff there exists a chart $(U,\varphi)$
of $M$ such that 
\begin{enumerate}
\item $x_{0}$, $y_{0}\in U$ 
\item $\varphi(x_{t})=\varphi(y_{t})+\text{{\rm o}}(t)$. 
\end{enumerate}
\noindent Moreover the previous conditions do not depend on the chart
$(U,\varphi)$. In particular if $X=U$ is an open set in $\R^{k}$,
then $x=y$ in $\ext{U}$ is simply equivalent to the limit relation
$x(t)=y(t)+\text{o}(t)$ as $t\to0^{+}$; hence if $i:U\hookrightarrow\R^{k}$
is the inclusion map, it's easy to prove that its Fermat extension
$\ext{i}:\ext{U}\freccia\ER^{k}$ is injective. We will always identify
$\ext{U}$ with $\ext{i}(\ext{U})$, so we simply write $\ext{U}\subseteq\ER^{k}$.
According to this identification, if $U$ is open in $\R^{k}$, we
can also prove that\begin{equation}
\ext{U}=\{x\in\ext{\R^{k}}\,|\,\st{x}\in U\}.\label{eq:FermatExtOf_U}\end{equation}
This property says that the preliminary definition of $\ext{U}$ given
in Definition \ref{def:extensionOfSubsetsOfR} is equivalent to the
previous, more general, Definition \ref{def:equalityInExt_X} of extension.
Using the previous equivalent way to express the relation $\sim$
on manifolds, we see that $(x,y)=(x',y')$ in $\ext{(M\times N)}$
iff $x=x'$ in $\ext{M}$ and $y=y'$ in $\ext{N}$. From this conclusion
and from Theorem \ref{thm:littleOhPolynomialsAndProductOfManifolds}
we prove that the following applications \begin{equation}
\alpha_{\sss MN}:=\alpha\,\,:\,\,([x]_{\sim},[y]_{\sim})\in\ext{M}\times\ext{N}\longmapsto[(x,y)]_{\sim}\in\ext{(M\times N)}\label{eq:alpha}\end{equation}
 \begin{equation}
\beta_{\sss MN}:=\beta\,\,:\,\,[z]_{\sim}\in\ext{(M\times N)}\longmapsto([z\cdot p_{M}]_{\sim},[z\cdot p_{N}]_{\sim})\in\ext{M}\times\ext{N}\label{eq:beta}\end{equation}
 (for clarity we have used the notation with the equivalence classes)
are well-defined bijections with $\alpha^{-1}=\beta$ (obviously $p_{\sss M},p_{\sss N}$
are the projections). We will use the first one of them in the following
section with the temporary notation $\langle p,x\rangle:=\alpha(p,x)$,
hence $f\langle p,x\rangle=f(\alpha(p,x))$ for $f:\ext{(M\times N)}\freccia Y$.
This simplifies our notations but permits to avoid the identification
of $\ext{M}\times\ext{N}$ with $\ext{(M\times N)}$ until we will
have proved that $\alpha$ and $\beta$ are arrows of the category
$\ECInfty$.

\section{The category of Fermat spaces}

\noindent Up to now every $\ext{X}$ is a simple set only. Now we
want to use the general passage from a category of the types of figures
$\F$ to its cartesian closure $\bar{\F}$ so as to put on any $\ext{X}$
a useful structure of $\bar{\F}$ space. Our aim is to obtain in this
way a new cartesian closed category $\bar{\F}=:\ECInfty$, called
the \emph{category of Fermat} \emph{spaces}, and a functor $\ext{(-)}:\CInfty\freccia\ECInfty$,
called the \emph{Fermat functor}. Therefore we have to choose $\F$,
that is we have to understand what can be the types of figures of
$\ext{X}$. It may seem very natural to take $\ext{g}:\ext{U}\freccia\ext{V}$
as arrow in $\F$ if $g:U\freccia V$ is in $\ORInfty$ (in \citet{Gio}
we followed this way). The first problem in this idea is that, e.g.
\[
\ER\xfrecciad{\ext{f}}\ER\then\ext{f}(0)=f(0)\in\R,\]
hence there cannot exist a constant function of the type $\ext{f}$
to a non-standard value, and so we cannot satisfy the closure of $\F$
with respect to generic constant functions (see the hypotheses about
the types of figures $\F$ in Section \ref{sec:hypothesisForCartesianClosure}).
But we can make further considerations about this problem so as to
better motivate the choice of $\F$. The first one is that we surely
want to have the possibility to lift maps%
\footnote{\noindent I.e. to consider their adjoint function using cartesian
closedness.%
} as simple as the sum between Fermat reals: \[
s:(p,q)\in\ER\times\ER\freccia p+q\in\ER.\]
Therefore, we have to choose $\F$ so that the map $s^{\wedge}(p):q\in\ER\freccia p+q\in\ER$
is an arrow of $\ECInfty$. Note that this map is neither constant
nor of the type $\ext{f}$ because $s^{\wedge}(p)(0)=p$ and $p$
could be a non standard Fermat real.

\noindent The second consideration is about the map $\alpha$ defined
in \eqref{eq:alpha}: if we want $\alpha$ to be an arrow of $\ECInfty$,
then in the following situation we have to obtain a $\ECInfty$ arrow
\begin{gather*}
\ER\times\ER\xfrecciad{p\times1_{\ER}}\ER\times\ER\xfrecciad{\alpha}\ext{(\R\times\R)}\xfrecciad{\ext{g}}\ER\\
(r,s)\longmapsto(p,s)\longmapsto\langle p,s\rangle\longmapsto\ext{g}\langle p,s\rangle,\end{gather*}
where $p\in\ER$ and $g\in\Cc^{\infty}(\R^{2},\R)$. The idea we shall
follow is exactly to take as arrows of $\F$ all the maps that locally
are of the form $\delta(s)=\ext{g}\langle p,s\rangle$, where $p\in\ext{(\R^{\sf p})}$
works as a parameter of $\ext{g}\langle-,-\rangle$. Obviously, in
this way $\delta$ could also be a constant map to a non standard
value (take as $g$ a projection). Frequently one can find maps of
the form $\ext{g}\langle p,-\rangle$ in informal calculations in
physics or geometry. Actually, they simply are $\Cc^{\infty}$ maps
with some fixed parameter $p$, which could be an infinitesimal distance
(e.g. in the potential of the electric dipole, see below), an infinitesimal
coefficient associated to a metric (like, e.g., in Einstein's formula
\eqref{eq:EinsteinInfinitesimal}), or a side $l:=s(a,-)$ of an infinitesimal
surface $s:[a,b]\times[c,d]\freccia\ER$, where $[a,b],[c,d]\subseteq D_{k}$.

\noindent Note the importance of the map $\alpha$ to perform passages
like the following\[
M\times N\xfrecciad{f}Y\e{in}\CInfty\]
\[
\ext(M\times N)\xfrecciad{\ext{f}}\ext{Y}\e{in}\ECInfty\]
\[
\ext{M}\times\ext{N}\xfrecciad{\ext{f}}\ext{Y}\e{in}\ECInfty\text{ (identification via }\alpha\text{)}\]
\[
\ext{N}\xfrecciad{\ext{f}^{\wedge}}\ext{Y}^{\ext{M}}\e{using\ cartesian\ closedness.}\text{ }\]
This motivates the choice of arrows in $\F$, but there is a second
problem about the choice of the objects of the category $\mathcal{F}$.
Take a manifold $M$ and an arrow $t:D\freccia\ext{M}$ in $\ECInfty$.
Even if we have not still defined formally what is the meaning of
this {}``arrow'', we want to think $t$ as a tangent vector applied
either to a standard point $t(0)\in M$ or to a non standard one,
$t(0)\in\ext{M}\setminus M$. Roughly speaking, this is the case if
we can write $t(h)=\ext{g}\langle p,h\rangle$ for every $h\in D$
and for some $g$, $p$. If we want to obtain this equality it is
useful to have two properties: the first one is that the identity
map over $D$, i.e. $1_{D}$, is a figure of $D$, i.e. $1_{D}\In{D}D$.
In this way, from the property $t:D\freccia\ext{M}$ of being an arrow
of $\ECInfty$ we can deduce that $t$ is a figure of $\ext{M}$ of
the type $D$, i.e. $t\In{D}\ext{M}$. The second property we would
like to obtain is to have maps of the form $\ext{g}\langle p,-\rangle:D\freccia\ext{M}$
as figures of $\ext{M}$. Of course, we can thus say that necessarily
$t=\ext{g}\langle p,-\rangle$ for some $g\in\Cc^{\infty}(\R^{\sf p},\R)$
and $p\in\R^{\sf p}$. Therefore, to obtain these properties, it would
be useful to have $D$ as an object of $\F$. But $D$ is not the
extension of a standard subset of $\R$, thus what will be the objects
of $\F$? We will take generic subsets $S$ of $\ext{(\R^{\sf s})}$
with the topology $\Top{S}$ generated by $\mathcal{U}=\ext{U}\cap S$,
for $U$ open in $\R^{\sf s}$ (in this case we will say that the
open set $\mathcal{U}$ \emph{is defined by} $U$ \emph{in} $S$).
In other words $A\in\Top{S}$ if and only if\begin{equation}
A=\bigcup\left\{ \ext{U}\cap S\subseteq A\,|\, U\text{ is open in }\R^{\sf s}\right\} .\label{eq:topologyInSERInfty}\end{equation}

These are the motivations to introduce the category of the types of
figures $\F$ by means of the following 
\begin{defn}
\noindent \label{def:arrowsOfSERInfty}We call $\SERInfty$ the category
whose objects are topological spaces $(S,\Top{S})$, with $S\subseteq\ext{(\R^{\sf s})}$
for some $\sf s\in\N$ which depends on $S$, and with the previous
topology $\Top{S}$. In the following we will frequently use the simplified
notation $S$ instead of the complete $(S,\Top{S})$.

If $S\subseteq\ext{(\R^{\sf s})}$ and $T\subseteq\ext{(\R^{\sf t})}$
then we say that \[
S\xfrecciad{f}T\e{in}\SERInfty\]
 iff $f$ maps $S$ in $T$ and for every $s\in S$ we can write \begin{equation}
f(x)=\ext{g}\langle p,x\rangle\quad\forall x\in\ext{V}\cap S\label{eq:localFormForArrowsOfSERn}\end{equation}
 for some \begin{align*}
{} & V\text{ open in }\R^{\sf s}\text{ such that }s\in\ext{V}\\
{} & p\in\ext{U},\text{ where }U\text{ is open in }\R^{\sf p}\\
{} & g\in\Cc^{\infty}(U\times V,\R^{\sf t}).\end{align*}
Moreover we will consider on $\SERInfty$ the forgetful functor given
by the inclusion $|-|:\SERInfty\hookrightarrow\text{{\bf Set}}$,
i.e. $|(S,\Top{S})|:=S$. The category $\SERInfty$ will be called
\emph{the category of subsets of $\ext{\R^{\infty}}$} (but note that
here $\infty$ indicates the class of regularity of the functions
we are considering).\end{defn}
\begin{rem*}
\ 
\begin{enumerate}
\item In other words \emph{locally }a $\CInfty$ function $f:S\freccia T$
between two types of figures $S\subseteq\ext{(\R^{\sf s})}$ and $T\subseteq\ext{(\R^{\sf t})}$
is constructed in the following way:

\begin{enumerate}
\item start with an ordinary standard function $g\in\Cc^{\infty}(U\times V,\R^{\sf t})$,
with $U$ open in $\R^{\sf p}$ and $V$ open in $\R^{\sf s}$. The
space $\R^{\sf p}$ has to be thought as a space of parameters for
the function $g$;
\item consider its Fermat extension obtaining $\ext{g}:\ext{(U\times V)}\freccia\ext{(\R^{\sf t})}$;
\item consider the composition $\ext{g}\circ\langle-,-\rangle:\ext{U}\times\ext{V}\freccia\ext{(\R^{\sf t})}$,
where $\langle-,-\rangle$ is the map $\alpha$ given by \eqref{eq:alpha};
\item fix a parameter $p\in\ext{U}$ as a first variable of the previous
composition, i.e. consider $\ext{g}\langle p,-\rangle:\ext{V}\freccia\ext{(\R^{\sf t})}$.
Locally, the map $f$ is of this form: $f=\ext{g}\langle p,-\rangle$.
\end{enumerate}
\item Because in the Definition \ref{def:arrowsOfSERInfty} we ask $s\in\ext{V}$
we have that $\mathcal{V}:=\ext{V}\cap S$ is a neighborhood of $s$
defined by $V$ in $S$ (see \eqref{eq:topologyInSERInfty}). Analogously
$\ext{U}$ is a neighborhood of the parameter $p$.
\end{enumerate}
\end{rem*}
To simplify the presentation, in case the context will be sufficiently
clear, we shall consider the coupling of variables%
\footnote{Note the use of a different font for the second variable in the pairing,
so that it will be easier to identify such pairings.%
} $(S,\sf s)$, $(T,\sf t)$, $(p,\sf p)$, $(q,\sf q)$ etc. in properties
of the form $S\subseteq\R^{\sf s}$, $T\subseteq\R^{\sf t}$, $p\in\ext{(\R^{\sf p})}$
or $q\in\ext{(\R^{\sf q})}$ respectively. In fact, in these cases
we have that the second variable in the pairing, e.g. the number $\sf s\in\N$
in the pairing $(S,\sf s)$, is uniquely determined by the first variable
$S$. E.g. the number $\sf p\in\N$ is uniquely determined by the
point $p\in\ext{(\R^{\sf p})}$. Therefore, if we denote by $\sigma(V)\in\N$
the unique $\sf v\in\N$ in a pairing $(V,\sf v)$, then any formula
of the form $\mathcal{P}(V,\sf v)$ can be interpreted as\[
{\sf v}=\sigma(V)\then\mathcal{P}(V,\sf v).\]

Now we have to prove that $\SERInfty$ verifies the hypothesis of
Section \ref{sec:hypothesisForCartesianClosure} about the category
of the types of figures. Firstly, we prove that $\SERInfty$ is indeed
a category. In the following proofs we will frequently use the properties\[
x\in\ext{U}\iff\st{x}\in U\]
\[
\st{\left(\ext{g}\langle p,x\rangle\right)}=g(\st{p},\st{x}).\]
The first one follows from \eqref{eq:FermatExtOf_U}, and the second
one can be proved directly:\[
\st{\left(\ext{g}\langle p,x\rangle\right)}=\left(g(p_{t},x_{t})\right)|_{t=0}=g(p_{0},x_{0})=g(\st{p},\st{x}).\]

\begin{thm}
\label{thm:SERInftyIsACategory}$\SERInfty$ is a category
\end{thm}
\noindent \textbf{Proof:} In this proof we will consider the coupling
of variables $(S,\sf s)$, $(T,\sf t)$, $(R,\sf r)$, $(p,\sf p)$
and $(q,\sf q)$. If we consider any $p\in\ER$ and the projection
$g:(r,s)\in\R\times\R^{\sf s}\mapsto s\in\R^{\sf s}$, then we have
that $\ext{g}\langle p,s\rangle_{t}=g(p_{t},s_{t})=s_{t}$, hence
$\ext{g}\langle p,s\rangle=s$ and this suffices to prove that the
identity $1_{S}$ for $S\in\SERn$ is always an arrow of $ $$\SERInfty$.

Now let us consider\[
S\xfreccia{f}T\xfreccia{g}R\e{in}\SERInfty\]
and a point $s\in S$. We have to prove that $f\circ g$ is again
an arrow of $\SERInfty$. Using self-evident notations we can assert
that\begin{equation}
f(x)=\ext{h}\langle p,x\rangle\quad\forall x\in\ext{V_{s}}\cap S\ni s\label{eq:SERnCategory_f_h}\end{equation}
\begin{equation}
g(y)=\ext{k}\langle q,y\rangle\quad\forall y\in\ext{V_{fs}}\cap T\ni f(s),\label{eq:SERnCategory_g_k}\end{equation}
where $h\in\Cc^{\infty}(U_{p}\times V_{s},\R^{\sf t})$ and $k\in\Cc^{\infty}(U_{q}\times V_{fs},\R^{\sf r})$.
Hence $U_{q}\times h^{-1}(V_{fs})$ is open in $\R^{\sf q}\times\R^{\sf p}\times\R^{\sf s}$.
But $\st{f(s)}=h(\st{p},\st{s})\in V_{fs}$ because $f(s)\in\ext{V_{fs}}$,
and $\st{q}\in U_{q}$, so\begin{equation}
(\st{q},\st{p},\st{s})\in U_{q}\times h^{-1}(V_{fs}).\label{eq:St_q_p_s}\end{equation}
Hence, we can find three open sets $A\subseteq\R^{\sf q}$, $B\subseteq\R^{\sf p}$
and $C\subseteq\R^{\sf s}$ such that$(\st{q},\st{p},\st{s})\in A\times B\times C\subseteq U_{q}\times h^{-1}(V_{fs})$
and we can correctly define\[
\rho:(x_{1},x_{2},y)\in A\times B\times C\mapsto k\left[x_{1},h(x_{2},y)\right]\in\R^{\sf r},\]
obtaining a map $\rho\in\Cc^{\infty}(A\times B\times C,\R^{\sf r})$;
this is the first step to prove that locally the composition $f(g(-))$
is of the form \eqref{eq:localFormForArrowsOfSERn}. The parameter
corresponding to this local form is $\langle q,p\rangle\in\ext{(A\times B)}$
because of \eqref{eq:alpha} and \eqref{eq:St_q_p_s}. The neighborhood
we are searching for this local equality is $\ext{C}\cap S\ni s$,
in fact let us take a generic $x\in\ext{C}\cap S$, then\begin{equation}
(\st{p},\st{x})\in B\times C\subseteq h^{-1}(V_{fs})\subseteq U_{p}\times V_{s}.\label{eq:SERnCategory_st_p_st_x}\end{equation}
Therefore, $\st{x}\in V_{s}$ and hence $x\in\ext{V_{s}}\cap S$ so
that we can use \eqref{eq:SERnCategory_f_h} obtaining $f(x)=\ext{h}\langle p,x\rangle$.
We can continue, saying that then $\st{f(x)}=h(\st{p},\st{x})\in V_{fs}$,
because of \eqref{eq:SERnCategory_st_p_st_x}, and hence $f(x)\in\ext{V_{fs}}\cap T$.
Now we can apply \eqref{eq:SERnCategory_g_k} with $y=f(x)$ obtaining
\[
g(f(x))=\ext{k}\langle q,fx\rangle=\ext{k}\langle q,\ext{h}\langle p,x\rangle\rangle=\ext{\rho}\langle q,p,x\rangle.\]
\qedWithFinalEq

To prove that $\SERInfty$ is a subcategory of the category $ $$\textbf{Top}$
of topological spaces, we need the following
\begin{thm}
\label{thm:SERInftyIsASubCategoryOfTop}If $f:S\freccia T$ in $\SERInfty$,
then $f$ is continuous with respect to the topologies $\Top{S}$
and $\Top{T}$.
\end{thm}
\noindent \textbf{Proof:} Take $A$ open in $T$ and $s\in f^{-1}(A)$;
we have to prove that, for some $W_{s}$ open in $\R^{\sf s}$, we
have $s\in\ext{W_{s}}\cap S\subseteq f^{-1}(A)$ (see \eqref{eq:topologyInSERInfty}).
From $f(s)\in A\in\Top{T}$ we have that $f(s)\in\ext{V_{fs}}\cap T\subseteq A$
for some open set $V_{fs}\subseteq\R^{\sf t}$. On the other hand,
from $s\in f^{-1}(A)\subseteq S$ and the Definition \ref{def:arrowsOfSERInfty}
of arrow in $\SERInfty$, it follows that in a neighborhood $\mathcal{V}_{s}:=\ext{V_{s}}\cap S$
of $s$ we can write the function $f$ as $f=\ext{g}\langle p,-\rangle$,
where $p\in\ext{U}\subseteq\ext{(\R^{\sf p})}$ is the usual parameter.
Diagrammatically the situation is as follow:\[
\xyR{40pt}\xyC{40pt}\xymatrix{S\ar[r]^{f} & T & V_{fs}\ni f(s)\ar@{_{(}->}[l]\\
s\in\mathcal{V}_{s}\ar@{_{(}->}[u]\ar[ru]_{\ext{g}\langle p,-\rangle}\ar@{^{(}->}[r] & \ext{(\R^{\sf s})}}
\]
Intuitively, the idea is to consider the standard part of $f$ and
to define the open set $W_{s}$ we searched for using the counter
image, along this standard part, of the open set $V_{fs}$. In fact,
let us define\[
W_{s}:=\left[g(\st{p},-)\right]^{-1}(V_{fs}),\]
then $W_{s}$ is open in $\R^{\sf s}$ and we have \begin{align}
s\in\ext{W_{s}} & \iff\st{s}\in W_{s}\nonumber \\
 & \iff g(\st{p},\st{s})\in V_{fs}\nonumber \\
 & \iff\ext{g}\langle p,s\rangle\in\ext{V_{fs}}\nonumber \\
 & \iff f(s)\in\ext{V_{fs}}.\label{eq:sInExt_W}\end{align}
The latter property $f(s)\in\ext{V_{fs}}$ is true, so we have that
$s\in\ext{W_{s}}\cap S$. It remains to prove that $\ext{W_{s}}\cap S\subseteq f^{-1}(A)$.
Let us take a point $x\in\ext{W_{s}}\cap S$, then $\st{x}\in W_{s}\subseteq V_{s}$,
and hence $g(\st{p},\st{x})\in V_{fs}$. So, $x\in\ext{V_{s}}$ and
$f(x)=\ext{g}\langle p,x\rangle\in\ext{V_{fs}}$. But $f(x)\in T$,
so $f(x)\in\ext{V_{fs}}\cap T\subseteq A$.\qedNoNewLine

In the following theorem we prove that the category $\SERInfty$ is
closed with respect to subspaces (with the induced topology) and the
corresponding inclusion:
\begin{thm}
\label{thm:SERInftyCategorySubspaces}Let $S\subseteq\ext{(\R^{\sf s})}$,
and $U\in\Top{S}$ be an open set, with $i:U\hookrightarrow S$ the
corresponding inclusion. Then we have:
\begin{enumerate}
\item $(U\prec\Top{S})\in\SERInfty$, that is the topology $\Top{U}$ defined
by \eqref{eq:topologyInSERInfty} coincides with the induced topology
$\Top{(U\prec S)}$.
\item \label{enu:SERInftyCategoryInclusionsAreArrows}The inclusion $i:U\freccia S$
is an arrow of $\SERInfty$.
\end{enumerate}
\end{thm}
\noindent \textbf{Proof : }By \eqref{eq:topologyInSERInfty}, if $A\in\Top{U}$
we have that $A$ is the union of $\ext{V}\cap U\subseteq A$ for
$V$ open in $\R^{\sf s}$. But $\ext{V}\cap U=\left(\ext{V}\cap S\right)\cap U$
because $U\subseteq S$. Therefore, $A$ is the union of sets of the
form $W\cap U$ with $W\in\Top{S}$, because $W:=\ext{V}\cap S\in\Top{S}$,
i.e. $A$ is open in the subspace $(U\prec\Top{S})$. Vice versa if
we can write $A=B\cap U$, where $B$ is open in $\Top{S}$, then
by \eqref{eq:topologyInSERInfty}\[
\forall s\in A\,\exists\, V\,\text{open in }\R^{\sf s}\pti s\in\ext{V}\cap S\subseteq B,\]
so $s\in A\subseteq U$ and hence $s\in\ext{V}\cap U\subseteq B\cap U=A$,
and this proves that $A\in\Top{U}$, thus $\Top{U}=\Top{(U\prec S)}$.

Property \ref{enu:SERInftyCategoryInclusionsAreArrows}. can be proved
following ideas similar to those used in Theorem \ref{thm:SERInftyIsACategory}
to prove that the identities $1_{S}$ are always arrows of the category
$\SERInfty$.\qedNoNewLine

Now we will prove the closure of $\SERInfty$ with respect to restriction
to open sets (see Hypothesis \ref{ass:basicHypothesesOn_F}):
\begin{thm}
Let $f:S\freccia T$ in $\SERInfty$, $U$ an open set in $S$ and
$V$an open set in $T,$ with $f(U)\subseteq V$. Then\[
f|_{U}:(U\prec S)\freccia(V\prec T)\e{in}\SERInfty.\]

\end{thm}
\noindent \textbf{Proof:} Recalling the Definition \ref{def:arrowsOfSERInfty}
of an arrow in $\SERInfty$, and using the fact that, by hypotheses
we already have that $f:S\freccia T$ in $\SERInfty$, we only have
to prove that equalities of the form \eqref{eq:localFormForArrowsOfSERn}
hold locally also with respect to the topology of $(U\prec S)$. Because
of the previous Theorem \ref{thm:SERInftyCategorySubspaces} we can
work with $\Top{U}$ instead of $\Top{(U\prec S)}$. Take $s\in U$,
since $U\subseteq S$ and $f:S\freccia T$ in $\SERInfty$, using
the usual notations we can write\begin{equation}
f(x)=\ext{g}\langle p,x\rangle\label{eq:SERnCategoryRestrictions_f_g}\end{equation}
for every $x\in\ext{V_{s}}\cap S\ni s$ and where $g\in\Cc^{\infty}(U_{p}\times V_{s})$.
Hence $s\in\ext{V_{s}}\cap U$, and because $\ext{V_{s}}\cap U\subseteq\ext{V_{s}}\cap S$
we have again the equality \eqref{eq:SERnCategoryRestrictions_f_g}
in the neighborhood $\ext{V_{s}}\cap U$ of $s$ in $\Top{U}$, and
this proves the conclusion.\qedNoNewLine

Since it is trivial to prove that $\SERInfty$ contains all the constant
maps (it suffices to take $g(p,x):=p$), to prove that $\SERInfty$
is a category of the types of figures, it remain to prove the sheaf
property:
\begin{thm}
\label{thm:SheafPropertyOfSERn}Let $H$, $K\in\SERInfty,$ and $(H_{i})_{i\in I}$
be an open cover of $H$ such that the map $f:H\freccia K$ verifies\begin{equation}
\forall i\in I\pti f|_{H_{i}}\in\SERInfty(H_{i},K).\label{eq:SheafSERn_f_restrictedTo_Hi}\end{equation}
Then\[
f:H\freccia K\e{in}\SERInfty.\]

\end{thm}
\textbf{Proof:} Take $s\in H$, then $s\in H_{i}$ for some $i\in I$,
and from \eqref{eq:SheafSERn_f_restrictedTo_Hi} it follows that we
can write\[
\left(f|_{H_{i}}\right)(x)=f(x)=\ext{g}\langle p,x\rangle\quad\forall x\in\ext{V_{s}}\cap H_{i}.\]
But $H_{i}$ is open in $H$ so that we can also say that $s\in\ext{V'_{s}}\cap H\subseteq H_{i}$
for some open set $V'_{s}$ of $\R^{\sf h}$. The new neighborhood
$\ext{(V_{s}\cap V'_{s})}\cap H$ of $s$ and the restriction $g|_{U_{p}\times(V_{s}\cap V'_{s})}$
verify that the function $f$ is locally of the form $f=\ext{g}\langle p,-\rangle$
in a neighborhood of $s$ in $H$.\qedNoNewLine

We have proved that $\SERInfty$ and the forgetful functor $|-|$
verify the hypotheses of Section \ref{sec:hypothesisForCartesianClosure}
about the category of the types of figures and hence we can define
\[
\ECInfty:=\overline{\SERInfty}.\]
 Each object of $\ECInfty$ will be called a \emph{Fermat space}.

We close this section with the following simple but useful result
that permits to obtain functions in $\SERInfty$ starting from ordinary
$\Cc^{\infty}$ functions.
\begin{thm}
Let $f\in\Cc^{\infty}(\R^{\sf k},\R^{\sf h})$ be a standard $\Cc^{\infty}$
function and $H\subseteq\ext{(\R^{\sf h})}$ and $K\subseteq\ext{(\R^{\sf k})}$
be subsets of Fermat reals. If the function $f$ verifies $\ext{f}|_{K}(K)\subseteq H$,
then\[
\ext{f}|_{K}:K\freccia H\e{in}\SERInfty.\]

\end{thm}
\textbf{Proof:} It suffices to define $g(x,y):=f(y)$ for $x\in\R$
and $y\in\R^{\sf k}$ to obtain that\[
\ext{g}\langle0,k\rangle_{t}=g(0,k_{t})=f(k_{t})=\ext{f}(k)_{t},\]
that is $\left(\ext{f}|_{K}\right)=\ext{g}\langle0,k\rangle\in H$
for every $k\in K$.\qedNoNewLine

\chapter{\label{cha:TheFermatFunctor}The Fermat functor}

\section{Putting a structure on the sets $\ext{X}$}

\noindent Now the problem is: what Fermat space could we associate
to sets like $\ext{X}$ or $D$? 
\begin{defn}
\noindent \label{def:FermatFunctor}Let $X\in\CInfty$, then for any
subset $Z\subseteq\ext{X}$ we call $\ext{(ZX)}$ the extended space
generated on $Z$ (see Section \ref{sec:DefinitionOfCartesianClosure})
by the following set of figures $d:T\freccia Z$ (where $T\subseteq\ext{(\R^{\sf t})}$
is a type of figure in $\SERInfty$) \begin{equation}
\begin{aligned}d\in\Ds{T}^{0}(Z)\DIff & \text{$d$ is constant or we can write}\\
{} & \text{$d=\ext{h}|_{T}$ for some $h\In{V}X$ such that $T\subseteq\ext{V}$.}\end{aligned}
\label{def: GeneratingSetForExtension}\end{equation}

\end{defn}
Thus in the non-trivial case we start from a standard figure $h\In{V}X$
of type $V\in\ORInfty$ such that $\ext{V}\supseteq T$; we extend
this figure obtaining $\ext{h}:\ext{V}\freccia\ext{X}$, and finally
the restriction $\ext{h}|_{\sss T}$ is a generating figure if it
maps $T$ in $Z$. This choice is very natural, and the adding of
the alternative {}``$d$ is constant'' in the previous disjunction
is due to the need to have all constant figures in a family of generating
figures.

Using this definition of $\ext{(ZX)}$, we set (with some abuses of
language) \begin{gather*}
\ext{X}:=\ext{(\ext{X}X)}\\
D:=\ext{(D\R)}\\
\ER:=\ext{(\ext{\R}\R)}\\
\ER^{k}:=\ext{(\ext{(\R^{k})}\R^{k})}\\
D_{k}:=\ext{(D_{k}\R^{k})}.\end{gather*}
 We will call $\ext{(ZX)}$ \emph{the Fermat space induced on $Z$
by $X\in\CInfty$}. We can now study the extension functor:
\begin{thm}
\label{thm:FermatIsAFunctor}Let $f\in\CInfty(X,Y)$ and $Z$ a subset
of $\ext{X}$ with $\ext{f}(Z)\subseteq W\subseteq\ext{Y}$, then
in $\ECInfty$ we have that \[
\ext{(ZX)}\xfrecciad{\ext{f}|_{Z}}\ext{(WY)}.\]
 Therefore $\ext{(-)}:\CInfty\freccia\ECInfty$ is a functor, called
the \emph{Fermat functor}.
\end{thm}
\noindent \textbf{Proof:} Take a figure $\delta\In{S}\ext{(ZX)}$
of type $S\in\SERInfty$ in the domain. We have to prove that $\delta\cdot\ext{f}|_{Z}$
locally factors through $\SERInfty$ and $\D^{0}(W)$ (see in Section
\ref{sec:DefinitionOfCartesianClosure} the definition of space generated
by a family of figures). Hence taking $s\in S$, since $\delta\In{S}\ext{(ZX)}$,
we can write $\delta|_{U}=f_{1}\cdot d$, where $U$ is an open neighborhood
of $s$, $f_{1}\in\SERInfty(U\prec S,T)$ and $d\in\Ds{T}^{0}(Z)$:\[
\xyR{45pt}\xyC{55pt}\xymatrix{S\ar[r]^{\delta} & Z\\
s\in U\ar[r]_{f_{1}}\ar@{_{(}->}[u]\ar[ru]^{\delta|_{U}} & T\ar[u]_{d}}
\]
We omit the trivial case $d$ constant, hence we can suppose, using
the same notations used in the Definition \ref{def:FermatFunctor},
to have $d=\ext{h}|_{T}:T\freccia Z$ with $h\In{V}X$. Therefore
\[
(\delta\cdot\ext{f}|_{Z})|_{U}=\delta|_{U}\cdot\ext{f}|_{Z}=f_{1}\cdot d\cdot\ext{f}|_{Z}=f_{1}\cdot\ext{h}|_{T}\cdot\ext{f}|_{Z}=f_{1}\cdot\ext{(hf)}|_{T}.\]
 But $hf\In{V}Y$ since $f\in\CInfty(X,Y)$ and $h\In{V}X$, so $(\delta\cdot\ext{f}|_{Z})|_{U}=f_{1}\cdot d_{1}$,
where $d_{1}:=\ext{(hf)}|_{T}\in\Ds{T}^{0}(W)$, which is the conclusion.
The other functorial properties, i.e. $\ext{\left(1_{X}\right)}=1_{\ext{X}}$
and $\ext{\left(f\cdot g\right)}=\ext{f}\cdot\ext{g}$, follow directly
from the definition of the Fermat extension $\ext{f}$ of $f\in\CInfty(X,Y)$.$\qedWithFinalEq$

\section{\label{sec:FermatFunctorAndProductOfManifolds}The Fermat functor
preserves product of manifolds}

We want to prove that the bijective applications $\alpha$ defined
in \ref{eq:alpha} and $\beta$ defined in \ref{eq:beta}, i.e.\begin{equation}
\alpha_{\sss MN}:([x]_{\sim},[y]_{\sim})\in\ext{M}\times\ext{N}\longmapsto[(x,y)]_{\sim}\in\ext{(M\times N)}\label{eq:alphaBis}\end{equation}
\begin{equation}
\beta_{\sss MN}:[z]_{\sim}\in\ext{(M\times N)}\longmapsto([z\cdot p_{M}]_{\sim},[z\cdot p_{N}]_{\sim})\in\ext{M}\times\ext{N}\label{eq:betaBis}\end{equation}
are arrows of $ECInfty$. Where it will be clear from the context,
we shall use the simplified notations $\alpha:=\alpha_{\sss MN}$
and $\beta:=\beta_{\sss MN}$. To simplify the proof we will use the
following preliminary results. The first one is a general property
of the cartesian closure $\bar{\F}$ of a category of figures $\F$
(see Chapter \ref{cha:theCartesianClosure}). 
\begin{lem}
\label{lem: ProductOfGeneratedObjects} Suppose that $\F$ admits
finite products $K\times J$ for every objects $K$, $J\in\F$, and
an isomorphism%
\footnote{Recall the definition of the embedding $\bar{(-)}:\F\freccia\bar{\F}$
given in Section \ref{sec:DefinitionOfCartesianClosure}.%
}\[
\xymatrix{\xyC{40pt}\gamma_{\sss{KJ}}:\overline{K\times J}\ar[r]\sb(0.57)\sim & \bar{K}\times\bar{J}}
\e{in}\bar{\F}.\]
Moreover, let $Z$, $X$, $Y\in\bar{\F}$ with $X$ and $Y$ generated
by $\D^{\sss X}$ and $\D^{\sss Y}$ respectively. Then we have\[
\xymatrix{\xyC{40pt}X\times Y\ar[r]\sp(0.57){\displaystyle f} & Z}
\e{in}\bar{\F}\]
if and only if for any $K$, $J\in\F$ and $d\in\Ds{K}^{\sss X}$,
$\delta\in\Ds{J}^{\sss Y}$ we have \[
\gamma_{\sss{KJ}}\cdot(d\times\delta)\cdot f\In{K\times J}Z.\]

\end{lem}
The second Lemma asserts that the category of figures $\F=\SERInfty$
verifies the hypotheses of the previous one.
\begin{lem}
\label{lem: ProductInSERn} The category $\SERInfty$ admits finite
products and the above mentioned isomorphisms $\gamma_{\sss KJ}$.
For $K\subseteq\ext{(\R^{\sf k})}$ and $J\subseteq\ext{(\R^{\sf j})}$
these are given by\[
K\times J=\left\{ \langle x,y\rangle\in\ext{(\R^{{\sf k}+{\sf j}})}\,|\, x\in K\ ,\ j\in J\right\} \]
\[
\gamma_{\sss KJ}:\langle x,y\rangle\in K\times J\longmapsto(x,y)\in K\times_{\text{{\rm s}}}J,\]
where we recall that $\langle x,y\rangle=\alpha_{\sss\R^{\sf k}\R^{\sf j}}(x,y)=[t\mapsto(x_{t},y_{t})]_{\sim}$
and where $K\times_{\text{{\rm s}}}J$ is the set theoretical product
of the subsets $K$ and $J$.

\noindent Moreover let $M$, $N$ be $\Cc^{\infty}$ manifolds, and
$h\In{V}M$, $l\In{V'}N$ with $K\subseteq\ext{V}$ and $J\subseteq\ext{V'}$,
then \[
\gamma_{\sss KJ}\cdot(\ext{h}|_{\sss K}\times\ext{l}|_{\sss J})\cdot\alpha_{\sss MN}=\ext{(h\times l)}|_{\sss K\times J}.\]

\end{lem}
\noindent The proofs of these lemmas are direct consequences of the
given definitions.
\begin{thm}
\label{thm:FermatPreservesProductOfMan}Let $M$, $N$ be $\Cc^{\infty}$
manifolds, then in $\ECInfty$ we have the isomorphism \[
\ext{(M\times N)}\simeq\ext{M}\times\ext{N}.\]

\end{thm}
\noindent \textbf{Proof:} Note that in the statement each manifold
is identified with the corresponding $\CInfty$ space $\bar{M}$.
Hence we mean $\ext{M}=\ext{\bar{M}}=\ext{(\ext{M}\bar{M})}$ (see
Definition \ref{def:FermatFunctor} for the notation $\ext{(ZX)}$).
To prove that $\alpha$ is a $\ECInfty$ arrow we can use Lemma \ref{lem: ProductOfGeneratedObjects},
because of Lemma \ref{lem: ProductInSERn} and considering that $\ext{M}$
and $\ext{N}$ are generated by $\D^{0}(\ext{M})$ and $\D^{0}(\ext{N})$.
Since these generating sets are defined using a disjunction (see \eqref{def: GeneratingSetForExtension})
we have to check four cases depending on $d\in\Ds{K}^{0}(\ext{M})$
and $\delta\in\Ds{J}^{0}(\ext{N})$. In the first case we have $d=\ext{h}|_{\sss K}\in\Ds{K}^{0}(\ext{M})$
and $\delta=\ext{l}|_{\sss J}\in\Ds{J}^{0}(\ext{N})$ (we are using
the same notations of the previous Lemma \ref{lem: ProductInSERn}).
Thus \[
\gamma_{\sss KJ}\cdot(d\times\delta)\cdot\alpha=\gamma_{\sss KJ}\cdot(\ext{h}|_{\sss K}\times\ext{l}|_{\sss J})\cdot\alpha=\ext{(h\times l)}|_{\sss K\times J}.\]
That is $\gamma_{\sss KJ}\cdot(d\times\delta)\cdot\alpha$ is a generating
element in $\ext{(M\times N)}$, and so it is also a figure. In the
second case let us suppose $\delta$ constant with value $n\in\ext{N}$
and $d=\ext{h}|_{\sss K}\in\Ds{K}^{0}(\ext{M})$. Take a chart $l:\R^{\sf p}\freccia U$
on $\st{n}=n_{0}\in U\subseteq N$ and let $W:=\ext{(\R^{\sf p})}$,
$p:=\ext{l^{-1}}(n)\in W$. Note that $\ext{l}(p)=n=\delta(-)$. We
have to prove that $\gamma_{\sss KJ}\cdot(d\times\delta)\cdot\alpha\In{K\times J}\ext{(M\times N)}$,
so let us start to calculate the map $\gamma_{\sss KJ}\cdot(d\times\delta)\cdot\alpha$
at a generic element $\langle k,j\rangle\in K\times J$. We have\begin{align}
\alpha\{(d\times\delta)[\gamma_{\sss{KJ}}(\langle k,j\rangle)]\} & =\alpha\left[(d\times\delta)(k,j)\right]\nonumber \\
 & =\alpha\left[d(k),n\right]\nonumber \\
 & =\alpha[\ext{h}(k),\ext{l}(p)]\nonumber \\
 & =\{\gamma_{\sss{KW}}\cdot[\ext{h}|_{\sss K}\times\ext{l}|_{\sss J}]\cdot\alpha\}\langle k,p\rangle\nonumber \\
 & =\ext{(h\times l)}|_{\sss K\times W}\langle k,p\rangle,\label{eq: AlphaBetaDelta}\end{align}
where we have used once again the equality of Lemma \ref{lem: ProductInSERn}.
Thus let us call $\tau$ the map $\tau:\langle k,j\rangle\in|K\times J|\mapsto\langle k,p\rangle\in|K\times W|$,
so that we can write \eqref{eq: AlphaBetaDelta} as \[
\gamma_{\sss{KJ}}\cdot(d\times\delta)\cdot\alpha=\tau\cdot\ext{(h\times l)}|_{\sss K\times W}.\]
But $\ext{(h\times l)}|_{\sss K\times W}$ is a generating figure
of $\ext{(M\times N)}$ and $\tau$ is an arrow of $\SERInfty$, and
this proves that $\gamma_{\sss{KJ}}\cdot(d\times\delta)\cdot\alpha\In{K\times J}\ext{(M\times N)}$.
The remaining cases are either trivial (both $d$ and $\delta$ constant)
or analogous to the latter one.

To prove that the map $\beta_{\sss MN}$ is an arrow of $\ECInfty$
is simpler. Indeed, take $d\In{H}\ext{(M\times N)}$ to prove that
$d\cdot\beta_{\sss MN}\In{H}\ext{M}\times\ext{N}$. Due to the universal
property of the product $\ext{M}\times\ext{N}$, it suffices to consider
the composition of this map $d\cdot\beta_{\sss MN}$ with the projections
of this product. But, if $p_{\sss M}:M\times N\freccia M$ is the
projection on $M$, then \[
\xymatrix{\xyC{45pt}\ext{M}\times\ext{N}\ar[r]^{{\displaystyle \alpha}{}_{\sss MN}} & \ext{(M\times N)}\ar[r]^{{\displaystyle {\ \ \ \ \ext{p_{\sss M}}}}}\ar[r] & \ext{M}}
\]
and $\ext{p_{\sss M}}(\alpha_{\sss MN}(x,y))=\ext{p_{\sss M}}(\langle x,y\rangle)=x$,
so $\alpha_{\sss MN}\cdot\ext{p_{\sss M}}$ is the projection of the
product $\ext{M}\times\ext{N}$ on $\ext{M}$. Therefore the conclusion
$d\cdot\beta_{\sss MN}\In{H}\ext{M}\times\ext{N}$ is equivalent to\[
d\cdot\beta_{\sss MN}\cdot\alpha_{\sss MN}\ext{p_{\sss M}}=d\cdot\ext{p_{\sss M}}\In{H}\ext{M}\]
\[
d\cdot\beta_{\sss MN}\cdot\alpha_{\sss MN}\ext{p_{\sss N}}=d\cdot\ext{p_{\sss N}}\In{H}\ext{M}\]
 which are true since $\ext{p_{\sss M}}$ and $\ext{p_{\sss N}}$
are arrows of $\ECInfty$.$\qedWithFinalEq$

In the following we shall always use the isomorphism $\alpha$ to
identify these spaces, hence we write $\ext{M}\times\ext{N}=\ext{(M\times N)}$,
e.g. $\ext{(\R^{d})}=\left(\ER\right)^{d}=:\ER^{d}$.

\subsection{Figures of Fermat spaces}

\label{sec: FiguresOfExtendedSpaces} In this section we want to understand
better the figures of the Fermat space $\ext{(ZX)}$; we will use
these results later, for example when we will study the embedding
of $\ManInfty$ into $\ECInfty$, or to prove some logical properties
of the Fermat functor.

From the general definition of $\bar{\F}$--space generated by a family
of figures $\D^{0}$ (see Section \ref{sec:DefinitionOfCartesianClosure}),
a figure $\delta\In{S}\ext{(ZX)}$, for $S\in\SERInfty$, can be locally
factored as $\delta|_{V}=f\cdot d$ through an arrow $f\in\SERInfty(V,T)$
and a generating function $d\in\Ds{T}^{0}(Z)$; here $V=V(s)$ is
an open neighborhood of the considered point $s\in S$, so that we
can always suppose $V$ to be of the form $V=\ext{B}\cap S$ (see
in \eqref{eq:topologyInSERInfty} the definition of topology for $S$).
Hence, either $\delta|_{V}$ is constant (if $d$ is constant) or
we can write $d=\ext{h}|_{T}$ and $f=\ext{g}(p,-)$ so that \[
\delta(x)=d[f(x)]=\ext{h}[\ext{g}(p,x)]=\ext{(gh)}(p,x)\quad\forall x\in V=\ext{B}\cap S,\]
where $A\times B$ is an open neighborhood of $(\st{p},\st{s})$.
Therefore we can write \[
\delta(x)=\ext{\gamma}(p,x)\quad\forall x\in\ext{B}\cap S,\]
with $\gamma:=g|_{A\times B}\cdot h\in\CInfty(A\times B,X)$. Thus
figures of $\ext{(ZX)}$ are locally necessarily either constant maps
or a natural generalization of the maps of $\SERInfty$, that is {}``extended
$\CInfty$ arrows $\ext{\gamma}(-,-)$ with a fixed parameter $\ext{\gamma}(p,-)$''.
Using the properties of $\ECInfty$ and of its arrow $\alpha_{\R^{\sf p}\R^{\sf s}}$
it is easy to prove that these conditions are sufficient too. Moreover
if $X=M$ is a manifold, the condition {}``$\delta|_{V}$ constant''
can be omitted. In fact if $\delta|_{\mathcal{V}}$ is constant with
value $m\in Z\subseteq\ext{M}$, then taking a chart $\varphi$ on
$\st{m}\in M$ we can write $\delta(x)=m=\ext{\gamma}(p,x)$, where
$p=\ext{\varphi}(m)$ and $\gamma(x,y)=\varphi^{-1}(x)$. We have
proved the following
\begin{thm}
\label{thm:figuresOfFermatSpaces}Let $X\in\CInfty$, $Z\subseteq\ext{X}$,
$S\subseteq\ER^{\sf s}$ and $\delta:S\freccia Z$. Then we have\[
\delta\In{S}\ext{(ZX)}\]
iff for every point $s\in S$ there exist an open set $B$ in $\R^{\sf s}$
such that $s\in\ext{B}$ and such that either\begin{equation}
\delta|_{\ext{B}\cap S}\text{ is constant},\label{eq:deltaIsConstant}\end{equation}
or we can write\[
\delta(x)=\ext{\gamma}\langle p,x\rangle\quad\forall x\in\ext{B}\cap S\]
for some\begin{align*}
 & p\in\ext{A},\text{ where }A\text{ is open in }\R^{\sf p}\\
 & \gamma\in\Cc^{\infty}(A\times B,X).\end{align*}
Moreover if $X=M$ is a manifold, condition \eqref{eq:deltaIsConstant}
can be omitted and there remains only the second alternative.$\qedWithFinalEq$
\end{thm}
Using this result we can prove several useful properties of the Fermat
functor. The following ones say that we can arrive at the same Fermat
space starting from several different constructions.
\begin{thm}
\label{thm:spacesPropertiesOfFermatFunctor}The Fermat functor has
the following properties:
\begin{enumerate}
\item \label{enu:twoWays}If $X\in\CInfty$ and $Z\subseteq|\ext{X}|$,
then $\ext{(ZX)}=(Z\prec\ext{X})$.
\item \label{enu:threeWays}If $S\subseteq|\ER^{\sf s}|$, then $\bar{S}=\ext{(S\R^{\sf s})}=(S\prec\ER^{\sf s})$.
\end{enumerate}
\end{thm}
\noindent E.g. if $f:\ER^{\sf s}\freccia\ext{X}$ is a $\ECInfty$
arrow, then we also have $f:\overline{\ER^{\sf s}}\freccia\ext{X}$
because $\ER^{\sf s}=\ext{(\ext{(\R^{\sf s})}\R^{\sf s})}=\overline{\ER^{\sf s}}$
and the previous property \emph{\ref{enu:threeWays}} holds. Therefore,
$f\In{\ER^{\sf s}}\ext{X}$ and locally we can write $f$ either as
a constant function or, with the usual notations, as $f(x)=\ext{\gamma}(p,x)$.
For functions $f:I\freccia\ext{X}$ defined on some set $I\subseteq D_{\infty}$
of infinitesimals which contains $0\in I$, these two alternatives
globally holds instead of only locally, because the set of infinitesimals
$I$ is contained in any open neighborhood of $0$.

\noindent \textbf{Proof:} To prove \emph{\ref{enu:twoWays}}. let
us consider a figure $\delta\In{S}(Z\prec\ext{X})$ of type $S\in\SERInfty$
and let $i:Z\hookrightarrow|\ext{X}|$ be the inclusion. We have to
prove that $\delta\In{S}\ext{(ZX)}$, and we will prove it locally,
that is using the sheaf property of the space $\ext{(ZX)}$. By the
definition of subspace, we have that $\delta\cdot i=\delta\In{S}\ext{X}=\ext{(|\ext{X}|X)}$,
so that for every $s\in S$ we can locally factor the figure $\delta$
through $\SERInfty$ and a generating figure $d\in\Ds{K}^{0}(|\ext{X}|)$,
i.e. $\delta|_{U}=f\cdot d$ for some open neighborhood $U$ of $s$
and some $f:(U\prec S)\freccia K$ in $\SERInfty$. If $d$ is constant,
then so is $\delta|_{U}$ and hence $\delta|_{U}\In{U}\ext{(ZX)}$.
Otherwise, we can write $d=\ext{h}|_{K}$ for some $h\In{V}X$, with
$V$ open in $\R^{\sf k}$ such that $K\subseteq\ext{V}$ (see Definition
\ref{def:FermatFunctor}). To prove that $\delta|_{U}\In{U}\ext{(ZX)}$
we exactly need to prove that the map $\delta|_{U}$ factors in the
same way, but with a generating figure $d'$ having values in $Z$
and not in the bigger $|\ext{X}|$ (like $d$ does). For this reason
we change the subset $K$ with the smaller $K':=f(U)\subseteq K\subseteq\ext{V}\subseteq\ER^{\sf k}$,
so $K'\in\SERInfty$, and we set $d':=\ext{h}|_{K'}$. The map $d'$
has values in $Z$, in fact for $x\in K'=f(U)$ we have $x=f(u)$
for some $u\in U$, and \[
d'(x)=\ext{h}(x)=\ext{h}(f(u))=d(f(u))=\delta(u)\in Z.\]
Hence $d'\in\Ds{K}^{0}(Z)$ and $\delta|_{U}(u)=d(f(u))=d'(f(u))$
for every $u\in U$, so $\delta|_{U}\In{U}\ext{(ZX)}$. We have proved
that\[
\forall s\in S\,\exists\, U\text{ open neighborhood of }s\text{ in }S\pti\delta|_{U}\In{U}\ext{(ZX)},\]
hence $\delta\In{S}\ext{(ZX)}$ from the sheaf property of the space
$\ext{(ZX)}\in\ECInfty$. For the opposite inclusion we only have
to make the opposite passage: from $d':K'\freccia Z$ with values
in $Z$ to $d:=d':K'\freccia|\ext{X}|$ with values in the bigger
$|\ext{X}|$, but this is trivial.

Because of the just proved property \emph{\ref{enu:twoWays}}., to
prove \emph{\ref{enu:threeWays}}. we have to verify only the equality
$\bar{S}=\ext{(S\R^{\sf s})}$, so take a figure $\delta\In{T}\bar{S}=(\SERInfty(-,S),S)$
and then $\delta\in\SERInfty(T,S)$. But $1_{\R^{\sf s}}\In{\R^{\sf s}}\R^{\sf s}$,
so $\ext{(1_{\R^{\sf s}})}|_{S}=1_{\ext{\R^{\sf s}}}|_{S}=1_{S}$,
so $1_{S}\in\Ds{S}^{0}(S)$ and $\delta=\delta\cdot1_{S}$ factors
through a map of $\SERInfty(T,S)$ ($\delta$ itself) and a generating
figure of $\Ds{S}^{0}(S)$, i.e. $\delta\In{T}\ext{(S\R^{\sf s})}$.
To prove the opposite inclusion, let us take $\delta\In{T}\ext{(S\R^{\sf s})}$,
then from Theorem \ref{thm:figuresOfFermatSpaces} we have that in
a suitable neighborhood $U$ of a given generic point $s\in T$ we
have that either $\delta|_{U}$ is constant, or we can write $\delta|_{U}=\ext{\gamma}\langle p,-\rangle|_{U}$
for some $\gamma\in\Cc^{\infty}(A\times B,\R^{\sf s})$. In both cases
we have that $\delta|_{U}\in\SERInfty(U,S)$, so $\delta|_{U}\In{U}\bar{S}$,
and the conclusion follows from the sheaf property of $\bar{S}$.$\qedWithFinalEq$

\section{\label{sec:TheEmbeddingOfManifoldsInFermatSpaces}The embedding of
manifolds in $\ECInfty$}

\label{sec: EmbeddingInExtendedSpaces} If we consider a $\CInfty$
space $X$, we have just seen that we have the possibility to associate
a Fermat space to any subset $Z\subseteq|\ext{X}|$. Thus if $X$
is separated we can put a structure of $\ECInfty$ space on the set
$|X|$ of standard points of $X$, by means of $\bar{X}:=\ext{(|X|X)}=(|X|\prec\ext{X})$.
Intuitively $X$ and $\bar{X}$ seem very similar, and in fact we
have 
\begin{thm}
\label{thm: EmbeddingOfSepInExtendedCn} Let $X$, $Y$ be $\CInfty$
separated spaces, then
\begin{enumerate}
\item \label{enu:separatedSpacesAreEmbedded}$\bar{X}=\bar{Y}\then X=Y$
\item \label{enu:arrowsBetweenSepSpacesAreEmbedded}$\bar{X}\xfrecciad{f}\bar{Y}\e{in}\ECInfty\ \iff\ X\xfrecciad{f}Y\e{in}\CInfty.$ 
\end{enumerate}
\noindent Hence $\CInfty$ separated spaces are fully embedded in
$\ECInfty$, and so is $\ManInfty$.
\end{thm}
\noindent \textbf{Proof:} The equality $\bar{X}=\bar{Y}$ implies
the equality of the support sets $|X|=|Y|$. We consider now a figure
$d\In{H}X$ of type $H$, where $H$ is an open set of $\R^{\sf h}$.
Taking the extension of $d$ and then the restriction to standard
points we obtain \begin{equation}
(H\prec\ext{\bar{H}})\xfrecciad{\ext{d}|_{\sss H}}(|X|\prec\ext{X})=\bar{X}=\bar{Y}.\end{equation}
But from Theorem \ref{thm:spacesPropertiesOfFermatFunctor} we have
$(H\prec\ext{\bar{H}})=(H\prec\ER^{\sf h})=\ext{(H\R^{\sf h})}=\bar{H}$,
hence \[
\ext{d}|_{\sss H}=d:\bar{H}\freccia\bar{Y}\e{in}\ECInfty\]
and so $d\In{H}\bar{Y}$. Therefore for every $s\in H$ either $d$
is constant in some open neighborhood $V$ of $s$, or, using the
usual notations, we can write \begin{equation}
d(x)=\ext{\gamma}(p,x)\quad\forall x\in\ext{B}\cap H=B\cap H,\label{eq:d_ext_p_x}\end{equation}
where $\ext{B}\cap H=B\cap H$ because $H\subseteq\R^{\sf h}$ is
made of standard point only. Let us note that the equality in \eqref{eq:d_ext_p_x}
has to be understood in the space $\ext{Y}$. Hence for every $x\in B\cap H$
we have that $\st{d(x)}\asymp\st{[\gamma(p,x)]}$ in $Y$, and so
we can write $d(x)=\gamma(p_{0},x)$ because $Y$ is separated and
$x\in B\cap H\subseteq\R^{\sf h}$ is standard. Therefore $d|_{B\cap H}$
is a $Y$-valued arrow of $\CInfty$ defined in a neighborhood of
the fixed $s$. The conclusion $d\In{H}Y$ thus follows from the sheaf
property of $Y$. Analogously we can prove the opposite inclusion,
so $X=Y$.

If we suppose that $f:\bar{X}\freccia\bar{Y}$ in $\ECInfty$, then
from the proof of \emph{\ref{enu:separatedSpacesAreEmbedded}}. we
have seen that if $d\In{H}X$ then $d\In{H}\bar{X}$. Hence $f(d)\In{H}\bar{Y}$.
But once again from the previous proof of \emph{\ref{enu:separatedSpacesAreEmbedded}}.
we have seen that this implies that $f(d)\In{H}Y$, and so $f:X\freccia Y$
in $\CInfty$.

To prove the opposite implication it suffices to extend $f$ so that
$\ext{f}:\ext{X}\freccia\ext{Y}$, to restrict it to standard points
only so that\[
\ext{f}|_{|X|}:(|X|\prec\ext{X})=\bar{X}\freccia(|Y|\prec\ext{Y})=\bar{Y},\]
and finally to consider that our spaces are separated so that $\ext{f}|_{|X|}=f$.$\qedWithFinalEq$

An immediate corollary of this theorem is that the extension functor
is another full embedding for separated spaces. 
\begin{cor}
\label{cor:FermatFunctorIsAnEmbedding}Let $X,Y$ be $\CInfty$ separated
spaces, then
\begin{enumerate}
\item \label{enu:FermatFunctorIsAnEmbedding}$\ext{X}=\ext{Y}\then X=Y$
\item \label{enu:fromECnToCn}If $\ext{X}\xfrecciad{f}\ext{Y}$ in $\ECInfty$
and $f(|X|)\subseteq|Y|$ then\[
X\xfrecciad{f|_{|X|}}Y\e{in}\CInfty\]

\item \label{enu:FermatFunctorIsAnEmbeddingOnArrows}$\ext{X}\xfrecciad{\ext{f}}\ext{Y}\e{in}\ECInfty\ \iff\ X\xfrecciad{f}Y\e{in}\CInfty$
\item \label{enu:FermatFunctorIsInjectiveOnArrows}If $f$, $g:X\freccia Y$
are $\CInfty$ functions, then\[
\ext{f}=\ext{g}\then f=g.\]

\end{enumerate}
\end{cor}
\noindent \textbf{Proof:} To prove \emph{\ref{enu:FermatFunctorIsAnEmbedding}}.
we start to prove that the support sets of $X$ and $Y$ are equal.
Indeed, if we take standard parts, since $\ext{X}=\ext{Y}$, we have
\[
\{\st{x}\,\,|\,\, x\in\ext{X}\}=|X|=\{\st{x}\,\,|\,\, x\in\ext{Y}\}=|Y|.\]
Hence $\bar{X}=(|X|\prec\ext{X})=(|Y|\prec\ext{Y})=\bar{Y}$ and the
conclusion follows from \emph{\ref{enu:separatedSpacesAreEmbedded}}.
of Theorem \ref{thm: EmbeddingOfSepInExtendedCn}.

To prove \emph{\ref{enu:fromECnToCn}}. let us take the restriction
of $f$ to $|X|\subseteq|\ext{X}|$, then $f|_{|X|}:\bar{X}=(|X|\prec\ext{X})\freccia(|Y|\prec\ext{Y})=\bar{Y}$
in $\ECInfty$, so the conclusion follows from \emph{\ref{enu:arrowsBetweenSepSpacesAreEmbedded}}.
of Theorem \ref{thm: EmbeddingOfSepInExtendedCn}. Property \emph{\ref{enu:FermatFunctorIsAnEmbeddingOnArrows}}.
follows from the just proved \emph{\ref{enu:fromECnToCn}}. considering
that $\ext{f}|_{|X|}=f$ and using Theorem \ref{thm:FermatIsAFunctor}.
The same idea of considering restrictions can be used to prove \emph{\ref{enu:FermatFunctorIsInjectiveOnArrows}}.$\qedWithFinalEq$

\section{The standard part functor cannot exist}

It is very natural to ask if it is possible to define a standard part
functor, that is a way to associate to \emph{every} Fermat space $X\in\ECInfty$
a space $\st{X}\in\CInfty$ intuitively corresponding to its {}``standard
points'' only. This application $\st{(-)}:\ECInfty\freccia\CInfty$
has to satisfy some expected properties, some example of which are
functoriality, its support set has to be included in the original
space, i.e. $|\st{X}|\subseteq|X|$, and we must also have examples
like $\st{D}=\left\{ 0\right\} $ and $\st{(\ER)}=\R$. Because, intuitively,
the Fermat extension $\ext{X}\in\ECInfty$ appears to be some kind
of completion of the standard space $X\in\CInfty$, we also expect
that the Fermat functor is the left adjoint of the standard part functor,
$\ext{(-)}\dashv\st{(-)}$. Indeed, it is natural to expect that this
adjunction is related to the following equivalence%
\footnote{The horizontal line indicating the logical equivalence between the
formula above and the formula below, similar to the notations in the
logical calculus of Gentzen, but where the line indicates logical
deduction of the formula below from the formula above.%
}\begin{equation}
\cfrac{\CInfty\vDash X\xfrecciad{\st{f}}\st{Y}}{\ECInfty\vDash\ext{X}\xfrecciad{f}Y}\label{eq:FermatLeftAdjointOfStandardPart}\end{equation}

If one tries to define this standard part space (and the corresponding
standard part map acting on arrows, i.e. $f\in\ECInfty(X,Y)\mapsto\st{f}\in\CInfty(\st{X},\st{Y})$),
then several difficulties arise.

For example, the first trivial point that has to be noted in the searching
for the definition of $\st{X}$, is that we want to have $|\st{X}|\subseteq|X|$,
that is the standard points have to be searched in the same Fermat
space $X\in\ECInfty$ from which we have started. For a generic space
$X\in\ECInfty$, that is in general not a space of the form $X=\ext{Y}$,
we do not have an easy way to associate to each point $x\in X$ another
point $s\in X$ making the role of its standard part. Because, on
the contrary, the definition of standard part is a trivial problem
in numerical spaces of the form $\ER^{d}$, the natural idea seems
to use, as it has been done several times in past definitions, observables
like $X\supseteq U\xfreccia{\phi}\ER^{d}$ and to reduce the problem
from the space $X$ to the numerical space $\ER^{d}$. But this idea
naturally leads to the problem of how it is possible to return back
from $\ER^{d}$ to $X$. Unfortunately, this seems solvable only for
spaces $X$ sufficiently similar to manifolds, where charts are invertible
observables (thus not for generic spaces $X$).

Moreover, we also have to consider examples like $X=\left\{ \text{d}t\right\} \subseteq D\setminus\left\{ 0\right\} $,
where it seems natural to expect that $\st{X}=\emptyset$, so the
searched map $x\mapsto\st{x}=s$ in general cannot be defined and
we have to restrict our aim to prove, whether this would be possible,
that for every $x\in X$ there exists at most one $s\in X$ corresponding
to its standard part.

Another idea could be to identify the standard points $s\in X$ as
those points that can be obtained as standard values of figures of
the form $\delta:\ext{U}\freccia X$, i.e. of point of the form $s=\delta(r)$
for $r\in U$. But the case of constant figures having non standard
values, like $\delta(u)=\diff{t}$, represent a counter example to
this intuition.

These are only few examples of unsuccessful attempts that can be tried
if one would like to define a standard part functor. The confirmation
that this is not a trivial goal is given by the following impossibility
results. For their proof we need some preliminary lemmas.

\subsection{Smooth functions with standard values}

The following result state that a function defined on the Fermat reals
and having standard values only, i.e. of the form $f:\ER\freccia\R$,
is necessarily the Fermat extension of its restriction $f|_{\R}$
to the standard points only.
\begin{lem}
\label{lem:functionsWithStandardValues}If $f:\ER\freccia\R$ is smooth
(i.e. it is an arrow of $\ECInfty$), then
\begin{enumerate}
\item \label{enu:restrictionToRealsIsSmooth}$f|_{\R}:\R\freccia\R$ is
smooth in $\CInfty$
\item \label{enu:extensionOfRestrictionToReals}$f=\ext{\left(f|_{\R}\right)}$.
\end{enumerate}
\end{lem}
\noindent \textbf{Proof:} To prove \emph{\ref{enu:restrictionToRealsIsSmooth}}.
we only have to consider the general Theorem \ref{thm:cartesianClosureAndRestrictionOfMaps}
about the restriction of maps. Indeed, since the map $f$ has only
values in $\R$, we have $f(\R)\subseteq\R$ and hence since $f:\ER\freccia\R$
in $\ECInfty$, we have\[
f|_{\R}:(\R\prec\ER)=\bar{\R}\freccia(\R\prec\R)=\bar{\R}\e{in}\ECInfty,\]
from which the conclusion \emph{\ref{enu:restrictionToRealsIsSmooth}}.
follows thanks to Theorem \ref{thm: EmbeddingOfSepInExtendedCn}.

To prove \emph{\ref{enu:extensionOfRestrictionToReals}}. we will
use Theorem \ref{thm:figuresOfFermatSpaces}. In fact, for every $x\in\ER$
we can write $f(y)=\gamma(p,y)$ for every $y\in\mathcal{V}$ in an
open neighborhood $\mathcal{V}$ of $x$. Possibly considering the
composition with a translation, we can suppose $\st{p}=\underbar{0}$
and hence $p\in D_{n}^{\sf p}\subseteq\ER^{\sf p}$ for some order
$n\in\N_{>0}$. Considering the infinitesimal Taylor's formula of
$\gamma$ of order $n$ with respect to the variable $p\in\ER^{\sf p}$,
we have\begin{equation}
f(y)=\gamma(\underbar{0}+p,y)=\sum_{|\alpha|\le n}\frac{p^{\alpha}}{\alpha!}\cdot\partial_{1}^{\alpha}\gamma(\underbar{0},y)\quad\forall y\in\mathcal{V},\label{eq:TaylorWithRespetTo_p}\end{equation}
where $\partial_{1}$ indicates the derivation with respect to the
first slot in $\gamma(-,-)$. But $f(y)\in\R$ and hence $\st{(f(y))}=f(y)$,
so the infinitesimal part of $f(y)$ is zero. From \eqref{eq:TaylorWithRespetTo_p}
we thus obtain\[
\sum_{\substack{|\alpha|\le n\\
\alpha\ne\underbar{0}}
}\frac{p^{\alpha}}{\alpha!}\cdot\partial_{1}^{\alpha}\gamma(\underbar{0},y)=0.\]
Therefore $f(y)=\gamma(\underbar{0},y)$ for every $y\in\mathcal{V}$
and hence $f(x)=\gamma(\underbar{0},x)=\ext{\left[\gamma(\underbar{0},-)\right]}(x)=\ext{f|_{\R}}(x)$.$\qedWithFinalEq$

From this lemma we obtain the following expected result:
\begin{cor}
\label{cor:smoothFromFermatToStarndarAreConstant} If $f:\ER\freccia\R$
is smooth, then $f$ is constant.
\end{cor}
\noindent \textbf{Proof:} From the previous Lemma \ref{lem:functionsWithStandardValues},
if $g:=f|_{\R}$, then $f=\ext{g}$, hence using the derivation formula
with $g$ we have\begin{equation}
\forall x\in\R\:\forall h\in D\pti f(x+h)=f(x)+h\cdot g'(x).\label{eq:derivationFormulaWithFunctionWithStdValues}\end{equation}
But $f(x+h)\in\R$ hence $f(x+h)=\st{f(x+h)}=\st{\left[f(x)+h\cdot g'(x)\right]}=\st{f(x)}=g(x)$,
so from \eqref{eq:derivationFormulaWithFunctionWithStdValues} we
obtain $f(x+h)=g(x)=g(x)+h\cdot g'(x)$ and hence $g'(x)=0$ and so
$g$ is constant because from Lemma \ref{lem:functionsWithStandardValues}
we have that $g:\R\freccia\R$ is smooth.$\qedWithFinalEq$

The most natural example of a function defined on $\ER$ but with
standard values is the standard part map $\st{(-)}:\ER\freccia\R$,
which of course is not constant, so we have the following
\begin{cor}
\label{cor:TheStandardPartIsNotSmooth}The standard part map $\st{(-)}:\ER\freccia\R$
is not smooth.$\qedWithFinalEq$
\end{cor}
As a consequence of this corollary we have that the standard part
functor cannot exists. We will prove this assertion in two ways:
\begin{thm}
\label{thm:StdFunctorDoesNotExist_1}Let $\bar{(-)}:\CInfty_{\text{\emph{sep}}}\freccia\ECInfty$
be the embedding of separated $\CInfty$-spaces into the category
$\ECInfty$ of Fermat spaces (see Section \ref{sec:TheEmbeddingOfManifoldsInFermatSpaces}).
Then, there does not exist a functor\[
\st{(-)}:\ECInfty\freccia\CInfty\]
with the following properties:
\begin{enumerate}
\item There exists a universal arrow of the form $(\eta,\ER):\R\xfreccia{\eta}\st{\!\left(\ER\right)}$.
\item In $\CInfty$ we have the isomorphism $\st{\bar{\R}}\simeq\R$.
\item The functor $\st{(-)}$ preserves terminal objects.
\end{enumerate}
\noindent Therefore, there does not exists a right adjoint of the
Fermat functor that satisfies the isomorphism $\st{\bar{\R}}\simeq\R$
in $\CInfty$ and preserves terminal objects.
\end{thm}
\noindent \textbf{Proof:} We proceed by reduction to the absurd, recalling
(see Appendix \ref{app:someNotionsOfCategoryTheory}) that such a
universal arrow has to verify\[
\ER\in\ECInfty\quad\text{(trivial)}\]
\[
\CInfty\vDash\R\xfrecciad{\eta}\st{\!\left(\ER\right)},\]
and has to be the co-simplest arrow among all arrows satisfying this
property, i.e. for every pair $(\mu,A)$ that verifies\begin{equation}
A\in\ECInfty\label{eq:1UniversalArrowStdFctr}\end{equation}
\begin{equation}
\CInfty\vDash\R\xfrecciad{\mu}\st{A},\label{eq:2UniversalArrowStdFctr}\end{equation}
there exists one and only one arrow $\phi$ such that\begin{equation}
\ECInfty\vDash\ER\xfrecciad{\phi}A\label{eq:phi_arrow1}\end{equation}
\begin{equation}
\CInfty\vDash\xyR{40pt}\xyC{40pt}\xymatrix{\R\ar[r]^{\eta}\ar[dr]_{\mu} & \st{\!\left(\ER\right)}\ar[d]^{\st{\phi}}\\
 & \st{A}}
\label{eq:phi_arrow2}\end{equation}
Let us set $A=\bar{\R}\in\ECInfty$ in \eqref{eq:1UniversalArrowStdFctr}
and \eqref{eq:2UniversalArrowStdFctr} and let $\mu:\R\freccia\st{\bar{\R}}$
be the $\CInfty$-isomorphism of the hypothesis $\st{\bar{\R}}\simeq\R$,
then by \eqref{eq:phi_arrow1} and \eqref{eq:phi_arrow2} we obtain
that $\phi:\ER\freccia\bar{\R}$ in $\ECInfty$ and $\eta\cdot\st{\phi}=\mu$.
From Corollary \ref{cor:smoothFromFermatToStarndarAreConstant} we
have that $\phi$ must be constant, that is there exist a value $r\in\R$
such that\[
\xyR{40pt}\xyC{40pt}\xymatrix{\ER\ar[r]^{\phi}\ar[dr]_{t} & \bar{\R}\\
 & \mathbf{1}\ar[u]_{r}}
\]
where $\mathbf{1}\in\ECInfty$ is the terminal object. Therefore $\phi=t\cdot r$
and hence $\st{\phi}=\st{t}\cdot\st{r}$ since $\st{(-)}$ is supposed
to be a functor. So the map $\st{\phi}$ factors through $\st{\mathbf{1}}$
which, by hypothesis, is the terminal object of $\CInfty$, hence
$\st{\phi}$ is constant too. But this is impossible because $\eta\cdot\st{\phi}=\mu$
and $\mu:\R\freccia\st{\bar{\R}}$ is an isomorphism.$\qedWithFinalEq$

Finally, we want to prove a similar conclusion starting from the equivalence
\eqref{eq:FermatLeftAdjointOfStandardPart}
\begin{thm}
The equivalence \eqref{eq:FermatLeftAdjointOfStandardPart} is false
for $n=\infty$, $X=Y=\R$ and $f=\st{(-)}:\ER\freccia\R$ the standard
part map if $\st{\bar{\R}}=\R$ and $\st{(\st{(-)}})=\st{(-)}$.
\end{thm}
\noindent \textbf{Proof:} Indeed from Corollary \ref{cor:TheStandardPartIsNotSmooth}
we know that the standard part map $f$ is not smooth, that is the
property $\ECInfty\vDash\ER\xfreccia{f}\R$ is false. On the other
hand, we have that for $X=Y=\R$ and $f=\st{(-)}$ the property $\CInfty\vDash X\xfreccia{\st{f}}\st{Y}$
becomes\[
\CInfty\vDash\R\xfrecciad{\st{(\st{(-)}})}\st{\bar{\R}}\]
that is, by the assumed hypotheses\[
\CInfty\vDash\R\xfrecciad{\st{(-)}}\R,\]
which is true because the standard part map is the identity on $\R$.$\qedWithFinalEq$

Analyzing the proofs of these theorems, we can see that the only possibility
to avoid this impossibility result is to change radically the definition
of the category of Fermat spaces $\ECInfty$ so as to include non
constant maps of the form $f:\ER\freccia\R$. This seems possible
thanks to the flexibility of the cartesian closure construction (Chapter
\ref{cha:theCartesianClosure}), but this idea has not been developed
in the present work.

\chapter{\label{cha:logicalPropertiesOfTheFermatFunctor}Logical properties
of the Fermat functor}

In this section we want to investigate some logical properties of
the Fermat functor, with the aim to arrive to a general transfer theorem.
We will see that there are strict connections between the Fermat functor
and intuitionistic logic.

\section{\label{sub:BasicLogicalPropertiesOfFermatFunctor}Basic logical properties
of the Fermat functor}

In this section we will start to investigate some basic logical properties
of the Fermat functor, i.e. the relationships between a given logical
operator (i.e. a propositional connective or a quantifier) and the
related preservation of the Fermat functor of that operator.

The first theorem establishes the relationships between the Fermat
functor and the preservation of implication.
\begin{thm}
\label{thm:FermatFunctorAndImplication}Let $X$, $Y\in\CInfty$ with
$|X|$ is open in $Y$ and such that $X\subseteq Y$ in $\CInfty$
(see Section \ref{sec:CategoricalPropertiesOfTheCartesianClosure}),
then $\ext{X}\subseteq\ext{Y}$ in $\ECInfty$.
\end{thm}
\noindent In other words, the Fermat functor preserves implication
if the antecedent is a property represented by an open set.

\noindent \textbf{Proof:} Let us first assume that $X\subseteq Y$
and recall that $X\subseteq Y$ means $|X|\subseteq|Y|$ and $X=(|X|\prec Y)$,
i.e. the space $X$ has exactly the structure induced by the superspace
$Y$ on one of its subsets. This is equivalent to the following two
properties:\begin{equation}
\forall\delta\pti\delta\In{H}X\then\delta\In{H}Y\label{eq:firstConditionForInclusionOfSubspaces}\end{equation}
\begin{equation}
\forall\delta\pti\delta:|H|\freccia|X|\e{,}\delta\cdot i\In{H}Y\then\delta\In{H}X,\label{eq:secondConditionForInclusionOfSubspaces}\end{equation}
where $i:|X|\hookrightarrow|Y|$ is the inclusion. Using the Fermat
functor we have that $\ext{i}:\ext{X}\freccia\ext{Y}$ in $\ECInfty$.
How does the map $\ext{i}$ act? If, to be more clear, we use the
notation $[x]_{X}:=[(x_{t})_{t}]_{\sim}$ with explicit use of equivalence
classes, then we have\[
\forall x\pti[x]_{X}\in|\ext{X}|\then\ext{i}\left([x]_{X}\right)=[x\cdot i]_{Y}=[x]_{Y}\in|\ext{Y}|,\]
hence $\ext{i}:[x]_{X}\mapsto[x]_{Y}$. We want to prove that this
map is injective. In fact, let us take $[x]_{X}$, $[y]_{X}\in\ext{X}$
such that $[x]_{Y}=[y]_{Y}$ and an observable $\psi:(V\prec X)\freccia K$
defined on the open set $V\in\Top{X}$. From the results of Section
\ref{sec:CategoricalPropertiesOfTheCartesianClosure} it follows that
$(V\prec X)=(V\prec(|X|\prec Y))=(V\prec Y)$ and also that $V\in\Top{Y}$
because, by hypothesis, $|X|$ is open in $Y$. Therefore $VK$ is
a zone of $Y$ too and hence $\psi:(V\prec Y)\freccia K$ is an observable
of $Y$. From the equality $[x]_{Y}=[y]_{Y}$ it follows\[
x_{0}\in V\iff y_{0}\in V\]
\[
x_{0}\in V\then\psi(x_{t})=\psi(y_{t})+o(t)\]
which proves that $[x]_{X}=[y]_{X}$, that is the map $\ext{i}$ is
injective. This injection is exactly the generalization of the identification
that permits to write $\ext{U}\subseteq\ER^{k}$ if $U$ is open in
$\R^{k}$ (see Section \ref{sec: extensions}). For these reasons
we simply write $|\ext{X}|\subseteq|\ext{Y}|$ identifying $|\ext{X}|$
with $\ext{i}(|\ext{X}|)\subseteq|\ext{Y}|$. Now we have to prove
that $\ext{X}\subseteq\ext{Y}$, i.e. $\ext{X}=(|\ext{X}|\prec\ext{Y})$,
i.e. $\ext{(|\ext{X}|X)}=(|\ext{X}|\prec\ext{Y})$ since $\ext{X}=\ext{(|\ext{X}|X)}$
by the Definition \ref{def:FermatFunctor} of Fermat functor. So,
let us first consider a generic figure $\delta\In{S}\ext{X}$ of type
$S\in\SERInfty$; using Theorem \ref{thm:figuresOfFermatSpaces} we
have that for every $s\in S$ there exists an open neighborhood $V=\ext{B}\cap S$
of $s$ in $S$ such that either $\delta|_{V}$ is constant or we
can write $\delta|_{V}=\ext{\gamma}(p,-)|_{V}$ for some $\gamma\in\CInfty(A\times B,X)$.
In the first case, trivially $\delta|_{V}\In{S}(|\ext{X}|\prec\ext{Y})$,
because any space always contains all constant figures. In the second
case, since $i\in\CInfty(X,Y)$ we have $\gamma\cdot i=\gamma\in\CInfty(A\times B,Y)$
and, once again from Theorem \ref{thm:figuresOfFermatSpaces}, we
obtain that $\delta|_{U}\In{U}(|\ext{X}|\prec\ext{Y})$. From the
sheaf property of the space $(|\ext{X}|\prec\ext{Y})$ the conclusion
$\delta\In{S}(|\ext{X}|\prec\ext{Y})$ follows.

\noindent Vice versa if $\delta\In{S}(|\ext{X}|\prec\ext{Y})$, then
$\delta\In{S}(|\ext{X}|\ext{Y})$ by Theorem \ref{thm:spacesPropertiesOfFermatFunctor}
so that, using again Theorem \ref{thm:figuresOfFermatSpaces} and
notations similar to those used above, we have that either $\delta|_{V}$
is constant or $\delta|_{V}=\ext{\gamma}(p,-)|_{V}$, but now with
$\gamma\in\CInfty(A\times B,Y)$. The first case is trivial. For the
second one, it suffices to restrict $\gamma$ so as to obtain a function
with values in $X$ instead of $Y$. But $|X|$ is open in $Y$ so
$\gamma^{-1}(|X|)$ is open in $A\times B$. Thus, we can find $C$
and $D$ open neighborhood of $\st{p}$ and $\st{s}$ respectively
such that $\mu:=\gamma|_{C\times D}\in\CInfty(C\times D,(|X|\prec Y))=\CInfty(C\times D,X)$,
the last equality following from $X\subseteq Y$. Of course $\delta|_{\ext{D}\cap S}=\ext{\mu}(p,-)|_{\ext{D}\cap S}$
and hence $\delta\In{S}\ext{X}$.$\qedNoNewLine$

The following theorem says that the Fermat functor takes open sets
to open sets.
\begin{thm}
\label{thm:FermatFunctorAndOpenSets} If $X\in\CInfty$ and $U$ is
open in $X$, then $\ext{U}$ is open in $\ext{X}$
\end{thm}
\noindent \textbf{Proof:} From the previous theorem we know that $|\ext{U}|\subseteq|\ext{X}|$.
Let us take a figure $d\In{S}\ext{X}$ of type $S\subseteq\ER^{\sf s}$;
to prove that $\ext{U}$ is open in $\ext{X}$ we have to prove that
$d^{-1}(\ext{U})$ is open in $S$, that is we have to prove that
$d^{-1}(\ext{U})$ is generated by sets of the form $\ext{C}\cap S$
for $C$ open in $\R^{\sf s}$. So, let us take a point $s\in d^{-1}(\ext{U})$,
once again from the characterization of the figures of $\ext{X}$
(Theorem \ref{thm:figuresOfFermatSpaces}), we have the existence
of an open neighborhood $V=\ext{B}\cap S$ of $s$ in $S$ such that
either $d|_{V}$ is constant, or we can write $d|_{V}=\ext{\gamma}(p,-)|_{V}$,
for $\gamma\in\CInfty(A\times B,X)$ and $p_{0}\in A$ open in $\R^{\sf p}$.
In the first trivial case we can take $C:=\R^{\sf s}$, so we can
consider the second one only. Because $d(s)\in\ext{U}$, we have that
$\st{d(s)}=\gamma(p_{0},s_{0})\in U$. Since $U$ is open in $X$,
we have that $\gamma^{-1}(U)$ is open in $A\times B$, so from $(p_{0},s_{0})\in\gamma^{-1}(U)$
we get the existence of two open sets $D$ and $C$, respectively
in $A\subseteq\R^{\sf p}$ and $B\subseteq\R^{\sf s}$, such that
$(p_{0},s_{0})\in D\times C\subseteq\gamma^{-1}(U)$. From this we
obtain that $s\in\ext{C}\cap S$, which is the first part of our conclusion.
But $C$ is open in $B$, so $\ext{C}\subseteq\ext{B}$ from the previous
theorem and hence $\ext{C}\cap S\subseteq\ext{B}\cap S=V$, and we
can write $d(x)=\ext{\gamma}(p,x)$ for every $x\in\ext{C}\cap S$.
Therefore $\st{d(x)}=\gamma(p_{0},x_{0})\in U$ because $(p_{0},x_{0})\in D\times C\subseteq\gamma^{-1}(U)$.
From $\st{d(x)}\in U$ we hence get $d(x)\in\ext{U}$ because $U$
is open, and hence we have also proved that $x\in d^{-1}(U)$ for
every $x\in\ext{C}\cap S$, which is the final part of our conclusion.$\qedNoNewLine$

\noindent From this theorem we also obtain the important conclusion
that the Fermat functor preserves open covers, i.e. if $\left(U_{\alpha}\right)_{\alpha\in A}$
is an open cover of $X\in\CInfty$, then $\left(\ext{U}_{\alpha}\right)_{\alpha\in A}$
is an open cover of $\ext{X}$.

The following theorem is the converse of the previous \ref{thm:FermatFunctorAndImplication}
in the case where the spaces are separated.
\begin{thm}
\label{thm:converseOfImplication}In the hypothesis of Theorem \ref{thm:FermatFunctorAndImplication},
if $X$ and $Y$are separated, then $\ext{X}\subseteq\ext{Y}$ in
$\ECInfty$ implies $X\subseteq Y$ in $\CInfty$.
\end{thm}
\noindent \textbf{Proof:} If $\delta\In{U}X$ is figure, then $\ext{\delta}:\ext{U}\freccia\ext{X}$
in $\ECInfty$ and hence $\ext{\delta}\In{\ext{U}}\ext{X}$. but $\ext{X}\subseteq\ext{Y}$,
so $\ext{\delta}\In{\ext{U}}\ext{Y}$. From Corollary \ref{cor:FermatFunctorIsAnEmbedding}
we thus have $\delta\In{U}Y$. It remains to prove condition \eqref{eq:secondConditionForInclusionOfSubspaces}.
If $\delta\cdot i\In{H}Y$, then $\ext{\delta}\cdot\ext{i}\In{\ext{H}}\ext{Y}$.
Recalling that $\ext{X}$ is always identified with $\ext{i}(\ext{X})$,
we can set $j:\ext{i}(\ext{X})\hookrightarrow\ext{Y}$ the inclusion
so that $\ext{\delta}\cdot\ext{i}=\ext{\delta}\cdot\ext{i}\cdot j\In{\ext{H}}\ext{Y}$
and hence $\ext{\delta}\cdot\ext{i}\In{\ext{H}}\ext{i}(\ext{X})$
since $\ext{i}(\ext{X})\subseteq\ext{Y}$. Using the identification\[
\ECInfty\vDash\ext{i}:\xymatrix{\ext{X}\ar[r]\sb(0.42)\sim & \ext{i}(\ext{X})}
\]
this means that $\ext{\delta}\In{\ext{H}}\ext{X}$ and hence $\delta\In{H}X$
from Corollary \ref{cor:FermatFunctorIsAnEmbedding}.$\qedNoNewLine$

From the preservation of the inclusion we can prove that if $X$ is
an open subspace of $Y$, then the operators $(-\prec\ext{X})$ and
$(-\prec\ext{Y})$ conduct to the same subspaces, i.e. we can change
the superspace $Y$ with any other open superspace $X$.
\begin{cor}
\label{cor:noDependencyFromSuperSet}If $X\subseteq Y$ in $\CInfty$,
$|X|$ is open in $Y$ and $Z\subseteq|\ext{X}|$ then $(Z\prec\ext{X})=(Z\prec\ext{Y})$.
\end{cor}
\noindent \textbf{Proof:} This is a trivial consequence of Corollary
\ref{cor:changeOfSuperSpace}. In fact since from Theorem \ref{thm:FermatFunctorAndImplication},
we have that $\ext{X}\subseteq\ext{Y}$ in $\ECInfty$ and hence we
can apply the cited Corollary \ref{cor:changeOfSuperSpace}.$\qedWithFinalEq$

From this result we can prove that the Fermat functor preserves also
counter images of open sets through $\CInfty$ functions.
\begin{thm}
\label{thm:counterImages}Let $f:X\freccia Y$ and $Z\subseteq Y$
in $\CInfty$, with $|Z|$ open in $Y$. Moreover define the spaces
$\ext{f}^{-1}(\ext{Z}):=(\ext{f}^{-1}(|\ext{Z}|)\prec\ext{X})\in\ECInfty$
and $f^{-1}(Z):=(f^{-1}(|Z|)\prec X)\in\CInfty$. Then we have the
equality\[
\ext{[f^{-1}(Z)]}=\ext{f}^{-1}(\ext{Z})\]
 as Fermat spaces.
\end{thm}
\noindent \textbf{Proof:} Let us start from the support sets of the
two spaces:\begin{align*}
x\in\ext{f}^{-1}(|\ext{Z}|)=\left(\ext{f}\right)^{-1}(|\ext{Z}|)\iff & \ext{f}(x)\in|\ext{Z}|\\
\iff & \forall^{0}t\pti f(x_{t})\in|Z|\end{align*}
On the other hand we have\begin{align*}
x\in\ext{\left[f^{-1}(|Z|)\right]} & \iff\forall^{0}t\pti x_{t}\in f^{-1}(|Z|)\\
 & \iff\forall^{0}t\pti f(x_{t})\in|Z|\end{align*}
Hence the support sets are equal. Now we have\begin{equation}
\ext{\left[f^{-1}(Z)\right]}=\ext{\left(\ext{\left[f^{-1}(|Z|)\right]}f^{-1}(Z)\right)}=\left(\ext{\left[f^{-1}(|Z|)\right]}\prec\ext{\left[f^{-1}(Z)\right]}\right),\label{eq:equalitiesForCounterImage}\end{equation}
the first equality following from the Definition \ref{def:FermatFunctor}
of Fermat functor, and the second one from Theorem \ref{thm:spacesPropertiesOfFermatFunctor}.
But $f^{-1}(|Z|)$ is open in $X$ because $|Z|$ is open in $Y$,
hence from the previous Corollary \ref{cor:noDependencyFromSuperSet}
we can change in \eqref{eq:equalitiesForCounterImage} the superspace
$\ext{\left[f^{-1}(Z)\right]}$ with $\ext{X}\supseteq\ext{\left[f^{-1}(Z)\right]}$,
hence\[
\ext{\left[f^{-1}(Z)\right]}=\left(\ext{\left[f^{-1}(|Z|)\right]}\prec\ext{X}\right)=\left(\ext{f}^{-1}(|\ext{Z}|)\prec\ext{X}\right)=\ext{f}^{-1}(\ext{Z}),\]
where we have used the equality of support sets, i.e. \[
\ext{\left[f^{-1}(|Z|)\right]}=\ext{f}^{-1}(|\ext{Z}|)\]
and the definition of the space $\ext{f}^{-1}(\ext{Z})$.$\qedNoNewLine$

Now we consider the relationships between the Fermat functor and the
other propositional connectives.
\begin{thm}
\label{thm:FermatFunctor_AndOrNegation}The Fermat functor preserves
intersections and unions of open sets and the intuitionistic negations,
i.e.
\begin{enumerate}
\item \label{enu:conjunctionAndDisjunction}If $A\subseteq X$ and $B\subseteq X$
in $\CInfty$ and $|A|$, $|B|$ are open in $X$, then\[
\ext{(A\cap_{\sss{X}}B)=\ext{A}\cap_{\sss{{}^{\bullet\!}X}}\ext{B}}\]
and\[
\ext{(A\cup_{\sss{X}}B)}=\ext{A}\cup\ext{B}\]
where, e.g. $A\cup_{\sss{X}}B:=(|A|\cup|B|\prec X)$, $\ext{A}\cap_{\sss{^{\bullet\!}X}}\ext{B}:=(|\ext{A}|\cap|\ext{B}|\prec\ext{X})$,
etc.
\item \label{enu:negation1}If $X\subseteq Y$ in $\CInfty$ and $|X|$
is open in $Y$, then\[
\ext{\left[\text{\emph{int}}_{Y}(Y\setminus X)\right]}\subseteq\text{\emph{int}}_{\ext{Y}}(\ext{Y}\setminus\ext{X}),\]
where $\text{\emph{int}}_{T}(S)$ is the interior of the set $S$
in the topological space $T$.
\item \label{enu:negation2}In the hypotheses of the previous item, if $X$
and $Y$ are separated and the topology of $\ext{Y}$ is generated
by open subsets of the form $\ext{B}$ with $B$ open in $Y$, i.e.
$A=\bigcup\left\{ \ext{B}\subseteq A\,|\, B\in\Top{Y}\right\} $ for
every $A\in\Top{\ext{Y}}$, then\[
\ext{\left[\text{\emph{int}}_{Y}(Y\setminus X)\right]}=\text{\emph{int}}_{\ext{Y}}(\ext{Y}\setminus\ext{X}),\]
i.e. in this case the Fermat functor preserves intuitionistic negations.
\end{enumerate}
\end{thm}
When the topology of a Fermat space of the form $\ext{Y}$ is generated
by open subsets of the form $\ext{B}$ with $B$ open in $Y$, we
will say that \emph{the topology of $\ext{Y}$ is $\ext{(-)}$-generated}.

\noindent \textbf{Proof:}

\noindent \emph{\ref{enu:conjunctionAndDisjunction}}.\quad{}We start
proving that the space $A\cap_{\sss{X}}B$ is the infimum of the spaces
$A$ and $B$ with respect to the partial order of inclusion between
$\CInfty$ spaces. In fact, because of Corollary \ref{cor:noDependencyFromSuperSet}
we have\[
A\cap_{\sss{X}}B=(|A|\cap|B|\prec X)=(|A|\cap|B|\prec A)=(|A|\cap|B|\prec B),\]
that is, $A\cap_{\sss{X}}B\subseteq A$ and $A\cap_{\sss{X}}B\subseteq B$.
Now, let us consider a space $C\in\CInfty$ such that $C\subseteq A$
and $C\subseteq B$, then $|C|\subseteq|A|\cap|B|$ and, e.g., $C=(|C|\prec A)=(|C|\prec(|A|\prec X))=(|C|\prec X)$.
But $|A\cap_{\sss{X}}B|$ is open in $X$ and we can hence apply Corollary
\ref{cor:noDependencyFromSuperSet} again, obtaining $C=(|C|\prec A\cap_{\sss{X}}B)$,
i.e. $C\subseteq A\cap_{\sss{X}}B$. Analogously we can prove that
$A\cup_{\sss{X}}B$ is the supremum of the spaces $A$ and $B$, or
the analogous properties in the category $\ECInfty$.

Therefore, from $A\cap_{\sss{X}}B\subseteq A$ and $A\cap_{\sss{X}}B\subseteq B$
we obtain $\ext{(A\cap_{\sss{X}}B)}\subseteq\ext{A}$ and $\ext{(A\cap_{\sss{X}}B)}\subseteq\ext{B}$
and hence $\ext{(A\cap_{\sss{X}}B)}\subseteq\ext{A}\cap_{\sss{^{\bullet\!}X}}\ext{B}$
because of the greatest lower bound property. Vice versa, if $\delta\In{S}\ext{A}\cap_{\sss{^{\bullet\!}X}}\ext{B}$
is a figure of type $S\subseteq\ER^{\sf s}$, then using the characterization
of the figures of a Fermat space, i.e. Theorem \ref{thm:figuresOfFermatSpaces},
we can say that for every $s\in S$ there exists an open neighborhood
$V=\ext{C}\cap S$ of $s$ such that either $\delta|_{V}$ is constant
or we can write $\delta|_{V}=\ext{\gamma}(p,-)|_{V}$ for $\gamma\in\CInfty(C\times D,X)$.
In the first case $\delta|_{V}\In{V}\ext{(A\cap_{\sss{X}}B)}$; in
the second one $\st{\delta(s)}=\gamma(p_{0},s_{0})\in A\cap_{\sss{X}}B$,
therefore we can find a sufficiently small neighborhood $E\times F$
of $(p_{0},s_{0})$ such that for $U:=\ext{F}\cap S$ we have $\delta|_{U}=\ext{\gamma}(p,-)|_{U}:U\freccia\ext{(A\cap_{\sss{X}}B)}$,
so that $\delta|_{U}\In{U}\ext{(A\cap_{\sss{X}}B)}$. The conclusion
$\delta\In{S}\ext{(A\cap_{\sss{X}}B)}$ follows from the sheaf property
of the space $\ext{(A\cap_{\sss{X}}B)}$. Analogously we can prove
that the Fermat functor preserves unions of $\CInfty$ spaces.

\medskip{}

\noindent \emph{\ref{enu:negation1}}.\quad{}Let us start proving
that $\text{int}_{Y}(Y\setminus X)$ verifies the expected lattice
properties. Being defined as a subspace of $Y$, we have \begin{equation}
\CInfty\vDash\text{int}_{Y}(Y\setminus X)\subseteq Y.\label{eq:setminusInclusionProp}\end{equation}
Moreover, because $|\text{int}_{Y}(Y\setminus X)|\subseteq|Y|\setminus|X|$,
we have that $|X|\cap|\text{int}_{Y}(Y\setminus X)|=\emptyset$, so\begin{equation}
\CInfty\vDash X\cap\text{int}_{Y}(Y\setminus X)=\emptyset.\label{eq:setminusIntersectionProp}\end{equation}
Now, we can prove that among the open subspaces of the space $Y$,
the subspace $\text{int}_{Y}(Y\setminus X)$ is the greatest one verifying
the previous properties \eqref{eq:setminusInclusionProp} and \eqref{eq:setminusIntersectionProp}.
Indeed if $A\in\CInfty$ is open in $Y$, i.e. $|A|\in\Top{Y}$, and
$A\subseteq Y$, $X\cap A=\emptyset$, considering its support set
we have $|A|\subseteq|Y|\setminus|X|$ and hence $|A|\subseteq|\text{int}_{Y}(Y\setminus X)|$
because $|A|$ is open in $Y$. From $A\subseteq Y$, and using Corollary
\ref{cor:noDependencyFromSuperSet} we also get\[
A=(|A|\prec Y)=(|A|\prec\text{int}_{Y}(Y\setminus X)),\]
that is $A\subseteq\text{int}_{Y}(Y\setminus X)$.

Applying the Fermat functor to the properties \eqref{eq:setminusInclusionProp}
and \eqref{eq:setminusIntersectionProp} we obtain $\ext{\left[\text{int}_{Y}(Y\setminus X)\right]}\subseteq\ext{Y}$
and $\ext{X}\cap\ext{\left[\text{int}_{Y}(Y\setminus X)\right]}=\emptyset$,
and hence\[
\ext{\left[\text{int}_{Y}(Y\setminus X)\right]}\subseteq\text{int}_{\ext{Y}}(\ext{Y}\setminus\ext{X})\]

\noindent \emph{\ref{enu:negation2}}.\quad{}To prove the opposite
inclusion, let us take a figure $\delta\In{S}\text{int}_{\ext{Y}}(\ext{Y}\setminus\ext{X})$
of type $S\subseteq\ER^{\sf s}$. Then, for every $s\in S$ we have
$\delta(s)\in\text{int}_{\ext{Y}}(\ext{Y}\setminus\ext{X})$, so that
$\delta(s)\in A\subseteq|\ext{Y}|\setminus|\ext{X}|$, with $A$ open
in $\ext{Y}$. But, by hypothesis, we can find an open set $B\in\Top{Y}$
such that $\delta(s)\in\ext{B}\subseteq A\subseteq|\ext{Y}|\setminus|\ext{X}|$,
and hence $B\subseteq\text{int}_{Y}(Y\setminus X)$ because $B\subseteq\ext{B}$
and $|X|\subseteq|\ext{X}|$ (all the spaces and their subspaces are
separated by hypothesis). Now, we can proceed in the usual way using
the characterization of the figures of a Fermat space (Theorem \ref{thm:figuresOfFermatSpaces}),
from which we get the existence of an open neighborhood $V=\ext{C}\cap S$
of $s$ such that either $\delta|_{V}$ is constant or we can write
$\delta|_{V}=\ext{\gamma}(p,-)|_{V}$ for $\gamma\in\CInfty(C\times D,Y)$.
In the first case $\delta|_{V}\In{V}\ext{\left[\text{int}_{Y}(Y\setminus X)\right]}$;
in the second one $\st{\delta(s)}=\gamma(p_{0},s_{0})\in B$, therefore
we can find a sufficiently small neighborhood $E\times F$ of $(p_{0},s_{0})$
such that for $U:=\ext{F}\cap S$ we have $\delta|_{U}=\ext{\gamma}(p,-)|_{U}:U\freccia\ext{B}$,
so that $\delta|_{U}\In{U}\ext{B}\subseteq\ext{\left[\text{int}_{Y}(Y\setminus X)\right]}$.
The conclusion $\delta\In{S}\ext{\left[\text{int}_{Y}(Y\setminus X)\right]}$
follows from the sheaf property of the space $\ext{\left[\text{int}_{Y}(Y\setminus X)\right]}$.$\qedNoNewLine$
\begin{defn}
\label{def:negation}If $X$, $Y\in\CInfty$ are separated space and
$|X|$ is open in $Y$, then we will use the notation\[
\neg_{\sss{Y}}X:=\left(\text{\emph{int}}_{Y}(Y\setminus X)\prec Y\right)\]
\[
\neg_{\sss{\ext{Y}}}\ext{X}:=\left(\text{\emph{int}}_{\ext{Y}}(\ext{Y}\setminus\ext{X})\prec\ext{Y}\right).\]
Moreover, if $A$, $B$ are open in $Y$, then we also set\[
A\Rightarrow_{\sss{Y}}B:=\neg_{\sss{Y}}A\cup_{\sss{Y}}B.\]

\end{defn}
Therefore, from the previous theorem we can say that\[
\ext{(\neg_{\sss{Y}}X)}=\neg_{\sss{\ext{Y}}}\ext{X}\]
\[
\ext{\left(A\Rightarrow_{\sss{Y}}B\right)}=\left(\ext{A}\Rightarrow_{\sss{\ext{Y}}}\ext{B}\right)\]

\noindent Let us note that the hypotheses of \emph{\ref{enu:negation2}}.
in the previous theorem are surely verified for $X$, $Y$ manifolds.

Finally, we have to consider the relationships between the Fermat
functor and the logical quantifiers.
\begin{defn}
Let $\F$ be a category of types of figures and $f:X\freccia Y$ be
an arrow of the cartesian closure $\bar{\F}$. Then for $Z\subseteq|X|$
we set\begin{equation}
\exists_{f}(Z):=\left(f(Z)\prec Y\right)\label{eq:existentialQuantifier}\end{equation}
\begin{equation}
\forall_{f}(Z):=(\text{\emph{int}}_{Y}\{y\in|Y|\,|\, f^{-1}(\{y\})\subseteq Z\}\prec Y)\label{eq:universalQuantifier}\end{equation}
\end{defn}
\begin{thm}
\label{thm:existentialQuantifier}Let $f:X\freccia Y$ be a $\CInfty$-map.
Moreover, let us suppose that
\begin{enumerate}
\item \noindent $Z$ is open in $X$,
\item \noindent $f$ is open with respect to the topologies $\Top{X}$ and
$\Top{Y}$,
\item $f|_{Z}:(Z\prec X)\freccia(f(Z)\prec Y)$ has a left%
\footnote{Let us note that here the word {}``left'' is with respect to the
composition of functions represented by the symbol $\left(f\cdot g\right)(x)=g(f(x))$
(that permits an easier reading of diagrams), so that it corresponds
to {}``right'' with respect to the notation with $\left(f\circ g\right)(x)=f(g(x))$.%
} inverse in $\CInfty$,
\item \noindent $X$, $Y$ are separated.
\end{enumerate}
Then we have\[
\ext{(\exists_{f}(Z))}=\existsext{f}(\ext{Z})\]

\noindent i.e., in these hypotheses, the Fermat functor preserves
existential quantifiers.

\end{thm}
\smallskip{}

\begin{thm}
\label{thm:universalQuantfier}Let $f:X\freccia Y$ be a $\CInfty$-map.
Moreover, let us suppose that
\begin{enumerate}
\item \noindent $Z$ is open in $X$,
\item the topology of $\ext{Y}$ is $\ext{(-)}$-generated,
\item \noindent $X$, $Y$ are separated.
\end{enumerate}
Then we have\[
\ext{(\forall(Z))}=\forallext{f}(\ext{Z})\]

\noindent i.e., in these hypotheses, the Fermat functor preserves
existential quantifiers.

\end{thm}
\noindent To motivate the definitions \eqref{eq:existentialQuantifier}
and \eqref{eq:universalQuantifier} we can consider as $f$ a projection
$p:A\times B\freccia B$ of a product, then for $Z\subseteq|A\times B|$
we have\[
|\exists_{p}(Z)|=p(Z)=\{b\,|\,\exists x\in Z:\ b=p(x)\}=\{b\in B\,|\,\exists a\in A:\ Z(a,b)\},\]
where we used $Z(a,b)$ for $(a,b)\in Z$. This justifies the definition
of $\exists_{f}$ as a generalization of this $\exists_{p}$.

Taking the difference $|Y|\setminus|\exists_{f}(Z)|$ we obtain\begin{align*}
|Y|\setminus|\exists_{f}(Z)| & =|Y|\setminus f(Z)=\{y\in Y\,|\,\neg(\exists x\in Z:\ y=f(x))\}=\\
 & =\{y\,|\,\forall x\in X:\ y=f(x)\Rightarrow x\notin Z\}=\\
 & =\{y\,|\, f^{-1}(\{y\})\subseteq X\setminus Z\}=\\
 & =|\forall_{f}(X\setminus Z)|\end{align*}
This justifies fully the definition of $\forall_{f}$ in the case
of classical logic. For example, in the case of a projection $p:A\times B\freccia B$,
for $Z\subseteq|A\times B|$ we have\[
|\forall_{p}(Z)|=\{b\,|\,\forall x\in X:\ b=p(x)\Rightarrow x\in Z\}=\{b\in B\,|\,\forall a\in A:\ Z(a,b)\}\]

\noindent In an intuitionistic context %
\footnote{We recall that in intuitionistic logic a quantifier cannot be defined
starting from the other one; the best result that it is possible to
obtain is that $[\forall x:\,\neg\phi(x)]\iff[\neg\exists x:\,\phi(x)]$,
where it is important to recall that, in general, $\neg\neg\phi(x)$
is not equivalent to $\phi(x)$ in intuitionistic logic (as it can
be guessed using topological considerations, because of the interior
operator, starting from our Definition \ref{def:negation} of negation).%
} the interpretation of a formula in a topological space must always
result in an open set (We recall that like the classical logic can
be interpreted in any boolean algebra of generic subsets of a given
superset, the intuitionistic logic can be interpreted in the Heyting
algebra of the open sets of any topological space (see e.g. \citet{Ra-Si,Sco})
and this motivates the use of the interior operator $\text{int}_{Y}$
in the definition \eqref{eq:universalQuantifier}. Finally, we recall
that the projection of a product is always an open map if on the product
space $A\times B$ we have the product topology, like in our context
if $A$ and B are manifolds (see Section \ref{sec:FermatFunctorAndProductOfManifolds},
\eqref{eq:productTopologyAndProductSpaces} and the final discussion
in Section \ref{sec:CategoricalPropertiesOfTheCartesianClosure}).
Moreover if $a\in A$, then $g:b\in B\freccia(a,b)\in A\times B$
is a left inverse of class $\CInfty$ of the projection $p$, so the
map $p$ verifies all the hypotheses of Theorem \ref{thm:existentialQuantifier}.

To prove this theorem we need the following two lemmas, which repeat
in our context well known results (see e.g. \citet{Tay}).
\begin{lem}
\label{lem:existIsLeftAdjoint}If $\F$ is a category of types of
figures, and $f:X\freccia Y$ in $\bar{\F}$, then we have:
\begin{enumerate}
\item If $A$, $A'$ are subspaces of $X$ (not necessarily open) with $A\subseteq A'$,
then $\exists_{f}(A)\subseteq\exists_{f}(A')$.
\item If $A\subseteq X$ and $B\subseteq Y$, then in the category $\bar{\F}$
we have the equivalence\begin{equation}
\cfrac{A\subseteq f^{-1}(B)}{\exists_{f}(A)\subseteq B}\label{eq:existsLeftAdjoint}\end{equation}
that is $\exists_{f}\dashv f^{-1}$ with respect to the order relation
$\subseteq$ between subspaces.
\end{enumerate}
\end{lem}
\smallskip{}

\begin{lem}
\label{lem:forallIsRightAdjoint}If $\F$ is a category of types of
figures, and $f:X\freccia Y$ in $\bar{\F}$, then we have:
\begin{enumerate}
\item If $A$, $A'$ are subspaces of $X$ (not necessarily open) with $A\subseteq A'$,
then $\forall_{f}(A)\subseteq\forall_{f}(A')$.
\item if $A\subseteq X$ and $B\subseteq Y$, then in the category $\bar{\F}$
we have the equivalence\begin{equation}
\cfrac{f^{-1}(B)\subseteq A}{B\subseteq\forall_{f}(A)}\label{eq:forallRightAdjoint}\end{equation}
that is $f^{-1}\dashv\forall_{f}$ with respect to the order relation
$\subseteq$ between subspaces.
\end{enumerate}
\end{lem}
\smallskip{}

\begin{lem}
\label{lem:preservationOfRestriction}If $f:X\freccia Y$ in $\CInfty$
and $Z\subseteq|X|$, then\[
\ext{\left(f|_{Z}\right)}=\ext{f}|_{\ext{Z}}\]

\end{lem}
\noindent \textbf{Proof:} Both the functions are defined in $\ext{Z}=\ext{(Z\prec X)}$,
so let $x\in\ext{Z}$, we have $\ext{\left(f|_{Z}\right)}(x)=\left(f(x_{t})\right)_{t\ge0}=\left(\ext{f}|_{\ext{Z}}\right)(x)$.$\qedNoNewLine$

\noindent \textbf{Proof of Lemmas \ref{lem:existIsLeftAdjoint} and
\ref{lem:forallIsRightAdjoint}:} let us assume that $A$ and $A'$
are subspaces of $X$ with $A\subseteq A'$. We recall that $\exists_{f}(A)=(f(A)\prec Y)$
and $\exists_{f}(A')=(f(A')\prec Y)$; but $|f(A)|\subseteq|f(A')|$,
we can hence apply Corollary \ref{cor:changeOfSuperSpace} to change
in $\exists_{f}(A)$ the superspace $Y$ with the superspace $(f(A')\prec Y)\subseteq Y$,
obtaining\[
\exists_{f}(A)=(f(A)\prec(f(A')\prec Y))=(f(A)\prec\exists_{f}(A')),\]
that is $\exists_{f}(A)\subseteq\exists_{f}(A')$.

Now let us assume that $A\subseteq f^{-1}(B)$, then $|f(A)|\subseteq B$
as sets so that, applying once again Corollary \ref{cor:changeOfSuperSpace}
we can change the superspace $Y$ in $\exists_{f}(A)=(f(A)\prec Y)$
with the superspace $B$ obtaining $\exists_{f}(A)=(f(A)\prec B)$,
that is the conclusion $\exists_{f}(A)\subseteq B$. Reversing this
deduction we can obtain a proof for the opposite implication. In a
similar way we can also prove the analogous properties of the universal
quantifier.$\qedNoNewLine$

\noindent \textbf{Proof of Theorem \ref{thm:existentialQuantifier}:}
The first idea is to use the uniqueness of the adjoints of $f^{-1}$,
that is the property that the spaces $\exists_{f}(A)$ and $\forall_{f}(A)$
are uniquely determined by the equivalences \eqref{eq:existsLeftAdjoint}
and \eqref{eq:forallRightAdjoint} respectively, and to use the preservation
of the relation $X\subseteq Y$ by the Fermat functor. Indeed, if
we suppose that $\ext{Z}\subseteq\ext{f}^{-1}(\ext{W})$, then we
also have $\ext{Z}\subseteq\ext{\left(f^{-1}(W)\right)}$ by the preservation
of counter images. By Theorem \ref{thm:converseOfImplication} this
implies $Z\subseteq f^{-1}(W)$ and hence $\exists_{f}(Z)\subseteq W$
by Lemma \ref{lem:existIsLeftAdjoint} and so $\ext{\left(\exists_{f}(Z)\right)}\subseteq\ext{W}$
applying the preservation of implications. At the same time, the hypothesis
$\ext{Z}\subseteq\ext{f}^{-1}(\ext{W})$ implies $\existsext{f}(\ext{Z})\subseteq\ext{W}$
since Lemma \ref{lem:existIsLeftAdjoint} is true for the category
$\ECInfty$ too. All these implications can be reversed in a direct
way using Theorem \ref{thm:converseOfImplication} and our hypothesis
that the spaces $X$ and $Y$ (and hence all their subspaces) are
separated. Therefore, we have the equivalences\begin{equation}
\cfrac{\ext{Z}\subseteq\ext{f}^{-1}(\ext{W})}{\cfrac{\existsext{f}(\ext{Z})\subseteq\ext{W}}{\ext{\left(\exists_{f}(Z)\right)}\subseteq\ext{W}}}\label{eq:3equivalencesExists}\end{equation}
In them, if we set $W:=\exists_{f}(Z)$, then the third one is trivially
true, and from the second one we obtain\begin{equation}
\existsext{f}(\ext{Z})\subseteq\ext{\left(\exists_{f}(Z)\right)}.\label{eq:existsFermatInclusion}\end{equation}
This part of the deduction cannot be reversed because, e.g., in \eqref{eq:3equivalencesExists}
instead of a generic subspace of $\ext{Y}$ we have a subspace of
the form $\ext{W}$ only. So, let us first recall that\[
\ext{(\exists_{f}(Z))}=\ext{(f(Z)\prec Y)}\]
\[
\existsext{f}(\ext{Z})=(\ext{f}(\ext{Z})\prec\ext{Y})\]
To prove the opposite relations of \eqref{eq:existsFermatInclusion}
we need to assume the existence of a left inverse $g$ of the restriction
$f|_{Z}$, i.e. a $\CInfty$-map $g:(f(Z)\prec Y)\freccia(Z\prec X)$
such that $g\cdot f|_{Z}=1_{f(Z)}$. Let us take a figure $\delta\In{S}\ext{(\exists_{f}(Z))}$
of type $S\subseteq\ER^{\sf s}$. Then\[
\xyR{40pt}\xyC{40pt}\xymatrix{\ECInfty\vDash S\ar[r]\sp(0.45)\delta & \ext{(f(Z)\prec Y)}\ar[r]\sp(0.53){\ext{g}} & \ext{(Z\prec X)}}
,\]
and hence $\delta\cdot\ext{g}\In{S}\ext{(Z\prec X)}$. Composing this
map with the restriction $\ext{(f|_{Z})}=\ext{f}|_{\ext{Z}}:\ext{(Z\prec X)}\freccia(\ext{f}(\ext{Z})\prec\ext{Y})$
we obtain\[
\delta\cdot\ext{g}\cdot\ext{f}|_{\ext{Z}}=\delta\cdot\ext{(g\cdot f|_{Z})}=\delta\In{S}(\ext{f}(\ext{Z})\prec\ext{Y})=\existsext{f}(\ext{Z}).\]
We have hence proved the first condition \eqref{eq:firstConditionForInclusionOfSubspaces}
to prove that $\ext{\left(\exists_{f}(Z)\right)}\subseteq\existsext{f}(\ext{Z})$.
This part of the deduction also proves that we have the relation $|\ext{(\exists_{f}(Z))}|\subseteq|\existsext{f}(\ext{Z})|$
between the corresponding support sets. Hence we can now prove the
second condition \eqref{eq:secondConditionForInclusionOfSubspaces};
let us consider a map $\delta:S\freccia|\ext{(\exists_{f}(Z))}|$
such that $\delta\cdot i\In{S}\existsext{f}(\ext{Z})$, where $i:|\ext{(\exists_{f}(Z))}|\hookrightarrow|\existsext{f}(\ext{Z})|$
is the inclusion map. So we have $\delta\cdot i=\delta\In{S}\existsext{f}(\ext{Z})$
and hence also $\delta\In{S}\ext{\left(\exists_{f}(Z)\right)}$ since
\eqref{eq:existsFermatInclusion}. This easily proves also the second
condition \eqref{eq:secondConditionForInclusionOfSubspaces} and hence
$\existsext{f}(\ext{Z})=\ext{\left(\exists_{f}(Z)\right)}$.$\qedNoNewLine$

\noindent \textbf{Proof of Theorem \ref{thm:universalQuantfier}:}
Analogously to how we did in the previous proof, we can proceed for
the universal quantifier obtaining the equivalences\begin{equation}
\cfrac{\ext{f}^{-1}(\ext{W})\subseteq\ext{Z}}{\cfrac{\ext{W}\subseteq\forallext{f}(\ext{Z})}{\ext{W}\subseteq\ext{\left(\forall_{f}(Z)\right)}}}\label{eq:3equivalencesForall}\end{equation}

\noindent from which we obtain\begin{equation}
\ext{\left(\forall_{f}(Z)\right)}\subseteq\forallext{f}(\ext{Z}).\label{eq:forallFermatInclusion}\end{equation}

\noindent Now, let us consider the opposite inclusion, recalling that\[
\ext{(\forall_{f}(Z))}=\ext{(\left[\text{int}_{Y}\left\{ y\,|\, f^{-1}\{y\}\subseteq Z\right\} \right]\prec Y)}\]
\[
\forallext{f}(\ext{Z})=(\text{int}_{\ext{Y}}\left\{ y\,|\,\ext{f}^{-1}\{y\}\subseteq\ext{Z}\right\} \prec\ext{Y}).\]
So let us consider a figure $\delta\In{S}\forallext{f}(\ext{Z})$
and a point $s\in S$, then\[
\delta(s)\in(\text{int}_{\ext{Y}}\left\{ y\,|\,\ext{f}^{-1}\{y\}\subseteq\ext{Z}\right\} \prec\ext{Y}).\]
Because, by hypothesis, the topology of $\ext{Y}$ is generated by
open sets of the form $\ext{U}$, $U\in\Top{X}$, by the definition
of interior we obtain\begin{equation}
\exists U\in\Top{X}:\ \delta(s)\in\ext{U}\subseteq\left\{ y\,|\,\ext{f}^{-1}\{y\}\subseteq\ext{Z}\right\} .\label{eq:forallInterior}\end{equation}
It is natural to expect that the property $\ext{f}^{-1}\{y\}\subseteq\ext{Z}$
can be extended to the whole set $\ext{U}$, indeed\begin{align*}
\forall x\in\ext{f}^{-1}(\ext{U}):\  & \ext{f}(x)\in\ext{U}\\
 & \ext{f}^{-1}\{\ext{f}x\}\subseteq\ext{Z}\quad\text{ by }\eqref{eq:forallInterior}\\
 & \text{but }x\in\ext{f}^{-1}\{\ext{f}x\}\\
 & \text{hence }x\in\ext{Z}.\end{align*}
Therefore we have $\ext{f}^{-1}(\ext{U})\subseteq\ext{Z}$, that is
$\ext{\left(f^{-1}(U)\right)}\subseteq\ext{Z}$, and hence $f^{-1}(U)\subseteq Z$
because we are considering separated spaces, and so $U\subseteq\forall_{f}(Z)$.
But $\delta(s)\in\ext{U}$, and setting $V:=\delta^{-1}(\ext{U})$
we obtain an open neighborhood of $s$ such that\[
\delta|_{V}:V\freccia\ext{U}.\]
Therefore $\delta|_{V}\In{V}\ext{U}\subseteq\ext{\left(\forall_{f}(Z)\right)}$.
The conclusion $\delta\In{S}\ext{\left(\forall_{f}(Z)\right)}$ follows
from the sheaf property of the space $\ext{\left(\forall_{f}(Z)\right)}$.
The second condition \eqref{eq:secondConditionForInclusionOfSubspaces}
can be proved analogously to what we already did above for the existential
quantifier.$\qedNoNewLine$

\section{The general transfer theorem}

In this section, for simplicity of notations, every arrow $f$ of
the categories $\CInfty$ and $\ECInfty$ is supposed to have unique
domain and codomain (they will be denoted by $\text{dom}(f)$ and
$\text{cod}(f)$ respectively; see Appendix \ref{app:someNotionsOfCategoryTheory}
for more details about this hypothesis, which at a first reading may
seem trivial).

\medskip{}

In the previous section, it has been underlined that the logical operators
defined above, like $A\cap_{\sss{Y}}B$ or $\neg_{\sss{Y}}B$, or
$\forall_{f}(A)$ take subspaces of a given space $Y$ to subspaces
of the same or of another space (like e.g. $f^{-1}(A)\subseteq X$
if $A\subseteq Y$). Therefore, we have now the possibility to compose
these operators to construct new spaces, like e.g. the following\begin{equation}
S:=\forall_{\eps}\left(A\Rightarrow_{\sss{Z}}\exists_{\delta}\left(B\cap_{\sss{Y}}\forall_{x}\left(C\Rightarrow_{\sss{X}}D\right)\right)\right)\in\CInfty,\label{eq:exampleLogicalSpaces}\end{equation}
where e.g.\begin{gather}
X\xfrecciad{x}Y\xfrecciad{\delta}Z\xfrecciad{\eps}W\label{eq:cond1}\\
C,D\subseteq X\label{eq:cond2}\\
B\subseteq Y\label{eq:cond3}\\
A\subseteq Z\label{eq:cond4}\end{gather}
In this section, we want to
\begin{enumerate}
\item define the family of formulae, like that used in \eqref{eq:exampleLogicalSpaces}
to define $S$, that permit to define spaces in $\CInfty$ or in $\ECInfty$
by means of logical operators;
\item show that to each formula there corresponds a suitable operator that
maps subspaces of $\CInfty$ into new subspaces of the same category;
\item define a $\ext{(-)}$-transform $\ext{\xi}$ of a formula $\xi$,
called the \emph{Fermat transform of }$\xi$. To the Fermat transform
$\ext{\xi}$ corresponds an operator acting on spaces of the category
$\ECInfty$;
\item find a way to associate to every formula $\phi$, a set of conditions
like \eqref{eq:cond1}, \eqref{eq:cond2}, \eqref{eq:cond3}, \eqref{eq:cond4}
and other suitable hypotheses that will permit to apply all the theorems
of the previous Section \ref{sub:BasicLogicalPropertiesOfFermatFunctor}.
Indeed, in the general transfer theorem we have to assume on superspaces,
subspaces and maps, all the hypotheses of the theorems of the previous
section, if we want that the Fermat functor preserves all the logical
operations;
\item prove that the operator corresponding to $\ext{\xi}$ is the Fermat
transform of the operator corresponding to the formula $\xi$, that
is the general transfer theorem.
\end{enumerate}
We will also include, in our formulae, the symbol of product because
in case of manifolds the Fermat functor preserves also this operation
(see Theorem \ref{thm:FermatPreservesProductOfMan}).
\begin{defn}
\label{def:formulasForSpaces}Let\[
\mathcal{S}:=\left\{ \corners{\times},\corners{\neg},\corners{\Rightarrow},\corners{\cap},\corners{\cup},\corners{\exists},\corners{\forall},\corners{^{-1}},\corners{(},\corners{)}\right\} \]
be a set of distinct elements called \emph{symbols}. An \emph{expression
in $\CInfty$} is a finite sequence of symbols in $\mathcal{S}$,
objects or arrows of $\CInfty$. Sequences of length 0 are admitted,
but those of length 1 are identified with the element itself. For
example the following\[
\ulcorner\neg_{\sss{Y}}A\urcorner:=(\ulcorner\neg\urcorner,Y,A)\]
\[
\ulcorner\exists_{f}(A)\urcorner:=(\ulcorner\exists\urcorner,f,\ulcorner(\urcorner,A,\ulcorner)\urcorner)\]
$ $are examples of expressions. We will use similar abbreviations
for other expressions like, e.g., $\ulcorner A\Rightarrow_{\sss{Y}}B\urcorner:=\left(A,\ulcorner\Rightarrow\urcorner,Y,B\right)$.

\noindent If $\phi$ and $\psi$ are expressions, then with the symbol
$\corners{(\phi\cap_{\chi}\psi)}$ we mean the n-tuple $(\corners{(},\phi,\corners{\cap},\chi,\psi,\corners{)})$.
We will use similar notations to construct expressions, like e.g.\[
\corners{\exists_{f}(\phi)}:=(\corners{\exists},f,\corners{(},\phi,\corners{)}).\]

\noindent We will denote with $\mathcal{L}^{+}(\CInfty)$ the intersection
of all the classes $L$ of expressions verifying
\begin{enumerate}
\item If $A\in\CInfty$, then $A\in L$
\item If $\phi$, $\chi$, $\psi\in L$, then\begin{equation}
\corners{(\phi\times\psi)}\ ,\ \corners{\neg_{\sss{\chi}}\phi}\ ,\ \corners{(\phi\Rightarrow_{\sss{\chi}}\psi)}\ ,\ \corners{(\phi\cap_{\sss{\chi}}\psi)}\ ,\ \corners{(\phi\cup_{\sss{\chi}}\psi)}\in L\label{eq:operationsForFormulas}\end{equation}

\item If $f$ is an arrow of $\CInfty$ and $\phi\in L$, then\[
\ulcorner\exists_{f}(\phi)\urcorner\e{,}\ulcorner\forall_{f}(\phi)\urcorner\e{,}\ulcorner f^{-1}(\phi)\urcorner\in L\]

\end{enumerate}
An analogous definition can be stated in the category $\ECInfty$,
and the related class of expressions will be denoted by $\mathcal{L}^{+}(\ECInfty$).

\end{defn}
\noindent As usual, see e.g. \citet{Mon}, we can prove the following
\begin{thm}
\label{thm:uniqueReadabilityForLPlus}If $\xi\in\mathcal{L}^{+}(\CInfty)$,
then one and only one of the following holds:
\begin{enumerate}
\item $\xi=A$ for some object $A\in\CInfty$ (expression of length 1);
\item $\xi=\corners{(\phi\times\psi)}$ for some $\phi$, $\psi\in\mathcal{L}^{+}(\CInfty)$;
\item $\xi=\corners{\neg_{\sss{\chi}}\psi}$ for some $\chi$, $\psi\in\mathcal{L}^{+}(\CInfty)$;
\item $\xi=\corners{(\phi\Rightarrow_{\chi}\psi)}$ for some $\phi,$ $\chi$,
$\psi\in\mathcal{L}^{+}(\CInfty)$;
\item $\xi=\corners{(\phi\cap_{\sss{\chi}}\psi)}$ for some $\phi,$ $\chi$,
$\psi\in\mathcal{L}^{+}(\CInfty)$;
\item $\xi=\corners{(\phi\cup_{\sss{\chi}}\psi)}$ for some $\phi,$ $\chi$,
$\psi\in\mathcal{L}^{+}(\CInfty)$;
\item $\xi=\corners{\exists_{f}(\phi)}$ for some $\phi\in\mathcal{L}^{+}(\CInfty)$
and some arrow $f$ of $\CInfty$;
\item $\xi=\corners{\forall_{f}(\phi)}$ for some $\phi\in\mathcal{L}^{+}(\CInfty)$
and some arrow $f$ of $\CInfty$;
\item $\xi=\corners{f^{-1}(\phi)}$ for some $\phi\in\mathcal{L}^{+}(\CInfty)$
and some arrow $f$ of $\CInfty$.
\end{enumerate}
Moreover, the expressions $\phi$, $\psi$, $\chi$, the object $A$
and the arrow $f$ asserted to exist are uniquely determined by $\xi$.

\end{thm}
\noindent Actually, the expressions of $\mathcal{L}^{+}(\CInfty)$
are not well formed formulae because we can consider in the set $\mathcal{L}^{+}(\CInfty)$
expressions like $\ulcorner A\cap_{\sss{X}}B\urcorner$, but with
$A$ and $B$ that are not subspaces of $X\in\CInfty$. Analogously,
an expression of the form $\ulcorner\exists_{f}(A)\urcorner$ is a
formula only if $f:X\freccia Y$ and $A\subseteq X$. This means that
we are dealing with a typed language and, e.g., the previous $\ulcorner\exists_{f}(A)\urcorner$
is a formula only if $A$ is of the form {}``subsets of the domain
of $f$''. In the following definition we will define what is this
type.
\begin{defn}
\label{def:typeOfAnExpression}If $\xi\in\mathcal{L}^{+}(\CInfty)$,
then the type $\tau(\xi)$ is defined recursively by the following
conditions:
\begin{enumerate}
\item If $\xi=A$ for some object $A\in\CInfty$, then $\tau(\xi):=A$.
\item If $\xi=\corners{(\phi\times\psi)}$ for some $\phi$, $\psi\in\mathcal{L}^{+}(\CInfty)$,
then $\tau(\xi):=\tau(\phi)\times\tau(\psi)$.
\item If $\xi=\corners{\neg_{\sss{\chi}}\psi}$ for some $\chi$, $\psi\in\mathcal{L}^{+}(\CInfty)$,
and if $\tau(\psi)\subseteq\tau(\chi)$ in $\CInfty$, then \[
\tau(\xi):=\neg_{\sss{\tau(\chi)}}\tau(\psi)\]

\item If $\xi=\corners{(\phi\Rightarrow_{\chi}\psi)}$ for some $\phi,$
$\chi$, $\psi\in\mathcal{L}^{+}(\CInfty)$, and if $\tau(\phi)\subseteq\tau(\chi)$
and $\tau(\psi)\subseteq\tau(\chi)$, then \[
\tau(\xi):=\tau(\phi)\Rightarrow_{\sss{\tau(\chi)}}\tau(\psi)\]

\item If $\xi=\corners{(\phi\cap_{\sss{\chi}}\psi)}$ for some $\phi,$
$\chi$, $\psi\in\mathcal{L}^{+}(\CInfty)$, and if $\tau(\phi)\subseteq\tau(\chi)$
and $\tau(\psi)\subseteq\tau(\chi)$, then \[
\tau(\xi):=\tau(\phi)\cap_{\sss{\tau(\chi)}}\tau(\psi)\]

\item If $\xi=\corners{(\phi\cup_{\sss{\chi}}\psi)}$ for some $\phi,$
$\chi$, $\psi\in\mathcal{L}^{+}(\CInfty)$, and if $\tau(\phi)\subseteq\tau(\chi)$
and $\tau(\psi)\subseteq\tau(\chi)$, then \[
\tau(\xi):=\tau(\phi)\cup_{\sss{\tau(\chi)}}\tau(\psi)\]

\item If $\xi=\corners{\exists_{f}(\phi)}$ for some $\phi\in\mathcal{L}^{+}(\CInfty)$
and some arrow $f:X\freccia Y$ of $\CInfty$, and if $\tau(\phi)\subseteq X$,
then \[
\tau(\xi):=\exists_{f}(\tau(\phi))\]

\item If $\xi=\corners{\forall_{f}(\phi)}$ for some $\phi\in\mathcal{L}^{+}(\CInfty)$
and some arrow $f:X\freccia Y$ of $\CInfty$, and if $\tau(\phi)\subseteq X$,
then \[
\tau(\xi):=\forall_{f}(\tau(\phi))\]

\item If $\xi=\corners{f^{-1}(\phi)}$ for some $\phi\in\mathcal{L}^{+}(\CInfty)$
and some arrow $f:X\freccia Y$ of $\CInfty$, and if $\tau(\phi)\subseteq Y$,
then \[
\tau(\xi):=f^{-1}(\tau(\phi))\]

\end{enumerate}
In all the other cases the type $\tau(\xi)$ is not defined. Analogously
we can define $\ext{\tau}(\xi)$, the type of expressions $\xi\in\mathcal{L}^{+}(\ECInfty)$
in the category of Fermat spaces.

\end{defn}
\noindent Let us note that e.g. when we say {}``If $\xi=\corners{(\phi\times\psi)}$
for some $\phi$, $\psi\in\mathcal{L}^{+}(\CInfty)$, then $\tau(\xi):=\tau(\phi)\times\tau(\psi)$'',
we implicitly mean {}``If $\xi=\corners{(\phi\times\psi)}$ for some
$\phi$, $\psi\in\mathcal{L}^{+}(\CInfty)$, then $\tau(\phi)$ \emph{and
}$\tau(\psi)$ \emph{are defined and }$\tau(\xi):=\tau(\phi)\times\tau(\psi)$''.

\noindent Now we can define the formulae of $\CInfty$ as the expressions
$\xi$ in $\mathcal{L}^{+}(\CInfty)$ for which the type $\tau(\xi)$
is defined:
\begin{defn}
\label{def:formulas}The set $\mathcal{L}(\CInfty)$ of \emph{formulae
in }$\CInfty$ is defined recursively by the following condition:
$\xi\in\mathcal{L}(\CInfty)$ if and only if one of the following
alternatives is true:
\begin{enumerate}
\item $\xi=A$ for some object $A\in\CInfty$;
\item $\xi=\corners{(\phi\times\psi)}$ for some $\phi$, $\psi\in\mathcal{L}(\CInfty)$;
\item If $\xi=\corners{\neg_{\sss{\chi}}\psi}$ for some $\chi$, $\psi\in\mathcal{L}(\CInfty)$,
then $\tau(\psi)\subseteq\tau(\chi)$ in $\CInfty$;
\item If $\xi=\corners{(\phi\Rightarrow_{\chi}\psi)}$ for some $\phi,$
$\chi$, $\psi\in\mathcal{L}(\CInfty)$, then $\tau(\phi)\subseteq\tau(\chi)$
and $\tau(\psi)\subseteq\tau(\chi)$;
\item If $\xi=\corners{(\phi\cap_{\sss{\chi}}\psi)}$ for some $\phi,$
$\chi$, $\psi\in\mathcal{L}(\CInfty)$, then $\tau(\phi)\subseteq\tau(\chi)$
and $\tau(\psi)\subseteq\tau(\chi)$;
\item If $\xi=\corners{(\phi\cup_{\sss{\chi}}\psi)}$ for some $\phi,$
$\chi$, $\psi\in\mathcal{L}(\CInfty)$, then $\tau(\phi)\subseteq\tau(\chi)$
and $\tau(\psi)\subseteq\tau(\chi)$;
\item If $\xi=\corners{\exists_{f}(\phi)}$ for some $\phi\in\mathcal{L}(\CInfty)$
and some arrow $f:X\freccia Y$ of $\CInfty$, then $\tau(\phi)\subseteq X$;
\item If $\xi=\corners{\forall_{f}(\phi)}$ for some $\phi\in\mathcal{L}(\CInfty)$
and some arrow $f:X\freccia Y$ of $\CInfty$, then $\tau(\phi)\subseteq X$;
\item If $\xi=\corners{f^{-1}(\phi)}$ for some $\phi\in\mathcal{L}(\CInfty)$
and some arrow $f:X\freccia Y$ of $\CInfty$, then $\tau(\phi)\subseteq Y$.
\end{enumerate}
\end{defn}
\noindent Therefore, if $\xi\in\mathcal{L}(\CInfty)$ is a formula,
then the type $\tau(\xi)$ is defined, and hence the type $\tau$
is an application\[
\tau:\mathcal{L}(\CInfty)\freccia\text{Obj}(\CInfty),\]
where $\text{Obj}(\CInfty)$ is the class of all the objects of the
category $\CInfty$. An analogous property can be stated for $\ECInfty$.
As usual, we can say that $\phi$ is a \emph{subformula} of $\xi$
if both $\xi$ and $\phi$ are formulae and $\xi=(\chi,\phi,\psi)$
for some expressions $\chi$ and $\psi$.

The condition that the type $\tau(\xi)$ is defined is exactly the
minimal condition for the formula $\xi$ of being meaningful. E.g.
for the formula\begin{equation}
\xi:=\corners{\forall_{\eps}\left(A\Rightarrow_{\sss{Z}}\exists_{\delta}\left(B\cap_{\sss{Y}}\forall_{x}\left(C\Rightarrow_{\sss{X}}D\right)\right)\right)}\label{eq:exampleOfFormula}\end{equation}
we have that the type $\tau(\xi)$ is defined if and only if all the
following conditions are true:\begin{gather*}
C\subseteq X\\
D\subseteq X\\
(C\Rightarrow_{\sss{X}}D)\subseteq\text{dom}(x)\\
B\subseteq Y\\
\forall_{x}(C\Rightarrow_{\sss{X}}D)\subseteq Y\\
B\cap_{\sss{Y}}\forall_{x}(C\Rightarrow_{\sss{X}}D)\subseteq\text{dom}(\delta)\\
A\subseteq Z\\
\exists_{\delta}(B\cap_{\sss{Y}}\forall_{x}(C\Rightarrow_{\sss{X}}D))\subseteq Z\\
\left[A\Rightarrow_{\sss{Z}}\exists_{\delta}(B\cap_{\sss{Y}}\forall_{x}(C\Rightarrow_{\sss{X}}D))\right]\subseteq\text{dom}(\eps).\end{gather*}

\noindent They are obviously more complicated, but more general, than
conditions \ref{eq:cond1}, \ref{eq:cond2}, \ref{eq:cond3} and \ref{eq:cond4}.
Nevertheless, the hypothesis that the type $\tau(\xi)$ is defined
(which, by Definition \ref{def:formulas}, is a consequence of the
condition that $\xi$ is a formula) is not everything we need to apply
all the theorems of Section \ref{sub:BasicLogicalPropertiesOfFermatFunctor}.
For example, to the previously listed conditions related to the formula
$\xi$ of \eqref{eq:exampleOfFormula}, we have to add hypotheses
like: {}``the spaces $X$, $Y$, $Z$ are separated and the topology
of their Fermat extension is $\ext{(-)}$-generated'', {}``the arrows
$x$, $\delta$, $\eps$, are open and with left inverse'' and {}``all
the subspaces appearing in the previous list of conditions are open
in the corresponding superspace''. We will introduce these types
of hypotheses directly in the statement of the general transfer theorem.

Now we can define the list of objects and arrows occurring in a formula
$\phi$. They are formally different from the free variables defined
for a logical formula, because they have to be thought of as all the
elements of the category $\CInfty$ (or $\ECInfty$) occurring in
the formula $\phi$. These objects and arrows will be the elements
that have to be $\ext{(-)}$-transformed in the general transfer theorem,
so e.g. in the formula $\corners{\exists_{f}(A)}$ the only object
is $A$ and the only arrow is $f$ (whereas in a logical formula of
the form $\exists\, f\,(A)$ the variable $f$ is not free).
\begin{defn}
\label{def:variablesInFormula}Let $\xi\in\mathcal{L}(\CInfty)$ be
a formula, then the \emph{list of objects} ${\rm ob}(\xi)$ and the
\emph{list of arrows} ${\rm ar}(\xi)$ are expressions defined recursively
by the following conditions:
\begin{enumerate}
\item If $\xi=\corners{A}$ for some object $A\in\CInfty$, then \begin{align*}
{\rm ob}(\xi): & =A\\
{\rm ar}(\xi): & =\emptyset.\end{align*}

\item If $\xi=\corners{(\phi\times\psi)}$ for some $\phi$, $\psi\in\mathcal{L}(\CInfty)$,
then \begin{align*}
{\rm ob}(\xi): & =({\rm ob}(\chi),{\rm ob}(\psi))\\
{\rm ar}(\xi): & =({\rm ar}(\chi),{\rm ar}(\psi)).\end{align*}

\item If $\xi=\corners{\neg_{\sss{\chi}}\psi}$ for some $\chi$, $\psi\in\mathcal{L}(\CInfty)$,
then \begin{align*}
{\rm ob}(\xi): & =({\rm ob}(\chi),{\rm ob}(\psi))\\
{\rm ar}(\xi): & =({\rm ar}(\chi),{\rm ar}(\psi)).\end{align*}

\item If $\xi=\corners{(\phi\Rightarrow_{\chi}\psi)}$ for some $\phi,$
$\chi$, $\psi\in\mathcal{L}(\CInfty)$, then \begin{align*}
{\rm ob}(\xi): & =({\rm ob}(\phi),{\rm ob}(\chi),{\rm ob}(\psi))\\
{\rm ar}(\xi): & =({\rm ar}(\phi),{\rm ar}(\chi),{\rm ar}(\psi)).\end{align*}

\item If $\xi=\corners{(\phi\cap_{\sss{\chi}}\psi)}$ for some $\phi,$
$\chi$, $\psi\in\mathcal{L}(\CInfty)$, then \begin{align*}
{\rm ob}(\xi): & =({\rm ob}(\phi),{\rm ob}(\chi),{\rm ob}(\psi))\\
{\rm ar}(\xi): & =({\rm ar}(\phi),{\rm ar}(\chi),{\rm ar}(\psi)).\end{align*}

\item If $\xi=\corners{(\phi\cup_{\sss{\chi}}\psi)}$ for some $\phi,$
$\chi$, $\psi\in\mathcal{L}(\CInfty)$, then \begin{align*}
{\rm ob}(\xi): & =({\rm ob}(\phi),{\rm ob}(\chi),{\rm ob}(\psi))\\
{\rm ar}(\xi): & =({\rm ar}(\phi),{\rm ar}(\chi),{\rm ar}(\psi)).\end{align*}

\item If $\xi=\corners{\exists_{f}(\phi)}$ for some $\phi\in\mathcal{L}(\CInfty)$
and some arrow $f:X\freccia Y$ of $\CInfty$, then \begin{align*}
{\rm ob}(\xi): & ={\rm ob}(\phi)\\
{\rm ar}(\xi): & =(f,{\rm ar}(\phi)).\end{align*}

\item If $\xi=\corners{\forall_{f}(\phi)}$ for some $\phi\in\mathcal{L}(\CInfty)$
and some arrow $f:X\freccia Y$ of $\CInfty$, then \begin{align*}
{\rm ob}(\xi): & ={\rm ob}(\phi)\\
{\rm ar}(\xi): & =(f,{\rm ar}(\phi)).\end{align*}

\item If $\xi=\corners{f^{-1}(\phi)}$ for some $\phi\in\mathcal{L}(\CInfty)$
and some arrow $f:X\freccia Y$ of $\CInfty$, then \begin{align*}
{\rm ob}(\xi): & ={\rm ob}(\phi)\\
{\rm ar}(\xi): & =(f,{\rm ar}(\phi)).\end{align*}

\end{enumerate}
\end{defn}
Now we can define the operator corresponding to a given formula $\xi\in\mathcal{L}(\CInfty)$
simply as the type $\tau(\xi)$ of the formula with the explicit indication
of objects and arrows occurring in the formula itself.
\begin{defn}
\label{def:operatorOfAFormula}If $\xi\in\mathcal{L}(\CInfty)$ is
a formula and ${\rm ob}(\xi)=:(A_{1}\ldots,A_{n})$, ${\rm ar}(\xi)=:(f_{1},\ldots,f_{m})$
are the lists of objects and arrows occurring in $\xi$, then\[
\omega_{\xi}(A_{1},\ldots,A_{n},f_{1},\ldots,f_{m}):=\tau(\xi)\]

\end{defn}
Finally, we can define the Fermat transform of a formula.
\begin{defn}
Let $\xi\in\mathcal{L}(\CInfty)$ be a formula, then the \emph{Fermat
transform }$\ext{\xi}$ is defined recursively by the following conditions:
\begin{enumerate}
\item If $\xi=\corners{A}$ for some object $A\in\CInfty$, then \[
\ext{\xi}:=\ext{A}\]

\item If $\xi=\corners{(\phi\times\psi)}$ for some $\phi$, $\psi\in\mathcal{L}(\CInfty)$,
then \[
\ext{\xi}:=\corners{(\ext{\phi}\times\ext{\psi})}\]

\item If $\xi=\corners{\neg_{\sss{\chi}}\psi}$ for some $\chi$, $\psi\in\mathcal{L}(\CInfty)$,
then \[
\ext{\xi}:=\corners{\neg_{\sss{\ext{\chi}}}\ext{\phi}}\]

\item If $\xi=\corners{(\phi\Rightarrow_{\chi}\psi)}$ for some $\phi,$
$\chi$, $\psi\in\mathcal{L}(\CInfty)$, then \[
\ext{\xi}:=\corners{(\ext{\phi}\Rightarrow_{\sss{\ext{\chi}}}\ext{\psi})}\]

\item If $\xi=\corners{(\phi\cap_{\sss{\chi}}\psi)}$ for some $\phi,$
$\psi\in\mathcal{L}(\CInfty)$, then \[
\ext{\xi}:=\corners{(\ext{\phi}\cap_{\sss{{}^{\bullet\!}\chi}}\ext{\psi})}\]

\item If $\xi=\corners{(\phi\cup_{\sss{\chi}}\psi)}$ for some $\phi,$
$\psi\in\mathcal{L}(\CInfty)$, then \[
\ext{\xi}:=\corners{(\ext{\phi}\cup_{\sss{{}^{\bullet\!}\chi}}\ext{\psi})}\]

\item If $\xi=\corners{\exists_{f}(\phi)}$ for some $\phi\in\mathcal{L}(\CInfty)$
and some arrow $f:X\freccia Y$ of $\CInfty$, then \[
\ext{\xi}:=\corners{\existsext{f}(\ext{\phi})}\]

\item If $\xi=\corners{\forall_{f}(\phi)}$ for some $\phi\in\mathcal{L}(\CInfty)$
and some arrow $f:X\freccia Y$ of $\CInfty$, then \[
\ext{\xi}:=\corners{\forallext{f}(\ext{\phi})}\]

\item If $\xi=\corners{f^{-1}(\phi)}$ for some $\phi\in\mathcal{L}(\CInfty)$
and some arrow $f:X\freccia Y$ of $\CInfty$, then \[
\ext{\xi}:=\corners{\ext{f}^{-1}(\ext{\phi})}\]

\end{enumerate}
\end{defn}
We can now state the general transfer theorem:
\begin{thm}
\noindent \label{thm:generalTransferTheorem}Let $\xi\in\mathcal{L}(\CInfty)$
be a formula in $\CInfty$ without occurrences of $\corners{\times}$,
and let ${\rm ob}(\xi)=:(A_{1},\ldots,A_{n})$, ${\rm ar}(\xi)=:(f_{1},\ldots,f_{m})$
be objects and arrows occurring in the formula $\xi$. Let us suppose
that for every $i=1,\ldots,m$ and every $j$, $k=1,\ldots,n$:
\begin{enumerate}
\item \label{enu:f_i_isOpenAndWithLeftInverse}$f_{i}:X_{i}\freccia Y_{i}$
is open and with left inverse.
\item \label{enu:subformulasAndOpen}Let $\phi$ and $\psi$ be subformulae
of $\xi$ and $Z$ be any space in the list $\tau(\psi)$, $X_{1},\ldots,X_{n}$,
$Y_{1},\ldots,Y_{n}$, then\[
\tau(\phi)\subseteq Z\then|\tau(\phi)|\text{ is open in }Z.\]

\item \label{enu:subformulasAndSeparated}Let $\phi$ be a subformula of
$\xi$, then the topology of $\ext{\tau(\phi)}$ is $\ext{(-)}$-generated.
\item \label{enu:separatedAndTopology}All the spaces $X_{i}$ are separated
and the topology of their Fermat extension is $\ext{(-)}$-generated.
\end{enumerate}
Then we have:\[
\ext{\left[\omega_{\xi}(A_{1},\ldots,A_{n},f_{1},\ldots,f_{m})\right]}=\omega_{{}^{\bullet\!}\xi}(\ext{A_{1}},\ldots,\ext{A_{n}},\ext{f_{1}},\ldots,\ext{f_{m}})\]

\end{thm}
\noindent For manifolds we can also include the product:
\begin{thm}
\noindent \label{thm:generalTransferTheoremForManifolds}Let $\xi\in\mathcal{L}(\CInfty)$
be a generic formula in $\CInfty$, and let ${\rm ob}(\xi)=:(A_{1},\ldots,A_{n})$,
${\rm ar}(\xi)=:(f_{1},\ldots,f_{m})$ be objects and arrows occurring
in the formula $\xi$. Let us suppose that for every $i=1,\ldots,m$
and every $j$, $k=1,\ldots,n$:
\begin{enumerate}
\item $f_{i}:X_{i}\freccia Y_{i}$ is open and with left inverse.
\item Let $\phi$ and $\psi$ be subformulae of $\xi$ and $Z$ be any space
in the list $\tau(\psi)$, $X_{1},\ldots,X_{n}$, $Y_{1},\ldots,Y_{n}$,
then\[
\tau(\phi)\subseteq Z\then|\tau(\phi)|\text{ is open in }Z.\]

\item Let $\phi$ be a subformula of $\xi$, then the topology of $\ext{\tau(\phi)}$
is $\ext{(-)}$-generated.
\item All the spaces $X_{i}$ are separated and the topology of their Fermat
extension is $\ext{(-)}$-generated.
\item If $\corners{(\phi\times\psi)}$ is a subformula of $\xi$, then $\tau(\phi)$
and $\tau(\psi)$ are manifolds.
\end{enumerate}
Then we have:\begin{equation}
\ext{\left[\omega_{\xi}(A_{1},\ldots,A_{n},f_{1},\ldots,f_{m})\right]}=\omega_{{}^{\bullet\!}\xi}(\ext{A_{1}},\ldots,\ext{A_{n}},\ext{f_{1}},\ldots,\ext{f_{m}})\label{eq:conclusionTrasferTheorem}\end{equation}

\end{thm}
\noindent \textbf{Proof of Theorem \ref{thm:generalTransferTheorem}
and Theorem \ref{thm:generalTransferTheoremForManifolds}:} We proceed
by induction on the length of the formula $\xi$. If $\xi$ is made
of one object only, i.e. $\xi=A$, then $\omega_{\xi}(A_{1},\ldots,A_{n},f_{1},\ldots,f_{m})=\tau(\xi)=A$,
$n=1$, $A_{1}=A$, $m=0$. Analogously\[
\omega_{{}^{\bullet\!}\xi}(\ext{A_{1}},\ldots,\ext{A_{n}},\ext{f_{1}},\ldots,\ext{f_{m}})=\ext{A}\]
 since $\ext{\xi}=\ext{A}$ and $\ext{\tau}(\ext{\xi})=\ext{A}$;
the conclusion is hence trivial.

Now suppose that the equality \eqref{eq:conclusionTrasferTheorem}
is true for every formula of length less than $N>0$ and that in the
formula $\xi$ occur $N$ symbols. Using the Definition \eqref{def:formulas}
we have to consider several cases depending on the form of $\xi$.
We will proceed for the case $\xi=\corners{(\phi\Rightarrow_{\sss{\chi}}\psi)}$,
$\xi=\corners{\exists_{f_{i}}(\phi)}$ and $\xi=\corners{f_{i}^{-1}(\phi)}$,
the other ones being analogous.

In the case $\xi=\corners{\exists_{f_{i}}(\phi)}$, from the Definition
\eqref{def:formulas} we get $\tau(\phi)\subseteq X_{i}$. Moreover,
$\phi$ is a subformula of $\xi$, hence from the hypothesis \emph{\ref{enu:subformulasAndOpen}}.
we obtain that $|\tau(\phi)|$ is open in $X_{i}$. We can thus apply
Theorem \ref{thm:existentialQuantifier} since $f_{i}$ is open by
hypotheses \emph{\ref{enu:f_i_isOpenAndWithLeftInverse}}., and we
obtain that\begin{equation}
\ext{[\tau(\xi)]}=\ext{[\exists_{f_{i}}(\tau(\phi)]}=\existsext{f_{i}}(\ext{[\tau(\phi)]}).\label{eq:GTT_preservationExists}\end{equation}
By induction hypotheses, we get\begin{equation}
\ext{[\tau(\phi)]}=\ext{\tau}(\ext{\phi})=\omega_{\ext{\phi}}(\ext{A_{r_{1}}},\ldots,\ext{A_{r_{a}}},\ext{f_{s_{1}}},\ldots,\ext{f_{s_{b}}})\label{eq:inductive1}\end{equation}
where ${\rm ob}(\phi)=(A_{r_{1}},\ldots,A_{r_{a}})$ and ${\rm ar}(\phi)=(f_{s_{1}},\ldots,f_{s_{b}})$
are objects and arrows occurring in $\phi$, hence $\{r_{1},\ldots,r_{a}\}\subseteq\{1,\ldots,n\}$
and $\{s_{1},\ldots,s_{b}\}\subseteq\{1,\ldots,m\}$. On the other
hand, $\ext{\xi}=\corners{\existsext{f_{i}}(\ext{\phi})}$ and hence\begin{equation}
\ext{\tau}(\ext{\xi})=\existsext{f_{i}}(\ext{\tau}(\ext{\phi})).\label{eq:typeOfFermatTransform2}\end{equation}
The conclusion for this case follows from \eqref{eq:GTT_preservationExists},
\eqref{eq:inductive1} and \eqref{eq:typeOfFermatTransform2}, indeed:\begin{align*}
\ext{\left[\omega_{\xi}(A_{1},\ldots,A_{n},f_{1},\ldots,f_{m})\right]} & =\ext{[\tau(\xi)]}\\
 & =\existsext{f_{i}}(\ext{[\tau(\phi)]})\\
 & =\existsext{f_{i}}(\omega_{\ext{\phi}}(\ext{A_{r_{1}}},\ldots,\ext{A_{r_{a}}},\ext{f_{s_{1}}},\ldots,\ext{f_{s_{b}}}))\\
 & =\ext{\tau}(\ext{\xi})\\
 & =\omega_{{}^{\bullet\!}\xi}(\ext{A_{1}},\ldots,\ext{A_{n}},\ext{f_{1}},\ldots,\ext{f_{m}}).\end{align*}

Let us note that the hypotheses that the topology of all the spaces
$\ext{X_{i}}$ is $\ext{(-)}$-generated must be used in the case
$\xi=\corners{\forall_{f_{i}}(\phi)}$.

In the case $\xi=\corners{(\phi\Rightarrow_{\sss{\chi}}\psi)}$, from
the Definition \eqref{def:formulas} we get $\tau(\phi)\subseteq\tau(\chi)$
and $\tau(\psi)\subseteq\tau(\chi)$. Moreover, $\phi$, $\chi$ and
$\psi$ are subformulae of $\xi$, hence from the hypothesis \emph{\ref{enu:subformulasAndOpen}}.
we obtain that both $|\tau(\phi)|$ and $|\tau(\psi)|$ are open in
$\tau(\chi)$, and from the hypothesis \emph{\ref{enu:subformulasAndSeparated}}.
we get that the topology of $\ext{\tau(\chi)}$ is $\ext{(-)}$-generated.
We can hence apply Theorem \ref{thm:FermatFunctor_AndOrNegation}
obtaining that\begin{equation}
\ext{[\tau(\xi)]}=\ext{[\tau(\phi)\Rightarrow_{\sss{\tau(\chi)}}\tau(\psi)]}=\ext{[\tau(\phi)]}\Rightarrow_{\sss{\ext{[\tau(\chi)}]}}\ext{[\tau(\psi)]}.\label{eq:GTT_preservationImplications}\end{equation}
But, by induction hypotheses we get equalities like \eqref{eq:inductive1},
i.e.:\begin{equation}
\ext{[\tau(\phi)]}=\ext{\tau}(\ext{\phi})=\omega_{\ext{\phi}}(\ext{A_{r_{1}}},\ldots,\ext{A_{r_{a}}},\ext{f_{s_{1}}},\ldots,\ext{f_{s_{b}}})\label{eq:inductive1bis}\end{equation}
\begin{equation}
\ext{[\tau(\chi)]}=\ext{\tau}(\ext{\chi})=\omega_{\ext{\chi}}(\ext{A_{t_{1}}},\ldots,\ext{A_{t_{c}}},\ext{f_{u_{1}}},\ldots,\ext{f_{u_{d}}})\label{eq:inductive2}\end{equation}
\begin{equation}
\ext{[\tau(\psi)]}=\ext{\tau}(\ext{\psi})=\omega_{\ext{\psi}}(\ext{A_{v_{1}}},\ldots,\ext{A_{v_{e}}},\ext{f_{w_{1}}},\ldots,\ext{f_{w_{h}}}),\label{eq:inductive3}\end{equation}
On the other hand, $\ext{\xi}=\corners{\ext{\phi}\Rightarrow_{\sss{\ext{\chi}}}\ext{\psi}}$
and hence\begin{equation}
\ext{\tau}(\ext{\xi})=\ext{\tau}(\ext{\phi})\Rightarrow_{\sss{\ext{\tau}(\ext{\chi})}}\ext{\tau}(\ext{\psi}).\label{eq:typeOfFermatTransform1}\end{equation}
The conclusion for the first case follows from \eqref{eq:GTT_preservationImplications},
\eqref{eq:inductive1bis}, \eqref{eq:inductive2}, \eqref{eq:inductive3}
and \eqref{eq:typeOfFermatTransform2}.

Finally, let us note that in the case $\xi=\corners{f_{i}^{-1}(\phi)}$
we have to use the hypotheses \emph{\ref{enu:subformulasAndOpen}}.
to prove that $|\tau(\phi)|$ is open in $Y_{i}$, but we do not need
any other hypotheses on the codomain space $Y_{i}$. For this reason
condition \emph{\ref{enu:separatedAndTopology}}. is stated for the
domain spaces $X_{i}$ only.

For Theorem \ref{thm:generalTransferTheoremForManifolds} we can proceed
in a similar way, using Theorem \ref{thm:FermatPreservesProductOfMan}
in case of formulae of type $\xi=\corners{(\phi\times\psi)}$.$\qedNoNewLine$

It is natural to expect that there would be some relationship between
our transfer theorem and a transfer theorem more similar to those
of NSA. The principal difference is that our transfer theorem, even
if it concerns formulae, it is used to construct spaces and the theorem
itself states an equality between spaces of $\ECInfty$. On the contrary,
the transfer theorem of NSA asserts an equivalence between two sentences.
Nevertheless, it seems possible to follow the following scheme:
\begin{enumerate}
\item Define the meaning of the sentence {}``\emph{the formula $\xi$ is
intuitionistically true in $\CInfty$}'' using the intuitionistic
interpretation of the propositional connectives and quantifiers in
this category. An analogous definition of intuitionistic validity
can be done in the category $\ECInfty$.
\item Define the $\ext{(-)}$-transform of a given formula $\xi$.
\item Prove that $\xi$ is intuitionistically true in $\CInfty$ if and
only if $\ext{\xi}$ is intuitionistically true in $\ECInfty$.
\end{enumerate}
This work is planned in future projects.

A specification is adequate here. Though the theory of Fermat reals
is compatible with classical logic, the previous theorems state that
the Fermat functor behaves really better if the logical formulae are
interpreted in open sets. This may seem in contraddiction with the
thread of the present work (see Section \ref{cha:IntroductionAndGeneralProblem}).
Indeed, we remember that one of the main aims of the present work
is to develop a sufficiently powerful theory of infinitesimal \emph{without
forcing the reader to learn a strong formal control of the mathematics
he/she is doing}, e.g. forcing the reader to learn to work in intuitionistic
logic\emph{.} Of course, this is not incompatible with the \emph{possibility}
to gain more if one is interested to have this type of strong formal
control, e.g. if one is already able to work in intuitionistic logic,
and the results of this section go exactly in this direction.

\part{\label{par:TheStartingOfANewTheory}The beginning of a new theory}

\chapter{\label{cha:CalculusOnOpenDomains}Calculus on open domains}

\section{Introduction}

We have defined and studied plenty of instruments that can be useful
to develop the differential and integral calculus of functions defined
on infinitesimal domains like $D_{k}$ or on bigger sets like the
extension $\ext{(a,b)}$ of a real interval. We can then start the
development of infinitesimal differential geometry, following, where
possible, the lines of SDG. But further development can be glimpsed
in the calculus of variations, because of cartesian closedness of
our categories, because of the possibility to use infinitesimal methods
and because of the properties of diffeological maps that, e.g., do
not require any compactness hypothesis on the domain of our functions%
\footnote{We can say that compactness assumptions are only required because
of the non adequacy of a tool like normed space (as our Chapter \ref{cha:TheCategoryCn}
and Section \ref{sec:DefinitionOfCartesianClosure} prove), in the
sense that nothing in the problem of defining smooth spaces and maps
forces us to introduce a norm.%
}. Of course, this could also open the possibility of several applications,
e.g. in general relativity or in continuum mechanics. Indeed, Fermat
reals can be considered as the first theory of infinitesimals having
a good intuitive interpretation and without the need to possess a
non trivial background of knowledge in formal logic to be understood
(see Appendix \ref{app:theoriesOfInfinitesimals}), and this characteristic
can be very useful for its diffusion among physicists, engineers and
even mathematicians.

But, exactly as SDG required tens of years to be developed, we have
to expect a comparable amount of time for the full development of
applications to the geometry of the approach we introduced here. At
the same time, Fermat reals seems sufficiently stable and with good
properties to permit us to state that such a development can be achieved.

In this chapter we want to introduce the basic theorems and ideas
that permits this further development. We shall prove all the theorems
which are useful for the development of the calculus both for $\ECInfty$
functions of the form $f:\ext{U}\freccia\ER^{d}$, where $U$ is open
in $\ER^{n}$, and for functions defined on infinitesimal sets, like
$g:D_{n}^{d}\freccia\ext{X}$. Subsequently we shall present a first
development of infinitesimal differential geometry, primarily in manifolds
and in spaces of smooth functions of the form $\ext{N}^{\ext{M}}$.

\bigskip{}

Using the Taylor's formula as stated in Theorem \ref{thm:TaylorForStndSmoothAtNonStndPoint},
we have a powerful instrument to manage derivatives of functions $\ext{f}$
obtained as extensions of ordinary smooth functions $f:\R^{d}\freccia\R^{u}$.
But this is not the case if $f:\ER^{d}\freccia\ER^{u}$ is a generic
$\ECInfty$ arrow, that is if we can write locally $f(x)=\ext{\alpha}(p,x)$,
where $p\in\ER^{n}$ and $g$ is smooth, because generally speaking
$f$ does not have standard derivatives $\partial_{j}f(x)\in\ER^{u}\setminus\R^{u}$.
Therefore, the problem arises how to define the derivatives of this
type of functions in our setting. On the one hand, we would like to
set e.g. $f'(x):=\ext{(\partial\alpha/\partial x)(p,x)}$ (if $d=u=1$,
for simplicity), and so the problem would become the independence
in this definition from both the function $g$ and the non standard
parameter $p$. For example, for functions defined on an infinitesimal
domain we can see that this problem of independence is not trivial.
Let us consider two first order infinitesimals $p$, $p'\in D$, $p\ne p'$.
Because the product of first order infinitesimals is always zero,
we have that the null function $f(x)=0$, for $x\in D$, can be written
both as $f(x)=p\cdot x=:\ext{\alpha}(p,x)$ and as $f(x)=p'\cdot x=\ext{\alpha}(p',x)$.
But $\ext{(\partial\alpha/\partial x)(p,x)}=p\ne p'=\ext{(\partial\alpha/\partial x)(p',x)}$.
For functions defined on an open set, this independence can be established,
using the method originally used by Fermat and studied by G.E. Reyes
(see \citet{Mo-Re}; see also \citet{Ber} and \citet{Sha} for analogous
ideas in a context different from that of SDG).

In all this section we will use the notation for intervals as subsets
of $\ER$, e.g. $[a,b):=\left\{ x\in\ER\,|\, a\le x<b\right\} $.
Notations of the type\[
[a,b)_{\R}:=\left\{ x\in\R\,|\, a\le x<b\right\} \]
will be used to specify that the interval has to be understood as
a subset of $\R$.

\section{\label{sec:theFermatMethod}The Fermat-Reyes method}

The method used by Fermat to calculate derivatives is to assume $h\ne0$,
to construct the incremental ratio\[
\frac{f(x+h)-f(x)}{h},\]
and then to set $h=0$ in the final result. This idea, which sounds
as inconsistent, can be perfectly understood if we think that the
incremental ratio can be extended with continuity at $h=0$ if the
function $f$ is differentiable at $x$. In our smooth context, we
need a theorem confirming the existence of a {}``smooth version''
of the incremental ratio. We firstly introduce the notion of segment
in an $n$-dimensional space $\ER^{n}$, that, as we will prove later,
for $n=1$ coincide with the notion of interval in $\ER$.
\begin{defn}
\label{def:segment}If $a$, $b\in\ER^{n}$, then \[
\overrightarrow{[a,b]}:=\left\{ a+s\cdot(b-a)\,|\, s\in[0,1]\right\} \]
is the segment of $\ER^{n}$ going from $a\in\ER^{n}$ to $b\in\ER^{n}$.\end{defn}
\begin{thm}
\label{thm:smoothIncrementalRatio}Let $U$ be an open set of $\R$,
and $f:\ext{U}\freccia\ER$ be a $\ECInfty$ function. Let us define
the \emph{thickening of $\ext{U}$ along the $x$-axis }by\emph{ }\[
\widetilde{\ext{U}}:=\left\{ (x,h)\,|\,\overrightarrow{[x,x+h]}\subseteq\ext{U}\right\} ,\]
Then $\widetilde{\ext{U}}$ is open in $\ER^{2}$ and there exists
one and only one $\ECInfty$ map $r:\widetilde{\ext{U}}\freccia\ER$
such that \[
f(x+h)=f(x)+h\cdot r(x,h)\quad\forall(x,h)\in\widetilde{\ext{U}}.\]
Hence we define $f^{\prime}(x):=r(x,0)\in\ER$ for every $x\in\ext{U}$.

\noindent Moreover if $f(x)=\ext{\alpha}(p,x)$, $\forall x\in\mathcal{V}\subseteq\ext{U}$
with $\alpha\in\Cc^{\infty}(A\times B,\R)$, then \[
f'(x)={}^{{}^{{}^{{\scriptstyle \bullet}}}}{\!\left(\frac{\partial\alpha}{\partial x}\right)}(p,x).\]

\end{thm}
\noindent We anticipate the proof of this theorem by the following
lemmas
\begin{lem}
\label{lem:thickeningIsOpen}Let $U$ be an open set of $\R^{n}$
and $v\in\ER^{n}$, then the \emph{thickening of $\ext{U}$ along
$v$} defined as\begin{equation}
\widetilde{\ext{U}_{v}}:=\left\{ (x,h)\in\ER^{n}\times\ER\,|\,\overrightarrow{[x,x+hv]}\subseteq\ext{U}\right\} \label{eq:thickeningGeneric}\end{equation}
is open in $\ER^{n}\times\ER$.
\end{lem}
\noindent \textbf{Proof:} Let us take a generic point $(x,h)\in\widetilde{\ext{U}_{v}}$;
we want to prove that $(x,h)\in\ext{(A\times B)}\subseteq\widetilde{\ext{U}_{v}}$
for some subsets $A$ of $\ER^{n}$ and $B$ of $\ER$. Because the
point $(x,h)$ is in the thickening, we have that\[
\forall s\in[0,1]\pti x+s\cdot hv\in\ext{U}.\]
Taking the standard parts we obtain\[
\forall s\in[0,1]_{\R}\pti\st{x}+s\cdot\st{h}\cdot\st{v}=:\phi(s)\in U.\]
The function $\phi:[0,1]_{\R}\freccia U$ is continuous and thus\[
\phi\left([0,1)_{\R}\right)=\overrightarrow{[\st{x},\st{x}+\st{h}\st{v}]}=:K\]
is compact in $\R^{n}$. But $K\subseteq U$ and $U$ is open, so
the distance of $K$ from the complement $\R^{n}\setminus U$ is strictly
positive; let us call $2a:=d\left(K,\R^{n}\setminus U\right)>0$ this
distance, so that for every $c\in K$ we have that\[
B_{a}(c):=\left\{ x\in\R^{n}\,|\, d(x,c)<a\right\} \subseteq U.\]
Now, set $A:=B_{a/2}(\st{x})$ and $B:=B_{b}(\st{h})$, where we have
fixed $b\in\R_{>0}$ such that $b\cdot\Vert\st{v}\Vert\le\frac{a}{2}$.
We have $x\in\ext{A}$ because $\st{x}\in A$ and $A$ is open; analogously
$h\in\ext{B}$ and thus $(x,h)\in\ext{A}\times\ext{B}=\ext{(A\times B)}$.
We have finally to prove that taking a generic point $(y,k)\in\ext{(A\times B)}$,
the whole segment $\overrightarrow{[y,y+kv]}$ is contained in $\ext{U}$;
so, let us take also a Fermat number $0\le s\le1$. Since $U$ is
open, to prove that $y+skv\in\ext{U}$ is equivalent to prove that
the standard part $y+skv$ is in $U$, i.e. that $\st{y}+\st{s}\st{k}\st{v}\in U$.
For, let us observe that\[
\Vert\st{y}+\st{s}\st{k}\st{v}-\st{x}-\st{s}\st{h}\st{v}\Vert\le\Vert\st{y}-\st{x}\Vert+|\st{s}|\cdot\Vert\st{v}\Vert\cdot|\st{k}-\st{h}|\le\frac{a}{2}+1\cdot\Vert\st{v}\Vert\cdot b\le a.\]
Therefore, $\st{y}+\st{s}\st{k}\st{v}\in B_{a}(c)\subseteq U$, where
$c=\st{x}+\st{s}\st{h}\st{v}\in K$ from our definition of the compact
set $K$.$\qedNoNewLine$
\begin{lem}
\label{lem:twoWaysToSeeIntervals}If $a$, $b\in\ER$, then\begin{align*}
a & <b\then\overrightarrow{[a,b]}=[a,b]\\
b & \le a\then\overrightarrow{[a,b]}=[b,a].\end{align*}

\end{lem}
\noindent \textbf{Proof:} We will prove the first implication, the
second being a simple consequence of the first one. To prove the inclusion
$\overrightarrow{[a,b]}\subseteq[a,b]$ take $x=a+s\cdot(b-a)$ with
$0\le s\le1$, then $0\le s\cdot(b-a)\le b-a$ because $b-a>0$. Adding
$a$ to these inequalities we get $a\le x\le b$. For the proof of
the opposite inclusion, let us consider $a\le x\le b$. If we prove
the inclusion for $a=0$ only, we can prove it in general: in fact,
$0\le x-a\le b-a$, so that if $[0,b-a]\subseteq\overrightarrow{[0,b-a]}$
we can derive the existence of $s\in[0,1]$ such that $x-a=0+s\cdot(b-a)$,
which is our conclusion. So, let us assume that $a=0$. If $\st{b}\ne0$,
then $b$ is invertible and it suffices to set $s:=\frac{x}{b}$ to
have the conclusion. Otherwise, $\st{b}=0$ and hence also $\st{x}=0$.
Let us consider the decompositions of $x$ and $b$\begin{align*}
x & =\sum_{i=1}^{k}\st{x_{i}}\cdot\diff{t_{\omega_{i}(x)}}\\
b & =\sum_{j=1}^{h}\st{b_{j}}\cdot\diff{t_{\omega_{j}(b)}}.\end{align*}
We have to find a number $s=\st{s}+\sum_{n=1}^{N}\st{s_{n}}\cdot\diff{t_{\omega_{n}(s)}}$
such that $s\cdot b=x$. It is interesting to note that the attempt
to find the solution $s\in[0,1]$ directly from these decompositions
and from the property $s\cdot b=x$ is not as easy as to find the
solution using directly little-oh polynomials. In fact $\forall^{0}t>0:\ b_{t}>0$
because $b>0$ and hence for $t>0$ sufficiently small, we can form
the ratio\begin{align}
\frac{x_{t}}{b_{t}} & =\frac{\sum_{i=1}^{k}\st{x_{i}}\cdot t^{\frac{1}{\omega_{i}(x)}}}{\sum_{j=1}^{h}\st{b_{j}}\cdot t^{\frac{1}{\omega_{j}(b)}}}\nonumber \\
 & =\frac{t^{\frac{1}{\omega_{1}(b)}}\cdot\sum_{i=1}^{k}\st{x_{i}}\cdot t^{\frac{1}{\omega_{i}(x)}-\frac{1}{\omega_{1}(b)}}}{t^{\frac{1}{\omega_{1}(b)}}\cdot\sum_{j=1}^{h}\st{b_{j}}\cdot t^{\frac{1}{\omega_{j}(b)}-\frac{1}{\omega_{1}(b)}}}.\label{eq:ratio_x_over_b_firstStep}\end{align}
Let us note that from Theorem \ref{thm:effectiveCriterionForTheOrder}
we can deduce that $\st{x_{1}}>0$ since $x>0$ and hence that $\omega(b)=\omega_{1}(b)>\omega(x)\ge\omega_{i}(x)$
because $x<b$. From \eqref{eq:ratio_x_over_b_firstStep} we have

\begin{align*}
\frac{x_{t}}{b_{t}} & =\frac{\sum_{i=1}^{k}\st{x_{i}}\cdot t^{\frac{1}{\omega_{i}(x)}-\frac{1}{\omega_{1}(b)}}}{\st{b_{1}}\cdot\left(1+\sum_{j=2}^{h}\frac{\st{b_{j}}}{\st{b_{1}}}\cdot t^{\frac{1}{\omega_{j}(b)}-\frac{1}{\omega_{1}(b)}}\right)}\\
 & =\frac{1}{\st{b_{1}}}\cdot\sum_{i=1}^{k}\st{x_{i}}\cdot t^{\frac{1}{\omega_{i}(x)}-\frac{1}{\omega_{1}(b)}}\cdot\sum_{k=0}^{+\infty}(-1)^{k}\cdot\left(\sum_{j=2}^{h}\frac{\st{b_{j}}}{\st{b_{1}}}\cdot t^{\frac{1}{\omega_{j}(b)}-\frac{1}{\omega_{1}(b)}}\right)^{k}.\end{align*}
Writing, for simplicity, $a\odot b:=\frac{a\cdot b}{a+b}$ we can
write the previous little-oh polynomial using the common notation
with $\diff{t}_{a}$:\begin{equation}
s_{t}:=\frac{x_{t}}{b_{t}}=\frac{1}{\st{b_{1}}}\cdot\sum_{i=1}^{k}\st{x_{i}}\cdot\diff{t_{\omega_{j}(x)\odot\h{.009}\omega_{1}(b)}}\cdot\sum_{k=0}^{+\infty}(-1)^{k}\cdot\left(\sum_{j=2}^{h}\frac{\st{b_{j}}}{\st{b_{1}}}\cdot\diff{t_{\omega_{j}(b)\odot\h{.009}\omega_{1}(b)}}\right)^{k}.\label{eq:DefinitionOf_s}\end{equation}
As usual, the series in this formula is really a finite sum, because
$D_{\infty}$ is an ideal of nilpotent infinitesimals. Going back
in these passages, it is quite easy to prove that $ $the previously
defined $s\in\ER$ verifies the desired equality $s\cdot b=x$. Moreover,
from Theorem \ref{thm:equivalentFormulationForOrderRelation} the
relations $0\le s\le1$ follow.$\qedNoNewLine$

\noindent It is interesting to make some considerations based on the
proof of this lemma. Indeed, we have just proved that in the Fermat
reals every equation of the form $a+x\cdot b=c$ with $a<c<a+b$ has
a solution%
\footnote{Let us note explicitly, that this is not in contradiction with the
non Archimedean property of $\ER$ (let $a=0$ and $b\in D_{\infty}$)
because of the inequalities that $c$ must verifies to have a solution.%
}. If $b$ is invertible, this is obvious and we have a unique solution.
If $b$ is a nilpotent infinitesimal, a possible solution is given
by a formula like \eqref{eq:DefinitionOf_s}, but we do not have uniqueness.
E.g. if $a=0$, $c=\diff{t_{2}}+\diff{t}$ and $b=\diff{t_{3}}$,
then $x=\diff{t_{6}}+\diff{t_{3/2}}$ is a solution of $a+x\cdot b=c$,
but $x+\diff{t}$ is another solution because $\diff{t}\cdot\diff{t_{a}}=0$
for every $a\ge1$. Among all the solutions in the case $b\in D_{\infty}$,
we can choose the simplest one, i.e. that {}``having no useless addends
in its decomposition'', that is such that\[
\frac{1}{\omega_{i}(x)}+\frac{1}{\omega(b)}\le1\]
for every addend $\st{x_{i}}\cdot\diff{t_{\omega_{i}(x)}}$ in the
decomposition of $x$. Otherwise, if for some $i$ we have the opposite
inequality, we can apply Lemma \ref{lem:producOfInfinitesimalAndIota}
with $k:=\omega_{i}(x)$ to have that $b\cdot x=b\cdot\iota_{k}x$,
i.e. we can delete some {}``useless addend'' considering $\iota_{k}x$
instead of $x$. We can thus understand that this algebraic problem
is strictly tied with the definition of derivative $f'(x)$, which
is the solution of the linear equation $f(x+h)=f(x)+h\cdot f'(x)$:
as we give an hint in Chapter \ref{cha:equalityUpTo_k-thOrder}, if
$f$ is defined only on an infinitesimal set like $D_{n}$, this equation
has not a unique solution and we can define the derivative $f'(x)$
only by considering {}``the simplest solution'', i.e. using a suitable
$\iota_{k}$. We will get back to the problem of defining $f'(x)$,
where the function $f$ is defined on an infinitesimal set, in the
next Section \ref{cha:CalculusOnInfinitesimalDomains}.

The uniqueness of the smooth incremental ratio stated in Theorem \ref{thm:smoothIncrementalRatio}
is tied with the following lemma, for the proof of which we decided
to introduce nilpotent paths (see Definition \ref{def:NilpotentFunctions})
instead of continuous paths at $t=0$, like in \citep{Gio}. We will
call this lemma \emph{the cancellation law of non-infinitesimal functions}.
\begin{lem}
\textbf{\emph{(Cancellation law of non-infinitesimal functions):}}\label{lem:cancellationOfNon-infFunctions}

\noindent Let $U$ be an open neighborhood of $0$ in $\R$, and let\[
f,\ g:\ext{U}\freccia\ER\]
be two $\ECInfty$ functions such that\[
\forall x\in\ext{U}:\ x\text{ is invertible}\then g(x)\text{ is invertible and }g(x)\cdot f(x)=0.\]
Then $f$ is the null function, i.e. $f=0$.
\end{lem}
\noindent \textbf{Proof:} We have that $f:\overline{\ext{U}}\freccia\ER$
and hence $f\In{\ext{U}}\ER$ and we can apply Theorem \ref{thm:figuresOfFermatSpaces}
at the point $0\in\ext{U}$ obtaining that the function $f$ can be
written as\[
f(x)=\ext{\alpha}(p,x)\quad\forall x\in\ext{B}\cap\ext{U}=\ext{(B\cap U)}=:\mathcal{V},\]
where $\alpha\in\CInfty(A\times B,\R)$, $p\in\ext{A}$, $A$ is an
open set of $\R^{\sf p}$ and $B$ is an open neighborhood of $0$
in $\R$. We can always assume that $\st{p}=0$ because, otherwise,
we can consider the standard smooth function $(y,x)\mapsto\alpha(y-\st{p},x)$.
We can thus write our main hypotheses as\begin{equation}
\forall x\in\mathcal{V}:\ x\text{ is invertible }\then\lim_{t\to0^{+}}\frac{g(x)_{t}\cdot\alpha(p_{t},x_{t})}{t}=0.\label{eq:mainHypothesesIn_V}\end{equation}
Let us provide some explanation about the notation $g(x)_{t}$ which
is a consequence of our notations concerning quotient sets: we have
that $g(x)\in\ER=\R_{o}[t]/\sim$, hence, avoiding the use of equivalence
classes in favor of the new notion of equality $\sim$ in $\R_{o}[t]$,
we have that $g(x)$ is a little-oh polynomial and hence $g(x):\R_{\ge0}\freccia\R$,
from which the notation $g(x)_{t}\in\R$ for $t\in\R_{\ge0}$ acquires
a clear meaning. We firstly want to prove that $\alpha(p_{t},x_{t})=o(t)$
for every $x\in\mathcal{V}$. Let us take a generic infinitesimal
$h\in D_{\infty}$ and choose a $k\in\N_{>0}$ such that%
\footnote{This passage is possible exactly because we are considering nilpotent
paths as elements of $\ER$.%
}\begin{align*}
h^{k} & =0\e{in}\ER\\
p^{k} & =\underline{0}\e{in}\ER^{\sf p},\end{align*}
and consider a generic non zero $r\in U\cap B\setminus\{0\}$. Then
$x:=h+r\in\ext{U}$, because $\st{x}=r\in U$, and $x$ is invertible
because its standard part is $r$ and $r\ne0$. From our hypothesis
we have that $g(x)$ is also invertible, i.e. $\st{g(x)}=g(x)_{0}=\lim_{t\to0^{+}}g(x)_{t}\ne0$,
and hence from \eqref{eq:mainHypothesesIn_V} we get\begin{equation}
\lim_{t\to0^{+}}\frac{\alpha(p_{t},h_{t}+r)}{t}=0.\label{eq:alphaIsALittleOh_x_invertible}\end{equation}
Because every invertible $x\in\mathcal{V}$ can be written as $x=h+r$
with $r\in\R_{\ne0}$ and $h\in D_{\infty}$, we have just proved
our conclusion for every $x\in\mathcal{V}$ which is invertible. Now
we have to prove \eqref{eq:alphaIsALittleOh_x_invertible} for $r=0$
too. Let us consider the Taylor's formula of order $k$ with the function
$\alpha$ at the point $(\underline{0},r)$ (which obviously is true
for $r=0$ too):\begin{multline}
\frac{\alpha(\underline{0}+p_{t},r+h_{t})}{t}=\frac{1}{t}\cdot\left[\sum_{\substack{q\in\N^{\sf p+1}\\
|q|\le k}
}\frac{\partial^{q}\alpha}{\partial(p,x)^{q}}(\underline{0},r)\cdot\frac{(p_{t},h_{t})^{q}}{q!}+\right.\\
\left.+\sum_{\substack{q\in\N^{\sf p+1}\\
|q|=k+1}
}\frac{\partial^{q}\alpha}{\partial(p,x)^{q}}(\xi_{t},\eta_{t})\cdot\frac{(p_{t},h_{t})^{q}}{q!}\right]\label{eq:TaylorForAlphaProofCancellationNonInfFunctions}\end{multline}
with $\xi_{t}\in(\underline{0},p_{t})$ and $\eta_{t}\in(r,r+h_{t})$.
But $h^{k}=0$ and $p^{k}=(p_{1},\ldots,p_{\sf p})^{k}=(p_{1}^{k},\ldots,p_{\sf p}^{k})=\underline{0}$,
hence $h$, $p_{i}\in D_{k}$. Moreover, if $|q|=k+1$, then\[
\sum_{i=1}^{\sf{p}+1}\frac{q_{i}}{k+1}=\frac{k+1}{k+1}=1,\]
so that from Corollary \ref{cor:suffCondToHaveProductZeroFrom_D_j_k}
we get\[
(p_{t},h_{t})^{q}=p_{1}(t)^{q_{1}}\cdot\ldots\cdot p_{\sf p}(t)^{q_{\sf p}}\cdot h(t)^{q_{\sf{p}+1}}=o(t).\]
Therefore, from \eqref{eq:alphaIsALittleOh_x_invertible} and \eqref{eq:TaylorForAlphaProofCancellationNonInfFunctions}
we obtain\begin{equation}
\lim_{t\to0^{+}}\sum_{\substack{q\in\N^{\sf p+1}\\
|q|\le k}
}\frac{\partial^{q}\alpha}{\partial(p,x)^{q}}(\underline{0},r)\cdot\frac{1}{q!}\cdot\frac{(p_{t},h_{t})^{q}}{t}=0\quad\forall r\in(U\cap B)_{\ne0}.\label{eq:limitFromTaylor}\end{equation}
Now, let $\{q_{1},\ldots,q_{N}\}$ be an enumeration of all the $q\in\N^{\sf{p}+1}$
such that $|q|\le k$, and for simplicity set\begin{align*}
b_{i}(r): & =\frac{\partial^{q_{i}}\alpha}{\partial(p,x)^{q_{i}}}(\underline{0},r)\cdot\frac{1}{q_{i}!}\quad\forall r\in U\cap B\\
s_{i}(t): & =\frac{(p_{t},h_{t})^{q_{i}}}{t}\quad\forall t\in\R_{\ge0},\end{align*}
so that we can write \eqref{eq:limitFromTaylor} as\begin{equation}
\forall r\in(U\cap B)_{\ne0}\pti\lim_{t\to0^{+}}\sum_{i=1}^{N}b_{i}(r)\cdot s_{i}(t)=0.\label{eq:limitIn_bi_and_si}\end{equation}
If all the functions $b_{i}$ are identically zero, then $b_{i}(\bar{r})=b_{i}(0)$
where $\bar{r}\in U\cap B\setminus\{0\}$, which always exists because
$U\cap B$ is open in $\R$. Therefore, \eqref{eq:limitIn_bi_and_si}
(and hence also \eqref{eq:limitFromTaylor}) is true for $r=0$ too.
Otherwise, taking a base of the subspace of $\CInfty(U\cap B,\R)$
generated by the smooth functions $b_{1},\ldots,b_{N}$ and expressing
all the $b_{i}$ in this base, we can suppose to have in \eqref{eq:limitIn_bi_and_si}
only linearly independent functions.

We can now use the following lemma:
\begin{lem}
\label{lem: IndependentFunctionsAndDeterminant}Let $U$ be an open
neighborhood of $0$ in $\R$ and $b_{1},\ldots,b_{N}:U\freccia\R$
be linearly independent functions continuous at $0$. Then we can
find\[
r_{1},\ldots,r_{N}\in U\setminus\{0\}\]
such that\[
\det\left[\begin{array}{ccc}
b_{1}(r_{1}) & \dots & b_{N}(r_{1})\\
\vdots & {} & \vdots\\
b_{1}(r_{N}) & \dots & b_{N}(r_{N})\end{array}\right]\ne0.\]

\end{lem}
\noindent From \eqref{eq:limitIn_bi_and_si} we can write\[
\lim_{t\to0^{+}}\left[\begin{array}{ccc}
b_{1}(r_{1}) & \dots & b_{N}(r_{1})\\
\vdots & {} & \vdots\\
b_{1}(r_{N}) & \dots & b_{N}(r_{N})\end{array}\right]\cdot\left[\begin{array}{c}
s_{1}(t)\\
\vdots\\
s_{N}(t)\end{array}\right]=\underline{0}\]
and hence from this lemma we can deduce that $s_{i}(t)\to0$ for $t\to0^{+}$.
Because these limits exist, we can take the limit for $r\to0$ of
\eqref{eq:limitIn_bi_and_si} and proceed in the following way\begin{align*}
\lim_{r\to0}\lim_{t\to0^{+}}\sum_{i=1}^{N}b_{i}(r)\cdot s_{i}(t) & =\sum_{i=1}^{N}\lim_{r\to0}b_{i}(r)\cdot\lim_{t\to0^{+}}s_{i}(t)\\
 & =\lim_{t\to0^{+}}\sum_{i=1}^{N}b_{i}(0)\cdot s_{i}(t)=0\end{align*}
(let us note that we do not exchange the limit signs). This proves
that \eqref{eq:limitIn_bi_and_si} is true for $r=0$ too. From \eqref{eq:TaylorForAlphaProofCancellationNonInfFunctions}
for $r=0$ we obtain\[
\lim_{t\to0^{+}}\frac{\alpha(p_{t},h_{t})}{t}=0.\]
This proves that $f(x)=0$ in $\ER$ for every $x\in\mathcal{V}$.
Finally, if $x\in\ext{U}\setminus\mathcal{V}$ then $\st{x}\ne0$
because otherwise we would have $x\in\mathcal{V}=\ext{(B\cap U)}$.
So $x$ is invertible and hence also $g(x)$ is invertible, so that
from $g(x)\cdot f(x)=0$ we can easily deduce $f(x)=0$ also in this
case.$\qedNoNewLine$

\noindent \textbf{Proof of Lemma \ref{lem: IndependentFunctionsAndDeterminant}:}
We prove the converse by induction on $N\ge2$, i.e. if all the determinants
cited in the statement are zero, then the functions $(b_{1},\ldots,b_{N})$
are linearly dependent. Let us suppose first that $N=2$ and that
all these determinants are zero, that is \begin{equation}
b_{1}(r)\cdot b_{2}(s)=b_{2}(r)\cdot b_{1}(s)\quad\forall r,s\in U_{\ne0}.\label{eq: N2HypAbs}\end{equation}
If the functions $b_{i}$, $i=1,2$, are both zero then they are trivially
linearly dependent, hence let us suppose, e.g., that $b_{1}(\bar{s})\ne0$
for some $\bar{s}\in U$. Due to the continuity of $b_{1}$ at $0$
we can suppose $\bar{s}\ne0$, hence from (\ref{eq: N2HypAbs}) \[
b_{2}(r)=b_{1}(r)\cdot\frac{b_{2}(\bar{s})}{b_{1}(\bar{s})}=:b_{1}(r)\cdot a\quad\forall r\in U_{\ne0}.\]
From the continuity of $b_{i}$ at $0$ we have that $b_{2}=b_{1}\cdot a$,
that is $(b_{1},b_{2})$ are linearly dependent.

\noindent Now suppose that the implication is true for any matrix
of $N$ functions and we prove the conclusion for matrices of order
$N+1$ too. By Laplace's formula with respect to the first row, for
every $r_{1},\ldots,r_{N+1}\in U_{\ne0}$ we have \begin{multline}
b_{1}(r_{1})\cdot\left|\begin{array}{ccc}
b_{2}(r_{2}) & \dots & b_{N+1}(r_{2})\\
\vdots & {} & \vdots\\
b_{2}(r_{N+1}) & \dots & b_{N+1}(r_{N+1})\end{array}\right|-\ldots+\\
+(-1)^{N+2}\cdot b_{N+1}(r_{1})\cdot\left|\begin{array}{ccc}
b_{1}(r_{2}) & \dots & b_{N}(r_{2})\\
\vdots & {} & \vdots\\
b_{1}(r_{N+1}) & \dots & b_{N}(r_{N+1})\end{array}\right|=0.\label{eq: LaplaceDevelopment}\end{multline}
Now we have two cases. Let $\alpha_{1}(r_{2},\ldots,r_{N+1})$ denote
the first determinant in the previous \eqref{eq: LaplaceDevelopment}.
If it is zero for any $r_{2},\ldots,r_{N+1}\in U_{\ne0}$, then by
the induction hypothesis $(b_{2},\ldots,b_{N+1})$ are linearly dependent,
hence the conclusion follows. Otherwise $\bar{\alpha}_{1}:=\alpha_{1}(\bar{r}_{2},\ldots,\bar{r}_{N+1})\ne0$
for some $\bar{r}_{2},\ldots,\bar{r}_{N+1}\in U_{\ne0}$. Then from
\eqref{eq: LaplaceDevelopment} it follows \[
b_{1}(r_{1})=b_{2}(r_{1})\cdot\frac{\alpha_{2}}{\bar{\alpha}_{1}}-\ldots-(-1)^{N+2}\cdot b_{N+1}(r_{1})\cdot\frac{\alpha_{N+1}}{\bar{\alpha}_{1}}\quad\forall r_{1}\in U_{\ne0},\]
where we used obvious notations for the other determinants in \eqref{eq: LaplaceDevelopment}.
From the continuity of $b_{i}$ the previous formula is true for $r_{1}=0$
too and this proves the conclusion.$\qedNoNewLine$

\noindent \textbf{Proof of Theorem \ref{thm:smoothIncrementalRatio}:}
We will define the function $r:\widetilde{\ext{U}}\freccia\ER$ patching
together smooth functions defined on open subsets covering $\widetilde{\ext{U}}$.
Therefore, we have to take a generic point $(x,h)\in\widetilde{\ext{U}}$,
to define the function $r$ on some open neighborhood of $(x,h)$
in $\widetilde{\ext{U}}$, and to prove that every two of such local
functions agree on the intersection of their domains.

\noindent As usual, we have that $f\In{\ext{U}}\ER$ and, since $x\in\ext{U}$,
we can write\begin{equation}
f|_{\mathcal{V}}=\ext{\alpha}(p,-)|_{\mathcal{V}},\label{eq:f_equal_alpha_proof_incremental_ratio}\end{equation}
where $\alpha\in\CInfty(\bar{U}\times\bar{V},\R)$, $\mathcal{V}:=\ext{\bar{V}}\cap\ext{U}=\ext{(\bar{V}\cap U)}$
is an open neighborhood of $x$ and $\ext{\bar{U}}$ is an open neighborhood
of $p\in\ER^{\sf p}$ defined by the open subset $\bar{U}$ of $\R^{\sf p}$.
Because $\widetilde{\ext{U}}$ is open in $\ER\times\ER=\ER^{2}$,
we can find two open subset $A$ and $B$ of $\R$ such that\[
(x,h)\in\ext{(A\times B)}\subseteq\widetilde{\ext{U}}\]
and such that\begin{equation}
a+s\cdot b\in\bar{V}\quad\forall a\in A,\ b\in B,\ s\in[0,1]_{\R}.\label{eq:interval_a_b_containedIn_V_bar}\end{equation}
Let us define\begin{equation}
\gamma(q,a,b):=\int_{0}^{1}\partial_{2}\alpha(q,a+s\cdot b)\,\diff{s}\quad\forall q\in\bar{U},\ a\in A,\ b\in B.\label{eq:defOfGamma_proofCancellationNonInfFunctions}\end{equation}
We have that $\gamma\in\CInfty(\bar{U}\times A\times B,\R)$, so that
if we define\[
r(a,b):=\ext{\gamma}(p,a,b)\quad\forall(a,b)\in\ext{(A\times B)},\]
then we have\begin{equation}
r\in\ECInfty(\ext{(A\times B)},\R)\label{eq:rSmooth}\end{equation}
\begin{equation}
\ext{(A\times B)}\text{ open neighborhood of }(u,h)\text{ in }\widetilde{\ext{U}}.\label{eq:A_Times_B_neighOf_x_h}\end{equation}
For every $(a,b)\in\ext{(A\times B)}$ we have\begin{align}
b_{t}\cdot r(a,b)_{t} & =\int_{0}^{1}\partial_{2}\alpha(p_{t},a_{t}+s\cdot b_{t})\cdot b_{t}\,\diff{s}\nonumber \\
 & =\int_{a_{t}}^{a_{t}+b_{t}}\partial_{2}\alpha(p_{t},y)\,\diff{y}\nonumber \\
 & =\alpha(p_{t},a_{t}+b_{t})-\alpha(p_{t},a_{t}).\label{eq:ratioPropertyIn_a_b}\end{align}
But from $(a,b)\in\ext{(A\times B)}=\ext{A}\times\ext{B}$ and \eqref{eq:interval_a_b_containedIn_V_bar}
it follows $\st{a}$, $\st{a}+\st{b}\in\bar{V}$, and hence also $a$,
$a+b\in\ext{\bar{V}}$. From the definition of thickening we also
have that $a$, $a+b\in\ext{U}$. We can thus use \eqref{eq:f_equal_alpha_proof_incremental_ratio}
at the points $a$, $b\in\mathcal{V}=\ext{\bar{V}}\cap\ext{U}$, so
that we can write \eqref{eq:ratioPropertyIn_a_b} as\begin{equation}
\forall(a,b)\in\ext{(A\times B)}\pti b\cdot r(a,b)=f(a+b)-f(a).\label{eq:ratioPropertyIn_a_b_with_f}\end{equation}
We have proved that for every $(x,h)\in\widetilde{\ext{U}}$ there
exist an open neighborhood $\ext{(A\times B)}$ of $(x,h)$ in $\widetilde{\ext{U}}$
and a smooth function $r\in\ECInfty(\ext{(A\times B)},\R)$ such that
\eqref{eq:ratioPropertyIn_a_b_with_f} holds.

If $\rho\in\ECInfty\left(\ext{(C\times D)},\ER\right)$ is another
such functions, then\[
\forall(x,h)\in\ext{(C\times D)}\cap\ext{(A\times B)}\pti h\cdot\left[r(x,h)-\rho(x,h)\right]=0,\]
so that for every $x\in\ext{C}\cap\ext{A}$ we have that\[
\forall h\in\ext{(D\times B)}\pti h\cdot\left[r(x,h)-\rho(x,h)\right]=0.\]
For Lemma \ref{lem:cancellationOfNon-infFunctions} applied with $g(h):=h$
and $f(h):=r(x,h)-\rho(x,h)$, we have $r(x,h)=\rho(x,h)$ for every
$(x,h)\in\ext{(C\times D)}\cap\ext{(A\times B)}$, which proves the
conclusion for the sheaf property of $\ER$. Finally, let us note
that from \eqref{eq:defOfGamma_proofCancellationNonInfFunctions}
for $b=0$ we obtain $r(a,0)=\partial_{2}\alpha(p,a)$, which is the
last part of the statement.$\qedNoNewLine$

Using this theorem, we can develop all the differential calculus for
non standard smooth functions of type $f:\ER^{n}\freccia\ER$. We
will see now the first steps of this development, underlining the
main differences with respect to \citet{Lav} and \citet{Mo-Re},
to which we refer as a guideline for a complete development.
\begin{defn}
\label{def:smoothIncrementalRatio}Let $U$ be an open subset of $\R$,
and $f:\ext{U}\freccia\ER$ a $\ECInfty$ function. Then
\begin{enumerate}
\item $f'[-]:\widetilde{\ext{U}}\freccia\ER$\label{enu:defSmoothIncrRatio1}
\item $f(x+h)=f(x)+h\cdot f'[x,h]\quad\forall(x,h)\in\widetilde{\ext{U}}$.\label{enu:defSmoothIncrRatio2}
\end{enumerate}
\noindent Moreover we will also set $f'(x):=f'[x,0]$ for every $x\in\ext{U}$.

\end{defn}
\noindent Let us note that the notation for the smooth incremental
ratio as a function uses square brackets like in $f'[-]$. For this
reason there is no way to confuse the smooth incremental ratio $f'[-]$
and its values $f'[x,h]$ with the corresponding derivative $f'$
and its values $f'(x)$.

First of all, from property \emph{\ref{enu:defSmoothIncrRatio1}}.
in the previous definition, it follows that\[
f':\ext{U}\freccia\ER.\]
The following theorem contains the first expected properties of the
derivative.
\begin{thm}
Let $U$ be an open subset of $\R$, and $f$, $g:\ext{U}\freccia\ER$
be smooth $\ECInfty$ functions. Finally, let us consider a Fermat
real $r\in\ER$. Then
\begin{enumerate}
\item $\left(f+g\right)'=f'+g'$
\item $\left(r\cdot f\right)'=r\cdot f'$
\item $\left(f\cdot g\right)'=f'\cdot g+f\cdot g'$
\item $\left(1_{\ER}\right)'=1$
\item $r'=0$
\end{enumerate}
\end{thm}
\noindent \textbf{Proof:} We report the proof essentially as a first
example to show how to use precisely the Fermat-Reyes method in our
context.

The first step is to prove, e.g., that $f+g$ is smooth in $\ECInfty$.
Looking at the diagram\[
\xyR{30pt}\xyC{30pt}\xymatrix{ & \ER\\
\ext{U}\ar[ur]^{f}\ar[dr]_{g}\ar[rr]^{\langle f,g\rangle} &  & \ER^{2}\ar[lu]_{p_{1}}\ar[ld]^{p_{2}}\ar[r]^{+} & \ER,\\
 & \ER}
\]
where $+:(r,s)\in\ER^{2}\mapsto r+s\in\ER$ is the sum of Fermat reals,
we can see that $f+g=\langle f,g\rangle\cdot+$ and hence it is smooth
because it can be expressed as a composition of smooth functions.
The proof that the sum $f+g$ is smooth, even if it is almost trivial,
can show us why it is very important to work in a cartesian closed
category like $\ECInfty$. We have, indeed, the possibility to consider
very general set theoretical operations like compositions or evaluations.

Now we have only to calculate $(f+g)(x+h)$ using the definition of
smooth incremental ratio and its uniqueness\begin{align*}
(f+g)(x+h) & =f(x+h)+g(x+h)\\
 & =f(x)+h\cdot f'[x,h]+g(x)+h\cdot g'[x,h]\\
 & =(f+g)(x)+h\cdot\left\{ f'[x,h]+g'[x,h]\right\} \quad\forall(x,h)\in\widetilde{\ext{U}.}\end{align*}
From the uniqueness of the smooth incremental ratio of $f+g$ we obtain
$(f+g)'[-]=f'[-]+g'[-]$ and thus the conclusion evaluating these
ratios at $h=0$.

As a further simple example, we consider only the derivative of the
product. The smoothness of $f\cdot g$ can be proved analogously to
what we have just done for the sum. Now, let us evaluate for every
$(x,h)\in\widetilde{\ext{U}}$\begin{align*}
(f\cdot g)(x+h) & =f(x+h)\cdot g(x+h)\\
 & =\left\{ f(x)+h\cdot f'[x,h]\right\} \cdot\left\{ g(x)+h\cdot g'[x,h]\right\} \\
 & =(f\cdot g)(x)+h\cdot\\
 & \phantom{=}\cdot\left\{ f(x)\cdot g'[x,h]+g(x)\cdot f'[x,h]+h^{2}\cdot f'[x,h]\cdot g'[x,h]\right\} .\end{align*}
From the uniqueness of the smooth incremental ratio of $f\cdot g$
we have thus\[
(f\cdot g)'[x,h]=f(x)\cdot g'[x,h]+g(x)\cdot f'[x,h]+h^{2}\cdot f'[x,h]\cdot g'[x,h],\]
which gives the conclusion setting $h=0$. The other properties can
be proved analogously.$\qedNoNewLine$

The next expected property that permits a deeper understanding of
the Fermat-Reyes method is the chain rule.
\begin{thm}
\label{thm:chainRule}If $U$ and $V$ are open subsets of $\R$ and\[
f:\ext{U}\freccia\ER\]
\[
g:\ext{V}\freccia\ext{U}\]
are $\ECInfty$ functions, then\[
(f\circ g)'=(f'\circ g)\cdot g'.\]

\end{thm}
\noindent We will give a proof of this theorem with the aim of explaining
in a general way the Fermat-Reyes method. We first need the following
\begin{lem}
\label{lem:compactnessSimpleCase}Let $U$ be an open subset of $\R^{\sf k}$,
$x\in\ext{U}$ and $v\in\ER^{\sf k}$. Then there exists\[
r\in\R_{>0}\]
such that\[
\forall h\in(-r,r)\pti(x,h)\in\widetilde{\ext{U}_{v}}.\]

\end{lem}
\noindent \textbf{Proof:} If $\st{v}=\underline{0}$, then for every
$s\in[0,1]$ and every $h\in\ER$ we have $\st{(x+shv)}=\st{x}\in U$,
hence $x+shv\in\ext{U}$, that is $\overrightarrow{[x,x+hv]}\subseteq\ext{U}$.
In this case we have thus $(x,h)\in\widetilde{\ext{U}_{v}}$ for every
$h\in\ER$.

Otherwise, if $\st{v}\ne\underline{0}$ then from $\st{x}\in U$ we
obtain\[
\exists\,\rho>0\pti B_{\rho}(\st{x})\subseteq U\]
because $U$ is open in $\R^{\sf k}$. Take as $r\in\R_{>0}$ any
real number verifying\[
0<r<\min\left(\rho,\frac{\rho}{\Vert\st{v}\Vert}\right).\]
For such an $r$, if $s\in[0,1]$ and $h\in(-r,r)$, then\begin{align}
\st{(x+shv)}=\st{x}+\st{s}\cdot\st{h}\cdot\st{v}\in B_{\rho}(\st{x})\iff & \Vert\st{s}\cdot\st{h}\cdot\st{v}\Vert<\rho\nonumber \\
\Longleftarrow\ \, & |\st{h}|\cdot\Vert\st{v}\Vert<\rho\label{eq:compactnessSimpleCaseLast}\end{align}
the last implication is due to the assumption that $s\in[0,1]$. But
\eqref{eq:compactnessSimpleCaseLast} holds because $|h|<r$ and hence
$\st{|h|}=|\st{h}|<r$ and $r\cdot\Vert\st{v}\Vert<\rho$ for the
definition of $r$.$\qedNoNewLine$

\noindent The next result works for the Fermat-Reyes methods like
a sort of {}``compactness principle'' analogous to the compactness
theorem of mathematical logic. It is the generalization to more than
just one open set $U$ of the previous lemma.
\begin{thm}
\textbf{\emph{(Compactness principle):}}\label{thm:compactnessPrinciple}

\noindent For $i=1,\ldots,n$, let $U^{i}$ be open sets of $\R^{\sf{k}_{i}}$,
$v\in\ER^{\sf{k}_{i}}$, $x_{i}\in\ext{U^{i}}$ and finally $a_{i}\in\ER$.
Then there exists\[
r\in\R_{>0}\]
such that\[
\forall i=1,\ldots,n\ \forall h\in(-r,r)\pti(x_{i},h\cdot a_{i})\in\widetilde{\ext{U_{v_{i}}^{i}}.}\]

\end{thm}
\noindent \textbf{Proof:} For every $x_{i}\in U^{i}$ we apply the
previous Lemma \ref{lem:compactnessSimpleCase} obtaining the existence
of $r_{i}\in\R_{>0}$ such that\begin{equation}
\forall k\in(-r_{i},r_{i})\pti(x_{i},k)\in\widetilde{\ext{U_{v_{i}}^{i}}}.\label{eq:1_compactnessPrinciple}\end{equation}
Now, let us set\[
r:=\min_{i:\st{a_{i}}\ne0}\frac{r_{i}}{|\st{a_{i}}|}\in\R_{>0},\]
then taking a generic $h\in(-r,r)$ we have\begin{equation}
-r<\st{h}<r.\label{eq:2_compactnessPrinciple}\end{equation}
If $\st{a_{i}}=0$, then trivially $-r_{i}<\st{h}\cdot\st{a_{i}}<r_{i}$
and hence $-r_{i}<h\cdot a_{i}<r_{i}$, so that from \eqref{eq:1_compactnessPrinciple}
we get the conclusion for this first case, i.e. $(x_{i},ha_{i})\in\widetilde{\ext{U_{v_{i}}^{i}}}$.

\noindent Otherwise, if $\st{a_{i}}\ne0$, then $r\le\frac{r_{i}}{|\st{a_{i}}|}$
and from \eqref{eq:2_compactnessPrinciple} we get $|\st{h}|<r\le\frac{r_{i}}{|\st{a_{i}}|}$
and hence $-r_{i}<ha_{i}<r$, and once again the conclusion follows
from \eqref{eq:1_compactnessPrinciple}.$\qedNoNewLine$

We can use this theorem in the following way:
\begin{enumerate}
\item every time in a proof we need a property of the form\begin{equation}
(x_{i},ha_{i})\in\widetilde{\ext{U_{i}}}\label{eq:Hp_i_FermatMethod}\end{equation}
we will assume {}``to have chosen $h$ so little that \eqref{eq:Hp_i_FermatMethod}
is verified''.
\item We derive the conclusion $\mathcal{A}(h)$ under $n$ of such hypothesis,
so that we have concretely deduced that\[
\left(\forall i=1,\ldots,n\pti(x_{i},ha_{i})\in\widetilde{\ext{U_{i}}}\right)\then\mathcal{A}(h).\]

\item At this point we can ap\-ply the com\-pact\-ness prin\-ci\-ple
obtaining\[
\exists\, r\in\R_{>0}\ \forall h\in(-r,r)\pti\mathcal{A}(h).\]

\item Usually the property $\mathcal{A}(h)$ is of the form\begin{equation}
\mathcal{A}(h)\iff h\cdot\tau(h)=h\cdot\sigma(h),\label{eq:1_explFermatMethod}\end{equation}
and hence we can deduce $\tau(h)=\sigma(h)$ for every $h\in(-r,r)$
from the cancellation law of non-infinitesimal functions, and in particular
$\tau(0)=\sigma(0)$. If the property $\mathcal{A}$ has the form
\eqref{eq:1_explFermatMethod}, then we can also suppose that $h$
is invertible because the cancellation law can be applied also in
this case. But at the end we will anyway set $h=0$, in perfect agreement
with the classical description of the Fermat method (see e.g. \citet{Bo-Fr-Ri,EBell,Edw}).
\end{enumerate}
Let us note that, as mentioned above, conceptually this way to proceed
reflects the same idea of the compactness theorem of mathematical
logic, because in every proof we can only have a finite number of
hypothesis of type \eqref{eq:Hp_i_FermatMethod}. Even if this method
does not involve explicitly infinitesimal methods, using it the final
proofs are very similar to those we would have if $h$ were an actual
infinitesimal, i.e. $h\in D_{\infty}$.

In the following proof we will concretely use this method.

\noindent \textbf{Proof of Theorem \ref{thm:chainRule}:} First of
all the composition\[
(-)\circ(-):\ext{U}^{\ext{V}}\times\ER^{\ext{U}}\freccia\ER^{\ext{V}}\]
is a smooth map of $\ECInfty$ and hence $f\circ g$ is smooth because
it can be written as a composition of smooth maps.

For a generic\begin{equation}
(x,h)\in\widetilde{\ext{V}}\label{eq:Hp1_ChainRule}\end{equation}
we can always write\[
(f\circ g)(x+h)=f\left[g(x+h)\right]=f\left[g(x)+h\cdot g'[x,h]\right]\]
because $x+h\in\ext{V}$ and hence $f\circ g$ is defined at $x+h$.
Now we would like to use the smooth incremental ratio of $f$ at the
point $g(x)$ with increment $h\cdot g'[x,h]$. For this end we assume\begin{equation}
(g(x),h\cdot g'[x,h])\in\widetilde{\ext{U}}\label{eq:Hp2_chainRule}\end{equation}
so that we can write\[
(f\circ g)(x+h)=f(gx)+h\cdot g'[x,h]\cdot f'\left[gx,h\cdot g'[x,h]\right].\]
Using the compactness principle and the cancellation law of non-infinitesimal
functions we get\[
\exists\, r\in\R_{>0}:\ \forall h\in(-r,r)\pti g'[x,h]\cdot f'\left[gx,h\cdot g'[x,h]\right]=(f\circ g)'[x,h],\]
and thus the conclusion for $h=0$.$\qedNoNewLine$

\noindent Let us note that these ideas, that do not use infinitesimal
methods, can be repeated in a standard context, with only slight modifications,
so that they represent an interesting alternative way to teach a significant
part of the calculus with strongly simpler proofs.

To realize a comparison with the Levi-Civita field (see Appendix \ref{app:theoriesOfInfinitesimals})
we now prove the inverse function theorem.
\begin{thm}
Let $U$ be an open subset of $\R$, $x$ a point in $\ext{U}$, and\[
f:\ext{U}\freccia\ER\]
a $\ECInfty$ map such that\[
f'(x)\text{ is invertible}.\]
Then there exist two open subsets $X$, $Y$ of $\R$ such that
\begin{enumerate}
\item $x\in\ext{X}$ and $f(x)\in\ext{Y}$, i.e. $\ext{X}$ and $\ext{Y}$
are open neighborhoods of $x$ and $f(x)$ respectively
\item $f|_{\ext{X}}:\ext{X}\freccia\ext{Y}$ is invertible and $\left(f|_{\ext{X}}\right)^{-1}:\ext{Y}\freccia\ext{X}$
is a $\ECInfty$ map
\item ${\displaystyle \left[\left(f|_{\ext{X}}\right)^{-1}\right]'(fx_{1})=\frac{1}{f'(x_{1})}}$
for every $x_{1}\in\ext{X}$
\end{enumerate}
\end{thm}
\noindent \textbf{Proof:} Because $x\in\ext{U}$ we can write $f|_{\mathcal{V}}=\ext{\alpha}(p,-)|_{\mathcal{V}}$,
where $\alpha\in\CInfty(A\times B,\R)$, $p\in\ext{A}$, $A$ is an
open set of $\R^{\sf p}$ and $B$ is an open subset of $\R$ such
that $x\in\ext{B}$ and finally $\mathcal{V}:=\ext{B}\cap\ext{U}$.
Considering $B\cap U$ instead of $B$ we can assume, for simplicity,
that $B\subseteq U$.

We have that%
\footnote{Because it is sufficiently clear from the context, we use here simplified
notations like $\partial_{2}\alpha(p,x)$ instead of $\ext{(\partial_{2}\alpha)}(p,x)$.%
} $f'(x)=\partial_{2}\alpha(p,x)$ is invertible, hence its standard
part is not zero\[
\st{f'(x)}=\partial_{2}\alpha(\st{p},\st{x})\in\R_{\ne0}.\]
Since $\alpha$ is smooth, we can find a neighborhood $C\times D\subseteq A\times B\subseteq A\times U$
of $(\st{p},\st{x})$ where $\partial_{2}\alpha(p_{1},x_{1})\ne0$
for every $(p_{1},x_{1})\in C\times D$. We can also assume to have
taken this neighborhood sufficiently small in order to have that\begin{equation}
\inf_{(p_{1},x_{1})\in C\times D}\left|\partial_{2}\alpha(p_{1},x_{1})\right|=:m>0.\label{eq:4_inverseFunctionTheorem}\end{equation}
By the standard implicit function theorem, we get an open neighborhood
$E\times X\subseteq C\times D$ of $(\st{p},\st{x})$, an open neighborhood
$Y$ of $\alpha(\st{p},\st{x})$ and a smooth function $\beta\in\CInfty(E\times Y,X)$
such that\begin{equation}
\forall p_{1}\in E\,\forall x_{1}\in X\pti\alpha(p_{1},x_{1})\in Y\label{eq:2_inverseFunctionTheorem}\end{equation}
\begin{equation}
\alpha\left[p_{1},\beta(p_{1},y_{1})\right]=y_{1}\quad\forall(p_{1},y_{1})\in E\times Y\label{eq:1_inverseFunctionTheorem}\end{equation}
\begin{equation}
\forall p_{1}\in E\,\forall y\in Y\ \exists!\, x\in X\pti\alpha(p_{1},x)=y.\label{eq:3_inverseFunctionTheorem}\end{equation}
We can assume that $X$ is connected. Let us define $g:=\ext{\beta}(p,-)$,
then $g\in\CInfty(\ext{Y},\ext{X})$. Moreover, $x\in\ext{X}$ and
$f(x)\in\ext{Y}$ because $\st{x}\in X$, $\st{f(x)}=\alpha(\st{p},\st{x})\in Y$
and $X$, $Y$ are open. From \eqref{eq:2_inverseFunctionTheorem},
if $x_{1}\in\ext{X}$, then $\st{f(x_{1})}=\alpha(\st{p},\st{x_{1}})\in Y$,
hence $f(x_{1})\in\ext{Y}$, so that $f$ maps $\ext{X}$ in $\ext{Y}$.
From \eqref{eq:1_inverseFunctionTheorem}, noting that $p\in\ext{E}$,
because $\st{p}\in E$, and that $X\subseteq D\subseteq B$, we obtain\[
\forall y\in\ext{Y}\pti f\left(g(y)\right)=\alpha\left[p,g(y)\right]=\alpha\left[p,\beta(p,y)\right]=y.\]
This proves that $g$ is a smooth left%
\footnote{With respect to the notation for the composition $(g\cdot f)(y)=f(g(y))$.%
} inverse of $f|_{\ext{X}}:\ext{X}\freccia\ext{Y}$, which is thus
surjective. If we prove that $f|_{\ext{X}}$ is injective, this left
inverse will also be the right inverse. So, let us suppose that $f(x_{1})=f(x_{2})$
in $\ext{Y}$ for $x_{1}$, $x_{2}\in\ext{X}$, i.e.\[
\lim_{t\to0^{+}}\frac{\alpha(p_{t},x_{1t})-\alpha(p_{t},x_{2t})}{t}=0.\]
But we can write\[
\alpha(p_{t},x_{1t})-\alpha(p_{t},x_{2t})=(x_{1t}-x_{2t})\cdot\partial_{2}\alpha(p_{t},\xi_{t})\quad\forall t\in\R_{>0}\]
for a suitable $\xi_{t}\in(x_{1t},x_{2t})$. Moreover, from \eqref{eq:4_inverseFunctionTheorem}
and from\[
\forall^{0}t>0\pti\xi_{t}\in(x_{1t},x_{2t})\subseteq X\subseteq D\e{and}p_{t}\in E\subseteq C\]
(here we are using the assumption that $X$ is connected), we have
that $\left|\partial_{2}\alpha(p_{t},\xi_{t})\right|\ge m$. Therefore\begin{align*}
\lim_{t\to0^{+}}\left|\frac{x_{1t}-x_{2t}}{t}\right| & =\lim_{t\to0^{+}}\left|\frac{\alpha(p_{t},x_{1t})-\alpha(p_{t},x_{2t})}{t\cdot\partial_{2}\alpha(p_{t},\xi_{t})}\right|\ge\\
 & \ge\lim_{t\to0^{+}}\left|\frac{\alpha(p_{t},x_{1t})-\alpha(p_{t},x_{2t})}{t\cdot m}\right|=0,\end{align*}
and this proves that $x_{1}=x_{2}$ in $\ext{X}$ and thus also that
$f|_{\ext{X}}:\ext{X}\freccia\ext{Y}$ is invertible with smooth inverse
given by $\left(f|_{\ext{X}}\right)^{-1}=g$.

Now we can use the Fermat-Reyes method to prove the formula for the
derivative of the inverse function. Let us consider a point $(x_{1},h)\in\widetilde{\ext{X}}$
in the thickening of $\ext{X}$, then $f(x_{1}+h)=f(x_{1})+h\cdot f'[x_{1},h]$.
Applying $g=\left(f|_{\ext{X}}\right)^{-1}$ to both sides of this
formula we obtain\[
x_{1}+h=g\left[fx_{1}+h\cdot f'[x_{1},h]\right].\]
It is natural, at this point, to try to use the smooth incremental
ratio of the smooth function $g$. For this end we have to assume
that\[
(fx_{1},h\cdot f'[x_{1},h])\in\widetilde{\ext{Y}}\]
so that we can write\[
x_{1}+h=g\left(f(x_{1})\right)+h\cdot f'[x_{1},h]\cdot g'\left[fx_{1},h\cdot f'[x,h]\right].\]
Because $g\left(f(x_{1})\right)=x_{1}$, we obtain the equality\[
h=h\cdot f'[x_{1},h]\cdot g'\left[fx_{1},h\cdot f'[x,h]\right].\]
From the compactness principle (Theorem \ref{thm:compactnessPrinciple})
and the cancellation law of non-infinitesimal functions (Lemma \ref{lem:cancellationOfNon-infFunctions})
we obtain\[
1=f'[x_{1},h]\cdot g'\left[fx_{1},h\cdot f'[x,h]\right],\]
from which the conclusion follows setting $h=0$.$\qedNoNewLine$

We have shown, using meaningful examples, that the Fermat-Reyes me\-th\-od
can be used to try a generalization of several results of differential
calculus to $\ECInfty$ functions of the form $f:\ext{U}\freccia\ER^{d}$,
with $U$ open in $\R^{n}$.

Indeed, this can be done for several theorems. We only list here the
main results that we have already proved, leaving a complete report
of them for a subsequent work. For most of them the proofs are very
similar to the analogous presented e.g. in \citet{Lav}:
\begin{enumerate}
\item \noindent the formula for the derivative of $\frac{1}{f(x)}$ if $f(x)\in\ER$
is invertible,
\item the notion of right and left derivatives, i.e. $f'_{+}(a)$ and $f'_{-}(b)$
for a $\ECInfty$ function of the form $f:[a,b]\freccia\ER^{d}$,
\item \noindent definition of higher order derivatives using higher order
smooth incremental ratios,
\item \noindent 1-dimensional Taylor's formula with integral rest (see the
next Section \ref{sec:IntegralCalculus} about the integral calculus),
\item \noindent uniqueness theorem for Taylor's formulas,
\item \noindent the functional operation of taking the derivative is smooth,
i.e. the map $f\in\ER^{\ext{U}}\mapsto f'\in\ER^{\ext{U}}$ is $\ECInfty$,
\item \noindent the functional operation of taking the smooth incremental
ratio is smooth, i.e. the map $f\in\ER^{\ext{U}}\mapsto f'[-]\in\ER^{\widetilde{\ext{U}}}$
is $\ECInfty$,
\item \noindent definition of partial derivatives using smooth partial incremental
ratio,
\item \noindent the functional operation of taking the partial derivative
and the smooth partial incremental ratio are smooth,
\item \noindent linearity of the map: $v\in\ER^{n}\mapsto\frac{\partial f}{\partial v}(x)\in\ER^{d}$,
\item \noindent definition of differentials of arbitrary order,
\item \noindent Euler-Schwarz theorem (differentials are symmetric),
\item \noindent $d$-dimensional chain rule,
\item \noindent several variables Taylor's formula with integral rest,
\item \noindent uniqueness of $d$-dimensional Taylor's formula,
\item \noindent majoration of differentials: $\Vert\diff{^{i}f}.h^{i}\Vert\le M\cdot\Vert h\Vert^{i}$
for every $h\in\ER^{d}$ and some positive constant $M$,
\item \noindent infinitesimal Taylor's formula for functions of the form
$f:\ext{U}\freccia\ER^{d}$ and $U$ open in $\R^{n}$. 
\end{enumerate}

\section{\label{sec:IntegralCalculus}Integral calculus}

\noindent It is now natural to study the existence of primitives of
generic smooth functions $f:[a,b]\freccia\ER$ and hence the existence
of an integration theory. We will tackle this problem firstly for
$a$, $b\in\R$, then for $a=-\infty$ and $b=+\infty$, and finally
for $a$, $b\in\ER$. Like in SDG, the problem is solved proving existence
and uniqueness of the simplest Cauchy initial value problem.

\noindent We firstly recall our notations for intervals, e.g. $(a,b]:=\{x\in\ER\,|\, a<x\le b\}$,
whereas if $a$, $b\in\R$, then $(a,b]_{\R}:=(a,b]\cap\R$. Using
Theorem \ref{thm:equivalentFormulationForOrderRelation} it is not
hard to prove that if $a$, $b\in\R$\begin{align*}
\ext{\left\{ (a,b)_{\R}\right\} } & =(a,b)\\
\ext{\left\{ [a,b]_{\R}\right\} } & \subsetneqq[a,b],\end{align*}
for example, $x_{t}:=a-t^{2}$ is equal to $a$ in $\ER$, and hence
it belongs to the interval $[a,b]$, but $x\notin\ext{\left\{ [a,b]_{\R}\right\} }$
because $x$ does not map $\R_{\ge0}$ into $[a,b]_{\R}$. We also
recall that there can be any order relationship between a Fermat number
$a\in\ER$ and its standard part: e.g. $a=\st{a}-\diff{t}<\st{a}$
whereas $a=\st{a}+\diff{t}>\st{a}$. For this reason, a general inclusion
relationship between the interval $(a,b)$ and the interval $(\st{a},\st{b})$
does not hold, even if $\st{a}<\st{x}<\st{b}$ implies $a<x<b$.

To solve the problem of existence and uniqueness of primitives, we
need two preliminary results. The first one is called by \citet{Bel}
the \emph{constancy principle}.
\begin{lem}
\noindent \label{lem:constancyPrinciple}Let $a$, $b\in\R$ with
$a<b$, and $f:(a,b)\freccia\ER$ a $\ECInfty$ function such that\[
f'(x)=0\quad\forall x\in(a,b).\]
Then $f$ is constant.
\end{lem}
\noindent \textbf{Proof:} Let $x$, $y\in(a,b)$ and $h:=y-x$. We
can suppose $h>0$, otherwise we can repeat the proof exchanging the
role of $x$ and $y$. So, we have that $\overrightarrow{[x,x+h]}=\overrightarrow{[x,y]}=[x,y]\subseteq(a,b)$
because $a<x<y<b$, therefore $(x,h)\in\widetilde{(a,b)}=\widetilde{\ext{(a,b)_{\R}}}$.
Using the smooth incremental ratio (Theorem \ref{thm:smoothIncrementalRatio})
we get\begin{equation}
f(y)=f(x)+h\cdot f'[x,h].\label{eq:2_constancyPrinciple}\end{equation}
As proved in Theorem \ref{thm:smoothIncrementalRatio}, we can always
find a smooth function $\alpha$ and a parameter $p\in\ER^{\sf p}$
such that\begin{equation}
f'[x,h]_{t}=\int_{0}^{1}\partial_{2}\alpha(p_{t},x_{t}+s\cdot h_{t})\,\diff{s}.\label{eq:1_constancyPrinciple}\end{equation}
But for every $s\in[0,1]_{\R}$ we have that\[
x+s\cdot h=x+s\cdot(y-x)\in[x,y]\subseteq(a,b),\]
so that $f'(x+s\cdot h)=0$, i.e. $\partial_{2}\alpha(p,x+s\cdot h)=0$
in $\ER$. Written in explicit form this means\[
\lim_{t\to0^{+}}\frac{\partial_{2}\alpha(p_{t},x_{t}+s\cdot h_{t})}{t}=0.\]
From \eqref{eq:1_constancyPrinciple} using dominated convergence
we have\begin{align*}
\lim_{t\to0^{+}}\frac{f'[x,h]_{t}}{t} & =\lim_{t\to0^{+}}\frac{1}{t}\cdot\int_{0}^{1}\partial_{2}\alpha(p_{t},x_{t}+s\cdot h_{t})\,\diff{s}\\
 & =\int_{0}^{1}\lim_{t\to0^{+}}\frac{\partial_{2}\alpha(p_{t},x_{t}+s\cdot h_{t})}{t}\,\diff{s}=0,\end{align*}
that is $f'[x,h]=0$ in $\ER$ and hence $f(y)=f(x)$ from \eqref{eq:2_constancyPrinciple}.$\qedNoNewLine$

The second preliminary result permits to extend the validity of an
equality from an open interval $(a,b)$ to its borders.
\begin{lem}
\noindent \label{lem:extesionOfEquationToTheBorders}Let $a$, $b\in\R$,
with $a<b$, $c\in\ER$, and $f:[a,b]\freccia\ER$ a $\ECInfty$ function
such that\[
\forall x\in(a,b)\pti f(x)=c.\]
Then\[
f(a)=f(b)=c\]

\end{lem}
\noindent \textbf{Proof:} We prove that $f(a)=c$, analogously we
can proceed for $f(b)=c$. Let us write the function $f$ as the parametrized
extension of a smooth function in a neighborhood of $x=a$:\begin{equation}
f(x)=\alpha(p,x)\quad\forall x\in\ext{V}\cap[a,b],\label{eq:3_extensionToTheBorders}\end{equation}
where $V$ is open in $\R$ and $a\in\ext{V}$. 

Let $\rho:\widetilde{\ext{V}}\freccia\ER$ be the incremental ratio
of $\alpha(p,-):\ext{V}\freccia\ER$:\begin{equation}
\alpha(p,x+h)=\alpha(p.x)+h\cdot\rho(x,h)\quad\forall(x,h)\in\widetilde{\ext{V}}.\label{eq:2_extensionToTheBorders}\end{equation}
Since $a=\st{a}\in V$ we can find a $\delta\in\R_{>0}$ such that
$(a-2\delta,a+2\delta)_{\R}\subseteq V$ and with $a+\delta<b$. Then
$[a,a+\delta]\subseteq(a-2\delta,a+2\delta)=\ext{(a-2\delta,a+2\delta)_{\R}}\subseteq\ext{V}$,
i.e. $(a,\delta)\in\widetilde{\ext{V}}$ and we can hence use the
previous \eqref{eq:2_extensionToTheBorders} and\eqref{eq:3_extensionToTheBorders}
obtaining\[
\alpha(p,a+\delta)=\alpha(p,a)+\delta\cdot\rho(a,\delta).\]
\begin{equation}
f(a+\delta)=f(a)+\delta\cdot\rho(a,\delta),\label{eq:4_extensionToTheBorders}\end{equation}
because $a+\delta<b$, $[a,a+\delta]\subseteq\ext{V}$.

But we know that it is always possible to take the smooth incremental
ratio $\rho$ so that\begin{equation}
\rho(a,\delta)_{t}=\int_{0}^{1}\partial_{2}\alpha(p_{t},a_{t}+s\cdot\delta)\,\diff{s}.\label{eq:1_extensionToTheBorders}\end{equation}
But $f$ is constant on $(a,b)$, so that for every $s\in[0,1]_{\R}$
we have $a+s\cdot\delta\in(a,b)$ and hence $f'(a+s\cdot\delta)=\partial_{2}\alpha(p,a+s\cdot\delta)=0$
in $\ER$, i.e.\[
\lim_{t\to0^{+}}\frac{\partial_{2}\alpha(p_{t},a_{t}+s\cdot\delta)}{t}=0\quad\forall s\in[0,1]_{\R}.\]
From this and from \eqref{eq:1_extensionToTheBorders}, using dominated
convergence we get\[
\lim_{t\to0^{+}}\frac{\rho(a,\delta)_{t}}{t}=0,\]
that is $\rho(a,\delta)=0$ in $\ER$. Finally, from this and from
\eqref{eq:4_extensionToTheBorders} we obtain the conclusion: $f(a+\delta)=f(a)=c$.$\qedNoNewLine$

We can now prove existence and uniqueness of primitives in the first
case of domains $[a,b]$ with real boundaries.
\begin{thm}
\label{thm:existenceOfIntegralForRealBoundaries}Let $a$, $b\in\R$
with $a<b$, $f:[a,b]\freccia\ER$ a $\ECInfty$ function and $u\in[a,b]$.
Then there exists one and only one $\ECInfty$ map\[
I:[a,b]\freccia\ER\]
such that\[
I'(x)=f(x)\quad\forall x\in(a,b)\]
\[
I(u)=0\]

\end{thm}
\noindent \textbf{Proof:} We can prove the existence assuming that
$u=a$; in fact, if $I'=f$ on $(a,b)$ and $I(a)=0$, then $J(x):=I(x)-I(u)$
verifies $J'=f$ on $(a,b)$ and $J(u)=0$.

For every $x\in[a,b]$ we can write\[
f|_{\mathcal{V}_{x}}=\alpha_{x}(p,-)|_{\mathcal{V}_{x}}\]
for suitable $p^{x}\in\ER^{\sf{p}^{x}}$, $U_{x}$ open subset of
$\R^{\sf{p}^{x}}$ such that $p^{x}\in\ext{U_{x}}$, $V_{x}$ open
in $\R$ such that $x\in\ext{V_{x}}\cap[a,b]=:\mathcal{V}_{x}$ and
$\alpha_{x}\in\CInfty(U_{x}\times V_{x},\R)$. We can assume that
the open sets $V_{x}$ are of the form $V_{x}=(x-\delta_{x},x+\delta_{x})_{\R}$
for a suitable $\delta_{x}>0$.

The idea is to patch together suitable integrals of the functions
$\alpha_{x}(p^{x},-)$. The problem in this idea is that we have to
realize the condition $I(a)=0$, which forces us to patch together
integrals which are {}``each one is the extension of the previous
one'', i.e. on the intersection of their domains two integrals must
have the same value at one point, so that we can prove they are equal
on the whole intersection. Moreover, we must use the compactness of
the interval $[a,b]_{\R}$ because, generally speaking, a smooth function
can be non integrable in an open set.

We have that $\left(V_{x}\right)_{x\in[a,b]_{\R}}$ is an open cover
of the real interval $[a,b]_{\R}$, thus we can cover $[a,b]_{\R}$
with a finite number of $V_{x}$, that is we can find $x_{1},\ldots,x_{n}\in[a,b]_{\R}$
such that $\left(V_{x_{i}}\right)_{i=1,\ldots,n}$ is an open cover
of $[a,b]_{\R}$. We will use simplified notations like $V_{i}:=V_{x_{i}}$,
,$\delta_{i}:=\delta_{x_{i}}$, $\alpha_{i}:=\alpha_{x_{i}}$, etc.

We can always suppose to have chosen the indexes $i=1,\ldots,n$ and
the amplitudes $\delta_{i}>0$ such that\[
a=x_{1}<x_{2}<\ldots<x_{n}=b\]
\[
x_{i}-\delta_{i}<x_{i+1}-\delta_{i+1}<x_{i}+\delta_{i}<x_{i+1}+\delta_{i+1}\quad\forall i=1,\ldots,n-1,\]
in this way the intervals $V_{i}=(x_{i}-\delta_{i},x_{i}+\delta_{i})_{\R}$
intersect in the sub-interval $(x_{i+1}-\delta_{i+1},x_{i}+\delta_{i})_{\R}$.

\noindent %
\begin{figure}[H]
\begin{centering}
\includegraphics[scale=0.68]{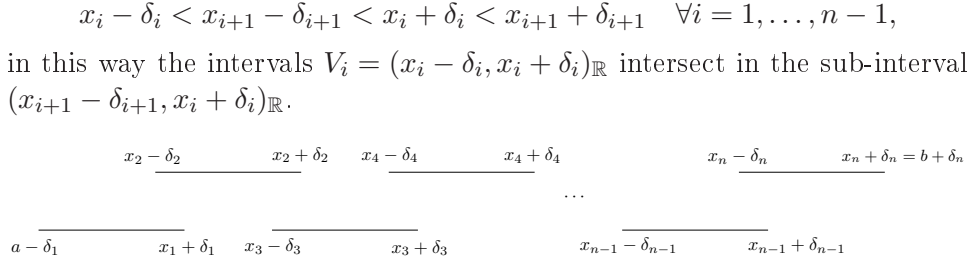}
\par\end{centering}

\caption{Intervals for the recursive definition of a primitive}

\end{figure}

For any $i=1,\ldots,n$ let us choose a point in this sub interval
$\bar{x}_{i}\in(x_{i+1}-\delta_{i+1},x_{i}+\delta_{i})_{\R}$ (these
are the points in the intersections of the domains of the integrals
we are going to define and mentioned in the previous intuitive sketch
of the ideas of this proof).

In fact, let us define, recursively:\[
I_{1}(x)_{t}:=\int_{a}^{x_{t}}\alpha_{1}(p_{t}^{1},s)\,\diff{s}\quad\forall x\in\ext{V_{1}}\]
\[
I_{k+1}(x)_{t}:=\int_{\bar{x}_{k}}^{x_{t}}\alpha_{k+1}(p_{t}^{k+1},s)\,\diff{s}+I_{k}(\bar{x}_{k})\quad\forall x\in\ext{V_{k+1}}\ \forall k=1,\ldots,n-1.\]
Every $I_{k}$ is a $\ECInfty$ function defined on $\ext{V_{k}}$,
and moreover\[
I'_{k}(x)=\alpha_{k}(p^{k},x)=f(x)\quad\forall x\in\ext{V_{k}}.\]
Therefore, in a generic point in the intersection \begin{align*}
\ext{V_{k}}\cap\ext{V_{k+1}} & =\ext{(V_{k}\cap V_{k+1})}=\ext{\left[(x_{k+1}-\delta_{k+1},x_{k}+\delta_{k})_{\R}\right]}\\
 & =(x_{k+1}-\delta_{k+1},x_{k}+\delta_{k})\end{align*}
 we have\begin{equation}
I'_{k}(x)=f(x)=I'_{k+1}(x)\label{eq:1_existenceOfIntegralForRealBoundaries}\end{equation}
\[
\left(I_{k}-I_{k+1}\right)'(x)=0\]
\[
\left(I_{k}-I_{k+1}\right)(\bar{x}_{k})=0,\]
so, from Theorem \ref{lem:constancyPrinciple} it follows $I_{k}=I_{k+1}$
on $(x_{k+1}-\delta_{k+1},x_{k}+\delta_{k})$. We can hence use the
sheaf property of the space $[a,b]$ with the open cover $\left(\ext{V_{k}}\cap[a,b]\right)_{i=1,\ldots,n}$
to patch together the functions $I_{k}|_{\ext{V_{k}}\cap[a,b]}$ obtaining
the map $I:[a,b]\freccia\ER$. This function satisfies the conclusion
of the statement because of \eqref{eq:1_existenceOfIntegralForRealBoundaries}
and because of the equalities $I(a)=I_{1}(a)=0$.

To prove the uniqueness, let us suppose that $J$ verifies $J'=f$
on $(a,b)$ and $J(u)=0$, then using again Theorem \ref{lem:constancyPrinciple}
we have that $(J-I)|_{(a,b)}$ is constant and equal to zero. Finally,
using Lemma \ref{lem:extesionOfEquationToTheBorders}, we can extend
this constancy to the whole closed interval $[a,b]$.$\qedNoNewLine$

The second case is for domains $[a,b]=\ER$.
\begin{thm}
\noindent If $f:\ER\freccia\ER$ is smooth and $u\in\ER$ then there
exists one and only one smooth $I:\ER\freccia\ER$ such that $I'=f$
and $I(u)=0$.
\end{thm}
\noindent \textbf{Proof:} For every $k\in\N_{>0}$ let us define\[
f_{k}:=f|_{[a-k,a+k]},\]
where $a:=\st{u}$. Moreover, define\[
I_{k}(x):=\int_{u}^{x}f_{k}\quad\forall x\in(a-k,a+k)\subseteq[a-k,a+k].\]
Therefore, we have that setting $V_{k}:=(a-k,a+k)=\ext{(a-k,a+k)_{\R}}$,
we obtain that $\left(V_{k}\right)_{k>0}$ is an open cover of $\ER$.
Moreover, $I_{k}'(a)=f_{k}(x)=f(x)$ for every $x\in V_{k}$, so that
$I_{k}$ and $I_{j}$ coincide in $V_{k}\cap V_{j}$. From the sheaf
property of $\ER$ we get\[
\exists!\, I:\ER\freccia\ER\text{ smooth}\pti I|_{V_{k}}=I_{k}\quad\forall k\in\N_{>0}.\]

Now, let us note that\[
\forall x,h\in\ER\ \exists\, k\in\N_{>0}\pti\overrightarrow{[x,x+h]}\subseteq V_{k},\]
then we also have\[
I(x+h)=I_{k}(x+h)=I(x)+h\cdot I'[x,h]=I_{k}(x)+h\cdot I'_{k}[x,h],\]
so that the smooth incremental ratios of $I$ and $I_{k}$ are equal,
i.e. $I'[x,h]=I'_{k}[x,h]$. Thus, $I'(x)=I'_{k}(x)=f(x)$, and finally
$I(u)=I_{1}(u)=0$.

This proves the existence part. The uniqueness follows from Lemma
\ref{lem:constancyPrinciple}.$\qedNoNewLine$

To extend Theorem \ref{thm:existenceOfIntegralForRealBoundaries}
to non standard boundaries $a$, $b\in\ER$ we need the following
result.
\begin{lem}
\label{lem:smoothFnctsOnClosedIntervalsCanBeExtended}Let $a$, $b\in\ER$
with $\st{a}<\st{b}$, and $f:[a,b]\freccia\ER$ a $\ECInfty$ function.
Then there exist $\delta\in\R_{>0}$ and a $\ECInfty$ function $\bar{f}:(a-\delta,b+\delta)\freccia\ER$
such that\[
\bar{f}|_{[a,b]}=f\]

\end{lem}
\noindent \textbf{Proof:} As usual, let us write the function $f$
as the parametrized extension of an ordinary smooth function in a
neighborhood of $x=a$:\begin{equation}
f(x)=\alpha(p,x)\quad\forall x\in\ext{V}\cap[a,b],\label{eq:1_smoothFnctsOnClosedIntervalsCanBeExtended}\end{equation}
where $V$ is open in $\R$ and $a\in\ext{V}$ so that $\st{a}\in V$.

We can make the same in a neighborhood of $x=b$:\begin{equation}
f(x)=\beta(q,x)\quad\forall x\in\ext{U}\cap[a,b],\label{eq:2_smoothFnctsOnClosedIntervalsCanBeExtended}\end{equation}
where $U$ is open in $\R$ and $b\in\ext{U}$ so that $\st{b}\in U$.

Because $U$ and $V$ are open subsets of $\R$, we can always suppose
to have $U=(\st{b}-\eta,\st{b}+\eta)_{\R}$ and $V=(\st{a}-\eta,\st{a}+\eta)_{\R}$
with $\eta\in\R_{>0}$ such that $\st{a}+\eta<\st{b}-\eta$ because
$\st{a}<\st{b}$. Therefore, we have\[
\ext{V}\cap(\st{a},\st{b})=(\st{a}-\eta,\st{a}+\eta)\cap(\st{a},\st{b})\subseteq\ext{V}\cap[a,b]\]
\[
\ext{U}\cap(\st{a},\st{b})=(\st{b}-\eta,\st{b}+\eta)\cap(\st{a},\st{b})\subseteq\ext{U}\cap[a,b]\]
\[
\ext{V}\cap\ext{U}=(\st{a}-\eta,\st{a}+\eta)\cap(\st{b}-\eta,\st{b}+\eta)=\emptyset,\]
so that any two of the following smooth functions\[
\alpha(p,-):\ext{V}\freccia\ER\]
\[
f|_{(\st{a},\st{b})}:(\st{a},\st{b})\freccia\ER\]
\[
\beta(q,-):\ext{U}\freccia\ER\]
are equal on the intersection of their domains for \eqref{eq:1_smoothFnctsOnClosedIntervalsCanBeExtended}
and \eqref{eq:2_smoothFnctsOnClosedIntervalsCanBeExtended}.

For the sheaf property of $(\st{a}-\eta,\st{a}+\eta)\cup(\st{a},\st{b})\cup(\st{b}-\eta,\st{b}+\eta)=(\st{a}-\eta,\st{b}+\eta)$
we thus have\[
\exists!\, g:(\st{a}-\eta,\st{b}+\eta)\freccia\ER\text{ smooth}\pti g|_{(\st{a},\st{b})}=f.\]
If we set $\delta:=\frac{\eta}{2}$ we have that $\st{a}-\eta<a-\delta<b+\delta<\st{b}+\eta$,
as we can verify considering the standard parts of all these numbers,
and hence $\bar{f}:=g|_{(a-\delta,b+\delta)}$ verifies\[
\bar{f}|_{(\st{a},\st{b})}=f.\]
Because $(\st{a},\st{b})\subseteq[a,b]$ we have to verify that the
function $\bar{f}$ and the function $f$ also coincide on $[a,b]\setminus(\st{a},\st{b})$.
We firstly note that $\st{a}\in V\subseteq\ext{V}$ and $\st{b}\subseteq U\subseteq\ext{U}$,
so we can apply \eqref{eq:1_smoothFnctsOnClosedIntervalsCanBeExtended}
and \eqref{eq:2_smoothFnctsOnClosedIntervalsCanBeExtended} at $x=\st{a}$
and $x=\st{b}$ too. Therefore, we have $\bar{f}|_{[\st{a},\st{b}]}=f$.
Secondly, if $x\in[a,b]\setminus(\st{a},\st{b})$, then either $a\le x\le\st{a}$
or $\st{b}\le x\le b$; we will deal with the first case, the second
being analogous. From these inequalities, it follows that $\st{x}=\st{a}$
so that from the infinitesimal Taylor's formula we get\begin{align*}
f(x) & =f\left[\st{x}+(x-\st{x})\right]=\sum_{i=0}^{n}\partial_{2}^{(i)}\alpha(p,\st{a})\cdot\frac{(x-\st{a})^{i}}{i!}=\\
 & =\bar{f}\left[\st{x}+(x-\st{x})\right]=\bar{f}(x).\end{align*}
$\qedWithFinalEq$
\begin{thm}
\label{thm:existenceOfIntegralForNonStdBoundaries}Let $a$, $b\in\ER$
with $\st{a}<\st{b}$, $f:[a,b]\freccia\ER$ a $\ECInfty$ function
and $u\in[a,b]$. Then there exists one and only one $\ECInfty$ map\[
I:[a,b]\freccia\ER\]
such that\[
I'(x)=f(x)\quad\forall x\in(a,b)\]
\[
I(u)=0\]

\end{thm}
\noindent \textbf{Proof:} From Lemma \ref{lem:smoothFnctsOnClosedIntervalsCanBeExtended}
there exist a $\delta\in\R_{>0}$ and a smooth function $\bar{f}:(a-\delta,b+\delta)\freccia\ER$
such that $\bar{f}|_{[a,b]}=f$.

But $a-\delta<\st{a}-\frac{\delta}{2}<\st{b}+\frac{\delta}{2}<b+\delta$,
so the interval with real boundaries $[\alpha,\beta]:=\left[\st{a}-\frac{\delta}{2},\st{b}+\frac{\delta}{2}\right]$
is contained in $(a-\delta,b+\delta)$. Finally $\alpha=\st{a}-\frac{\delta}{2}<a\le u\le b<\st{b}+\frac{\delta}{2}=\beta$
and we can thus apply Theorem \ref{thm:existenceOfIntegralForRealBoundaries}
obtaining existence and uniqueness of the primitive $J$ of the function
$\bar{f}|_{[\alpha,\beta]}$ such that $J(u)=0$. But $[a,b]\subseteq[\alpha,\beta]$
and hence $I:=J|_{[a,b]}$ verifies the existence part of the conclusion.
The uniqueness part follows, in the usual way, from Lemma \ref{lem:constancyPrinciple}
and Lemma \ref{lem:extesionOfEquationToTheBorders}.$\qedNoNewLine$

We can now define
\begin{defn}
\label{def:integral}Let $a$, $b\in\ER$ with $\st{a}<\st{b}$. Moreover,
let us consider a $\ECInfty$ function $f:[a,b]\freccia\ER$ and a
point $u\in[a,b]$. Then
\begin{enumerate}
\item ${\displaystyle \int_{u}^{(-)}f:=\int_{u}^{(-)}f(s)\,\diff{s}:[a,b]\freccia\ER}$
\item ${\displaystyle \int_{u}^{u}f=0}$
\item ${\displaystyle \forall x\in(a,b)\pti\left(\int_{u}^{(-)}f\right)'(x)=\frac{\diff{}}{\diff{x}}\int_{u}^{x}f(s)\,\diff{s}=f(x)}$
\end{enumerate}
\end{defn}
\noindent It is important to note that in this way we obtain a generalization
of the usual notion of integral. Indeed, for $a$, $b$, $u\in\R$
with $a<u<b$, let us consider a standard smooth function\[
f\in\CInfty([a,b]_{\R},\ER).\]
Let us extend smoothly this function on an open interval $(a-\delta,b+\delta)_{\R}$
with $\delta\in\R_{>0}$, so that outside $[a,b]_{\R}$ the extension
of $f$ is constant. Let us consider the smooth function\[
I:={}^{{}^{{}^{{}^{{\scriptstyle \bullet}}}}}\left(\int_{u}^{(-)}f(s)\diff{s}\right):(a-\delta,b+\delta)\freccia\ER,\]
where here the integral symbol has to be understood as the classical
Riemann integral on $\R$. Now we have that\[
[a,b]\subseteq(a-\delta,b+\delta)\]
so that we can consider the restriction $I|_{[a,b]}$. It is not hard
to prove that this restriction verifies all the properties of the
previous Definition \ref{def:integral} for the function $\ext{f}$,
but at the same time, because it is the extension of a classical integral,
it also verifies\[
\forall x\in[a,b]_{\R}\pti I(x)=\int_{u}^{x}f(s)\diff{s}\in\R.\]

These theorems can be used to try a generalization of several results
of integral (smooth!) calculus to $\ECInfty$ functions of the form
$f:\prod_{i=1}^{n}[a_{i},b_{i}]\freccia\ER^{d}$, with $\st{a_{i}}<\st{b_{i}}$.

Indeed, this can be done for several theorems. We only list here the
main results that we have already proved, leaving a complete report
of them for a subsequent work. For most of them the proofs are very
similar to the analogous presented e.g. in \citet{Lav}:
\begin{enumerate}
\item property of linearity of integrals,
\item fundamental theorem of calculus,
\item integration by parts formula,
\item formulas of the form $\int_{u}^{v}f+\int_{v}^{w}f=\int_{u}^{w}f$,
$\int_{u}^{v}f=-\int_{v}^{u}f$,
\item integration formula by change of variable,
\item derivation under the integral sign,
\item smoothness of the function $(f,u,v)\in\ER^{[a,b]}\times[a,b]^{2}\mapsto\int_{u}^{v}f\in\ER$,
\item majorization of integrals: if $|fx|\le M$ for every $x\in[a,b]$,
then $|\int_{a}^{b}f|\le M\cdot(b-a)$,
\item majorization of $d$-dimensional integrals $\Vert\int_{a}^{b}f\Vert\le M\cdot\sqrt{d}\cdot(b-a)$
if $\Vert f(x)\Vert\le M$ for every $x\in[a,b]$ and $f:[a,b]\freccia\ER^{d}$,
\item Fubini theorem for double integrals.\end{enumerate}
\begin{example*}
$\phantom{a}$
\begin{enumerate}
\item \noindent \textbf{Divergence and curl.} Classically the $\text{div}\vec{A}(x)$
is the density of the flux of $\vec{A}\in\CInfty(U,\R^{3})$ through
an {}``infinitesimal parallelepiped'' centered at $x\in U\subseteq\R^{3}$.
To formalize this concept we take three vectors $\vec{h}_{1}$, $\vec{h}_{2}$,
$\vec{h}_{3}\in\ER^{3}$ and express them with respect to a fixed
base $\vec{e}_{1}$, $\vec{e}_{2}$, $\vec{e}_{3}\in\R^{3}$: \[
\vec{h}_{i}=k_{i}^{1}\cdot\vec{e}_{1}+k_{i}^{2}\cdot\vec{e}_{2}+k_{i}^{3}\cdot\vec{e}_{3}\ee{where}k_{i}^{j}\in\ER.\]
We say that $P:=(x,\vec{h}_{1},\vec{h}_{2},\vec{h}_{3})$ is a \emph{(first
order) infinitesimal parallelepiped} if \begin{align*}
{} & x\in\R^{3}\\
{} & \forall i,j,k=1,2,3\pti k_{i}^{1}\cdot k_{i}^{2}\cdot k_{i}^{3}\in D.\end{align*}
The flux of the vector field $\vec{A}$ through such a parallelepiped
(toward the outer) is by definition the sum of the fluxes through
every {}``face'' \begin{align*}
\int_{P}\vec{A}\boldsymbol{\cdot}\vec{n}\,\diff{S} & :=\int_{0}^{1}\diff{t}\int_{0}^{1}\vec{A}(x+t\vec{h}_{1}+s\vec{h}_{2})\boldsymbol{\cdot}\vec{h}_{2}\times\vec{h}_{1}\,\diff{s}+\\
{} & \int_{0}^{1}\diff{t}\int_{0}^{1}\vec{A}(x+\vec{h}_{3}+t\vec{h}_{1}+s\vec{h}_{2})\boldsymbol{\cdot}\vec{h}_{1}\times\vec{h}_{2}\,\diff{s}+\ldots,\end{align*}
where the \ldots{} indicate similar terms for the other faces of
the parallelepiped. Let us note that e.g. the function $s\mapsto\vec{A}(x+t\vec{h}_{1}+s\vec{h}_{2})$
is a $\ECInfty$ arrow of type $\ext{\alpha}(p,s)$, where here the
parameter is $p=(x,t,\vec{h}_{1},\vec{h}_{2})\in\ER^{10}$. We have
hence concrete examples of $\ECInfty$ functions to which we can apply
the results of the previous sections. Now, it is easy to prove that
if $\vec{A}\in\CInfty(U;\R^{3})$ and ${\rm Vol}(\vec{h}_{1},\vec{h}_{2},\vec{h}_{3})$,
i.e. the oriented volume of the infinitesimal parallelepiped $P=(x,\vec{h}_{1},\vec{h}_{2},\vec{h}_{3})$,
is not zero, then the following ratio between first order infinitesimals
exists and is independent by $\vec{h}_{1}$, $\vec{h}_{2}$, $\vec{h}_{3}$:
\[
{\rm div}\vec{A}(x):=\frac{1}{{\rm Vol}(\vec{h}_{1},\vec{h}_{2},\vec{h}_{3})}\cdot\int_{P}\vec{A}\boldsymbol{\cdot}\vec{n}\,\diff{S}.\]
To define the curl of a vector field $\vec{A}\in\ECInfty(U,\R^{3})$
we can say that $C:=(x,\vec{h}_{1},\vec{h}_{2})$ is a \emph{(first
order) infinitesimal cycle }if \[
x\in U\ee{and}\forall p,q=1,2,3\pti\sum_{i,j=1}^{3}\vert k_{i}^{p}\cdot k_{j}^{q}\vert\in D.\]
The circulation of the vector field $\vec{A}$ on this cycle $C$
is defined as the sum of the {}``line integrals'' on every {}``side'':
\[
\int_{C}\vec{A}\boldsymbol{\cdot}\vec{t}\,\diff{l}:=\int_{0}^{1}\vec{A}(x+t\vec{h}_{1})\boldsymbol{\cdot}\vec{h}_{1}\,\diff{t}+\int_{0}^{1}\vec{A}(x+\vec{h}_{1}+t\vec{h}_{2})\boldsymbol{\cdot}\vec{h}_{2}\,\diff{t}+\ldots,\]
where \ldots{} indicates similar terms for the otehr side of the
cycle $C$. Once again, using exactly the calculations frequently
done in elementary courses of physics, one can prove that there exists
one and only one vector, $\text{{\rm curl}}\,\vec{A}(x)\in\R^{3}$,
such that \[
\int_{C}\vec{A}\boldsymbol{\cdot}\vec{t}\,\diff{l}=\text{{\rm curl}}\vec{A}(x)\boldsymbol{\cdot}\vec{h}_{1}\times\vec{h}_{2}\]
for every infinitesimal cycle $C=(x,\vec{h}_{1},\vec{h}_{2})$, representing
thus the (vector) density of the circulation of $\vec{A}$.
\item \textbf{Limits in $\mathbf{\ER}$}. Because the theory of Fermat reals
is not an alternative way for the foundation of calculus, but a rigorous
way to have at disposal infinitesimal methods, there is no need to
think that the notion of limit expressed by Weierstrass' $\eps-\delta$'s
is conceptually incompatible with our use of infinitesimals. A similar
approach is already used, e.g. in the study of the Levi-Civita field
(see Appendix \ref{app:theoriesOfInfinitesimals} and references therein).
Therefore, we can introduce the following:

\begin{defn}
\label{def:limit}Let $f:U\freccia\ER$ be a $\ECInfty$ function
defined in $U\subseteq\ER$ and $l$, $\bar{x}\in\ER$ be two Fermat
reals. Then we say that $l$ \emph{is the limit of $f(x)$ for $x\to\bar{x}$}
if and only if\[
\forall\eps\in\ER_{>0}\ \exists\,\delta\in\ER_{>0}:\ \forall x\in U\pti0<\left|x-\bar{x}\right|<\delta\Rightarrow\left|f(x)-l\right|<\eps\]

\end{defn}
\noindent Analogously we can define the right and the left limit.

Using the total order on $\ER$ and replicating the standard proof,
we can prove that this limit, if it exists, is unique.
\begin{thm}
In the hypothesis of the previous Definition \ref{def:limit}, there
exists at most one $l\in\ER$ such that $l$ is the limit of $f(x)$
for $x\to\bar{x}$. In this case, we will use the notation ${\displaystyle l=\lim_{x\to\bar{x}}f(x)}$.
\end{thm}
If $f:(a,b)\freccia\ER$, $\st{a}<\st{b}$, and $\st{\bar{x}}\in(\st{a},\st{b})$,
we want to prove that ${\displaystyle \lim_{x\to\bar{x}}f(x)=f(\bar{x})}$.
Let us consider a generic $\eps>0$, we want to find a $\delta\in D_{>0}$.
From the inequalities $0<|x-\bar{x}|<\delta$ it follows that $x-\bar{x}\in D$
and hence from the first order Taylor's formula\[
\left|f(x)-f(\bar{x})\right|=\left|f'(\bar{x})\cdot(x-\bar{x})\right|\le\left|f'(\bar{x})\right|\cdot\delta.\]
If $f'(\bar{x})\in D_{\infty}$, then $\left|f'(\bar{x})\right|\cdot\delta=0<\eps$
because $\delta\in D$ is a first order infinitesimal. Otherwise,
$f'(\bar{x})$ is invertible and it suffices to fix $\delta$ such
that\[
\delta<\frac{\eps}{\left|f'(\bar{x})\right|},\]
e.g.\[
\delta:=\min\left\{ \diff{t},\frac{\eps}{2\left|f'(\bar{x})\right|}\right\} \in D_{>0}.\]
This expected result (even if the topology we considered on the Fermat
reals has not been defined as the one induced by the absolute value,
but the natural topology induced by the smooth figures of $\ER$;
see the Definition \ref{def:topologyOnTheCartesianClosure}) says
us that in the context of $\ECInfty$ functions, the notion of limit
is interesting only at the border points $\bar{x}=a$ or $\bar{x}=b$
on which the function $f$ is not defined. From this point of view,
lemmas \ref{lem:extesionOfEquationToTheBorders} and \ref{lem:smoothFnctsOnClosedIntervalsCanBeExtended}
represent possible substitutes of the notion of limit in $\ER$.\end{enumerate}
\end{example*}

\chapter{\label{cha:CalculusOnInfinitesimalDomains}Calculus on infinitesimal
domains}

It is natural to expect that we cannot restrict our differential calculus
to smooth functions defined on open sets, but that we have to extend
the notion of derivatives to functions defined on infinitesimal sets,
e.g. $0\in I\subseteq D_{\infty}$.

As we prompted above, the infinitesimal Taylor's formula does not
uniquely identifies the derivatives appearing in its addends, so that
we must use the map $\iota_{k}$ to consider the simplest numbers
that verify a given Taylor's formula.

\section{\label{sec:TheGeneralizedTaylorFormula}The generalized Taylor's
formula}

In this section we want to prove the Taylor's formula for functions
defined on an infinitesimal domain, like e.g. $f:D_{\alpha}\freccia\ext{X}$,
with $\alpha\in\R_{>0}$ and $X\in\CInfty$. The possibility to prove
the following theorems has been the first motivation to choose little-oh
polynomials instead of the more general nilpotent functions (like
in \citet{Gio3}) to define $\ER$. A stronger algebraic control on
the properties of nilpotent infinitesimals, and better order properties
have been the second motivation.

We start proving some preliminary results that permit to affirm that
if $f(0)\in\ext{U}$, where $U$ is open in the space $X$, then $f(h)\in\ext{U}$
for every $h$ in the infinitesimal domain of $f$.
\begin{lem}
\label{lem:stdPartOf_x_inUand_x_inExt_U}Let $X$ be a $\CInfty$
space, $U$ an open subset of $X$ and $x\in\ext{X}$, then\[
\st{x}\in U\then x\in\ext{U}.\]

\end{lem}
\noindent Let us note that we have already frequently used the analogous
of this result for spaces of the form $X=\R^{s}$, but in this particular
case the notion of little-oh polynomial does not depend on observables
but only on the norm of $\R^{s}$. For this reason, in this particular
situation, the passage from $x\in\ext{X}$ to $x\in\ext{U}$ is trivial.

\noindent \textbf{Proof:} Because $\st{x}=x(0)\in U\in\Top{U}$ and
$x$ is continuous at $t=0$, we have that locally $x$ has values
in $U$, i.e.\[
\forall^{0}t:\ x_{t}\in U.\]
So, because all the properties we are considering are local, we can
assume that $x:\R_{\ge0}\freccia U$, i.e. $x$ has globally values
in $U$. To prove that $x\in\ext{U}$ it remains to prove that%
\footnote{Let us recall the general definition of the set of all the little-oh
polynomials in the space $X\in\CInfty$, i.e. the Definition \ref{def:LittleOhPolynomialsRd}.%
} $x\in U_{o}[t]$, where we are meaning $U=(U\prec X)$, that is on
the subset $U$ the structure induced by the superspace $X$. So,
let us consider a zone $VK$ of $(U\prec X)$ such that $x(0)\in V$
and an observable $\phi\In{VK}(U\prec X)$. We have that $V\in\Top{(U\prec X)}\subseteq\Top{X}$
because $U$ is open in $X$, and hence $VK$ is a zone of the space
$X$ too. Moreover\[
(V\prec(U\prec X))=(V\prec X)\xfrecciad{\phi}K\]
and hence $\phi$ is also an observable of $X$ such that $\phi\In{VK}X$.
But, by hypotheses $x\in\ext{X}$, so that $x\in X_{o}[t]$ and hence
$\phi\circ x\in\R_{o}^{\sf k}[t]$, which is the conclusion.$\qedNoNewLine$

The main aim of this section is to prove an infinitesimal Taylor's
formula for functions of the form $f:D_{\alpha}\freccia\ext{X}$ through
the composition with observables $\varphi\InUp{UK}X$. Precisely we
want to consider a $\ECInfty$ function $f:D_{\alpha}\freccia\ext{X}$
with $f(0)\in\ext{U}:=\ext{(U\prec X)}$ (in general, the function
$f$ will not be the extension of a classical one, that is $f$ is
not necessarily of the form $f=\ext{g}|_{D}$) and we will prove the
Taylor's formula for the function $\ext{\varphi}(f(-)):D_{\alpha}\freccia K\subseteq\ER^{\sf k}$.
First of all, we prove that this composition is well defined, that
is the following theorem holds:
\begin{thm}
\label{thm: DandOpenGeneralized} Let $X$ be a $\Cn$ space and let
$U\in\Top{X}$ be an open set. Let us consider an infinitesimal set
$I\subseteq D_{\infty}^{d}$, with $d\in\N_{>0}$ and containing the
null vector: $\underline{0}\in I$. Finally, let $f:I\freccia\ext{X}$
be a $\ECInfty$ function with $f(\underline{0})\in\ext{U}$. Then
$f(h)\in\ext{U}$ for every $h\in I$.
\end{thm}
\textbf{Proof:} From the hypothesis on $f$ it follows that $f\In{I}\ext{X}$
because $I=\ext{(I\R)}=\overline{I}$ (see Theorem \ref{thm:spacesPropertiesOfFermatFunctor}).
Hence, since $\underline{0}\in I$, by the results of Section \ref{sec: FiguresOfExtendedSpaces},
we can globally say that either $f$ is constant, and the proof is
trivial, or we can write the equality $f(h)=\ext{\gamma}(p,h)$ in
$\ext{X}$ for every $h\in I$. For the sake of clarity let $y:=f(h)$,
thus taking standard parts we get \begin{equation}
\st{y}\asymp\st{[\ext{\gamma}(p,h)]}=\gamma(p_{0},0)=\st{[\ext{\gamma}(p,0)]}\asymp\st{f(0)},\label{eq: StOf_y_and_StOff0}\end{equation}
that is $\st{y}$ and $\st{f(0)}$ are identified in $X$ (see Definition
\ref{def: xIdentifiedy} for the definition of the relation $x\asymp y$).
But $f(0)\in\ext{U}$, hence $\st{f(0)}\in U$ and $\st{y}\in U$
from \eqref{eq: StOf_y_and_StOff0}. Finally, $y=f(h)\in\ext{X}$
and hence $y=f(h)\in\ext{U}$ because of the previous Lemma \ref{lem:stdPartOf_x_inUand_x_inExt_U}.$\qedNoNewLine$

\noindent We will state both the $1$-dim Taylor's formula and the
$d$-dimensional one, because the first case can be stated in a considerably
simpler way.
\begin{thm}
\label{thm: GeneralizedDerivationFormula} Let $X$ be a $\CInfty$
space, $\alpha\in\R_{>0}$, $U\in\Top{X}$ be an open set of $X$
and\[
f:D_{\alpha}\xfrecciad{}\ext{X}\e{with}f(0)\in\ext{U}\]
\[
\phi:\ext{U}\xfrecciad{}\ext{\R^{\sf k}}\]
be $\ECInfty$ maps. Define $k_{j}\in\R$ such that\[
k_{0}:=0\]
\[
\frac{1}{k_{j}}+\frac{j}{\alpha+1}=1\quad\forall j=1,\ldots,[\alpha]=:n.\]
Then there exists one and only one $n$-tuple $(m_{1},\ldots,m_{n})$
such that
\begin{enumerate}
\item $m_{j}\in\ER_{k_{j}}^{\sf k}$ for every $j=1,\ldots,n$
\item $\phi[f(h)]=\sum\limits _{j=0}^{n}\frac{h^{j}}{j!}\cdot m_{j}\quad\forall h\in D_{\alpha}$.
\end{enumerate}
\end{thm}
The more general statement, with infinitesimal increments taken in
a product of ideals of different order, i.e. $h\in D_{\alpha_{1}}\times\ldots\times D_{\alpha_{d}}$,
is the following%
\footnote{Recall the Definition \ref{def:quotientOfn-tuples} for the definition
of the term $\frac{j}{\alpha+1}$, where $\alpha$, $j\in\N^{d}$.%
}
\begin{thm}
\label{thm:GeneralizedDerivationFormula_d-dim}Let $X$ be a $\CInfty$
space, $\alpha_{1},\ldots,\alpha_{d}\in\R_{>0}$, $U\in\Top{X}$ an
open subset of $X$ and\[
f:D_{\alpha_{1}}\times\cdots\times D_{\alpha_{d}}\xfrecciad{}\ext{X}\e{with}f(\underline{0})\in\ext{U}\]
\[
\phi:\ext{U}\xfrecciad{}\ER^{\sf k}\]
be $\ECInfty$ maps. Define $k_{j}\in\R$ such that\[
k_{\underline{0}}:=0\]
\[
\frac{1}{k_{j}}+\frac{j}{\alpha+1}=1\quad\forall j\in\N^{d}:\ 0<\frac{j}{\alpha+1}<1.\]
Then there exists one and only one\[
m:\left\{ j\in\N^{d}\,|\,\frac{j}{\alpha+1}<1\right\} \xfrecciad{}\ER^{\sf k}\]
such that
\begin{enumerate}
\item $m_{j}\in\ER_{k_{j}}^{\sf k}$ for every $j\in\N^{d}$ such that $\frac{j}{\alpha+1}<1$
\item $\phi\left[f(h)\right]=\sum\limits _{\substack{j\in\N^{d}\\
\frac{j}{\alpha+1}<1}
}\frac{h^{j}}{j!}\cdot m_{j}\quad\forall h\in D_{\alpha_{1}}\times\cdots\times D_{\alpha_{d}}$.
\end{enumerate}
\end{thm}
\noindent \textbf{Proof:} The domain of our function is\begin{align*}
D_{\alpha_{1}}\times\cdots\times D_{\alpha_{d}} & =\ext{(D_{\alpha_{1}}\R)\times\cdots\times\ext{(D_{\alpha_{d}}\R)}}=\\
 & =\overline{D_{\alpha_{1}}}\times\cdots\times\overline{D_{\alpha_{d}}}=\overline{D_{\alpha_{1}}\times\cdots\times D_{\alpha_{d}}}\end{align*}
where we have used Theorem \ref{thm:spacesPropertiesOfFermatFunctor}
for the second equality and Lemma \ref{lem: ProductInSERn} for the
latter equality. Setting $I:=D_{\alpha_{1}}\times\cdots\times D_{\alpha_{d}}$
for the sake of simplicity, we have thus\[
f:\overline{I}\freccia\ext{X}\e{and\ hence}f\In{I}\ext{X}.\]
From the results of Section \ref{sec: FiguresOfExtendedSpaces}, since
$I$ is an infinitesimal set containing $\underline{0}$, we get that
either $f$ is constant or we can write\begin{equation}
f=\ext{\gamma}(p,-)|_{I}\label{eq:secondCaseTaylorFormula}\end{equation}
for some $p\in\ext{A}$, $A$ is open in $\R^{\sf p}$, and some $\gamma\in\CInfty(A\times B,X)$
with $I\subset\ext{B}$ and $B$ open in $\R^{d}$. The case $f$
constant is trivial because it suffices to set $m_{j}:=0$ for $j\ne\underline{0}$,
$m_{\underline{0}}:=\phi\left[f(\underline{0})\right]$ to have the
existence part and to apply Corollary \ref{cor:uniquenessInGeneralTaylorFormulae}
for the uniqueness part. In the second case \eqref{eq:secondCaseTaylorFormula}
our aim is, of course, to use the composition $\phi\circ\ext{\gamma}(p,-)$,
so that now we would like to find where this composition is defined
and to prove that its domain contains the previous infinitesimal set
$I$. We have that $\eta:=\ext{\gamma}(p,-):\ext{B}\freccia\ext{X}$
in $\ECInfty$, hence, since $\ext{U}\in\Top{\ext{X}}$, we also get
that $\eta^{-1}(\ext{U})$ is open in $\ext{B}$ and hence it is also
open in $\ER^{d}$ because $B$ is open in $\R^{d}$. But we have
that $\underline{0}\in\eta^{-1}(\ext{U})$ if and only if $\eta(\underline{0})=\ext{\gamma}(p,\underline{0})=f(\underline{0})\in\ext{U}$
which is true by hypothesis. Thus, since $\eta^{-1}(\ext{U})$ is
open in $\ER^{d}$ we obtain that\[
\exists\, B_{1}\text{ open in }\R^{d}:\quad\underline{0}\in\ext{B_{1}}\subseteq\eta^{-1}(\ext{U})\subseteq\ext{B}.\]
So we are in the following situation\[
\xyC{40pt}\xymatrix{\ext{B_{1}}\ar[r]^{\eta|_{\ext{B_{1}}}} & \ext{U}\ar[r]^{\phi} & \ER^{\sf k}}
,\]
and hence the composition $\phi\circ\ext{\gamma}(p,-)|_{\ext{B_{1}}}=:\psi$
is defined in $\ext{B_{1}}$ which, being open and containing $\underline{0}$,
it also contains the infinitesimal set $I$ (Lemma \ref{lem:stdPartOf_x_inUand_x_inExt_U}).
But the idea is to use the Taylor's formula for standard smooth functions,
i.e. Theorem \ref{thm:TaylorForStndSmoothAtNonStndPoint}, and we
do not know whether the function $\psi$ is the extension of an ordinary
standard function. So, we have to note that $\ext{B_{1}}=\overline{\ext{B_{1}}}$
and hence $\psi\In{\ext{B_{1}}}\ER^{\sf k}$, so that we can apply
once again Theorem \ref{thm:figuresOfFermatSpaces} obtaining that
locally, in a neighborhood of $\underline{0}\in\ext{B_{1}}$, we can
express the figure $\psi$ as $\psi=\ext{\delta}(q,-)$ for a suitable
$\delta\in\CInfty(C\times E,\R^{\sf k})$, with $I\subseteq\ext{E}$
(the case $\psi$ constant can be dealt as seen above). Therefore
we have \[
\psi(x)=\ext{\delta}(q,x)=\phi\left[\ext{\gamma}(p,x)\right]=\phi\left[f(x)\right]\quad\forall x\in\ext{E}\cap\ext{B_{1}}=\ext{(E\cap B)}.\]
To the standard smooth function $\delta\in\CInfty(C\times E,\R^{\sf k})$
we can apply Theorem \ref{thm:TaylorForStndSmoothAtNonStndPoint}
at the non-standard point $\left(q,\underline{0}\right)\in\ext{C}\times\ext{E}=\ext{(C\times E)}$
with infinitesimal increment $(q,\underline{0}+h)$, $h\in I$; we
obtain\[
\forall h\in I:\quad\psi(h)=\left[\ext{\delta}(q,h)\right]=\phi\left[f(h)\right]=\sum\limits _{\substack{j\in\N^{d}\\
\frac{j}{\alpha+1}<1}
}\frac{h^{j}}{j!}\cdot\frac{\partial^{|j|}\delta}{\partial x^{j}}(\underline{0}).\]
Now it suffices to apply Corollary \ref{cor:uniquenessInGeneralTaylorFormulae}
to obtain the conclusion.$\qedNoNewLine$

\noindent Analogously we can state and prove a Taylor's formula for
functions $f:D_{\infty}^{d}\freccia\ext{X}$, with coefficients $m_{\underline{0}}\in\ER_{0}^{\sf k}$
and $m_{j}\in\ER_{1}^{\sf k}$.
\begin{defn}
\label{def:derivativesOfFunctionsDefinedOnInfSets}In the hypothesis
of the previous Theorem \ref{thm:GeneralizedDerivationFormula_d-dim},
we set\[
\partial\phi(f):\left\{ j\in\N^{d}\,|\,\frac{j}{\alpha+1}<1\right\} \xfrecciad{}\ER^{\sf k}\]
such that:
\begin{enumerate}
\item $\partial_{j}\phi(f)\in\ER_{k_{j}}^{\sf k}$ for every $j\in\N^{d}$
such that $\frac{j}{\alpha+1}<1$
\item $\phi\left[f(h)\right]=\sum\limits _{\substack{j\in\N^{d}\\
\frac{j}{\alpha+1}<1}
}\frac{h^{j}}{j!}\cdot\partial_{j}\phi(f)\quad\forall h\in D_{\alpha_{1}}\times\cdots\times D_{\alpha_{d}}$.
\end{enumerate}
In the case $X=\R^{\sf k}$ and $\phi=1_{\R^{\sf k}}$ we will use
the simplified notations\[
\partial_{j}f:=\partial_{j}f(\underline{0}):=\partial_{j}\phi(f)\]
\[
f^{(n)}(0):=\partial_{n}f(0)\quad\text{if}\quad f:D_{\alpha}\freccia\ER\ \text{and}\ n<\alpha+1.\]

\end{defn}
\noindent Let us note that using these notations we have that $\partial_{j}\phi(f)=\partial_{j}(\phi\circ f)$.

\noindent For example if $f:D\freccia\ER$ is smooth, then we have\[
f'(0)\in\ER_{2}\e{and}\forall h\in D:\quad f(h)=f(0)+h\cdot f'(0),\]
with $f'(0)$ uniquely determined by this property. Using this notation
we have that $f\mapsto f'(0)$ is a derivation up to second order
infinitesimals
\begin{thm}
Let $f$, $g:D\freccia\ER$ and $r\in\ER$, then
\begin{enumerate}
\item $\left(f+g\right)'(0)=f'(0)+g'(0)$
\item $\left(r\cdot f\right)'(0)=_{2}r\cdot f'(0)$ and if $r\in\R$, then
$\left(r\cdot f\right)'(0)=r\cdot f'(0)$
\item \label{enu:productAndDerivative}$\left(f\cdot g\right)'(0)=_{2}f'(0)\cdot g(0)+f(0)\cdot g'(0)$
\end{enumerate}
In other words the map $f\in\ER^{D}\mapsto\ER_{\sss{=_{2}}}$ is a
derivation (see Theorem \ref{thm:ER_k_andAlgebraicOperations} for
the definition of $\ER_{\sss{=_{k}}}$).

\end{thm}
\noindent \textbf{Proof:} We use the notations of the proof of Theorem
\ref{thm:GeneralizedDerivationFormula_d-dim} and we prove property
\emph{\ref{enu:productAndDerivative}}., the others being similar.
Thus we can write\begin{align}
f & =\ext{\gamma}(p,-)|_{D}\e{,}g=\ext{\eta}(q,-)|_{D}\nonumber \\
f'(0) & =\iota_{2}\left[\partial_{2}\gamma(p,0)\right]\e{,}g'(0)=\iota_{2}\left[\partial_{2}\eta(q,0)\right],\label{eq:firstDerivativeFormula_thmDerivation}\end{align}
where $\partial_{2}$ means partial derivative with respect to the
second slot. Therefore, recalling Theorem \ref{thm:ER_k_andAlgebraicOperations}
about the properties of $=_{k}$, we have\begin{align}
\left(f\cdot g\right)'(0) & =_{2}\partial_{2}\left(\gamma(p,-)\cdot\eta(q,-)\right)(0)=_{2}\nonumber \\
 & =_{2}\partial_{2}\gamma(p,0)\cdot\eta(q,0)+\gamma(p,0)\cdot\partial_{2}\eta(q,0).\label{eq:stepFormulaAboutProductOfDerivative}\end{align}
But from \eqref{eq:firstDerivativeFormula_thmDerivation} we have
that $f'(0)=_{2}\partial_{2}\gamma(p,0)$, $g'(0)=_{2}\partial_{2}\eta(q,0)$
and $=_{2}$ is a congruence relation with respect to ring operations
(Theorem \ref{thm:ER_k_andAlgebraicOperations}), hence from \eqref{eq:stepFormulaAboutProductOfDerivative}
we obtain the conclusion.$\qedNoNewLine$

It is important now to make some considerations about the meaning
of the derivative $\partial_{j}f(\underline{0})$, for $f:D_{n}^{d}\freccia\ER$,
with respect to the order of infinitesimals $n\in\N_{>0}$. We have
already hinted in Section \ref{sec:theFermatMethod} to the fact that
the best properties for derivatives can be proved using the Fermat
method for functions $f:V\freccia\ER$ defined in a neighborhood $V$
of the point we are interested to, e.g. $\underline{0}\in\ext{U}\subseteq V$,
with $U$ open in $\R^{d}$. But if we start from a function of the
form $f:D_{n}^{d}\freccia\ER$ defined on an infinitesimal set, then,
roughly speaking, the domain $D_{n}^{d}$ is {}``too small to give
sufficient information'' for the definition of $\partial_{j}f(\underline{0})$
using the Fermat method. Indeed, we do not have as domain a full neighborhood
to uniquely determine the smooth incremental ratio of $f$. The Taylor's
formula determines the derivatives $\partial_{j}f(\underline{0})$
in the set $\ER_{k_{j}}^{\sf k}$ and hence forces us to work up to
$k_{j}$-th order infinitesimals, i.e. using the congruence%
\footnote{We have to note that $k_{j}$, defined in the statement of Theorem
\ref{thm:GeneralizedDerivationFormula_d-dim}, really depends on the
order $n$, thus if we need to distinguish two situations with two
orders, we will use the more complete notation $k_{j}(n)$.%
} $=_{k_{j}}$. As a further proof of these informal considerations,
it seems plausible to expect that the larger is the order $n\in\N_{>0}$
the larger is the {}``information'' we have at disposal. More precisely,
the situation we want to analyze is the following: if we take a smooth
function $f:D_{m}^{d}\freccia\ext{X}$ and $n<m$, then what is the
relationship between $\partial_{j}\phi(f)$ and $\partial_{j}\phi\left(f|_{D_{n}^{d}}\right)$?
The answer is: it results $k_{j}(n)>k_{j}(m)$ and $\partial_{j}\phi(f)$
is equal to $\partial_{j}\phi\left(f|_{D_{n}^{d}}\right)$ up to infinitesimals
of order $k_{j}(n)$.
\begin{thm}
Let $X\in\CInfty$ be a smooth space. Let us consider $n$, $m$,
$d\in\N_{>0}\cup\{+\infty\}$, $j\in\N^{d}$ with $n<m$ and $1\le|j|\le n$.
Moreover, let us consider an open set $U\in\Top{X}$ and smooth maps
of the form\begin{align*}
f:D_{m}^{d}\freccia\ext{X}\e{with}f(\underline{0}) & \in\ext{U}\end{align*}
\[
\phi:\ext{U}\freccia\ER^{\sf k}.\]
Finally, let\[
\frac{1}{k_{j}(p)}+\frac{j}{p+1}=1\quad\forall p\in\N_{>0}\cup\{+\infty\}.\]
Then we have
\begin{enumerate}
\item \label{enu:1ThmAboutDerivativesOnInfSetsAndRestriction}$k_{j}(n)>k_{j}(m)$
\item \label{enu:2ThmAboutDerivativesOnInfSetsAndRestriction}$\partial_{j}\phi(f)=_{k_{j}(n)}\partial_{j}\phi\left(f|_{D_{n}^{d}}\right)$
\item \label{enu:3ThmAboutDerivativesOnInfSetsAndRestriction}$k_{j}(m)<\omega\left[\partial_{j}\phi(f)-\partial_{j}\phi\left(f|_{D_{n}^{d}}\right)\right]\le k_{j}(n)$
and hence\[
\partial_{j}\phi(f)-\partial_{j}\phi\left(f|_{D_{n}^{d}}\right)\in D_{k_{j}(n)}\]

\end{enumerate}
\end{thm}
\noindent \textbf{Proof:} To prove \emph{\ref{enu:1ThmAboutDerivativesOnInfSetsAndRestriction}}.
it suffices to note that $k_{j}(p)=\frac{p+1}{p+1-|j|}$, because
$\frac{j}{p+1}=\frac{|j|}{p+1}$, and that for $a>b$ the real function
$x\mapsto\frac{x+a}{x+b}$ has a derivative $\frac{b-a}{(x+b)^{2}}$
which is negative for every $x$.

\noindent To prove \emph{\ref{enu:2ThmAboutDerivativesOnInfSetsAndRestriction}}.
from Definition \ref{def:derivativesOfFunctionsDefinedOnInfSets}
we have\[
\forall h\in D_{m}^{d}:\quad\phi\left[f(h)\right]=\sum\limits _{|j|\le m}\frac{h^{j}}{j!}\cdot\partial_{j}\phi(f).\]
But $n<m$ so $D_{n}^{d}\subseteq D_{m}^{d}$ and thus from Corollary
\ref{cor:suffCondToHaveProductZeroMulti-indexes} we have\[
\forall h\in D_{n}^{d}:\quad\phi\left[f(h)\right]=\sum\limits _{|j|\le n}\frac{h^{j}}{j!}\cdot\partial_{j}\phi(f).\]
Now using Theorem \ref{thm:generalCancellationLaw} we obtain $h^{j}\cdot\partial_{j}\phi(f)=h^{j}\cdot\iota_{k_{j}(n)}\left[\partial_{j}\phi(f)\right]$
for every $h\in D_{n}^{d}$ and substituting we get\[
\forall h\in D_{n}^{d}:\quad\phi\left[f(h)\right]=\sum\limits _{|j|\le n}\frac{h^{j}}{j!}\cdot\iota_{k_{j}(n)}\left[\partial_{j}\phi(f)\right],\]
hence from the uniqueness in Taylor's formula we obtain $\iota_{k_{j}(n)}\left[\partial_{j}\phi(f)\right]=\partial_{j}\phi\left(f|_{D_{n}^{d}}\right)=\iota_{k_{j}(n)}\left[\partial_{j}\phi\left(f|_{D_{n}^{d}}\right)\right]$,
because $\iota_{k}(x)=x$ if $x\in\ER_{k}$, and this proves \emph{\ref{enu:2ThmAboutDerivativesOnInfSetsAndRestriction}}.

Property \emph{\ref{enu:3ThmAboutDerivativesOnInfSetsAndRestriction}}.
follows directly from Theorem \ref{thm:orderAndIota_k}.$\qedNoNewLine$

Much in the same way as the Fermat method provides a very useful instrument
to derive the calculus for functions defined on open sets, the previous
theorems show us that the derivatives $\partial_{j}\phi(f)$ provides
a useful instrument to study functions defined on infinitesimal sets.
The following theorem states that equality of derivatives through
observables implies identity of the functions:
\begin{thm}
\label{thm:equalityOfDerivativesImpliesEqualityEverywhere}Let $X\in\CInfty$
be a smooth space and $n$, $d\in\N_{>0}$. Let us consider two smooth
functions\[
f,\ g:D_{n}^{d}\freccia\ext{X}\e{with}f(\underline{0})=g(\underline{0}).\]
Moreover, let us assume that the derivatives of these functions are
equal, i.e.\[
\partial_{j}\ext{\phi}(f)=\partial_{j}\ext{\phi}(g)\quad\forall j\in\N^{d}:\ 1\le|j|\le n\]
for every observable $\phi:U\freccia\R^{\sf k}$ of the space $X$
with $f(\underline{0})\in\ext{U}$. Then\[
f=g\]

\end{thm}
\noindent Let us note that this theorem, which is a consequence of
our definition of equality in $\ext{X}$ using observables (see Definition
\ref{def:equalityInExt_X}), is not trivial, because in our context
we do not have charts on our spaces $X\in\CInfty$.

\noindent \textbf{Proof:} Take $h\in D_{n}^{d}$, we have to prove
that $y:=f(h)$ and $z:=g(h)$ are equal in $\ext{X}$. Using typical
notations and neglecting some details, we can say that $y=f(h)=\ext{\gamma}(p,h)$
in $\ext{X}$ so that\begin{align*}
\st{y} & =\st{\left[\ext{\gamma}(p,h)\right]}=\gamma(\st{p},\underline{0})=\\
 & =\st{\left[\ext{\gamma(p,\underline{0})}\right]}=\st{f}(\underline{0}).\end{align*}
But $f(\underline{0})=g(\underline{0})$ in $\ext{X}$ by hypothesis,
so $\st{y}=\st{f}(\underline{0})\asymp\st{g}(\underline{0})=\st{z}$.
Now let us take an observable $\phi:U\freccia\R^{\sf k}$ of $X$,
we have to prove that\begin{equation}
y_{0}\in U\ \iff\ z_{0}\in U\label{eq:firstThingToProve}\end{equation}
\begin{equation}
y_{0}\in U\ \Longrightarrow\ \phi(y_{t})=\phi(z_{t})+o(t).\label{eq:secondThingToProve}\end{equation}
The first one follows directly from the identification $\st{y}\asymp\st{z}$.
For the second one, if $y_{0}\in U$, then $y_{0}=\st{y}=\st{f}(\underline{0})$,
thus $\st{f}(\underline{0})\in U$ and hence $f(\underline{0})\in\ext{U}$
from Lemma \ref{lem:stdPartOf_x_inUand_x_inExt_U}. We can thus apply
our hypotheses to obtain the equality of derivatives $\partial_{j}\ext{\phi}(f)=\partial_{j}\ext{\phi}(g)$
for every $j\in\N^{d}$ with $1\le|j|\le n$. Using Taylor's formula\[
\ext{\phi}(y)=\ext{\phi}\left[f(h)\right]=\sum\limits _{|j|\le n}\frac{h^{j}}{j!}\cdot\partial_{j}\ext{\phi}(f)=\sum\limits _{|j|\le n}\frac{h^{j}}{j!}\cdot\partial_{j}\ext{\phi}(g)=\ext{\phi}\left[g(h)\right]=\ext{\phi}(z),\]
this equality being in $\ER^{\sf k}$, i.e. $\phi(y_{t})=\phi(z_{t})+o(t)$,
which is the conclusion stated in the theorem.$\qedNoNewLine$

There is the possibility to connect the methods developed for the
differential calculus of function defined on open sets (see the previous
Section \ref{cha:CalculusOnOpenDomains}) with the differential calculus
of smooth functions defined on infinitesimal sets. Indeed, the following
results prove that functions of the form $f:S\freccia\ER^{n}$ can
be seen locally as {}``infinitesimal polynomials with smooth coefficients''.
\begin{thm}
\label{thm:infPolynomialsWithSmoothCoefficients}Let $S\subseteq\ER^{\sf s}$
and $f:S\freccia\ER^{n}$ a map (in $\Set$). Then it results that\begin{equation}
f:S\freccia\ER^{n}\text{ is smooth in }\ECInfty\label{eq:HypoInfPolyWithSmoothCoeff}\end{equation}
 if and only if for every $x\in S$ we can write\begin{equation}
f(y)=\sum_{\substack{|q|\le k\\
q\in\N^{d}}
}a_{q}(y)\cdot p^{q}\quad\forall y\in\ext{V}\cap S,\label{eq:thesisInfPolyWithSmoothCoeff}\end{equation}
for suitable:
\begin{enumerate}
\item $d$, $k\in\N$
\item $p\in D_{k}^{d}$
\item $V$ open subset of $\R^{\sf s}$ such that $x\in\ext{V}$
\item $(a_{q})_{\substack{|q|\le k\\
q\in\N^{d}}
}$ family of $\CInfty(V,\R^{n})$.
\end{enumerate}
\end{thm}
In other words, every smooth function $f:S\freccia\ER^{n}$ can be
constructed locally starting from some {}``infinitesimal parameters''\[
p_{1},\ldots,p_{d}\in D_{k}\]
and from ordinary smooth functions\[
a_{q}\in\CInfty(V,\R^{n})\]
and using polynomial operation only with $p_{1}$, ..., $p_{d}$ and
with coefficients $a_{q}(-)$. Roughly speaking, we can say that they
are {}``infinitesimal polynomials with smooth coefficients. The polynomials
variables act as parameters only''. By the sheaf property, here {}``locally''
means that this construction using infinitesimal polynomials has to
be done in a neighborhood of each point $x\in S$, but in such a way
to have equal polynomials on intersecting neighborhoods.

If $f:I\freccia\ER^{n}$ with $\underline{0}\in I\subseteq D_{\infty}^{\sf s}$,
then for $x=\underline{0}$ we can write \eqref{eq:thesisInfPolyWithSmoothCoeff}
globally as\begin{equation}
f(h)=\sum_{\substack{|q|\le k\\
q\in\N^{d}}
}a_{q}(h)\cdot p^{q}\quad\forall h\in I,\label{eq:globalInfPolyWithSmoothCoeffs}\end{equation}
because $I\subseteq\ext{V}$. From \eqref{eq:globalInfPolyWithSmoothCoeffs}
we obtain\[
\partial_{j}f(h_{1})=\iota_{k_{j}}\left(\sum_{\substack{|q|\le k\\
q\in\N^{d}}
}\partial_{j}a_{q}(h_{1})\cdot p^{q}\quad\forall h_{1}\in I\right),\]
where $j\in\N^{d}$, $k_{j}\in\R$ and $\alpha_{1}$, ..., $\alpha_{d}\in\N_{>0}$
are such that $0<\frac{j}{\alpha+1}<1$ and $\frac{1}{k_{j}}+\frac{j}{\alpha+1}=1$.
Therefore from Theorem \ref{thm:generalCancellationLaw} we get\begin{equation}
\forall h\in D_{\alpha_{1}}\times\cdots\times D_{\alpha_{d}}\pti h\cdot\partial_{j}f(h_{1})=h\cdot\sum_{\substack{|q|\le k\\
q\in\N^{d}}
}\partial_{j}a_{q}(h_{1})\cdot p^{q}.\label{eq:derivativeWithout_iota_k}\end{equation}
All this permits to use the results about the differential calculus
of functions like $a_{q}\in\CInfty(V,\ER^{n})$, defined on open sets,
to functions defined on infinitesimal sets. Moreover, equalities of
the form \eqref{eq:derivativeWithout_iota_k} permit to avoid the
use of the map $\iota_{k}:\ER^{n}\freccia\ER_{\sss{=_{k}}}^{n}$.

\noindent \textbf{Proof:} The implication \eqref{eq:thesisInfPolyWithSmoothCoeff}$\Rightarrow$
\eqref{eq:HypoInfPolyWithSmoothCoeff} follows directly from Theorem
\ref{thm:figuresOfFermatSpaces}. For the opposite implication, let
us write $f|_{\mathcal{V}}=\ext{\alpha}(\pi,-)|_{\mathcal{V}}$, as
usual, in a neighborhood of $x\in\ext{V}\cap S$, for $\alpha\in\CInfty(U\times V,\R^{n})$
and where $\pi\in\ext{U}\subseteq\ER^{d}$ works as the usual non
standard parameter. Set $r:=\st{\pi}$ and $p:=\pi-r$ so that $p\in D_{k}^{d}$
for some $k\in\N_{>0}$. Using the infinitesimal Taylor's formula
we get\[
f(y)=\alpha(\pi,y)=\alpha(r+p,y)=\sum_{|q|\le k}\partial_{q}\alpha(r,y)\cdot\frac{p^{q}}{q!},\]
from which we have the conclusion setting $a_{q}:=\frac{1}{q!}\cdot\partial_{q}\alpha(r,-)$.$\qedNoNewLine$

\noindent Taking an enumeration of all these multi-indexes $q\in\N^{d}$,
i.e.\[
\{q_{1},\ldots,q_{N}\}=\left\{ q\in\N^{d}\,:\,|q|\le k\ ,\ p^{q}\ne0\right\} \setminus\{\underline{0}\}\]
\[
q_{i}\ne q_{j}\e{if}i\ne j,\]
then we can write the infinitesimal polynomial \eqref{eq:thesisInfPolyWithSmoothCoeff}
in a simpler way, even if it hide the powers $p^{q}$ of the infinitesimal
parameter $p\in D_{k}^{d}$. In fact, using this enumeration we can
write\[
f(y)=\sum_{i=1}^{N}a_{q_{i}}(y)\cdot p_{1}^{q_{i1}}\cdot\ldots\cdot p_{d}^{q_{id}}+a_{\underline{0}}(y).\]
It suffices to set $\pi_{i}:=p_{1}^{q_{i1}}\cdot\ldots\cdot p_{d}^{q_{id}}$,
$b_{i}:=a_{q_{i}}$ and $b_{0}:=a_{\underline{0}}$ to have\begin{equation}
f(y)=b_{0}(y)+\sum_{i=1}^{N}b_{i}(y)\cdot\pi_{i}\quad\forall y\in\ext{V}\cap S.\label{eq:infPolyWithSmoothCoeffsSimpleForm}\end{equation}

\noindent As usual, both the smooth functions $b_{i}$ and the infinitesimal
parameters $\pi_{i}$ are not uniquely determined by formulas of the
form \eqref{eq:infPolyWithSmoothCoeffsSimpleForm}.

\section{Smoothness of derivatives}

In our smooth context, it is important that the definition of derivative
for our non standard smooth functions always produces a smooth operator.
On the other hand, it is natural to expect, exactly as for the standard
part map (see Corollary \ref{cor:TheStandardPartIsNotSmooth}) that
every function $\iota_{k}:\ER\freccia\ER$ mapping a Fermat real $x\in\ER$
to {}``$x$ up to $k$-th order infinitesimals'', i.e. $\iota_{k}(x)$,
cannot be smooth. If this is so, then also the first derivative cannot
be smooth being thought as a function of $f$, i.e. $f\in\ER^{D}\mapsto f'(0)\in\ER$.
Let us consider, for example, the following function defined on the
infinitesimal set $D_{2}$ of second order infinitesimals and depending
on the parameter $p\in\ER$:\[
f_{p}(h):=\frac{1}{2}(h+p)^{2}\quad\forall h\in D_{2}.\]
We have $f_{p}(h)=\frac{1}{2}p^{2}+hp+\frac{h^{2}}{2}$, so that from
the Taylor's formula we have\begin{align*}
f_{p}'(0) & =\iota_{\frac{3}{2}}(p)\e{in\ fact:}\frac{n+1}{n+1-j}=\frac{2+1}{2+1-1}=\frac{3}{2}\\
f_{p}''(0) & =\iota_{3}(1)=1\e{in\ fact:}\frac{n+1}{n+1-j}=\frac{2+1}{2+1-2}=3.\end{align*}
So, if the map\[
f\in\ER^{D_{2}}\mapsto f'(0)\in\ER\]
is smooth, also the map\[
p\in\ER\mapsto f_{p}\in\ER^{D_{2}}\mapsto f'_{p}(0)=\iota_{\frac{3}{2}}(p)\in\ER\]
would be smooth. Therefore, the smoothness of the maps $\iota_{k}:\ER\freccia\ER$
is strictly tied with the smoothness of the derivatives. We will see
that these maps are not smooth, as a simple consequence of the following
general result.
\begin{thm}
\label{thm:smoothFnctsThatTakeStndToStnd}Let $M$, $N$ be manifolds,
and $f:\ext{N}\freccia\ext{M}$ be a $\ECInfty$ function. Then\[
f(N)\subseteq M\then f=\ext{\left(f|_{N}\right)}.\]

\end{thm}
\noindent In other words, a $\ECInfty$ function between extended
manifolds that takes standard points to standard points, can be realized
as the extension of an ordinary smooth function.

\noindent \textbf{Proof:} Let us consider a generic point $n_{1}\in\ext{N}$.
We want to prove that $f(n_{1})=\ext{\left(f|_{N}\right)(n_{1})}$.
Let $n:=\st{n_{1}}\in N$ so $f(n)\in M$ and we can consider a chart
$(U,\phi)$ on $n\in N$ and another one $(V,\psi)$ on $f(n)$. We
can assume, for simplicity, $\phi(U)=\R^{\sf n}$, $\psi(V)=\R^{\sf m}$,
$\phi(n)=\underline{0}$ and the open set $U$ sufficiently small
so that $f(\ext{U})\subseteq\ext{V}$. Diagrammatically, in $\ECInfty$,
the situation is the following\[
\xyC{40pt}\xymatrix{\ER^{\sf n}\ar[r]^{\ext{\phi}^{-1}} & \ext{U}\ar[r]^{f|_{\ext{U}}} & \ext{V}\ar[r]^{\ext{\psi}} & \ER^{\sf m}}
.\]
To the smooth function $\ext{\phi}^{-1}\cdot f|_{\ext{U}}\cdot\ext{\psi}:\ER^{\sf n}\freccia\ER^{\sf m}$
we can apply Theorem \ref{thm:figuresOfFermatSpaces} obtaining that
in a neighborhood of $\underline{0}$ we can write \[
\psi\left[f\left(\phi^{-1}y\right)\right]=\ext{\gamma}(p,y),\]
where $\gamma\in\CInfty(A\times B,\R^{\sf m})$, $p\in\ext{A}$, $A$
is open in $\R^{\sf p}$ and $B$ is an open neighborhood of $\underline{0}$
in $\R^{\sf n}$. Setting $r:=\st{p}\in\R^{\sf p}$ and $h:=p-r\in D_{k}^{\sf p}$
for some $k\in\N_{>0}$, we get\[
\psi\left[f\left(\phi^{-1}y\right)\right]=\ext{\gamma}(r+h,y)=\sum\limits _{\substack{q\in\N^{\sf p}\\
|q|\le k}
}\partial_{1}^{q}\gamma(r,y)\cdot\frac{h^{q}}{q!},\]
where $\partial_{1}$ means the derivative with respect to the first
slot of $\gamma(-,-)$. Setting for simplicity $a_{i}(y):=\frac{1}{q_{i}!}\cdot\partial_{1}^{q_{i}}\gamma(r,y)$,
$h_{i}:=h^{q_{i}}$, where $\{q_{1},\ldots,q_{N}\}=\{q\in\N^{\sf p}\,:\,|q|\le k\}$,
$q_{i}\ne q_{j}$ if $i\ne j$, we can write\begin{equation}
\psi\left[f\left(\phi^{-1}y\right)\right]=a_{0}(y)+\sum\limits _{i=1}^{N}a_{i}(y)\cdot h_{i}\quad\forall y\in\ext{B},\label{eq:psi_f_phi_to_minus1_AsLinearCombinationOfSmoothFncs}\end{equation}
where the functions $a_{i}\in\CInfty(B,\R^{\sf m})$ are standard
smooth maps. We can suppose $(a_{1},\ldots,a_{N})$ linearly independent
in the real vector space $\CInfty(B,\R^{\sf m})$, because otherwise
we can select among them a basis and express the other functions as
a linear combination of the basis. Now, let us evaluate the standard
part of \eqref{eq:psi_f_phi_to_minus1_AsLinearCombinationOfSmoothFncs}
at a generic standard point $r\in B$:\begin{equation}
\forall r\in B:\quad\st{\left\{ \psi\left[f\left(\phi^{-1}r\right)\right]\right\} }=\st{a_{0}}(r)+\sum\limits _{i=1}^{N}\st{a_{i}}(r)\cdot\st{h_{i}}=\st{a_{0}}(r),\label{eq:psi_f_phi_to_minus1_atStndPoints}\end{equation}
because $\st{h_{i}}=0$ since $h_{i}=h^{q_{i}}=(p-\st{p})^{q_{i}}\in D_{\infty}$.
But $r\in B\subseteq\R^{\sf n}$, so $\phi^{-1}(r)\in U\subseteq\ext{U}$
and $f\left(\phi^{-1}r\right)\in M\cap\ext{V}=V$ because $f(N)\subseteq M$.
Thus, $\psi\left[f\left(\phi^{-1}r\right)\right]\in\R^{\sf m}$ and
hence from \eqref{eq:psi_f_phi_to_minus1_atStndPoints} we get \begin{equation}
\psi\left[f\left(\phi^{-1}r\right)\right]=a_{0}(r)\quad\forall r\in B.\label{eq:psi_f_phi_at_r}\end{equation}
Therefore\begin{equation}
\sum\limits _{i=1}^{N}a_{i}(r)\cdot h_{i}=\underline{0}\quad\forall r\in B.\label{eq:systemWith_a_i_And_h_i}\end{equation}
The functions $a_{i}$ are continuous and linearly independent, so
that from Lemma \ref{lem: IndependentFunctionsAndDeterminant} we
can find $r_{1},\ldots,r_{N}\in B$ such that\begin{align*}
\det\left[\begin{array}{ccc}
a_{1}(r_{1}) & \ldots & a_{N}(r_{1})\\
\vdots\\
a_{1}(r_{N}) & \ldots & a_{N}(r_{N})\end{array}\right]\ne0 & .\end{align*}
We can write \eqref{eq:systemWith_a_i_And_h_i} as\[
\left[\begin{array}{ccc}
a_{1}(r_{1}) & \ldots & a_{N}(r_{1})\\
\vdots\\
a_{1}(r_{N}) & \ldots & a_{N}(r_{N})\end{array}\right]\cdot\left[\begin{array}{c}
h_{1}\\
\vdots\\
h_{N}\end{array}\right]=\underline{0},\]
and therefore we obtain that $h_{i}=0$ for every $i=1,\ldots,N$,
and\[
\psi\left[f\left(\phi^{-1}y\right)\right]=a_{0}(y)\]
for every $y\in\ext{B}$ from \eqref{eq:psi_f_phi_to_minus1_AsLinearCombinationOfSmoothFncs},
which has to be understood as an abbreviation of\begin{equation}
\ext{\psi}\left[f\left(\ext{\phi}^{-1}y\right)\right]=\ext{a_{0}}(y)\quad\forall y\in\ext{B}.\label{eq:psi_f_phi_isStnd}\end{equation}
Now $n_{1}\asymp n$ because $\st{n_{1}}=n$, so that $\ext{\phi}(n_{1})\asymp\phi(n)=\underline{0}\in\ext{B}$
and thus also $\ext{\phi}(n_{1})\in\ext{B}$ from the definition of
$\asymp$. We can thus apply \eqref{eq:psi_f_phi_isStnd} with $y=\ext{\phi}(n_{1})$
obtaining that $\ext{\psi}\left[f(n_{1})\right]=\ext{a_{0}}\left(\ext{\phi}n_{1}\right)$,
and hence\begin{equation}
f(n_{1})=\ext{\psi}^{-1}\left[\ext{a_{0}}\left(\ext{\phi}n_{1}\right)\right]=\ext{\left(\psi^{-1}\circ a_{0}\circ\phi\right)|_{\phi^{-1}(B)}}(n_{1}).\label{eq:f(n_1)}\end{equation}
Finally, we must prove that $\left(\psi^{-1}\circ a_{0}\circ\phi\right)(x)=\left(f|_{N}\right)(x)$
in an open neighborhood of $n$, but from \eqref{eq:psi_f_phi_at_r}
and taking a generic $x\in\phi^{-1}(B)\subseteq N$ we get $ $\begin{align*}
\psi\left[f(x)\right] & =a_{0}\left[\phi(x)\right]\\
f(x) & =\left(f|_{N}\right)(x)=\psi^{-1}\left[a_{0}\left(\phi x\right)\right]=\left(\psi^{-1}\circ a_{0}\circ\phi\right)|_{\phi^{-1}(B)}(x),\end{align*}
 and therefore $f(n_{1})=\ext{\left(f|_{N}\right)}(n_{1})$ from \eqref{eq:f(n_1)}.$\qedNoNewLine$

\noindent From this general result it follows
\begin{cor}
\label{cor:i_k_IsNotSmooth}Let $k\in\R_{\ge1}$, then the function\begin{equation}
\iota_{k}:\ER\freccia\ER\label{eq:i_k_IsNotSmoothMapWithWrongCod}\end{equation}
is not smooth.
\end{cor}
\textbf{Proof:} In fact $\iota_{k}(r)=r$ for every $r\in\R$, so
from the previous Theorem \ref{thm:smoothFnctsThatTakeStndToStnd}
we have\[
\iota_{k}=\ext{\left(\iota_{k}|_{\R}\right)}=\ext{1_{\R}}=1_{\ER}\]
if $\iota_{k}$ is smooth. But this is impossible because, e.g. $\iota_{k}(\diff{t_{k}})=0\ne\diff{t_{k}}=1_{\ER}(\diff{t_{k}})$.$\qedNoNewLine$

This negative result will be counteracted in two ways: in the first
one we will prove that, in spite of this corollary, any map of the
form\[
(f,\phi,h)\mapsto h^{j}\cdot\partial_{j}\phi(f)\]
is smooth; the second one says that the negative result is due to
the choice of a wrong codomain in \eqref{eq:i_k_IsNotSmoothMapWithWrongCod}.
\begin{thm}
\label{thm:nDerInTermsOfFnctValues1-dim}Let $X\in\CInfty$ be a smooth
space, $n\in\N_{>0}$ and $U\in\Top{X}$ an open set of $X$. Let
us consider the $\ECInfty$-maps\[
f:D_{n}\freccia\ext{X}\e{with}f(0)\in\ext{U}\]
\[
\phi:\ext{U}\freccia\ER^{\sf k}.\]
Finally, let\[
A=\left[\begin{array}{cccc}
1 & 1 & \ldots & 1\\
2 & 2^{2} & \ldots & 2^{n}\\
\vdots\\
n & n^{2} & \ldots & n^{n}\end{array}\right]^{-1}\]
be the inverse of the submatrix $V(1;1)$ obtained deleting the first
row and the first column of the Vandermonde matrix $V$ determined
by $(0,1,2,\ldots,n)$.

\noindent Then for every $j=1,\ldots,n$ we have\begin{equation}
h^{j}\cdot\partial_{j}\phi(f)=\sum_{i=1}^{n}j!\cdot a_{ij}\cdot\left\{ \phi\left[f(i\cdot h)\right]-\phi\left[f(0)\right]\right\} \quad\forall h\in D_{n},\label{eq:smoothDiffFinalFormula}\end{equation}
hence the function\begin{equation}
(f,\phi,h)\in\ext{U}^{D_{n}}\times\left(\ER^{\sf k}\right)^{\ext{U}}\times D_{n}\mapsto h^{j}\cdot\partial_{j}\phi(f)\in\ER^{\sf k}\label{eq:smoothDiffFinalFormula2}\end{equation}
is smooth in $\ECInfty$.
\end{thm}
\textbf{Proof:} From the infinitesimal Taylor's formula, Theorem \ref{thm: GeneralizedDerivationFormula},
we have\begin{equation}
\phi\left[f(h)\right]=\phi\left[f(0)\right]+\sum_{j=1}^{n}\frac{h^{j}}{j!}\cdot\partial_{j}\phi(f)\quad\forall h\in D_{n}.\label{eq:TaylorForSmoothDiff}\end{equation}
Let $x_{j}:=\frac{h^{j}}{j!}\cdot\partial_{j}\phi(f)$ for $j=1,\ldots,n$
and for a fixed $h\in D_{n}$ let $h$, $2h$, $3h$, ..., $n\cdot h$
in \eqref{eq:TaylorForSmoothDiff}. We obtain\begin{align*}
\phi(fh)-\phi(f0) & =x_{1}+\ldots+x_{n}\\
\phi\left[f(2h)\right]-\phi(f0) & =2x_{1}+2^{2}x_{2}+\ldots+2^{n}x_{n}\\
 & \ldots\\
\phi\left[f(nh)\right]-\phi(f0) & =nx+n^{2}x_{2}+\ldots+n^{n}x_{n}.\end{align*}
So that we can write this system of equations in $x_{1},\ldots,x_{n}$
as\[
\left[\begin{array}{cccc}
1 & 1 & \ldots & 1\\
2 & 2^{2} & \ldots & 2^{n}\\
\vdots\\
n & n^{2} & \ldots & n^{n}\end{array}\right]\cdot\left[\begin{array}{c}
x_{1}\\
\vdots\\
x_{n}\end{array}\right]=\left[\begin{array}{c}
\phi(fh)-\phi(f0)\\
\vdots\\
\phi\left[f(nh)\right]-\phi(f0)\end{array}\right],\]
from which the first part \eqref{eq:smoothDiffFinalFormula} of the
conclusion follows.

The second part follows noting that the right hand side of \eqref{eq:smoothDiffFinalFormula}
gives the function \eqref{eq:smoothDiffFinalFormula2} as a composition
of smooth functions (among which we have to consider some evaluations,
that in a cartesian closed category are always smooth; see Section
\ref{sec:ExamplesOfCnSpacesAndFunctions}).$\qedNoNewLine$

\noindent Analogously, we can prove a corresponding result in the
$d$-dimensional case:
\begin{thm}
Let $X\in\CInfty$ be a smooth space, $n$, $d\in\N_{>0}$ and $U\in\Top{X}$
an open set. Let us consider the $\ECInfty$-maps\[
f:D_{n}^{d}\freccia\ext{X}\e{with}f(\underline{0})\in\ext{U}\]
\[
\phi:\ext{U}\freccia\ER^{\sf k}.\]
Then for every $j\in\N^{d}$ such that $|j|\le n$, the map\[
(f,\phi,h)\in\ext{U}^{D_{n}^{d}}\times\left(\ER^{\sf k}\right)^{\ext{U}}\times D_{n}^{d}\mapsto h^{j}\cdot\partial_{j}\phi(f)\in\ER^{\sf k}\]
is a $\ECInfty$ map.
\end{thm}
The second solution of the negative result of Corollary \ref{cor:i_k_IsNotSmooth}
is to admit that the codomain of the map $\iota_{k}$ is not correct,
but we have to change it as follows\[
\iota_{k}:\ER\freccia\ER_{\sss{=_{k}}}.\]
Indeed, we have
\begin{thm}
\label{thm:i_k_IsSmooth}Let $k\in\R_{\ge1}$, then the map \[
\iota_{k}:\ER\freccia\ER_{\sss{=_{k}}}\]
is a $\ECInfty$ map. We recall here that $\ER_{\sss{=_{k}}}$ is
the quotient set $\ER/=_{k}$.
\end{thm}
\noindent \textbf{Proof:} For clarity, we will use the notations with
the equivalence classes, so that here the map $\iota_{k}$ as to be
understood as defined by\begin{equation}
{\displaystyle \iota_{k}(m):=\left[\st{m}+\sum_{\substack{i=1\\
\omega_{i}(m)>k}
}^{N}\st{m_{i}}\cdot\diff{t_{\omega_{i}(m)}}\right]_{=_{k}}\quad\forall m\in\ER}.\label{eq:iota_k-definitionWithEqualityUpTo_k}\end{equation}
But using the notations with the equivalence classes, the equality\[
\st{u}+\sum_{\substack{i=1\\
\omega_{i}(u)>k}
}^{N}\st{u_{i}}\cdot\diff{t_{\omega_{i}(u)}}=_{k}u\quad\forall u\in\ER\]
now gives\[
\iota_{k}(u)=[u]_{\sss{=_{k}}}\quad\forall u\in\ER,\]
and hence the map $\iota_{k}$ defined in \eqref{eq:iota_k-definitionWithEqualityUpTo_k}
is simply the projection onto the quotient set $\ER_{\sss{=_{k}}}$
which is always smooth by the co-completeness of the category $\ECInfty$
(see Theorem \ref{thm: limit in Cn}).$\qedNoNewLine$

To clarify further the relationships between $\ER_{k}$ and $\ER_{\sss{=_{k}}}$
we also prove the following
\begin{thm}
Let $k\in\R_{\ge1}$, then the map\[
i_{k}:x\in\ER_{\sss{=_{k}}}\mapsto\iota_{k}(x)\in\ER_{k}\]
is not smooth in $\ECInfty$.
\end{thm}
\noindent Even if we have these negative results, the map $\iota_{k}:\ER\freccia\ER_{k}$
as defined in the Definition \ref{def:upTo_k-thOrderInfinitesimals}
has not to be forgotten: if we only need algebraic properties like
those expressed in results like those of Chapter \ref{cha:equalityUpTo_k-thOrder},
then we do not need the whole map $\iota_{k}$ but only terms of the
form $\iota_{k}(m)$, and in this case, we can think $\iota_{k}(m)\in\ER_{k}$.
If, instead, we need to prove smoothness of derivatives, then we need
the map $\iota_{k}$ thought with codomain: $\iota_{k}:\ER\freccia\ER_{\sss{=_{k}}}$.
This little bit of notational confusion disappears completely once
we specify domains and codomains of the map $\iota_{k}$ we are considering.

\noindent \textbf{Proof:} Firstly, the map $i_{k}:\ER_{\sss{=_{k}}}\freccia\ER_{k}$
is well defined because the definition of $x=_{k}y$ is exactly $\iota_{k}(x)=\iota_{k}(y)$
(see Definition \ref{def:upTo_k-thOrderInfinitesimals}).

Now, let us take a figure $\delta\In{H}\ER_{\sss{=_{k}}}$ on the
quotient set $\ER_{\sss{=_{k}}}$. By Theorem \ref{thm: limit in Cn}
this means that for every $h\in H$ we can find a neighborhood $U$
of $h$ in $H$ and a figure $\alpha\In{U}\ER$ such that $\delta|_{U}=\alpha\cdot[-]_{=_{k}}$.
Let us consider the composition $\delta\cdot i_{k}:H\freccia\ER_{k}$
in the neighborhood $U$:\[
i_{k}\left(\delta(u)\right)=i_{k}\left\{ \left[\alpha(u)\right]_{=_{k}}\right\} =\iota_{k}\left[\alpha(u)\right].\]
So, taking $\alpha=1_{\ER}$, from this we would obtain that if $i_{k}$
were smooth, then also $\iota_{k}(u)$ would be smooth in $u$. In
other words, it would be smooth if considered as a map from $\ER$
to $\ER_{k}\hookrightarrow\ER$, but we already know that this does
not hold from Corollary \ref{cor:i_k_IsNotSmooth}.$\qedNoNewLine$

Now we have to understand with respect to what variables we have to
mean that {}``derivatives are smooth functions'', because we are
considering functions defined on infinitesimal sets. The natural answer
is given by the following
\begin{defn}
\label{def:derivativesOfFncsDefinedOnInfSets_x}Let $X\in\CInfty$
be a smooth space, $d\in\N_{>0}$, $n\in\N_{>0}\cup\{+\infty\}$,
and $U\in\Top{X}$ an open set. Let us consider the $\ECInfty$ maps\[
\phi:\ext{U}\freccia\ER^{\sf k}\]
\[
f:V\freccia\ext{X}\e{with}x\in V\subseteq\ER^{d}\e{and}f(x)\in\ext{U}.\]
Moreover, let us suppose that $V$ verifies\[
\forall h\in D_{n}^{d}\pti x+h\in V,\]
so that we can define $f_{x}:h\in D_{n}^{d}\mapsto f(x+h)\in\ext{X}$.
Then for every multi-index $j\in\N^{d}$ with $1\le|j|\le n$ we define\[
\partial_{j}\phi(f)_{x}:=\partial_{j}\phi\left(f_{x}\right).\]
As usual, if $X=\R^{\sf k}$ and $\phi=1_{\R^{\sf k}}$, we will use
the simplified notations $\partial_{j}f_{x}:=\partial_{j}f(x):=\partial_{j}\phi(f)_{x}$
and $f^{(j)}(x):=\partial_{j}f(x)$ if $k=1$ and $1\le j\le n$.
\end{defn}
Therefore the derivative $\partial_{j}\phi(f)_{x}$ is characterized
by the Taylor's formula\[
\forall h\in D_{n}^{d}\pti\phi\left[f(x+h)\right]=\sum_{\substack{j\in\N^{d}\\
|j|\le n}
}\frac{h^{j}}{j!}\cdot\partial_{j}\phi(f)_{x}\]
and by the conditions $\partial_{j}\phi(f)_{x}\in\ER_{\sss{=_{k_{j}}}}$
for every multi-index $j$. As above, with these notations we have
that $\partial_{j}\phi(f)_{x}=\partial_{j}(\phi\circ f)_{x}$.
\begin{example*}
\noindent Let us consider $f:D\freccia\ER$, we want to find $f'(i)$
for $i\in D$. By Definition \ref{def:derivativesOfFncsDefinedOnInfSets_x}
we have\begin{equation}
f(i+h)=f(i)+h\cdot f'(i)\quad\forall h\in D,\label{eq:exampleDerivativeGenericPointInfDom1}\end{equation}
with $f'(i)\in\ER_{\sss{=_{2}}}$. But $i+h\in D$ so that we also
get\begin{equation}
f(i+h)=f(0)+(i+h)\cdot f'(0)\label{eq:eq:exampleDerivativeGenericPointInfDom2}\end{equation}
with $f'(0)\in\ER_{\sss{=_{2}}}$. On the other hand from $i\in D$
we also have $f(i)=f(0)+i\cdot f'(0)$, and substituting in \eqref{eq:exampleDerivativeGenericPointInfDom1}
we get $f(i+h)=f(0)+i\cdot f'(0)+h\cdot f'(i)$. From this and from
\eqref{eq:eq:exampleDerivativeGenericPointInfDom2} we finally obtain\[
h\cdot f'(0)=h\cdot f'(i)\quad\forall h\in D,\]
that is (see Theorem \ref{thm:generalCancellationLaw}) $f'(0)=_{2}f'(i)$,
i.e. $f'(0)=f'(i)$ because both derivatives are in $\ER_{\sss{=_{2}}}$.
This confirms an intuitive result, i.e. that every smooth function
$f:D\freccia\ER$ is a straight line, and hence it has constant derivative.
\end{example*}
\noindent Using the notation of the Definition \ref{def:derivativesOfFncsDefinedOnInfSets_x}
we can now state the following
\begin{thm}
Let $X$, $d$, $n$, $U$, $V$ and $j$ as in the hypothesis of
Definition \ref{def:derivativesOfFncsDefinedOnInfSets_x}, then the
function\[
\partial_{j}:(f,\phi,x)\in\ext{U}^{V}\times\left(\ER^{\sf k}\right)^{\ext{U}}\times V\mapsto\partial_{j}\phi(f)_{x}\in\ER_{\sss{=_{k_{j}}}}^{\sf k}\]
is a $\ECInfty$ map.
\end{thm}
\noindent \textbf{Proof:} Let us consider figures $\alpha\In{A}\ext{U}^{V}$,
$\beta\In{B}\left(\ER^{\sf k}\right)^{\ext{U}}$ and $\gamma\In{C}V$,
we have to prove that $(\alpha\times\beta\times\gamma)\cdot\partial_{j}\In{A\times B\times C}\ER_{\sss{=_{k_{j}}}}^{\sf k}$.
Due to cartesian closedness we have that $\alpha^{\vee}:\bar{A}\times V\freccia\ext{U}$
and $\beta^{\vee}:\bar{B}\times\ext{U}\freccia\ER^{\sf k}$ are smooth.
In the following we will always use identifications of type $\bar{A}\times\bar{B}=\overline{A\times B}$,
based on the isomorphism of Lemma \ref{lem: ProductOfGeneratedObjects}
and Lemma \ref{lem: ProductInSERn}. We proceed locally, that is using
the sheaf property of the space $\ER_{\sss{=_{k}}}^{\sf k}$, so let
us fix generic $(a,b,c)\in A\times B\times C$. First of all, we have
to note that the function \[
\Theta:(a_{1},b_{1},x)\in A\times B\times V\mapsto\beta^{\vee}\left\{ b_{1},\alpha^{\vee}\left[a_{1},x\right]\right\} \in\ER^{\sf k},\]
is smooth, being the composition of smooth functions. Let, for simplicity,
$\mathcal{V}:=A\times B\times V$. We have that $\Theta:\overline{A\times B\times V}\freccia\ER^{\sf k}$
in $\ECInfty$, and hence we have the figure $\Theta\In{A\times B\times V}\ER^{\sf k}$.
Let us apply to this figure Theorem \ref{thm:figuresOfFermatSpaces}
at the point $(a,b,\gamma(c))\in\mathcal{V}$, obtaining that in an
open neighborhood $\ext{\mathcal{U}}\cap\mathcal{V}$ of $(a,b,\gamma(c))$
generated by the open set $\mathcal{U}$ of $\R^{\sf a}\times\R^{\sf b}\times\R^{d}$
we can write $\Theta(a_{1},b_{1},x)=\ext{\delta}(p,a_{1},b_{1},x)$
for every $(a_{1},b_{1},x)\in\ext{\mathcal{U}}\cap\mathcal{V}$, where
$\delta\in\CInfty(\mathcal{A}\times\mathcal{U},\R^{\sf k})$, $p\in\ext{\mathcal{A}}$
and $\mathcal{A}$ is open in $\R^{\sf p}$. Note that, being $\R^{\sf k}$
a manifold, we do not have the classical alternative {}``$\Theta$
is locally constant or $\Theta=\ext{\delta}(p,-)$'', see the above
cited theorem. Roughly speaking, to obtain the map $(\alpha\times\beta\times\gamma)\cdot\partial_{j}$
of our conclusion we have to derive the function $\delta$ with respect
to the fourth variable $x$, to compose the result with the mapping
$\iota_{k_{j}}:\ER^{\sf k}\freccia\ER_{\sss{=_{k_{j}}}}^{k}$ and
finally to compose the final result with the figure $\gamma:C\freccia V$.
To formalize this reasoning, we start from the open set $\mathcal{U}$
of $\R^{\sf a}\times\R^{\sf b}\times\R^{d}$. Therefore, $p_{3}(\mathcal{U})$
is open in $\R^{d}$, where $p_{3}:\R^{\sf a}\times\R^{\sf b}\times\R^{d}\freccia\R^{d}$
is the projection onto the third space. Hence $p_{3}(\mathcal{U})\cap V$
is open in $V$ and $\gamma^{-1}\left[p_{3}(\mathcal{U})\cap V\right]=:\mathcal{C}$
is open in $C$ because $\gamma:C\freccia V$, being smooth, is continuous.
On the other hand, $p_{12}(\mathcal{U})$ is open in $\R^{\sf a}\times\R^{\sf b}$,
where $p_{12}:\R^{\sf a}\times\R^{\sf b}\times\R^{d}\freccia\R^{\sf a}\times\R^{\sf b}$
is the projection onto the first two factors. Therefore, $p_{12}(\mathcal{U})\cap(A\times B)=:\mathcal{D}$
is open in $A\times B$ and hence $\mathcal{D}\times\mathcal{C}$
is open in $A\times B\times C$. In this open set we will realize
the above mentioned compositions. Indeed, first of all we have that
$(a,b,\gamma(c))\in\mathcal{U}$ and hence $(a,b)\in\mathcal{D}=p_{12}(\mathcal{U})\cap(A\times B)$;
moreover, $\gamma(c)\in p_{3}(\mathcal{U})\cap V$ and hence $(a,b,c)\in\mathcal{D}\times\mathcal{C}$.
Now, for a generic $(a_{1},b_{1},c_{1})\in\mathcal{D}\times\mathcal{C}$
we have\begin{align*}
\left[(\alpha\times\beta\times\gamma)\cdot\partial_{j}\right](a_{1},b_{1},c_{1}) & =\partial_{j}\beta(b_{1})\left(\alpha(a_{1})\right)_{\gamma(c_{1})}\\
 & =\partial_{j}\left(\beta(b_{1})\circ\alpha(a_{1})\right)_{\gamma(c_{1})}\\
 & =\partial_{j}\left(\beta^{\vee}(b_{1},-)\circ\alpha^{\vee}(a_{1},-)\right)_{\gamma(c_{1)}}\\
 & =\partial_{j}\left(\Theta(a_{1},b_{1},-)\right)_{\gamma(c_{1})}\\
 & =\iota_{k_{j}}\left[\ext{\left(\partial_{4}\delta\right)}(p,a_{1},b_{1},\gamma(c_{1}))\right].\end{align*}
This proves that we can express $(\alpha\times\beta\times\gamma)\cdot\partial_{j}$
on the open neighborhood $\mathcal{D}\times\mathcal{C}$ of $(a,b,c)$
as a composition of smooth maps (so, here we are using Theorem \ref{thm:i_k_IsSmooth}
about the smoothness of the map $\iota_{k_{j}}$) and hence the conclusion
follows from the sheaf property of the space $\ER_{\sss{=_{k_{j}}}}^{\sf k}\in\ECInfty$.$\qedNoNewLine$

Using this result, or the analogous for derivative of smooth functions
defined on open sets, we can easily extend the Taylor's formula to
vector spaces of smooth functions of the form $\ER^{Z}$.
\begin{thm}
\label{thm:TaylorForFunctionsSpaces}Let $Z$ be a $\CInfty$ space,
$\alpha_{1},\ldots,\alpha_{d}\in\R_{>0}$, and\[
f:D_{\alpha_{1}}\times\cdots\times D_{\alpha_{d}}\xfrecciad{}\ER^{Z}\]
be a $\ECInfty$ map. Define $k_{j}\in\R$ such that\[
k_{\underline{0}}:=0\]
\[
\frac{1}{k_{j}}+\frac{j}{\alpha+1}=1\quad\forall j\in J:=\left\{ j\in\N^{d}\,|\,\frac{j}{\alpha+1}<1\right\} \pti j\ne\underline{0}.\]
Then there exists one and only one family of smooth functions\[
m:J\xfrecciad{}\sum_{j\in J}\ER_{\sss{=_{k_{j}}}}^{Z}\]
such that
\begin{enumerate}
\item $m_{j}(z)\in\ER_{\sss{=_{k_{j}}}}^{\sf k}$ for every $j\in\N^{d}$
such that $\frac{j}{\alpha+1}<1$ and every $z\in Z$
\item $f(h)=\sum\limits _{\substack{j\in\N^{d}\\
\frac{j}{\alpha+1}<1}
}\frac{h^{j}}{j!}\cdot m_{j}\quad\forall h\in D_{\alpha_{1}}\times\cdots\times D_{\alpha_{d}}$ in the vector space $\ER^{Z}$.
\end{enumerate}
\end{thm}
\noindent \textbf{Proof:} For cartesian closedness, the adjoint of
the map $f$ is smooth:\[
f^{\vee}:Z\times D_{\alpha_{1}}\times\cdots\times D_{\alpha_{d}}\xfrecciad{}\ER.\]
Let us indicate this map, for simplicity, again with $f(-,-)$. Then
it suffices to consider the smooth functions\[
m_{j}:z\in Z\mapsto f(z,-)\in\ER^{D_{\alpha_{1}}\times\cdots\times D_{\alpha_{d}}}\mapsto\partial_{j}f(z,-)\in\ER_{\sss{=_{k_{j}}}}\]
to obtain the desired pointwise equality\[
f(z,h)=\sum\limits _{\substack{j\in\N^{d}\\
\frac{j}{\alpha+1}<1}
}\frac{h^{j}}{j!}\cdot m_{j}(z)\quad\forall z\in Z\]
from the infinitesimal Taylor's formula (Theorem \ref{thm:GeneralizedDerivationFormula_d-dim}).
The uniqueness part follows from the corresponding uniqueness of the
cited theorem.$\qedNoNewLine$

Using analogous ideas, we can also extend Theorem \ref{thm:infPolynomialsWithSmoothCoefficients}
to functions of the form $f:S\freccia\ER^{Z}$, where $Z$ is a generic
$\ECInfty$ space; the corresponding statement can be easily obtain
simply replacing in Theorem \ref{thm:infPolynomialsWithSmoothCoefficients}
the space $\ER^{n}$ with the space $\ER^{Z}$.

\chapter{\label{cha:Infinitesimal-differential-geometry}Infinitesimal differential
geometry}

The use of nilpotent infinitesimals permits to develop many concepts
of differential geometry in an intrinsic way, without being forced
to use coordinates. In this way the use of charts becomes specific
of suitable areas of differential geometry, e.g. where one strictly
needs some solution in a finite neighborhood and not in an infinitesimal
one only (e.g. this is the case for the inverse function theorem).

We recall that we named this kind of intrinsic geometry \emph{infinitesimal
differential geometry}.

\noindent The possibility to avoid coordinates using infinitesimal
neighborhoods instead, permits to perform some generalizations to
more abstract spaces, like spaces of mappings. Even if the categories
$\CInfty$ and $\ECInfty$ are very big and not very much can be said
about generic objects, in this section we shall see that the best
properties can be formulated for a restricted class of extended spaces,
the \emph{infinitesimally linear }ones, to which spaces of mappings
between manifolds belong to.

All this section takes strong inspiration from the corresponding part
of SDG, in the sense that all the statements of theorems and definitions
have a strict analogue in SDG. Our reference for proofs not depending
on the model presented in this work but substantially identical to
those in SDG is \citet{Lav}.

\section{Tangent spaces and vector fields}

We start from the fundamental idea of tangent vector. It is natural
to define a tangent vector to a space $X\in\ECInfty$ as an arrow
(in $\ECInfty$) of type $t:D\freccia X$. Therefore ${\rm T}X:=X^{D}=\ECInfty(D,X)$
with projection $\pi:t\in{\rm T}X\mapsto t(0)\in X$ is the tangent
bundle of $X$.

We can also define the differential of an application $f:X\freccia Y$
in $\ECInfty$ simply by composition\[
\diff{f}:t\in\text{T}X=X^{D}\freccia f\circ t\in\text{T}Y=Y^{D}.\]

\noindent In the following we will also use the notations\[
\text{T}_{x}X:=\left(\left\{ t\in\text{T}X\,|\, t(0)=x\right\} \prec\text{T}M\right)\in\ECInfty\]
\[
\diff{f}_{x}:t\in\text{T}_{x}X\mapsto\diff{f}_{x}[t]:=f\circ t\in\text{T}_{f(x)}Y\]
for the tangent space at the point $x\in X$ and for the differential
of the application $f:X\freccia Y$ at the point $x$.

Note that using the absolute value it is also possible to consider
{}``boundary tangent vectors'' taking $|D|:=\{\;|h|:\, h\in D\}$
instead of $D$, for example at the initial point of a curve or at
a point in the boundary of a closed set. In the following, $M\in\ManInfty$
will always be a $d$-dimensional smooth manifold and we will use
the simplified notation ${\rm T}M$ for ${\rm T}(\ext{M})$.

\noindent It is important to note that with this definition of tangent
vector we obtain a generalization of the classical notion. In fact,
in general we have that $t(0)\in\ext{M}$ and $\phi'(t):=\partial_{1}\phi(t)\in\ER_{\sss{=_{2}}}^{d}$
if $\phi$ is a chart on $\st{t(0)}\in M$. In other words, a tangent
vector $t:D\freccia\ext{M}$ can be applied to a non standard point
or have a non standard speed. If we want to study classical tangent
vectors we have to consider the following $\CInfty$ object
\begin{defn}
\noindent \label{def: StndTangentFunctor} We call ${\rm T}_{{\rm st}}M$
the $\CInfty$ object with support set \[
\left|{\rm T}_{{\rm st}}M\right|:=\{\ext{f}|_{D}\,:\, f\in\CInfty(\R,M)\},\]
and with figures of type $U$ (open in $\R^{\sf u}$) given by the
substructure induced by $\text{T}M$, i.e. \[
d\In{U}{\rm T}_{{\rm st}}M\DIff d:U\freccia\left|{\rm T}_{{\rm st}}M\right|\e{and}\ECInfty\vDash d\cdot i\In{\bar{U}}{\rm T}M,\]
where $i:\left|{\rm T}_{{\rm st}}M\right|\hookrightarrow{\rm T}M$
is the inclusion.
\end{defn}
\noindent That is in ${\rm T}_{{\rm st}}M$ we consider only tangent
vectors of the form $t=\ext{f}|_{D}$, i.e. obtained as extension
of ordinary smooth functions $f:\R\freccia M$, and we take as figures
of type $U\subseteq\R^{\sf u}$ the functions $d$ with values in
${\rm T}_{{\rm st}}M$ which in the category $\ECInfty$ verify $d^{\vee}:\bar{U}\times D\freccia\ext{M}$.
Note that, intuitively speaking, $d$ takes a standard element $u\in U\subseteq\R^{k}$
to the standard element $d(u)\in{\rm T}_{{\rm st}}M$.
\begin{thm}
\noindent \label{thm: StandardTangentVector} Let $t\in{\rm T}M$
be a tangent vector and $(U,\phi)$ a chart of $M$ on $\st{t(0)}$.
Then \[
t\in{\rm T}_{{\rm st}}M\iff t(0)\in M\e{and}\phi'(t)\in\R^{d}.\]

\end{thm}
\noindent \textbf{Proof:} If $t=\ext{f}|_{D}\in\text{T}_{\text{st}}M$
then $t(0)=f(0)=\st{t(0)}\in M$ and $\phi'(t)=\iota_{2}\left[\left(\phi\circ f\right)'(0)\right]\in\R^{d}$
because\[
\CInfty\vDash V:=f^{-1}(U)\xfrecciad{f|_{V}}U\xfrecciad{\phi}\R^{d}\]
and hence $\left(\phi\circ f\right)'(0)\in\R^{d}$. Vice versa if
$t(0)\in M$ and $\phi'(t)\in\R^{d}$, then applying the generalized
derivation formula (Theorem \ref{thm: GeneralizedDerivationFormula})
we obtain $\ext{\varphi}(t(h))=\ext{\phi}(t(0))+h\cdot\phi'(t)$ for
any $h\in D$. But $\ext{\varphi}(t(0))=\varphi(t(0))$ because $t(0)\in M$.
Hence setting $a:=\phi(t(0))\in\R^{d}$ and $b:=\phi'(t)\in\R^{d}$
we can define\[
f(s):=\phi^{-1}(a+s\cdot b)\in M\quad\forall s:\ |s|<r,\]
where $r\in\R_{>0}$ has been taken such that $B_{r\cdot\Vert b\Vert}(a)\subseteq U$.
The standard smooth function $f:(-r,r)\freccia M$ can be defined
on the whole of $\R$ in any way that preserves its smoothness.

We have that$ $\[
t(h)=\ext{\varphi}^{-1}(\ext{\phi(t(h))})=\ext{\phi}^{-1}(a+h\cdot b)=:\ext{f}|_{D}(h)\quad\forall h\in D,\]
and this proves that $t\in\text{T}_{\text{st}}M$ is a standard tangent
vector.$\qedNoNewLine$

\noindent In the following result we prove that the definition of
standard tangent vector $t\in{\rm T}_{{\rm st}}M$ is equivalent to
the classical one. 
\begin{thm}
\noindent In the category $\CInfty$ the object ${\rm T}_{{\rm st}}M$
is isomorphic to the usual tangent bundle of $M$
\end{thm}
\noindent \textbf{Proof:} We have to prove that ${\rm T}_{{\rm st}}^{m}:=\{t\in{\rm T}_{{\rm st}}M\,|\, t(0)=m\}\simeq{\rm T}_{m}$
where here ${\rm T}_{m}:=\{f\in\Cc^{\infty}(\R,M)\,|\, f(0)=m\}/\sim$
is the usual tangent space of $M$ at $m\in M$. Note that ${\rm T}_{m}\in\CInfty$
because of completeness and co-completeness.

\noindent Firstly we prove that \begin{align}
\alpha:\quad[f]_{\sim}\in{\rm T}_{m}\quad & \mapsto\quad\frac{\diff{(\varphi\circ f)}}{\diff t}(0)\in\R^{d}\nonumber \\
\alpha^{-1}:\quad v\in\R^{d}\quad & \mapsto\quad[r\mapsto\varphi^{-1}(\varphi m+r\cdot v)]_{\sim}\in{\rm T}_{m}\label{eq:1_standardTangentSpace}\end{align}
are arrows of $\CInfty$, where $\varphi:V\freccia\R^{d}$ is a chart
on $m$ with $\varphi(V)=\R^{d}$.

\noindent Secondly we prove that \begin{align*}
\beta:\quad t\in{\rm T}_{{\rm st}}^{m}\quad & \mapsto\quad\varphi'(t)\in\R^{d}\\
\beta^{-1}:\quad v\in\R^{d}\quad & \mapsto\quad\ext{[r\mapsto\varphi^{-1}(\varphi m+r\cdot v)]}|_{D}\in{\rm T}_{{\rm st}}^{m}\end{align*}
are arrows of $\CInfty$.

To prove the smoothness of $\alpha:\text{T}_{m}\freccia\R^{d}$ in
$\CInfty$, let us take a figure $d\In{H}\text{T}_{m}$, where $H$
is open in $\R^{\sf h}$. For the sheaf property of $\text{T}_{m}$
to prove that $d\cdot\alpha:H\freccia\R^{d}$ is smooth we can proceed
proving that it is locally smooth.

From the definition of figures of the quotient space $\text{T}_{m}$
(see Theorem \ref{thm: limit in Cn}), for every $h\in H$ there exist
an open neighborhood $U$ of $h$ in $H$ and a smooth function $\delta\in\CInfty(\R,M^{\R})$
such that \[
\forall u\in U\pti\delta(u)(0)=m\]
\[
d|_{U}=\delta\cdot[-]_{\sim},\]
where $[-]_{\sim}:\{f\in\Cc^{\infty}(\R,M)\,|\, f(0)=m\}\freccia\text{T}_{m}$
is the canonical projection map of the quotient set $\text{T}_{m}$.

Thus, we have\[
\left(d\cdot\alpha\right)|_{U}=d|_{U}\cdot\alpha:u\in U\mapsto\alpha\left(du\right)=\alpha\left[\delta u\right]_{\sim}=\frac{\diff{(\phi\circ\delta u)}}{\diff{t}}(0).\]
But \[
\forall r\in\R\pti(\phi\circ\delta u)(r)=\phi\left[\delta(u)(r)\right]=\phi\left[\delta^{\vee}(u,r)\right],\]
and for the cartesian closedness of $\CInfty$ we have that $\delta^{\vee}:U\times\R\freccia M$
is smooth. Therefore\[
(d\cdot\alpha)|_{U}:u\in U\mapsto\diff{}_{m}\phi\left[\partial_{2}\delta^{\vee}(u,0)\right]\in\R^{d},\]
where $\diff{}_{m}\phi$ is the differential of $\phi$ at the point
$m\in M$. We thus have that $(d\cdot\alpha)|_{U}=\diff{}_{m}\phi\left[\partial_{2}\delta^{\vee}(-,0)\right]\in\CInfty(U,\R^{d})$
which proves that $d\cdot\alpha$ is locally smooth.

Now, let us consider the inverse $\alpha^{-1}$ defined in \eqref{eq:1_standardTangentSpace}.
This map is exactly the composition of the adjoint $\tilde{\alpha}^{\wedge}$
of the smooth map\[
\tilde{\alpha}:(v,r)\in\R^{d}\times\R\mapsto\phi^{-1}\left(\phi m+r\cdot v\right)\in M\]
with the canonical projection:\begin{align*}
\xymatrix{\xyR{50pt}\xyC{50pt}\R^{d}\ar[r]\sp(0.33){\tilde{\alpha}^{\wedge}}\ar[rd]_{\alpha^{-1}} & \left\{ f\in M^{\R}\,|\, f(0)=m\right\} \ar[d]^{[-]_{\sim}}\\
 & \text{T}_{m}}
\end{align*}
and hence it is smooth because of the type of figures we have on a
quotient set, see Theorem \ref{thm: limit in Cn}.

We now prove that $\beta:\text{T}_{\text{st}}^{m}M\freccia\R^{d}$
is smooth. If $d\In{U}{\rm T}_{{\rm st}}^{m}$ is a figure in $\ECInfty$
of type $U$, where $U$ is an open set of $\R^{\sf u}$, then $d^{\vee}:\bar{U}\times D\freccia\ext{M}$
in $\ECInfty$. But $\bar{U}\times D=\bar{U}\times\bar{D}=\overline{U\times D}$
hence $d^{\vee}\In{U\times D}\ext{M}$. Thus, for every $u\in U$
we can locally write $d^{\vee}|_{\mathcal{V}}=\ext{\gamma}(p,-,-)|_{\mathcal{V}}$
where $\mathcal{V}:=\ext{(A\times B)}\cap(U\times D)$ is an open
neighborhood of $(u,0)$ defined by $A\times B$ in $U\times D$,
$\gamma\in\Cc^{\infty}(W\times A\times B,M)$ is an standard smooth
function, and $p\in\ext{W}$, where $W$ is open in $\ER^{\sf p}$.
But\[
\mathcal{V}=\ext{(A\times B)}\cap(U\times D)=(\ext{A}\cap U)\times(\ext{B}\cap D)=(A\cap U)\times D\]
because $U\subseteq\R^{\sf u}$. Now, we have\begin{align*}
\beta[d(x)] & =\varphi'[d(x)]\\
 & =\phi'\left[d^{\vee}(x,-)\right]\\
 & =\iota_{2}\left\{ \frac{\diff{}}{\diff{r}}\{\varphi[\gamma(p,x,r)]\}|_{r=0}\right\} \\
 & =\iota_{2}\left\{ \diff{}_{m}\phi\left[\partial_{3}\gamma(p,x,0)\right]\right\} \quad\forall x\in A\cap U.\end{align*}
But $\st{\left\{ \beta[d(x)]\right\} }=\beta[d(x)]$ because $\beta:\text{T}_{\text{st}}\freccia\R^{d}$
and hence \begin{align}
\beta[d(x)] & ={}^{{}^{{}\circ}}\left[\iota_{2}\left\{ \diff{}_{m}\phi\left[\partial_{3}\gamma(p,x,0)\right]\right\} \right]\nonumber \\
 & ={}^{{}^{{}\circ}}\left[\diff{}_{m}\phi\left[\partial_{3}\gamma(p,x,0)\right]\right]\nonumber \\
 & =\diff{}_{m}\phi\left[\partial_{3}\gamma(p_{0},x,0)\right]\quad\forall x\in A\cap U,\label{eq:2_standardTangentSpace}\end{align}
so that $(d\cdot\beta)|_{A\cap U}=\diff{}_{m}\phi\left[\partial_{3}\gamma(p_{0},-,0)\right]\in\CInfty(A\cap U,\R^{d})$
is an ordinary smooth function. Note the importance to have as $U$
a standard open set in the last passage of \eqref{eq:2_standardTangentSpace},
and this represents a further strong motivation for the definition
we gave for ${\rm T}_{{\rm st}}M$.

To prove the regularity of $\beta^{-1}$ we consider the map\[
\tilde{\beta}:v\in\R^{d}\mapsto\left[r\in\R\mapsto\phi^{-1}(\phi m+r\cdot v)\in M\right]\in\CInfty(\R,M).\]
Then we have\[
\CInfty\vDash\R^{d}\xfrecciad{\tilde{\beta}}M^{\R}\]
\[
\CInfty\vDash\R^{d}\times\R\xfrecciad{\tilde{\beta}^{\vee}}M\]
\[
\ECInfty\vDash\ER^{d}\times\ER\xfrecciad{\ext{\tilde{\beta}^{\vee}}}\ext{M}\]
\[
\ECInfty\vDash\R^{d}\times D\xfrecciad{\bar{\beta}}\ext{M}\e{where}\bar{\beta}:={}^{{}^{{}{\scriptstyle \bullet}}}\left(\tilde{\beta}^{\vee}\right)|_{\R\times D}\]
\[
\ECInfty\vDash\R^{d}\xfrecciad{\bar{\beta}^{\wedge}}\ext{M}^{D}.\]
This map is actually $\beta^{-1}$, in fact\begin{align*}
\forall v\in\R^{d}\pti\bar{\beta}^{\wedge}(v):h\in D\mapsto\bar{\beta}(v,h) & ={}^{{}^{{}{\scriptstyle \bullet}}}\left(\tilde{\beta}^{\vee}\right)(v,h)\\
 & =\left(\tilde{\beta}^{\vee}(v,h_{t})\right)_{t\ge0}\\
 & =\left(\tilde{\beta}(v)(h_{t})\right)_{t\ge0}\\
 & =\beta^{-1}(v)(h).\end{align*}
Therefore $\beta^{-1}:\R^{d}\freccia\text{T}M$ is smooth in $\ECInfty$.
But, finally, $\beta^{-1}$ is actually with values in $\text{T}_{\text{st}}^{m}$
because $\beta^{-1}(v)(0)=m$ for every $v\in\R^{d}$, so that $\beta^{-1}=\beta^{-1}\cdot i\In{\bar{\R^{d}}}\text{T}M$
where $i:\text{T}_{\text{st}}^{m}\hookrightarrow\text{T}M$ is the
inclusion. We have thus prove that $\beta^{-1}$ is a figure of type
$\R^{d}$ of the $\CInfty$ space $\text{T}_{\text{st}}^{m}$ and
hence it is also smooth in this category, which is the conclusion.$\qedNoNewLine$

For any object $X\in\ECInfty$ the multiplication of a tangent vector
$t$ by a scalar $r\in\ER$ can be defined simply {}``increasing
its speed'' by a factor $r$: \begin{equation}
(r\cdot t)(h):=t(r\cdot h)\quad\forall h\in D.\label{eq:definitionOfProductByScalar}\end{equation}
But, as we have already noted, in the category $\ECInfty$ we have
spaces with singular points too, like algebraic curves with double
points. For this reason, we cannot define the sum of tangent vectors
for every smooth space $X\in\ECInfty$, but we need to introduce a
class of objects in which this operation is possible.

\noindent The following definition simply states that in these spaces
there always exists the infinitesimal parallelogram generated by a
finite number of given vectors at the same point $m$.
\begin{defn}
\noindent \label{def:inf-linear}Let $X\in\ECInfty$, then we say
that $X$ is \emph{infinitesimally linear}, or simply \emph{inf-linear},
\emph{at the point }$m\in X$ if and only if the following conditions
are fulfilled
\begin{enumerate}
\item \noindent for any $k\in\N_{>1}$ and for any $t_{1},\ldots,t_{k}\in{\rm T}_{m}X$,
there exists one and only one $p:D^{k}\freccia X$ in $\ECInfty$
such that \[
\forall i=1,\ldots,k\pti p(0,\ptind^{i-1},0,h,0,\ldots,0)=t_{i}(h)\quad\forall h\in D.\]
We will call the map $p$ \emph{the infinitesimal parallelogram generated
by$ $} $t_{1},\ldots,t_{k}$.
\item The application\[
(-)+_{m}\ldots+_{m}(-):(t_{1},\ldots,t_{k})\in\left(\text{T}_{m}X\right)^{k}\mapsto p\in X^{D^{k}}\]
that associates to the $k$ tangent vectors at $m\in X$ the infinitesimal
parallelogram $p$, is $\ECInfty$-smooth.
\end{enumerate}
\noindent Moreover, we will simply say that $X$ is \emph{inf-linear}
if it is inf-linear at each point $m\in X$ and if the application\[
m\in X\mapsto(-)+_{m}\ldots+_{m}(-)\in\sum_{m\in X}\left(X^{D^{k}}\right)^{\left(\text{T}_{m}X\right)^{k}}\]
is $\ECInfty$-smooth.

\end{defn}
\begin{figure}[h]
\noindent \begin{centering}
\includegraphics[scale=0.6]{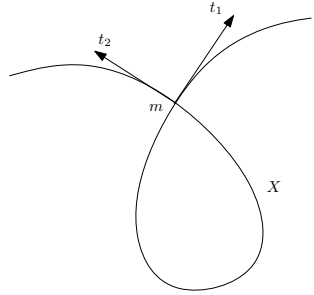}
\par\end{centering}

\caption{An example of space which is not inf-linear at $m\in X$.}

\end{figure}

\noindent The following theorem gives meaningful examples of inf-linear
objects.
\begin{thm}
\noindent \label{thm:manifoldsAndExpOfManifoldsAreInfLinear}The extension
of any manifold $\ext{M}$ is inf-linear at every point $m\in\ext{M}$.
If $M_{i}\in\ManInfty$ for $i=1,\ldots,s$ then the exponential object
\begin{equation}
\ext{M_{1}}^{{\ext{M_{2}}^{\,\dots}}^{\ext{M_{s}}}}\simeq\ext{M_{1}}^{\ext{(M_{2}\times\dots\times M_{{\scriptstyle s}})}}\label{eq:isoOfFunctionSpaces}\end{equation}
 is also inf-linear at every point.
\end{thm}
\noindent The importance of the isomorphism \eqref{eq:isoOfFunctionSpaces}
lies in the fact that complex spaces like\[
\ext{M_{1}}^{{\ext{M_{2}}^{\,\dots}}^{\ext{M_{s}}}}=\ECInfty(\ext{M_{s}},\ECInfty(\ext{M_{s-1}},\cdots,\ECInfty(\ext{M_{2}},\ext{M_{1}})\cdots)\]
are now no more difficult to handle than classical spaces of mappings
like $\ext{M}^{\ext{N}}=\ECInfty(\ext{N},\ext{M})$. Let us note explicitly
that this isomorphism is a consequence of cartesian closedness and
of the preservation of products of manifolds of the Fermat functor.

\noindent \textbf{Proof:} Given any chart $(U,\varphi)$ on $\st{m}\in M$
we can define the infinitesimal parallelogram $p$ as \begin{equation}
p(h_{1},\ldots,h_{k}):=\ext{\varphi}^{-1}\left(\ext{\varphi}(m)+\sum_{i=1}^{k}h_{i}\cdot\varphi'(t_{i})\right)\quad\forall h_{1},\ldots,h_{k}\in D.\label{eq:defOfInfParallelogram}\end{equation}
If fact if $\tau(h):=p(0,\ptind^{i-1},0,h,0,\ldots,0)$ then $\varphi(\tau(h))=\varphi(m)+h\cdot\varphi'(t_{i})$
for every $h\in D$; this implies that $t(0)=\tau(0)$ and $\varphi'(\tau)=\varphi'(t_{i})$,
hence $t_{i}=\tau$. To prove the uniqueness of the parallelogram
generated by $t_{1},\ldots,t_{k}\in\text{T}_{m}M$, let us consider
that if $p:D^{k}\freccia\ext{M}$ is such that $p(0,\ptind^{i-1},0,h,0,\ldots,0)=t_{i}(h)$
for every tangent vector $t_{i}$ and every $h\in D$, then \[
\varphi[p(0,\ptind^{i-1},0,h,0,\ldots,0)]=\varphi[t_{i}(h)]=\varphi(m)+h\cdot\varphi'(t_{i})\]
and so \[
\varphi[p(h_{1},\ldots,p_{k})]=\varphi(m)+\sum_{i=1}^{k}h_{i}\cdot\varphi'(t_{i})\]
from the first order infinitesimal Taylor's formula, so that we obtain
again the definition \eqref{eq:defOfInfParallelogram}, and this proves
the uniqueness part.

Considering the exponential object, because \[
\ext{M_{1}}^{{\ext{M_{2}}^{\dots}}^{\ext{M_{s}}}}\simeq\ext{M_{1}}^{\ext{M_{2}}\times\dots\times\ext{M_{{\scriptstyle s}}}}\simeq\ext{M_{1}}^{\ext{(M_{2}\times\dots\times M_{{\scriptstyle s}})}},\]
it suffices to prove the conclusion for $s=2$. First of all we note
that, because of the previously proved uniqueness, the definition
\ref{eq:defOfInfParallelogram} of the infinitesimal parallelogram
does not depend on the chart $\varphi$ on $\st{m}$. Now let $t_{1},\ldots,t_{k}$
be $k$ tangent vectors at $f\in\ext{N}^{\ext{M}}$. We shall define
their parallelogram%
\footnote{For simplicity, in this proof we will use implicitly the cartesian
closedness property.%
} $p:\ext{M}\freccia\ext{N}^{D^{k}}$ patching together smooth functions
defined on open subsets, and using the sheaf property of $\ext{N}^{D^{k}}$.
Indeed, for every $m\in\ext{M}$ we can find a chart $(U_{m},\varphi_{m})$
of $N$ on $\st{f(m)}$ with $\varphi_{m}(U_{m})=\R^{n}$. Now $m\in V_{m}:=f^{-1}(\ext{U_{m}})\in\Top{\ext{M}}$
and for every $x\in V_{m}$ we have $t_{i}^{\vee}(0,x)=f(x)\in\ext{U_{m}}$.
Hence $t_{i}^{\vee}(h,x)\in\ext{U_{m}}$ for any $h\in D$ by Theorem
\ref{thm: DandOpenGeneralized}. Therefore we can define \begin{equation}
p_{m}^{\vee}(x,h):=\varphi_{m}^{-1}\left\{ \sum_{i=1}^{k}\varphi_{m}[t_{i}^{\vee}(h^{i},x)]-(k-1)\cdot\varphi_{m}(fx)\right\} \enskip\,\forall x\in V_{m},\forall h\in D^{k}\label{eq:1_expOfManifoldsAreInf-Linear}\end{equation}
 and we have that $p_{m}^{\vee}:(V_{m}\prec\ext{M})\times D^{k}\freccia\ext{N}$
is smooth, because it is a composition of smooth functions. We claim
that if $x\in V_{m}\cap V_{m^{\prime}}$ then $p_{m}^{\vee}(x,-)=p_{m^{\prime}}^{\vee}(x,-)$,
in fact from the generalized Taylor's formula we have $\varphi_{m}[t_{i}^{\vee}(h^{i},x)]=\varphi_{m}(fx)+h^{i}\cdot\varphi_{m}^{\prime}[t_{i}^{\vee}(-,x)]$
and hence substituting in \eqref{eq:1_expOfManifoldsAreInf-Linear}
we can write \begin{align}
p_{m}^{\vee}(x,h) & =\varphi_{m}^{-1}\left\{ k\varphi_{m}(fx)+\sum_{i=1}^{k}h^{i}\cdot\varphi_{m}^{\prime}[t_{i}^{\vee}(-,x)]-k\varphi_{m}(fx)+\varphi_{m}(fx)\right\} \label{eq: ParallelogramFunctionsSpace}\\
 & =\varphi_{m}^{-1}\left\{ \varphi_{m}(fx)+\sum_{i=1}^{k}h^{i}\cdot\varphi_{m}^{\prime}[t_{i}^{\vee}(-,x)]\right\} \quad\forall x\in V_{m},\forall h\in D^{k}.\nonumber \end{align}
 But $(U_{m},\varphi_{m})$ is a chart on $\st{f(x)}$, so $p_{m}^{\vee}(x,-)$
is the infinitesimal parallelogram generated by the tangent vectors
$t_{i}^{\vee}(-,x)$ at $f(x)$, and we know that \eqref{eq: ParallelogramFunctionsSpace}
does not depend on $\varphi_{m}$, so $p_{m}=p_{m^{\prime}}$. For
the sheaf property of $\ext{N}^{D^{k}}$ we thus have the existence
of a smooth $p:\ext{M}\freccia\ext{N}^{D^{k}}$ such that\[
\forall m\in\ext{M}\pti p|_{V_{m}}=p_{m}.\]
From this and from \eqref{eq: ParallelogramFunctionsSpace} it is
also easy to prove that $p:D^{k}\freccia\ext{N}^{\ext{M}}$ verifies
the desired properties. Uniqueness follows noting that $p^{\vee}(m,-)$
is the infinitesimal parallelogram generated by $t_{i}^{\vee}(-,m)$.
From \eqref{eq:1_expOfManifoldsAreInf-Linear} it also follow easily
that the map $(m,t_{1},\ldots,t_{k})\mapsto p$ is smooth because
it is given by the composition of smooth maps.$\qedNoNewLine$

Another important family of inf-linear spaces is given by the following
\begin{thm}
\label{thm:euclideanSpacesAreInf-Linear}Let $X$ be an inf-linear
space and $Z\in\ECInfty$ be another smooth space. Then the space\[
X^{Z}\]
is inf-linear.
\end{thm}
\noindent \textbf{Proof:} Let $t_{1},\ldots,t_{k}:D\freccia X^{Z}$
be $k\in\N_{>1}$ tangent vectors at the point $m\in X^{Z}$. Because
of cartesian closedness the adjoint maps $t_{i}^{\vee}:Z\times D\freccia X$
are smooth; we will simply denote them with the initial symbol $t_{i}$
again. Finally, for every $z\in Z$, because $X$ is inf-linear, we
know that the map\[
(-)+_{m(z)}\ldots+_{m(z)}(-):\left(\text{T}_{m(z)}X\right)^{k}\freccia X^{D^{k}}\]
is smooth in $z\in Z$ because it is composed by smooth functions.

Then the adjoint of the map: \[
p(z,h_{1},\ldots,h_{k}):=\left[t_{1}(z,-)+_{m(z)}\ldots+_{m(z)}t_{k}(z,-)\right](h_{1},\ldots,h_{k})\]
verifies the desired properties.$\qedNoNewLine$

If $X$ is inf-linear at $x\in X$ then we can define the sum of tangent
vectors $t_{1},t_{2}\in{\rm T}_{x}X$ simply taking the diagonal of
the parallelogram $p$ generated by these vectors \begin{equation}
(t_{1}+t_{2})(h):=p(h,h)\quad\forall h\in D.\label{eq:definitionOfSumTangentVectors}\end{equation}
With these operations ${\rm T}_{x}X$ becomes a $\ER$-module:
\begin{thm}
If $X$ is inf-linear at the point $x\in X$, then with respect to
the sum defined in \eqref{eq:definitionOfSumTangentVectors} and the
product by scalar defined by \eqref{eq:definitionOfProductByScalar},
the tangent space $\text{T}_{x}X$ is a $\ER$-module.
\end{thm}
\noindent \textbf{Proof:} We only prove that the sum is associative.
Analogously, one can prove the other axioms of module, . Let us consider
the tangent vectors $t_{1}$, $t_{2}$, $t_{3}\in\text{T}_{x}X$,
and denote by $p_{12}$ the infinitesimal parallelogram generated
by $t_{1}$ and $t_{2}$, by $p_{12,3}$ the parallelogram generated
by $t_{1}+t_{2}$ and by $t_{3}$, and analogously for the symbols
$p_{23}$ and $p_{1,23}$. Then $p_{12,3}$ is characterized by the
properties\[
p_{12,3}:D^{2}\freccia X\]
\[
p_{12,3}(h,0)=(t_{1}+t_{2})(h)=p_{12}(h,h)\quad\forall h\in D\]
\[
p_{12,3}(0,h)=t_{3}(h)\quad\forall h\in D.\]
Now, let $l:D^{3}\freccia X$ be the parallelogram generated by all
the three vectors. Then the map $l(-,-,0)$ verifies\[
l(-,-,0):D^{2}\freccia X\]
\[
l(h,0,0)=t_{1}(h)\e{and}l(0,h,0)=t_{2}(h)\quad\forall h\in D,\]
so $l(-,-,0)=p_{12}$. Now let us consider the application\[
\lambda:(h,k)\in D^{2}\mapsto l(h,h,k)\in X.\]
It is smooth as a composition of smooth maps and verifies\[
\lambda(h,0)=l(h,h,0)=p_{12}(h,h)=(t_{1}+t_{2})(h)\quad\forall h\in D\]
\[
\lambda(0,k)=l(0,0,k)=t_{3}(k)\quad\forall k\in D.\]
Therefore, $p_{12,3}=\lambda$ and $\left((t_{1}+t_{2})+t_{3}\right)(h)=p_{12,3}(h,h)=\lambda(h,h)=l(h,h,h)$.
Analogously we can prove that $\left(t_{1}+(t_{2}+t_{3})\right)(h)=l(h,h,h)$,
that is, we get the conclusion.$\qedNoNewLine$

It is now quite easy to prove that the differential at a point is
linear
\begin{thm}
If $f:X\freccia Y$ is $\ECInfty$ smooth, the space $X$ is inf-linear
at the point $x\in X$, and the space $Y$ is inf-linear at the point
$f(x)\in Y$, then the differential\[
\diff{f}_{x}:\text{T}_{x}X\freccia\text{T}_{f(x)}Y\]
is linear.
\end{thm}
\noindent \textbf{Proof:} Let $r\in\ER$ and $t\in\text{T}_{x}X$,
we first prove homogeneity\begin{align*}
\diff{f}_{x}[r\cdot t](h) & =f\left(\left(r\cdot t\right)(h)\right)=f\left(t(r\cdot h)\right)\quad\forall h\in D.\end{align*}
On the other hand\[
\left(r\cdot\diff{f}_{x}[t]\right)(h)=\diff{f}_{x}(r\cdot h)=f\left(t(r\cdot h)\right)\quad\forall h\in D,\]
and therefore $\diff{f}_{x}[r\cdot t]=r\cdot\diff{f}_{x}[t]$.

To prove additivity, let $p$ be the infinitesimal parallelogram generated
by $t_{1}$, $t_{2}\in\text{T}_{x}X$ and $l$ the parallelogram generated
by $\diff{f}_{x}[t_{1}]$, $\diff{f}_{x}[t_{2}]\in\text{T}_{fx}Y$.
We have\begin{equation}
\diff{f}_{x}[t_{1}+t_{2}](h)=f\left(\left(t_{1}+t_{2}\right)(h)\right)=f\left(p(h,h)\right)\quad\forall h\in D.\label{eq:1_differentialIsLinear}\end{equation}
On the other hand, we obviously have\begin{equation}
\left(\diff{f}_{x}[t_{1}]+\diff{f}_{x}[t_{2}]\right)(h)=l(h,h)\quad\forall h\in D.\label{eq:2_differentialIsLinear}\end{equation}
But the smooth map\[
f\left(p(-,-)\right):D^{2}\freccia X\]
verifies\[
f\left(p(h,0)\right)=f\left(t_{1}(h)\right)=\diff{f}_{x}[t_{1}](h)\]
\[
f\left(p(0,h)\right)=f\left(t_{2}(h)\right)=\diff{f}_{x}[t_{2}](h),\]
and therefore $l=f\left(p(-,-)\right)$. From this and \eqref{eq:1_differentialIsLinear},
\eqref{eq:2_differentialIsLinear} we get the conclusion.$\qedNoNewLine$

In the case $X=\ER^{d}$ and $Y=\ER^{n}$ we have\begin{align*}
\diff{f}_{x}[t](h) & =f\left(t(h)\right)\\
 & =f\left(t(0)+h\cdot t'(0)\right)\\
 & =f\left(t(0)\right)+h\cdot t'(0)\cdot f'\left(t(0)\right)\\
 & =f(x)+h\cdot t'(0)\cdot f'(x).\end{align*}
The differential $\diff{f}_{x}[t]\in\text{T}_{f(x)}Y$ is thus uniquely
determined by the linear function $h\in D\mapsto h\cdot t'(0)\cdot f'(x)$
and hence it is uniquely determined by the vector of Fermat reals
$\iota_{2}\left[f'(x)\right]\in\ER_{2}^{n}$, as expected.

Vector fields on a generic object $X\in\ECn$ are naturally defined
as $\ECInfty$ maps of the form \[
V:X\freccia\hbox{T}X\e{such\ that}V(x)(0)=x\quad\forall x\in X.\]
In the case of manifolds, $X=\ext{M}$, this implies that $V(m)(0)\in M$
for every $m\in M$, we therefore introduce the following condition
to characterize the standard vector fields:
\begin{defn}
\label{def:tangentVectorWithStdSpeed}If $X\in\ECInfty$ is a Fermat
space and $t\in\text{T}_{x}X$ is a tangent vector at $x\in X$, then
we say that $t$\emph{ has standard speed} if and only if for every
observable $\phi\InUp{UK}X$ with $x\in U$ and $K\subseteq\ER^{\sf k}$
we have\begin{equation}
(\phi\circ t)'(0)\in\R^{\sf k}.\label{eq:1_tangentVectorWithStdSpeed}\end{equation}

\end{defn}
\noindent As usual, if $X=\ext{M}$ is a manifold, this condition
is equivalent to saying that there exists a chart on the point $\st{x}\in M$
such that condition \eqref{eq:1_tangentVectorWithStdSpeed} holds.
Using Theorem \eqref{thm: StandardTangentVector} we have the following
equivalence:\begin{equation}
\forall m\in M\pti V(m)\text{ has standard speed}\label{eq:V_hasStdSpeed}\end{equation}
if and only if \[
V|_{M}:(M\prec\ext{M})\freccia(\{\ext{f}|_{D}:\, f\in\Cn(\R,M)\}\prec{\rm T}M).\]
From this, using the definition of arrow in $\CInfty$ and the embedding
Theorem \ref{thm: EmbeddingOfSepInExtendedCn}, it follows that \eqref{eq:V_hasStdSpeed}
holds if and only if \[
V|_{M}:M\freccia{\rm T}_{{\rm st}}(M)\text{ in }\Cn,\]
that is we obtain the standard notion of vector field on $M$ because
of Theorem \ref{thm: StandardTangentVector}.

Vice versa if we have \[
W:M\freccia{\rm T}_{{\rm st}}(M)\text{ in }\Cn\]
then we can extend it to $\ext{M}$ obtaining a vector field verifying
condition \eqref{eq:V_hasStdSpeed}. In fact, for fixed $m\in\ext{M}$
and $h\in D$ we can choose a chart $(U,x)$ on $\st{m}$ and we can
write \[
W|_{U}=\sum_{i=1}^{d}A_{i}\cdot\frac{\partial}{\partial x_{i}},\]
with $A_{i}\in\CInfty(U,\R)$. But $m\in\ext{U}$ because $\st{m}\in U$
and hence we can define \[
\tilde{W}(m,h):=\sum_{i=1}^{d}\ext{A_{i}}(m)\cdot\frac{\partial}{\partial x_{i}}(m)(h)\quad\forall h\in D.\]
This definition does not depend on the chart $(U,x)$ and, because
of the sheaf property of $\ext{M}$ it provides a $\ECInfty$ function
\[
\tilde{W}:\ext{M}\times D\freccia\ext{M}\e{such\ that}\tilde{W}(m,0)=m\]
and with $(\tilde{W}^{\wedge})|_{M}=V$, that is verifying condition
\eqref{eq:V_hasStdSpeed} of standard speed.

Finally we can easily see that any vector field can be identified
equivalently with an infinitesimal transformation of the space into
itself. In fact, using cartesian closedness we have \[
V\in(X^{D})^{X}\simeq X^{X\times D}\simeq X^{D\times X}\simeq(X^{X})^{D}.\]
If $W$ corresponds to $V$ in this isomorphism then $W:D\freccia X^{X}$
and $V(x)(0)=x$ is equivalent to say that $W(0)=1_{X}$, that is
$W$ is the tangent vector at $1_{X}$ to the space of transformations
$X^{X}$, that is an infinitesimal path traced from $1_{X}$.

\section{Infinitesimal integral curves}

To the notion of vector field there is naturally associated the notion
of integral curve. In our context we are interested to define this
concept in infinitesimal terms, i.e. for curves defined on an infinitesimal
set.
\begin{defn}
\label{def:inf-integralCurve}Let $X\in\ECInfty$ be a smooth space,
$V:X\freccia\text{T}X$ a vector field on $X$ and $x\in X$ a point
in it. Then we say that $\gamma$\emph{ is the (inf-)integral curve
of $V$ at }$x$ if and only if
\begin{enumerate}
\item $\gamma:D_{\infty}\freccia X$ is smooth
\item $\gamma(0)=x$
\item $\gamma(t+h)=V\left[\gamma(t)\right](h)$ for every $t\in D_{\infty}$
and every $h\in D$.
\end{enumerate}
Moreover, we say that the vector field $V$ is \emph{inf-complete}
if and only if
\begin{enumerate}
\item \label{enu:1_defInf-Complete}$\forall x\in X\,\exists!\,\gamma_{x}\in X^{D_{\infty}}:\ \gamma_{x}$
is the integral curve of $V$ at $x$
\item \label{enu:2_defInf-Complete}The map associating to each point $x\in X$
the corresponding integral curve $\gamma_{x}$\[
x\in X\mapsto\gamma_{x}\in X^{D_{\infty}}\]
is smooth.
\end{enumerate}
\end{defn}
\begin{figure}[h]
\noindent \begin{centering}
\includegraphics[scale=0.8]{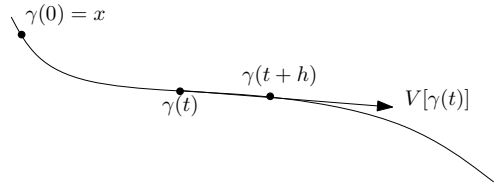}
\par\end{centering}

\caption{Explanation of the definition of integral curve}

\end{figure}

Let us note explicitly the methodological analogy among the Definition
\ref{def:inf-linear} of inf-linear space and the previous definition
of inf-complete vector field. These definitions are indeed divided
into two parts: in the first one we have that the predicate we are
defining depends on some parameter (the point $m$ in the definition
of inf-linear at $m$ and the point $x$ in the definition of integral
curve). To each value of this parameter there corresponds a unique
smooth function defined on an infinitesimal object ($D^{k}$ in the
definition of inf-linearity and $D_{\infty}$ in the definition of
integral curve). Let us note that this uniqueness is possible only
because the object is defined on an infinitesimal space. In the second
part of the definition we extend the predicate to a global object
(the whole space $X$ in the definition of inf-linearity and the vector
field in the definition of inf-completeness) universally quantifying
over these parameters, that is requiring that the first predicate
holds for every possible value of the parameters. To every universal
quantification there corresponds a further smoothness condition about
the function that to each parameter assigns the corresponding unique
infinitesimal function. The same requirement has been used in the
definition of inf-linearity at a given point $m$, where the universal
quantification is over every $k$-tuples of tangent vectors. This
method, which in some sense is implicit in SDG where every function
defined in intuitionistic logic is smooth, can be used to transpose
several definitions of SDG to our infinitesimal differential geometry.

First of all, we have to prove that the notion of inf-integral curve
generalizes, in some way, the classical notion. For simplicity let
$X=\ER^{d}$, the same reasoning can be applied to the case of more
general manifolds, because all the notions we are using are local.
So, let $V:\ER^{d}\freccia\left(\ER^{d}\right)^{D}$ be a standard
vector field, then from what we have just seen above, we know that
we can find a smooth function\[
\bar{V}:\ER^{d}\freccia\ER^{d}\]
such that\[
V(x)(h)=x+h\cdot\bar{V}(x)\quad\forall x\in\ER^{d}\ \forall h\in D.\]
Therefore, if $\gamma:D_{\infty}\freccia\ER^{d}$ is an integral curve
of $V$ at $x\in\ER^{d}$, by the Definition \ref{def:inf-integralCurve}
we get\[
\gamma(t+h)=V\left[\gamma(t)\right](h)\]
\[
\gamma(t)+h\cdot\gamma'(t)=\gamma(t)+h\cdot\bar{V}\left[\gamma(t)\right]\]
\[
h\cdot\gamma'(t)=h\cdot\bar{V}\left[\gamma(t)\right]\quad\forall t\in D_{\infty}\ \forall h\in D.\]
This implies\[
\gamma'(t)=_{2}\bar{V}\left[\gamma(t)\right]\quad\forall t\in D_{\infty}.\]
So we have the classical notion of integral curve up to second order
infinitesimals. Now, let $\eta:(-\delta,+\delta)_{\R}\freccia\R^{d}$,
$\delta\in\R_{>0}$, be a standard integral curve of $V'$, i.e.\begin{equation}
\eta'(t)=\bar{V}\left[\eta(t)\right]\quad\forall t\in(-\delta,+\delta)_{\R}.\label{eq:1_classicalAndInfNotionOfIntegralCurve}\end{equation}
Then extending $\eta$ to $(-\delta,+\delta)\subseteq\ER$ and using
the elementary transfer theorem (Theorem \ref{thm:elementaryTransferTheorem})
we obtain that the equality \eqref{eq:1_classicalAndInfNotionOfIntegralCurve}
holds also for every $t\in(-\delta,+\delta)$, and hence it holds
also in $D_{\infty}$:\[
\eta'(t)=\bar{V}\left[\eta(t)\right]\quad\forall t\in D_{\infty},\]
and thus $\eta|_{D_{\infty}}:D_{\infty}\freccia\ER^{d}$ is an inf-integral
curve of the vector field $V$. Any two of these standard integral
curves, let us say $\eta_{1}$ and $\eta_{2}$, agree in some neighborhood
$\mathcal{U}$ of $t=0$ if $\eta_{1}(0)=\eta_{2}(0)$. Therefore,
the corresponding inf-integral curves coincide on the whole $D_{\infty}\subseteq\mathcal{U}$:\[
\eta_{1}|_{D_{\infty}}=\eta_{2}|_{D_{\infty}}.\]
For this reason in Definition \ref{def:inf-integralCurve}, we say
that $\gamma$ is \emph{the }inf-integral curve of $V$ at the point
$x$.

The next step is to prove that spaces of mappings between manifolds
always verify the just introduced definition.
\begin{thm}
\label{thm:vectorFieldsOnSpacesOfMappingBetweenManifoldsAreComplete}Every
vector field $V$ in spaces of the form $X=\ext{M}$ or $X=\ext{M}^{\ext{N}}$,
where $M$ and $N$ are manifolds and where $N$ admits partitions
of unity, is inf-complete.
\end{thm}
\noindent \textbf{Proof:} The first part of the statement, i.e. the
case $X=\ext{M}$, is really a particular case of the second one where
one takes as $N=\{*\}$ any $0$-dimensional manifold. So, let us
prove only the second part of the statement. Moreover, to simplify
the notations, we will simply use the symbols $M$ and $N$ to indicate
the extensions $\ext{M}$ and $\ext{N}$.

Our vector field is a smooth map of the form\[
V:M^{N}\freccia\left(M^{N}\right)^{D}.\]
Moreover, let us consider a point $\mu\in M^{N}$. We have to prove
that there exists one and only one inf-integral curve $\gamma:D_{\infty}\freccia M^{N}$
passing from $\mu$ at $t=0$. We will construct $\gamma$ using the
sheaf property of the space $N$.

Because $N$ admits partitions of unity, we can consider a standard
open cover $\left(U_{n}\right)_{n\in N}$ of $N$ such that the closure
$\bar{U}_{n}=:K_{n}$ is compact and such that the partition of unity
$\left(\rho_{n}\right)_{n\in N}$ is subordinate to the open cover
$\left(U_{n}\right)_{n\in N}$. From cartesian closedness we can think
of $V$ as a map of the form\[
V:M^{N}\times D\times N\freccia M.\]
Using the partition of unity $\left(\rho_{n}\right)_{n\in N}$ every
smooth map $f:K_{n}\freccia M$ can be extended to a smooth map defined
on the whole $N$. Moreover, this extension, which essentially is
the multiplication by a cut-off function, can be defined as a smooth
application\[
\chi_{n}:M^{K_{n}}\freccia M^{N}\quad\forall n\in N\]
such that\[
\chi_{n}\left(f\right)|_{K_{n}}=f\quad\forall f\in M^{K_{n}}\ \forall n\in N.\]
For each $n\in N$ we can hence define\[
W_{n}:(f,h,x)\in M^{K_{n}}\times D\times K_{n}\mapsto V\left[\chi_{n}\left(f\right),h,x\right]\in M\]
obtaining a family $\left(W_{n}\right)_{n\in N}$ of smooth functions.

From cartesian closedness, these functions can be thought as\[
W_{n}:M^{K_{n}}\times D\freccia M^{K_{n}}.\]
But here $M^{K_{n}}$ is a Banach manifold because $K_{n}$ is compact,
and hence we can apply the standard local existence of integral curves
for the vector field $W_{n}$ in Banach spaces obtaining the existence
of a smooth map $\gamma_{n}:D_{\infty}\freccia M^{K_{n}}$ such that\begin{equation}
\begin{cases}
\gamma_{n}(t+h)=W_{n}\left[\gamma_{n}(t),h\right] & \forall t\in D_{\infty}\ \forall h\in D\\
\gamma_{n}(0)=\mu|_{K_{n}}\end{cases}\label{eq:initialValueProblemInBS}\end{equation}
It is not hard to prove that $\gamma_{n}$ and $\gamma_{m}$ agree
on $K_{n}\cap K_{m}$ because they verify the same initial value problem. 

From the sheaf property of the space $N$ there exist one and only
one smooth function $\gamma:D_{\infty}\times N\freccia M$ such that\[
\gamma|_{D_{\infty}\times U_{n}}=\gamma_{n}|_{U_{n}}\quad\forall n\in N.\]
From \eqref{eq:initialValueProblemInBS} it hence follows that $\gamma$
is the integral curve of $V$ at $\mu$ we searched for. Condition
\emph{\ref{enu:2_defInf-Complete}}. of the definition of inf-completeness
of $V$ follows from the classical theorem of smooth dependence from
the initial conditions (see e.g. \citet{Ab-Ma-Ra}).$\qedNoNewLine$

\section{\label{sec:IdeasForTheCalculusOfVariations}Ideas for the calculus
of variations}

In this section we want to show the flexibility of our theory proving
a very general form of the Euler-Lagrange equation. Even if the result
holds for lagrangians defined on very general spaces, the proof uses
infinitesimal methods and, when specified in the space $\ER$, is
essentially identical to the one sometimes presented in classical
courses of physics using informal infinitesimal argumentations.

We start with the notion of minimum of a functional
\begin{defn}
\label{def:minimunOFAFunctional}Let $Y\in\ECInfty$ be a Fermat space,
$\mu\in Y$ a point in it, and $J:Y\freccia\ER$ a $\ECInfty$ function.

Then we say that $J$\emph{ has a minimum at }$\mu$ if and only if\[
\forall\tau\in\text{T}_{\mu}Y\pti J\left[\tau(h)\right]\ge J(\mu)\quad\forall h\in D.\]
In other words, the value $J(\mu)$ has to be minimum along every
tangent vector of $Y$ starting from $\mu$.
\end{defn}
The first positive characteristic of our approach is that in this
definition of minimum we have used tangent vectors $\tau:D\freccia Y$
at $\mu\in Y$ instead of some notion of neighborhood of $\mu\in Y$
(like in the classical approach, see e.g. \citet{Ge-Fo}).

The total order of $\ER$ seems essential in the proof of the following
\begin{thm}
\label{thm:necessaryConditionToHaveAMinimum}Let $Y\in\ECInfty$ be
a Fermat space, $\mu\in Y$ a point in it, and $J:Y\freccia\ER$ a
$\ECInfty$ function. Moreover, let us suppose that\[
J\text{ has a minimum at }\mu.\]
Then\begin{equation}
\forall\tau\in\text{T}_{\mu}Y\pti\diff{J}_{\mu}[\tau]=\underline{0}\label{eq:LagrangeDiffZero}\end{equation}

\end{thm}
\noindent \textbf{Proof:} Firstly, let us note that\[
\xyC{40pt}\xymatrix{D\ar[r]^{\tau} & Y\ar[r]^{J} & \ER}
,\]
so that $J\circ\tau$ is smooth and we can apply the Taylor's formula
(Theorem \ref{thm: GeneralizedDerivationFormula})\begin{equation}
\forall h\in D\pti\diff{J}_{\mu}[\tau]=J\left[\tau(h)\right]=J\left[\tau(0)\right]+h\cdot(J\circ\tau)'(0)=J(\mu)+h\cdot(J\circ\tau)'(0),\label{eq:1_necessaryConditionToHaveAMinimum}\end{equation}
where $(J\circ\tau)'(0)\in\ER_{2}$ so that its order verifies\[
\omega\left[(J\circ\tau)'(0)\right]=:b>2.\]

If $\left(J\circ\tau\right)'(0)>0$, then we could set $a:=\frac{4b}{3b-2}$
and $h:=-\diff{t_{a}}$. It is easy to check that $1\le a<2$ because
$b>2$, so $h\in D_{<0}$ and we have $h\cdot(J\circ\tau)'(0)\le0$.
But\[
\frac{1}{\omega(h)}+\frac{1}{\omega\left[(J\circ\tau)'(0)\right]}=\frac{1}{a}+\frac{1}{b}=\frac{3b+2}{4b}<1\]
because $b>2$, so it is $h\cdot(J\circ\tau)'(0)<0$. But then, from
\eqref{eq:1_necessaryConditionToHaveAMinimum} we would have $J\left[\tau(h)\right]<J(\mu)$
in contradiction with the hypothesis that $J$ has a minimum at $\mu$.
Analogously, we can prove that it cannot be that $\left(J\circ\tau\right)'(0)<0$
and thus we obtain\[
\left(J\circ\tau\right)'(0)=0\]
from the trichotomy law. From \eqref{eq:1_necessaryConditionToHaveAMinimum}
it follows that $\diff{J}[\tau]=J(\mu)$, that is $\diff{J}[\tau]$
is the null tangent vector.$\qedNoNewLine$

Let us note explicitly the importance, in the previous proof, of the
possibility to construct an infinitesimal $h\in D$ having the desired
properties with respect to the order relation, e.g. $h<0$, and of
a suitable order $\omega(h)$ so that to assure that the product $h\cdot(J\circ\tau)'(0)$
is not zero.

The functionals we are interested in are of the form\begin{equation}
J(\eta)=\int_{a}^{b}L\left[t,\eta(t),\diff{\eta}_{t}\right]\diff{t}\quad\forall\eta\in X^{[a,b]}=:Y,\label{eq:functional}\end{equation}
where\[
a,b\in\R\e{with}\st{a}<\st{b}\]
\[
L:X\times\text{T}X\freccia\ER\e{in}\ECInfty\]
and where we recall that $\diff{\eta}_{t}$ is the differential of
$\eta:[a,b]\freccia X$ at the point $t\in[a,b]$, i.e. the map $\diff{\eta}_{t}:\tau\in\text{T}_{t}[a,b]\freccia\diff{\eta}_{t}[\tau]=\tau\cdot\eta\in\text{T}_{\eta(t)}X$;
moreover, we recall that $\text{T}_{y}Y=\left(\left\{ t\in\text{T}Y\,|\, t(0)=y\right\} \prec\text{T}Y\right)$
and that $\text{T}Y=Y^{D}$.

Concretely, the proof works if we can apply a Taylor's formula to
$J\left[\tau(h)\right]$ and if we also have a vector space structure
on the tangent space $\text{T}_{x}X$ (for the derivation by parts
formula), so that interesting cases are $X=\ER^{d}$ or, more generally,
any inf-linear vector space of the form $X=\ER^{\ext{M}}$, $M$ being
a generic smooth manifold (see Theorem \ref{thm:TaylorForFunctionsSpaces}
and Theorem \ref{thm:euclideanSpacesAreInf-Linear}).

Let us note that, due to cartesian closedeness of $\ECInfty$, the
notion of smooth map both for the Lagrangian $L$ and for the functional
$J$ does not present any problem even for spaces of functions like\[
X=\ext{\R}^{{\ext{M_{1}}^{\,\dots}}^{\ext{M_{s}}}}\simeq\ext{\R}^{\ext{(M_{1}\times\dots\times M_{{\scriptstyle s}})}}.\]

\noindent For these reasons, in the following we will assume\[
X=\ER^{\ext{M}}\e{,}M\text{ manifold}\]
so that our functional \eqref{eq:functional} if a map of the form
$J:Y\freccia\ER$ where%
\footnote{In the following we will use implicitly the cartesian closedeness,
without changing notation from a map to its adjoint.%
}\[
Y:=\left(\ER^{\ext{M}}\right)^{[a,b]}=\ER^{[a,b]\times\ext{M}}\]

We want to prove the Euler-Lagrange equations for a standard La\-gran\-gian
at a standard point $\mu\in Y$, so let us firstly assume that $J$
has a minimum at a standard function\[
\mu:[a,b]\freccia\ER^{\ext{M}},\]
\begin{equation}
\forall t\in(a,b)_{\R}\pti\mu(M)\subseteq\R,\label{eq:5_Lagrange}\end{equation}
(recalling Theorem \ref{thm:smoothFnctsThatTakeStndToStnd}).

Secondly, let us assume that the Lagrangian $L$ gets standard values
at $\mu$, i.e.\begin{equation}
\forall t\in(a,b)_{\R}\pti L\left[t,\mu(t),\diff{\mu}_{t}\right]\in\R.\label{eq:6_Lagrange}\end{equation}
To prove the Euler-Lagrange equations in a space of the form $X=\ER^{\ext{M}}$
(that, we recall, in general is not a Banach space) we will use infinitesimal
methods, ensuing the following thread of thoughts.

Let us start considering a tangent vector $\tau\in\text{T}_{\mu}Y$,
i.e. a function\[
\tau:D\freccia\ER^{[a,b]\times\ext{M}}.\]
Because of cartesian closedness, we can think of $\tau$ as a map
from $[a,b]\times D$ into $\ER^{\ext{M}}$. Using Taylor's formula
in the space $\ER^{\ext{M}}$ (see Theorem \ref{thm:TaylorForFunctionsSpaces})
we can write\begin{equation}
\forall t\in[a,b]\,\forall h\in D\pti\tau(t,h)=\mu(t)+h\cdot\nu(t),\label{eq:1_Lagrange}\end{equation}
where $\nu:=\tau'(0):[a,b]\freccia\ER_{\sss{=_{2}}}^{\ext{M}}$. Because
Euler-Lagrange equations are a necessary condition that follows from
\eqref{eq:LagrangeDiffZero}, let us assume that the derivative $\nu$
of our tangent vector $\tau$ is a standard smooth function, i.e.
let us assume that\[
\nu:(a,b)\freccia\ER^{\ext{M}}=X\]
\begin{equation}
\forall t\in(a,b)_{\R}\pti\nu(M)\subseteq\R\label{eq:8_Lagrange}\end{equation}
and that verifies \eqref{eq:1_Lagrange} on the open set $(a,b)$.

For a generic $h\in D$, let us calculate\begin{align*}
J\left[\tau(h)\right] & =\int_{a}^{b}L\left[t,\tau(h,t),\frac{\partial\tau}{\partial t}(h,t)\right]\diff{t}\\
 & =\int_{a}^{b}L\left[t,\mu(t)+h\cdot\nu(t),\mu'(t)+h\cdot\nu'(t)\right]\diff{t}.\end{align*}
We use the first order Taylor's formula firstly with respect to the
second variable and after with respect to the third variable (traditionally
indicated with $q$ and $\dot{q}$ respectively) obtaining\begin{align*}
J\left[\tau(h)\right] & =\int_{a}^{b}\left\{ L\left[t,\mu(t),\mu'(t)+h\cdot\nu'(t)\right]+\right.\\
 & \phantom{=\qquad\ \ }\left.+h\cdot\diff{_{2}L}\left[t,\mu(t),\mu'(t)+h\cdot\nu'(t)\right]\ldotp\nu(t)\diff{t}\right\} =\\
 & =\int_{a}^{b}\left\{ L\left[t,\mu(t),\mu'(t)\right]+h\cdot\diff{_{3}L}\left[t,\mu(t),\mu'(t)\right]\ldotp\nu'(t)+\right.\\
 & \phantom{=\qquad\ \ }+\left.h\cdot\diff{_{2}L}\left[t,\mu(t),\mu'(t)\right]\ldotp\nu(t)+h^{2}\cdot T\right\} \diff{t},\end{align*}
where we have used the notation $\diff{_{i}L}[t,q,\dot{q}]\ldotp v$
for the differential of the Lagrangian with respect to its $i$-th
argument at the point $(t,q,\dot{q})$ and applied to the tangent
vector $v$, and where $T$ is a term containing the second derivative
of $L$, but non influencing our calculation because it is multiplied
by $h^{2}=0$. Therefore, we have\begin{multline*}
J\left[\tau(h)\right]=J\left[\tau(0)\right]+h\cdot\int_{a}^{b}\left\{ \diff{_{3}L}\left[t,\mu(t),\mu'(t)\right]\ldotp\nu'(t)+\right.,\\
\left.+\diff{_{2}L}\left[\mu(t),\mu'(t)\right]\ldotp\nu(t)\diff{t}\right\} \end{multline*}
that is\[
\diff{J}_{\mu}[\tau](h)=h\cdot\int_{a}^{b}\left\{ \diff{_{3}L}\left[t,\mu(t),\mu'(t)\right]\ldotp\nu'(t)+\diff{_{2}L}\left[t,\mu(t),\mu'(t)\right]\ldotp\nu(t)\right\} \diff{t}.\]
But for Theorem \ref{thm:necessaryConditionToHaveAMinimum} we have
$\diff{J}_{\mu}[\tau]=\underline{0}$, that is\[
h\cdot\int_{a}^{b}\left\{ \diff{_{3}L}\left[t,\mu(t),\mu'(t)\right]\ldotp\nu'(t)+\diff{_{2}L}\left[t,\mu(t),\mu'(t)\right]\ldotp\nu(t)\right\} \diff{t}=0.\]
We will not delete now the factor $h\in D$ from this equation because
this would imply that the integral is equal to zero only up to second
order infinitesimals, but we will continue to take this factor for
another step, where we will use the hypothesis about the standard
nature of both functions $\mu$, $\nu$ and of the Lagrangian $L$
(see equations \eqref{eq:5_Lagrange}, \eqref{eq:6_Lagrange} and
\eqref{eq:8_Lagrange}).

Now we can apply the integration by part formula to the term\[
\diff{_{3}L}[t,\mu(t),\mu'(t)]\ldotp\nu'(t)\]
and with the bilinear form $\beta(\delta,v):=\delta\ldotp v$, where
$\delta$ is a smooth linear functional, i.e. $\delta\in\text{Lin}(\ER^{\ext{M}},\ER)$,
and where $v\in\ER^{\ext{M}}$. We obtain\begin{multline}
0=h\cdot\left[\diff{_{3}L}\left[t,\mu(t),\mu'(t)\right]\ldotp\nu(t)\right]_{a}^{b}-\\
-h\cdot\int_{a}^{b}\left\{ \diff{_{2}L}\left[t,\mu(t),\mu'(t)\right]-\frac{\diff{}}{\diff{t}}\diff{_{3}L}\left[t,\mu(t),\mu'(t)\right]\right\} \ldotp\nu(t)\diff{t}\label{eq:10_Lagrange}\end{multline}

Restricting to the case where\begin{equation}
\nu(a)=\nu(b)=\underline{0}\in\ER^{\ext{M}}\label{eq:9_Lagrange}\end{equation}
we obtain that necessarily\[
h\cdot\int_{a}^{b}\left\{ \diff{_{2}L}\left[t,\mu(t),\mu'(t)\right]-\frac{\diff{}}{\diff{t}}\diff{_{3}L}\left[t,\mu(t),\mu'(t)\right]\right\} \ldotp\nu(t)\diff{t}=0\]
 holds, because $\diff{_{3}L}[x].(-)$ is linear. In this equality
we note that the integrated function is a standard function, because
of our hypothesis \eqref{eq:5_Lagrange}, \eqref{eq:6_Lagrange} and
\eqref{eq:8_Lagrange}, so the integral itself is a standard real
and we can delete the first order infinitesimal factor $h$ because
of Theorem \ref{thm:firstCancellationLaw}:\[
\int_{a}^{b}\left\{ \diff{_{2}L}\left[t,\mu(t),\mu'(t)\right]-\frac{\diff{}}{\diff{t}}\diff{_{3}L}\left[t,\mu(t),\mu'(t)\right]\right\} \ldotp\nu(t)\diff{t}=0.\]

The usual proof of the so called fundamental lemma of calculus of
variations, which uses a continuous function for $\nu$ and not a
smooth one, can be easily substituted by a formally identical argumentation,
but with a smooth function of the form\[
\beta(\bar{t},\bar{h},\delta_{1},\delta_{2},x):=b\left(\frac{x-\bar{t}+\delta_{1}}{\delta_{1}}\right)\cdot b\left(\frac{\bar{t}+\bar{h}+\delta_{2}-x}{\delta_{2}}\right)\quad\forall x\in\R.\]
where $b\in\CInfty(\R,\R)$ is a standard smooth bump function, i.e.\begin{equation}
b(t)=0\ \ \forall t\le0\e{and}b(s)=1\ \ \forall s\ge1,\label{eq:2_Lagrange}\end{equation}
and where $(\bar{t},\bar{h},\delta_{1},\delta_{2})$ are real parameters.

\begin{figure}[h]
\noindent \begin{centering}
\includegraphics[scale=0.35]{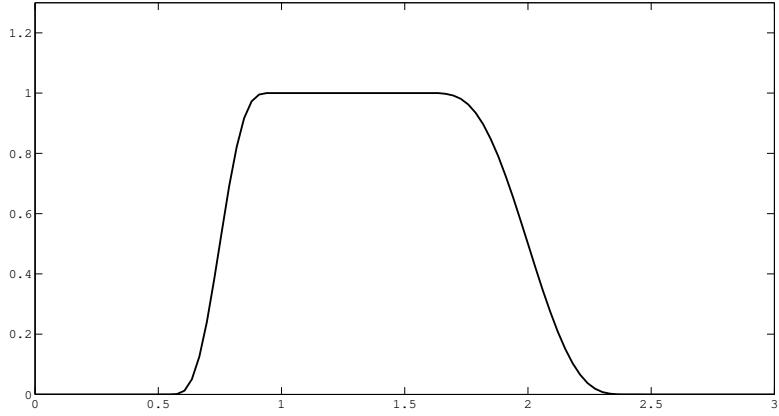}
\par\end{centering}

\caption{An example of function $\beta(\bar{t},h,\delta_{1},\delta_{2},-)$
for $\bar{t}-\delta_{1}=0.5$, $\bar{t}=1$, $\bar{t}+h=1.5$, $\bar{t}+h+\delta_{2}=2.5$.}

\end{figure}

From the smooth version of the fundamental lemma and from \eqref{eq:10_Lagrange}
we obtain the conclusion:\[
\diff{_{2}L}\left[t,\mu(t),\mu'(t)\right]=\frac{\diff{}}{\diff{t}}\diff{_{3}L}\left[t,\mu(t),\mu'(t)\right]\quad\forall\bar{t}\in(a,b)_{\R}.\]

\begin{thm}
Let $M$ be a manifold and set for simplicity $X=\ER^{\ext{M}}$.
Let us consider a smooth map\[
L:X\times\text{T}X\freccia\ER\]
and an interval $[a,b]$ with $\st{a}<\st{b}$. Let us define the
functional\begin{equation}
J(\eta)=\int_{a}^{b}L\left[t,\eta(t),\diff{\eta}_{t}\right]\diff{t}\quad\forall\eta\in X^{[a,b]}\label{eq:1_statementLagrange}\end{equation}
and assume that $J$ has a minimum at the point $\mu\in X^{[a,b]}$
such that\[
\forall t\in(a,b)_{\R}\pti\mu(M)\subseteq\R\]
\[
\forall t\in(a,b)_{\R}\pti L\left[t,\mu(t),\diff{\mu}_{t}\right]\in\R.\]
Then we have\[
\diff{_{2}L}\left[t,\mu(t),\mu'(t)\right]=\frac{\diff{}}{\diff{t}}\diff{_{3}L}\left[t,\mu(t),\mu'(t)\right]\quad\forall t\in(a,b)_{\R}.\]

\end{thm}
$\qedWithFinalEq$

We have to admit that the proof we gave of the Euler-Lagrange equation
in the space $X=\ER^{\ext{M}}$ could be elaborated further and presented
in a more clear way, e.g. clarifying better some passages, like the
identification of the tangent space $\text{T}X$ with the space of
the derivatives of the form $\mu'(r)\in\ER^{\ext{M}}$ (i.e. the identification
of the differential $\diff{\mu}_{r}[\tau](h)=\mu(r)+h\cdot\tau'(0)\cdot\mu'(r)\in\ER^{\ext{M}}$
with the element $\mu'(r)$ of the $\ER$-module $\ER^{\ext{M}}$).
However, in our opinion already in the present form it has positive
features:
\begin{enumerate}
\item The general notion of differential $\diff{J}$ of a function $J:Y\freccia Z$
between two inf-linear spaces $Y$, $Z\in\ECInfty$ can be used to
define the notion of minimum of a functional, without any need to
define norms on function spaces.
\item Functionals of the form \eqref{eq:1_statementLagrange} are smooth
even if the domain can be of the form $X^{[a,b]}$, with\[
X=\ext{\R}^{{\ext{M_{1}}^{\,\dots}}^{\ext{M_{s}}}}\]
and without any compactness hypothesis on the manifolds $M_{1},\ldots,M_{s}$.
\item The proof is formally the usual one used in the situation where $X=\R^{d}$,
but our smooth framework is more appropriate, e.g. because of cartesian
closedness and completeness and co-completeness.
\end{enumerate}

\chapter{\label{cha:FurtherDevelopments}Further developments}

Several ideas can be developed starting from this foundations of the
theory of Fermat reals we provided in the present work. Some are systematic,
with high feasibility; some other are, at the present stage, only
sketches of ideas. In the next sections we should present some of
them, with no aim to be exhaustive in their presentation.

\section{\label{sec:FirstOrderInfinitesimalsWhoseProductIsNotZero}First order
infinitesimals whose product is not zero}

We have seen (see Theorem \ref{thm:orderImpliesProductOfFirstOrderInfinitesimalsIsZero}
and the related discussion) that it is impossible to have good properties
for the order relation of the ground ring and at the same time to
have the existence of two first order infinitesimals whose product
is not zero. On the other hand we have had to develop the notion of
equality up to $k$-th order infinitesimals (Chapter \ref{cha:equalityUpTo_k-thOrder})
to bypass this algebraic problem, first of all in connection with
its relationships with Taylor's formula for functions defined on infinitesimal
domains (Section \ref{cha:CalculusOnInfinitesimalDomains}). In the
present work we have seen that a total order can be very useful. For
example, our geometrical representation of Fermat reals is strongly
based on the trichotomy law, and we have also seen that the possibility
to have a total order can be very useful in some proofs (see Section
\ref{sec:IdeasForTheCalculusOfVariations}). On the other hand, the
possibility to have two first order infinitesimals whose product is
not zero, opens, like in SDG, the possibility to prove a general cancellation
law of the form\[
\left(\forall h\in D:\ h\cdot m=h\cdot n\right)\then m=n,\]
and hence to avoid the use of the equality up to a $k$-th order infinitesimal.

The ideal solution would be to keep all the results we have shown
in the present work and, at the same time, to have the possibility
to consider pairs of first order infinitesimals whose product is not
necessarily zero. An idea, inspired by rings like\[
\R[t,s]/\langle t^{2}=0,s^{2}=0\rangle,\]
we can try to explore, can be roughly stated saying that {}``two
first order infinitesimals $(h_{t})_{t}$ and $(k_{s})_{s}$ have
a non zero product $(h_{t}\cdot k_{s})_{t,s}$ if they depend on two
independent variables $t$ and $s$''. A possible formalization of
this idea can be sketched in the following way.

Firstly let us fix a way to embed a space of type $\R^{n}$ into $\R^{m}$
if $n<m$, e.g.\begin{equation}
(r_{1},\ldots,r_{n})\in\R^{n}\mapsto(r_{1},\ldots,r_{n},0,\ptind^{m-n},0)\in\R^{m}.\label{eq:1_firstOrderInfWitNoZeroProduct}\end{equation}
Then, instead of little-oh polynomials $x:\R_{\ge0}\freccia\R$, let
us consider maps of the form\[
y:r\in\R_{\ge0}^{n}\mapsto x\left(t_{i}^{n}(r)\right)\in\R,\]
where $x\in\R_{0}[t]$ is a usual little-oh polynomial and where $t_{i}^{n}:r\in\R_{\ge0}^{n}\mapsto r_{i}\in\R_{\ge0}$
is the projection onto the $i$-th component. In this case we say
that $y$ \emph{depends on the variable} $t_{i}^{n}$ or, where there
is no confusion, simply on the variable $t_{i}$. Therefore, our map
$y$ can now be written as\[
x_{t_{i}}=r+\sum\limits _{j=1}^{k}\alpha_{j}\cdot t_{i}^{a_{j}}+o(t_{i})\quad\text{ as }\quad t_{i}\to0^{+},\]
where the limit has to be understood along the directed set\[
(\R^{n},\le)\]
\[
(r_{1},\ldots,r_{n})\le(s_{1},\ldots,s_{n})\DIff r_{i}\le s_{i}.\]
But if we sum this map $y$ with a map $z$ that depends on the variable
$t_{j}^{m}$, what do we obtain? Intuitively, a map which is a function
of the two variables $t_{i}^{m}$ and $t_{j}^{m}$ if we firstly embed
$\R^{n}$ into $\R^{m}$ using \eqref{eq:1_firstOrderInfWitNoZeroProduct}.
Thus, more generally, we have to consider maps of the form\begin{equation}
x_{t_{i_{1}}\cdots t_{i_{v}}}=r+\sum_{j=1}^{k}\alpha_{j}\cdot t_{i_{1}}^{a_{1j}}\cdot\ldots\cdot t_{i_{v}}^{a_{vj}}+o(t_{i_{1}})+\ldots+o(t_{i_{v}}),\label{eq:2_firstOrderInfWithNoZeroProduct}\end{equation}
$t_{i_{1}},\ldots,t_{i_{v}}$ being all the variables from which the
map $x$ depends on. In \eqref{eq:2_firstOrderInfWithNoZeroProduct}
the limit has to be mean along the directed set\[
(\R^{m},\le)\]
\[
(r_{1},\ldots,r_{m})\le(s_{1},\ldots,s_{m})\DIff r_{i_{1}}\le s_{i_{1}},\ldots,r_{i_{v}}\le s_{i_{v}}.\]
More precisely, with a writing like\[
\mathcal{P}(o(\varphi_{1}),\ldots,o(\varphi_{n}))\]
where $\mathcal{P}$ is a generic property and $\varphi_{1},\ldots,\varphi_{n}$
are free variables in $\mathcal{P}$ for functions in the space $\R^{\R_{\ge0}^{m}}$,
we mean\[
\exists w_{1},\ldots,w_{n}\in\R^{\R_{\ge0}^{m}}\pti\begin{cases}
\forall^{0}t\pti\mathcal{P}\left(w_{1}(t),\ldots,w_{n}(t)\right)\\
w_{i}=o(\varphi_{i})\quad\forall i=1,\ldots,n\end{cases}\]
And an {}``equality'' of the type $w=o(\varphi)$, as usual in our
context, means\[
\exists\lim_{t\to0^{+}}\frac{w(t)}{\varphi(t)}\in\R\text{ and }w(0)=0.\]
The analogue of the equality in $\ER$ (i.e. the equivalence relation
introduced in Definition \ref{def:equalityInFermatReals}) is now
that $x\sim y$ if and only if\[
x_{t_{i_{1}}\ldots t_{i_{v}}}=y_{t_{i_{1}}\ldots t_{i_{v}}}+o(t_{i_{1}})+\ldots+o(t_{i_{v}})\text{ as }t_{i_{k}}\to0^{+}\ \forall k,\]
where $t_{i_{1}},\ldots,t_{i_{v}}$ are all the variables from which
the maps $x$ and $y$ depend on.

This idea seems positive for two reasons: firstly, if we define a
new Fermat reals ring in this way, considering only the subring of
all the maps $\R_{o}[t_{i}]$ which only depend on one variable $t_{i}$,
we obtain a ring $\ER[t_{i}]$ isomorphic to the present $\ER$. This
means that we are not loosing all the results we have proved in the
present work.

Secondly, let us consider $h_{t_{1}}:=t_{1}$ and $k_{t_{2}}:=t_{2}$,
then we have that $h^{2}\sim0$ and $k^{2}\sim0$, but if we were
to have $h\cdot k\sim0$, then we would get\[
t_{1}\cdot t_{2}=o(t_{1})+o(t_{2})\]
\[
t_{2}=\frac{o(t_{1})}{t_{1}}+\frac{o(t_{2})}{t_{1}}.\]
But the left hand side of this equality goes to zero for $t_{1}\to0^{+}$
and $t_{2}\to0^{+}$, whereas on the right hand side the limit\[
\lim_{\substack{t_{1}\to0^{+}\\
t_{2}\to0^{+}}
}\frac{o(t_{2})}{t_{1}}\]
does not exist. We therefore have indeed an example of two first order
infinitesimals whose product is not zero.

Of course, from Theorem \ref{thm:orderImpliesProductOfFirstOrderInfinitesimalsIsZero}
it follows that every subring $\ER[t_{i}]$ is totally ordered, but
the whole ring cannot be totally ordered.

\section{Relationships with Topos theory}

It is possible to define a meaningful notion of powerset diffeology
(see \citet{Igl}) defined on the powerset $\mathcal{P}(X)$ of any
diffeological space $X$. Let us recall that any diffeological space
is also a $\CInfty$ space. Therefore, we can try to see whether there
is some relation between this powerset diffeology and the powerset
object as defined in Topos theory. In case of a positive answer, this
would imply that our category $\CInfty$ is a Topos. It would start
thus the possibility to consider its internal language, almost surely
in intuitionistic logic, to describe the objects of $\CInfty$. Independently
from the results related to the powerset diffeology, we can try to
see whether an axiomatic approach to $\CInfty$ and $\ECInfty$ can
be developed. This could be useful for those readers who are interested
in the study of infinitesimal differential geometry without being
forced to consider the whole construction of $\CInfty$ and $\ECInfty$.
Almost surely, this axiomatic description can be introduced in classical
logic instead of intuitionistic logic. Indeed, such an axiomatic description
can be developed without considering the internal language of a Topos,
and hence without assuming that all our objects and maps are smooth.
Of course, we would need an axiom that permits to construct a general
family of smooth functions starting from smooth functions, and a {}``starting
point'' for this construction, like the assumption that all the standard
smooth functions are arrows of the category $\CInfty$.

\section{A transfer theorem for sentences}

We have seen the proof of a transfer theorem for the construction
of $\ECInfty$ spaces using logical formulas and the preservation
properties of the Fermat functor $\ext{(-)}$ (see Chapters \ref{cha:TheFermatFunctor}
and \ref{cha:logicalPropertiesOfTheFermatFunctor}). As already stated
at the end of Chapter \ref{cha:logicalPropertiesOfTheFermatFunctor},
differently from our situation, the transfer theorem of NSA asserts
an equivalence between two sentences. Nevertheless, it seems possible
to follow the following scheme:
\begin{enumerate}
\item Define the meaning of the sentence {}``\emph{the formula $\xi$ is
intuitionistically true in $\CInfty$}'' using the intuitionistic
interpretation of the propositional connectives and quantifiers in
this category. An analogous definition of intuitionistic validity
can be done in the category $\ECInfty$.
\item Define the $\ext{(-)}$-transform of a given formula $\xi$.
\item Prove that $\xi$ is intuitionistically true in $\CInfty$ if and
only if $\ext{\xi}$ is intuitionistically true in $\ECInfty$.
\end{enumerate}

\section{Two general theorems for two very used techniques}

We used several times two techniques in our proofs. The first one
is usually a way to speed up several proofs saying {}``the considered
function is smooth because it can be expressed as a composition of
smooth functions''. Among these functions we have also to consider
set theoretical operations like those listed in Section \ref{sec:ExamplesOfCnSpacesAndFunctions}
or those related to cartesian closedness. It would be useful to define
generally which logical terms can be obtained in this way and to prove
a general theorem that roughly states that every function given by
local formulas {}``smooth in each variable'' is indeed smooth in
$\ECInfty$. This theorem can be assumed as an axiom in the above
mentioned axiomatic description of $\CInfty$ and $\ECInfty$ and
it would substitute very well the global hypotheses to work in intuitionistic
logic (where every function can be assumed to be smooth). In other
words, instead of saying: {}``because we are working in a Topos of
smooth spaces and functions, every space and every functions we will
define using intuitionistic logic is smooth'', we can say: {}``because
of the categorical property of our categories and because the considered
function $f$ is locally smooth in every variable, the considered
space $X$ and the function $f$ are smooth''.

Another very useful technique we have used is based on the local form
of figures of Fermat spaces (see Theorem \ref{thm:figuresOfFermatSpaces}).
It would be useful, even if it seems not easy to find the corresponding
statement, to prove a general theorem that permits to transfer a property
that is {}``locally true and valid for smooth function of the form
$\ext{\alpha}(p,-)$'' to a property that is {}``globally true''
for function that are locally of the form $\ext{\alpha}(p,-)$.

\section{Infinitesimal differential geometry}

After a verification of the idea presented in Section \ref{sec:FirstOrderInfinitesimalsWhoseProductIsNotZero},
it would be natural to present a development of infinitesimal differential
geometry along the lines already presented in SDG (see e.g. \citet{Lav}).
As we have already said several times, frequently the proofs and the
definitions given in SDG can be easily reformulated in the context
of Fermat spaces, so that the development of this idea sometimes coincides
with the formal repetition in our context of those proofs. On the
other hand, the property that the product of two first order infinitesimal
is always zero, which is one of the most important differences between
our theory and SDG, forces us to find a completely new thread of ideas.
In contrast to SDG, in our context the study of the relationships
between classical results on manifolds and our infinitesimal version
is usually a not hard task, whereas in SDG these relationships must
always pass through the construction of a suitable topos and a corresponding
non trivial embedding of a class of standard smooth manifolds (see
Section \ref{sec:SDG} and e.g. \citet{Mo-Re} for more details).

\section{Automatic differentiation }

Like in the Levi-Civita field (see Section \ref{sec:LeviCivitaField})
using Fermat reals we have all the instrument to try a computer implementation
of algorithms for automatic differentiation. Even if in the present
work we have concentrated ourselves in developing a {}``smooth framework'',
it is not hard to prove the following result
\begin{thm}
\label{thm:TaylorForCnFunctions}Let $f:\R\freccia\R$, $x\in\R$
and $n\in\N_{>0}$, then\[
f\text{ is of class }\Cn\text{ at the point }x\]
if and only if the following conditions are verified
\begin{enumerate}
\item $f$ is locally Lipschitz in a neighborhood of $x$
\item $\exists\, m_{1},\ldots,m_{n}\in\R\ \forall h\in D_{n}\pti f(x+h)=\sum\limits _{j=0}^{n}\frac{h^{j}}{j!}\cdot m_{j}$
\end{enumerate}
\end{thm}
This permits to reproduce in the context of Fermat reals several applications
of the Levi-Civita field in the frame of automatic differentiation
theory (see \citet{Sha,Berz-et-al,Berz1} and Section \ref{sec:LeviCivitaField}).

\section{Calculus of variations}

We sketched in Section \ref{sec:IdeasForTheCalculusOfVariations}
some ideas that our framework can give in the context of the calculus
of variations. In SDG this topic has been approached in \citep{Bu-He}
and \citep{Nis}. It is thus natural to try to reformulate in our
context these results and in general to study whether the possibility
to consider exponential spaces in the category $\ECInfty$ can lead
to some more general results, or at least to have a more natural approach
to some classical results. Indeed, we have shown that the use of infinitesimal
methods can be useful both to define well-known notion of calculus
of variations without being forced to introduce a norm, and hence
without assuming a corresponding compactness hypotheses. On the other
hand, we have also shown that these infinitesimal methods can also
be very useful to generalize in spaces of mappings the Euler-Lagrange
equations. What other results are generalizable in this type of spaces?
What other notions can be defined using tangent vectors like in Definition
\ref{def:minimunOFAFunctional} without considering a neighborhood
generated by a norm instead?

\section{Infinitesimal calculus with distributions}

In the present work, every space and function we have considered is
smooth. This can be useful in a context like infinitesimal differential
geometry, but it is obviously a limitation if one needs to apply infinitesimal
methods in contexts with non smooth functions. A possibility is to
extend the theory developing an infinitesimal calculus for distributions.
Definitions in our framework of the space of all the distributions
given by families of smooth functions with a suitable equivalence
relation (like in \citep{An-Mi-Si} or in \citep{Col2}, where non
linear polynomial operations on distributions can also be considered)
are the most promising ones for this type of generalization of the
Fermat reals to a non smooth context.

\section{Stochastic infinitesimals}

Let $(\Omega,\mathcal{F},P)$ be a probability space and let us consider
stochastic processes indexed by $t\in\R_{\ge0}$. With the symbol
${\displaystyle \mathop{\freccia}^{P}}$ we will denote the convergence
in probability.

Using the notion of little-oh relation for stochastic processes, i.e.\[
X_{t}=o_{P}\left(Y_{t}\right)\DIff\exists\left(Q_{t}\right)_{t}\text{ stochastic process}\pti\begin{cases}
X_{t}=Y_{t}\cdot Q_{t}\\
Q_{t}\,{\displaystyle \mathop{\freccia}^{P}}\,0\end{cases}\]
we can try to consider suitable stochastic processes $(X_{t})_{t\ge0}$
instead of little-oh polynomials used in the present definition of
$\ER$. What will we obtain in this way? Does the corresponding algebraic
structure permit to prove in a rigorous and formally identical way
informal equalities like\[
\diff{B}(t)=\sqrt{\diff{t}}\]
for a Brownian motion $B$? Let us note that the square root is not
smooth at the origin, and hence the term $\sqrt{\diff{t}}$ has to
be understood in a suitable way. For example we can denote by $\sqrt{k}$,
for $k\in D_{\infty}$ infinitesimal, the simplest $h\in D_{2\omega(k)}$
such that $h^{2}=k$. Here {}``simplest'' means that the decomposition
of $h$ does not contains first order infinitesimals, i.e. $h\in\ER_{1}$.
Let us immediately note that this notion of square root does not verifies
$\sqrt{h^{2}}=|h|$. E.g. if $h^{2}=0$, then $\sqrt{h^{2}}=0$ (the
simplest number whose square is zero is the zero itself), whereas
it can be $|h|\ne0$, so $\sqrt{h^{2}}\ne|h|$.

In this context it is possible to conceive the possibility to develop
a differential geometry extending a manifold using such stochastic
infinitesimals.

From the point of view of cartesian closedness, this possibility is
tied with the one of defining interesting probability measurea on
the space $\Omega_{1}^{\Omega_{2}}$ of measurable mappings between
two given probability spaces, without any particular assumption about
the topology%
\footnote{Usually this problem is solved in the case of complete separable metric
spaces or locally compact linear spaces or in case of spaces which
are Borel-isomorphic to a Borel subset of $\R$. Recall that, due
to cartesian closedness, each one of the spaces $\Omega_{i}$ can
be itself an exponential object of the form $\Omega_{3}^{\Omega_{4}}$
and and so on, so that it is not natural to make strong topological
assumptions (see Chapter \ref{cha:ApproachesToDiffGeomOfInfDimSpaces}).%
} of the spaces $\Omega_{i}$. For this a combination of ideas of integrals
in infinite dimensional spaces (see e.g. \citet{Schw} and \citet{Ito}
and references therein) and generalized Riemann integral (see e.g.
\citep{Mul1,Mul2,Ku-Sw} and references therein) could be useful.

\section{Infinite numbers and nilpotent infinitesimal}

In every field the property $h^{2}=0$ implies $h=0$, so its seems
impossible to make infinities and nilpotent infinitesimals to coexist.
But with the usual properties, also the existence of the square root
would be incompatible with the existence of non zero nilpotent infinitesimals,
but we have just seen in fact that some meaningful notion of square
root is indeed possible. Of course, not all the usual property of
this square root can be maintained in the extension from the real
field $\R$ to the Fermat ring $\ER$. On the other hand, infinities
and nilpotent infinitesimals do coexist in standard analysis, and
in our theory we have a good dialectic between potential infinitesimals
and actual infinitesimals in $\ER$. These are the motivations to
try to make coexist these two types of extended numbers in the same
structure. The problem is what property cannot be extended from $\R$
to $\ER$? Is the corresponding formalism sufficiently natural to
work with? Let us present some more concrete ideas in this direction.

If we wish to introduce infinities in the ring $\ER$, we will have
the problem of the meaning of products of the form $h\cdot H$, where
$h$ is an infinitesimals and $H$ is an infinite number. But, unlike
NSA where the solution is only formal in case of non convergent sequences
$(h_{n}\cdot H_{n})_{n\in\N}$, here we want to follow the way used
in standard analysis: {}``a product of the form $0\cdot\infty$ can
be anything: $0$, $\infty$, $r\ne0$ or nothing in case it does
not converge''. Based on this informal motivation, we can understand
that the property we have to criticize is\[
\left(x=x'\e{and}y=y'\right)\then x\cdot y=x'\cdot y',\]
because, if we want to have infinitesimals and infinities in $\ER$,
we cannot multiply freely two numbers in this ring and to hope to
always obtain a meaningful result. E.g. we can try to obtain sufficient
conditions of the form: we can multiply $x$ and $y$ in case both
are finite or if $x-x'$ goes to zero more quickly than the order
with which $y$ goes to infinite and vice versa''. E.g. if we define\[
x=_{n}y\DIff x_{t}=y_{t}+o\left(t^{n}\right)\e{as}t\to0^{+},\]
then we have not one equality only, but a family of equalities, one
for each $n\in\N_{>0}$. On the one hand, this is positive because
the subring of finite numbers with the equality $=_{1}$ is exactly
the present ring of Fermat reals. On the other hand we can prove the
following:
\begin{thm}
Let $x$, $y$, $x'$, $y':\R_{\ge0}\freccia\R$ be maps, and let
$p$, $m$, $n\in\N_{>0}$, $a$, $b\in\R$ be numbers that verify\[
p\le n\e{and}p\le m\]
\begin{equation}
\forall^{0}t\pti\left|x_{t}\cdot t^{m-p}\right|<a\e{and}\left|y'_{t}\cdot t^{n-p}\right|<b.\label{eq:conditionForProductWithInfinities}\end{equation}
Then we have\[
x=_{n}x'\e{and}y=_{m}y'\then x\cdot y=_{p}x'\cdot y'\]

\end{thm}
\textbf{Proof:} We have that\begin{align*}
xy-x'y' & =xy-xy'+xy'-x'y'\\
 & =x\cdot(y-y')+(x-x')\cdot y'\\
 & =x\cdot o_{1}\left(t^{m}\right)+o_{2}\left(t^{n}\right)y'.\end{align*}
sBut\[
\frac{x_{t}\cdot o_{1}\left(t^{m}\right)}{t^{p}}=\frac{x_{t}\cdot t^{m-p}\cdot o_{1}\left(t^{m}\right)}{t^{m}}\to0\]
because, by hypotheses, $\left|x_{t}\cdot t^{m-p}\right|$ is bounded
from above. Analogously we can deal with the term $o_{2}\left(t^{n}\right)\cdot y'$
and hence we have the conclusion.$\qedNoNewLine$

Condition \eqref{eq:conditionForProductWithInfinities} says that
the numbers $x$ and $y$ cannot be infinities {}``too large'',
and hence includes the case where both $x$ and $y$ are finite. But
if we have $x=_{1}y$ and $z$ infinite, then $z=_{m}z$ for every
$m\in\N_{>0}$, but $\left|z_{t}\cdot t^{1-1}\right|=|z_{t}|$ which
is unbounded and hence we cannot use the previous theorem to deduce
that $x\cdot z=y\cdot z$. In other words, in this structure we cannot
multiply an equality of the form $=_{1}$ with an infinite number.
This make it possible to have the coexistence of $h^{2}=_{1}0$ with
the existence of the inverse of the nilpotent $h$, i.e. a number
$k$ such that $k\cdot h=_{m}1$ for every $m\in\N_{>0}$.

As mentioned above, the feasibility of this simple idea is tied with
the possibility to create a sufficiently flexible formalism to deal
with nilpotent infinitesimals and infinite numbers at the same time,
using a family of equalities $=_{m}$. The first aims to test this
construction are of course tied with the possibility to describe Riemann
integral sums using our infinitesimals and infinities and the possibility
to define at least some $\delta$ Dirac like distributions.

\begin{appendix}
\renewcommand{\chaptermark}[1]{\markboth {\appendixname\ \thechapter. \ #1}{}}

\part{Appendices}

\chapter{\label{app:someNotionsOfCategoryTheory}Some notions of category
theory}

This appendix recalls those (more or less) standard definitions and
basic results which are used in the present work. It also aims at
clarifying the notations of category theory we use in this work, but
it is not meant as an introduction to the subject. For this reason,
no proofs and no intuitive interpretations, nor a sufficient amount
of examples, are given; they can be found in several standard textbooks
on category theory (see e.g. \citet{Ad-He-St,Ar-Ma,Mac}).

All the definitions and theorems we will state are framed in the set
theory $\text{\textbf{NBG}}$ of von Neumann-Bernay-G$\text{\textipa{ö}}$del,
where, in some cases, we can add the axiom about the existence of
Grothendieck universes.

\section{Categories}
\begin{defn}
\label{def:category}A category $\mathbf{C}$ is a structure of the
form\[
\mathbf{C}=\left((-)\xfreccia{(-)}(-),1_{(-)},\cdot,\mathcal{O},\mathcal{A}\right),\]
where $\mathcal{O}$ and $\mathcal{A}$ are classes, called respectively
the class of \emph{objects} and the class of \emph{arrows or morphisms}
of $\mathbf{C}$. The relation\[
(-)\xfreccia{(-)}(-)\subseteq\mathcal{O}\times\mathcal{A}\times\mathcal{O}\]
is called the \emph{arrow relation of $\mathbf{C}$}. The function\[
1_{(-)}:\mathcal{O}\freccia\mathcal{A}\]
assigns an arrow $1_{A}$, called the \emph{identity} of the object
$A$. Finally the function\[
\cdot:\left\{ (f,g)\in\mathcal{A}\times\mathcal{A}\,|\,\exists\, f\cdot g\right\} \freccia\mathcal{A}\]
is called the \emph{composition of the arrows $f$ and $g$ with respect
to the objects $A$, $B$ and $C$}. The predicate $\exists\, f\cdot g$
will be defined by the following conditions, which hold for every
object $A$, $B$, $C$, $D\in\mathcal{O}$ and every arrow $f$,
$g$, $h\in\mathcal{A}$:
\begin{enumerate}
\item $\exists\, f\cdot g\iff\exists\, A,B,C\pti A\xfrecciad{f}B\xfrecciad{g}C$,
i.e. the composition $f\cdot g$ is defined if the arrow $f$ takes
some object $A$ into $B$ and the arrow $g$ takes $B$ into $C$.
\item $A\xfrecciad{f}B\xfrecciad{g}C\then A\xfrecciad{f\cdot g}C$
\item $A\xfrecciad{f}B\xfrecciad{g}C\xfrecciad{h}D\then f\cdot(g\cdot h)=(f\cdot g)\cdot h$,
i.e. the composition is associative
\item $A\in\mathcal{O}\then A\xfrecciad{1_{A}}A$
\item $A\xfrecciad{f}B\xfrecciad{g}C\then f\cdot1_{B}=f\e{and}1_{B}\cdot g=g$
\end{enumerate}
\end{defn}
The following notation\[
\mathbf{C}(A,B):=\left\{ f\in\mathcal{A}\,|\, A\xfreccia{f}B\right\} \]
is also very used. Let us note that generally speaking $\mathbf{C}(A,B)$,
usually called \emph{hom-set}, is a class and not a set. In case it
is a set and not a proper class, then the category is said to be \emph{locally
small.} If the classes of objects and that of arrows of a category
are sets and not proper classes, then the category is called \emph{small}.

It is possible to prove that for every object $A$ of $\mathbf{C}$
there exists one and only one arrow $u$ such that\[
A\xfrecciad{u}A\]
\[
B\xfrecciad{f}A\xfrecciad{g}C\then f\cdot u=f\e{and}u\cdot g=g.\]
From the definition of category, this arrow is $u=1_{A}$, so that
the notion of identity can be defined starting from the arrow relation
and the composition map. For this reason, in defining a category we
have no need to specify the definition of the identity $1_{A}$.

In the present work, unless it is differently specified, we will not
assume that if $A\xfrecciad{f}B$, then the objects $A$ and $B$
are uniquely determined by the arrow $f$. On the contrary, if the
property\[
A\xfrecciad{f}B\e{and}A'\xfrecciad{f}B'\then A=A'\e{and}B=B'\]
holds, then we say that the category $\mathbf{C}$ \emph{has domains
and codomains }and we can define the domain and codomain maps:\[
\text{dom}:\mathcal{A}\freccia\mathcal{O}\]
\[
\text{cod}:\mathcal{A}\freccia\mathcal{O}\]
\[
\text{dom}(f)\xfrecciad{f}\text{cod}(f)\quad\forall f\in\mathcal{A}.\]

In case we have to consider more than one category, we will use notations
like\[
\mathbf{C}\vDash A\xfrecciad{f}B\]
\[
\mathbf{C}\vDash f\cdot g=h.\]
Moreover, we will also use the notations\[
\text{Obj}(\mathbf{C}):=\mathcal{O\e{,}\text{Arr}(\mathbf{C}):=\mathcal{A}}\]
\[
A\in\mathbf{C}\DIff A\in\text{Obj}(\mathbf{C})\]
In almost all the examples of categories considered in the present
work, the objects are sets with some additional structure and the
morphisms are maps between the underlying sets that preserve this
structure. So we have the category $\Set$ of all sets, the category
$\text{\textbf{Grp}}$ of all groups, the category $\ManInfty$ of
smooth manifolds, etc. Let us note that every set is a category with
only identity arrows.

An example used in this work that is not a category of sets with a
structure is given by the category corresponding to a preorder. Indeed,
let $(P,\le)$ be a preordered set; let us fix any element $*\in\Set$
(it is not important what element concretely we choose, e.g. it can
be $*=0\in\R$) and define\[
\mathcal{O}:=P\e{,}\mathcal{A}:=\{*\}\]
\[
x\xfreccia{f}y\DIff x,y\in P\e{,}x\le y\e{,}f=*\]
\[
x\xfreccia{f}y\xfreccia{g}z\then f\cdot g:=*.\]
It is easy to prove that in this way we obtain a category.
\begin{defn}
\label{def:duality-subcategory}Let $\mathbf{C}$ be a category, then
$\mathbf{C}^{\text{op}}$ is the category obtained {}``reversing
the direction of all the arrows'', i.e.\[
\text{\emph{Obj}}\left(\mathbf{C}^{\text{op}}\right):=\text{\emph{Obj}}(\mathbf{C})\e{,}\text{\emph{Arr}}\left(\mathbf{C}^{\text{op}}\right):=\text{\emph{Arr}}(\mathbf{C})\]
\[
\mathbf{C}^{\text{\emph{op}}}\vDash A\xfreccia{f}B\DIff\mathbf{C}\vDash B\xfreccia{f}A\]
\[
\mathbf{C}\vDash f\cdot g=h\then\mathbf{C}^{\text{\emph{op}}}\vDash g\cdot f:=h\]
Moreover, if $\mathbf{D}$ is another category, we say that $\mathbf{D}$
\emph{is a subcategory of }$\mathbf{C}$ if and only if the following
conditions hold:
\begin{enumerate}
\item $\text{\emph{Obj}}(\mathbf{D})\subseteq\text{\emph{Obj}}(\mathbf{C})$
\item $\mathbf{D}\vDash A\xfreccia{f}B\then\mathbf{C}\vDash A\xfreccia{f}B\quad\forall A,B,f$
\item $\mathbf{D}\vDash f\cdot g=h\then\mathbf{C}\vDash f\cdot g=h.$
\end{enumerate}
\end{defn}
For two categories, the \emph{product category} is defined by\begin{align*}
\text{Obj}(\mathbf{C}\times\mathbf{D}) & :=\text{Obj}(\mathbf{C})\times\text{Obj}(\mathbf{D})\\
\text{Arr}(\mathbf{C}\times\mathbf{D}) & :=\text{Arr}(\mathbf{C})\times\text{Arr}(\mathbf{D})\end{align*}
\[
\mathbf{C}\times\mathbf{D}\vDash(C_{1},D_{1})\xfreccia{(f,g)}(C_{2},D_{2})\DIff\begin{cases}
\mathbf{C}\vDash C_{1}\xfreccia{f}C_{2}\\
\mathbf{D}\vDash D_{1}\xfreccia{g}D_{2}\end{cases}\]
\[
\mathbf{C}\times\mathbf{D}\vDash(c,d)\cdot(\gamma,\delta)=(f,g)\DIff\begin{cases}
\mathbf{C}\vDash c\cdot\gamma=f\\
\mathbf{D}\vDash d\cdot\delta=g\end{cases}\]

In Chapter \ref{cha:theCartesianClosure} we mention at the notion
of \emph{Grothendieck universe}, which is defined as follows.
\begin{defn}
\label{def:GrothedieckUniverse}We say that the class $\mathcal{U}$
is a \emph{Grothendieck universe} if and only if the following conditions
hold:
\begin{enumerate}
\item $x\in\mathcal{U}\e{and}y\in x\then y\in\mathcal{U}$
\item $x,y\in\mathcal{U}\then\{x,y\}\in\mathcal{U}\e{and}(x,y)\in\mathcal{U}$
\item $x\in\mathcal{U}\then\left\{ y\,|\, y\subseteq x\right\} \in\mathcal{U}$
\item If $(x_{i})_{i\in I}$ is a family of elements of $\mathcal{U}$ and
if $I\in\mathcal{U}$, then $\bigcup_{i\in I}x_{i}\in\mathcal{U}$
\item $x$, $y\in\mathcal{U}$ and $f:x\freccia y$ is a map between these
sets, then $f\in\mathcal{U}$
\item $\N\in\mathcal{U}$, i.e. the set of natural numbers belongs to the
universe $\mathcal{U}$
\end{enumerate}
\end{defn}
In other words in a Grothendieck universe all the usual constructions
of set theory are possible. A supplementary axiom of set theory that
one may need when using category theory is\begin{equation}
\forall x\,\exists\,\mathcal{U}\text{ Grothendieck universe}\pti x\in\mathcal{U},\label{eq:GrothendieckAxiom}\end{equation}
that is, every class is an element of a suitable universe. The theory
$\text{\textbf{NBG}}$ changes radically if we assume this axiom.
E.g. all our categories can be defined in a given fixed universe (obtaining
in this way classes of that universe), but if we need to consider
e.g. $\ManInfty_{\sss{\,\mathcal{U}}}$ as an element of another class,
then we can consider another Grothendieck universe $\mathcal{U}_{2}$
that contains $\ManInfty_{\sss{\,\mathcal{U}}}$ as an element. In
this way $\ManInfty_{\sss{\,\mathcal{U}}}$ is now a set, and not
a proper class, in the universe $\mathcal{U}_{2}$. All our construction
do not depend on the axiom \eqref{eq:GrothendieckAxiom}.

\section{Functors}
\begin{defn}
\label{def:Functor}Let $\mathbf{C}$ and $\mathbf{D}$ be two categories,
then a functor\[
F:\mathbf{C}\freccia\mathbf{D}\]
is a pair $F=(F_{\text{o}},F_{\text{a}})$ of maps\begin{align*}
F_{\text{o}} & :\text{\emph{Obj}}(\mathbf{C})\freccia\text{\emph{Obj}}(\mathbf{D})\\
F_{\text{a}} & :\mathcal{AR}(\mathbf{C})\freccia\text{\emph{Arr}}(\mathbf{D}),\end{align*}
where\[
\mathcal{AR}(\mathbf{C}):=\left\{ (A,f,B)\,|\,\mathbf{C}\vDash A\xfreccia{f}B\right\} .\]
Because we will always use different symbols for objects and arrows,
and because it should be from the context what domain and codomain
we are considering, we will simply use the notations\begin{align*}
F(A) & :=F_{\text{o}}(A)\quad\forall A\in\text{\emph{Obj}}(\mathbf{C})\\
F(f) & :=F_{\text{a}}(A,f,B)\quad\forall(A,f,B)\in\mathcal{AR}(\mathbf{C}).\end{align*}
Moreover, the following conditions must hold:
\begin{enumerate}
\item F$\left(1_{A}\right)=1_{F(A)}$ for every object $A\in\mathbf{C}$,
i.e. the functor preserves the identity maps.
\item $\mathbf{C}\vDash A\xfrecciad{f}B\then\mathbf{D}\vDash F(A)\xfrecciad{F(f)}F(B)$,
i.e. the functor preserves the arrow relation.
\item $\mathbf{C}\vDash A\xfrecciad{f}B\xfrecciad{g}C\then\mathbf{D}\vDash F(f\cdot g)=F(f)\cdot F(g)$,
i.e. the functor preserves the composition of arrows.
\end{enumerate}
Finally, a functor of the form\[
F:\mathbf{C}^{\text{\emph{op}}}\freccia\mathbf{D}\]
is called a \emph{contravariant} functor.

\end{defn}
\begin{example*}
Let $\mathbf{P}$ and $\mathbf{Q}$ be the categories induced by two
preordered sets $(P,\le)$ and $(Q,\preceq)$ respectively. Then,
only the preservation of the arrow relation is non trivial in this
case, and a functor $f:\mathbf{P}\freccia\mathbf{Q}$ preserves this
relation if and only if\[
x\le y\then f(x)\preceq f(y)\quad\forall x,y\in P,\]
that is the functors correspond to order preserving morphisms.
\end{example*}
If the category $\mathbf{C}$ is locally small, then we can consider
the functor\[
\mathbf{C}(-,-):\mathbf{C}^{\text{op}}\times\mathbf{C}\freccia\Set\]
called the \emph{hom-functor} of $\mathbf{C}$ defined on the objects
$(A,B)$ as the hom-set $\mathbf{C}(A,B)\in\Set$, and on arrows $(A,B)\xfreccia{(f,g)}(C,D)$
as\[
\mathbf{C}(f,g):h\in\mathbf{C}(A,B)\mapsto f\cdot h\cdot g\in\mathbf{C}(C,D).\]

Functors $F:\mathbf{C}\freccia\mathbf{D}$ and $G:\mathbf{D}\freccia\mathbf{E}$
can be composed by considering the composition of the corresponding
maps acting on objects and arrows.
\begin{defn}
\label{def:forgetful-full-faithful-concreteCat}A functor $F:\mathbf{C}\freccia\mathbf{D}$
is called \emph{faithful }(resp. \emph{full}) if and only if for any
two objects $A$, $B\in\mathbf{C}$, the mapping\begin{equation}
f\in\mathbf{C}(A,B)\mapsto F(f)\in\mathbf{D}(FA,FB)\label{eq:mapFunctorFF}\end{equation}
is injective (resp. surjective). A full and faithful functor is called
an \emph{embedding.}

A category $\mathbf{C}$ with a faithful functor $F:\mathbf{C}\freccia\mathbf{D}$
is called a \emph{concrete category based on $\mathbf{D}$}.
\end{defn}
All the categories of sets with a suitable structure and the corresponding
morphisms are concrete categories based on $\Set$. The corresponding
faithful functor associate to each pair $(S,\mathcal{S})$ made of
a set $S$ with the structure $\mathcal{S}$ the underlying set $S\in\Set$,
and to each morphisms the corresponding map between the underlying
sets.
\begin{defn}
\label{def:naturalTransformation}Given two functors $F$, $G:\mathbf{C}\freccia\mathbf{D}$
taking the same domain category $\mathbf{C}$ to the same codomain
category $\mathbf{D}$, we say that $\tau:F\freccia D$ \emph{is a
natural transformation} if and only if
\begin{enumerate}
\item $\tau:\text{\emph{Obj}}(\mathbf{C})\freccia\text{\emph{Arr}}(\mathbf{D})$.
Usually the notation $\tau_{A}:=\tau(A)$ is used.
\item If $A\in\mathbf{C}$, then $\mathbf{D}\vDash F(A)\xfrecciad{\tau_{A}}G(A)$
\item If $\mathbf{C}\vDash A\xfreccia{f}B$, then the following diagram
commutes\[
\xyR{40pt}\xyC{40pt}\xymatrix{F(A)\ar[r]^{{\displaystyle \tau_{A}}}\ar[d]_{{\displaystyle F(f)}} & G(A)\ar[d]^{{\displaystyle G(f)}}\\
F(B)\ar[r]^{{\displaystyle \tau_{B}}} & G(B)}
\]

\end{enumerate}
\end{defn}
If the categories $\mathbf{C}$ and $\mathbf{D}$ are small (in some
universe), then taking as objects all the functors $F:\mathbf{C}\freccia\mathbf{D}$,
as arrows the natural transformations between these functors and with
the composition of natural transformation defined by\[
(\tau\cdot\sigma)_{A}:=\tau_{A}\cdot\sigma_{A},\]
we obtain a category indicated by the symbol $\mathbf{D}^{\mathbf{C}}$.
In this category we can thus say when two functors are isomorphic.
In particular a functor $F:\mathbf{C}\freccia\Set$ is called \emph{representable}
if\[
\exists\, A\in\mathbf{C}\pti\Set^{\mathbf{C}}\vDash F\simeq\mathbf{C}(A,-),\]
such an isomorphism is called a \emph{representation}.

\section{Limits and colimits}
\begin{defn}
\label{def:cone}Let $\mathbf{C}$, $\mathbf{I}$ be two categories
and $F:\mathbf{I}\freccia\mathbf{C}$ a functor, then we say that\[
\left(V\xfrecciad{f_{i}}F(i)\right)_{i\in\mathbf{I}}\text{ is a cone with base }F\]
if and only if:
\begin{enumerate}
\item $f:\text{\emph{Obj}}(\mathbf{I})\freccia\text{\emph{Arr}}(\mathbf{C})$.
We will use the notation $f_{i}:=f(i)$ for $i\in\mathbf{I}$.
\item $V\in\mathbf{C}$
\item $\forall i\in\mathbf{I}\pti\mathbf{C}\vDash V\xfrecciad{f_{i}}F(i)$
\item If $\mathbf{I}\vDash i\xfrecciad{h}j$, then in the category $\mathbf{C}$
the following diagram commutes\[
\xyR{40pt}\xyC{20pt}\xymatrix{ & V\ar[dl]_{{\displaystyle f_{i}}}\ar[dr]^{{\displaystyle f_{j}}}\\
F(i)\ar[rr]_{{\displaystyle F(h)}} &  & F(j)}
\]

\end{enumerate}
\end{defn}
A universal cone with base $F$ is called a limit of $F$:
\begin{defn}
\label{def:categorialLimit}In the previous hypothesis, we say that\[
\left(L\xfrecciad{\mu_{i}}F(i)\right)_{i\in\mathbf{I}}\text{ is a limit of }F\]
if and only if:
\begin{enumerate}
\item $\left(L\xfrecciad{\mu_{i}}F(i)\right)_{i\in\mathbf{I}}\text{ is a cone with base }F$
\item If $\left(V\xfrecciad{f_{i}}F(i)\right)_{i\in\mathbf{I}}$ is another
cone with base $F$, then there exists one and only one morphism $\phi$
such that, in the category $\mathbf{C}$, the following conditions
hold

\begin{enumerate}
\item $V\xfrecciad{\phi}L$
\item For every $i\in\mathbf{I}$, we have\[
\xyR{40pt}\xyC{40pt}\xymatrix{L\ar[r]^{{\displaystyle \mu_{i}}} & F(i)\\
V\ar[u]^{{\displaystyle \phi}}\ar[ru]_{{\displaystyle f_{i}}}}
\]

\end{enumerate}
\end{enumerate}
\end{defn}
The notions of cocone and of colimit are dual with respect to these,
so that the analogous definition can be obtained by simply reversing
the directions of all the arrows. It is possible to prove that if
a limit exists, it is unique up to isomorphisms in $\mathbf{C}$.
For these reasons, if the limit exists, we will denote the corresponding
object $L$ by\[
\lim_{i\in\mathbf{I}}F(i).\]
Analogously the colimit will be denoted by\[
{\displaystyle \mathop{\text{colim}}_{i\in\mathbf{I}}F(i)}.\]

A category $\mathbf{C}$ is said\emph{ }to be \emph{complete} if every
functor $F:\mathbf{I}\freccia\mathbf{C}$ defined in a small category
$\mathbf{I}$ admits a limit; whereas it is said to be \emph{cocomplete}
if each one of such functor admits a colimit.
\begin{example*}
$\ $\end{example*}
\begin{enumerate}
\item If $\mathbf{I}=\{0,1\}$, then the limit $\left(P\xfrecciad{p_{i}}F(i)\right)_{i\in\{0,1\}}$
of $F$ is given by an object $P\in\mathbf{C}$ and two morphisms\[
\xyC{40pt}\xymatrix{F(1) & P\ar[r]\sp(0.45){\displaystyle p_{0}}\ar[l]\sb(0.4){\displaystyle p_{1}} & F(0)}
\]
which verify the universal property: if $\left(V\xfrecciad{f_{i}}F(i)\right)_{i\in\{0,1\}}$is
another pair of morphisms of this form, then there exists one and
only one arrow in $\mathbf{C}$\[
\phi:V\freccia P\]
such that\begin{equation}
\xyR{40pt}\xyC{50pt}\xymatrix{F(0) & P\ar[r]\sp(0.45){\displaystyle p_{1}}\ar[l]\sb(0.4){\displaystyle p_{0}} & F(1)\\
 & V\ar[ur]_{{\displaystyle f_{1}}}\ar[ul]^{{\displaystyle f_{0}}}\ar[u]_{{\displaystyle \phi}}}
\label{eq:universalPropertyOfProduct}\end{equation}
Therefore, in this special case the notion of limit of $F$ gives
the usual notion of product of the objects $F(0)$, $F(1)\in\mathbf{C}$.
In the present work, the unique morphism $\phi$ that verifies \eqref{eq:universalPropertyOfProduct}
is denoted by $\langle f_{1},f_{2}\rangle$. With the notion of cocone
and the same index category $\mathbf{I}=\{0,1\}$ we obtain the usual
notion of sum of two objects.
\item If $\mathbf{I}$ is the category generated by the graph\[
\xyC{40pt}\xymatrix{2\ar[r]^{{\displaystyle a}} & 0 & 1\ar[l]_{{\displaystyle b}}}
\]
then the notion of limit corresponds to the notion of pull-back of
the diagram\[
\xyR{50pt}\xyC{50pt}\xymatrix{ & F(1)\ar[d]^{{\displaystyle F(b)}}\\
F(2)\ar[r]_{{\displaystyle F(a)}} & F(0)}
\]

\item If $\mathbf{I}$ is the category generated by the graph\[
\xyC{40pt}\xymatrix{0\ar@<3pt>[r]^{{\displaystyle a}}\ar@<-3pt>[r]_{{\displaystyle b}} & 1}
\]
then the notion of limit corresponds to that of \emph{equalizer} of
the diagram\[
\xyC{60pt}\xymatrix{F(0)\ar@<3pt>[r]^{{\displaystyle F(a)}}\ar@<-3pt>[r]_{{\displaystyle F(b)}} & F(1)}
\]
that is an arrow $E\xfrecciad{e}F(0)$ such that $e\cdot F(a)=e\cdot F(b)$
which is universal among all the arrows that verify these relations.
\end{enumerate}
In case of concrete categories the notion of limit can be simplified
using the notion of lifting.
\begin{defn}
\label{def:lifting}Let $\mathbf{C}$ be a concrete category based
on $\mathbf{D}$ with faithful functor $U:\mathbf{C}\freccia\mathbf{D}$.
We will use the notation $U^{-1}(f)$ every time $f$ is in the image
set of the functor $U$. Let $I\in\Set$. Then, we say that\[
\left(C\xfrecciad{\gamma_{i}}C_{i}\right)_{i\in I}\text{ is a \emph{lifting} of }\left(D\xfrecciad{\delta_{i}}D_{i}\right)_{i\in I}\]
if and only if:
\begin{enumerate}
\item $U\left(C\xfrecciad{\gamma_{i}}C_{i}\right)=D\xfrecciad{\delta_{i}}D_{i}\quad\forall i\in I$
\item If $\mathbf{D}\vDash U(A)\xfrecciad{\phi}U(C)$ and for every $i\in I$
we have\[
\mathbf{C}\vDash A\xfrecciad{U^{-1}(\phi\cdot\delta_{i})}C_{i}\]
then\[
\mathbf{C}\vDash A\xfrecciad{U^{-1}(\phi)}C\]

\end{enumerate}
\end{defn}
The theorem which connects the two concepts is the following.
\begin{thm}
\label{thm:relationBetweenLimitsAndLifting}Under the previous hypothesis
of Definition \ref{def:lifting}, let us consider a functor $F:I\freccia\mathbf{C}$
and let $\left(D\xfrecciad{\delta_{i}}U(F(i))\right)_{i\in I}$ be
the limit of $U\circ F$ in the category $\mathbf{D}$. Finally, let
us suppose that\[
\left(C\xfrecciad{\gamma_{i}}F(i)\right)_{i\in I}\text{ is a lifting of }\left(D\xfrecciad{\delta_{i}}U(F(i))\right)_{i\in I}.\]
Then\[
\left(C\xfrecciad{\gamma_{i}}F(i)\right)_{i\in I}\text{ is the limit of }F\]

\end{thm}

\section{\label{sec:TheYonedaEmbedding}The Yoneda embedding}

Every object $A$ of a locally small category $\mathbf{C}$ defines
a contravariant functor\[
\text{Y}(A):=\mathbf{C}(-,A):\mathbf{C}^{\text{op}}\freccia\Set.\]
This map $\text{Y}$ can be extended to the arrow of $\mathbf{C}$.
Indeed, every morphism $f:A\freccia B$ in $\mathbf{C}$ induces a
natural transformation\[
\text{Y}(f):=\mathbf{C}(-,f):\text{Y}(A)\freccia\text{Y}(B),\]
so that, at the end we obtain a functor\[
\text{Y}:\mathbf{C}\freccia\Set^{\mathbf{C}^{\text{op}}}\]
called the \emph{Yoneda embedding}. The name is justified by the following
two results. To state the first one of them, we will use the following
language to express a bijection
\begin{defn}
Let $\mathcal{A}(x)$, $\mathcal{B}(y)$ and $\mathcal{C}(x,y)$ be
three property in the free variables $x$ and $y$. Then with the
statement\[
\text{\emph{To give }}x:\, A(x)\text{\emph{ is equivalent to give }}y:\,\mathcal{B}(y)\text{\emph{ so that }}\mathcal{C}(x,y)\text{\emph{ holds}}\]
we mean:
\begin{enumerate}
\item $\forall x:\,\mathcal{A}(x)\Rightarrow\exists!\, y:\,\mathcal{B}(y)\ \text{and}\ \mathcal{C}(x,y)$
\item $\forall y:\,\mathcal{B}(y)\Rightarrow\exists!\, x:\,\mathcal{A}(x)\ \text{and}\ \mathcal{C}(x,y)$
\end{enumerate}
In other words, these properties define a bijection and the property
$\mathcal{C}(x,y)$ acts as the formula connecting the objects $x$
and the objects $y$.

\end{defn}
\begin{thm}
Let $\mathbf{C}$ be a locally small category, $F:\mathbf{C}^{\text{\emph{op}}}\freccia\Set$
a functor and $C\in\mathbf{C}$. Then to give a natural transformation
$\tau$:\begin{equation}
\tau:\text{Y}(C)\freccia F\label{eq:YonedaNaturalTransformation}\end{equation}
is equivalent to give an element $s$:\begin{equation}
s\in F(C)\label{eq:YonedaElement}\end{equation}
so that the following properties hold:
\begin{enumerate}
\item $s=\tau_{C}\left(1_{C}\right)$\label{enu:1_Yoneda}
\item $\tau_{A}(g)=F(g)(s)\quad\forall A\in\mathbf{C}\,\forall g\in\mathbf{C}(A,C).$\label{enu:2_Yoneda}
\end{enumerate}
\end{thm}
\noindent As a consequence of this theorem we have the following result,
which is cited at Chapter \ref{cha:ApproachesToDiffGeomOfInfDimSpaces}
of the present work.
\begin{cor}
The Yoneda embedding is a full and faithful functor.
\end{cor}

\section{Universal arrows and adjoints}
\begin{defn}
\label{def:universalArrow}Let $G:\mathbf{D}\freccia\mathbf{C}$ be
a functor and $C\in\mathbf{C}$, then we say that\[
C\xfrecciad{\eta}G(D)\text{ is a universal arrow}\]
if and only if:
\begin{enumerate}
\item $D\in\mathbf{D}$
\item $\mathbf{C}\vDash C\xfrecciad{\eta}G(D)$
\item The pair $(D,\eta)$ is $G$\emph{-couniversal}%
\footnote{Let us note explicitly the inconsistency between the property of co-universality
(i.e. the unique morphism $\phi$ \emph{starts} from the couniversal
object $C$) and the name {}``universal arrow''. This inconsistency
in the name, even if it creates a little bit of confusion, is well
established in the practice of category theory.%
}\emph{ }among all the pairs which satisfy the previous two conditions,
i.e. if $D_{1}\in\mathbf{D}$ and $\mathbf{C}\vDash C\xfrecciad{\eta_{1}}G(D_{1})$,
then there exists one and only one $\mathbf{D}$-morphism $\phi$
such that

\begin{enumerate}
\item $\mathbf{D}\vDash D\xfrecciad{\phi}D_{1}$
\item $\mathbf{C}\vDash\xyR{40pt}\xyC{40pt}\xymatrix{C\ar[r]^{{\displaystyle \eta}}\ar[dr]_{{\displaystyle \eta_{1}}} & G(D)\ar[d]^{{\displaystyle G(\phi)}}\\
 & G(D_{1})}
$
\end{enumerate}
\end{enumerate}
\end{defn}
The notion of \emph{couniversal arrow} is dual with respect to that
of universal arrow.
\begin{defn}
\label{def:adjointFunctors}Let $\xyC{40pt}\xymatrix{\mathbf{C}\ar@<3pt>[r]^{{\displaystyle F}} & \mathbf{D}\ar@<3pt>[l]^{{\displaystyle G}}}
$ be a pair of functors with opposite directions, then we write\[
F\dashv G\text{ with unit }\eta,\]
and we read it \emph{$G$ is right adjoint of $F$ with unit $\eta$},
if and only if:
\begin{enumerate}
\item $1_{\mathbf{C}}\xfrecciad{\eta}F\cdot G$, i.e. $\eta$ is a natural
transformation from the identity functor $1_{\mathbf{C}}$ to the
composition $F\cdot G=G\circ F$.
\item $C\xfrecciad{\eta_{C}}G(F(C))$ is a universal arrow.
\end{enumerate}
\end{defn}
In case of locally small categories, the notion of pair of adjoint
functors can be reformulated in the following way
\begin{thm}
If $\xyC{40pt}\xymatrix{\mathbf{C}\ar@<3pt>[r]^{{\displaystyle F}} & \mathbf{D}\ar@<3pt>[l]^{{\displaystyle G}}}
$ and $\mathbf{C}$, $\mathbf{D}$ are locally small, then to give
$\eta$:\[
F\dashv G\text{ with unit }\eta,\]
is equivalent to give a natural transformation $\vartheta$:\[
\vartheta:\mathbf{D}\left(F(-),-\right)\xrightarrow[\sim]{\phantom{|\qquad|}}\mathbf{C}\left(-,G(-)\right)\]
so that it results\[
\vartheta_{CD}(\psi)=\eta_{C}\cdot G(\phi)\]
for every $c\in\mathbf{C}$, $D\in\mathbf{D}$ and $\psi\in\mathbf{D}(FC,D)$.
\end{thm}
In the particular case where the categories $\mathbf{C}$ and $\mathbf{D}$
are generated by preordered sets $(C,\le)$ and $(D,\preceq)$, a
pair of adjoints $F\dashv G$ correspond to a pair of order preserving
morphisms such that\[
F(c)\preceq d\iff c\le G(d)\]
(a so called \emph{Galois connection}). In case of concrete categories
based on $\Set$, the notion of cartesian closedness is fully presented
in Chapter \ref{cha:ApproachesToDiffGeomOfInfDimSpaces}. In more
abstract categories, it is defined in the following way.
\begin{defn}
\label{def:cartesianClosedCat}We say that $\left(\mathbf{C},\times,T,\pi,\eps,h\right)$
is a cartesian closed category if and only if:
\begin{enumerate}
\item $\mathbf{C}$ is locally small
\item For every objects $A$, $B\in\mathbf{C}$, the diagram\[
\xyC{50pt}\xymatrix{A & A\times B\ar[r]^{{\displaystyle \pi_{AB}^{2}}}\ar[l]_{{\displaystyle \pi_{AB}^{1}}} & B}
\]
is a product
\item $T$ is a terminal object, i.e. for every $A\in\mathbf{C}$ there
exists one and only one morphism $t$ such that\[
\mathbf{C}\vDash A\xfrecciad{t}T\]

\item For every $A\in\mathbf{C}$\[
(-)\times A\dashv h(A,-)\text{ with counit }\eps^{A}\]

\end{enumerate}
\end{defn}

\chapter[Other theories of infinitesimals]{\label{app:theoriesOfInfinitesimals}A comparison with other theories
of infinitesimals}

It is not easy to clarify in a few pages the relationships between
our theory of Fermat reals and other, more developed and well established
theories of actual infinitesimals. Nevertheless, in this chapter we
want to sketch a first comparison, mostly underlining the conceptual
differences instead of the technical ones, hoping in this way to clarify
the foundational and philosophical choices we made in the present
work.

Our focus will fall on the most studied theories like NSA, SDG and
surreal numbers, or on constructions having analogies with our Fermat
reals like Weil functors and the Levi-Civita field, but we will not
dedicate a section to more algebraic theories whose first aim is not
to develop properties of infinitesimals or infinities and related
applications, but instead to construct a general framework for the
study of fields extending the reals (like formal power series or super-real
fields). In the case of surreal numbers and the Levi-Civita field
we will also give a short presentation of the topic.

A general distinction criterion to classify a theory of infinitesimals
is the possibility to establish a dialogueue between potential infinitesimals
and actual infinitesimals. On the one hand of this dialogue there
are potential infinitesimals, represented by some kind of functions
$i:E\freccia\R$ defined on a directed set $(E,\le)$, like sequences
$i:\N\freccia\R$ or functions defined on a subset $E$ of $\R$,
and such that\begin{equation}
\lim_{(E,\le)}i=0.\label{eq:limitZeroPotentialInfinitesimals}\end{equation}
Classical example are, of course, $i(n)=\frac{1}{n}$ for $n\in\N_{>0}$
and $i(t)=t$ for $t\in\R_{\ge0}$. On the other hand, there are actual
infinitesimals as elements $d\in R$ of a suitable ring $R$. The
dialogue can be realized, if any, in several ways, using e.g. the
standard part and the limit \eqref{eq:limitZeroPotentialInfinitesimals},
or through some connection between the order relation defined on $R$
and the order of the directed set $(E,\le)$, or through the ring
operations of $R$ and pointwise operations on the set of potential
infinitesimals. From our point of view, it is very natural to see
this dialogue as an advantage, if the theory permits this possibility.
First of all, it is a dialogue between two different, but from several
aspects equivalent, instruments to formalize natural phenomena and
mathematical problems, and hence it seems natural to expect a close
relation between them. Secondly, this dialogue can remarkably increase
our intuition on actual infinitesimals and can suggest further generalizations.
For example, in the context of Fermat reals, it seems very natural
to try a generalization taking some stochastic processes $(x_{t})_{t\in\R_{\ge0}}$
instead of our little-oh polynomials, creating in this way {}``stochastic
infinitesimals''.

Theories with a, more or less strong, dialogue between potential infinitesimals
and actual infinitesimals are: NSA, the theory of surreal numbers
and our theory of Fermat reals.

This dialogue, and hence the consequent generalizations or intuitions,
are more difficult in formal algebraic approaches to infinitesimals.
Very roughly, these approaches can be summarized following the spirit
of J. Conway's citation on pag. \pageref{quo:ConwayOnSDG}: if one
needs some kind of infinitesimal $d$, add this new symbol to $\R$
and impose to it the properties you need. In this class of theories
we can inscribe all the other theories: SDG, Weil functors, differential
geometry over a base ring, and Levi-Civita field. They can be thought
of as theories generated by two different elementary ideas: the ring
of dual numbers $\R[\eps]/\left\langle \eps^{2}=0\right\rangle $
(firstly generalized by the strongly stimulating and influential article
\citet{Wei}) and the fields of formal power series. The distinction
between these two classes of theories, those that try a dialogue with
potential infinitesimals and those approaching formally the problem,
is essentially philosophical and at the end choosing one of them rather
than the other one is more of a personal opinion than a rational choice.
First of all, the distinction is not always so crisp, and (non constructive)
NSA represents a case where the above mentioned dialogue cannot always
be performed. Moreover, it is also surely important to note that formal
theories of infinitesimals are able to reach a great formal power
and flexibility, and sometimes through them a sort of a posteriori
intuition about actual infinitesimals can be gained.

\section{Nonstandard Analysis}

A basic request in the construction of NSA is to extend the real field
by a larger field $\hyperR\supseteq\R$. As a consequence of this
request, in NSA every non zero infinitesimal is invertible and so
we cannot have non trivial nilpotent elements (in a field $h^{2}=0$
always implies $h=0$). On the contrary, in the theory of Fermat reals
we aim at obtaining a ring extending the reals, and, as a result of
our choices, we cannot have non-nilpotent infinitesimals, in particular
they cannot be invertible. In the present work, our first aim was
to obtain a meaningful theory from the point of view of the intuitive
interpretation, to the disadvantage of some formal property, only
partially inherited from the real field. Vice versa every construction
in NSA has, as one of its primary aims, to obtain the inheritance
of \emph{all} \emph{the properties }of the reals through the transfer
principle. This way of thinking conducts NSA towards the necessity
to extend every function $f:\R\freccia\R$, e.g. $f=\sin$, from $\R$
to ${}^{*}\R$, and to the property that any sequence of standard
reals $(x_{n})_{n\in\N}\in\R^{\N}$, even the more strange, e.g. $\left(\sin(n)\right)_{n\in\N}$,
represents one and only one hyperreal.

Of course, in the present work we followed a completely different
way: to define the ring of Fermat reals $\ER$ we restrict our construction
to the use of little-oh polynomials $x\in\R_{o}[t]$ only, and therefore
we can extend only smooth functions from $\R$ to $\ER$. Obviously,
our purpose is to develop infinitesimal instruments for smooth differential
geometry only, and we have not the aim of developing an alternative
foundation for all mathematics, like NSA does. In exchange, not every
property is transferred to $\ER$, e.g. our, as presently developed,
is not a meaningful framework where to talk of a $\Cc^{1}$ but not
smooth function $f:\ER\freccia\ER$.

In NSA, this attention to formally inherit every property of the reals
implies that on the one hand we have the greatest formal strength,
but on the other hand we need a higher formal control and sometimes
we lose the intuitive point of view. We can argue for the truth of
this assertion from two points of view: the first one is connected
with the necessity to use a form of the axiom of choice to construct
the non principal ultrafilter needed to define $\hyperR$. In the
second one, we will study more formally the classical motivation used
to introduce $\hyperR$: two sequences of reals are equivalent if
they agree almost everywhere on a {}``large'' set.

It is rather interesting to recall here that the work of \citet{Sc-La}
predates by a few year the construction of $\hyperR$ by A. Robinson.
In \citet{Sc-La} using the the filter of co-finite sets and not an
ultrafilter, a ring extending the real field $\R$ and containing
infinitesimals and infinities is constructed. This work has been of
great inspiration for subsequent works in constructive non-standard
analysis like \citet{Pal1,Pal2,Pal3}, where a field extending the
reals is developed constructively, with a related transfer theorem,
but without a standard part map. Because of their construtive nature,
in these works, no use of the axiom of choice is made.

To study the relationships between the axiom of choice and the hyperreals,
we start from \citet{Conn}, where the author argued that in NSA it
is impossible to give an example of nonstandard infinitesimal, even
{}``to name'' it. More precisely, A. Connes asserts that to any
infinitesimal $e\in\hyperR_{\ne0}$ it is possible to associate, in
a canonical way, a non Lebesgue-measurable subset of $(0,1)$. The
following result of \citet{Sol}
\begin{thm}
\label{thm:SolovayAndLebesgue}There exists a model of the Zermelo-Fraenkel
theory of sets without axiom of choice (\textbf{ZF}), in which every
subset of $\R$ is Lebesgue measurable.
\end{thm}
\noindent would show us the impossibility, in the point of view of
A. Connes, to give an example of infinitesimal in NSA. These affirmations,
not proved in \citet{Conn}, can be formalized using the following
results:
\begin{thm}
\label{thm:infinitesimalsAndUltrafilters}Let $e\in\hyperR_{\ne0}$
be an infinitesimal, and set\[
\mathcal{U}_{e}:=\left\{ X\subseteq\N\,|\,\left[\frac{1}{e}\right]\in{}^{*}X\right\} ,\]
where $[x]$ is the integer part of the hyperreal $x$. Then $\mathcal{U}_{e}$
is an ultrafilter on $\N$ containing the filter of all co-finite
sets.
\end{thm}
\noindent \textbf{Proof:} Directly from the definitions of ultrafilter
and from the properties of the operator ${}^{*}(-)$.$\qedNoNewLine$

\noindent The second result we need is due to \citet{Sie} and does
not need the axiom of choice to be proved:
\begin{thm}
\label{thm:Sierpinski} Let $f:\mathcal{P}(\N)\freccia\{0,1\}$ be
a finitely additive measure defined on every subset of $\N$. For
each $x\in(0,1)$, let \[
x=\frac{1}{2^{n_{1}}}+\frac{1}{2^{n_{2}}}+\frac{1}{2^{n_{3}}}+\dots\]
be the binary representation of $x$, and set\[
\phi(x):=f\mbox{\ensuremath{\left(\left\{  n_{1},n_{2},n_{3},\ldots\right\}  \right)}}.\]
Then the function $\phi:(0,1)\freccia\{0,1\}$ is not Lebesgue measurable
and hence $\phi^{-1}(\{1\})$ is a non Lebesgue-measurable subset
of $(0,1)$.$\qedNoNewLine$
\end{thm}
\noindent Using these results the sentence of A. Connes is now more
clear: to any $e\in\hyperR_{\ne0}$ infinitesimal we can associate
the ultrafilter $\mathcal{U}_{e}$ on $\N$; to this ultrafilter we
can associate the finitely additive measure $f_{e}(S):=1$ if $S\in\mathcal{U}_{e}$
and $f_{e}(S):=0$ if $S\notin\mathcal{U}_{e}$; to this measure we
can finally associate the non Lebesgue-measurable subset of $(0,1)$
given by $S_{e}:=\phi_{e}^{-1}(\{1\})$, where $\phi_{e}$ is defined
as in Theorem \ref{thm:Sierpinski}. The association $e\mapsto S_{e}$
is canonical in the sense that it does not depend on the axiom of
choice. But the result of Solovay, i.e. Theorem \ref{thm:SolovayAndLebesgue},
proves that it is impossible to construct a non Lebesgue-measurable
set without using some form of the axiom of choice, so the association
$e\mapsto S_{e}$ shows the impossibility to define $\hyperR$ without
some form of this axiom%
\footnote{Let us note explicitly, that Theorem \ref{thm:SolovayAndLebesgue}
refers to the full version of the axiom of choice. Indeed, it is well
known, see e.g. \citet{Al-Fe-Ho-Li} and references therein, that
the existence of an ultrafilter on $\N$ is less stronger than the
full axiom of choice. Roughly speaking, we have just proved that if
we are able to construct the hyperreal field $\hyperR$, then some
form of the axiom of choice must holds, not necessarily the full one.%
}. This is the technical result. Whether this can be interpreted as
{}``it is impossible to give an example of infinitesimal in NSA''
or not, it depends on how one means the words {}``to give an example''.
It seems indeed, undeniable that if one accepts the axiom of choice
and $\mathcal{U}$ is an ultrafilter on $\N$ containing the filter
of co-finite sets, then the hyperreal\[
e:=\left[\left(\frac{1}{n}\right)_{n\in\N}\right]_{\mathcal{U}}\in\hyperR\]
is an example of infinitesimal.

The last example seems a typical solution to several problems of NSA
related to the existence of ultrafilters, and can be synthesized in
the sentence {}``the ultrapower construction is intuitively clear
once the ultrafilter is fixed''. For example, an ultrafilter $\mathcal{U}$
on $\N$ containing the filter of co-finite sets is frequently presented
as a possible notion of {}``large sets of natural numbers'' and
the basic equivalence relation\[
(x_{n})_{n}\sim(y_{n})_{n}\DIff\left\{ n\in\N\,|\, x_{n}=y_{n}\right\} \in\mathcal{U}\]
is hence interpreted as {}``the two sequences of real numbers are
almost everywhere equal, i.e. they agree on a large set (with respect
to the notion of large sets given by $\mathcal{U}$)''. We want to
show now that this intuition is not always correct, despite of the
{}``natural'' choice of the ultrafilter $\mathcal{U}$.

To compare two elements of an ultrafilter, i.e. two infinite subsets
of $\N$ we will use the notion of \emph{natural density} (also called
\emph{asymptotic density}, see e.g. \citet{Ten}):
\begin{defn}
\label{def:naturalDensity}If $A\subseteq\N$ and $n\in\N$, we will
set $A_{\le n}:=\{a\in A\,|\, a\le n\}$. Now let $A$, $B$ be subsets
of $\N$, we will say that \emph{there exists the (natural) density
of $A$ with respect to $B$} iff there exists the limit\[
\rho(A,B):=\lim_{n\to+\infty}\frac{\text{\emph{card}}(A_{\le n})}{\text{\emph{card}}(B_{\le n})}\in\R\cup\{+\infty\}.\]
The set of pairs $(A,B)$ for which the density $\rho(A,B)$ is defined
will be denoted by $\mathcal{D}$.
\end{defn}
\noindent For example if $P:=\{2n\,|\, n\in\N\}$ is the set of even
numbers, then $\rho(P,\N)=\frac{1}{2}$, that is the set of even number
is dense $\frac{1}{2}$ with respect to the set of all the natural
numbers.

The notion of natural density has the following properties:
\begin{thm}
\label{thm:propNaturalDensity}Let $A$ and $B$ be subsets of $\N$,
then we have:
\begin{enumerate}
\item $\rho(A,B)=\frac{\text{\emph{card}}(A)}{\text{\emph{card}}(B)}$ if
$A$ and $B$ are both finite.
\item $\rho(A,B)=0$ if $A$ is finite and $B$ is infinite; vice versa
$\rho(A,B)=+\infty$.
\item $\rho(A,B)\le1$ if $A\subseteq B$ and $(A,B)\in\mathcal{D}$.
\item $\rho(-,B)$ is finitely additive.
\item $\rho(m+A,m+B)=\rho(A,B)$ if $(A,B)\in\mathcal{D}$, i.e. the natural
density is translation invariant.
\item \label{enu:multiplesAndDensity}$\rho\left(\left\{ h\cdot n\,|\, n\in\N\right\} ,\N\right)=\frac{1}{h}$
if $h\in\N_{\ne0}$.
\item If $(A,B)$, $(C,D)\in\mathcal{D}$, then the following implications
are true:

\begin{enumerate}
\item $A\cap C=\emptyset\then(A\cup C,B)\in\mathcal{D}$
\item $A\subseteq B\then(B\setminus A,B)\in\mathcal{D}$
\item $A\cup C=B\then(A\cap C,B)\in\mathcal{D}$
\end{enumerate}
\end{enumerate}
\end{thm}
\noindent \textbf{Proof:} see \citet{Ten} and references therein.$\qedNoNewLine$

Our first aim is to generalize the conclusion \emph{\ref{enu:multiplesAndDensity}}.
of this theorem and secondly to prove that given an infinite element
$P\in\mathcal{U}$ of a fixed ultrafilter, we can always find in the
ultrafilter a subset $S\subseteq P$ having density $\frac{1}{2}$
with respect to $P$. This means, intuitively, that an ultrafilter
is closed not only with respect to supersets, but also with respect
to suitable subsets.
\begin{lem}
\label{lem:increasingSequenceAndDensity}Let $b:\N\freccia\N$ be
a strictly increasing sequence of natural numbers, and set for simplicity
of notations\[
B:=\left\{ b_{n}\,|\, n\in\N\right\} .\]
Then we have\[
\rho\left(\left\{ b_{h\cdot n}\,|\, n\in\N\right\} ,B\right)=\frac{1}{h}\quad\forall h\in\N_{\ne0}.\]

\end{lem}
\textbf{Proof:} Let $\text{int}(r)$ be the integer part of the real
$r\in\R$, i.e. the greatest integer number greater or equal to $r$,
and for simplicity of notations set $B_{h}:=\left\{ b_{h\cdot n}\,|\, n\in\N\right\} $.
We first want to prove that \[
\text{card}(B_{h})_{\le n}=\text{int}\left(\frac{\text{card}(B_{\le n})-1}{h}\right)+1.\]
Indeed, since $b$ is strictly increasing, we have\begin{equation}
\text{card}(B_{\le n})=\max\left\{ k\,|\, b_{k}\le n\right\} +1=:K+1\label{eq:firstEqIncrDensity}\end{equation}
\begin{equation}
\text{card}(B_{h})_{\le n}=\max\left\{ k\,|\, b_{h\cdot k}\le n\right\} +1=:H+1,\label{eq:secondEqIncrDensity}\end{equation}
so that we want to prove that $H+1=\text{int}\left(\frac{K}{h}\right)+1$,
i.e. that $H=\text{int}\left(\frac{K}{h}\right)$. In fact from \eqref{eq:secondEqIncrDensity}
we have $b_{h\cdot H}\le n$ and hence $h\cdot H\le K$ from \eqref{eq:firstEqIncrDensity},
i.e. $H\le\frac{K}{h}$ and $H\le\text{int}\left(\frac{K}{h}\right)$.
To prove the opposite, let us consider a generic integer $m\le\frac{K}{h}$
and let us prove that $m\le H$. In fact, since $b$ is increasing
we have $b_{h\cdot m}\le b_{K}$ and $b_{K}\le n$ from \eqref{eq:firstEqIncrDensity}.
Hence $b_{h\cdot m}\le n$ and from \eqref{eq:secondEqIncrDensity}
we obtain the conclusion $m\le H$.

Now we can evaluate the limit (in the sense that this limit exists
if and only if any one of the limits in this series of equalities
exists)\[
\lim_{n\to+\infty}\frac{\text{card}(B_{h})_{\le n}}{\text{card}(B_{\le n})}=\lim_{n\to+\infty}\frac{1}{\text{card}(B_{\le n})}\cdot\left\{ \text{int}\left(\frac{\text{card}(B_{\le n})-1}{h}\right)+1\right\} .\]
Let, for simplicity, $\beta_{n}:=\text{card}(B_{\le n})$ and note
that $\beta_{n}\to+\infty$ because $b$ is strictly increasing. Finally,
let $\text{frac}(r):=r-\text{int}(r)\in[0,1)$ be the fractional part
of the generic real $r\in\R$. With these notations, our limit becomes\begin{multline*}
\lim_{n\to+\infty}\frac{1}{\beta_{n}}\cdot\left\{ \frac{\beta_{n}-1}{h}-\text{frac}\left(\frac{\beta_{n}-1}{h}\right)+1\right\} =\\
=\lim_{n\to+\infty}\left\{ \frac{1}{h}\cdot\left(1-\frac{1}{\beta_{n}}\right)-\frac{1}{\beta_{n}}\cdot\text{frac}\left(\frac{\beta_{n}-1}{h}\right)+\frac{1}{\beta_{n}}\right\} =\frac{1}{h}\end{multline*}
 since $\beta_{n}\to+\infty$ and the fractional part is limited.$\qedNoNewLine$

Now we can prove that if $P$ is an infinite element of a given ultrafilter
$\mathcal{U}$, then in $\mathcal{U}$ we can also find a subset of
$P$ with one half of the elements of $P$.
\begin{lem}
\label{lem:ultrafilterClosedwrtOneHalfSubsets}Let $\mathcal{U}$
be an ultrafilter on $\N$, $P\in\mathcal{U}$ an infinite element
of the ultrafilter and finally $n\in\N_{\ne0}$. Then we can always
find an $S\in\mathcal{U}$ such that
\begin{enumerate}
\item $S\subseteq P$
\item Either $\rho(S,P)=\frac{1}{n}$ or $\rho(S,P)=1-\frac{1}{n}$.
\end{enumerate}
Therefore we have\[
\forall P\in\mathcal{U}:\ P\text{ infinite}\then\exists\, S\in\mathcal{U}:\ S\subseteq P\text{ and }\rho(S,P)=\frac{1}{2}.\]

\end{lem}
\noindent \textbf{Proof:} Since $P$ is infinite, setting\begin{align*}
p_{0}: & =\min(P)\\
p_{n+1}: & =\min\left(P\setminus\left\{ p_{0},\ldots,p_{n}\right\} \right)\end{align*}
we obtain a strictly increasing sequence of natural numbers. Setting
$S':=\left\{ p_{n\cdot k}\,|\, k\in\N\right\} $ from Lemma \ref{lem:increasingSequenceAndDensity}
we have $\rho(S',P)=\frac{1}{n}$. Therefore, if $S'\in\mathcal{U}$,
we have the conclusion for $S:=S'$. Otherwise, we have $\N\setminus S'\in\mathcal{U}$,
so that setting $S:=(\N\setminus S')\cap P=P\setminus S'$ we obtain
$S\in\mathcal{U}$ and\[
\rho(S,P)=\rho(P\setminus S',P)=1-\rho(S',P)=1-\frac{1}{n}.\]
The second part of the conclusion follows setting $n=2$.$\qedNoNewLine$

\noindent Now we only have to apply recursively this lemma to obtain
that in any ultrafilter we can always find elements with arbitrary
small density:
\begin{thm}
\label{thm:closurewrtArbitrarySmallDensity}Let $\mathcal{U}$ be
an ultrafilter on $\N$ and $P\in\mathcal{U}$ with $P$ infinite,
then we can find a sequence $(P_{n})_{n}$ of elements of $\mathcal{U}$
such that
\begin{enumerate}
\item $P_{0}=P$
\item For every $n\in\N$\[
P_{n+1}\subseteq P_{n}\]
\[
\rho(P_{n+1},P_{n})=\frac{1}{2}\]
\[
\rho(P_{n},P)=\frac{1}{2^{n}}.\]

\end{enumerate}
Therefore in any ultrafilter we can always find elements of arbitrary
small density, i.e.\[
\forall\eps>0\ \exists\, S\in\mathcal{U}:\ \rho(S,\N)<\eps.\]

\end{thm}
\noindent \textbf{Proof:} Set $P_{0}:=P$ and apply recursively the
previous Lemma \ref{lem:ultrafilterClosedwrtOneHalfSubsets} (note
that following the proof of this lemma, we can affirm that we are
not applying here the axiom of countable choice) we obtain\[
P_{n}\in\mathcal{U}\e{,}P_{n+1}\subseteq P_{n}\e{,}\rho(P_{n+1},P_{n})=\frac{1}{2}.\]
But\[
\frac{\text{card}(P_{n})_{\le k}}{\text{card}(P_{0})_{\le k}}=\frac{\text{card}(P_{n})_{\le k}}{\text{card}(P_{n-1})_{\le k}}\cdot\frac{\text{card}(P_{n-1})_{\le k}}{\text{card}(P_{n-2})_{\le k}}\cdot\ldots\cdot\frac{\text{card}(P_{1})_{\le k}}{\text{card}(P_{0})_{\le k}}.\]
 Therefore, for $k\to+\infty$ we obtain $\rho(P_{n},P)=\frac{1}{2^{n}}$.
The final sentence of the statement follows from the previous one
if $P:=\N$ and if we take $n$ such that $2^{-n}<\eps$.$\qedNoNewLine$

In the precise sense given by this theorem, we can hence affirm that
in any ultrafilter on $\N$ we can always find also {}``arbitrary
small'' sets. For example, if we set $\eps:=10^{-100}$, we can find
$S\in\mathcal{U}$ with density $\rho(S,\N)<10^{-100}$. The characteristic
function of $S$\[
x_{n}:=\begin{cases}
1 & \text{ if }n\in S\\
0 & \text{ if }n\notin S\end{cases}\quad\forall n\in\N\]
generates, modulo $\mathcal{U}$, an hyperreal $y:=[(x_{n})_{n}]_{\mathcal{U}}\in\hyperR$
equal to $1$ but (with respect to the density $\rho(-,\N)$) almost
always equal to $0$. Finally, the set $S$ of indexes $n\in\N$ where
$x_{n}=1$ has a density strongly lower with respect to the set $ $$\N\setminus S$
of indexes where $x_{n}=0$, in fact \begin{align*}
\rho(S,\N\setminus S) & =\lim_{n\to+\infty}\frac{\text{card}(S_{\le n})}{n+1}\cdot\frac{n+1}{\text{card}(\N\setminus S)_{\le n}}=\\
 & =\frac{\rho(S,\N)}{\rho(\N\setminus S,\N)}\le\frac{10^{-100}}{1-10^{-100}}.\end{align*}

\noindent See also Chapter \ref{cha:IntroductionAndGeneralProblem},
where we already compared NSA with the basic aims of the present work
on Fermat reals.

\section{Synthetic differential geometry}

We have already mentioned, several times, the relationships between
Fermat reals and SDG, and we have already presented very briefly the
main ideas for the construction of a model in SDG (see Section \ref{sec:SDG}).
For these reasons, here we essentially summarize and underline the
differences between the two theories.

There are many analogies between SDG and Fermat reals, so that sometimes
the proofs of several theorems remain almost unchanged. But the differences
are so important that, in spite of the similarities, these theories
can be said to describe {}``different kind of infinitesimals''.

\noindent We have already noted (see Section \ref{sec:Infinitesimals-and-order-properties})
that one of the most important differences is that for the Fermat
reals we have $h\cdot k=0$ if $h^{2}=k^{2}=0$, whereas this is not
the case for SDG, where first order infinitesimals $h,k\in\Delta:=\{d\,|\, d^{2}=0\}$
with $h\cdot k$ not necessarily equal zero, sometimes play an important
role. Note that, as shown in the proof of Schwarz theorem using infinitesimals
(see Section \ref{sec:someElementaryExample}), to bypass this difference,
sometimes completely new ideas are required (to compare our proof
with that of SDG, see e.g. \citet{Koc,Lav}). Because of these diversities,
in our derivation formula we are forced to state $\exists!\, m\in\R_{2}$
and not $\exists!\, m\in\ER$ (see \ref{sec:TheGeneralizedTaylorFormula}).
This is essentially the only important difference between this formula
and the Kock-Lawvere axiom. Indeed to differentiate a generic smooth
map $f:\ER\freccia\ER$ we need the Fermat method (see Section \ref{sec:theFermatMethod})
i.e. the notion of {}``smooth incremental ratio''.

Another point of view regarding the relationships between Fermat reals
and SDG concerns models of SDG. As we hint in Section \ref{sec:SDG},
these models are topos of not simple construction, so that we are
almost compelled to work with the internal language of the topos itself,
that is in intuitionistic logic. If on the one hand this implies that
{}``all our spaces and functions are smooth'', and so we do not
have to prove this, e.g. after every definition, on the other hand
it requires a more strong formal control of the Mathematics one is
doing.

Everyone can be in agreement or not with the assertion whether it
is difficult or easy to learn to work in intuitionistic logic and
after to translate the results using topos based models. Anyway we
think undeniable that the formal beauty achieved by SDG can hardly
be reached using a theory based on classical logic. It suffices to
say, as a simple example, that to prove the infinitesimal linearity
of $M^{N}$ (starting from $M$, $N$ generic inf-linear spaces),
it suffices to fix $n\in N$, to note that $t_{i}(-,n)$ are tangent
vectors at $f(n)$, to consider their parallelogram $p(-,n)$, and
automatically, thanks to the use of intuitionistic logic, $p$ is
smooth without any need to use directly the sheaf property to prove
it. See our Theorem \ref{thm:manifoldsAndExpOfManifoldsAreInfLinear}
to compare this proof with the proof of the analogous statement in
our context.

\noindent On the other hand, if we need a partition of unity, we are
forced to assume a suitable axiom for the existence of bump functions
(whose definition, in the models, necessarily uses the law of the
excluded middle).

Indeed, we think that, as we hint in Chapter \ref{cha:logicalPropertiesOfTheFermatFunctor},
the best properties of the theory of Fermat reals can be obtained
using an {}``intuitionistic interpretation''. We can say that also
this theory {}``proves'' that the best logic to deal with nilpotent
infinitesimals in differential geometry is the intuitionistic one
and not the classical one. All the efforts done in the present work
can be framed into an attempt to obtain a sufficiently simple model
of nilpotent infinitesimals, having a strong intuitive interpretation
but, at the same time, without forcing the reader to switch to intuitionistic
logic. Indeed, we think that the best result in the theory of Fermat
reals would be to prove that the category of smooth spaces $\CInfty$
and that of Fermat spaces $\ECInfty$ are really topoi: in this way
the reader working in this theory would have the possibility to use
the internal language of these topoi, in intuitionistic logic, and
at the same time a sufficiently simple model to work directly in classical
logic or to interpret the results obtained using the internal language.
We plan to achieve some steps in this direction in future works.

Moreover, from the intuitive, classical, point of view, SDG sometimes
presents counter-intuitive properties. For example, it is a little
strange to think that we do not have {}``examples'' of infinitesimals
in SDG (it is only possible to prove that $\neg\neg\exists d\in\Delta$),
so that, e.g., we cannot construct a physical theory containing a
fixed infinitesimal parameter; moreover any $d\in\Delta$ is at the
same time negative $d\le0$ and positive $d\ge0$; finally the definition
of the Lie brackets using $h\cdot k$ for $h,k\in\Delta$, i.e. \[
[X,Y]_{h\cdot k}=Y_{-k}\circ X_{-h}\circ Y_{k}\circ X_{h},\]
is very far to the usual definitions given on manifolds.

\section{Weil functors}

Weil functors (in the following WF; see \citet{Ko-Mi-Sl} and \citet{Kr-Mi2})
represent a way to introduce some kind of useful infinitesimal method
without the need to possess a non-trivial background in mathematical
logic. The construction of WF does not achieve the construction of
a whole {}``infinitesimal universe'', like in the theory of Fermat
reals or in NSA and SDG, but it defines functors $T_{A}:\ManInfty\freccia\ManInfty$,
related to certain geometrical constructions of interest, starting
from a \emph{Weil algebra}. A Weil algebra is a real commutative algebra
with unit of the form $A=\R\cdot1\oplus N$, where $N$ is a finite
dimensional ideal of nilpotent elements. The flexibility of its input
$A$ gives a corresponding flexibility to the construction of these
functors. But, generally speaking, if one changes the geometrical
problem, one has also to change the algebra $A$ and so the corresponding
functor $T_{A}$. E.g., if $A=\R[x]/\langle x^{2}\rangle$, then $T_{A}$
is the ordinary tangent bundle functor, whereas if $B=R[x,y]/\langle x^{2},y^{2}\rangle$,
then $T_{B}=T_{A}\circ T_{A}$ is the second tangent bundle. The definition
of a WF starting from a generic Weil algebra $A$ is very long, and
we refer the reader e.g. to \citealt{Kr-Mi,Kr-Mi2,Ko-Mi-Sl}. Note
that, in the previous example, $x$, $y\in B$ verify $x^{2}=y^{2}=0$
but $x\cdot y\ne0$. This provides us the first difference between
WF and Fermat reals. In fact $\ER=\R\cdot1\oplus D_{\infty}$ and
$\dim_{\R}D_{\infty}=\infty$, so that using the infinitesimals of
$\ER$ we can generate a large family of Weil algebras, e.g. any $A=\R\cdot1\oplus N\subset\R\cdot1\oplus D_{k}$
(which represents $k$-th order infinitesimal Taylor's formulas) where
$N$ is an $\R$-finite dimensional ideal of infinitesimals taken
in $D_{k}$. On the other hand, not every algebra can be generated
in this way, e.g. the previous $B=R[x,y]/\langle x^{2},y^{2}\rangle$.
But using exponential objects of $\CInfty$ and $\ECInfty$ we can
give a simple infinitesimal representation of a large class of WF.
For $\alpha_{1},\ldots,\alpha_{c}\in\N^{n}$, $c\ge n$, let \[
D_{k}^{\alpha}:=\left\{ h\in D_{k_{1}}\times\ldots\times D_{k_{n}}\,|\, h^{\alpha_{i}}=0\ \ \forall i=1,\ldots,c\right\} .\]
E.g. if $k_{1}=(3,0)$, $k_{2}=(0,2)$ and $\alpha=(1,1)$, then $D_{k}^{\alpha}=\{(h,k)\in D_{3}\times D_{2}\,|\, h\cdot k=0\}$.
To any infinitesimal object $D_{k}^{\alpha}$ there is associated
a corresponding Taylor's formula: let $f=\ext{g}|_{D_{k}^{\alpha}}$,
with $g\in\Cc^{\infty}(\R^{n},\R)$, then \begin{equation}
f(h)=\sum_{r\in\iota(\alpha)}\frac{h^{r}}{r!}\cdot m_{r}\quad\forall h\in D_{k}^{\alpha}.\label{eq:TaylorInD_k^alpha}\end{equation}
Here $\iota(\alpha):=\{r\in\N^{n}\,|\,\exists h\in D_{k}^{\alpha}:h^{r}\ne0,|r|\le k\}$
is the set of multi-indexes $r\in\N^{n}$ corresponding to a non zero
power $h^{r}$, and $k:=\max(k_{1},\dots,k_{n})$. The coefficients
$m_{r}=\frac{\partial^{r}g}{\partial x^{r}}(0)\in\R$ are uniquely
determined by the formula \eqref{eq:TaylorInD_k^alpha}. We can therefore
proceed generalizing the definition \ref{def: StndTangentFunctor}
of standard tangent functor. 
\begin{defn}
If $M\in\ManInfty$ is a manifold, we call $M^{D_{k}^{\alpha}}$ the
$\CInfty$ object with support set \[
\left|M^{D_{k}^{\alpha}}\right|:=\{\ext{f}|_{D_{k}^{\alpha}}\,:\, f\in\CInfty(\R^{n},M)\},\]
 and with generalized elements of type $U$ (open in $\R^{u}$) defined
by: \[
d\In{U}M^{D_{k}^{\alpha}}\DIff d:U\freccia\left|M^{D_{k}^{\alpha}}\right|\e{and}d\cdot i\In{\bar{U}}\ext{M}^{D_{k}^{\alpha}},\]
 where $i:\left|M^{D_{k}^{\alpha}}\right|\hookrightarrow\ext{M}^{D_{k}^{\alpha}}$
is the inclusion.
\end{defn}
Let us note explicitly that writing $M^{D_{k}^{\alpha}}$ we are doing
an abuse of notation because this is not an exponential object. We
can extend this definition to the arrows of $\ManInfty$ by setting
$f^{D_{k}^{\alpha}}(t):=t\cdot f\in N^{D_{k}^{\alpha}}$, where $t\in M^{D_{k}^{\alpha}}$
and $f\in\ManInfty(M,N)$. With these definitions we obtain a product
preserving functor $(-)^{D_{k}^{\alpha}}:\ManInfty\freccia\ManInfty$.
Finally we have a natural transformation $e_{0}:(-)^{D_{k}^{\alpha}}\freccia1_{\ManInfty}$
defined by evaluation at $0\in\R^{n}$: $e_{0}(M)(t):=t(0)$. The
functor $(-)^{D_{k}^{\alpha}}$ and the natural transformation $e_{0}$
verify the {}``locality condition'' of Theorem 1.36.1 in \citet{Ko-Mi-Sl}:
if $U$ is open in $M$ and $i:U\hookrightarrow M$ is the inclusion,
then $U^{D_{k}^{\alpha}}=e_{0}(M)^{-1}(U)$ and $i^{D_{k}^{\alpha}}$
is the inclusion of $U^{D_{k}^{\alpha}}$ in $M^{D_{k}^{\alpha}}$.
We can thus apply the above cited theorem to obtain that $(-)^{D_{k}^{\alpha}}$
is a Weil functor, whose algebra is \[
\text{Al}\left((-)^{D_{k}^{\alpha}}\right)=\R^{D_{k}^{\alpha}}.\]

Not every Weil functor has this simple infinitesimal representation.
E.g., the second tangent bundle $(-)^{D}\circ(-)^{D}$ is not of type
$(-)^{D_{k}^{\alpha}}$; indeed it is easy to prove that the only
possible candidate could be $D_{k}^{\alpha}=D\times D$, but $(\R^{D})^{D}$
is a four dimensional manifold, whereas $\R^{D\times D}$ has dimension
three. We do not have this kind of problems with the functor $(-)^{D_{k}^{\alpha}}=\ECInfty(D_{k}^{\alpha},-):\ECInfty\freccia\ECInfty$
which generalizes the previous one as well as ${\rm T}M=\ext{M}^{D}$
generalizes the standard tangent functor. In fact because of cartesian
closedness we have \[
\left(X^{D_{k}^{\alpha}}\right)^{D_{h}^{\beta}}\simeq X^{D_{k}^{\alpha}\times D_{h}^{\beta}}\]
 and $D_{k}^{\alpha}\times D_{h}^{\beta}$ is again of type $D_{k}^{\alpha}$.

Summarizing, we can affirm that WF permit to consider nilpotent infinitesimals
which are more algebraic and hence more general than those occurring
in the Fermat reals. The typical example is the WF $T_{B}$ for $B=R[x,y]/\langle x^{2},y^{2}\rangle$,
corresponding to the second tangent bundle. On the other hand, WF
do not permit to consider an extension of the real field with the
addition of new infinitesimal points (like in our framework, where
we have the extension from $\R$ to $\ER$), and hence they do not
permit to consider properties like order between infinitesimals, an
extension functor analogous of the Fermat functor and the related
properties, like the transfer theorem, tangent vectors as infinitesimal
curves, infinitesimal parallelograms to add tangent vectors, infinitesimal
fluxes, and so on. This implies that with the WF we do not have {}``a
framework with the possibility to extend standard spaces adding infinitesimals'',
but we are forced to consider a new WF for every geometrical construction
we are considering. Finally, the general definition of WF works on
the category of smooth manifolds modelled on convenient vector spaces,
because it needs the existence of charts (see \citet{Kr-Mi,Kr-Mi2,Ko-Mi-Sl}),
and we already mentioned (see Section \ref{sec:TheConvenientVectorSpacesSettings})
that this category is not cartesian closed. Therefore, WF cannot be
defined for spaces like $N^{M}$, where $M$ is a non-compact manifold.
On the contrary, we have seen (see Chapter \ref{cha:Infinitesimal-differential-geometry})
that some results of infinitesimal differential geometry can be obtained
also for spaces of the form $\ext{N}^{\ext{M}}$, where $M$ is a
generic manifold.

Finally, a recent approach similar in essence to Weil functors is
\emph{differential geometry over a general base ring}, see \citet{Ber}
and references therein. The basic idea is to develop, as far as possible,
all the topics of differential geometry not dealing with integration
theory, in the framework of manifolds modelled over a generic topological
module $V$ over a topological ring $\mathbb{K}$. This of course,
includes ordinary finite dimensional real or complex manifolds, but
also infinite dimensional manifolds modelled on Banach spaces and
even on the hyper-vector spaces $\hyperR^{n}$. One of the basic results
is that in this way the tangent functor $TM$ becomes a manifold over
the scalar extension $V\oplus\eps V$, i.e. over the module of all
the expressions of the form $u+\eps v$ over the ring $\mathbb{K}[\eps]$
of dual numbers over $ $$\mathbb{K}$, i.e. $\mathbb{K}[\eps]:=\mathbb{K}\oplus\eps\mathbb{K}:=\mathbb{K}[x]/\left(x^{2}\right)$.
The process can be iterated obtaining that the double tangent bundle
$T^{2}M$ is a manifold over $V\oplus\eps_{1}V\oplus\eps_{2}V\oplus\eps_{1}\eps_{2}V$,
which is a module over $\mathbb{K}[\eps_{1},\eps_{2}]:=\mathbb{K}[x,y]/\left(x^{2},y^{2}\right)$.
Analogous results are available for $T^{k}M$ and for the jet bundle
$J^{k}M$. The theory is appealing for its generality and for the
possibility to obtain in a simple way a context with formal infinitesimals.
This construction does not deal with cartesian closedness and hence
generic spaces like $\ManInfty(M,N)$ cannot be considered.

\section{Surreal numbers}

Surreal numbers has been introduced by J.H. Conway and presented in
\citet{Knu} and in \citet{Con2}%
\footnote{Really, the same field of numbers has been predate by \citet{Cue}
(in Spanish) and \citet{Harz} (in German).%
}. One of the most surprising features of surreal numbers is that starting
from a simple set of rules it is possible to construct a rich algebraic
structure containing the real numbers as well as infinite and infinitesimals,
but also all the ordinal numbers, the hyperreals of NSA, the Levi-Civita
field and the field of rational functions. Indeed, in a precise sense
we will see later, the ordered field $\text{\textbf{No}}$ of surreal
numbers is the largest possible ordered field or, in other words,
the above mentioned simple rules for the construction of surreal numbers,
represent the most general way to obtain a notion of number culminating
in an ordered field.

There are two basic ideas to introduce surreal numbers: the first
is to have the possibility to construct numbers in a transfinite-recursive
way using a notion analogous to that of Dedekind cut (called \emph{Conway
cut}). If we have a totally ordered set $(N,<)$, a Conway cut is
simply a pair $(L,R)$ of subsets $L$, $R\subseteq N$ such that
\begin{equation}
\forall l\in L\,\forall r\in R\pti l<r,\label{eq:L_lessThan_R_Def}\end{equation}
in this case we will simply write $L<R$. This is exactly the notion
of Dedekind cut without the condition that the subsets $L$, $R$
have to be contiguous (i.e. without the condition that $\forall\eps>0\,\exists\, l\in L\,\exists\, r\in R:\ |l-r|<\eps$).
Exactly because we do not have this further condition, we need another
condition for a pair $(L,R)$ to identify a unique {}``number''.
Indeed, the second idea, intuitively stated, is that every Conway
cut identifies uniquely \emph{the simplest number} $x$ between $L$
and $R$:\begin{equation}
\forall l\in L\,\forall r\in R\pti l<x<r.\label{eq:simplicityWRT_options}\end{equation}

\noindent We can intuitively represent a Conway cut and the associated
simplest number in the following way\[
\xymatrix{\xyC{50pt}{}\ar@{-}[r]^{L} & x\ar@{-}[r]^{R} & {}}
\]
A little more formally, the class $\text{\textbf{No}}$ of surreal
numbers is introduced by Conway using a suitable set of rules. We
can think at these rules as axioms defining a suitable structure $(\text{\textbf{No}},\le,\{-\,|\,-\})$.
In the following, as usual, $x<y$ means $x\le y$ and $x\ne y$.
\begin{description}
\item [{Construction}] If $L$, $R\subseteq\text{\textbf{No}}$ and $L<R$,
then $\{L\,|\, R\}\in\text{\textbf{No}}$, that is starting from a
Conway cut $(L,R)$ we can construct a surreal with $\{L\,|\, R\}\in\text{\textbf{No}}$.
\item [{Surjectivity}] If $x\in\text{\textbf{No}}$, then there exist $L$,
$R\subseteq\text{\textbf{No}}$ such that $L<R$ and $x=\{L\,|\, R\}$,
that is all surreal numbers can be constructed starting from a Conway
cut.
\item [{Inequality}] If $x=\{L_{x}\,|\, R_{x}\}$ and $y=\{L_{y}\,|\, R_{y}\}$
are well defined%
\footnote{That is $L_{x}<R_{x}$ and $L_{y}<R_{y}$. Let us note that using
a notation like $x=\{L_{x}\,|\, R_{x}\}$ we do not mean that a number
$x\in\text{\textbf{No}}$ uniquely determines the subsets $L_{x}$
and $R_{x}$.%
}, then $x\le y$ if and only if $L_{x}<\{y\}$ and $\{x\}<R_{y}$,
i.e. $l_{x}<y$ and $x<r_{y}$ for every $l_{x}\in L_{x}$ and every
$r_{y}\in R_{y}$. This rule can be represented in the following way\[
\xymatrix{\xyC{50pt}{}\ar@{-}[r]^{L_{x}} & x\ar@{-}[r]^{R_{x}} & {\ldots}\ar@{-}[r]^{L_{y}} & {y}\ar@{-}[r]^{R_{y}} & {}}
\]

\item [{Equality}] If $x$, $y\in\text{\textbf{No}}$, then $x=y$ if and
only if $x\le y$ and $y\le x$, that is equality between surreal
numbers can be defined starting from the order relation $\le$.
\end{description}
Starting from these simple rules/axioms we can already construct several
meaningful examples of numbers in $\text{\textbf{No}}$. From the
definition \eqref{eq:L_lessThan_R_Def} of $L<R$ we see that always
$L<\emptyset$ and $\emptyset<R$ for every $L$, $R\subseteq\text{\textbf{No}}$.
So we have $\emptyset<\emptyset$ and from the Construction rule $\{\emptyset\,|\,\emptyset\}\in\text{\textbf{No}}$.
Therefore, $\text{\textbf{No}}$ is not empty and we can iterate the
process. For simplicity, we will write $\{|\}:=\{\emptyset\,|\,\emptyset\}$,
$\{L\,|\}:=\{L\,|\,\emptyset\}$, $\{|\, R\}:=\{\emptyset\,|\, R\}$,
$\{x_{1},\ldots,x_{n}\,|\, R\}:=\left\{ \{x_{1},\ldots,x_{n}\}\,|\, R\right\} $
and $\{L\,|\, x_{1},\ldots,x_{n}\}:=\left\{ L\,|\,\{x_{1},\ldots,x_{n}\}\right\} $.
Hence we have, e.g.,\begin{align*}
\{|\}\in\text{\textbf{No}}\\
x\in\text{\textbf{No}} & \then\{|x\},\ \{x|\}\in\text{\textbf{No}}.\end{align*}
But the understanding of the class $\text{\textbf{No}}$ has a great
improvement if we introduce the above mentioned interpretation of
simplicity. Conway's idea is that a number $x\in\text{\textbf{No}}$
is simpler than $y\in\text{\textbf{No}}$ if $x$ is defined before
$y$ in the previous iterative process (using Conway's terminology:
$x$ was born before $y$). So, $\{|\}$ is the simplest number%
\footnote{Of course, at this stage of developement and using this not-strictly
formal point of view, our use of the notion of {}``simplicity''
is only informal and it is natural to ask for a more formal definition,
considering, moreover, its uniqueness. This will be done in the next
section.%
} and this justify the definition $\{|\}=:0$. On the next step, we
have e.g. $\{|0\}$ and $\{0|\}$ which have the same degree of simplicity
(because they are both defined in the second step of the iterative
process). To interpret these numbers we have hence to use the idea
of simplicity expressed in \eqref{eq:simplicityWRT_options}: $\{0|\}$
is the simplest number greater than $0$ and $\{|0\}$ is the simplest
number less than $0$. This justify the definition $1:=\{0|\}$ and
$-1:=\{|0\}$. Up to isomorphisms we can hence affirm that $\N$,
$\mathbb{Z}\subseteq\text{\textbf{No}}$. Another meaningful example
based on this interpretation is the number $\{0\,|\,1\}$ which has
to be thought as the simplest number between $0$ and $1$, i.e. $\frac{1}{2}:=\{0\,|\,1\}$.
From the Inequality rule we can prove that $\{|\, x\}\le0\le\{x\,|\}$
for every $x\in\text{\textbf{No}}$, and that $1\nleqslant0$, hence
$0<1$ follows from the Equality rule. Analogously one can prove that
$\ldots<-3<-2<-1<0<1<2<3<\ldots$ Moreover, we can also easily see
that e.g. $\{-1\,|\,1\}=0$ so we note that different subsets $L$,
$R\subseteq\text{\textbf{No}}$ can define the same number $x=\{L\,|,R\}$.

\noindent But now we also have that $\N\subseteq\text{\textbf{No}}$
and hence we can form the number $\N+1:=\{\N\,|\}$ and this catch
a glimpse of the possibility to extend all this using transfinite
induction.

Instead of further proceeding with Conway's approach to $\text{\textbf{No}}$
we want to sketch his point of view to the foundational questions
arising from his construction. These ideas are precisely stated in
the Appendix to Part Zero of \citep{Con2}. The mainstream's approach
to a topic like $\text{\textbf{No}}$, where one must use e.g. transfinite
induction, is to fix a formal theory of sets (like Zermelo-Fraenkel
\textbf{ZF} theory or von Neumann-Bernay-G{\"o}del \textbf{NBG} theory)
and to formalize every construction inside that theory. Sometimes,
this formalization can conduct to a theory far from the original intuition,
because different formalizations are possible of a given informal
construction. Let us note explicitly that here J. Conway's term of
judgment is the searching for a good dialectic between informal constructions
and their formal counterpart, which has been a leading thread of all
the present work. Conway's foundational point of view can be summarized
citing, as in \citep{Con2}:
\begin{quotation}
It seems to us, however, that mathematics has now reached the stage
where formalization within some particular axiomatic theory is irrelevant,
even for foundational studies. It should be possible to specify conditions
on a mathematical theory which would suffice for embeddability within
\textbf{ZF} (supplemented by additional axioms of infinity if necessary),
but which do not otherwise restrict the possible constructions in
that theory.
\end{quotation}

\subsubsection{The Cuesta Dutari definition of $\text{\textbf{No}}$}

\noindent From our point of view, one has the best result if there
exists a formalization respecting the above mentioned good formal-informal
dialectic, and hence a great effort has to be dedicated to the searching
for this best formalization, if any, before assuming a point of view
so general and radical like J. Conway's one. This type of formalization
for surreal numbers is possible and is given by \citep{Cue}. In this
appendix we only sketch the first steps of this development; for a
complete treatment, see \citep{All}. The set theory we will consider
is \textbf{NBG}.
\begin{defn}
Let $(T,\le)$ a totally ordered set, then we say that $(L,R)$\emph{
is a Cuesta Dutari cut in }$T$ iff
\begin{enumerate}
\item $L$, $M\subseteq T$
\item $L<R$
\item \label{enu:CD_unionIs_T}$L\cup R=T$.
\end{enumerate}
Moreover, we will denote by $\text{\emph{CD}}(T)$ the class of all
Cuesta Dutari cuts in $T$.

\end{defn}
\noindent Essentially a Cuesta Dutari cuts is a Conway cut with the
additional condition $L\cup R=T$. Since $(\emptyset,T)$ and $(T,\emptyset)$
are always Cuesta Dutari cuts, we have that $\text{CD}(T)$ is never
empty.

If we think than each Cuesta Dutari cut identifies a new number, the
union $T\cup\text{CD}(T)$ can be thought of as a \emph{completion}
of the totally ordered set $(T,\le)$ if we can extend the order relation
$\le$ to this union. This is done in the following
\begin{defn}
Let $(T,\le)$ be a totally ordered set, then on the \emph{Cuesta
Dutari completion $\chi(T):=T\cup\text{CD}(T)$ of $T$} we define
the order relation:
\begin{enumerate}
\item If $x$, $y\in T$ then we will say that $x$ is less than or equal
to $y$ iff $x\le y$ in $T$. Because of this first case, the order
relation on $\chi(T)$ will be denoted again by the symbol $\le$.
\item If $x\in T$ and $y=(L,R)\in\text{\emph{CD}}(T)$, then:

\begin{enumerate}
\item $x\in L\then x<y$
\item $x\in R\then y<x$
\end{enumerate}
\item If $x=(L_{x},R_{x}),$ $y=(L_{y},R_{y})\in\chi(T)$, then $x<y$ iff
$L_{x}\subset L_{y}$.
\end{enumerate}
\end{defn}
It is indeed possible to prove (see \citet{All}) that $(\chi(T),\le)$
is a totally ordered set. For example if we take $t$, $\tau\in T$
with $t<\tau$, we can consider the cut $c=((-\infty,t],[\tau,+\infty))$
and we have $t<c<\tau$. If $x=(L,R)\in\text{CD}(T)$, then $L<\{x\}<R$
and, as a further example, $(\emptyset,T)$ is the least element of
$\chi(T)$, whereas $(T,\emptyset)$ is the greatest element.

So, how can we form $0$ using Cuesta Dutari cuts? We do not have
to think at the Cuesta Dutari completion as a final completion starting
from a single given ordered set $(T,\le)$ but, instead, as a tool
for a transfinite-recursive construction:
\begin{defn}
Let $\text{\textbf{\emph{On}}}$ be the class of all ordinals, we
define by transfinite recursion the family $\left(T_{\alpha}\right)_{\alpha\in\text{\textbf{\emph{On}}}}$
of ordered sets given by:
\begin{enumerate}
\item $T_{0}$ is the empty set ordered with the empty relation,
\item For every $\beta\in\text{\textbf{\emph{On}}}$:

\begin{enumerate}
\item If $\alpha+1=\beta$, then $T_{\beta}:=\chi(T_{\alpha})$
\item If $\beta$ is a non-zero limit ordinal, then $T_{\beta}:=\bigcup_{\alpha<\beta}T_{\alpha}$.
\end{enumerate}
\end{enumerate}
Finally we set $\text{\textbf{\emph{No}}}:=\bigcup_{\alpha\in\text{\textbf{\emph{On}}}}T_{\alpha}$.

\end{defn}
So, e.g., $0:=(\emptyset,\emptyset)\in T_{1}=\chi(T_{0})\subset\text{\textbf{No}}$.
The ordinal index $\alpha$ in the previous transfinite recursive
definition gives the notion of simplicity of a number $x\in\text{\textbf{No}}$,
that is its \emph{birthday} using Conway's terminology.
\begin{defn}
If $x\in\text{\emph{\textbf{No}}}$, we define the \emph{birth-order}
function by\[
b(x):=\min_{\le}\left\{ \alpha\in\text{\emph{\textbf{No}}}\,|\, x\in T_{\alpha}\right\} ,\]
where $\le$ is the order relation defined on $\textbf{\emph{No}}$.
\end{defn}
So we have, e.g., $b(0)=0$, $b(1)=b(-1)=1$, $b(2)=b\left(\frac{1}{2}\right)=b\left(-\frac{1}{2}\right)=2$.

At this point, the most important result is Conway's simplicity theorem:
it states that every Conway cut $(L,R)$ determines uniquely the simplest
surreal number filling the gap between the subsets $L$ and $R$:
\begin{thm}
\label{thm:ConwaySimplicityTheorem}Let $L$, $ $$R\subseteq\text{\emph{\textbf{No}}}$
with $L<R$, then there exist one and only one $x\in\text{\emph{\textbf{No}}}$
such that:
\begin{enumerate}
\item $L<\{x\}<R$
\item If $y\in\text{\emph{\textbf{No}}}\setminus\{x\}$ verifies $L<\{y\}<R$,
then $b(x)<b(y)$.
\end{enumerate}
This unique $x$ will be denoted by $\{L\,|\, R\}$.

\end{thm}
\noindent For a proof see \citet{All}.

\subsubsection{Ring operations on $\text{\textbf{No}}$}

On the class of surreal numbers we can define all the field operations,
the simplest one being subtraction:\[
x=\{L\,|\, R\}\then-x:=\{-R\,|\,-L\}.\]

\noindent The definition of sum and product requires some motivations.
Let $x=\{L_{x}\,|\, R_{x}\}$ and $y=\{L_{y}\,|\, R_{y}\}$ be surreal
numbers, so that\begin{align}
l_{x} & <x<r_{x}\label{eq:inequalities_x_surreal}\\
l_{y} & <y<r_{y}\label{eq:inequalities_y_surreal}\end{align}
for every $l_{x}\in L_{x}$, $r_{x}\in R_{x}$, $l_{y}\in L_{y}$
and $r_{y}\in R_{y}$. Since we want $\text{\textbf{No}}$ to be an
ordered group with respect to addition, we must have:\begin{align*}
l_{x}+y & <x+y\\
x+l_{y} & <x+y\\
x+y & <r_{x}+y\\
x+y & <x+r_{y}.\end{align*}
Hence Conway defines $x+y$ as the simplest number verifying these
inequalities, i.e. using transfinite recursion we can define\[
x+y:=\left\{ (L_{x}+y)\cup(x+L_{y})\,|\,(R_{x}+y)\cup(x+R_{y})\right\} ,\]
where, e.g., $L_{x}+y:=\{l_{x}+y\,:\, l_{x}\in L_{x}\}$. Analogously
we can proceed to justify the definition of product. From \eqref{eq:inequalities_x_surreal}
and \eqref{eq:inequalities_y_surreal}, in the hypothesis that $\text{\textbf{No}}$
be an ordered group under multiplication, we must have that $x-l_{x}$,
$y-l_{y}$, $r_{x}-x$ and $r_{y}-y$ are all greater than zero. Taking
all the products of these terms involving an $x$ and a $y$ we have\begin{align*}
0 & <(x-l_{x})\cdot(y-l_{y})=xy-l_{x}y-xl_{y}+l_{x}l_{y}\\
0 & <(r_{x}-x)\cdot(r_{y}-y)=xy-r_{x}y-xr_{y}+r_{x}r_{y}\\
0 & <(x-l_{x})\cdot(r_{y}-y)=-xy+l_{x}y+xr_{y}-l_{x}r_{y}\\
0 & <(r_{x}-x)\cdot(y-l_{y})=-xy+r_{x}y+xl_{y}+r_{x}l_{y}.\end{align*}
As a consequence, from these we get inequalities bounding $xy$:\begin{align*}
l_{x}y+xl_{y}-l_{x}l_{y} & <xy<l_{x}y+xr_{y}-l_{x}r_{y}\\
r_{x}y+xr_{y}-r_{x}r_{y} & <xy<r_{x}y+xl_{y}+r_{x}l_{y}.\end{align*}
We can hence define (once again by transfinite recursion):\begin{align*}
L_{x\cdot y}: & =(L_{x}y+xL_{y}-L_{x}L_{y})\cup(R_{x}y+xR_{y}-R_{x}R_{y})\\
R_{x\cdot y}: & =(L_{x}y+xR_{y}-L_{x}R_{y})\cup(R_{x}y+xL_{y}+R_{x}L_{y})\\
x\cdot y: & =\left\{ L_{x\cdot y}\,|\, R_{x\cdot y}\right\} ,\end{align*}
where e.g. $L_{x}y:=\{l_{x}\cdot y\,:\, l_{x}\in L_{x}\}$ and $L_{x}L_{y}:=\{l_{x}\cdot l_{y}\,|\, l_{x}\in L_{x}\,,\, l_{y}\in L_{y}\}$.Using
these definitions we can prove that $\text{\textbf{No}}$ verifies
the axioms of an ordered field.

\subsubsection{Examples of surreal numbers}

As we already sketched, up to isomorphism we have $n=\{0,1,2,\ldots,n-1\,|\}\in\N\subseteq\text{\textbf{No}}$,
$-n=\{|\,-n+1,-n+2,\ldots,-2,-1,0\}\in\mathbb{Z}\subseteq\text{\textbf{No}}$,
but also $\omega:=\{\N\,|\}$. It results $n<\omega$ for every $n\in\N$
and hence $\text{\textbf{No}}$ is a non-Archimedean field. Moreover,
because we have an ordered field containing the integers, we also
have $\mathbb{Q}\subseteq\text{\textbf{No}}$, i.e. all the rationals
can be seen as surreal numbers. Finally, using Dedekind cuts we can
also identify $\R$ with a subfield of $\text{\textbf{No}}$.

But using transfinite induction we can also define\begin{align*}
\phi(0): & =0\\
\phi(\beta): & =\left\{ \{\phi(\alpha)\,:\,\alpha<\beta\}\,|\right\} \in\text{\textbf{No}}\quad\forall\beta\in\text{\textbf{On}},\end{align*}
and we can prove that $\phi$ is an order-preserving map from the
class of all ordinals $\text{\textbf{On}}$ into $\text{\textbf{No}}$,
with birthday function verifying $b(\phi(\beta))=\beta$. This also
proves that $\text{\textbf{No}}$ is a proper class and not a set
because $\phi$, being order preserving, is injective. For this reason
usually one says that $\text{\textbf{No}}$ is a Field, with the capital
initial letter to underline that it verifies the axioms of a field,
but its support set is a proper class.

We have already seen that in $\text{\textbf{No}}$ we have infinities
like $\omega=\{\N\,|\}$, but we can also easily construct infinitesimals
like\[
\eps:=\left\{ 0\,|\,\left\{ \frac{1}{n}\,:\, n\in\N_{>0}\right\} \right\} ,\]
which is strictly positive but smaller than any strictly positive
real.

All these examples can conduct us toward the conjecture that the class
$\text{\textbf{No}}$ is some kind of {}``universal'' field containing
every possible extension of the real field. Indeed we have the following
theorem (see \citet{Con2} for the proof; see also \citet{Ehrl} for
a more general and systematic treatment)
\begin{thm}
\label{thm:NoIsUniversal}The field $\text{\emph{\textbf{No}}}$ verifies
the following properties
\begin{enumerate}
\item \label{enu:NoIsAField}$\text{\emph{\textbf{No}}}$ is an ordered
Field
\item \label{enu:subfieldsAndNo}If \begin{align*}
A\text{ is an ordered subfield of \emph{\textbf{No}}}\\
A\text{ is an ordered subfield of }B,\end{align*}
with $|A|$ and $|B|$ sets and not proper classes, then there exist\begin{gather*}
B'\text{ ordered subfield of \emph{\textbf{No}}}\\
f:B\xrightarrow[\sim]{\phantom{\longrightarrow}}B'\text{ isormophism of ordered fields}\end{gather*}
such that $f|_{A}=1_{A}$.
\end{enumerate}
Moreover, if $F$ verifies these properties \ref{enu:NoIsAField}.
and \ref{enu:subfieldsAndNo} (like $\text{\emph{\textbf{No}}}$ does).
then $F\simeq\text{\emph{\textbf{No}}}$ as an ordered field.

\end{thm}
We can represent the situation in the following way: if we have (as
diagram of morphisms between ordered fields)\[
\xyR{45pt}\xyC{55pt}\xymatrix{A\ar@{_{(}->}[r]\ar@{^{(}->}[d] & \text{\textbf{No}}\\
B}
\]
then we can complete it with the commutative diagram\[
\xyR{45pt}\xyC{55pt}\xymatrix{A\ar@{_{(}->}[r]\ar@{^{(}->}[d] & \text{\textbf{No}}\\
B\ar[r]_{\sim}^{f} & B'\ar@{_{(}->}[u]}
\]
From this point of view the field of surreal numbers is remarkably
inclusive%
\footnote{From this point of view the name {}``surreal numbers'' is less meaningful
than the original Conway's {}``numbers'' without any adjective.%
}. For example applying the previous theorem with $A=\R$ and $B=\hyperR$,
we obtain that $\text{\textbf{No}}$ contains, up to isomorphism,
the hyperreals of NSA.

\subsubsection{Comparison with Fermat reals}

The first comparison between surreals and Fermat reals comes from
the previous Theorem \eqref{thm:NoIsUniversal} which cannot be applied
to the ring $\ER$. More trivially, the existence of non-zero nilpotent
infinitesimals is not compatible with field axioms.

Moreover, the construction of $\text{\textbf{No}}$ is deeply based
on order properties and produces a single numeric field and not a
category of extended spaces, including manifolds, like our $\ECInfty$.

The field $\text{\textbf{No}}$ has many remarkable properties, it
is a real closed field, there is the possibility to define exponential
and logarithm and even a notion of Riemann integral (see \citet{For}).
On the other hand, like any other non Archimedean ordered field, $\text{\textbf{No}}$
is totally disconnected, therefore we have examples of functions differentiable
on an interval with everywhere zero derivative which are not constant,
we do not have the uniqueness of the primitive of a continuous function
and we do not have uniqueness in the simplest initial value problem:
$y'(x)=0$, $y(0)=0$. This cannot be directly compared with our results
regarding the development of the calculus on the Fermat reals (see
Chapters \ref{cha:CalculusOnOpenDomains} and \ref{cha:CalculusOnInfinitesimalDomains})
because our results are applicable to smooth functions only and not
to a lower degree of differentiability.

From a methodological point of view, as we have already sketched above,
Conway's construction seems to be based on the search of a theory
with strong intuitive meaning, essentially due to Conway's simplicity
theorem \eqref{thm:ConwaySimplicityTheorem}. Formalization like \citet{Cue}
and \citet{All} permit to obtain a good dialectic between formal
theory and intuitive interpretation, which is also the leading design
of the present work.

\section{\label{sec:LeviCivitaField}Levi-Civita field}

The Levi-Civita field (from now on: LCF) originally appeared in \citet{Lev1}
and \citet{Lev2}, but it was subsequently rediscovered by \citet{Ost},
\citet{Ned}, \citet{Berz1} and \citet{Berz2} (to whom, together
with K. Shamseddine, we can attribute the modern development of the
topic). For an account of Levi-Civita's work see also \citet{Lau}.
For a detailed work in this topic and the proofs of the theorems we
will state in this section, see e.g. \citet{Sha}. Because of the
several analogies between our Fermat reals and the LCF, we will introduce
this topic with a certain level of detail.

To motivate the introduction of the LCF we start from the idea to
add to the real field $\R$ a single new invertible infinitesimal
number $d$ to which we want to apply all the ring operations but
also arbitrary roots. Hence, we would like to be able to form numbers
like $3+d+2d^{2}$, but also like $d^{-1}$, $d^{1/2}$, $d^{-1/2}$,
$2-d^{3}+4d^{3/5}-\frac{1}{2}d^{-2/3}$. This can be easily obtained
using purely algebraic methods, e.g. considering formal power series
of the form\begin{equation}
x=\sum_{q\in\mathbb{Q}}x_{q}\cdot d^{q}.\label{eq:formalPowerSeries}\end{equation}
There is no problem in defining the sum pointwise\[
x+y:=\sum_{q\in\mathbb{Q}}x_{q}d^{q}+\sum_{q\in\mathbb{Q}}y_{q}d^{q}:=\sum_{q\in\mathbb{Q}}\left(x_{q}+y_{q}\right)d^{q},\]
 but we can recognize a first limitation considering the product,
which is defined, as usual for formal power series, as\begin{equation}
x\cdot y:=\sum_{q\in\mathbb{Q}}\left(\sum_{r+s=q}x_{r}\cdot y_{s}\right)\cdot d^{q}.\label{eq:definitionOfProductFormalPowerSeries}\end{equation}
Indeed, the sum $\sum_{r+s=q}x_{r}\cdot y_{s}$ can have an infinite
number of addends, depending on how many non-zero coefficients $x_{r}$
and $y_{s}$ we have in the factors $x$ and $y$, and hence the previous
definition of product $x\cdot y$ can be meaningless for generic formal
power series $x$, $y$. Because we want that the definition \eqref{eq:definitionOfProductFormalPowerSeries}
works for every pair $x$, $y$, we must introduce a limitation on
the number of coefficients in our formal power series. In other words,
we must limit the number of non zero coefficients in the formal series.
For example we can have\[
\bar{x}=d^{-3}+d^{-2}+d^{-1}+1+d+d^{2}+d^{3}+d^{4}+d^{5}+\ldots,\]
or\[
\bar{y}=\ldots+d^{-5}+d^{-4}+d^{-3}+d^{-2}+d^{-1}+1+d+d^{2}.\]
More generally, the equation $r+s=q$ may have infinitely many solutions
if there is an accumulation point for the indexes $s\in\mathbb{Q}$
such that $y_{s}\ne0$. In the LCF the choice fall on power series
with a finite number of exponents {}``on the left'', i.e. such that\begin{equation}
\forall q\in\mathbb{Q}:\ \text{card}\left\{ r\in\mathbb{Q}\,|\, x_{r}\ne0\,,\, r\le q\right\} \text{ is finite}.\label{eq:left-finiteCondition}\end{equation}
From our point of view, which is not near to the formal point of view
expressed in the construction of the LCF - remember that at present
we do not have a notion of convergence for our series - it seems hard
to motivate this choice instead of the limitation {}``on the right''.
Moreover, let us note explicitly that if we do not want to introduce
limitations on the exponents we consider in our power series, then
we are forced to say that the product $x\cdot y$ is not always defined,
but only for those pairs $x$, $y$ such that the sum $\sum_{r+s=q}x_{r}\cdot y_{s}$
converges with respect to some notion of convergence%
\footnote{This notion may be trivial, e.g. if we consider only those pairs for
which for every $q\in\mathbb{Q}$, there is only a finite number of
exponents $r$, $s$ such that $r+s=q$.%
}. This may seem strange from an algebraic point of view, but it can
be considered more common in the calculus, where, e.g. in the standard
Schwartz's theory of distribution, the product of two distributions
is not always defined and can be considered in some cases only (see
e.g. \citet{Col2}), or where the set of convergent or divergent real
sequences is not closed with respect to pointwise product.

Leaving the intuitive motivations to arrive to a more formal mathematics,
we can introduce our formal power series thinking of the corresponding
definition for polynomials: identifying a polynomial $a_{0}+a_{1}x+a_{2}x^{2}+\ldots+a_{n}x^{n}$
with the $n$-tuple of its coefficients $(a_{0},a_{1},a_{2},\ldots,a_{n})$
is equivalent to say that a polynomial is a function of the form $a:\left\{ 0,1,2,\ldots,n\right\} \freccia\R$.
Analogously we can define
\begin{defn}
\label{def:LCF}The support set $\mathcal{R}$ of the LCF is the set
of all the functions $x:\mathbb{Q}\freccia\R$ with \emph{left-finite
support}, i.e. such that\[
\forall q\in\mathbb{Q}:\ \text{\emph{card}}\left\{ r\in\mathbb{Q}\,|\, x_{r}\ne0\ ,\ r\le q\right\} \text{ is finite}.\]

\end{defn}
So, the value%
\footnote{The notation with square brackets $x[q]$ permits to avoid confusion
when one consider functions defined on the LCF $\mathcal{R}$.%
} $x(q)=:x[q]$ has to be thought as the coefficient of the addend
$x[q]\cdot d^{q}$. The ring operations are defined for $q\in\mathbb{Q}$
as\begin{align*}
(x+y)[q]: & =x[q]+y[q]\\
(x\cdot y)[q]: & =\sum_{\substack{r,\, s\in\mathbb{Q}\\
r+s=q}
}x[r]\cdot y[s],\end{align*}
and we can verify that $(\mathcal{R},+,\cdot)$ becomes a field (see
e.g. \citet{Sha}). In a formal power series like \eqref{eq:formalPowerSeries}
the leading term $x_{m}\cdot d^{m}$ with the lower value $m\in\mathbb{Q}$
of the exponent determines the behavior of the number from several
points of view, e.g. with respect to order. For this reason we introduce
the following notations:
\begin{defn}
\label{def:orderAndOthers}Let $x$, $y\in\mathcal{R}$, then
\begin{enumerate}
\item $\text{\emph{supp}}(x):=\left\{ q\in\mathbb{Q}\,|\, x[q]\ne0\right\} $ 
\item $\lambda(x):=\min\left(\text{\emph{supp}}(x)\right)$ for $x\ne0$
and $\lambda(0):=+\infty$. The term $\lambda(x)$ is called \emph{order
of magnitude.}
\item $x\sim y\DIff\lambda(x)=\lambda(y)$. This relation is called \emph{agreement
of order of magnitude.}
\item $x\approx y\DIff\lambda(x)=\lambda(y)\e{and}x\left[\lambda(x)\right]=y\left[\lambda(y)\right]$
\item $x=_{r}y\DIff\forall q\in\mathbb{Q}_{\le r}:\ x[q]=y[q]$
\end{enumerate}
\end{defn}
If $r>0$, the relation $=_{r}$ is the analogous of our equality
up to $k$-th order infinitesimals (see Chapter \ref{cha:equalityUpTo_k-thOrder}):
if $x=_{r}y$, then $x-y$ is given by sum of infinitesimals $a_{q}d^{q}$
of order $d>r$. E.g. we can expect to have $1+d+d^{2}-2d^{3}=_{2}1+d+d^{2}+4d^{5}$
or $d^{n+1}=_{n}0$. So, it appears sufficiently clear that, even
if we do not have a ring with nilpotent elements, the equivalence
relation $=_{r}$ can supply a possible alternative language.

In the LCF we can prove the existence of roots:
\begin{thm}
\label{thm:rootsInLCF}Let $x\in\mathcal{R}$ be non zero and $n\in\N_{>0}$.
If $n$ is even and the leading term $x[\lambda(x)]>0$, then $x$
has two $n$-th roots in $\mathcal{R}$. If $n$ is even and $x[\lambda(x)]<0$,
then $x$ has no $n$-th roots in $\mathcal{R}$. If $n$ is odd,
then $x$ has a unique $n$-th root in $\mathcal{R}$.
\end{thm}
Exactly like in the study of polynomials, the particular number $d[q]:=1$
for $q=1$ and $d[q]:=0$ otherwise, works as the independent variable
in our formal power series, and the equality \eqref{eq:formalPowerSeries}
can now be proved for every $x\in\mathcal{R}$ because $\left(d^{r}\right)[q]=1$
if $q=r$ and $\left(d^{r}\right)[q]=0$ otherwise. Let us note explicitly
that if $r=\frac{p}{q}\in\mathbb{Q}$ with $p$, $q\in\mathbb{Z}$,
$q>0$, then $d^{r}=\sqrt[q]{d^{p}}$, so we need the previous Theorem
\ref{thm:rootsInLCF}.

Obviously, the embedding of the reals is given by $r\in\R\mapsto r[-]\in\mathcal{R}$,
where $r[0]=r$ and $r[q]=0$ otherwise, but it is now also clear
that formal Laurent series (and hence also D. Tall's superreal numbers,
see e.g. \citet{Tall}), i.e. numbers of the form\[
x=\sum_{k=-N}^{+\infty}x_{k}\cdot d^{k}\]
are embedded in the LCF.

Essential for the development of the LCF as an ordered field but also
for the different notions of continuity and differentiability of functions
$f:\mathcal{R}\freccia\mathcal{R}$ is the order relation. As hinted
above, we can define the order relation by comparison of the leading
terms
\begin{defn}
\label{def:orderRelationLCF}If $x$, $y\in\mathcal{R}$, we define\begin{align*}
x & >0\DIff x[\lambda(x)]>0\\
x & >y\DIff x-y>0\end{align*}

\end{defn}
With this relation the LCF becomes a totally ordered field extending
the real field.

In the framework of Fermat reals, the natural topology is the final
one with respect to which any figure is continuous (see Section \ref{sec:DefinitionOfCartesianClosure}).
As we mentioned ibidem, for Fermat reals, and more generally for Fermat
spaces, the topology is a byproduct of the diffeological structure.
Using this structure, we have a natural way to define smooth maps
between $\CInfty$-spaces and hence to extend these maps to the corresponding
Fermat spaces, without any particular focusing on the topology. In
the LCF, it is not so clear what functions $f:\R\freccia\R$ can be
extended to the whole $\mathcal{R}$ and hence the approach is different
and mimics the classical approach of calculus. The next step is hence
to use the order relation to define a corresponding order topology.
\begin{defn}
\label{def:orderTopologyLCF}Because the order relation on $\mathcal{R}$
is total, we can define the absolute value in the usual way\[
|x|:=\begin{cases}
x & \text{if }x\ge0\\
-x & \text{if }x<0\end{cases}\]
Moreover, we say that a subset $U\subseteq\mathcal{R}$ is \emph{open
in the order topology} iff\[
\forall u\in U\ \exists\,\delta\in\mathcal{R}_{>0}:\quad\left\{ x\in\mathcal{R}\,:\,|x-u|<\delta\right\} \subseteq U.\]

\end{defn}
For example the sequence $(d^{n})_{n\in\N}$ converges to the zero
sequence in the order topology. Using the same idea, that is the formal
analogy with the reals $\R$, we can define continuity, differentiability
and convergence of sequences
\begin{defn}
\label{def:topContLCF}Let $D\subseteq\mathcal{R}$ and $f:D\freccia\mathcal{R}$,
then we say that $f$ \emph{is topologically continuous at $x_{0}\in D$}
iff \[
\forall\eps\in\mathcal{R}_{>0}\ \exists\,\delta\in\mathcal{R}_{>0}:\ \forall x\in D:\ |x-x_{0}|<\delta\Rightarrow|f(x)-f(x_{0})|<\eps.\]

\end{defn}
\medskip{}

\begin{defn}
\label{def:topDiffLCF}Let $D\subseteq\mathcal{R}$ and $f:D\freccia\mathcal{R}$,
then we say that $f$ \emph{is topologically differentiable at $x_{0}\in D$}
iff there exists a number $l\in\mathcal{R}$ such that\[
\forall\eps\in\mathcal{R}_{>0}\ \exists\,\delta\in\mathcal{R}_{>0}\,\forall x\in D:\ 0<|x-x_{0}|<\delta\Rightarrow\left|\frac{f(x)-f(x_{0})}{x-x_{0}}-l\right|<\eps.\]

\end{defn}
\medskip{}

\begin{defn}
\label{def:convergenceOfsequencesLCF}Let $s:\N\freccia\mathcal{R}$,
then we say that $s$\emph{ converges strongly to $s\in\mathcal{R}$
iff\[
\forall\eps\in\mathcal{R}_{>0}\ \exists\, N\in\N:\ \forall n\in\N:\ n\ge N\Rightarrow|s_{n}-s|<\eps.\]
}
\end{defn}
It is interesting that now, using this notion of convergence, we can
associate to our formal power series a notion of convergence:
\begin{thm}
\label{thm:convergenceFormalPowerSeriesLCF}Let $x\in\mathcal{R}$
and define recursively\begin{align*}
q_{0}: & =\lambda(x)\end{align*}
\[
q_{k+1}:=\begin{cases}
\min\left(\text{\emph{supp}}(x)\setminus\left\{ q_{0},q_{1},\ldots,q_{k}\right\} \right) & \text{if }\text{\emph{supp}}(x)\supset\left\{ q_{0},q_{1},\ldots,q_{k}\right\} \\
0 & \text{otherwise}\end{cases}\]
then the sequence $n\in\N\mapsto\sum_{k=0}^{n}x[q_{k}]\cdot d^{q_{k}}\in\mathcal{R}$
converges strongly to $x$, so that we can write\[
x=\sum_{k=0}^{+\infty}x[q_{k}]\cdot d^{q_{k}}.\]

\end{thm}
Note that this theorem does not realize the above mentioned dialogue
between potential infinitesimals and actual infinitesimals because,
trivially, in its statement there is no mention of any such potential
infinitesimal, i.e. of a function of the form $i:E\freccia\R$ convergent
to zero, instead in the statement convergence is understood in the
LCF sense.

We can now give some motivations for the choice of the domain $\mathbb{Q}$
for the elements of the LCF $\mathcal{R}$: why is there, in the definition
of $\mathcal{R}$, the field $\mathbb{Q}$ instead of $\R$ like in
our Fermat reals%
\footnote{Recall e.g. that we can consider infinitesimal $\diff t_{a}$ for
every real number $a\ge1$.%
}? The answer can be anticipated saying that $\mathbb{Q}$ is the simplest
domain to obtain some of the desired properties. To render this statement
more precise we need the notion of skeleton group.
\begin{defn}
\label{def:skeletonGroupLCF}Let $(F,+,\cdot,<)$ be a totally ordered
field and define in it the absolute value as usual. Let $a$, $b\in F_{\ne0}$,
then we say\[
a\ll b\DIff\forall n\in\N:\ n\cdot a<b\]
and we will read it \emph{$a$ is infinitely smaller than $b$.} Moreover,
we will say\[
a\sim b\DIff\neg(|a|\ll|b|)\e{and}\neg(|b|\ll|a|).\]
The relation $\sim$ is an equivalence relation, and we will denote
by\[
S_{F}:=\left\{ [a]_{\sim}\,|\, a\in F_{\ne0}\right\} \]
the set of all its equivalence classes. Moreover, it is possible to
prove that the following definitions are correct:\begin{align*}
[a]_{\sim}\cdot[b]_{\sim}: & =[a\cdot b]_{\sim}\\
{}[a]_{\sim}^{-1}: & =[a^{-1}]_{\sim}\end{align*}
\[
[a]_{\sim}<[b]_{\sim}\DIff\forall n\in\N:\ n|a|<|b|.\]

\end{defn}
It is possible to prove that $(S_{F},\cdot,<)$ is a totally ordered
group, called the \emph{skeleton group of $F$}. This notion is naturally
tied with the notion of non-Archimedean field. Indeed, the skeleton
group of the real field is trivial $S_{\R}=\{[1]_{\sim}\}$, but it
is not so for non-Archimedean fields, as stated in the following
\begin{thm}
Let $F$ be a totally ordered non-Archimedean field, then\[
\mathbb{Z}\subseteq S_{F}.\]
Moreover, if $F$ admits roots of positive elements, then\[
\mathbb{Q}\subseteq S_{F}.\]

\end{thm}
This motivates why we take $\mathbb{Q}$ as domain of our functions
$x\in\mathcal{R}$, that is as exponents of $d$ in our formal power
series: it is the smallest set of exponents that permits to have a
non-Archimedean field and roots of positive elements.

But the idea to follow formal analogies to define continuity and differentiability
(see Definitions \ref{def:topContLCF} and \ref{def:topDiffLCF})
presents several problems (not characteristic of the LCF, but common
to every non-Archimedean totally ordered field): e.g. the function
$f:[0,1]\freccia\mathcal{R}$ defined by\[
f(x):=\begin{cases}
0 & \text{if }x\text{ is infinitely small}\\
1 & \text{if }x\text{ is finite}\end{cases}\]
is topologically continuous and topologically differentiable, but
it does not assume the value $d\in[0,1]$ even though $f(0)<d<f(1)$,
hence it does not verifies the intermediate value theorem. Moreover,
$f'(x)=0$ for all $x\in[0,1]$ but $f$ is not constant on $[0,1]$.
Therefore, also the simplest initial value problem $y'=0$, $y(0)=0$
does not have a unique solution. This is due to the fact that infinitesimals
are totally disconnected from finite numbers because $d\ll r$ for
every $r\in\R_{>0}$, and this is a general problem of non-Archimedean
fields. On the other hand, as we have seen in Chapter \ref{cha:Infinitesimal-differential-geometry},
in the context of Fermat reals, we do not encounter this type of problems.

The solution adopted in the framework of the LCF is to change the
notion of topological continuity introducing a Lipschitz condition:
\begin{defn}
\label{def:continuityLCF}Let $a<b$ be given in $\mathcal{R}$, let
$I\in\left\{ (a,b),[a,b],(a,b],[a,b]\right\} $ be an interval of
$\mathcal{R}$ and $f:I\freccia\mathcal{R}$, then we say that $f$\emph{
is continuous in $I$} iff\[
\exists\, M\in\mathcal{R}:\ \forall x\in I:\ \left|f(x)-f(y)\right|\le M\cdot\left|x-y\right|.\]

\end{defn}
This is very interesting for our comparison because a Lipschitz condition
is essential both for our Fermat reals (see e.g. the Definitions \ref{def:extensionOfFunctions}
and \ref{def:NilpotentFunctions}, but also \citet{Gio3} where the
Lipschitz condition is even more needed) and also for \citet{Fr-Kr,Kr-Mi}.
The corresponding definition of differentiability recalls our approach
with smooth incremental ratios (see Section \ref{sec:theFermatMethod}):
\begin{defn}
\label{def:diffLCF}Under the same hypothesis as in the previous definition
we say that $f$\emph{ is differentiable on $I$} iff there exists
a continuous function\[
f'[-,-]:I\times I\freccia\mathcal{R}\]
such that\[
\forall x,y\in I:\ f(y)=f(x)+f'[x,y]\cdot(y-x).\]
As usual we will set $f'(x):=f'[x,x]$.
\end{defn}
This definition is also very interesting, because, even if the approach
of the LCF is a formal approach to infinitesimals, some of the chosen
solutions are very similar to those adopted in non formal approaches
to infinitesimals. For example, note the requirement that $f'[-,-]$
should be continuous, and hence it should verify a Lipschitz condition,
in analogy to the Lipschitz differentiability introduced in \citet{Fr-Kr}.

A natural problem dealing with LCF is what class of functions $f:\R\freccia\R$
can be extended to a meaningful subclass of $\mathcal{R}$ with the
possibility to generalize to them some properties, like the intermediate
value theorem, an inverse function theorem, the maximum theorem, the
mean value theorem, Rolle's theorem, the existence of primitive functions,
or the constancy principle. Because of the left-finiteness of our
formal power series $x\in\mathcal{R}$, the most natural class of
functions are those locally expandable in power series (the term \emph{analytic}
is used for real functions only in the context of the LCF). The solution
adopted in \citet{Berz2} and \citet{Sha} (see also \citet{Sh-Be}
and references therein for a more recent article) is to define a notion
of convergence of power series with coefficients in $\mathcal{R}$,
to prove for them the above mentioned theorems, and hence to show
that standard power series in $\R$ are included as a particular case
of this notion of convergence in $\mathcal{R}$. It is also interesting
to note that this concept of convergence is not the one derived from
the formal analogies with the real case (see e.g. the Definition \ref{def:topContLCF})
but it is rather derived from a family of seminorms. For more details
on this development, see the above mentioned references.

The left-finiteness of the Levi-Civita numbers permits to represent
them on a computer. Indeed, for every $r\in\mathbb{Q}$ the amount
of information we have to store in the power series \eqref{eq:formalPowerSeries},
up to the terms $x_{q}d^{q}$ with $q\le r$, is finite and we can
represent all these numbers with the usual precision available in
a computer. Therefore, using the equivalence relation $=_{r}$ (see
Definition \ref{def:orderAndOthers}) we can implement a calculus
of Levi-Civita numbers on a computer. Obviously, this is possible
for rare cases only in other theories of infinitesimals%
\footnote{We only mentioned here that the surreal numbers has been implemented
in the computer based proof assistant Coq, see \citet{Mam}.%
} and it founds useful applications in automatic differentiation (see
e.g. \citet{Berz1}). Here the problem is to find a computer algorithm
to calculate the derivatives, up to a given order, of computationally
complicated functions like (see \citet{Sha})\[
f(x)=\frac{\sin(x^{3}+2x+1)+\frac{3+\cos(\sin(\ln|1+x|))}{\exp\left(\tanh\left(\sinh\left(\cosh\left(\frac{\sin(\cos(\tan(\exp(x))))}{\cos(\sin(\exp(\tan(x+2))))}\right)\right)\right)\right)}}{2+\sin(\sinh(\cos(\tan^{-1}(\ln(\exp(x)+x^{2}+3)))))},\]
obtained by composition of elementary functions like $\sin$, $\cos$,
$\exp$, ..., the Heaviside function\[
H(x):=\begin{cases}
1 & \text{if }x\ge0\\
0 & \text{otherwise}\end{cases}\]
and of the field operations. These are called \emph{computer functions},
and can be extended to a suitable subset of $\mathcal{R}$ using their
expansion in power series. The property that permits to compute these
derivatives is expressed in the following theorem and it presents
strong analogies with the calculus in our Fermat reals:
\begin{thm}
\label{thm:computingDerivativesInLCF}Let $f$ be a computer function
continuous at $x_{0}\in\R$ and extendable to $x_{0}\pm d$. Then
$f$ is $n$-times differentiable at $x_{0}$ if and only if there
exist real numbers $\alpha_{1},\ldots,\alpha_{n}$ such that\begin{align*}
f(x_{0}+d)=_{n}\sum_{j=0}^{n}\frac{\alpha_{j}}{j!}\cdot d^{j}\end{align*}
\[
f(x_{0}-d)=_{n}\sum_{j=0}^{n}(-1)^{j}\frac{\alpha_{j}}{j!}\cdot d^{j}.\]
Moreover, in this case we have $f^{(j)}(x_{0})=\alpha_{j}$ for $j=0,\ldots,n$.
\end{thm}
In \citet{Sha} one can find several examples of computation of derivatives
using these formulas, and of non smooth functions whose regularity
is proved using this theorem. A software, called COSY INFINITY, has
also been created, which is suitable for the computation of derivatives
of functions using the LCF (see \citet{Berz-et-al}, \citet{Sha}
and references therein also for the comparison with other methods
of computation of derivatives).

\subsubsection{Comparison with Fermat reals}

We have tried to introduce the LCF with a certain detail, due to the
many analogies that one can see between the LCF theory and our own
theory, even if the two approaches are very different from the philosophical
point of view. We can underline several points where the LCF can be
considered a better framework with respect to Fermat reals, and several
others where we can state the opposite:

\let\savelabelitemi=\labelitemi

\renewcommand\labelitemi[0]{\smiley}
\begin{itemize}
\item The LCF is a field and not a ring. This is surely reassuring for some
readers, even if we have mentioned in connection with some deep problems
that the theory of non-Archimedean fields find in the development
of the calculus. On the other hand, the availability of infinities
can be very useful to express e.g. the Riemann integral as an infinite
sum or to define Dirac delta functions as ordinary functions, like
in NSA.
\item In the LCF we have the existence of roots for every positive numbers,
in particular also for infinitesimals. Of course, this is incompatible
with a ring containing nilpotent elements and indeed, using Fermat
reals, we are able to define roots only for invertible numbers. We
hinted in Section \ref{cha:FurtherDevelopments} to the possibility
to define the square root of an infinitesimal Fermat real $h\in D_{\infty}$
as the simplest $k\in D_{2\cdot\omega(h)}$ such that $k^{2}=h$,
but this notion, even if useful, does \emph{not} verify the usual
equality $\sqrt{k^{2}}=|k|$.
\item In the present work, we chose to develop a theory of smooth functions
only, so as to obtain the simplest results useful for smooth differential
geometry. There is the possibility to extend some of our results to
functions which only belong to $\Cn$, keeping present some of the
ideas used in \citet{Gio3}. An example in this direction is given
by Theorem \ref{thm:TaylorForCnFunctions}. But at present, the theory
of Fermat reals and Fermat spaces is not developed in this direction.
The possibility to define continuity and differentiability in the
theory of the LCF is hence interesting. Because the theory of LCF
is not a theory of smooth functions only, we have the possibility
to prove a useful theorem like \ref{thm:computingDerivativesInLCF},
even if that theorem is applicable only to computer functions expandable
to $x_{0}\pm d$, i.e. to a class smaller than the one considered
in the previous analogous Theorem \ref{thm:TaylorForCnFunctions}.
\end{itemize}
\renewcommand\labelitemi[0]{\savelabelitemi}

On the other hand we have:

\noindent \let\savelabelitemi=\labelitemi

\noindent \renewcommand\labelitemi[0]{\frownie}
\begin{itemize}
\item Until now, the theory of LCF permits to extend the real field only,
and not the general case of smooth manifolds, like in the case of
Fermat reals.
\item The calculus with nilpotent infinitesimals seems easier, for smooth
functions, with respect to the use of the equivalence relation $=_{r}$.
As we mentioned above, because on the right the power series of $\mathcal{R}$
are not necessarily finite, the functions that naturally extends from
$\R$ to the LCF are the analytical ones. So, we have the methodological
contradiction that the LCF permits to develop a meaningful notion
of continuity and differentiability, but at the same time, because
of the form of the formal power series considered in the LCF, the
best results are for functions locally expandable in power series
and not for a lower degree of regularity. At the same time each Fermat
number need only a finite number of reals to be stored in a computer,
without any need of the equivalence relation $=_{r}$, and hence $\ER$
can also be implemented in a computer.
\item Because of the formal approach to infinitesimals, the intuitive meaning
of $\mathcal{R}$ as connected to potential infinitesimals of $\R$
is missing. For example, at the best of our knowledge, there is no
idea about how it would be possible to extend a given function $f:\R\freccia\R$
to infinitely large numbers in $\mathcal{R}$.
\end{itemize}
\renewcommand\labelitemi[0]{\savelabelitemi}

\end{appendix}

\bibliographystyle{plainnat} \clearpage{}\addcontentsline{toc}{chapter}{\bibname}\bibliographystyle{plainnat}
\bibliography{mybib}

\end{document}